\documentclass[12pt]{amsbook}
\usepackage{amssymb}
\usepackage[all]{xy}

\usepackage[colorlinks=true,linkcolor=magenta,citecolor=magenta]{hyperref}
\makeindex

\newtheorem{theorem}{Theorem}[chapter]
\newtheorem{answer}[theorem]{Answer}
\newtheorem{cat}[theorem]{Cat}
\newtheorem{challenge}[theorem]{Challenge}
\newtheorem{claim}[theorem]{Claim}
\newtheorem{comment}[theorem]{Comment}
\newtheorem{conclusion}[theorem]{Conclusion}

\newtheorem{definition}[theorem]{Definition}
\newtheorem{elephant}[theorem]{Elephant}
\newtheorem{example}[theorem]{Example}
\newtheorem{exercise}[theorem]{Exercise}
\newtheorem{fact}[theorem]{Fact}
\newtheorem{goat}[theorem]{Goat}
\newtheorem{lemma}[theorem]{Lemma}
\newtheorem{liste}[theorem]{List}
\newtheorem{magpie}[theorem]{Magpie}
\newtheorem{plan}[theorem]{Plan}
\newtheorem{poetry}[theorem]{Poetry}
\newtheorem{principle}[theorem]{Principle}

\newtheorem{proposal}[theorem]{Proposal}
\newtheorem{proposition}[theorem]{Proposition}
\newtheorem{question}[theorem]{Question}

\newtheorem{rat}[theorem]{Rat}
\newtheorem{remark}[theorem]{Remark}

\newtheorem{scenario}[theorem]{Scenario}
\newtheorem{speculation}[theorem]{Speculation}
\newtheorem{tiger}[theorem]{Tiger}

\begin{document}

\title{Normal random variables}

\author{Teo Banica}
\address{Department of Mathematics, University of Cergy-Pontoise, F-95000 Cergy-Pontoise, France. {\tt teo.banica@gmail.com}}

\subjclass[2010]{60C05}
\keywords{Central limit, Normal variable}

\begin{abstract}
This is an advanced introduction to the various types of normal random variables, with most of the needed preliminaries included. We first discuss the probability basics, standard central limits, and the theory of the usual, real normal variables, with examples, illustrations and numerous formulae. Then we go on a similar discussion regarding the complex normal variables, and the Rayleigh variables too. We then move to arbitrary dimensions, with a discussion regarding the Gaussian vectors, and related probability laws, featuring some functional analysis, and geometry and physics too. Finally, we provide an introduction to the quantum versions of the normal variables, and notably to those coming from free probability and random matrices.
\end{abstract}

\maketitle

\chapter*{Preface}

What is a normal variable? Good question, depending on the type of measurements that you make, and in answer, we can have here real normal variables, then complex normal variables, with these being discrete or continuous, then general vector normal variables, of various types, and even some quantum versions of these. In short, many things to be learned, and here are a few things that can be said, about this:

\bigskip

(1) The simplest random variables are those discrete, of type $f:X\to\mathbb N$. In this setting, the simplest variables are the Bernoulli and geometric ones, and by suitably summing independent such variables, we are led to the Poisson law, and its versions:
$$p_1=\frac{1}{e}\sum_{k=0}^\infty\frac{\delta_k}{k!}$$

(2) In the continuous case, there are several possible choices for the simplest variables, and in view of this, the best is to just say that normality means that $f:X\to\mathbb R$ appears as a suitably rescaled sum of centered, independent, identically distributed variables. This leads to the normal law, which is as follows, and its numerous versions:
$$g_1=\frac{1}{\sqrt{2\pi}}\,e^{-x^2/2}dx$$

To be more precise, among the versions of $g_1$ we have various shifts and rescalings, then complex versions, then higher dimensional versions, then lenghts and squared lengths of these vector versions, and finally interpolations too, using the gamma function.

\bigskip

(3) And with this, end of the story? Not really, because in certain contexts, such as the random matrix one, the independence used above must be replaced by freeness. And when doing so, the normal law gets replaced by the Wigner semicircle law:
$$\gamma_1=\frac{1}{2\pi}\sqrt{4-x^2}\,dx$$

(4) As for the Poisson law, its free version turns to have a continuous density, and is the Marchenko-Pastur law, which is given by the following formula:
$$\pi_1=\frac{1}{2\pi}\sqrt{4x^{-1}-1}\,dx$$

Summarizing, many things to be learned, of both basic and advanced type. 

\bigskip

This book is an introduction to this, normal random variables of all types, with most of the needed preliminaries included. The book is organized in 4 parts, as follows:

\bigskip

I - We first discuss the probability basics, standard central limits, and the theory of the usual, real normal variables, with examples, illustrations and many formulae. 

\bigskip

II - Then we go on a similar discussion regarding the complex normal variables, and some related variables as well, such as the Rayleigh variables and their squares.

\bigskip

III - We then move to arbitrary dimensions, with a discussion of the Gaussian vectors, and related laws, featuring some functional analysis, and geometry and physics too. 

\bigskip

IV - Finally, we discuss the quantum versions of the central limits and normal variables, and notably those coming from free probability and random matrices.

\bigskip

In the hope that you will find this book useful. The material in Parts I, II, III is of course very standard, but we have tried to be as complete as possible, which each law, once accepted to the ``normal galaxy'', being carefully investigated, with as many formulae about it as we could find, ranging from very simple to fairly complicated. As for Part IV, free probability basics, this is again standard material, but I had to make some work and improvements here and there, in relation with some missing laws and formulae.

\bigskip

This book, although quite advanced, is readable as such. For some help with the basics you have my book on discrete variables \cite{ba1}, which is elementary. For a more standard introduction to probability you have my book on continuous variables \cite{ba2}. But if you're not afraid of natural science, and some missing proofs, go directly with this one.

\bigskip

Finally, at the technical level, we will use the moment method, for pretty much everything. And many thanks here to my cats, for some help with the computations.

\bigskip

\

{\em Cergy, September 2026}

\smallskip

{\em Teo Banica}

\baselineskip=15.95pt
\tableofcontents
\baselineskip=14pt

\part{Normal variables}

\ \vskip50mm

\begin{center}
{\em Perdono, perdono, perdono

Io soffro piu ancora di te

Perdono, perdono, perdono

Il male l'ho fatto piu a me}
\end{center}

\chapter{Random variables}

\section*{1a. Coins and dice}

Welcome to probability, which will be of discrete type, to start with. Inspired by what happens with simple games, from the real life, let us have as starting point:

\begin{definition}
A discrete probability space is a set $X$, usually finite or countable, whose elements $x\in X$ are called events, together with a function
$$P:X\to[0,\infty)$$
called probability function, which is subject to the condition
$$\sum_{x\in X}P(x)=1$$
telling us that the overall probability for something to happen is $1$.
\end{definition}

Which sounds quite neat, this is definitely a rigorous mathematical definition, that we can build our theory upon. Here are a few comments, in relation with this:

\bigskip

(1) As a first comment, our condition $\sum_{x\in X}P(x)=1$ perfectly makes sense, and this even if $X$ is uncountable, because the sum of positive numbers is always defined, as a number in $[0,\infty]$, and this no matter how many positive numbers we have.

\bigskip

(2) As a second comment, we have chosen in the above not to assume that $X$ is finite or countable, and this for instance because we want to be able to regard any probability function on $\mathbb N$ as a probability function on $\mathbb R$, by setting $P(x)=0$ for $x\notin\mathbb N$. 

\bigskip

(3) As a third comment, standing as a complement to Definition 1.1, once we have a probability function $P:X\to[0,\infty)$, given a  set of events $Y\subset X$, we can say that the probability for the events of type $y\in Y$ to happen is $P(Y)=\sum_{y\in Y}P(y)$.

\bigskip

As our first result now, making the link with our intuition from real life, we have:

\begin{proposition}
In the simplest case, that where $X$ is finite and $P$ is constant,
$$P(Y)=\frac{{\rm number\ of\ times}\ Y\ {\rm happens}}{{\rm total\ number\ of\ possibilities}}$$
for any set of events $Y\subset X$.
\end{proposition}

\begin{proof}
Assume indeed that our probability function $P:X\to[0,\infty)$ is constant. But then, $\sum_{x\in X}P(x)=1$ implies that $X$ must be finite, and that we must have:
$$P(x)=\frac{1}{|X|}\quad,\quad\forall x\in X$$

Now by summing this over the elements $y\in Y$ of a subset $Y\subset X$, we obtain:
$$P(Y)=\frac{|Y|}{|X|}\quad,\quad\forall Y\subset X$$

We are therefore led to the conclusion in the statement.
\end{proof}

Time now for some examples? The basic examples will come of course from basic games, and at the beginning of everything, we certainly have flipping coins:

\begin{example}
Flipping coins.
\end{example}

Here things are simple and clear, because when you flip a coin the corresponding discrete probability space, together with its probability measure, is as follows:
$$X=\big\{{\rm heads},\,{\rm tails}\big\}\quad,\quad P({\rm heads})=P({\rm tails})=\frac{1}{2}$$

In the case where the coin is biased, as to land on heads with probability $2/3$, and on tails with probability $1/3$, the corresponding probability space is as follows:
$$X=\big\{{\rm heads},\,{\rm tails}\big\}\quad,\quad P({\rm heads})=\frac{2}{3}\quad,\quad P({\rm tails})=\frac{1}{3}$$

More generally, given any number $p\in[0,1]$, we have an abstract probability space as follows, where we have replaced heads and tails by win and lose:
$$X=\big\{{\rm win},\,{\rm lose}\big\}\quad,\quad P({\rm win})=p\quad,\quad P({\rm lose})=1-p$$

Things become even more interesting when flipping a coin, biased or not, several times in a row. In the case of a biased coin, landing on heads with probability $p_1$ and on tails with probability $p_2=1-p_2$, thrown $n$ times, the probability space is:
$$X=\big\{1,2\}^n\quad,\quad P(i_1\ldots i_n)=p_{i_1}\ldots p_{i_n}\quad,\quad p_i\geq0\ ,\ p_1+p_2=1$$

Finally, in relation with this latter example, let us check that the sum 1 condition in Definition 1.1 is indeed satisfied. But this is indeed the case, as shown by:
$$\sum_{i\in X}P(i)
=\sum_{i_1,\ldots,i_n}p_{i_1}\ldots p_{i_n}
=\sum_{i_1}p_{i_1}\ldots\sum_{i_n}p_{i_n}
=1$$

Summarizing, coins reasonably understood, or at least the mathematical formalism of the coin game, reasonably understood. However, in what regards the mathematics itself of the coin game, many interesting things can be said. More on this later.

\bigskip

Moving on, at a more advanced level we can talk about rolling dice, as follows:

\begin{example}
Rolling dice.
\end{example}

Again, things here are simple and clear, because when you throw a die the corresponding probability space, together with its probability measure, is as follows:
$$X=\big\{1,\ldots,6\big\}\quad,\quad P(i)=\frac{1}{6}\ ,\ \forall i$$

As before with the coins, we can further complicate this by assuming that the die is biased, say landing on face $i$ with probability $p_i\in[0,1]$. In this case the corresponding probability space, together with its probability measure, is as follows:
$$X=\big\{1,\ldots,6\big\}\quad,\quad P(i)=p_i\quad,\quad p_i\geq0\ ,\ \sum_ip_i=1$$

Also as before with coins, things become more interesting when throwing a die several times in a row, or equivalently, when throwing several identical dice at the same time. In this latter case, with $n$ identically biased dice, the probability space is as follows:
$$X=\big\{1,\ldots,6\big\}^n\quad,\quad P(i_1\ldots i_n)=p_{i_1}\ldots p_{i_n}\quad,\quad p_i\geq0\ ,\ \sum_ip_i=1$$

Observe that the sum 1 condition in Definition 1.1 is indeed satisfied, and with this proving that our dice modeling is bug-free, due to the following computation:
$$\sum_{i\in X}P(i)
=\sum_{i_1,\ldots,i_n}p_{i_1}\ldots p_{i_n}
=\sum_{i_1}p_{i_1}\ldots\sum_{i_n}p_{i_n}
=1$$

In fact, this computation is identical to the one that we did before, for coins.

\bigskip

Getting back now to theory, in the general context of Definition 1.1, we can see that what we have there is very close to the biased die, from Example 1.4. Indeed, in the general context of Definition 1.1, we can say that what happens is that we have a die with $|X|$ faces, which is biased such that it lands on face $i$ with probability $P(i)$. 

\bigskip

Which is something quite interesting, theoretically speaking, in relation with the discrete probability theory that we want to develop in this chapter, because we have now some intuition on all this. So, as a conclusion, let us record this finding as follows:

\begin{principle}
Discrete probability can be understood as being about throwing a general die, having an arbitrary number of faces, and which is arbitrarily biased too:
$$\xymatrix@R=7pt@C=7pt{
\ar@{-}[rr]\ar@{-}[dd]&&\ar@{-}[dd]&&\ar@{-}[rr]\ar@{-}[dd]&&\ar@{-}[dd]&&\ar@{-}[rr]\ar@{-}[dd]&&\ar@{-}[dd]&&\ar@{-}[rr]\ar@{-}[dd]&&\ar@{-}[dd]\\
&p_x&&&&p_y&&&&p_z&&&&p_t&&&&\ldots\\
\ar@{-}[rr]&&&&\ar@{-}[rr]&&&&\ar@{-}[rr]&&&&\ar@{-}[rr]&&}$$
To be more precise, when writing our probability space as $X=\{x,y,z,t,\ldots\}$, with the corresponding probabilities being $p_x,p_y,p_z,p_t,\ldots\,$, the picture is the one above.
\end{principle}

Moving on, in order to do some mathematics, in the context of Definition 1.1, we need a random variable $f:X\to\mathbb R$, and the mathematics will consist in computing the expectation of this variable, $E(f)\in\mathbb R$. Let us axiomatize this situation as follows:

\index{random variable}
\index{average gain}
\index{expectation}

\begin{definition}
A random variable on a probability space $X$ is a function
$$f:X\to\mathbb R$$
and the expectation of such a random variable is the quantity
$$E(f)=\sum_{x\in X}f(x)P(x)$$
which is best thought as being the average gain, when the game is played.
\end{definition}

Here the word ``game'' refers to the probability space interpretation from Principle 1.5. Indeed, in that context, with our discrete set of events $X$ being thought of as corresponding to a generalized die, and by thinking of $f:X\to\mathbb R$ as representing some sort of money, the above quantity $E(f)$ is what we win, on average, when playing the game.

\bigskip

Getting back now to the examples, as our first result, we have:

\begin{theorem}
When flipping a usual coin $n$ times in a row:
\begin{enumerate}
\item $E({\rm heads})=n/2$.

\item $E({\rm tails})=n/2$.

\item $E(a\times{\rm heads}+b\times{\rm tails})=n(a+b)/2$.
\end{enumerate} 
\end{theorem}

\begin{proof}
Let us first fix some mathematical notation, in the spirit of what we have in Definition 1.6. Consider the variables in the statement, namely:
$$f:\big\{{\rm heads},\,{\rm tails}\big\}^n\to\mathbb R\quad,\quad f={\rm number\ of\ heads}$$ 
$$g:\big\{{\rm heads},\,{\rm tails}\big\}^n\to\mathbb R\quad,\quad g={\rm number\ of\ tails}$$ 
$$h:\big\{{\rm heads},\,{\rm tails}\big\}^n\to\mathbb R\quad,\quad h=a\times{\rm heads}+b\times{\rm tails}$$ 

With a bit of probabilistic know-how, we can compute the expectations of these variables, by using several methods, which are all instructive, as follows:

\medskip

(1) \underline{Symmetry proof}. Observe that, by using the obvious symmetry between heads and tails, at the level of the corresponding expectations, we have:
$$E(f)=E(g)$$

On the other hand, we also have $f+g=n$, which gives the following formula:
$$E(f)+E(g)=n$$

Thus $E(f)=E(g)=n/2$, and then $E(h)=n(a+b)/2$ follows by linearity.

\medskip

(2) \underline{Additive proof}. When flipping our coin several times in a row, the expectations will obviously sum up. Thus, if we denote by $f_n,g_n,h_n$ our variables, we have:
$$E(f_n)=nE(f_1)\quad,\quad E(g_n)=nE(g_1)\quad,\quad E(h_n)=nE(h_1)$$

But together with $E(f_1)=E(g_1)=1/2$ and $E(h_1)=(a+b)/2$, we are done.

\medskip

(3) \underline{Algebraic proof}. Here is as well a third proof, by brainless computation:
\begin{eqnarray*}
E(f)
&=&\sum_{x\in\{{\rm heads,\, tails}\}^n}\#({\rm heads}\in x)\times\frac{1}{2^n}\\
&=&\sum_{s=1}^n\binom{n}{s}\times s\times\frac{1}{2^n}\\
&=&\frac{n}{2^n}\sum_{t=0}^{n-1}\binom{n-1}{t}\\
&=&\frac{n}{2}
\end{eqnarray*}

Similarly, we have $E(g)=n/2$, and then $E(h)=n(a+b)/2$ comes by linearity.
\end{proof}

Getting now to the more general case of a biased coin, we have here:

\index{biased coin}

\begin{theorem}
When flipping a biased coin $n$ times in a row:
\begin{enumerate}
\item $E({\rm heads})=np$.

\item $E({\rm tails})=n(1-p)$.

\item $E(a\times{\rm heads}+b\times{\rm tails})=n(ap+b(1-p))$.
\end{enumerate} 
\end{theorem}

\begin{proof}
As before with the usual coins, we can denote by $f,g,h$ our random variables, and we have 3 possible approaches here, the situation with them being as follows:

\medskip

(1) \underline{Symmetry proof}. Normally, the idea here would be to exploit the symmetry between heads and tails, with the goal of reaching to the following formula:
$$(1-p)E(f)=pE(g)$$

Indeed, together with $E(f)+E(g)=n$, coming from $f+g=n$, this would finish the proof. However, the symmetry study is not obvious at $p\neq1/2$. So, wrong way. 

\medskip

(2) \underline{Additive proof}. As before, when flipping our coin several times in a row, the expectations will sum up. Thus, if we denote by $f_n,g_n,h_n$ our variables, we have:
$$E(f_n)=nE(f_1)\quad,\quad E(g_n)=nE(g_1)\quad,\quad E(h_n)=nE(h_1)$$

But together with $E(f_1)=p$, $E(g_1)=1-p$, $E(h_1)=ap+b(1-p)$, we are done.

\medskip

(3) \underline{Algebraic proof}. Here is as well the algebraic proof, your old man's favorite:
\begin{eqnarray*}
E(f)
&=&\sum_{s=1}^n\binom{n}{s}sp^s(1-p)^{n-s}\\
&=&np\sum_{t=0}^{n-1}\binom{n-1}{t}p^t(1-p)^{n-t-1}\\
&=&np
\end{eqnarray*}

Similarly, $E(g)=n(1-p)$, and then $E(h)=n(ap+b(1-p))$ comes by linearity.
\end{proof}

Next, come the dice. Here there are several possible natural choices for our random variable, and regarding the essentials, in the unbiased die case, we have:

\begin{theorem}
When rolling a usual die $n$ times in a row:
\begin{enumerate}
\item For $f_i=\#i$, we have $E(f_i)=n/6$.

\item For $g_i=i\#i$, we have $E(g_i)=ni/6$.

\item For $h={\rm sum}$, we have $E(h)=3.5n$.
\end{enumerate} 
\end{theorem}

\begin{proof}
With respect to our study for the usual coins, the situation is as follows:

\medskip

(1) \underline{Symmetry proof}. This obviously works, exactly as for unbiased coins.

\medskip

(2) \underline{Additive proof}. Again this works, obviously, as for the unbiased coins.

\medskip

(3) \underline{Algebraic proof}. Things are a bit more tricky here, depending on the exact variables that we have in mind, the situation being as follows:

\medskip

-- Normally, we can compute indeed $E(f_i)$ algebrically, a bit as before for the coins, by using binomial coefficients and some summing, and I will leave this as an exercise for you, and then pass to $E(g_i)$, and then to $E(h)$, by linearity.

\medskip

-- However, if we want to compute for instance $E(h)$ directly, algebrically, we are a bit in trouble, because, as you can convince yourself by doing some computations, finding the formula of the distribution of the sum $h$ is a non-trivial question.
\end{proof}

Finally, in the case of a biased die, the result is as follows:

\index{biased die}

\begin{theorem}
When rolling a biased die $n$ times in a row:
\begin{enumerate}
\item For $f_i=\#i$, we have $E(f_i)=np_i$.

\item For $g_i=i\#i$, we have $E(g_i)=nip_i$.

\item For $h={\rm sum}$, we have $E(h)=n\sum_iip_i$.
\end{enumerate} 
\end{theorem}

\begin{proof}
This is something hybrid between Theorems 1.8 and 1.9:

\medskip

(1) \underline{Symmetry proof}. This won't quite work for a biased die, for the technical reasons explained in the proof of Theorem 1.8, dealing with a biased coin.

\medskip

(2) \underline{Additive proof}. This works, as usual, because if we denote by $f_i^n,g_i^n,h^n$ our variables, it is clear that we have the following formulae:
$$E(f_i^n)=nE(f_i^1)\quad,\quad E(g_i^n)=nE(g_i^1)\quad,\quad E(h^n)=nE(h^1)$$

Thus, together with $E(f_i^1)=p_i$, $E(g_i^1)=ip_i$, $E(h^1)=\sum_iip_i$, we are done.

\medskip

(3) \underline{Algebraic proof}. This won't work well for a biased die, for the technical reasons explained in the proof of Theorem 1.9, dealing with a usual die.
\end{proof}

And we will end our study of coins and dice with this. So, was this a success or not? Well, we certainly managed to compute all our expectations, but what was quite frustrating was that our arsenal of methods was shrinking on the way. No good.

\bigskip

So, what to do? Leaving the debugging of symmetry methods to our pure mathematics colleagues, let us formulate the following challenge, for us probabilists:

\begin{challenge}
What is the distribution of the sum, when rolling a die $n$ times?
\end{challenge}

And with this being now on our to-do list, we will be back to it, as soon as possible. But of course, feel free to do some computations, this is how math is best learned.

\section*{1b. General theory}

We have already seen some good illustrations for Definition 1.6, so time now to get into more delicate aspects. Imagine that you want to set up some sort of business, with your variable $f:X\to\mathbb R$. You are of course mostly interested in the expectation $E(f)\in\mathbb R$, but passed that, the way this expectation comes in matters too. For instance:

\medskip

-- When your variable is constant, $f=c$, you certainly have $E(f)=c$, and your business will run smoothly, with not so many surprises on the way.

\medskip

-- On the opposite, for a complicated variable satisfying $E(f)=c$, your business will be more bumpy, with wins and loses on the way, depending on your skills.

\medskip

In short, and extrapolating now from business to mathematics, physics, chemistry and everything else, we must complement Definition 1.6 with something finer, regarding the ``quality'' of the expectation $E(f)\in\mathbb R$ appearing there. And the first thought here, which is the correct one, goes to the following definition, of the variance of our variable:

\index{variance}

\begin{definition}
The variance of a variable $f:X\to\mathbb R$ is the quantity
\begin{eqnarray*}
V(f)
&=&E\left([f-E(f)]^2\right)\\
&=&E(f^2)-E(f)^2
\end{eqnarray*}
intuitively measuring how far is $f$ from a constant variable.
\end{definition}

As a first observation, which is in tune with what we were saying in the above, the variance is 0 precisely when the variable is constant, equal to its expectation:
\begin{eqnarray*}
V(f)=0
&\iff&f={\rm constant}\\
&\iff&f=E(f)
\end{eqnarray*}

As a complement now to Definition 1.12, or rather as an alternative to it, which can bring some intuition too, with this being a matter of taste, let us formulate as well:

\index{standard deviation}

\begin{definition}
The standard deviation of a variable $f:X\to\mathbb R$ is the quantity
$$\sigma(f)=\sqrt{V(f)}$$
intuitively measuring as well how far is $f$ from a constant variable.
\end{definition}

Moving on, but still regarding expectations and variances of variables, as constructed above, at the abstract level, following Markov, we have the following result:

\index{Markov inequality}

\begin{proposition}
We have the following Markov inequality,
$$P\big(|f|\geq a\big)\leq\frac{E(f)}{a}$$
valid for any random variable $f:X\to\mathbb R$, and any $a>0$.
\end{proposition}

\begin{proof}
This is something trivial, coming from definitions. Indeed, we have:
$$E(f)
\geq\sum_{|f(x)|\geq a}aP(x)
=aP\big(|f|\geq a\big)$$

Thus, we are led to the Markov inequality in the statement.
\end{proof}

Next, following Chebycheff, we have the following key estimate:

\index{Chebycheff inequality}

\begin{theorem}
We have the following Chebycheff inequality,
$$P\big(|g-E|\geq b\big)=\frac{V}{b^2}$$
valid for any random variable $g:X\to\mathbb R$, having mean $E$ and variance $V$.
\end{theorem}

\begin{proof}
Given a random variable $g:X\to\mathbb R$, having mean $E$ and variance $V$, we can use the Markov inequality above with $f=(g-E)^2$, and we obtain in this way:
$$P\big(|g-E|^2\geq b^2\big)
\leq\frac{E((g-E)^2)}{b^2}
=\frac{V}{b^2}$$

Thus, we are led to the Chebycheff inequality in the statement.
\end{proof}

Getting back now to the real-life considerations from the beginning of this section, good that we have the notion of variance, for understanding how to run our business, but let us not stop here. For a total control of your business, be that of financial, mathematical, physical or chemical type, you will certainly want to know more about your variable $f:X\to\mathbb R$. Which leads us into general moments, constructed as follows:

\index{moments}

\begin{definition}
The moments of a variable $f:X\to\mathbb R$ are the numbers
$$M_k=E(f^k)$$
which satisfy $M_0=1$, then $M_1=E(f)$, and then $V(f)=M_2-M_1^2$.
\end{definition}

And with this, good news, we have now all the needed tools in our bag for doing some good business, in what follows. To put things in a very compacted way:

\medskip

-- $M_0$ is about foundations.

\medskip

-- $M_1$ is about running some business.

\medskip

-- $M_2$ is about running that business well.

\medskip

-- $M_3$ and higher are advanced level, about ruining all the competing businesses. 

\medskip

As a next piece of discussion, still regarding the moments, we can formulate the following version of Definition 1.16, making a more clear link with the variance:

\index{central moments}

\begin{definition}
The central moments of a variable $f:X\to\mathbb R$ are the numbers
$$M_k'=E((f-E)^k)$$
with $E=E(f)$, which satisfy $M_0'=1$, $M_1'=0$, $M_2'=V(f)$.
\end{definition}

And I will leave it to you, to think a bit at this, and decide whether you prefer this over Definition 1.16. Along the same lines, we can in fact do better, as follows:

\index{normalized central moments}

\begin{definition}
The normalized central moments of $f:X\to\mathbb R$ are the numbers
$$M''_k=E\left(\left(\frac{f-E}{\sigma}\right)^k\right)$$
with $\sigma=\sigma(f)$, which satisfy $M''_0=1$, $M''_1=0$, $M''_2=1$.
\end{definition}

Here we assume of course $\sigma>0$, which is the same as saying that $f:X\to\mathbb R$ is not constant. Of particular interest are the cases $k=3,4$, where we can formulate:

\index{skewness}
\index{kurtosis}

\begin{definition}
The third and fourth normalized central moments
$$\gamma=M_3''\quad,\quad 
\kappa=M_4''$$
are called skewness and kurtosis of the variable $f:X\to\mathbb R$.
\end{definition}

As a conclusion to all this, a random variable $f:X\to\mathbb R$ is quite reasonably described by the data $(E,V,\gamma,\kappa)$, whose knowledge is the same as that of $M_1,M_2,M_3,M_4$. And for more delicate questions, we still have the higher moments, $M_k$ with $k\geq5$.

\bigskip

Done with the theory? Not really, because we still have to formulate one more key result, which is actually the most important one, in basic probability, as follows:

\index{law of variable}

\begin{theorem}
If we define the law of $f:X\to\mathbb R$ as being the real measure
$$\mu=\sum_{x\in X}P(x)\delta_{f(x)}$$
then $\mu(\mathbb R)=1$, and we have the following formula, for any function $\phi:\mathbb R\to\mathbb R$,
$$E(\phi(f))=\int_\mathbb R\phi(y)d\mu(y)$$
with the usual convention that each Dirac mass integrates up to $1$.
\end{theorem} 

\begin{proof}
To start with, we can certainly talk about the law $\mu$, as being the linear combination of Dirac masses in the statement. Observe that, alternatively, we have:
$$\mu=\sum_{y\in\mathbb R}P(f=y)\delta_y$$

As yet another alternative definition for $\mu$, we have the following formula, with $\nu$ being the probability measure on $X$, and with $*$ being the push-forward operation:
$$\mu=f_*\nu$$

Getting now to the expectation formula in the statement, this is clear, as follows:
$$E(\phi(f))
=\sum_{x\in X}P(x)\phi(f(x))
=\int_\mathbb R\phi(y)d\mu(y)$$

We are therefore led to the conclusions in the statement.
\end{proof}

Next, quite remarkably, the sequence of moments uniquely determines the law:

\begin{theorem}
The sequence of moments of a variable $f:X\to\mathbb R$,
$$M_k=E(f^k)$$
uniquely determines the law of the variable.
\end{theorem}

\begin{proof}
By Theorem 1.20, we have the following formula, for the moments:
$$\mu=\sum_i\lambda_i\delta_{x_i}\implies M_k=\sum_i\lambda_ix_i^k$$

But it is then a standard question to recover the numbers $\lambda_i,x_i\in\mathbb R$, and so the measure $\mu$, out of the sequence of numbers $M_k$. Indeed, assuming that the numbers $x_i$ are $0<x_1<\ldots<x_n$ for simplifying, with $k\to\infty$ we have the following estimate:
$$M_k\sim\lambda_nx_n^k$$

Thus, we got the parameters $\lambda_n,x_n\in\mathbb R$ of our measure $\mu$, and by substracting them and doing an obvious recurrence, we get the other parameters $\lambda_i,x_i\in\mathbb R$ as well.
\end{proof}

Switching topics, inspired now by what happens for coins and dice, let us formulate: 

\index{independence}
\index{independent variables}

\begin{definition}
We say that two variables $f,g:X\to\mathbb R$ are independent when
$$P(f=x,g=y)=P(f=x)P(g=y)$$
happens, for any $x,y\in\mathbb R$.
\end{definition}

As mentioned, this is something quite intuitive, inspired by what happens for coins and dice. More on this, and on other examples available, later. For the moment, let us develop some theory. As our first result regarding independence, we have:

\index{mixed moments}

\begin{theorem}
Assuming that $f,g:X\to\mathbb R$ are independent, we have
$$E(f^kg^l)=E(f^k)E(g^l)$$
and the converse holds, in the sense that this implies the independence of $f,g$.
\end{theorem}

\begin{proof}
This is something very standard, the idea being as follows:

\medskip

(1) In one sense, we have the following computation, for the mixed moments:
\begin{eqnarray*}
E(f^kg^l)
&=&\sum_{xy}x^ky^lP(f=x,g=y)\\
&=&\sum_{xy}x^ky^lP(f=x)P(g=y)\\
&=&\sum_xx^kP(f=x)\sum_yy^lP(g=y)\\
&=&E(f^k)E(g^l)
\end{eqnarray*}

(2) Conversely, the moment condition $E(f^kg^l)=E(f^k)E(g^l)$ reformulates as:
$$\sum_{xy}c(x,y)x^ky^l=0\quad,\quad 
c(x,y)=P(f=x,g=y)-P(f=x)P(g=y)$$

So, let us examine this. Our variables $f,g:X\to\mathbb R$ being discrete, let $x_1<\ldots<x_n$ and $y_1<\ldots<y_m$ be their images. With $c_{ij}=c(x_i,x_j)$, our formula reads:
$$\sum_{ij}c_{ij}x_i^ky_j^l=0\quad,\quad\forall k,l\in\mathbb N$$

Assuming now $x_i>0$ for simplifying, we have, with $k\to\infty$, and with $l$ fixed:
$$\sum_{ij}c_{ij}x_i^ky_j^l\sim x_n^k\sum_jc_{nj}y_j^l$$

Thus $\sum_jc_{nj}y_j^l=0$, so $c_{nj}=0$. But with this, again with $k\to\infty$, and $l$ fixed:
$$\sum_{ij}c_{ij}x_i^ky_j^l\sim x_{n-1}^k\sum_jc_{n-1,j}y_j^l$$

And so on, the idea being that we obtain by recurrence $c_{ij}=0$, as desired.
\end{proof}

Getting now to sums of independent variables, we have here:

\begin{theorem}
Assuming that $f,g:X\to\mathbb R$ are independent, we have
$$\mu_{f+g}=\mu_f*\mu_g$$
where $*$ is the convolution operation, defined by $\delta_x*\delta_y=\delta_{x+y}$ and linearity.
\end{theorem}

\begin{proof}
We have indeed the following straightforward computation:
\begin{eqnarray*}
\mu_{f+g}
&=&\sum_{x\in\mathbb R}P(f+g=x)\delta_x\\
&=&\sum_{y,z\in\mathbb R}P(f=y,g=z)\delta_{y+z}\\
&=&\sum_{y,z\in\mathbb R}P(f=y)P(g=z)\delta_y*\delta_z\\
&=&\left(\sum_{y\in\mathbb R}P(f=y)\delta_y\right)*\left(\sum_{z\in\mathbb R}P(g=z)\delta_z\right)\\
&=&\mu_f*\mu_g
\end{eqnarray*}

Thus, we are led to the conclusion in the statement.
\end{proof}

As a first application of all this, we can now solve Challenge 1.11, as follows:

\begin{theorem}
When rolling a die $n$ times, the distribution of the sum is
$$P(s)=\frac{1}{6^n}\sum_{k=0}^{[\frac{s-n}{6}]}(-1)^k\binom{n}{k}\binom{s-6k-1}{n-1}$$
with this coming from independence, and the binomial formula for negative exponents.
\end{theorem}

\begin{proof}
This is something quite tricky, the idea being as follows:

\medskip

(1) According to our independence theory above, the law of the sum is:
$$\mu_n=\frac{1}{6^n}(\delta_1+\delta_2+\delta_3+\delta_4+\delta_5+\delta_6)^{*n}$$

(2) Equivalently, when thinking a bit, we have the following formula for the law of the sum, and I will leave it to you, to clarify the details here:
$$\sum_{s=n}^{6n}P(s)x^s=\frac{1}{6^n}(x+x^2+x^3+x^4+x^5+x^6)^n$$

(3) So, let us compute the function on the right. This function is given by:
$$f(x)
=\frac{x^n}{6^n}\cdot\frac{(1-x^6)^n}{(1-x)^n}$$

(4) Next, let us recall that the generalized binomial formula at negative integer exponents takes the following form, and with this being in fact something which is elementary, easy to establish by recurrence on $n$, by using the Pascal formula for binomials:
$$\frac{1}{(1-x)^n}=\sum_{l=0}^\infty\binom{n+l-1}{n-1}x^l$$

(5) Now by using the usual binomial formula for the numerator $(1-x^6)^n$, and the generalized  binomial formula for the fraction $1/(1-x)^n$, we obtain:
\begin{eqnarray*}
f(x)
&=&\frac{x^n}{6^n}\sum_{k=0}^n(-1)^k\binom{n}{k}x^{6k}\sum_{l=0}^\infty\binom{n+l-1}{n-1}x^l\\
&=&\frac{1}{6^n}\sum_{k=0}^n\sum_{l=0}^\infty(-1)^k\binom{n}{k}\binom{n+l-1}{n-1}x^{n+6k+l}\\
&=&\frac{1}{6^n}\sum_{s=0}^\infty\sum_{k=0}^n(-1)^k\binom{n}{k}\binom{s-6k-1}{n-1}x^s
\end{eqnarray*}

(6) Finally, remember that $f$ was a polynomial of degree $6n$, divisible by $x^n$. Thus, the true range of the first summing index is $s=n,\ldots,6n$. Also, from $s=n+6k+l$ in the above computation, we get $s\geq n+6k$, and so $k\leq[(s-n)/6]$, integer part. Thus, our formula in (5), written by ignoring zero coefficients appearing from cancellation, is:
$$f(x)=\frac{1}{6^n}\sum_{s=n}^{6n}\sum_{k=0}^{[\frac{s-n}{6}]}(-1)^k\binom{n}{k}\binom{s-6k-1}{n-1}x^s$$

But with this, we are led via (2) to the formula in the statement.
\end{proof}

\section*{1c. Binomial laws}

With Challenge 1.11 solved, what is next? Simpler things, I guess, and getting now to variables and laws that we can study, at the beginning of everything, we have:

\index{Bernoulli law}
\index{binomial law}

\begin{theorem}
When playing with a biased coin, your winning law is
$$b_p=(1-p)\delta_0+p\delta_1$$
called Bernoulli law of parameter $p\in[0,1]$. When playing $n$ times, your winning law is
$$b_{np}=\sum_{s=0}^n\binom{n}{s}p^s(1-p)^{n-s}\delta_s$$
called binomial law of parameter $p\in[0,1]$. We have the formula $b_{np}=b_p^{*n}$.
\end{theorem}

\begin{proof}
This is something coming from our independence technology, as follows:

\medskip

(1) The first assertion, regarding the Bernoulli law $b_p$, is more of a definition, for this law. Next, when playing $n$ times, in order to have $s$ wins, we must select the $s$ occurrences of heads, among the $n$ throws, and there are $\binom{n}{s}$ choices here. And then, we must sum over these $\binom{n}{s}$ choices the common probability $p^s(1-p)^{n-s}$, for $s$ heads and $n-s$ tails. Thus, we are led to the formula for the binomial law $b_{np}$ in the statement.

\medskip

(2) The formula $b_{np}=b_p^{*n}$ comes from the fact that the Bernoulli laws $b_p$ produce the binomial laws $b_{np}$, by iterating the game $n$ times, via the independence of the throws. Equivalently, this latter formula comes from the following formal computation:
\begin{eqnarray*}
b_p^{*n}
&=&\big((1-p)\delta_0+p\delta_1\big)^{*n}\\
&=&\sum_{s=0}^np^s(1-p)^{n-s}\binom{n}{s}\delta_0^{*(n-s)}*\delta_1^{*s}\\
&=&\sum_{s=0}^np^s(1-p)^{n-s}\binom{n}{s}\delta_s\\
&=&b_{np}
\end{eqnarray*}

Thus, we are led to the conclusions in the statement.
\end{proof}

Getting now to the study of the binomial laws, we first have the following result:

\index{mean of binomial law}
\index{variance of binomial law}
\index{moments of binomial law}

\begin{theorem}
The binomial law $b_{np}$ has the following properties:
\begin{enumerate}
\item The mean is $E=np$.

\item The variance is $V=np(1-p)$.
\end{enumerate}
\end{theorem}

\begin{proof}
Regarding the mean of a binomial variable $f:X\to\mathbb R$, we have:
\begin{eqnarray*}
E(f)
&=&\sum_{s=1}^ns\binom{n}{s}p^s(1-p)^{n-s}\\
&=&np\sum_{s=1}^n\frac{(n-1)!}{(s-1)!(n-s)!}\,p^{s-1}(1-p)^{n-s}\\
&=&np\sum_{t=0}^{n-1}\binom{n-1}{t}p^t(1-p)^{n-t-1}\\
&=&np(p+1-p)^{n-1}\\
&=&np
\end{eqnarray*}

Coming next, with the same trick, we can compute the following quantity:
\begin{eqnarray*}
E(f^2)-E(f)
&=&\sum_{s=2}^n(s^2-s)\binom{n}{s}p^s(1-p)^{n-s}\\
&=&n(n-1)p^2\sum_{s=2}^n\frac{(n-2)!}{(s-2)!(n-s)!}p^{s-2}(1-p)^{n-s}\\
&=&n(n-1)p^2\sum_{t=0}^{n-2}\binom{n-2}{t}p^t(1-p)^{n-t-2}\\
&=&n(n-1)p^2(p+1-p)^{n-2}\\
&=&n(n-1)p^2
\end{eqnarray*}

Now by using the above formulae, we conclude that the variance is given by:
\begin{eqnarray*}
V(f)
&=&[E(f^2)-E(f)]+[E(f)-E(f)^2]\\
&=&n(n-1)p^2+np-(np)^2\\
&=&np(1-p)
\end{eqnarray*}

We are therefore led to the conclusions in the statement.
\end{proof}

Regarding now the higher moments, we have the following result, about them:

\index{number of blocks}
\index{Bell numbers}
\index{Stirling numbers}

\begin{theorem}
The moments of the binomial law $b_{np}$ are given by
$$M_k=\sum_{\pi\in P(k)}\frac{n!}{(n-|\pi|)!}\,p^{|\pi|}$$
where $P(k)$ is the set of partitions of $\{1,\ldots,k\}$, and $|.|$ is the number of blocks.
\end{theorem}

\begin{proof}
This is something very standard, the idea being as follows:

\medskip

(1) Some numerics first. At $k=1$ we only have one partition, namely $|$\,, and the formula in the statement holds indeed, as shown by the following computation:
$$M_1=\sum_{\pi\in P(1)}\frac{n!}{(n-|\pi|)!}\,p^{|\pi|}
=\frac{n!}{(n-1)!}\,p^1
=np$$

At $k=2$ now, we have two partitions, namely $|\,|$ and $\sqcap$, and the moment is:
$$M_2=\frac{n!}{(n-2)!}\,p^2+\frac{n!}{(n-1)!}\,p^1
=n(n-1)p^2+np$$

At $k=3$ the partitions are $|\ |\ |$\,, $\sqcap\ |$\,, $\sqcap\hskip-3.2mm{\ }_|$\,\,, $|\ \sqcap$\,, $\sqcap\hskip-0.7mm\sqcap$, and the moment is:
$$M_3=n(n-1)(n-2)p^3+3n(n-1)p^2+np$$

At $k=4$ we have 15 partitions, with the list of relevant partitions, and their multiplicities up to permutations, which is self-explanatory, being as follows:
$$|\ |\ |\ |\times1\quad,\quad \sqcap\ |\ |\times 6\quad,\quad \sqcap\sqcap\times 3
\quad,\quad \sqcap\hskip-1.6mm\sqcap|\times4\quad,\quad \sqcap\hskip-1.6mm\sqcap\hskip-1.6mm\sqcap\times1$$

Now by doing the computation, as before, the formula of the moment is:
$$M_4=n(n-1)(n-2)(n-3)p^4+6n(n-1)(n-2)p^3+7n(n-1)p^2+np$$

(2) In general now, in order to prove the result, we can use the trick from the proof of Theorem 1.27. Indeed, assume that we managed to find a formula as follows:
$$s^k=\sum_{r=1}^kc_r\cdot s(s-1)\ldots(s-r+1)$$

We have then the following computation, based on this formula:
\begin{eqnarray*}
M_k
&=&\sum_{r=1}^kc_r\sum_{s=r}^ns(s-1)\ldots(s-r+1)\binom{n}{s}p^s(1-p)^{n-s}\\
&=&\sum_{r=1}^kc_r\cdot n(n-1)\ldots(n-r+1)\sum_{s=r}^n\binom{n-r}{s-r}p^s(1-p)^{n-s}\\
&=&\sum_{r=1}^kc_r\cdot n(n-1)\ldots(n-r+1)\cdot p^r
\end{eqnarray*}

(3) The problem is now, how to find that magic formula from (2)? And the answer here comes from partitions. Indeed, with the standard convention that for a multi-index $i=(i_1,\ldots,i_k)$, its kernel is the partition $\ker i\in P(k)$ joining equal indices, we have:
\begin{eqnarray*}
s^k
&=&\sum_{i_1=1}^s\ldots\sum_{i_k=1}^s1\\
&=&\sum_{\pi\in P(k)}\#\left\{(i_1,\ldots,i_k)\in\{1,\ldots,s\}^k\Big|\ker i=\pi\right\}\\
&=&\sum_{\pi\in P(k)}\frac{s!}{(s-|\pi|)!}
\end{eqnarray*}

Thus, good news, we have the magic formula that we need as data in (2), and according now to that computation in (2), we are led to the following formula, for the moment:
$$M_k=\sum_{\pi\in P(k)}\frac{n!}{(n-|\pi|)!}\,p^{|\pi|}$$

Thus, we have indeed the formula in the statement.
\end{proof}

Let us end this discussion with the following result, capturing the essentials:

\index{mean}
\index{variance}
\index{skewness}
\index{kurtosis}

\begin{theorem}
The mean, variance, skewness and kurtosis of $b_{np}$ are
$$E=np\quad,\quad V=npq\quad,\quad\gamma=\frac{q-p}{\sqrt{npq}}\quad,\quad\kappa=3+\frac{1-6pq}{npq}$$
with the convention $q=1-p$.
\end{theorem}

\begin{proof}
We know the first two formulae, from Theorem 1.27. Next, according to Definition 1.17, the central moments of a variable $f:X\to\mathbb R$ are given by:
\begin{eqnarray*}
M_k'
&=&E((f-E)^k)\\
&=&\sum_{r=0}^k(-1)^r\binom{k}{r}E^rM_{k-r}\\
&=&\sum_{r=0}^{k-2}(-1)^r\binom{k}{r}E^rM_{k-r}+(-1)^{k-1}(k-1)E^k
\end{eqnarray*}

At $k=3$, by using the formulae from the proof of Theorem 1.28, we obtain:
\begin{eqnarray*}
M_3'
&=&M_3-3EM_2+2E^3\\
&=&[n(n-1)(n-2)p^3+3n(n-1)p^2+np]-3np(n(n-1)p^2+np)+2(np)^3\\
&=&np(1-p)(1-2p)
\end{eqnarray*}

At $k=4$, again by using the formulae from the proof of Theorem 1.28, we have:
\begin{eqnarray*}
&&M_4'\\
&=&M_4-4EM_3+6E^2M_2-3E^4\\
&=&n(n-1)(n-2)(n-3)p^4+6n(n-1)(n-2)p^3+7n(n-1)p^2+np\\
&&\!\!-4np[n(n-1)(n-2)p^3+3n(n-1)p^2+np]+6(np)^2(n(n-1)p^2+np)-3(np)^4\\
&=&np(1-p)(1+(3n-6)p(1-p))
\end{eqnarray*}

Observe now that with the convention $q=1-p$, our formulae for $M_3',M_4'$ read:
$$M_3'=npq(q-p)\quad,\quad M_4'=npq(1+(3n-6)pq)$$

Now by dividing respectively by $\sigma^3=V\sqrt{V}$ and $\sigma^4=V^2$, we obtain the formulae in the statement for the skewness $\gamma$ and the kurtosis $\kappa$, as desired.
\end{proof}

\section*{1d. Poisson laws}

We would like to discuss now the central objects in discrete probability, which are the Poisson laws $p_t$, appearing via the Poisson Limit Theorem (PLT). Let us start with:

\index{Poisson law}
\index{Poisson variable}

\begin{definition}
The Poisson law of parameter $1$ is the following measure,
$$p_1=\frac{1}{e}\sum_{k\geq0}\frac{\delta_k}{k!}$$
and the Poisson law of parameter $t>0$ is the following measure,
$$p_t=e^{-t}\sum_{k\geq0}\frac{t^k}{k!}\,\delta_k$$
with the letter ``p'' standing for Poisson.
\end{definition}

We will see in the moment why the above laws appear a bit everywhere, in discrete contexts, the reasons for this coming from the Poisson Limit Theorem. In the meantime, let us first develop some straightforward theory. We first have the following result:

\index{mean of Poisson law}
\index{variance of Poisson law}

\begin{proposition}
The mean and variance of $p_t$ are given by:
$$E=t\quad,\quad V=t$$
In particular for the Poisson law $p_1$ we have $E=1,V=1$.
\end{proposition}

\begin{proof}
Regarding the mean, this can be computed as follows:
$$E
=e^{-t}\sum_{k\geq0}\frac{t^k}{k!}\cdot k
=e^{-t}\sum_{k\geq1}\frac{t^k}{(k-1)!}
=te^{-t}\sum_{l\geq0}\frac{t^l}{l!}
=t$$

For the variance, we can first compute the second moment, as follows:
\begin{eqnarray*}
M_2
&=&e^{-t}\sum_{k\geq0}\frac{t^k}{k!}\cdot k^2\\
&=&e^{-t}\sum_{k\geq1}\frac{t^kk}{(k-1)!}\\
&=&e^{-t}\sum_{l\geq0}\frac{t^{l+1}(l+1)}{l!}\\
&=&te^{-t}\sum_{l\geq0}\frac{t^ll}{l!}+te^{-t}\sum_{l\geq0}\frac{t^l}{l!}\\
&=&te^{-t}\sum_{l\geq1}\frac{t^l}{(l-1)!}+t\\
&=&t^2+t
\end{eqnarray*}

Thus the variance is given by $V=(t^2+t)-t^2=t$, as claimed.
\end{proof}

Next, we have the following conceptual result, in relation with convolution:

\index{convolution semigroup}

\begin{theorem}
We have the following formula, for any $s,t>0$,
$$p_s*p_t=p_{s+t}$$
so the Poisson laws form a convolution semigroup.
\end{theorem}

\begin{proof}
By using $\delta_k*\delta_l=\delta_{k+l}$ and the binomial formula, we obtain:
\begin{eqnarray*}
p_s*p_t
&=&e^{-s}\sum_{k\geq0}\frac{s^k}{k!}\,\delta_k*e^{-t}\sum_{l\geq0}\frac{t^l}{l!}\,\delta_l\\
&=&e^{-s-t}\sum_{n\geq0}\delta_n\sum_{k+l=n}\frac{s^kt^l}{k!l!}\\
&=&e^{-s-t}\sum_{n\geq0}\frac{\delta_n}{n!}\sum_{k+l=n}\frac{n!}{k!l!}s^kt^l\\\
&=&e^{-s-t}\sum_{n\geq0}\frac{(s+t)^n}{n!}\,\delta_n\\
&=&p_{s+t}
\end{eqnarray*}

Thus, we are led to the conclusion in the statement.
\end{proof}

Getting now to what we wanted to do, Poisson limits, it is convenient to use, as our main tool for dealing with independence, the Fourier transform, coming via:

\index{independence}
\index{Fourier transform}

\begin{theorem}
Assuming that $f,g:X\to\mathbb R$ are independent, we have
$$F_{f+g}=F_fF_g$$
where $F_f(x)=E(e^{ixf})$ is the Fourier transform.
\end{theorem}

\begin{proof}
We have indeed the following computation, using $\mu_{f+g}=\mu_f*\mu_g$:
\begin{eqnarray*}
F_{f+g}(x)
&=&\int_Xe^{ixz}d\mu_{f+g}(z)\\
&=&\int_Xe^{ixz}d(\mu_f*\mu_g)(z)\\
&=&\int_{X\times X}e^{ix(z+t)}d\mu_f(z)d\mu_g(t)\\
&=&\int_Xe^{ixz}d\mu_f(z)\int_Xe^{ixt}d\mu_g(t)\\
&=&F_f(x)F_g(x)
\end{eqnarray*}

Thus, we are led to the conclusion in the statement.
\end{proof}

Regarding now the Fourier transform computation for $p_t$, this is as follows:

\index{Fourier transform}

\begin{proposition}
The Fourier transform of $p_t$ is given by
$$F_{p_t}(x)=\exp\left((e^{ix}-1)t\right)$$
for any $t>0$, with its logarithm being linear in $t$, as it should.
\end{proposition}

\begin{proof}
We have indeed the following computation, which gives the result:
\begin{eqnarray*}
F_{p_t}(x)
&=&e^{-t}\sum_{k\geq0}\frac{t^k}{k!}\,e^{ikx}\\
&=&e^{-t}\sum_{k\geq0}\frac{(e^{ix}t)^k}{k!}\\
&=&\exp(-t)\exp(e^{ix}t)\\
&=&\exp\left((e^{ix}-1)t\right)
\end{eqnarray*}

As for the last assertion, this is just a remark, to be related to Theorem 1.32.
\end{proof}

Good news, we can now establish the Poisson Limit Theorem, as follows:

\index{PLT}
\index{Poisson Limit Theorem}
\index{Bernoulli laws}
\index{Poisson limit}

\begin{theorem}[PLT]
We have the following convergence, in moments,
$$\left(\left(1-\frac{t}{n}\right)\delta_0+\frac{t}{n}\,\delta_1\right)^{*n}\to p_t$$
for any $t>0$. Equivalently, $b_{np}\to p_t$, when $p=t/n$ with $t>0$ fixed.
\end{theorem}

\begin{proof}
Let us denote by $\mu_n$ the Bernoulli law under the convolution sign:
$$\mu_n=\left(1-\frac{t}{n}\right)\delta_0+\frac{t}{n}\,\delta_1$$

We have the following computation, for the Fourier transform of the limit: 
\begin{eqnarray*}
F_{\delta_r}(x)=e^{irx}
&\implies&F_{\mu_n}(x)=\left(1-\frac{t}{n}\right)+\frac{t}{n}\,e^{ix}\\
&\implies&F_{\mu_n^{*n}}(x)=\left(\left(1-\frac{t}{n}\right)+\frac{t}{n}\,e^{ix}\right)^n\\
&\implies&F_{\mu_n^{*n}}(x)=\left(1+\frac{(e^{ix}-1)t}{n}\right)^n\\
&\implies&F(x)=\exp\left((e^{ix}-1)t\right)
\end{eqnarray*}

Thus, we obtain indeed the Fourier transform of $p_t$, as desired.
\end{proof}

In practice, Theorem 1.35 leads to many occurrences of the Poisson laws, in relation with continuous phenomena from the real life. For instance, we have:

\begin{theorem}
The number of chewing gum pieces on a single tile of a sidewalk
$$\xymatrix@R=5pt@C=7pt{
\ar@{-}[rrrrrrrrrr]\ar@{-}[dddddd]&&\ar@{-}[dddddd]&&\ar@{-}[dddddd]&&\ar@{-}[dddddd]&&\ar@{-}[dddddd]&&\ar@{-}[dddddd]\\
&\bullet&&&&\bullet&&\ \,\bullet\bullet\,\ &&&&\ldots\\
\ar@{-}[rrrrrrrrrr]&&&&&&&&&&\\
&&&\,\bullet\bullet\bullet\,&&&&\bullet&&\ \,\bullet\bullet\,\ &&\ldots\\
\ar@{-}[rrrrrrrrrr]&&&&&&&&&&\\
&\ \ \, \bullet\ \ \,&&&&\ \,\bullet\bullet\,\ &&\bullet&&&&\ldots\\
\ar@{-}[rrrrrrrrrr]&&&&&&&&&&
}$$
follows a Poisson law.
\end{theorem}

\begin{proof}
Please don't throw your chewing gum on the sidewalk. Instead, you can deduce this from Theorem 1.35, with some thinking here being an excellent exercise.
\end{proof}

Regarding now the moments of the Poisson laws, the result here is as follows:

\index{number of blocks}
\index{Stirling numbers}

\begin{theorem}
The moments of $p_t$ with $t>0$ are given by
$$M_k(p_t)=\sum_{\pi\in P(k)}t^{|\pi|}$$
where $|.|$ is the number of blocks.
\end{theorem}

\begin{proof}
We have the following recurrence formula for the moments:
\begin{eqnarray*}
M_{k+1}
&=&e^{-t}\sum_{n\geq1}\frac{t^nn^{k+1}}{n!}\\
&=&e^{-t}\sum_{m\geq0}\frac{t^{m+1}(m+1)^k}{m!}\\
&=&e^{-t}\sum_{m\geq0}\frac{t^{m+1}m^k}{m!}\left(1+\frac{1}{m}\right)^k\\
&=&e^{-t}\sum_{m\geq0}\frac{t^{m+1}m^k}{m!}\sum_{s=0}^k\binom{k}{s}m^{-s}\\
&=&\sum_{s=0}^k\binom{k}{s}\cdot e^{-t}\sum_{m\geq0}\frac{t^{m+1}m^{k-s}}{m!}\\
&=&t\sum_{s=0}^k\binom{k}{s}M_{k-s}
\end{eqnarray*}

But the numbers in the statement are easily seen to satisfy the same recurrence, and the initial values match too, so we are led to the conclusion in the statement.
\end{proof}

Let us end this discussion with the following result, capturing the essentials:

\begin{theorem}
For the Poisson law $p_t$ of parameter $t>0$,
$$E=t\quad,\quad  V=t\quad ,\quad \gamma=\frac{1}{\sqrt{t}}\quad ,\quad \kappa=3+\frac{1}{t}$$
are the mean, variance, skewness and kurtosis.
\end{theorem}

\begin{proof}
Regarding the formulae of $E,V$, we know these from Proposition 1.31. As for the formulae of $\gamma,\kappa$, these follow from Theorem 1.37. Indeed, that formula gives:
$$M_3=t+3t^2+t^3$$
$$M_4=t+7t^2+6t^3+t^4$$

By passing now to central moments, as in the proof of Theorem 1.29, we obtain:
$$M_3'=t$$
$$M_4'=t+3t^2$$

Thus, by normalizing, we are led to the formulae of $\gamma$ and $\kappa$ in the statement.
\end{proof}

\section*{1e. Exercises}

This was a standard introduction to probability, and as exercises, we have:

\begin{exercise}
Do some computations for card games, that we forgot to talk about.
\end{exercise}

\begin{exercise}
Further meditate at the symmetry question, for biased coins and dice.
\end{exercise}

\begin{exercise}
Think a bit at measures, and at the push-forward formula $\mu=f_*\nu$.
\end{exercise}

\begin{exercise}
Learn the binomial formula with arbitrary real exponents.
\end{exercise}

\begin{exercise}
Learn about geometric laws, and negative binomial laws.
\end{exercise}

\begin{exercise}
Learn as well about hypergeometric laws, positive and negative.
\end{exercise}

\begin{exercise}
Learn a bit about the Bell numbers, $B_k=|P(k)|=M_k(p_1)$.
\end{exercise}

\begin{exercise}
Read about compound Poisson laws, and the compound PLT.
\end{exercise}

As bonus exercise, read as well about Catalan numbers. We will meet them later.

\chapter{Normal variables}

\section*{2a. Variables, laws}

We have seen so far that some substantial discrete probability theory can be developed, based on some very simple axioms, which were basically as follows:

\begin{fact}
The axioms for discrete probability theory are as follows:
\begin{enumerate}
\item A probability space is a discrete measured space $(X,\nu)$ of mass one.

\item A random variable on $X$ is a real function $f:X\to\mathbb R$.

\item The expectation of $f:X\to\mathbb R$ is its integral, $E(f)=\int_Xf(x)d\nu(x)$.

\item The law of $f$ is the real probability measure $\mu=f_*\nu$, push-forward of $\nu$ by $f$.
\end{enumerate}
\end{fact}

We would like to talk now, as a continuation of this, about continuous probability. And here, thinking at measures, as a first surprise that awaits us, we have:

\begin{proposition}
In the continuous setting, some sets are not measurable.
\end{proposition}

\begin{proof}
Consider the unit circle $\mathbb T$, measured as usual, the total mass being $\mu(\mathbb T)=2\pi$, and then consider the group $G$ consisting of the rational rotations of $\mathbb T$:
$$G=\left\{z\to e^{2r\pi i}z\Big|r\in\mathbb Q\cap[0,1)\right\}$$

Now if we choose $A\subset\mathbb T$ containing exactly one point from each orbit, we have:
$$\mathbb T=\bigsqcup_{r\in\mathbb Q\cap[0,1)}e^{2r\pi i}A$$  

But this shows that $A$ is not measurable. Indeed, with the convention that we list our group as $\mathbb Q\cap[0,1)=\{r_1,r_2,r_3,\ldots\}$, we have two possible cases, as follows:

\medskip

(1) Assuming $\mu(A)=c>0$, we would get from this a contradiction, as follows:
$$\mu(\mathbb T)=\sum_{n=1}^\infty\mu(e^{2r_n\pi i}A)=\sum_{n=1}^\infty c=c\cdot\infty=\infty$$

(2) Assuming $\mu(A)=0$, we would get a contradiction too, because for $\varepsilon>0$:
$$\mu(\mathbb T)=\sum_{n=1}^\infty\mu(e^{2r_n\pi i}A)<\sum_{n=1}^\infty\frac{\varepsilon}{2^n}=\varepsilon$$

Thus, our set $A$ cannot be measurable, and we will have to live with that.
\end{proof}

Summarizing, as a first task in our axiomatization work, we must start with a set $X$, and come up with axioms for the measurable subsets $Y\subset X$. And here, we have:

\index{measurable set}

\begin{definition}
A measurable space is a set $X$, given with a set of subsets $M\subset P(X)$, called measurable sets, which form an algebra, in the sense that:
\begin{enumerate}
\item $\emptyset,X\in M$.

\item $E\in M\implies E^c\in M$.

\item $M$ is stable under countable unions and intersections.
\end{enumerate}
\end{definition}

Observe that some of the axioms above are redundant, because assuming (2) we have $\emptyset\in M\iff X\in M$, which can help in verifying (1). The same goes for (3), with only one of the conditions there being in need to be verified, due to $\left(\bigcup_iE_i\right)^c=\bigcap_iE_i^c$.

\bigskip

At the level of examples, these usually come via the following simple fact:

\begin{proposition}
Given a set of subsets $S\subset P(X)$, there is a smallest algebra
$$M=\bar{S}$$
containing it, called algebra generated by $X$.
\end{proposition}

\begin{proof}
This can be viewed in two possible ways, as follows:

\medskip

(1) In order to construct $M=\bar{S}$ we have to start with $S$, then add $\emptyset,X$ to it, along with the complements $E^c$ of all the sets $E\in S$, and then take countable unions and intersections of such sets. And, what we get in this way is indeed an algebra.

\medskip

(2) Alternatively, we can define $M=\bar{S}$ as being the intersection of all algebras containing $S$, and with the remark that we have at least one such algebra, namely $P(X)$ itself. It is then clear from definitions that $M$ is an algebra, as desired.
\end{proof}

Getting now to the concrete examples of measurable spaces, we have:

\index{Borel set}

\begin{proposition}
Any metric space $X$ is automatically a measurable space, with the algebra of measurable sets being 
$$B=\bar{O}$$
that is, the smallest algebra containing the open sets, called Borel algebra of $X$.
\end{proposition}

\begin{proof}
This is indeed something self-explanatory, and trivial, based on the construction in Proposition 2.4. Observe that the Borel sets include all open sets, all closed sets, as well as all countable unions of closed sets, and all countable intersections of open sets. As an example here, in the case $X=\mathbb R$, all kinds of intervals are Borel sets:
$$(a,b),\ [a,b],\ (a,b],\ [a,b)\in B$$

Indeed, the first interval is open, and the second one is closed, so these are certainly Borel. As for the third and fourth intervals, these appear as countable unions of closed intervals, or as countable intersections of open intervals, so they are Borel too.
\end{proof}

Getting back now to the general case, as a complement to Definition 2.3, we have:

\begin{definition}
Given a measurable space $(X,M)$, a measure on it is a function $\mu:M\to[0,\infty]$ which is countably additive, in the sense that
$$\mu\left(\bigcup_{i=1}^\infty E_i\right)=\sum_{i=1}^\infty\mu(E_i)$$
for any countable family of disjoint measurable sets $E_i\in M$. In this case, we say that $(X,M,\mu)$ is a measured space.
\end{definition}

Observe that we did not assume that our measures are finite, and this for including spaces like $X=\mathbb R$. In the case where $\mu$ happens to be bounded, up to a rescaling we can assume $\mu(X)=1$, and we say in this case that we have a probability measure on $X$.

\bigskip

Good news, we can now quickly axiomatize probability theory, as folows:

\begin{definition}
The axioms for probability theory are as follows:
\begin{enumerate}
\item A probability space is a measured space $(X,M,\nu)$ of mass one, $\nu(X)=1$.

\item A random variable on $X$ is a real measurable function $f:X\to\mathbb R$.

\item The expectation of $f:X\to\mathbb R$ is its integral, $E(f)=\int_Xf(x)d\nu(x)$.

\item The law of $f$ is the real probability measure $\mu=f_*\nu$, push-forward of $\nu$ by $f$.
\end{enumerate}
\end{definition}

In short, what we are doing here is to extend our previous discrete probability theory, axiomatized as in Fact 2.1, by using our measure theory knowledge, and with the main difference coming from the need of the set $M\subset P(X)$ of measurable sets.

\bigskip

Observe that, as another new phenomenon appearing, in this setting, the expectation $E(f)=\int_Xf(x)d\nu(x)$ can now take infinite values, $\pm\infty$, or even be undefined, due to the lack of convergence of the integral. So, in a word, always careful with this.

\bigskip

In practice now, Definition 2.7 remains something a bit abstract, and in what regards the random variables and their laws, which is what matters the most, we have:

\begin{theorem}
Given a random variable $f:X\to\mathbb R$ having law $\mu$, we have the following formula, for any measurable function $\phi:\mathbb R\to\mathbb R$,
$$E(\phi(f))=\int_\mathbb R\phi(x)d\mu(x)$$
and this can stand, via the Riesz theorem, as an alternative definition for $\mu$.
\end{theorem}

\begin{proof}
This is something quite self-explanatory, the idea being as follows:

\medskip

(1) To start with, the formula in the statement for the expectation $E(\phi(f))$ comes indeed from our definition of the law as a push-forward measure, $\mu=f_*\nu$.

\medskip

(2) As for the last assertion, this comes from the Riesz theorem, stating that a measure $\mu$ is uniquely determined by the integration over it. See Rudin \cite{rud}.
\end{proof}

As in the discrete probability setting, the law of a variable $f:X\to\mathbb R$ is not the only main quantity related to $f$ that we can talk about. As a rival to it, we have:

\index{CDF}
\index{cumulative distribution function}

\begin{definition}
Given a random variable $f:X\to\mathbb R$, the function
$$F(x)=P(f\leq x)$$
is called its cumulative distribution function (CDF).
\end{definition}

Observe that the cumulative distribution function $F:\mathbb R\to[0,1]$ must be by definition increasing, and also, that it must have the following limiting properties:
$$F(-\infty)=0\quad,\quad F(\infty)=1$$

With a bit more care, by looking into continuity properties as well, we are led to the following result, regarding the CDF of the general random variables $f:X\to\mathbb R$:

\index{right continuous}

\begin{theorem}
The cumulative distribution function $F:\mathbb R\to[0,1]$ must be
\begin{enumerate}
\item Increasing,

\item Right continuous,

\item Such that $F(-\infty)=0$,

\item Such that $F(\infty)=1$,
\end{enumerate}
and conversely, any function satisfying these properties appears as a CDF.
\end{theorem}

\begin{proof}
Here the first assertion is quite clear, and I will leave it to you to check the right continuity property, after remembering what that property exactly is. As for the second assertion, this is something which is clear too, good exercise for you.
\end{proof}

Moving on, as something which is new, of genuine continuous nature, we have:

\begin{theorem}
Assume that a random variable $f:X\to\mathbb R$ has density $\varphi:\mathbb R\to\mathbb R$, in the sense that its law is given by $d\mu(x)=\varphi(x)dx$.
\begin{enumerate}
\item The density must be positive, $\varphi\geq0$, and of mass one, $\int_\mathbb R\varphi(x)dx=1$.

\item $\mu(Z)=\int_Z\varphi(x)dx$, which can stand as an alternative definition for $\varphi$.

\item $E(\phi(f))=\int_\mathbb R\phi(x)\varphi(x)dx$, which can stand as well as a definition for $\varphi$.

\item The CDF is $F(x)=\int_{-\infty}^x\varphi(t)dt$, and conversely, we have $\varphi(x)=F'(x)$.
\end{enumerate}
\end{theorem}

\begin{proof}
As before with other such things, this is quite self-explanatory, as follows:

\medskip

(1) This is something clear, coming from $d\mu(x)=\varphi(x)dx$.

\medskip

(2) This is again clear, also coming from $d\mu(x)=\varphi(x)dx$.

\medskip

(3) This comes indeed via the general theory from Theorem 2.8.

\medskip

(4) This comes from the fundamental theorem of calculus.
\end{proof}

We can in fact merge Theorem 2.11 with our previous theory from chapter 1, and we are led in this way to the following result, summarizing our knowledge so far:

\index{generalized density}

\begin{theorem}
Assume that a random variable $f:X\to\mathbb R$ has generalized density, in the sense that its law is given by $d\mu(x)=\varphi(x)dx+\sum_ic_i\delta_{x_i}$.
\begin{enumerate}
\item The density must satisfy $\varphi\geq0$, $c_i>0$ and $\int_\mathbb R\varphi(x)dx+\sum_ic_i=1$.

\item $\mu(Z)=\int_Z\varphi(x)dx+\sum_{x_i\in Z}c_i$, with this uniquely determining the density.

\item $E(\phi(f))=\int_\mathbb R\phi(x)\varphi(x)dx+\sum_ic_i\phi(x_i)$, uniquely determining the density too.

\item The CDF is $F(x)=\int_{-\infty}^x\varphi(t)dt$, and conversely, we have $\varphi(x)=F'(x)$.
\end{enumerate}
\end{theorem}

\begin{proof}
Again, this is something quite self-explanatory, as follows:

\medskip

(1) This is clear, of course under the assumption that the points $x_i$ are distinct.

\medskip

(2) This is something which is clear, from what we know.

\medskip

(3) This is again something which is clear, from what we know.

\medskip

(4) This comes with the convention that jumps differentiate into Dirac masses.
\end{proof}

Quite nice all this, we are learning new things here, with respect to what we previously knew from chapter 1, this is first class theory, for sure. However, no idea if you noticed, but there is an elephant in the room. So, let's listen to what this elephant has to say:

\begin{elephant}
Not all measures on $\mathbb R$ have a generalized density. And this because big can be small, and vice versa, I mean size is something relative.
\end{elephant}

Okay, thanks elephant, and by the way no worries regarding your size, both me and cat we love you a lot, somehow you're to me what I am to cat, and so does rat, and with my calculus teacher salary I can certainly feed the four of us, no problems with that.

\bigskip

In practice now, getting now to what elephant was saying, now that he's gone for a walk, besides the fact that I will have to talk with cat and rat, see how we can find some cash, that actually sounds quite terrifying. Here is a mathematical result, about this:

\index{Cantor set}

\begin{proposition}
Exotic probability measures on $\mathbb R$ exist, with a standard example being the Cantor distribution, supported by the Cantor set $K\subset\mathbb R$.
\end{proposition}

\begin{proof}
Consider indeed the Cantor set $K\subset\mathbb R$, which is something quite tricky, appearing by definition as the intersection of the following unions of intervals:
$$\xymatrix@R=10pt@C=10pt{
K_0:&0\ar@{-}[rrrrrrrrr]&&&&&&&&&1\\
K_1:&0\ar@{-}[rrr]&&&\frac{1}{3}&&&\frac{2}{3}\ar@{-}[rrr]&&&1\\
K_2:&0\ar@{-}[r]&\frac{1}{9}&\frac{2}{9}\ar@{-}[r]&\frac{1}{3}&&&\frac{2}{3}\ar@{-}[r]&\frac{7}{9}&\frac{8}{9}\ar@{-}[r]&1\\
&&&&&\vdots&\vdots}$$

As a first observation, the measure of this Cantor set $K$ is given by:
$$\lambda(K)
=1-\frac{1}{3}-\frac{2}{9}-\frac{4}{27}-\ldots
=0$$

The point now is that we can construct a natural probability measure on $K$, by stating that $\mu$ must be uniform with respect to the 2 intervals of $K_1$, then with respect to the 4 intervals of $K_2$, and so on. And, exercise for you, to learn more about all this.
\end{proof}

Well, so what to do? In the lack of a better idea, I will ask the rat. Our fellow is as usual, just a few feet away, hiding behind the radiator, and here is what he says:

\begin{rat}
Peace upon you, and the truth can be found in Rudin. Read about the Lebesgue measure, and Radon-Nikodym. By the way, don't forget independence.
\end{rat}

Thanks rat, this sounds wise. So, leaving now to you reader the pleasure to consult Rudin \cite{rud}, that's best done when you're young, and is part of the standard initiation to mathematics, like these tribal rites that we used to have in the good old days, let me get now to independence. We can start our study here with the following definition:

\index{independence}
\index{independent variables}

\begin{definition}
We say that two variables $f,g:X\to\mathbb R$ are independent when
$$P(f\in I,g\in J)=P(f\in I)P(g\in J)$$
happens, for any two intervals $I,J\subset\mathbb R$.
\end{definition}

As a first observation, in the discrete case this is the same as our notion of independence from chapter 2, which was asking for the following equality, for any $x,y\in\mathbb R$:
$$P(f=x,g=y)=P(f=x)P(g=y)$$

Indeed, the equivalence comes as a particular case of the following statement:

\begin{proposition}
Two variables $f,g:X\to\mathbb R$ are independent precisely when
$$P(f\in K,g\in L)=P(f\in K)P(g\in L)$$
happens, for any two measurable sets $K,L\subset\mathbb R$.
\end{proposition}

\begin{proof}
In the case where our measurable sets $K,L\subset\mathbb R$ are disjoint unions of intervals, $K=\sqcup_rI_r$ and $L=\sqcup_sJ_s$, we have indeed the following computation:
\begin{eqnarray*}
P(f\in K,g\in L)
&=&\sum_r\sum_sP(f\in I_r,g\in J_s)\\
&=&\sum_r\sum_sP(f\in I_r)P(g\in J_s)\\
&=&\sum_rP(f\in I_r)\sum_sP(g\in J_s)\\
&=&P(f\in K)P(g\in L)
\end{eqnarray*}

As for the general case, this follows from this, via an approximation argument.
\end{proof}

Importantly, Proposition 2.17 allows us to pass to expectations, as follows:

\begin{proposition}
Two variables $f,g:X\to\mathbb R$ are independent precisely when
$$E(\phi(f)\psi(g))=E(\phi(f))E(\psi(g))$$
happens, for any two measurable functions $\phi,\psi:\mathbb R\to\mathbb R$.
\end{proposition}

\begin{proof}
We know from Proposition 2.17 that the independence of $f,g:X\to\mathbb R$ is equivalent to the following condition, for any two measurable sets $K,L\subset\mathbb R$:
$$P(f\in K,g\in L)=P(f\in K)P(g\in L)$$

Now observe that in terms of characteristic functions, this condition reads:
$$E(\chi_K(f)\chi_L(g))=E(\chi_K(f))E(\chi_L(g))$$

But with this, we are done, by approximating $\phi,\psi:\mathbb R\to\mathbb R$ by step functions.
\end{proof}

As a third and main result now on the subject, generalizing what we knew from chapter 1, in the discrete case, we have the moment characterization of independence:

\index{mixed moments}

\begin{theorem}
Assuming that $f,g:X\to\mathbb R$ are independent, we have
$$E(f^kg^l)=E(f^k)E(g^l)$$
and the converse holds, in the sense that this implies the independence of $f,g$.
\end{theorem}

\begin{proof}
This folllows indeed from Proposition 2.18, by approximating the continuous functions, and so the measurable functions too, by polynomials. We will leave the details here, including some learning about polynomial approximation, as an exercise. 
\end{proof}

Still following the material from chapter 1, as a next task, we must talk about convolution of measures. In the present continuous case, this can be done as follows:

\index{convolution}
\index{semigroup}
\index{group}

\begin{definition}
Given a space $X$ with a sum operation $+$, we can define the convolution of any two probability measures on it by the following formula:
$$(\mu*\nu)(E)=(\mu\times\nu)\left((x,y)\Big|x+y\in E\right)$$
Equivalently, in terms of the associated integrals, we must have 
$$\int_Xf(x)d(\mu*\nu)(x)=\int_{X\times X}f(x+y)d\mu(x)d\nu(y)$$
for any measurable function $f:X\to\mathbb R$.
\end{definition}

As a first observation, this fits with our previous notion of convolution of discrete measures, from chapter 1, because the convolution of two Dirac masses is given by:
$$(\delta_x*\delta_y)(E)=(\delta_x\times\delta_y)\left((x,y)\Big|x+y\in E\right)=\delta_{x+y}(E)$$

In general now, observe that our $*$ operation in Definition 2.20 is indeed well-defined, the total mass being 1. Also, the equivalent formulation, in terms of integrals, is clear by using characteristic functions, as shown by the following computation:
\begin{eqnarray*}
(\mu*\nu)(E)
&=&\int_X\chi_E(x)d(\mu*\nu)(x)\\
&=&\int_{X\times X}\chi_E(x+y)d\mu(x)d\nu(y)\\
&=&(\mu\times\nu)\left((x,y)\Big|x+y\in E\right)
\end{eqnarray*}

In relation now with the notion of independence, we have the following key result, extending what we previously knew from chapter 1, in the discrete case:

\begin{theorem}
Assuming that $f,g:X\to\mathbb R$ are independent, we have
$$\mu_{f+g}=\mu_f*\mu_g$$
where $*$ is the convolution of real probability measures.
\end{theorem}

\begin{proof}
There are several possible proofs here, and with our present knowledge, which is rather functional analysis oriented, the simplest is to use moments. Indeed, we have the following computation, using the formula for mixed moments in Theorem 2.19:
\begin{eqnarray*}
M_k(f+g)
&=&E((f+g)^k)\\
&=&\sum_r\binom{k}{r}E(f^rg^{k-r})\\
&=&\sum_r\binom{k}{r}M_r(f)M_{k-r}(g)
\end{eqnarray*}

On the other hand, if $h:X\to\mathbb R$ is a random variable following the law $\mu_f*\mu_g$, the corresponding moments are given by the following formula:
\begin{eqnarray*}
M_k(h)
&=&\int_Xx^kd(\mu_f*\mu_g)(x)\\
&=&\int_{X\times X}(x+y)^kd\mu_f(x)d\mu_g(y)\\
&=&\sum_r\binom{k}{r}\int_Xx^rd\mu_f(x)\int_Xy^{k-r}d\mu_g(y)\\
&=&\sum_r\binom{k}{r}M_r(f)M_{k-r}(g)
\end{eqnarray*}

Thus, our laws have the same moments. But with this in hand, it follows from the general theory from Theorem 2.8 that our laws must coincide, as stated.
\end{proof}

Still in relation with the basic theory of independence, here is now a second result, coming as a continuation of Theorem 2.21, which is something more advanced:

\index{independence}
\index{Fourier transform}

\begin{theorem}
Assuming that $f,g:X\to\mathbb R$ are independent, we have
$$F_{f+g}=F_fF_g$$
where $F_f(x)=E(e^{ixf})$ is the Fourier transform.
\end{theorem}

\begin{proof}
We have the following computation, using Theorem 2.21:
\begin{eqnarray*}
F_{f+g}(x)
&=&\int_Xe^{ixz}d\mu_{f+g}(z)\\
&=&\int_Xe^{ixz}d(\mu_f*\mu_g)(z)\\
&=&\int_{X\times X}e^{ix(z+t)}d\mu_f(z)d\mu_g(t)\\
&=&\int_Xe^{ixz}d\mu_f(z)\int_Xe^{ixt}d\mu_g(t)\\
&=&F_f(x)F_g(x)
\end{eqnarray*}

Thus, we are led to the conclusion in the statement.
\end{proof}

And with this, end of our discussion regarding.. but wait, cat is here, back from her hunt, let's see what she has to say about this, difficulties of measure theory:

\begin{cat}
In our daily feline work, we impose regularity conditions on the CDF.
\end{cat}

Well, I must admit that this sounds really interesting. So thanks cat, and next time I teach a probability class, I will think about this, modifying a bit the axioms, along the lines that you suggest. In what concerns us, what we have in the above will do.

\section*{2b. Large numbers}

Time now for some concrete mathematics, based on the above? The question that we would like to solve, which looks like something quite fundamental, is as follows:

\begin{question}
Given a sequence of independent identically distributed (i.i.d.) variables $f_1,f_2,f_3,\ldots\,$, what can we say about their partial sums
$$S_n=f_1+\ldots+f_n$$
with $n\to\infty$? Do we have convergence, under which normalization, and in which exact sense? What are the mean, variance of the limit? Is the limiting law continuous?
\end{question}

So, this will be our question, which is something quite subtle. More in detail now, here are a few comments about this, meant to clarify what our question exactly says:

\bigskip

(1) As a first observation, the PLT, namely $b_{np}\to p_t$ with $p=t/n$, does not fall in this class, due to the conditioning $p=t/n$ used for the i.i.d. variables appearing there. And so we are here, right from the beginning, into uncharted territory.

\bigskip

(2) As a second observation, due to our i.i.d. assumption on our variables, by obvious mean and variance reasons, the above sums $S_n$ cannot converge as such, without a suitable normalization, or an extra assumption on our variables. More on this later.

\bigskip

(3) Finally, we must talk about the meaning of the convergence too. Indeed, even if our variables $f_1,f_2,f_3,\ldots\,$ are discrete, their sums $S_n$, suitably normalized, can converge to bell-shaped curves. So, we need some convergence theory for such situations.

\bigskip

Quite interesting all this, hope you agree with me. So, let us start with the beginning, which as usual in probability, means study of the expectation $E$. We have here:

\begin{question}[soft version]
Given a sequence of i.i.d. variables $f_1,f_2,f_3,\ldots\,$, intuition tells us that we should have a convergence as follows,
$$\frac{f_1+\ldots+f_n}{n}\to E$$
 with $E=E(f_i)$ being the common expectation of our variables. So, do we really have this, and in which exact sense, in what regards the convergence?
\end{question}

Which sounds very good and concrete, as a starting question, so let us get now to work. Obviously, we must first talk about convergence of random variables. And here, we have several possible notions. As a first notion, which is quite natural, we have:

\index{convergence in law}
\index{convergence in moments}

\begin{definition}
Given variables $f_1,f_2,f_3,\ldots\,$, we say that $f_n\to f$ in law when 
$$E(\varphi(f_n))\to E(\varphi(f))$$
for any continuous function $\varphi:\mathbb R\to\mathbb R$.
\end{definition}

As a first observation, by linearity we can assume that we are dealing with the power functions $\varphi(x)=x^k$, so the following condition must be satisfied, for any $k\in\mathbb N$:
$$M_k(f_n)\to M_k(f)$$

Alternatively, again by linearity, and by using some standard analysis, we can assume that we are dealing with characteristic functions of intervals, which leads to:

\begin{proposition}
We have $f_n\to f$ in law precisely when
$$P(f_n\leq x)\to P(f\leq x)$$
at all continuity points of the function $x\to P(f\leq x)$.
\end{proposition}

\begin{proof}
This is indeed very standard, coming from definitions, as explained above, by using linearity, and approximation by step functions, say good exercise for you.
\end{proof}

Moving on, as a second notion of convergence for random variables, we have:

\index{convergence in probability}

\begin{definition}
We say that we have $f_n\to f$ in probability when 
$$P(|f_n-f|\geq\varepsilon)\to 0$$
for any $\varepsilon>0$.
\end{definition}

As before with the convergence in law, there are many illustrations for this notion, and further things that can be said. Before everything, however, let us formulate:

\begin{theorem}
We have the following implication,
$$\begin{pmatrix}f_n\to f\\ {\rm in\ probability}\end{pmatrix}\implies
\begin{pmatrix}f_n\to f\\ {\rm in\ law}\end{pmatrix}$$
and the reverse implication does not necessarily hold.
\end{theorem}

\begin{proof}
As a technical ingredient, we will need the following estimate, valid for any two random variables $f,g:X\to\mathbb R$, and any two numbers $a\in\mathbb R$ and $\varepsilon>0$:
\begin{eqnarray*}
P(g\leq a)
&=&P(g\leq a,f\leq a+\varepsilon)+P(g\leq a,f>a+\varepsilon)\\
&\leq&P(f\leq a+\varepsilon)+P(g\leq a,f>a+\varepsilon)\\
&=&P(f\leq a+\varepsilon)+P(g-f\leq a-f<-\varepsilon)\\
&\leq&P(f\leq a+\varepsilon)+P(g-f<-\varepsilon)\\
&\leq&P(f\leq a+\varepsilon)+P(|g-f|>\varepsilon)
\end{eqnarray*}

In order to prove now the result, the most convenient is to use the criterion for convergence in law from Proposition 2.27. So, let $a\in\mathbb R$ be a continuity point of the function $x\to P(f\leq x)$, and pick $\varepsilon>0$. By using the estimate found above, we have:
$$P(f_n\leq a)\leq P(f\leq a+\varepsilon)+P(|f_n-f|>\varepsilon)$$
$$P(f\leq a-\varepsilon)\leq P(f_n\leq a)+P(|f_n-f|>\varepsilon)$$

We conclude from this that we have the following estimate:
\begin{eqnarray*}
P(f\leq a+\varepsilon)-P(|f_n-f|>\varepsilon)
&\leq&P(f_n\leq a)\\
&\leq&P(f\leq a+\varepsilon)+P(|f_n-f|>\varepsilon)
\end{eqnarray*}

Now with $n\to\infty$ we obtain from this, in terms of $F(a)=P(x\leq a)$:
$$F(a-\varepsilon)\leq\lim_{n\to\infty}P(f_n\leq a)\leq F(a+\varepsilon)$$

Since $F$ was assumed to be continuous at $a$, with $\varepsilon\to 0$ this gives:
$$\lim_{n\to\infty}P(f_n\leq a)=P(f\leq a)$$

Thus we have indeed the convergence $f_n\to f$ in law, as stated. As for the counterexamples at the end, which are elementary, we will leave them as an exercise.
\end{proof}

\index{convergence almost surely}
\index{a.s.}
\index{almost surely}
\index{convergence a.s.}
\index{measure zero}

Moving on, as a third and final notion of convergence, we have:

\begin{definition}
We say that we have $f_n\to f$ almost surely when 
$$f_n(x)\to f(x)$$
outside a measure zero set.
\end{definition}

Again, there are many illustrations for this notion, and further things that can be said, about this. At the theoretical level, further buidling on Theorem 2.29, we have:

\begin{theorem}
We have the following implications,
$$\begin{pmatrix}f_n\to f\\ {\rm almost\ surely}\end{pmatrix}\implies
\begin{pmatrix}f_n\to f\\ {\rm in\ probability}\end{pmatrix}\implies
\begin{pmatrix}f_n\to f\\ {\rm in\ law}\end{pmatrix}$$
and the reverse implications do not necessarily hold.
\end{theorem}

\begin{proof}
The first implication is clear, because the measure zero set in Definition 2.30 can be assumed to be $\emptyset$. The second implication is something that we already know, from Theorem 2.29. As for the counterexamples, we will leave them as an exercise.
\end{proof}

Getting now to what we wanted to do, namely convergence of an average of i.i.d. variables towards their common mean, we have the following result, about this:

\index{law of large numbers}
\index{large numbers}

\begin{theorem}[Weak law of large numbers]
Given i.i.d. variables $f_1,f_2,f_3,\ldots\,$, with common mean $E(f_i)=E$, we have the following convergence in probability,
$$\frac{f_1+\ldots+f_n}{n}\to E$$
and as a consequence, this convergence holds as well in law.
\end{theorem}

\begin{proof}
This is something very standard, the idea being as follows:

\medskip

(1) Let us first establish the weaker result, namely the convergence in law, for our averages. Consider the sequence of averages of our variables:
$$g_n=\frac{f_1+\ldots+f_n}{n}$$

Now let us look at the Fourier transforms of our variables. We have, for any $i$, the following estimate, with $E=E(f_i)$ being the common expectation of our variables:
$$F_{f_i}(t)\simeq 1+iEt$$

By using independence, we deduce from this that we have, with $n\to\infty$:
$$F_{g_n}(t)
=\left[F_{f_i}\left(\frac{t}{n}\right)\right]^n
\simeq\left[1+iE\,\frac{t}{n}\right]^n
\simeq e^{iEt}$$

Now since the limiting function that we found $e^{iEt}$ is the Fourier transform of the constant variable $E$, we conclude from this that we have $g_n\to E$ in law.

\medskip

(2) In order to establish now the convergence in probability, as stated, we recall from chapter 1 that given a random variable $h:X\to\mathbb R$, having mean $E$ and variance $V$, the Chebycheff inequality states that we have, for any $a>0$:
$$P\big(|h-E|\geq a\big)\leq\frac{V}{a^2}$$

To be more precise, this was indeed something established in chapter 1, when doing discrete probability at that time, and the proof is identical in the general, continuous case. In order to prove now our result, consider the following sequence of variables:
$$h_n=\frac{f_1+\ldots+f_n}{n}-E$$

The mean and variance of these variables are then as follows, with $v=V(f_i)$:
$$E(h_n)=0\quad,\quad V(h_n)=\frac{v}{n}$$

Indeed, the mean formula $E(h_n)=0$ is clear, and by using this, we have as well:
\begin{eqnarray*}
V(h_n)
&=&E(h_n^2)\\
&=&E\left(\frac{\sum_{ij}f_if_j}{n^2}-\frac{2E\sum_if_i}{n}+E^2\right)\\
&=&E\left(\frac{\sum_if_i^2}{n^2}+\frac{\sum_{i\neq j}f_if_j}{n^2}-\frac{2E\sum_if_i}{n}+E^2\right)\\
&=&\frac{\sum_iE(f_i^2)}{n^2}+\frac{\sum_{i\neq j}E(f_i)E(f_j)}{n^2}-\frac{2E\sum_iE(f_i)}{n}+E^2\\
&=&\frac{v+E^2}{n}+n(n-1)\cdot\frac{E^2}{n^2}-2E^2+E^2\\
&=&\frac{v}{n}+\left(n+n(n-1)-n^2\right)\frac{E^2}{n^2}\\
&=&\frac{v}{n}
\end{eqnarray*}

Now let us apply the Chebycheff inequality above, to the variable $h_n$. This gives:
$$P(|h_n|\geq\varepsilon)\leq\frac{v/n}{\varepsilon^2}$$

Thus, we have the following inequality, in terms of the original variables $f_i$:
$$P\left(\left|\frac{f_1+\ldots+f_n}{n}-E\right|\geq\varepsilon\right)\leq\frac{v}{\varepsilon^2 n}$$

But with $n\to\infty$ this gives, as desired, the convergence in probability.
\end{proof}

The above result is not the end of the story with the law of large numbers, because, quite remarkably, the convergence there holds almost surely. In order to explain this, which is something quite technical, we will need some preliminaries. Let us start with:

\index{i.o.}
\index{infinitely often}

\begin{definition}
Given a sequence of events $A_1,A_2,A_3,\ldots\,$, we set
$$\limsup_nA_n=\bigcap_{m\in\mathbb N}\bigcup_{n=m}^\infty A_n$$
with the limit notation coming from the decreasing intersection. Equivalently, we set
$$\limsup_nA_n=\left\{w\in X\Big|w\in A_n\ {\rm i.o.}\right\}$$
with i.o. standing for infinitely often.
\end{definition}

To be more precise here, regarding the terminology that we use, given a sequence of events $A_1,A_2,A_3,\ldots\,$, we have a decreasing sequence of events, as follows:
$$\bigcup_{n=1}^\infty A_n\supset\bigcup_{n=2}^\infty A_n\supset\bigcup_{n=3}^\infty A_n\supset\ldots$$

Thus, we can consider the intersection of this sequence, and in analogy with what we know about numbers and limits, it makes sense to denote this intersection as above. Now with this notion in hand, we can state and prove a key technical result, as follows:

\index{Borel-Cantelli}

\begin{lemma}[Borel-Cantelli]
Given events $A_1,A_2,A_3,\ldots\,$, we have
$$\sum_{n=1}^\infty P(A_n)<\infty\implies P(w\in A_n\ {\rm i.o.})=0$$
with i.o. standing as usual for infinitely often.
\end{lemma}

\begin{proof}
We use the interpretation of $P(w\in A_n\ {\rm i.o.})$ coming from Definition 2.33. In that context, since the intersection there is decreasing, we have, for any $m\in\mathbb N$:
$$\limsup_nA_n\subset\bigcup_{n=m}^\infty A_n$$

We conclude from this that we have the following estimate, for any $m\in\mathbb N$:
$$P(\limsup_nA_n)\leq\sum_{n=m}^\infty P(A_n)$$

On the other hand, the assumption in the statement is that the series $\sum_{n=1}^\infty P(A_n)$ converges, and this can be stated in the following way:
$$\lim_{m\to\infty}\sum_{n=m}^\infty P(A_n)=0$$

Thus $P(\limsup_nA_n)$ is arbitrarily small, and we conclude that we have:
$$P(\limsup_nA_n)=0$$

But this is the same as saying that $P(w\in A_n\ {\rm i.o.})=0$, as stated.
\end{proof}

As a second technical ingredient, we will need the following straightforward generalization of the Chebycheff inequality from chapter 1, featuring an exponent $p$:

\index{Chebycheff inequality}

\begin{proposition}
We have the following Chebycheff inequality,
$$P\big(|g|\geq a\big)\leq\frac{E(g^p)}{a^p}$$
valid for any random variable $g:X\to\mathbb R$, and any $p\in\mathbb N$.
\end{proposition}

\begin{proof}
We have indeed the following computation, $\varphi$ being the density:
\begin{eqnarray*}
E(|g|^p)
&=&\int_\mathbb R|x|^p\varphi(x)dx\\
&\geq&a^p\int_{|x|\geq a}\varphi(x)dx\\
&=&a^pP(|g|\geq a)
\end{eqnarray*}

Thus, we are led to the Chebycheff inequality in the statement.
\end{proof}

We can now improve the weak law of large numbers, into something final, namely:

\index{strong law of large numbers}

\begin{theorem}[Strong law of large numbers]
Given i.i.d. variables $f_1,f_2,f_3,\ldots\,$, with common mean $E(f_i)=E$, we have the convergence
$$\frac{f_1+\ldots+f_n}{n}\to E$$
happening almost surely.
\end{theorem}

\begin{proof}
This is something quite technical, the idea being as follows:

\medskip

(1) By replacing the variables $f_i$ by their centered versions $f_i-E$, we can assume $E=0$. Consider now the sums of our variables, which have also zero mean:
$$g_n=f_1+\ldots+f_n$$

We would like to prove that the following convergence is almost sure:
$$\frac{g_n}{n}\to0$$

(2) In order to do this, consider an arbitrary event $w\in X$ such that:
$$\lim_{n\to\infty}\frac{g_n(w)}{n}\neq0$$

In this situation, there exists $\varepsilon>0$ such that, for infinitely many $n$:
$$\left|\frac{g_n(w)}{n}\right|>\varepsilon$$

Thus, in order to prove our theorem, we must show that we have, for any $\varepsilon>0$:
$$P(|g_n|>\varepsilon n\ {\rm i.o.})=0$$

(3) But for this latter purpose, we can use the Borel-Cantelli lemma, applied to:
$$A_n=\left\{w\in X\Big||g_n(w)|\geq n\varepsilon\right\}$$

Indeed, let us first verify that the convergence assumption in the Borel-Cantelli lemma is indeed satisfied. For this purpose, we will use the generalized Chebycheff inequality from Proposition 2.35, with exponent $p=4$. By using our assumption $E=0$, we have the following computation, with $M_2,M_4$ being the second and fourth moments of $f_i$:
\begin{eqnarray*}
E(g_n^4)
&=&E\left(\sum_{ijkl}f_if_jf_kf_k\right)\\
&=&E\left(\sum_if_i^4\right)+6E\left(\sum_{i<j}f_i^2f_j^2\right)\\
&=&nM_4+6\binom{n}{2}M_2^2\\
&=&nM_4+3n(n-1)M_2^2
\end{eqnarray*}

But this shows that for $n>>0$, we have a constant $C$ such that:
$$E(g_n^4)<Cn^2$$

(4) Thus, the generalized Chebycheff inequality from Proposition 2.35 applies to our situation, with exponent $p=4$, and gives the following estimate:
$$P(|g_n|\geq n\varepsilon)\leq\frac{E(g_n^4)}{(n\varepsilon)^4}\leq\frac{C}{\varepsilon^4n^2}$$

And with this, we are almost there. Indeed, if we denote by $N\in\mathbb N$ the smallest integer where the inequality found in (3) holds, we have:
$$\sum_{n\geq N}P(|g_n|\geq n\varepsilon)\leq\sum_{n\geq N}\frac{C}{\varepsilon^4n^2}<\infty$$

Thus, the Borel-Cantelli lemma applies, and shows that we have:
$$P(|g_n|>\varepsilon n\ {\rm i.o.})=0$$

But this ends the proof of our theorem, as explained in (2) above.
\end{proof}

\section*{2c. Central limits}

At a more advanced level now, on the same topics, you have certainly heard about bell-shaped curves, and perhaps even observed them in physics or chemistry class, because any routine measurement leads to such curves. Mathematically, here is the question that we would like to solve, coming as a refined version of our previous Question 2.24:

\index{CLT}
\index{Central Limit theorem}
\index{normal variable}
\index{Gaussian variable}

\begin{question}
Given i.i.d. random variables $f_1,f_2,f_3,\ldots\,$, assumed to be centered, and having common variance $t>0$, do we have
$$\frac{f_1+\ldots+f_n}{\sqrt{n}}\sim g_t$$
in the $n\to\infty$ limit, for some bell-shaped density $g_t$? And, what is the formula of $g_t$? 
\end{question}

Obviously, this is something quite tricky, notably with the above $1/\sqrt{n}$ factor needing some explanations. So, here is the idea. In order for our sums $S_n=f_1+\ldots+f_n$ to have a chance to converge, we must arrange our normalizations as for the expectation $E(S_n)$ not to blow up, and as for the variance $V(S_n)$ not to blow up either. And here:

\bigskip

(1) In general, the expectation requirement needs a $1/n$ normalization factor, and with this, we are into the law of large numbers, away from bell-shaped curves.

\bigskip

(2) So, we must waive the expectation requirement, and the way of doing this is by assuming that our variables are centered, $E(f_i)=0$, as said in Question 2.37. 

\bigskip

(3) Now with this done, we still have to find the correct normalization factor, making the variance $V(S_n)$ not to blow up. And here, we have the following computation:
$$V(S_n)
=\sum_iE(f_i^2)+\sum_{i\neq j}E(f_i)E(f_j)
=\sum_iE(f_i^2)
=nt$$

Thus the good normalization factor is indeed $1/\sqrt{n}$, as stated in Question 2.37.

\bigskip

In order to answer now Question 2.37, the straightforward approach would be via Fourier, and then the moment method for recapturing the limiting density. In practice, I mean trust me here, I did these computations, this leads to the following laws:

\index{normal law}

\begin{definition}
The normal law of parameter $1$ is the following measure:
$$g_1=\frac{1}{\sqrt{2\pi}}\,e^{-x^2/2}dx$$
More generally, the normal law of parameter $t>0$ is the following measure:
$$g_t=\frac{1}{\sqrt{2\pi t}}\,e^{-x^2/2t}dx$$
These are also called Gaussian distributions, with ``g'' standing for Gauss.
\end{definition}

Observe that the above laws have indeed mass 1, as they should, as shown by:
$$\int_\mathbb R e^{-x^2/2t}dx
=\int_\mathbb R e^{-y^2}\sqrt{2t}\,dy
=\sqrt{2t}\int_\mathbb R e^{-y^2}dy
=\sqrt{2\pi t}$$

Generally speaking, the normal laws appear as bit everywhere, in real life. The reasons behind this phenomenon come from the Central Limit Theorem (CLT), that we will explain in a moment, after developing some general theory. As a first result, we have:

\begin{proposition}
We have the variance formula
$$V(g_t)=t$$
valid for any $t>0$.
\end{proposition}

\begin{proof}
The first moment is 0, because our normal law $g_t$ is centered. As for the second moment, this can be computed by partial integration, as follows:
\begin{eqnarray*}
M_2
&=&\frac{1}{\sqrt{2\pi t}}\int_\mathbb Rx^2e^{-x^2/2t}dx\\
&=&\frac{1}{\sqrt{2\pi t}}\int_\mathbb R(tx)\left(-e^{-x^2/2t}\right)'dx\\
&=&\frac{1}{\sqrt{2\pi t}}\int_\mathbb Rte^{-x^2/2t}dx\\
&=&t
\end{eqnarray*}

We conclude from this that the variance is $V=M_2=t$, as stated.
\end{proof}

We can in fact compute all moments, by using the same trick, and we obtain:

\begin{theorem}
The even moments of the normal law are the numbers
$$M_k(g_t)=t^{k/2}\times k!!$$
where $k!!=(k-1)(k-3)(k-5)\ldots\,$, and the odd moments vanish. 
\end{theorem}

\begin{proof}
We have the following computation, valid for any integer $k\in\mathbb N$:
\begin{eqnarray*}
M_k
&=&\frac{1}{\sqrt{2\pi t}}\int_\mathbb Ry^ke^{-y^2/2t}dy\\
&=&\frac{1}{\sqrt{2\pi t}}\int_\mathbb R(ty^{k-1})\left(-e^{-y^2/2t}\right)'dy\\
&=&\frac{1}{\sqrt{2\pi t}}\int_\mathbb Rt(k-1)y^{k-2}e^{-y^2/2t}dy\\
&=&t(k-1)\times\frac{1}{\sqrt{2\pi t}}\int_\mathbb Ry^{k-2}e^{-y^2/2t}dy\\
&=&t(k-1)M_{k-2}
\end{eqnarray*}

Thus by recurrence, we are led to the formula in the statement.
\end{proof}

Regarding now the Fourier transform computation, this is as follows:

\begin{theorem}
We have the following formula, valid for any $t>0$:
$$F_{g_t}(x)=e^{-tx^2/2}$$
In particular, the normal laws satisfy $g_s*g_t=g_{s+t}$, for any $s,t>0$.
\end{theorem}

\begin{proof}
The Fourier transform formula can be established as follows:
\begin{eqnarray*}
F_{g_t}(x)
&=&\frac{1}{\sqrt{2\pi t}}\int_\mathbb Re^{-y^2/2t+ixy}dy\\
&=&\frac{1}{\sqrt{2\pi t}}\int_\mathbb Re^{-(y/\sqrt{2t}-\sqrt{t/2}ix)^2-tx^2/2}dy\\
&=&\frac{1}{\sqrt{2\pi t}}\int_\mathbb Re^{-z^2-tx^2/2}\sqrt{2t}dz\\
&=&\frac{1}{\sqrt{\pi}}e^{-tx^2/2}\int_\mathbb Re^{-z^2}dz\\
&=&\frac{1}{\sqrt{\pi}}e^{-tx^2/2}\cdot\sqrt{\pi}\\
&=&e^{-tx^2/2}
\end{eqnarray*}

As for the last assertion, this follows from the fact that $\log F_{g_t}$ is linear in $t$.
\end{proof}

We are now ready to state and prove the Central Limit Theorem, as follows:

\index{CLT}
\index{central limit}

\begin{theorem}[CLT]
Given i.i.d. random variables $f_1,f_2,f_3,\ldots\,$, assumed to be centered, and having common variance $t>0$, we have, with $n\to\infty$, in moments,
$$\frac{f_1+\ldots+f_n}{\sqrt{n}}\sim g_t$$
where $g_t$ is the Gaussian law of parameter $t$.
\end{theorem}

\begin{proof}
In terms of moments, the Fourier transform is given by:
$$F_f(x)
=E\left(\sum_{k=0}^\infty\frac{(ixf)^k}{k!}\right)
=\sum_{k=0}^\infty\frac{(ix)^kE(f^k)}{k!}
=\sum_{k=0}^\infty\frac{i^kM_k(f)}{k!}\,x^k$$

We conclude that the Fourier transform of the variable in the statement is:
$$F(x)
=\left[F_f\left(\frac{x}{\sqrt{n}}\right)\right]^n
\simeq\left[1-\frac{tx^2}{2n}\right]^n
\simeq e^{-tx^2/2}$$

But this latter function being the Fourier transform of $g_t$, we obtain the result.
\end{proof}

With this discussed, what is next? Normally, some more mathematics for the normal laws, in relation with their moments, as we usually do, when meeting a new law. And we first have here, coming as an interesting reformulation of Theorem 2.40:

\begin{theorem}
The moments of the normal law are the numbers
$$M_k(g_t)=t^{k/2}|P_2(k)|$$
where $P_2(k)$ is the set of pairings of $\{1,\ldots,k\}$.
\end{theorem}

\begin{proof}
Let us count the pairings of $\{1,\ldots,k\}$. In order to have such a pairing, we must pair $1$ with one of the numbers $2,\ldots,k$, and then use a pairing of the remaining $k-2$ numbers. Thus, we have the following recurrence formula:
$$|P_2(k)|=(k-1)|P_2(k-2)|$$

As for the initial data, this is $P_1=0$, $P_2=1$. Thus, we are led to the result.
\end{proof}

We are not done yet with moments, here being one more version of what we have:

\begin{theorem}
The moments of the normal law are the numbers
$$M_k(g_t)=\sum_{\pi\in P_2(k)}t^{|\pi|}$$
where $P_2(k)$ is the set of pairings of $\{1,\ldots,k\}$, and $|.|$ is the number of blocks.
\end{theorem}

\begin{proof}
This follows indeed from Theorem 2.43, because the number of blocks of a pairing of $\{1,\ldots,k\}$ is trivially $k/2$, independently of the pairing.
\end{proof}

The above result is quite interesting, making a link with the formula of the moments of the Poisson law from chapter 1, which was very similar, as follows:
$$M_k(p_t)=\sum_{\pi\in P(k)}t^{|\pi|}$$

And, more on such things later. Now back to our regular business, various moment manipulations, as our next and final result on the subject, let us formulate:

\begin{theorem}
The even normalized central moments of the normal law are
$$M_k''(g_t)=k!!$$
with $k!!=(k-1)(k-3)\ldots\,$ as usual, and the odd such moments vanish. Also,
$$E=0\quad,\quad V=t\quad,\quad \gamma=0\quad,\quad \kappa=3$$
are the expectation, variance, skewness and kurtosis.
\end{theorem}

\begin{proof}
The normal laws being centered, their central moments equal the plain ones, $M_k'=M_k$. Since the variance is $V=t$, the normalized central moments are given by:
$$M_k''=\frac{M'_k}{t^{k/2}}=\frac{M_k}{t^{k/2}}=\frac{\delta_{2|k}t^{k/2}\times k!!}{t^{k/2}}=\delta_{2|k}k!!$$

Finally, at $k=3$ this gives $\gamma=0$, and at $k=4$ this gives $\kappa=4!!=3$.
\end{proof}

\section*{2d. Shifted bells}

The laws $g_t$ that we studied in the previous section appeared there in a purely mathematical way, from central limits, and as the word ``central'' indicates, they are of course centered. However, this was theory, and in practice, say coming via various real-life measurements, the Gaussian bells quite often appear shifted by a parameter $a\in\mathbb R$, I mean with symmetry axis there at $a$. So, we must study these versions too, and we have:

\index{shifted normal law}

\begin{theorem}
The normal law $g_t$, shifted to the right by $a\in\mathbb R$, has density
$$g_t^a=\frac{1}{\sqrt{2\pi t}}\,e^{-(x-a)^2/2t}\,dx$$
and the basic properties of this shifted law $g_t^a$ are as follows:
\begin{enumerate}
\item We have $f\sim g_t\implies f+a\sim g_t^a$.

\item The mean is $E=a$, the variance is $V=t$.

\item The moments are $M_k=\sum_{r=0}^{[k/2]}\binom{k}{2r}(2r)!!t^ra^{k-2r}$.

\item The moments satisfy $M_k=aM_{k-1}+(k-1)tM_{k-2}$.

\item The central moments are $M_k'=\delta_{2|k}t^{k/2}k!!$.

\item The normalized central moments are $M_k''=\delta_{2|k}k!!$.

\item The skewness is $\gamma=0$, the kurtosis is $\kappa=3$.

\item The Fourier transform is $F(x)=\exp(iax-tx^2/2)$.

\item We have the convolution formula $g_s^a*g_t^b=g_{s+t}^{a+b}$.
\end{enumerate}
\end{theorem}

\begin{proof}
These are things that we know well at $a=0$, from the previous section, and the general case basically follows from this, with a few computations added. To start with, let us draw a picture, that will help us understanding what is going on:
$$
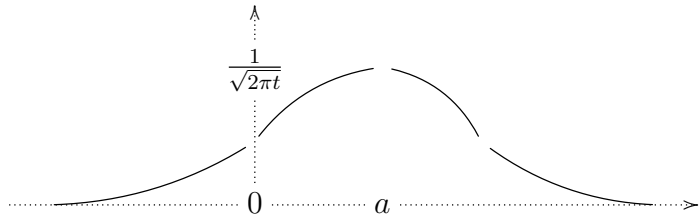
\xymatrix@R=13pt@C=10pt{
&&&&&\\
&&&&&\frac{1}{\sqrt{2\pi t}}\ar@.[u]&&\ar@/^/@{-}[drr]\ar@/_/@{-}[dll]&\\
&&&&&\ar@/^/@{-}[dllll]&&&&\ar@/_/@{-}[drrrr]&\\
\ar@{.}[rrrrr]&&&&&0\ar@{.}[uu]\ar@{.}[rr]&&a\ar@.[rrrrrrr]&&&&&&&}$$

(1) This is actually the definition of $g_t^a$, and the density formula in the statement comes from the following computation, assuming $f\sim g_t$, valid for any $\varphi:\mathbb R\to\mathbb R$:
\begin{eqnarray*}
E(\varphi(f+a))
&=&\frac{1}{\sqrt{2\pi t}}\int_\mathbb R\varphi(x+a)e^{-x^2/2t}\,dx\\
&=&\frac{1}{\sqrt{2\pi t}}\int_\mathbb R\varphi(x)e^{-(x-a)^2/2t}\,dx
\end{eqnarray*}

Indeed, the density of $f+a\sim g_t^a$ follows to be $1/\sqrt{2\pi t}\,e^{-(x-a)^2/2t}\,dx$, as claimed.

\medskip

(2) Regarding the mean, that is clearly $E=a$ on the picture, because the symmetry axis lies there, at $a\in\mathbb R$. As for the variance, by partial integration, we have:
\begin{eqnarray*}
M_2-aM_1
&=&\frac{1}{\sqrt{2\pi t}}\int_\mathbb R(x^2-ax)e^{-(x-a)^2/2t}dx\\
&=&\frac{1}{\sqrt{2\pi t}}\int_\mathbb R(tx)\left(-e^{-(x-a)^2/2t}\right)'dx\\
&=&\frac{1}{\sqrt{2\pi t}}\int_\mathbb Rte^{-(x-a)^2/2t}dx\\
&=&t
\end{eqnarray*}

Thus $M_2=a^2+t$, and the variance is $V=(a^2+t)-a^2=t$, as stated.

\medskip

(3) We have indeed the following computation, using a variable $f\sim g_t$:
\begin{eqnarray*}
M_k
&=&E((f+a)^k)\\
&=&\sum_{r=0}^{[k/2]}\binom{k}{2r}E(f^{2r})a^{k-2r}\\
&=&\sum_{r=0}^{[k/2]}\binom{k}{2r}(2r)!!t^ra^{k-2r}
\end{eqnarray*}

(4) By using the same trick as in (2), we obtain the following recurrence:
\begin{eqnarray*}
M_k-aM_{k-1}
&=&\frac{1}{\sqrt{2\pi t}}\int_\mathbb R(x^k-ax^{k-1})e^{-(x-a)^2/2t}dx\\
&=&\frac{1}{\sqrt{2\pi t}}\int_\mathbb R(tx^{k-1})\left(-e^{-(x-a)^2/2t}\right)'dx\\
&=&\frac{1}{\sqrt{2\pi t}}\int_\mathbb R(k-1)tx^{k-2}\cdot e^{-(x-a)^2/2t}dx\\
&=&(k-1)tM_{k-2}
\end{eqnarray*}

In practice, here is the list of the first few moments, obtained in this way:
$$M_1=a$$
$$M_2=a^2+t$$
$$M_3=a^3+3at$$
$$M_4=a^4+6a^2t+3t^2$$
$$M_5=a^5+10a^3t+15at^2$$
$$M_6=a^6+15a^4t+45a^2t^2+15t^3$$
$$M_7=a^7+21a^5t+105a^3t^2+105at^3$$
$$M_8=a^8+28a^6t+210a^4t^2+420a^2t^3+105t^4$$

(5) Let us look now at the central moments. And here, fortunately, we obtain the same moments as for $g_t$, as shown by the following computation:
\begin{eqnarray*}
M_k'(g_t^a)
&=&\frac{1}{\sqrt{2\pi t}}\int_\mathbb R(x-a)^ke^{-(x-a)^2/2t}dx\\
&=&\frac{1}{\sqrt{2\pi t}}\int_\mathbb Rx^ke^{-x^2/2t}dx\\
&=&M_k(g_t)
\end{eqnarray*}

Thus, by using now the formula for the moments of $g_t$, we have:
$$M_k'(g_t^a)=\delta_{2|k}t^{k/2}k!!$$

(6) Next, since the central moments are the same for $g_t^a$ and $g_t$, and since the variances are also the same, the normalized central moments must all coincide:
$$M_k''(g^t_a)=M_k''(g_t)$$

(7) In particular, we deduce from this formula that the skewness is $\gamma=0$, and that the kurtosis is $\kappa=3$, exactly as for the laws $g_t$ before, as claimed.

\medskip

(8) Regarding now the Fourier transform computation, we have:
\begin{eqnarray*}
F_{g_t^a}(x)
&=&\frac{1}{\sqrt{2\pi t}}\int_\mathbb Re^{-(y-a)^2/2t+ixy}dy\\
&=&\frac{1}{\sqrt{2\pi t}}\int_\mathbb Re^{-y^2/2t+ix(y+a)}dy\\
&=&\frac{e^{iax}}{\sqrt{2\pi t}}\int_\mathbb Re^{-y^2/2t+ixy}dy\\
&=&e^{iax}F_{g_t}(x)
\end{eqnarray*}

Now by using the formula of the Fourier transform of $g_t$, we get:
$$F_{g_t^a}(x)=e^{iax}e^{-tx^2/2}=e^{iax-tx^2/2}$$
 
(9) Observe now that this Fourier transform that we just computed satisfies:
\begin{eqnarray*}
F_{g_s^a}(x)F_{g_t^b}(x)
&=&e^{iax-sx^2/2}e^{ibx-tx^2/2}\\
&=&e^{i(a+b)x-(s+t)x^2/2}\\
&=&F_{g_{s+t}^{a+b}}(x)
\end{eqnarray*}

Thus, at the level of the corresponding laws, we have the following formula:
$$g_s^a*g_t^b=g_{s+t}^{a+b}$$

And with this, done, we are led to the conclusions in the statement.
\end{proof}

In practice now, many other things can be said about the normal laws $g_t$ and their shifted versions $g_t^a$, and we will be back to this, in what follows. Needless to say, such computations are very useful, when dealing with real-life probability questions.

\bigskip

As a conclusion now to what we did, a bit philosophical, focusing on the difficulties that we met, rather than on our successes, let us formulate things as follows:

\begin{conclusion}
The theory of $g_t^a$ is quite straightforward, save for:
\begin{enumerate}
\item The computation of the CDF, which amounts in fighting with $\int e^{-x^2}$.

\item The moments at $a\neq0$, most likely a hypergeometric function business.
\end{enumerate}
\end{conclusion}

And with of course, exercise for you to learn more, about all this. In fact, thinking a bit, for various real-life problems involving the laws $g_t^a$, the CDF is what you need.

\section*{2e. Exercises}

This chapter was a quite standard introduction to continuous phenomena in probability, and as exercises on all this, chosen rather computational, we have:

\begin{exercise}
Learn if needed the proof of the Gauss formula, $\int_\mathbb Re^{-x^2}dx=\sqrt{\pi}$.
\end{exercise}

\begin{exercise}
Learn as well the Fresnel formula, $\int_\mathbb R e^{it^2}dt=\sqrt{\pi i}$.
\end{exercise}

\begin{exercise}
Come across the normal laws via Fourier and the moment method.
\end{exercise}

\begin{exercise}
Learn more about the precise convergence in the CLT.
\end{exercise}

\begin{exercise}
Try to unify our moment formulae for $p_t$ and $g_t$.
\end{exercise}

\begin{exercise}
Further experiment with the moments of $g_t^a$.
\end{exercise}

\begin{exercise}
Prove $g_s^a*g_t^b=g_{s+t}^{a+b}$ with bare hands, without Fourier.
\end{exercise}

\begin{exercise}
Study a bit the approximation question, for the CDF of $g_t$.
\end{exercise}

As bonus exercise, in relation with the axiomatics, have a look at some standard books in measure theory or probability, such as Durrett \cite{dur}, Feller \cite{fel} or Rudin \cite{rud}.

\chapter{Advanced formulae}

\section*{3a. The seven laws}

We have seen some interesting theory for the Poisson laws $p_t$, which are the central laws in discrete probability, and for the normal laws $g_t$, which are the central laws in continuous probability. In this chapter we would like to develop some more theory for $p_t,g_t$, with formulae regarding a number of more specialized objects associated to them, namely Cauchy and Stieltjes transforms, Hankel determinants and orthogonal polynomials.

\bigskip

This will be something quite technical. To start with, the results about $p_t,g_t$ that we would like to talk about do not come alone, and are best understood in a larger context, that of a family of 7 main probability laws, according to the following principle:

\begin{principle}
The advanced theory of $p_t,g_t$ is best understood in the context of
$$\{p_t,g_t,e_t,u_t,\gamma_t,\alpha_t,\pi_t\}$$
the family formed by the main 7 probability laws.
\end{principle}

So, let us first talk about these laws. In order to do so, the best is to think in terms of support. The Poisson law $p_t$ is main law supported by $\mathbb N$, and the normal law $g_t$ is the main law supported by $\mathbb R$. As a next addition to our family, we have the exponential law $e_t$, which is the main law supported by $[0,\infty)$, whose basic theory is as follows:

\index{exponential law}

\begin{theorem}
The exponential law $e_t$ of parameter $t>0$, having density
$$e_t=te^{-tx}\,dx$$
on $[0,\infty)$, has the following properties:
\begin{enumerate}
\item The mean is $E=1/t$, the variance is $V=1/t^2$.

\item The moments are $M_k=k!/t^k$.

\item The central moments are $M_k'=k!/t^k\sum_{s=0}^k(-1)^s/s!$.

\item The normalized central moments are $M_k''=\sum_{s=0}^k(-1)^sk!/s!$.

\item The skewness is $\gamma=2$, the kurtosis is $\kappa=9$.

\item The Fourier transform is $F(x)=t/(t-ix)$.
\end{enumerate}
\end{theorem}

\begin{proof}
It is beyond our scope, here in this book, to have a detailed look at $e_t$, in the spirit of what we did before for $p_t,g_t$, with this being rather an auxiliary law, in our context. However, here is a brief survey of what can be said about $e_t$, going beyond the moment and Fourier essentials in the statement, that we will need in what follows:

\medskip

(1) To start with, we can certainly talk about $e_t=te^{-tx}\,dx$ with $t>0$, supported on $[0,\infty)$, with the mass one property coming from $\int_0^\infty e^{-x}dx=1$.

\medskip

(2) In practice, this law appears in many contexts. For instance in the context of particle decay, if we denote by $t>0$ the decay rate, that is, the probability per unit time that the particle will disintegrate, then the probability decay function is $e_t$.

\medskip

(3) Mathematically now, as a first observation, the CDF of the exponential law $e_t$ is computable, and given by a very simple formula, namely $P(f\leq x)=1-e^{-tx}$.

\medskip

(4) Getting now to the moments, by partial integration these are given by the formula in the statement, namely $M_k=k!/t^k$, as shown by the following computation: 
\begin{eqnarray*}
M_k
&=&t\int_0^\infty x^ke^{-tx}dx\\
&=&t\int_0^\infty kx^{k-1}\cdot\frac{e^{-tx}}{t}\,dx\\
&=&\frac{k}{t}\cdot t\int_0^\infty x^{k-1}e^{-tx}dx\\
&=&\frac{k}{t}\cdot M_{k-1}
\end{eqnarray*}

In particular the mean is $E=1/t$, and the variance is $V=(2/t^2-1/t^2)=1/t^2$.

\medskip

(5) Regarding the central moments, which in general are meant to simplify the moment formulae, in our case the moments are given by a simple formula, and there is nothing to simplify, and the central moments are in fact given by a quite complicated formula, namely $M_k'=k!/t^k\sum_{s=0}^k(-1)^s/s!$. An interesting peculiarity of the exponential law.

\medskip

(6) Next, the normalized central moments are given by $M_k''=\sum_{s=0}^k(-1)^sk!/s!$, and in particular, the skewness is $\gamma=2$, and the kurtosis is $\kappa=9$.

\medskip

(7) Regarding the Fourier transform, this is $F(x)=t/(t-ix)$. Observe that, unlike for $p_t,g_t$, the function $\log F$ is not linear in $t$, that is, $e_s*e_t\neq e_{s+t}$. However, this time in analogy with $p_t,g_t$, it is possible to prove that $e_t$ is infinitely divisible.

\medskip

(8) Finally, in relation with $P(f>x)=e^{-tx}$, an exponential variable $f:X\to\mathbb R$ experiences loss of memory, $P(f>s+t|f>s)=P(f>t)$. In fact, the exponential laws are the only continuous laws having this property. As for the discrete laws having this property, these are the geometric laws $c_p$ with $p\in[0,1]$, given by $P(s)=(1-p)^sp$.

\medskip

(9) So, this was the story with the exponential laws, get to know them better, and you have here for instance my book \cite{ba2}, where all this is explained in detail.
\end{proof}

Moving on, let us go back to the support considerations made after Principle 3.1. Now that we have main laws $p_t,g_t,e_t$ supported respectively by $\mathbb N,\mathbb R,[0,\infty)$, we need a compactly supported law too, in our collection, and as obvious choice here, we have:

\index{uniform law}

\begin{theorem}
The uniform law $u_t$ of parameter $t>0$, having density
$$u_t=\frac{1}{t}\,dx$$
on $[0,t]$, has the following properties:
\begin{enumerate}
\item The mean is $E=t/2$, the variance is $V=t^2/12$.

\item The moments are $M_k=t^k/(k+1)$.

\item The central moments are $M_k'=\delta_{2|k}(t/2)^k/(k+1)$.

\item The normalized central moments are $M_k''=\delta_{2|k}3^{k/2}/(k+1)$.

\item The skewness is $\gamma=0$, the kurtosis is $\kappa=9/5$.

\item The Fourier transform is $F(x)=(e^{itx}-1)/itx$.
\end{enumerate}
\end{theorem}

\begin{proof}
We can certainly talk about $u_t$, as a probability measure, the mean is $E=t/2$, then the moments are $M_k=t^k/(k+1)$, and so the variance is $V=t^2/3-t^2/4=t^2/12$. The central moments are easy to compute too, given by $M_k'=\delta_{2|k}(t/2)^k/(k+1)$, and by dividing by powers of $\sigma=t/(2\sqrt{3})$ we obtain $M_k''=\delta_{2|k}3^{k/2}/(k+1)$. In particular we have $\gamma=0$ and $\kappa=9/5$. Finally, the formula $F(x)=(e^{itx}-1)/itx$ is clear too.
\end{proof}

What is next? I don't know about you, but personally I find that in our family $\{p_t.g_t,e_t,u_t\}$ the weak member is $u_t$, I mean, that ain't no serious probability law. So, we need more compactly supported laws, and here, as an obvious choice, we have:

\index{semicircle law}
\index{Wigner law}

\begin{theorem}
The semicircle law $\gamma_t$ of parameter $t>0$, having density
$$\gamma_t=\frac{1}{2\pi t}\sqrt{4t-x^2}\,dx$$
on $[-2\sqrt{t},2\sqrt{t}]$, and called Wigner law, has the following properties:
\begin{enumerate}
\item The mean is $E=0$, the variance is $V=t$.

\item The even moments are $M_{2k}=\frac{t^k}{k+1}\binom{2k}{k}$.

\item The even central moments are $M_{2k}'=\frac{t^k}{k+1}\binom{2k}{k}$.

\item The even normalized central moments are $M_{2k}''=\frac{1}{k+1}\binom{2k}{k}$.

\item The skewness is $\gamma=0$, the kurtosis is $\kappa=2$.

\item The moment generating function is $M(z)=(1-\sqrt{1-4tz^2})/2tz^2$.
\end{enumerate}
\end{theorem}

\begin{proof}
Many things can be said here, the idea being as follows:

\medskip

(1) To start with, we can talk about the semicircle law over any interval $[a,b]$, the formula being as follows, based on the fact that the area of the unit circle is $\pi$:
$$\gamma_{ab}=\frac{8}{\pi(b-a)^2}\sqrt{(x-a)(b-x)}dx$$

(2) The problem is, how to choose a winner $\gamma_1$ among these laws. And then, once this winner found, how to correctly parametrize the family $\{\gamma_t|t>0\}$ it comes from.

\medskip

(3) In answer, let us count loops on $\mathbb N$, based at $0$. There are no loops of odd length, and the number of loops of length $2k$ is the $k$-th Catalan number $\frac{1}{k+1}\binom{2k}{k}$. But these numbers are precisely the moments of the semicircle law on $[-2,2]$, namely:
$$\gamma_1=\frac{1}{2\pi}\sqrt{4-x^2}\,dx$$

So, based on this loop business, this will be our winner. It is of course possible to say more here, by talking about random walks on $\mathbb N$, but let us not get further into this.

\medskip

(4) Alternatively, $\gamma_1$ appears as such from free probability, as being the free analogue of the normal law $g_1$, or from random matrices, as being the asymptotic law of the Wigner matrices, or from Lie groups, as being the main character law of $SU_2$, or from quantum groups, as being the main character law of $O_N^+$ with $N\geq2$. More on all this later.

\medskip

(5) Regarding now the correct scaling, the random walk approach in (3) does not really provide a clue, and nor does the $SU_2$ approach in (4). However, the free probability, random matrix and quantum group approaches in (4) do provide an answer, which fortunately is a common answer, leading to the laws in the statement, namely: 
$$\gamma_t=\frac{1}{2\pi t}\sqrt{4t-x^2}\,dx$$

(6) Next, we already talked about the moments of $\gamma_1$, as being $M_{2k}=\frac{1}{k+1}\binom{2k}{k}$, and it follows that for $\gamma_t$ we have $M_{2k}=\frac{t^k}{k+1}\binom{2k}{k}$. In particular,  $E=0$ and $V=t$.

\medskip

(7) Our laws being centered we have $M_{2k}'=M_{2k}$, and then by dividing by powers of $\sigma=\sqrt{t}$ we obtain $M_{2k}''=\frac{1}{k+1}\binom{2k}{k}$. In particular, we have $\gamma=0$ and $\kappa=2$.

\medskip

(8) Next, the Fourier transform of $\gamma_t$ is something quite complicated, and a bit irrelevant too, because as explained in (4,5) our measures are of ``free'' nature, and so having nothing to do with Fourier, which deals with independence, in the classical sense.

\medskip

(9) This being said, as a matter of doing some related mathematics, we can try to compute instead the moment generating function $M(z)=\sum_{k\geq0}M_kz^k$. And here, by using the generalized binomial formula with exponent $p=1/2$, we have:
$$M(z)=\sum_{k=0}^\infty\frac{1}{k+1}\binom{2k}{k}(tz^2)^k
=\frac{1-\sqrt{1-4tz^2}}{2tz^2}$$

(10) So, this was for the story with the Wigner semicircle laws, save for some advanced phenomenology, mentioned in (4,5), to be discussed later in this book, and for some standard computations involving loops on $\mathbb N$, Catalan numbers $C_k=\frac{1}{k+1}\binom{2k}{k}$, the semicircle density, and the formula of $\sqrt{1-x}$, that you can find for instance in my book \cite{ba2}.
\end{proof}

With this discussed, what is next? Comparing I guess our two compactly supported measures, $u_t$ and $\gamma_t$, and here, there is something interesting to say, as follows:

\index{beta distributions}
\index{beta function}

\begin{comment}
Up to shifts, scalings and mass 1 normalizations, our laws $u_t,\gamma_t$ are
$$u=1\quad,\quad\gamma=\sqrt{x(1-x)}$$
on $[0,1]$, and this suggests looking at the following laws on $[0,1]$, with $a,b>-1$:
$$\mu=x^a(1-x)^b$$
For instance, we can try to complement $u,\gamma$, coming from $a=b=0$ and $a=b=\frac{1}{2}$, with 
$$\alpha=\frac{1}{\sqrt{x(1-x)}}\quad,\quad\pi=\sqrt{\frac{x}{1-x}}$$
coming from the next simplest parameters, namely $a=b=-\frac{1}{2}$, and $a=\frac{1}{2},b=-\frac{1}{2}$.
\end{comment}

Obviously, this is something quite subjective, but in our situation now, lost in the mysteries of probability theory, a bit of subjectivity and common sense is what we need. Here are a few more comments, on all this, further clarifying the above:

\bigskip

(1) To start with, in view of the subtleties with the parametrization of the semicircle laws, which retrospectively bring some potential questions in relation with our parametrization of the uniform laws too, we decided to give up with the parametrization problematics, by having everything on $[0,1]$. More on this later, when we'll be wiser.

\bigskip

(2) Next, when looking at $\mu=x^a(1-x)^b$, the integral converges at $a,b>-1$, but is generically not computable. To be more precise, with $a=c-1$, $b=d-1$, there is a formula for this integral, called beta function of $b,d$, in terms of the gamma function:
$$\int_0^1x^{c-1}(1-x)^{d-1}dx=\frac{\Gamma(c)\Gamma(d)}{\Gamma(c+d)}$$

And the point is that the gamma function itself is generically not computable. And more on such things, which can be quite complicated, later in this book.

\bigskip

(3) In practice, we are left with looking for simple values of $a,b>-1$, that we can use. And here, with $a,b\in\{-\frac{1}{2},0,\frac{1}{2}\}$, leaving aside $u,\gamma$, that we already have, and leaving aside as well the cases $a=0,b\neq0$ and $a\neq 0,b=0$, leading nowhere, we are left, up to symmetry, with $a=b=-\frac{1}{2}$ producing $\alpha$, and $a=\frac{1}{2},b=-\frac{1}{2}$ producing $\pi$.

\bigskip

So, this was for the story with Comment 3.5, hope that I convinced you a bit, and that I get some credit for my literary skills. Getting to work now, let us first study $\alpha$ and see what we get, is that bound for interesting theorems, or for the trash can. And fortunately, we are led in this way to an interesting theorem, as follows:

\index{arcsine law}

\begin{theorem}
The arcsine law $\alpha_t$ of parameter $t>0$, having density
$$\alpha_t=\frac{1}{\pi\sqrt{4t-x^2}}\,dx$$
on $[-2\sqrt{t},2\sqrt{t}]$, has the following properties:
\begin{enumerate}
\item The mean is $E=0$, the variance is $V=2t$.

\item The even moments are $M_{2k}=t^k\binom{2k}{k}$.

\item The even central moments are $M_{2k}'=t^k\binom{2k}{k}$.

\item The even normalized central moments are $M_{2k}''=\frac{1}{2^k}\binom{2k}{k}$.

\item The skewness is $\gamma=0$, the kurtosis is $\kappa=3/2$.

\item The moment generating function is $M(z)=1/\sqrt{1-4tz^2}$.
\end{enumerate}
\end{theorem}

\begin{proof}
The story here is a bit similar to that of the semicircle laws, as follows:

\medskip

(1) To start with, we can talk about the arcsine law over any interval $[a,b]$:
$$\alpha_{ab}=\frac{1}{\pi\sqrt{(x-a)(b-x)}}\,dx$$

Indeed, the mass 1 property comes from the following computation, using the linear change of variables $x\to y$ mapping $[a,b]\to[-1,1]$, which makes $a,b$ dissapear:
$$\int_a^b\frac{1}{\sqrt{(x-a)(b-x)}}\,dx=\int_{-1}^1\frac{1}{\sqrt{1-y^2}}\,dy=\Big[\arcsin y\Big]_{-1}^1=\pi$$

(2) As before with the semicircle laws, we must choose a winner $\alpha_1$ among these laws. Now for this purpose, let us count loops on $\mathbb Z$, based at $0$. There are no loops of odd length, and the number of loops of length $2k$ is the $k$-th middle binomial $\binom{2k}{k}$. But these numbers are precisely the moments of the arcsine law on $[-2,2]$, namely:
$$\alpha_1=\frac{1}{\pi\sqrt{4-x^2}}\,dx$$

(3) So, based on this loop business, this will be our winner. It is of course possible to say more here, by talking about random walks on $\mathbb Z$, and also about certain Lie groups, but let us not further get into this. As for the parametric versions of our law $\alpha_1$, inspired by what happens for the semicircle laws, we will dilate the support by $\sqrt{t}$:
$$\alpha_t=\frac{1}{\pi\sqrt{4t-x^2}}\,dx$$

(4) Next, we already talked about the moments of $\alpha_1$, as being $M_{2k}=\binom{2k}{k}$, and it follows that for $\alpha_t$ we have $M_{2k}=t^k\binom{2k}{k}$. In particular,  $E=0$ and $V=2t$.

\medskip

(5) Our laws being centered we have $M_{2k}'=M_{2k}$, and then by dividing by powers of $\sigma=\sqrt{2t}$ we obtain $M_{2k}''=\frac{1}{2^k}\binom{2k}{k}$. In particular, we have $\gamma=0$ and $\kappa=3/2$.

\medskip

(6) Finally, for the same reasons as for semicircle laws, namely ``underlying freeness'', that we will learn about later in this book, we will be not interested in the Fourier transform. Instead, we can compute moment generating function $M(z)=\sum_{k\geq0}M_kz^k$. And here, by using the binomial formula with exponent $p=-1/2$, we have:
$$M(z)=\sum_{k=0}^\infty\binom{2k}{k}(tz^2)^k
=\frac{1}{\sqrt{1-4tz^2}}$$

(7) So, this was for the story with the arcsine laws, save for some standard computations involving loops on $\mathbb Z$, central binomial coefficients $D_k=\binom{2k}{k}$, the arcsine density, and the formula of $1/\sqrt{1-x}$, that you can find for instance in my book \cite{ba2}.
\end{proof}

As a last task, let us investigate the law $\pi$ from Comment 3.5. We have here:

\index{Marchenko-Pastur law}

\begin{theorem}
The Marchenko-Pastur law $\pi_t$ of parameter $t>0$, having density
$$\pi_t=\max(1-t,0)\delta_0+\frac{\sqrt{4t-(x-1-t)^2}}{2\pi x}\,dx$$
on $[(1-\sqrt{t})^2,(1+\sqrt{t})^2]$, has the following properties:
\begin{enumerate}
\item At $t=1$, the density is $\pi_1=\frac{1}{2\pi}\sqrt{4x^{-1}-1}\,dx$ on $[0,4]$.

\item At $t=1$, the moments are $M_k=\frac{1}{k+1}\binom{2k}{k}$.

\item The mean is $E=t$, the variance is $V=t$.

\item The moments are $M_k=\sum_{\pi\in NC(k)}t^{|\pi|}$.

\item The skewness is $\gamma=1/\sqrt{t}$, the kurtosis is $\kappa=2+1/t$.

\item $M(z)=(1-t)/2+(1-\sqrt{1-2(1+t)z+(1-t)^2z^2})/2z$.
\end{enumerate}
\end{theorem}

\begin{proof}
This is again a bit similar to what we did before for the semicircle law, but with some significant new twists on the way, the idea being as follows:

\medskip

(1) As per our general philosophy, coming from Comment 3.5, we can talk about measures of type $\pi_{ab}$ over any interval $[a,b]$. And what we have to do is to select a winner $\pi_1$, and then find its correct parametric versions $\{\pi_t|t>0\}$, presumably among $\{\pi_{ab}\}$.

\medskip

(2) In order to find the winner, we will look, as before in the context of the semicircle laws, at the loops on $\mathbb N$, based at $0$. We already know that there are no loops of odd length, and the number of loops of length $2k$ is the $k$-th Catalan number $\frac{1}{k+1}\binom{2k}{k}$. 

\medskip

(3) And here comes the point. Since there are no loops of odd length, we can totally ignore them, I mean ignore the 0 values coming from their count, and look for the measure having the Catalan numbers $\frac{1}{k+1}\binom{2k}{k}$ as moments. And this is the following law:
$$\pi_1=\frac{1}{2\pi}\sqrt{4x^{-1}-1}\,dx$$

(4) So, based on this loop business, which appears as a variation of what we did before for the semicircle laws, this will be our winner. It is of course possible to say more here, by talking about random walks on $\mathbb N$, but let us not get further into this. Instead, let us record the following fact, which can stand as an alternative definition for $\pi_1$:
$$f\sim\gamma_1\implies f^2\sim\pi_1$$

(5) Alternatively, $\pi_1$ appears from free probability, as free analogue of the Poisson law $p_1$, or from random matrices, as asymptotic law of the Wishart matrices, or from Lie groups, as the law of $\chi+1$, with $\chi$ being the main character law of $SO_3$, or from quantum groups, as the law of the main character of $S_N^+$ with $N\geq4$. More on all this later.

\medskip

(6) Getting now to parametric versions, things here are quite complicated, with an atom at $0$ appearing at $t<1$, for some reasons which are not easy to explain. The idea with all this is that the density in the statement comes from free probability, with $\pi_t$ being the free analogue of $p_t$, or from Wishart matrices, or from quantum groups.

\medskip

(7) To be more precise here, let us first look at moments. Our previous experience with $g_t\to\gamma_t$ suggests that ``when liberating, the crossing partitions dissapear''. Which is something that happens indeed, ultimately coming from the fact that the crossing $\slash\hskip-2.1mm\backslash$ corresponds to commutation, $ab=ba$, as we will later learn in this book. 

\medskip

(8) Now in our context, regarding the operation $p_t\to\pi_t$, the situation is as follows:
$$M_k(p_1)=|P(k)|\quad,\quad M_k(\pi_1)=|NC(k)|$$
$$M_k(p_t)=\sum_{\pi\in P(k)}t^{|\pi|}
\quad,\quad M_k(\pi_t)=?$$

But with this in hand, it is a no-brainer to conjecture the following formula:
$$M_k(\pi_t)=\sum_{\pi\in NC(k)}t^{|\pi|}$$

And the point is that this is what happens, indeed, with free probability naturally leading to a measure $\pi_t$ having these moments, and the density in the statement.

\medskip

(9) Next, since $P(k)=NC(k)$ at $k\leq 3$, we have $E=t$, $V=t$, $\gamma=1/\sqrt{t}$, exactly as for $p_t$. As for the kurtosis, here we must substract a $1$ factor from the kurtosis $3+1/t$ computed for $p_t$, coming from $P(4)=NC(4)\cup\{\cap\!\!\cap\}$, and we obtain $\kappa=2+1/t$.

\medskip

(10) As before with the semicircle and arcsine, we will be not interested in Fourier. Instead, we can compute the moment generating function. And here, at $t=1$:
$$M(z)=\sum_{k=0}^\infty\frac{1}{k+1}\binom{2k}{k}z^k
=\frac{1-\sqrt{1-4z}}{2z}$$

As for the case $t\neq1$, this is something more complicated. We will be back to this.

\medskip

(11) So, this was for the story with the Marchenko-Pastur laws, the idea being that we reasonably understood $\pi_1$, save for some standard computations, that you can find for instance in my book \cite{ba2}, and that $\pi_t$ with $t\neq1$ remains on our to-do list, for later.
\end{proof}

And with this, good news, we have our 7 laws, their summary being as follows:

\begin{theorem}
The main 7 laws in probability $p_1,g_1,e_1,u_1,\gamma_1,\alpha_1,\pi_1$ are
\begin{enumerate}
\item Poisson: $p_1=\frac{1}{e}\sum_{k\geq0}\frac{\delta_k}{k!}$, support $\mathbb N$, moments $M_k=|P(k)|$.

\item Normal: $g_1=\frac{1}{\sqrt{2\pi}}\,e^{-x^2/2}dx$, support $\mathbb R$, moments $M_{2k}=(2k)!!$.

\item Exponential: $e_1=e^{-x}dx$, support $[0,\infty)$, moments $M_k=k!$.

\item Uniform: $u_1=1$, support $[0,1]$, moments $M_k=\frac{1}{k+1}$.

\item Semicircle: $\gamma_1=\frac{1}{2\pi}\sqrt{4-x^2}\,dx$, support $[-2,2]$, moments $M_{2k}=\frac{1}{k+1}\binom{2k}{k}$.

\item Arcsine: $\alpha_1=\frac{1}{\pi\sqrt{4-x^2}}\,dx$, support $[-2,2]$, moments $M_{2k}=\binom{2k}{k}$.

\item Marchenko-Pastur: $\pi_1=\frac{1}{2\pi}\sqrt{4x^{-1}-1}\,dx$, support $[0,4]$, moments $M_k=\frac{1}{k+1}\binom{2k}{k}$.
\end{enumerate}
and these laws have parametric versions $p_t,g_t,e_t,u_t,\gamma_t,\alpha_t,\pi_t$, with a similar theory.
\end{theorem}

\begin{proof}
This is indeed a summary of what we have, with the $(E,V,\gamma,\kappa)$ problematics being left aside, and with the functionals transforms $(F,M)$ being left aside too. And with all this coming with two important comments, as follows:

\medskip

(1) As a rival to our arcsine law $\alpha_1$ on $[-2,2]$ we have the arcsine law $\alpha_1'=\frac{1}{\pi\sqrt{x(4-x)}}\,dx$ on $[0,4]$, with these being related by $f\sim\alpha_1\implies f^2\sim\alpha_1'$, a bit like the semicircle and Marchenko-Pastur laws. Also, the scaling of these arcsine laws $\alpha_1,\alpha_1'$ is a non-trivial business too, because there are several natural ways of modifying the support. Thus, what we called above $\alpha_1$ and $\alpha_t$ remains subject to discussion, and modifications.

\medskip

(2) As a natural question, don't we have also a nice law supported by $\mathbb Z$? In answer, yes, we have the Bessel law $p^2_t$, but this is a bit difficult to introduce and justify with bare hands, and we will talk about it in chapter 4, when doing group theory. Moreover, we will see later in this book that this law has a free analogue, the free Bessel law $\pi^2_t$. Thus, our family of 7 laws will soon increase to 8 laws, and then later, to 9 laws.
\end{proof}

And with this, end of our discussion regarding the 7 main laws. We are now ready for lots of computations, for these laws, and hang on, tough material to come.

\section*{3b. Cauchy, Stieltjes}

Time for some theory. As a starting point for the considerations in this section, and in the remainder of this chapter, we have the following fundamental question:

\begin{question}
How to recover a real probability measure $\mu$ out of its sequence of moments $M_0,M_1,M_2,M_3,\ldots$?
\end{question}

To be more precise, we first met this question in chapter 2, when talking CLT. At that time, I told you that the limiting measure can be recaptured via Fourier and moments, the point being that the Fourier transform of a central limit must be $F(x)=e^{-tx^2/2}$, which in turn shows that the moments must be $M_{2k}=t^k(2k)!!$, and with this in hand, we are led to the question of finding the law $g_t$ having these moments.

\bigskip

We also met Question 3.9 in the previous section, in relation with the laws $\gamma_t,\alpha_t,\pi_t$ investigated there, the point being that what comes out of concrete probability, meaning random walks, or more complicated things like free probability, random matrices or quantum groups, is not exactly the density, but rather the sequence of moments.

\bigskip

In answer now, we certainly know how to solve Question 3.9 by cheating. I mean if you have a good candidate $\mu$, found by whatever means, and no holds barred here, just don't kill someone, compute its moments, and if these are $M_0,M_1,M_2,M_3,\ldots\,$, job done. Simple and bright, and this is indeed what we did before, for $g_t,\gamma_t,\alpha_t,\pi_t$.

\bigskip

So, can we solve Question 3.9 without cheating? Normally, we have here:

\begin{answer}
We must convert the sequence of moments $M_0,M_1,M_2,M_3,\ldots$ into a functional transform $E_\mu(z)$, and then recover $\mu$ from $E_\mu(z)$ via analysis.
\end{answer}

Which brings us into functional transforms. We already know about the Fourier transform $F(x)$ and the moment generating function $M(z)$, and our collection of tools, 2 items so far, can be enlarged with two related transforms, $N(y)$ and $G(\xi)$, as follows:

\index{moment generating function}
\index{Cauchy transform}

\begin{definition}
Associated to a random variable $f:X\to\mathbb R$ are the transforms
$$F(x)=E(e^{ixf})=\sum_{k=0}^\infty\frac{(ix)^k}{k!}\,M_k$$
$$N(y)=E(e^{yf})=\sum_{k=0}^\infty\frac{y^k}{k!}\,M_k$$
called Fourier transform and normalized moment generating function, as well as
$$M(z)=E\left(\frac{1}{1-zf}\right)=\sum_{k=0}^\infty M_kz^k$$
$$G(\xi)=E\left(\frac{1}{\xi-f}\right)=\sum_{k=0}^\infty M_k\xi^{-k-1}$$
called moment generating function, and Cauchy transform.
\end{definition}

Here we assume of course that the expectations converge. Observe that one can easily pass $F\leftrightarrow N$ and $M\leftrightarrow G$. However, in what regards the passage $F/N\leftrightarrow M/G$, this can be something quite complicated, and even undoable, as we will soon discover. 

\bigskip

In practice now, in order to understand how these transforms work, nothing better than working out some examples. For our 7 main laws, the formulae are as follows:

\begin{theorem}
The basic transforms of the main 7 laws are as follows,
$$F_{p_t}(x)=\exp((e^{ix}-1)t)\ ,\ N_{p_t}(y)=\exp((e^y-1)t)\ ,\ M_{p_t}(z)=?\ ,\ G_{p_t}(\xi)=?$$
$$F_{g_t}(x)=e^{-tx^2/2}\ ,\ N_{g_t}(y)=e^{ty^2/2}\ ,\ M_{g_t}(z)=?\ ,\ G_{g_t}(\xi)=?$$
$$F_{e_t}(x)=\frac{t}{t-ix}\ ,\ N_{e_t}(y)=\frac{t}{t-y}\ ,\ M_{e_t}(z)=?\ ,\ G_{e_t}(\xi)=?$$
$$F_{u_t}(x)=\frac{e^{itx}-1}{itx},\,N_{u_t}(y)=\frac{e^{ty}-1}{ty},\,M_{u_t}(z)=-\frac{\log(1-tz)}{tz},\,G_{u_t}(\xi)=-\frac{\log(1-t\xi^{-1})}{t}$$
$$F_{\gamma_t}(x)=?\ ,\ N_{\gamma_t}(y)=?\ ,\ M_{\gamma_t}(z)=\frac{1-\sqrt{1-4tz^2}}{2tz^2}\ ,\ G_{\gamma_t}(\xi)=\frac{\xi-\sqrt{\xi^2-4t}}{2t}$$
$$F_{\alpha_t}(x)=?\ ,\ N_{\alpha_t}(y)=?\ ,\ M_{\alpha_t}(z)=\frac{1}{\sqrt{1-4tz^2}}\ ,\ G_{\alpha_t}(\xi)=\frac{1}{\sqrt{\xi^2-4t}}$$
$$F_{\pi_t}(x)=?\ ,\ N_{\pi_t}(y)=?\ ,\ M_{\pi_1}(z)=\frac{1-\sqrt{1-4z}}{2z}\ ,\ G_{\pi_1}(\xi)=\frac{\xi-\sqrt{\xi^2-4}}{2\xi}$$
with ? standing for ``stay away from that'', and with $M_{\pi_t},G_{\pi_t}$ being more complicated.
\end{theorem}

\begin{proof}
We know all $F,M$ formulae, and these produce the $N,G$ formulae too, via $N(y)=F(-iy)$ and $G(\xi)=\xi^{-1}M(\xi^{-1})$, except for $M/G$ for $u_t$. But here, we have:
$$M_{u_t}(\xi)=\frac{1}{t}\int_0^t\frac{1}{1-zx}\,dx=\frac{1}{t}\Big[-\frac{\log(1-zx)}{z}\Big]_0^t=-\frac{\log(1-tz)}{tz}$$

Alternatively, we can compute the Cauchy transform via moments, as follows:
$$G_{u_t}(\xi)=\sum_{k=0}^\infty\frac{t^k}{k+1}\,\xi^{-k-1}=\frac{1}{t}\sum_{k=0}^\infty\frac{(t\xi^{-1})^{k+1}}{k+1}=-\frac{\log(1-t\xi^{-1})}{t}$$

Thus, one way or another, we are led to all formulae in the statement.
\end{proof}

With this discussed, time now to answer Question 3.9? The result here, called Stieltjes inversion, is something quite tricky, using complex analysis, as follows:

\index{Cauchy transform}
\index{Stieltjes inversion}
\index{density of measure}
\index{sequence of moments}
\index{moment problem}

\begin{theorem}
The density of a real probability measure $\mu$ can be recaptured from the sequence of moments $(M_k)_{k\geq0}$ via the Stieltjes inversion formula
$$d\mu (x)=\lim_{s\searrow 0}-\frac{1}{\pi}\,Im\left(G(x+is)\right)\cdot dx$$
with the function on the right, given in terms of moments by
$$G(\xi)=\xi^{-1}+M_1\xi^{-2}+M_2\xi^{-3}+\ldots$$
being the Cauchy transform of the measure $\mu$.
\end{theorem}

\begin{proof}
The Cauchy transform of our measure $\mu$ is given by:
\begin{eqnarray*}
G(\xi)
&=&\xi^{-1}\sum_{k=0}^\infty M_k\xi^{-k}\\\
&=&\int_\mathbb R\frac{\xi^{-1}}{1-\xi^{-1}y}\,d\mu(y)\\
&=&\int_\mathbb R\frac{1}{\xi-y}\,d\mu(y)
\end{eqnarray*}

Now with $\xi=x+is$, we obtain the following formula:
\begin{eqnarray*}
Im(G(x+is))
&=&\int_\mathbb RIm\left(\frac{1}{x-y+is}\right)d\mu(y)\\
&=&\int_\mathbb R\frac{1}{2i}\left(\frac{1}{x-y+is}-\frac{1}{x-y-is}\right)d\mu(y)\\
&=&-\int_\mathbb R\frac{s}{(x-y)^2+s^2}\,d\mu(y)
\end{eqnarray*}

By integrating over $[a,b]$ we obtain, with the change of variables $x=y+sz$:
\begin{eqnarray*}
\int_a^bIm(G(x+is))dx
&=&-\int_\mathbb R\int_a^b\frac{s}{(x-y)^2+s^2}\,dx\,d\mu(y)\\
&=&-\int_\mathbb R\int_{(a-y)/s}^{(b-y)/s}\frac{s}{(sz)^2+s^2}\,s\,dz\,d\mu(y)\\
&=&-\int_\mathbb R\int_{(a-y)/s}^{(b-y)/s}\frac{1}{1+z^2}\,dz\,d\mu(y)\\
&=&-\int_\mathbb R\left(\arctan\frac{b-y}{s}-\arctan\frac{a-y}{s}\right)d\mu(y)
\end{eqnarray*}

Now observe that with $s\searrow0$ we have:
$$\lim_{s\searrow0}\left(\arctan\frac{b-y}{s}-\arctan\frac{a-y}{s}\right)
=\begin{cases}
\frac{\pi}{2}-\frac{\pi}{2}=0& (y<a)\\
\frac{\pi}{2}-0=\frac{\pi}{2}& (y=a)\\
\frac{\pi}{2}-(-\frac{\pi}{2})=\pi& (a<y<b)\\
0-(-\frac{\pi}{2})=\frac{\pi}{2}& (y=b)\\
-\frac{\pi}{2}-(-\frac{\pi}{2})=0& (y>b)
\end{cases}$$

We therefore obtain the following formula:
$$\lim_{s\searrow0}\int_a^bIm(G(x+is))dx=-\pi\left(\mu(a,b)+\frac{\mu(a)+\mu(b)}{2}\right)$$

Now when $\mu$ is continuous, this leads to the formula in the statement.
\end{proof}

As a comment on this, Theorem 3.13 assumes of course that the measure $\mu$ to be found has indeed a density, the point being that, in the opposite case, the density found there will not integrate up to 1. However, we can deal with atoms too, as shown by:

\index{atoms}

\begin{remark}
Stieltjes inversion can deal with atoms too, according to:
$$\mu=\nu+\sum_ic_i\delta_{a_i}\implies G_\mu(\xi)=G_{\nu}(\xi)+\sum_i\frac{c_i}{\xi-a_i}$$
That is, the poles of the Cauchy transform $G_\mu$ correspond to the atoms of $\mu$.
\end{remark}

With this discussed, time now for some examples? For the main 7 probability laws that we have, the applications of Stieltjes inversion range from trivial to undoable, and with the results that we find ranging from useless to very interesting, as follows:

\begin{theorem}
In the main 7 situations, Stieljes inversion applies as follows:
\begin{enumerate}
\item Poisson: $M_k=\sum_{\pi\in P(k)}t^{|\pi|}$ has troubles producing $p_t=e^{-t}\sum_{k\geq0}\frac{t^k\delta_k}{k!}$.

\item Normal: $M_{2k}=t^k(2k)!!$ has troubles producing $g_t=\frac{1}{\sqrt{2\pi t}}\,e^{-x^2/2t}dx$.

\item Exponential: $M_k=k!/t^k$ has troubles producing $e_t=te^{-tx}\,dx$. 

\item Uniform: $M_k=t^k/(k+1)$ produces indeed $u_t=\frac{1}{t}\,dx$.

\item Semicircle: $M_{2k}=\frac{t^k}{k+1}\binom{2k}{k}$ produces indeed $\gamma_t=\frac{1}{2\pi t}\sqrt{4t-x^2}\,dx$.

\item Arcsine: $M_{2k}=t^k\binom{2k}{k}$ produces indeed $\alpha_t=\frac{1}{\pi\sqrt{4t-x^2}}\,dx$.

\item Marchenko-Pastur: we get indeed $\pi_t=\max(1-t,0)\delta_0+\frac{\sqrt{4t-(x-1-t)^2}}{2\pi x}\,dx$.
\end{enumerate}
\end{theorem}

\begin{proof}
As mentioned, all this is quite tricky, the idea being as follows:

\medskip

(1) Assume indeed that, in your theoretical physics computations, you came upon the numbers $M_k=\sum_{\pi\in P(k)}t^{|\pi|}$, and are looking for the measure $\mu$ having these as moments. In this case Stieltjes inversion cannot really help you, because $M(z)=\sum_{k\geq 0}M_kz^k$ is not computable. However, as a fix, you can try to find equations for $M(z)$ and then $G(\xi)$, and use these in the context of Stieltjes inversion. And, good luck here.

\medskip

(2) For $M_{2k}=t^k(2k)!!$ same situation as for Poisson, $M(z)$ being not computable, and in fact even worse, because, while some work on Poisson could eventually lead you to the conclusion that $G(\xi)$ consists of poles only, here there is no such trick.

\medskip

(3) Same situation for $M_k=k!/t^k$, with the function $M(z)$ being not computable.

\medskip

(4) This works. So, assume that in your physics computations, you came upon the numbers $M_k=t^k/(k+1)$, and are looking for the measure $\mu$ having these as moments. Well, you ask your math buddy, who will first compute the Cauchy transform:
$$G(\xi)=\sum_{k\geq0}\frac{t^k\xi^{-k-1}}{k+1}=-\frac{\log(1-t\xi^{-1})}{t}$$

Next, the guy will set $\xi=x+is$, and look at the following quantity:
$$e^{-tG(x+is)}=1-\frac{t}{x+is}$$

And then, guy will argue that with $s\searrow 0$, the right term goes to $1-t/x\in\mathbb R$, so we must have $-tIm(G(x+is))\to k\pi$ with $k\in\mathbb Z$, and more specifically, by looking at what happens to the real part too, with $k=1$. Thus, by Stieltjes inversion:
$$d\mu (x)=\lim_{s\searrow 0}-\frac{1}{\pi}\,Im\left(G(x+is)\right)\cdot dx=\frac{1}{t}$$

Which sounds very good, you can now write a joint paper with your math friend.

\medskip

(5) This works too. Indeed, with $M_{2k}=\frac{t^k}{k+1}\binom{2k}{k}$, the Cauchy transform is: 
$$G(\xi)
=\xi^{-1}\sum_{k=0}^\infty\frac{1}{k+1}\binom{2k}{k}(t\xi^{-2})^k
=\frac{\xi-\sqrt{\xi^2-4t}}{2t}$$

Now let us apply Theorem 3.13. The study here goes as follows:

\medskip

-- According to the general philosophy of the Stieltjes formula, the first term, namely $\xi/2t$, which is ``trivial'', will not contribute to the density. 

\medskip

-- As for the second term, which is something non-trivial, this will contribute to the density, the rule here being that the square root $\sqrt{\xi^2-4t}$ will be replaced by the ``dual'' square root $\sqrt{4t-x^2}\,dx$, and that we have to multiply everything by $-1/\pi$. 

\medskip

-- As a conclusion, by Stieltjes inversion we obtain the following density:
$$d\mu(x)
=-\frac{1}{\pi}\cdot-\frac{\sqrt{4t-x^2}}{2t}\,dx
=\frac{1}{2\pi t}\sqrt{4t-x^2}dx$$

(6) This works too. Indeed, with $M_{2k}=t^k\binom{2k}{k}$, the Cauchy transform is:
$$G(\xi)
=\xi^{-1}\sum_{k=0}^\infty\binom{2k}{k}(t\xi^{-2})^k
=\frac{1}{\sqrt{\xi^2-4t}}$$

Thus, by Stieltjes inversion we obtain the following density, as stated:
$$d\mu(x)
=-\frac{1}{\pi}\cdot-\frac{1}{\sqrt{4t-x^2}}\,dx
=\frac{1}{\pi\sqrt{4t-x^2}}\,dx$$

(6') Still talking arcsine, let us do as well the computation for the ``rival'' situation from the proof of Theorem 3.8, namely $M_k=t^k\binom{2k}{k}$. The Cauchy transform is:
$$G(\xi)
=\xi^{-1}\sum_{k=0}^\infty\binom{2k}{k}(t\xi^{-1})^k
=\frac{1}{\sqrt{\xi^2-4t\xi}}$$

Thus, by Stieltjes inversion we obtain the following density, which at $t=1$ is the ``rival'' main arcsine law $\alpha_1'=\frac{1}{\pi\sqrt{x(4-x)}}\,dx$ on $[0,4]$, from the proof of Theorem 3.8:
$$d\mu(x)
=-\frac{1}{\pi}\cdot-\frac{1}{\sqrt{4tx-x^2}}\,dx
=\frac{1}{\pi\sqrt{4tx-x^2}}\,dx$$

(7) We have kept the best for the end. The Marchenko-Pastur situation, coming from free probability, random matrices and quantum groups, is as follows:
$$M_k=\sum_{\pi\in NC(k)}t^{|\pi|}$$

In the case $t=1$ we have $M_k=|NC(k)|=\frac{1}{k+1}\binom{2k}{k}$, leading right away to:
$$M(z)=\frac{1-\sqrt{1-4z}}{2z}\quad ,\quad G(\xi)=\frac{1-\sqrt{1-4\xi^{-1}}}{2}$$

Thus, by Stieltjes inversion we obtain the following density, as stated:
$$d\mu(x)
=-\frac{1}{\pi}\cdot-\frac{\sqrt{4x^{-1}-1}}{2}\,dx
=\frac{1}{2\pi}\sqrt{4x^{-1}-1}\,dx$$

In the case $t\neq1$ things are more complicated, and the moment formula above leads, say via some free cumulant know-how, to the following formula for the $M$ series:
$$M(z)=\frac{1-t}{2}+\frac{1-\sqrt{1-2(1+t)z+(1-t)^2z^2}}{2z}$$

Equivalently, the Cauchy transform is given by the following formula:
$$G(\xi)=\frac{1-t+\xi-\sqrt{(\xi-1-t)^2-4t}}{2\xi}$$

Thus, by Stieltjes inversion we obtain the following density, as stated:
$$d\mu(x)
=\max(1-t,0)\delta_0+\frac{\sqrt{4t-(x-1-t)^2}}{2\pi x}\,dx$$

And more on this, Marchenko-Pastur at $t\neq1$, later in this book.
\end{proof}

\section*{3c. Hankel determinants}

In view of what we have, namely bad results for the laws $p_t,g_t$ that we are officially interested in, then bad results for $e_t$ too, and then good results for $u_t,\gamma_t,\alpha_t,\pi_t$, I am pretty much sure that you currently suspect me to have introduced our 7 main probability laws, just as a matter of having a story to tell with Stieltjes inversion, and freeness.

\bigskip

Error. We will discover in the remainder of this chapter that the continuation of the story with Stieltjes inversion, featuring Hankel determinants and orthogonal polynomials, will put our family $p_t,g_t,e_t,u_t,\gamma_t,\alpha_t,\pi_t$ back together, and in a spectacular way.

\bigskip

As a starting point, Stieltjes inversion does not fully solve the moment problem, because we still have the question of understanding when a sequence of numbers $(M_k)_{k\geq0}$ can be the moments of a measure $\mu$. For instance we certainly must have $M_0=1$, but we must have as well the following inequality, corresponding to $V\geq0$:
$$M_2\geq M_1^2$$

In answer now, we have the following result, complementing Theorem 3.13:

\index{Hankel determinant}
\index{moment problem}

\begin{theorem}
A sequence of numbers $M_0,M_1,M_2,M_3,\ldots\in\mathbb R$, with $M_0=1$, is the series of moments of a real probability measure $\mu$ precisely when:
$$\begin{vmatrix}M_0\end{vmatrix}\geq0\quad,\quad 
\begin{vmatrix}
M_0&M_1\\
M_1&M_2
\end{vmatrix}\geq0\quad,\quad 
\begin{vmatrix}
M_0&M_1&M_2\\
M_1&M_2&M_3\\
M_2&M_3&M_4
\end{vmatrix}\geq0\quad,\quad 
\ldots$$
That is, the associated Hankel determinants must be all positive.
\end{theorem}

\begin{proof}
This is something a bit more advanced, the idea being as follows:

\medskip

(1) Some numerics first. We know that the first two Hankel determinants are $H_1=1$ and $H_2=V$, both positive. Regarding the third Hankel determinant, this is given by:
$$H_3=\begin{vmatrix}
1&M_1&M_2\\
M_1&M_2&M_3\\
M_2&M_3&M_4
\end{vmatrix}=(M_4-M_2^2)V-(M_3-M_1M_2)^2$$

Which does not look obviously positive, so let us do a check for the binomial law $b_{np}$. And here, after 1 hour of computations, we obtain, as desired:
$$H_3(b_{np})=2n^2(n-1)p^3(1-p)^3\geq0$$

Getting now to $H_4$, things get much worse, and as a wise decision, let us assume that our measure is centered, $M_{2k+1}=0$. In this case $H_4$ is easily computable, given by:
$$H_4=\begin{vmatrix}
1&0&M_2&0\\
0&M_2&0&M_4\\
M_2&0&M_4&0\\
0&M_4&0&M_6
\end{vmatrix}=(M_2M_6-M_4^2)(M_4-M_2^2)$$

Now with input from $g_t$, namely $M_2=t,M_4=3t^2,M_6=15t^3$, we obtain:
$$H_1(g_t)=1\quad,\quad H_2(g_t)=t\quad,\quad H_3(g_t)=2t^3\quad,\quad H_4(g_t)=12t^6$$

Which looks good, all positive numbers, but as an overall conclusion to this, the verification of $H_k\geq0$ by brute-force computation of $H_k$ is not a good method.

\medskip

(2) Fortunately, linear algebra comes to the rescue. Recall indeed from there the positivity theorem of Sylvester, stating that a matrix is positive precisely when its principal minors are all positive. In view of this, the conditions $H_k\geq0$, taken altogether, precisely mean that the corresponding matrices, say $A_k$, are all positive, $A_k\geq0$. 

\medskip

(3) Now observe that for each of these matrices $A_k$, given a vector $c\in\mathbb C^N$, we have:
$$<Ac,c>=\sum_{i,j=1}^kM_{i+j}c_i\bar{c_j}$$

Thus, as an overall conclusion, our theorem states that $(M_k)_{k\geq0}$ are the moments of a measure $\mu$ precisely when we have, for any $k\in\mathbb N$, and any $c\in\mathbb C^N$:
$$\sum_{i,j=1}^nM_{i+j}c_i\bar{c}_j\geq0$$

(4) But in one sense this is elementary, coming from the following computation:
$$\sum_{i,j=1}^nM_{i+j}c_i\bar{c}_j
=\int_\mathbb R\sum_{i,j=1}^nc_i\bar{c}_jx^{i+j}d\mu(x)
=\int_\mathbb R\left|\sum_{i=1}^nc_ix^i\right|^2d\mu(x)$$

(5) As for the other sense, this is something more delicate, requiring some functional analysis study. So, exercise for you, to learn more about all this.
\end{proof}

Getting now to our main 7 probability measures, we have, quite remarkably, the following beautiful and mysterious result regarding them:

\begin{theorem}
The Hankel determinants for the main $7$ laws are as follows:
\begin{enumerate}
\item Poisson, $M_k=|P(k)|$, here $H_k=\prod_{i=1}^{k-1}i!$.

\item Normal, $M_{2k}=(2k)!!$, here $H_k=\prod_{i=1}^{k-1}i!$.

\item Exponential, $M_k=k!$, here $H_k=\prod_{i=1}^{k-1}i!^2$. 

\item Uniform, $M_k=1/(k+1)$, here $H_k^{-1}=k!\prod_{i=1}^{2k-1}\binom{i}{[i/2]}$.

\item Semicircle, $M_{2k}=\frac{1}{k+1}\binom{2k}{k}$, here $H_k=1$.

\item Arcsine, $M_{2k}=\binom{2k}{k}$, here $H_k=2^{k-1}$.

\item Marchenko-Pastur, $M_k=\frac{1}{k+1}\binom{2k}{k}$, here $H_k=1$.
\end{enumerate}
\end{theorem}

\begin{proof}
This is something tough, that we will not attempt to prove here. This being said, we can do some verifications, and formulate some comments, as follows:

\medskip

(1) Some numerics first. At $k=2$ we obtain indeed the correct variances:
$$H_2\ =\ 1\ ,\ 1\ ,\ 1\ ,\ \frac{1}{12}\ ,\ 1\ ,\ 2\ ,\ 1$$

At $k=3$ we obtain the following numbers, which are the correct ones too, as you can check by using the formula of $H_3$ from the proof of Theorem 3.16:
$$H_3\ =\ 2\ ,\ 2\ ,\ 4\ ,\ \frac{1}{2160}\ ,\ 1\ ,\ 4\ ,\ 1$$

At $k=4$ we obtain the following numbers, correct too, good exercise for you:
$$H_4\ =\ 12\ ,\ 12\ ,\ 144\ ,\ \frac{1}{24192000}\ ,\ 1\ ,\ 8\ ,\ 1$$

(2) And so on. In general, however, there are several proofs, all complicated. So, come back here after reading the next section, on orthogonal polynomials, which can help.

\medskip

(3) As a comment now, for the rival arcsine law $\alpha_1'$, supported by $[0,4]$, and whose moments are $M_k=\binom{2k}{k}$, we obtain $H_k=2^{k-1}$, exactly as for the usual arcsine law $\alpha_1$. Which is not suprising, the passage $\alpha_1\to\alpha_1'$ being similar to $\gamma_1\to\pi_1$.

\medskip

(4) Finally, in what regards the parametric versions $p_t,g_t,e_t,u_t,\gamma_t,\alpha_t,\pi_t$, in some cases we can get away with homogeneity, while in other cases the problem becomes substantially more complicated. And of course, exercise for you, to learn more about all this.
\end{proof}

\section*{3d. Orthogonal polynomials} 

We would like to end this chapter with a discussion on orthogonal polynomials, which are advanced objects, related to the above, and whose values for the main 7 laws will be something spectacular again, in the spirit of Theorem 3.17. Let us start with:

\index{orthonormal basis}
\index{Gram-Schmidt}
\index{separable space}
\index{Hilbert space}
\index{orthogonal basis}

\begin{theorem}
Any Hilbert space $H$ has an orthonormal basis $\{e_i\}_{i\in I}$, which is by definition a set of vectors whose span is dense in $H$, and which satisfy
$$<e_i,e_j>=\delta_{ij}$$
with $\delta$ being a Kronecker symbol. The cardinality $|I|$ of the index set, which can be finite, countable, or uncountable, depends only on $H$, and is called dimension of $H$. We have
$$H\simeq l^2(I)$$
in the obvious way, mapping $\sum\lambda_ie_i\to(\lambda_i)$. The Hilbert spaces with $\dim H=|I|$ being countable, such as $l^2(\mathbb N)$, are all isomorphic, and are called separable.
\end{theorem}

\begin{proof}
This is certainly something that you know, coming from Gram-Schmidt, with of course due attention to infinite dimensionality, needing a few minor fixes.
\end{proof}

According to Theorem 3.18, there is only one separable Hilbert space, up to isomorphism. There are many interesting things that can be said, about this magic and unique Hilbert space. As a first such result, which is quite theoretical, we have:

\index{orthogonal polynomials}
\index{Weierstrass basis}

\begin{theorem}
The following happen, in relation with separability:
\begin{enumerate}
\item The Hilbert space $H=L^2[-1,1]$ is separable, with orthonormal basis coming by applying Gram-Schmidt to the basis $\{x^k\}_{k\in\mathbb N}$, coming from Weierstrass.

\item In fact, any $H=L^2(\mathbb R,\mu)$, with $d\mu(x)=f(x)dx$, is separable, and the same happens in higher dimensions, for $H=L^2(\mathbb R^N,\mu)$, with $d\mu(x)=f(x)dx$.

\item More generally, given a separable abstract measured space $X$, the associated Hilbert space of square-summable functions $H=L^2(X)$ is separable.
\end{enumerate}
\end{theorem}
 
\begin{proof}
Many things can be said here, the idea being as follows:

\medskip

(1) The fact that $H=L^2[-1,1]$ is separable is clear indeed from the Weierstrass approximation theorem, which provides us with the algebraic basis $g_k=x^k$, which can be orthogonalized via the Gram-Schmidt procedure, as explained in Theorem 3.18.

\medskip

(2) Regarding now more general spaces, of type $H=L^2(\mathbb R,\mu)$, we can use here the same argument, after modifying if needed our measure $\mu$, in order for the functions $g_k=x^k$ to be indeed square-summable. As for higher dimensions, the situation here is similar, because we can use here the multivariable polynomials $g_k(x)=x_1^{k_1}\ldots x_N^{k_N}$.

\medskip

(3) Concerning the last assertion, regarding the general spaces of type $H=L^2(X)$, which generalizes all this, this comes as a consequence of general measure theory, and we will leave some learning, and working out the details here, as an instructive exercise.

\medskip

(4) Finally, let us mention that all this is just the tip of the iceberg. For instance for the unit circle $X=\mathbb T$ the good orthonormal basis is the Fourier one, $\{z^n\}_{n\in\mathbb Z}$, and with this being unrelated to the constructions in (1,2,3). Many things to be learned here.
\end{proof}

At a more concrete level now, Theorem 3.19 suggests formulating:

\index{orthogonal polynomials}
\index{Gram-Schmidt}

\begin{definition}
The orthogonal polynomials with respect to $d\mu(x)=f(x)dx$ are polynomials $P_k\in\mathbb R[x]$ of degree $k\in\mathbb N$, which are orthogonal inside $H=L^2(\mathbb R,\mu)$:
$$\int_\mathbb RP_k(x)P_l(x)f(x)dx=0\quad,\quad\forall k\neq l$$
Equivalently, these orthogonal polynomials $\{P_k\}_{k\in\mathbb N}$, which are each unique modulo scalars, appear from the Weierstrass basis $\{x^k\}_{k\in\mathbb N}$, by doing Gram-Schmidt.
\end{definition}

As a first observation, the orthogonal polynomials exist indeed for any real measure $d\mu(x)=f(x)dx$, because we can obtain them from the monomials $x^k$ via Gram-Schmidt, as indicated above. It is possible to be a bit more explicit here, as follows:

\begin{theorem}
The orthogonal polynomials with respect to $\mu$ are given by
$$P_k=c_k\begin{vmatrix}
M_0&M_1&\ldots&M_k\\
M_1&M_2&\ldots&M_{k+1}\\
\vdots&\vdots&&\vdots\\
M_{k-1}&M_k&\ldots&M_{2k-1}\\
1&x&\ldots&x^k
\end{vmatrix}$$
where $M_k=\int_\mathbb Rx^kd\mu(x)$ are the moments of $\mu$, and $c_k\in\mathbb R^*$ can be any numbers.
\end{theorem}

\begin{proof}
Let us first see what happens at small values of $k\in\mathbb N$. At $k=0$ our formula holds indeed, due to $M_0=1$. At $k=1$, using again $M_0=1$, the formula is as follows:
$$P_1=c_1\begin{vmatrix}M_0&M_1\\ 1&x\end{vmatrix}=c_1(x-M_1)$$

But this is again the good formula, because the degree is 1, and we have:
\begin{eqnarray*}
<1,P_1>
&=&c_1<1,x-M_1>\\
&=&c_1(<1,x>-<1,M_1>)\\
&=&c_1(M_1-M_1)\\
&=&0
\end{eqnarray*}

At $k=2$ now, things get more complicated, with the formula being as follows:
$$P_2=c_2\begin{vmatrix}
M_0&M_1&M_2\\
M_1&M_2&M_3\\
1&x&x^2
\end{vmatrix}$$

However, no need for big computations here, in order to check the orthogonality, because by using the fact that $x^k$ integrates up to $M_k$, we obtain:
$$<1,P_2>=\int_\mathbb RP_2(x)d\mu(x)=c_2\begin{vmatrix}
M_0&M_1&M_2\\
M_1&M_2&M_3\\
M_0&M_1&M_2
\end{vmatrix}=0$$

Similarly, again by using the fact that $x^k$ integrates up to $M_k$, we have as well:
$$<x,P_2>=\int_\mathbb RxP_2(x)d\mu(x)=c_2\begin{vmatrix}
M_0&M_1&M_2\\
M_1&M_2&M_3\\
M_1&M_2&M_3
\end{vmatrix}=0$$

Thus, result proved at $k=0,1,2$, and the proof in general is similar.
\end{proof}

In practice now, all this leads us to a lot of interesting combinatorics, and countless things can be said. For the simplest measured space $X\subset\mathbb R$, which is the interval $[-1,1]$, with its uniform measure, the orthogonal basis problem can be solved as follows:

\index{Legendre polynomials}
\index{Rodrigues formula}
\index{Legendre equation}

\begin{theorem}
The orthogonal polynomials for $L^2[-1,1]$, subject to
$$\int_{-1}^1P_k(x)P_l(x)\,dx=\delta_{kl}$$
and called Legendre polynomials, satisfy the following differential equation,
$$(1-x^2)P_k''(x)-2xP_k'(x)+k(k+1)P_k(x)=0$$
which is the Legendre equation from physics. Moreover, we have the formula
$$(k+1)P_{k+1}(x)=(2k+1)xP_k(x)-kP_{k-1}(x)$$
called Bonnet recurrence formula, as well as the formula
$$P_k(x)=\frac{1}{2^kk!}\cdot\frac{d^k}{dx^k}\left(1-x^2\right)^k$$
called Rodrigues formula for the Legendre polynomials.
\end{theorem}

\begin{proof}
Many things going on here, the idea being as follows:

\medskip

(1) The first assertion is clear, because the Gram-Schmidt procedure applied to the Weierstrass basis $\{x^k\}$ can only lead to a certain family of polynomials $\{P_k\}$, with each $P_k$ being of degree $k$, and also unique, if we assume that it has positive leading coefficient, with this $\pm$ choice being needed, as usual, at each step of Gram-Schmidt.

\medskip

(2) In order to have now an idea about these beasts, here are the first few of them, which can be obtained say via a straightforward application of Gram-Schmidt:
\begin{eqnarray*}
P_0&=&1\\
P_1&=&x\\
P_2&=&(3x^2-1)/2\\
P_3&=&(5x^3-3x)/2\\
P_4&=&(35x^4-30x^2+3)/8
\end{eqnarray*}

(3) Now thinking about what Gram-Schmidt does, this is certainly something by recurrence. And examining the recurrence leads to the Legendre equation, as stated. As for the Bonnet recurrence formula, the story here is similar.

\medskip

(4) Regarding the Rodrigues formula, by uniqueness no need to try to understand where this formula comes from, and we have two choices here, either by verifying that $\{P_k\}$ is orthonormal, or by verifying the Legendre equation. And both methods work.

\medskip

(5) In relation with the 7 laws, the shifted Legendre polynomials $P_k(2x-1)$ solve the problem for $u_1$, uniform measure on $[0,1]$. So one done, 6 to go. By the way, the choice $[-1,1]$ in the statement is very standard, coming from physics, the Legendre polynomials defined in this way being able to solve the hydrogen atom. See Griffiths \cite{gri}.
\end{proof}

The above result is just the tip of the iceberg, and as a continuation, we have:

\index{Jacobi polynomials}
\index{Chebycheff polynomials}
\index{Gegenbauer polynomials}

\begin{theorem}
The orthogonal polynomials for $L^2[-1,1]$, with measure
$$d\mu(x)=(1-x)^a(1+x)^bdx$$
called Jacobi polynomials, satisfy as well a degree $2$ equation, namely
$$(1-x^2)P_k''(x)+(b-a-(a+b+2)x)P_k'(x)+k(k+a+b+1)P_k(x)=0$$
as well as an order $2$ recurrence relation, and are given by the following formula:
$$P_k(x)=\frac{(-1)^k}{2^kk!}(1-x)^{-a}(1+x)^{-b}\frac{d^k}{dx^k}\left[(1-x)^a(1+x)^b(1-x^2)^k\right]$$
As main examples, at $a=b=0$ we recover the Legendre polynomials, and at $a,b=\pm\frac{1}{2}$ we recover the four types of Chebycheff polynomials, from trigonometry. 
\end{theorem}

\begin{proof}
Again, many things going on here, the idea being as follows:

\medskip

(1) To start with, in what regards the precise statement, the order 2 recurrence relation mentioned there is something quite complicated, as follows:
\begin{eqnarray*}
&&2k(k+a+b)(2k+a+b-2)P_k(x)\\
&=&(2k+a+b-1)\left[(2k+a+b)(2k+a+b-2)x+a^2-b^2\right]P_{k-1}(x)\\
&-&2(k+a-1)(k+b-1)(2k+a+b)P_{k-2}(x)
\end{eqnarray*}

(2) Regarding now the main particular cases of the Jacobi polynomials, these are the Gegenbauer polynomials, appearing at $a=b$, with as particular cases:

\medskip

-- The Legendre polynomials from Theorem 3.22, appearing at $a=b=0$.

\medskip

-- The Chebycheff polynomials of the first kind, $T_k(\cos t)=\cos(kt)$, at $a=b=-\frac{1}{2}$.

\medskip

-- And of the second kind too, $U_k(\cos t)\sin t=\sin((k+1)t)$, appearing at $a=b=\frac{1}{2}$.

\medskip

(3) Passed Gegenbauer, at $a=-\frac{1}{2},b=\frac{1}{2}$ we have the Chebycheff polynomials of the third kind, $V_k(\cos t)\cos(t/2)=\cos((k+1/2)t)$, and at $a=\frac{1}{2},b=-\frac{1}{2}$ we have the Chebycheff polynomials of the fourth kind, $W_k(\cos t)\sin(t/2)=\sin((k+1/2)t)$.

\medskip

(4) In relation with probability, the Jacobi polynomials solve the orthogonalization problem for the beta distributions from Comment 3.5, up to some normalizations. In particular, in relation with the 7 laws, $u_1,\alpha_1,\gamma_1,\pi_1$ done, and 3 more to go.

\medskip

(5) Finally, regarding the proof, the statement itself appears as a generalization of Theorem 3.22, which corresponds to the particular case $a=b=0$, and the proof is quite similar. We will leave learning more about all this as an exercise.
\end{proof}

Getting now to other spaces $X\subset\mathbb R$, we first have the following result, dealing with $e_1$, which complements well Theorem 3.22, for the needs of basic quantum mechanics:

\index{Laguerre polynomials}

\begin{theorem}
The orthogonal polynomials for $L^2[0,\infty)$, with scalar product
$$<f,g>=\int_0^\infty f(x)g(x)e^{-x}\,dx$$
are the Laguerre polynomials $\{P_k\}$, satisfying the following differential equation,
$$xP_k''(x)+(1-x)P_k'(x)+kP_k(x)=0$$
as well as the following order $2$ recurrence relation,
$$(k+1)P_{k+1}(x)=(2k+1-x)P_k(x)-kP_{k-1}(x)$$
and which are given by the following formula,
$$P_k(x)=\frac{e^x}{k!}\cdot\frac{d^k}{dx^k}\left(e^{-x}x^k\right)$$
called Rodrigues formula for the Laguerre polynomials.
\end{theorem}

\begin{proof}
The story here is very similar to that of the Legendre and Jacobi polynomials, and many further things can be said here, with exercise for you to learn a bit about all this. Let us record as well a few numeric values, for the Laguerre polynomials:
\begin{eqnarray*}
P_0&=&1\\
P_1&=&1-x\\
P_2&=&(x^2-4x+2)/2\\
P_3&=&(-x^3+9x^2-18x+6)/6\\
P_4&=&(x^4-16x^3+72x^2-96x+24)/24
\end{eqnarray*}

Finally, for the story to be complete, no discussion about the Laguerre polynomials would be complete without a word about their use, in quantum mechanics. So, have a look at Griffiths \cite{gri}, see how these beasts produce the wave functions of hydrogen.
\end{proof}

Finally, regarding the space $X=\mathbb R$ itself, we have here the following result:

\index{Hermite polynomials}

\begin{theorem}
The orthogonal polynomials for $L^2(\mathbb R)$, with scalar product
$$<f,g>=\int_0^\infty f(x)g(x)e^{-x^2}\,dx$$
are the Hermite polynomials $\{P_k\}$, satisfying the following differential equation,
$$P_k''(x)-2xP_k'(x)+P_k(x)=0$$
as well as the following order $2$ recurrence relation,
$$P_{k+1}(x)=2xP_k(x)-2kP_{k-1}(x)$$
and which are given by the following formula,
$$P_k(x)=(-1)^k e^{x^2}\cdot\frac{d^k}{dx^k}\big(e^{-x^2}\big)$$
called Rodrigues formula for the Hermite polynomials.
\end{theorem}

\begin{proof}
As before, the story here is quite similar to that of the Legendre and other orthogonal polynomials, and exercise for you to learn a bit about all this. Let us record as well a few numeric values, for the Hermite polynomials:
\begin{eqnarray*}
P_0&=&1\\
P_1&=&2x\\
P_2&=&4x^2-2\\
P_3&=&8x^3-12x\\
P_4&=&16x^4-48x^2+12
\end{eqnarray*}

Finally, do not forget to work out some rescalings too, with the density in the statement $e^{-x^2}$ replaced by the more familiar $e^{-x^2/2}$, corresponding to the normal law $g_1$.
\end{proof}

And with this, good news, end of the story with the orthogonal polynomials, at least at the  introductory level, and this due to the following fact, which is something quite technical, and that we will not attempt to prove, or even explain in detail here:

\begin{fact}
From an abstract point of view, coming from degree $2$ equations, and Rodrigues formulae for the solutions, there are only three types of ``classical'' orthogonal polynomials, namely the Jacobi, Laguerre and Hermite ones, discussed above.
\end{fact}

And isn't this amazing, on par with the combinatorial beauties from Theorem 3.17. Simply put, up to some manipulations on our 7 laws, namely merging all the beta ones, and forgetting about Poisson, that is all. And with this, coming by theorem.

\bigskip

The continuation of the story, however, is more complicated, as follows:

\index{Charlier polynomials}
\index{Askey scheme}

\begin{fact}
The orthogonal polynomials for $p_1$ are the Charlier polynomials, which are quite complicated. The unification with Jacobi is done via the Askey polynomials, advanced level. On top of this, we have the Askey-Wilson polynomials, expert level.
\end{fact}

And we will end with this. Take it easy I guess, we learned many things in this chapter, regarding the 7 laws, but this is not the only way of viewing things. We will soon discover that another viewpoint, called easiness, allows us to reconcile $p_t,g_t,\gamma_t,\pi_t$.

\section*{3e. Exercises}

We had a lot of combinatorics in this chapter, and as exercises, we have:

\begin{exercise}
Learn more about $e_t,u_t,\gamma_t$, all basic distributions.
\end{exercise}

\begin{exercise}
Further meditate on $\alpha_t,\pi_t$, at the $t=1$ value, and scaling.
\end{exercise}

\begin{exercise}
Try to compute the missing transforms, for the 7 laws.
\end{exercise}

\begin{exercise}
What is the measure having $E_k=\binom{k}{[k/2]}$ as moments?
\end{exercise}

\begin{exercise}
Clarify the general theory of Hankel determinants.
\end{exercise}

\begin{exercise}
Do some Hankel computations, for the main 7 laws.
\end{exercise}

\begin{exercise}
Learn more about Chebycheff and other orthogonal polynomials.
\end{exercise}

\begin{exercise}
In view of what we saw, what is the simplest measure?
\end{exercise}

As bonus exercise, read more about Hilbert spaces, and about operators too.

\chapter{Groups, easiness}

\section*{4a. Poisson, revised}

We have seen in the previous chapter that, in relation with several advanced aspects, the Poisson and normal laws $p_t,g_t$ are best seen as being part of the family $\{p_t,g_t,e_t,u_t,\gamma_t,\alpha_t,\pi_t\}$, with $e_t,u_t$ being the exponential and uniform laws, and with $\gamma_t,\alpha_t,\pi_t$ being the semicircle, arcsine and Marchenko-Pastur laws. Informally, and up to rescalings, and up to an atom issue for $\pi_t$ at $t\in(0,1)$, our family can be described as ``Poisson, normal, exponential and basic beta''. And with this being certainly great.

\bigskip

As an issue, however, at the truly advanced level, all this leads us into the Askey scheme for orthogonal polynomials \cite{awi}, where the normal law $g_t$, corresponding to the Hermite polynomials, lies at the bottom, which is very nice, while the Poisson law $p_t$, corresponding to the Charlier polynomials, lies higher up, which is not ideal. 

\bigskip

So, what to do? Come up with an alternative approach to all this, according to:

\begin{principle}
There are two ways of conceptually viewing $p_t,g_t$:
\begin{enumerate}
\item As part of the family $\{p_t,g_t,e_t,u_t,\gamma_t,\alpha_t,\pi_t\}$, with $e_t,u_t,\gamma_t,\alpha_t,\pi_t$ being the exponential, uniform, semicircle, arcsine and Marchenko-Pastur laws. And with this being part of the Askey scheme for orthogonal polynomials.

\item As part of the family $\{p_t,g_t,p^2_t,g_t^t,\pi_t,\gamma_t,\pi^2_t,\gamma_t^t\}$, with $p^2_t,g_t^t,\pi_t,\gamma_t,\pi^2_t,\gamma_t^t$ being the Bessel, shifted normal, Marchenko-Pastur, semicircle, free Bessel and shifted semicircle laws. And with this coming from the ``easiness'' philosophy.
\end{enumerate}
\end{principle}

Which sounds quite good, guess I'll just have to explain you what easiness is, and then provide some details on the two newcomers, $p^2_t,\pi^2_t$, and we are ready to go. With this meaning theory and computations for (2), and then comparison with (1).

\bigskip

In answer, patience. Easiness is something that can be formally defined in terms of partitions and moment combinatorics, the key formulae here being as follows:
$$M_k(p_t)=\sum_{\pi\in P(k)}t^{|\pi|}\quad,\quad M_k(g_t)=\sum_{\pi\in P_2(k)}t^{|\pi|}$$ 

Indeed, we know that we have similar formulae for $\pi_t,\gamma_t$, involving $NC(k),NC_2(k)$. Which makes it pretty much clear that we can develop a theory, based on this.

\bigskip

However, and here comes my point, for best results and applications, and further extensions too, and so on, it is better to do things the slow way, according to:

\begin{principle}
Easiness comes from groups, with $p_t,g_t$ coming from $S_N,O_N$. And this is how it is best learned, with group theory knowledge making you powerful.
\end{principle}

But probably enough talking and principles, with black belt mathematical power at stake, believe me here, let us get interested in groups, and do some work. We will be first interested in the symmetric group $S_N$, and more specifically, in counting the permutations $\sigma\in S_N$ having no fixed points, which are called derangements. And here, we have:

\index{derangement}

\begin{theorem}
The probability for a random $\sigma\in S_N$ to be a derangement is
$$P\simeq\frac{1}{e}$$
in the $N\to\infty$ limit.
\end{theorem}

\begin{proof}
Consider the sets $S_N^i=\{\sigma\in S_N|\sigma(i)=i\}$. According to the inclusion-exclusion principle, the probability that we are interested in is given by:
\begin{eqnarray*}
P
&=&\frac{1}{N!}\left(|S_N|-\sum_i|S_N^i|+\sum_{i<j}|S_N^i\cap S_N^j|-\sum_{i<j<k}|S_N^i\cap S_N^j\cap S_N^k|+\ldots\right)\\
&=&\frac{1}{N!}\sum_{k=0}^N(-1)^k\sum_{i_1<\ldots<i_k}|S_N^{i_1}\cap\ldots\cap S_N^{i_k}|\\
&=&\frac{1}{N!}\sum_{k=0}^N(-1)^k\binom{N}{k}(N-k)!\\
&=&\sum_{k=0}^N\frac{(-1)^k}{k!}
\end{eqnarray*}

Thus, we are led to the conclusion in the statement.
\end{proof}

With a bit more work, this leads to a combinatorial model for the Poisson law $p_1$:

\index{fixed points}

\begin{theorem}[upgrade]
The number of fixed points of permutations,
$$\chi:S_N\to\mathbb N\quad,\quad \chi(\sigma)=\#\left\{i\in\{1,\ldots,N\}\Big|\sigma(i)=i\right\}$$
follows with $N\to\infty$ the Poisson law $p_1$.
\end{theorem}

\begin{proof}
We have to count the permutations $\sigma\in S_N$ having exactly $k$ points. Since having such a permutation amounts in choosing $k$ points among $1,\ldots,N$, and then permuting the $N-k$ points left, without fixed points allowed, we have:
$$\#\left\{\sigma\in S_N\Big|\chi(\sigma)=k\right\}
=\binom{N}{k}\#\left\{\sigma\in S_{N-k}\Big|\chi(\sigma)=0\right\}$$

Now by dividing everything by $N!$, we obtain from this the following formula:
$$\frac{\#\left\{\sigma\in S_N\Big|\chi(\sigma)=k\right\}}{N!}=\frac{1}{k!}\times\frac{\#\left\{\sigma\in S_{N-k}\Big|\chi(\sigma)=0\right\}}{(N-k)!}$$

By using now the computation at $k=0$, from Theorem 4.3, it follows that with $N\to\infty$ we have the following estimate, $\chi$ being the number of fixed points:
$$P(\chi=k)
\simeq\frac{1}{k!}\cdot P(\chi=0)
\simeq\frac{1}{k!}\cdot\frac{1}{e}$$

Thus, we are led to the conclusion in the statement.
\end{proof}

More generally now, and quite remarkably, we have in fact the following result:

\index{symmetric group}
\index{Poisson law}

\begin{theorem}[upgrade]
The number of partial fixed points of permutations,
$$\chi_t:S_N\to\mathbb N\quad,\quad \chi(\sigma)=\#\left\{i\in\{1,\ldots,[tN]\}\Big|\sigma(i)=i\right\}$$
follows with $N\to\infty$ the Poisson law $p_t$.
\end{theorem}

\begin{proof}
This is an upgrade of Theorem 4.4, the idea being as follows:

\medskip

(1) Consider first, as in the proof of Theorem 4.3, the sets $S_N^i=\{\sigma\in S_N|\sigma(i)=i\}$. As before in the proof of Theorem 4.3, we obtain by inclusion-exclusion that:
\begin{eqnarray*}
P(\chi_t=0)
&=&\frac{1}{N!}\sum_{k=0}^{[tN]}(-1)^k\sum_{i_1<\ldots<i_k<[tN]}|S_N^{i_1}\cap\ldots\cap S_N^{i_k}|\\
&=&\frac{1}{N!}\sum_{k=0}^{[tN]}(-1)^k\binom{[tN]}{k}(N-k)!\\
&=&\sum_{k=0}^{[tN]}\frac{(-1)^k}{k!}\cdot\frac{[tN]!(N-k)!}{N!([tN]-k)!}
\end{eqnarray*}

But with $N\to\infty$, we obtain from this the following estimate, as desired:
$$P(\chi_t=0)
\simeq\sum_{k=0}^{[tN]}\frac{(-1)^k}{k!}\cdot t^k
=\sum_{k=0}^{[tN]}\frac{(-t)^k}{k!}
\simeq e^{-t}$$

(2) More generally now, by counting the permutations $\sigma\in S_N$ having exactly $k$ fixed points among $1,\ldots,[tN]$, as in the proof of Theorem 4.4, our claim is that we get:
$$P(\chi_t=k)\simeq\frac{t^k}{k!e^t}$$

Indeed, we already know from (1) that this formula holds at $k=0$. In general now, we have to count the permutations $\sigma\in S_N$ having exactly $k$ fixed points among $1,\ldots,[tN]$. Since having such a permutation amounts in choosing $k$ points among $1,\ldots,[tN]$, and then permuting the $N-k$ points left, without fixed points among $1,\ldots,[tN]$ allowed, we obtain the following formula, where $s\in(0,1]$ is such that $[s(N-k)]=[tN]-k$:
$$\#\left\{\sigma\in S_N\Big|\chi_t(\sigma)=k\right\}
=\binom{[tN]}{k}\#\left\{\sigma\in S_{N-k}\Big|\chi_s(\sigma)=0\right\}$$

Now by dividing everything by $N!$, we obtain from this the following formula:
$$\frac{\#\left\{\sigma\in S_N\Big|\chi_t(\sigma)=k\right\}}{N!}=\frac{1}{k!}\times\frac{[tN]!(N-k)!}{N!([tN]-k)!}\times\frac{\#\left\{\sigma\in S_{N-k}\Big|\chi_s(\sigma)=0\right\}}{(N-k)!}$$

By using now the computation at $k=0$, that we have from (1), we obtain:
\begin{eqnarray*}
P(\chi_t=k)
&\simeq&\frac{1}{k!}\times\frac{[tN]!(N-k)!}{N!([tN]-k)!}\cdot P(\chi_s=0)\\
&\simeq&\frac{t^k}{k!}\cdot P(\chi_s=0)\\
&\simeq&\frac{t^k}{k!}\cdot\frac{1}{e^s}
\end{eqnarray*}

Now recall that $s\in(0,1]$ was chosen such that $[s(N-k)]=[tN]-k$. Thus with $N\to\infty$ we have $s=t$, so we obtain $P(\chi_t=k)\simeq t^k/(k!e^t)$, which gives the result.
\end{proof}

And with this, end of the story? You must be kidding. We have indeed:

\index{truncated character}

\begin{theorem}[upgrade]
For the symmetric group $S_N\subset_uO_N$, the truncated character
$$\chi_t:G\to\mathbb R\quad,\quad 
\chi_t=\sum_{i=1}^{[tN]}u_{ii}$$
follows with $N\to\infty$ the Poisson law $p_t$.
\end{theorem}

\begin{proof}
Let us view indeed $S_N$ as group of permutations of the $N$ coordinate axes of $\mathbb R^N$. The formula of the action is then $u(\sigma):e_j\to e_{\sigma(j)}$, which in matrix terms reads:
$$u_{ij}(\sigma)=\delta_{\sigma(j)i}$$

Thus, the sums of diagonal coordinates $u_{ii}$ count the fixed points, and with this observation in hand, Theorem 4.5 reformulates as in the statement.
\end{proof}

Next, now that we have a nice final statement, let us upgrade as well the proof:

\index{polynomial integrals}

\begin{theorem}[upgrade]
The polynomial integrals over $S_N\subset_uO_N$ are given by
$$\int_{S_N}u_{i_1j_1}\ldots u_{i_kj_k}=\begin{cases}
\frac{(N-|\ker i|)!}{N!}&{\rm if}\ \ker i=\ker j\\
0&{\rm otherwise}
\end{cases}$$
and by using this and summing, we obtain $\chi_t\sim p_t$ in the $N\to\infty$ limit.
\end{theorem}

\begin{proof}
We have two assertions here, the idea being as follows:

\medskip

(1) According to our convention for $S_N\subset_uO_N$, the integrals in the statement are:
$$\int_{S_N}u_{i_1j_1}\ldots u_{i_kj_k}=\frac{1}{N!}\#\left\{\sigma\in S_N\Big|\sigma(j_1)=i_1,\ldots,\sigma(j_k)=i_k\right\}$$

Now observe that the existence of $\sigma\in S_N$ as above requires $i_m=i_n\iff j_m=j_n$. Thus, the above integral vanishes when the following condition is satisfied:
$$\ker i\neq\ker j$$

Regarding now the case $\ker i=\ker j$, if we denote by $b\in\{1,\ldots,k\}$ the number of blocks of this partition $\ker i=\ker j$, we have $N-b$ points to be sent bijectively to $N-b$ points, and so $(N-b)!$ solutions, and the integral is $\frac{(N-b)!}{N!}$, as claimed.

\medskip

(2) For the second assertion, with $S_{kb}$ being the Stirling numbers, counting the partitions of $\{1,\ldots,k\}$ having exactly $b$ blocks, we have the following computation:
\begin{eqnarray*}
\int_{S_N}\chi_t^k
&=&\sum_{i_1,\ldots,i_k=1}^{[tN]}\int_{S_N}u_{i_1i_1}\ldots u_{i_ki_k}\\
&=&\sum_{\pi\in P(k)}\frac{[tN]!}{([tN]-|\pi|!)}\cdot\frac{(N-|\pi|!)}{N!}\\
&=&\sum_{b=1}^{[tN]}\frac{[tN]!}{([tN]-b)!}\cdot\frac{(N-b)!}{N!}\cdot S_{kb}
\end{eqnarray*}

In particular with $N\to\infty$ we obtain the following formula:
$$\lim_{N\to\infty}\int_{S_N}\chi_t^k=\sum_{b=1}^kS_{kb}t^b$$

But this is the $k$-th moment of the Poisson law $p_t$, and so we are done.
\end{proof}

So, are we done with this? Not really, because here is something even better:

\begin{theorem}[upgrade]
We have the following law formula,
$$law(\chi_t)=\frac{[tN]!}{N!}\sum_{p=0}^{[tN]}\frac{(N-p)!}{([tN]-p)!}
\cdot\frac{\left(\delta_1-\delta_0\right)^{*p}}{p!}$$ 
which shows in particular that $\chi_t\sim p_t$ in the $N\to\infty$ limit.
\end{theorem}

\begin{proof}
We have two assertions here, the idea being as follows:

\medskip

(1) The law formula follows either from the computations in the proof of Theorem 4.5, or from those in the proof of Theorem 4.7, by extracting from there the last formula, just before letting $N\to\infty$, and further processing it. We will leave this as an exercise.

\medskip

(2) The coefficients in the law formula can be estimated via Stirling, as follows:
$$c_p
\simeq\frac{(tN)^{tN}}{N^N}\cdot\frac{(N-p)^{N-p}}{(tN-p)^{tN-p}}
\simeq t^p$$

We conclude that the asymptotic law is given by the following formula:
$$law(\chi_t)\simeq \sum_{p=0}^st^p\cdot\frac{\left(\delta_1-\delta_0\right)^{*p}}{p!}$$

But on the right we have $p_t$, say via $F(x)=\exp(t(e^{ix}-1))$, and we are done.
\end{proof}

And we will end our study of derangements with this, guess we are now experts in counting. By the way, in case you got addicted to our series of upgrades, I have news for you, there will be one more, based on the easiness property of $S_N$. Coming soon. 

\section*{4b. Weingarten formula}

Getting now to the orthogonal group $O_N$, we would like to prove that $\chi_t\sim g_t$ with $N\to\infty$, which would stand as a good complement to our $\chi_t\sim p_t$ formula for $S_N$, and would allow us to have some business started, along the lines of Principles 4.1 and 4.2.

\bigskip

However, and coming as bad news, while in the case of $S_N$ we were totally at ease, with several counting methods being available, in the case of $O_N$, which is a continuous group, none of these discrete counting methods will apply, and we are in the dark.

\bigskip

So, in order to get started, we need a theorem regarding the integration over $O_N$, telling us what that integration is. And here is the theorem that we will need:

\index{Haar integration}

\begin{theorem}
The integration over a compact group $G\subset O_N$ can be constructed by starting with any faithful positive unital linear form $\varphi\in C(G)^*$, and setting:
$$\int_G=\lim_{n\to\infty}\frac{1}{n}\sum_{k=1}^n\varphi^{*k}$$
Moreover, for any representation $v:G\to O_n$ we have the following formula,
$$\left(\int_Gv_{ij}\right)_{ij}=P$$
where $P$ is the orthogonal projection onto $Fix(v)=\left\{\xi\in\mathbb C^n\big|v(g)\xi=\xi,\forall g\in G\right\}$.
\end{theorem}

\begin{proof}
This is something standard, which can be done in 3 steps, as follows:

\medskip

(1) Given $\varphi\in C(G)^*$, our claim is that the following converges, for any $f\in C(G)$:
$$\int_\varphi f=\lim_{n\to\infty}\frac{1}{n}\sum_{k=1}^n\varphi^{*k}(f)$$

Indeed, by linearity we can assume that $f$ is the coefficient of certain representation $v:G\to O_n$. But in this case, an elementary computation gives the following formula, with $P_\varphi\geq P$ being the orthogonal projection onto the $1$-eigenspace of $[\varphi(v_{ij})]_{ij}$:
$$\left(\int_\varphi v_{ij}\right)_{ij}=P_\varphi$$

(2) The point now is that when $\varphi\in C(G)^*$ is faithful, by using a standard positivity trick, we can prove that we have $P_\varphi=P$. Assume indeed $P_\varphi\xi=\xi$, and let us set:
$$f=\sum_i\left(\sum_jv_{ij}\xi_j-\xi_i\right)\overline{\left(\sum_kv_{ik}\xi_k-\xi_i\right)}$$

A straightforward computation shows then that $\varphi(f)=0$, and so $f=0$, as desired.

\medskip

(3) Thus, we proved our claim, and with this in hand, the left and right invariance of $\int_G=\int_\varphi$ is clear on coefficients, and so in general, and this gives all the assertions. For details here you can check my book \cite{ba1}, where all this is explained in detail.
\end{proof}

Getting back now to $O_N$, having Theorem 4.9 in hand is not the end of our troubles, because we still have to know how to use the formulae there. So, here is my proposal:

\begin{proposal}
We will first attempt to reprove $\chi_t\sim p_t$ for $S_N$, by using the technology from Theorem 4.9, even if this is something a bit unnatural, which might take some time. Then, once this understood, we will prove $\chi_t\sim g_t$ for $O_N$ too.
\end{proposal}

So, proposal accepted I hope, please cancel all your meetings until the holidays, turn off all electric machinery on you, and here we go again with $S_N$, somehow with the aim of inventing the most complicated proof ever for $\chi_t\sim p_t$. Let us start with:

\begin{proposition}
Given a compact group $G\subset_uO_N$, we have
$$\int_G\chi^k=\dim\left(Fix(u^{\otimes k})\right)$$
where $u^{\otimes k}:G\to O_{N^k}$ are the tensor powers of $u:G\to O_N$.
\end{proposition}

\begin{proof}
In the context of Theorem 4.9, by applying the trace there we obtain:
$$\int_GTr(v)=Tr(P)=\dim(Im(P))=\dim(Fix(v))$$

Now for the representation $v=u^{\otimes k}$, this gives the formula in the statement.
\end{proof}

Summarizing, in order to compute the character law for $G\subset_uO_N$, we must find bases for the spaces $Fix(u^{\otimes k})$. And getting now to the case of $G=S_N$, fortunately these spaces $Fix(u^{\otimes k})$ have a very simple description, in terms of partitions, as shown by:

\index{easiness}

\begin{theorem}
The symmetric group $S_N\subset_uO_N$ has the easiness property
$$Hom(u^{\otimes k},u^{\otimes l})=span\left(T_\pi\Big|\pi\in P(k,l)\right)$$
with $P(k,l)$ being the set of partitions between $k$ points and $l$ points, and
$$T_\pi(e_{i_1}\otimes\ldots\otimes e_{i_k})=\sum_{j_1\ldots j_l}\delta_\pi
\begin{pmatrix}i_1&\ldots& i_k\\ j_1&\ldots& j_l\end{pmatrix}e_{j_1}\otimes\ldots\otimes e_{j_l}$$
where $\delta_\pi=1$ if the indices match, and $\delta_\pi=0$ otherwise.
\end{theorem}

\begin{proof}
This is something worth discussing in detail, the idea being as follows:

\medskip

(1) To start with, forgetting about symmetric groups, and being a bit philosophers, what makes a compact group $G\subset_uO_N$ easy? Having 1 element, or being cyclic, or abelian, you would say. However, this is something naive. Indeed, at the advanced level, as we know from Proposition 4.11, this means that the spaces $Fix(u^{\otimes k})$ must be easy.

\medskip

(2) Next, observe that the representations $(u^{\otimes k})_{k\geq0}$ form a category, with the arrows between them being the intertwiners, which are as follows, with $H=\mathbb C^N$:
$$Hom(u^{\otimes k},u^{\otimes l})=\left\{T\in\mathcal L(H^{\otimes k},H^{\otimes l})\Big|Tg^{\otimes k}=g^{\otimes l}T,\forall g\in G\right\}$$

So, in practice, it is this collection of linear spaces that you want to be easy.

\medskip

(3) Then, what should this latter easiness property exactly mean? In answer, we would like the intertwiners $T$ to be as simple as possible, and some good candidates here are the linear maps in the statement, associated to the partitions $\pi\in P(k,l)$.

\medskip

(4) To be more precise, let us denote by $P(k,l)$ the set of usual partitions, between an upper row of $k$ points, and a lower row of $l$ points. As an example, here is a partition in $P(3,3)$, with two blocks, represented by strings, in the obvious way:
$$
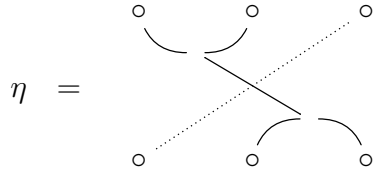
\xymatrix@R=2mm@C=4mm{\\ \\ \eta\ \ =\\ \\}\ \ \ 
\xymatrix@R=2mm@C=4mm{
\circ\ar@/_/@{-}[dr]&&\circ&&\circ\ar@{.}[ddddllll]\\
&\ar@/_/@{-}[ur]\ar@{-}[ddrr]\\
\\
&&&\ar@/^/@{-}[dr]\\
\circ&&\circ\ar@/^/@{-}[ur]&&\circ}$$

(5) Now given $\pi\in P(k,l)$ and multi-indices $i=(i_1,\ldots,i_k)$ and $j=(j_1,\ldots,j_l)$, we can put $i,j$ on the legs of $\pi$, in the obvious way. Then, if the indices fit, meaning that all strings of $\pi$ join equal indices of $i,j$, we set $\delta_\pi\binom{i}{j}=1$. Otherwise, we set $\delta_\pi\binom{i}{j}=0$. As an example, for the above partition $\eta\in P(3,3)$, we have the following formula:
$$\delta_\eta\begin{pmatrix}a&b&c\\d&e&f\end{pmatrix}=\delta_{abef}\delta_{cd}$$

(6) With this done, we can talk about the linear maps $T_\pi$ in the statement. As an example here, for the above partition $\eta\in P(3,3)$, the formula is as follows:
$$T_\eta(e_a\otimes e_b\otimes e_c)=\delta_{ab}e_c\otimes e_a\otimes e_a$$

(7) And with this, we can now axiomatize easiness. Indeed, we can say that a closed subgroup $G\subset_uO_N$ is easy when we have equalities of linear spaces as follows, for any two integers $k,l\in\mathbb N$, with $D(k,l)\subset P(k,l)$ being certain sets of partitions:
$$Hom(u^{\otimes k},u^{\otimes l})=span\left(T_\pi\Big|\pi\in D(k,l)\right)$$

(8) Getting to work now, our theorem says that $S_N\subset_uO_N$ is easy, coming from $D(k,l)=P(k,l)$, which in practice means that we have, for any $k,l\in\mathbb N$:
$$Hom(u^{\otimes k},u^{\otimes l})=span\left(T_\pi\Big|\pi\in P(k,l)\right)$$

(9) Let us first prove $\supset$. Given $\pi\in P(k,l)$ and $\sigma\in S_N$, we have indeed:
\begin{eqnarray*}
T_\pi\sigma^{\otimes k}(e_{i_1}\otimes\ldots\otimes e_{i_k})
&=&\sum_{j_1\ldots j_l}\delta_\pi
\begin{pmatrix}\sigma(i_1)&\ldots&\sigma(i_k)\\ j_1&\ldots& j_l\end{pmatrix}e_{j_1}\otimes\ldots\otimes e_{j_l}\\
&=&\sum_{j_1\ldots j_l}\delta_\pi
\begin{pmatrix}\sigma(i_1)&\ldots&\sigma(i_k)\\ \sigma(j_1)&\ldots&\sigma(j_l)\end{pmatrix}e_{\sigma(j_1)}\otimes\ldots\otimes e_{\sigma(j_l)}\\
&=&\sum_{j_1\ldots j_l}\delta_\pi
\begin{pmatrix}i_1&\ldots& i_k\\ j_1&\ldots& j_l\end{pmatrix}
e_{\sigma(j_1)}\otimes\ldots\otimes e_{\sigma(j_l)}\\
&=&\sigma^{\otimes l}T_\pi(e_{i_1}\otimes\ldots\otimes e_{i_k})
\end{eqnarray*}

(10) In order to prove now the reverse inclusion $\subset$, consider an arbitrary linear map $T:(\mathbb C^N)^{\otimes k}\to(\mathbb C^N)^{\otimes l}$, written as follows, with $\lambda\binom{i_1\ldots i_k}{j_1\ldots j_l}\in\mathbb C$ being certain scalars:
$$T(e_{i_1}\otimes\ldots\otimes e_{i_k})=\sum_{j_1\ldots j_l}\lambda
\begin{pmatrix}i_1&\ldots& i_k\\ j_1&\ldots& j_l\end{pmatrix}e_{j_1}\otimes\ldots\otimes e_{j_l}$$

Given a permutation $\sigma\in S_N$, by reasoning as before, we can see that we have:
$$T\sigma^{\otimes k}=\sigma^{\otimes l}T
\ \iff\ \lambda\begin{pmatrix}\sigma(i_1)&\ldots&\sigma(i_k)\\ \sigma(j_1)&\ldots&\sigma(j_l)\end{pmatrix}=\lambda\begin{pmatrix}i_1&\ldots& i_k\\ j_1&\ldots& j_l\end{pmatrix},\ \forall i,j$$

We conclude that we have $T\in Hom(u^{\otimes k},u^{\otimes l})$ precisely when the following happens:
$$\ker\begin{pmatrix}i_1&\ldots& i_k\\ j_1&\ldots& j_l\end{pmatrix}
=\ker\begin{pmatrix}i'_1&\ldots& i'_k\\ j'_1&\ldots& j'_l\end{pmatrix}\implies
\lambda\begin{pmatrix}i_1&\ldots& i_k\\ j_1&\ldots& j_l\end{pmatrix}
=\lambda\begin{pmatrix}i_1&\ldots& i_k\\ j_1&\ldots& j_l\end{pmatrix}$$

Thus, we are led to the conclusion that $\lambda:\{1,\ldots,N\}^{k+l}\to\mathbb C$ must come from a function $\varphi:P(k,l)\to\mathbb C$, via a formula of type $\lambda(x)=\varphi(\ker x)$, and it follows that the inclusion $\supset$ that we established in (9) is indeed an equality, as desired.
\end{proof}

With this done, it remains to see if this technology is any good for our probability questions. And here, according to Proposition 4.11, and then to easiness, we have:
$$\int_{S_N}\chi_1^k=\dim\left(Fix(u^{\otimes k})\right)
=\dim\left[span\left(T_\pi\Big|\pi\in P(k)\right)\right]$$

Thus, we are done with $t=1$, save for the fact that we do not know yet that the vectors on the right are linearly independent with $N\to\infty$. But this comes from:

\index{Lindst\"om formula}
\index{Gram determinant}

\begin{theorem}[Lindst\"om]
The Gram matrix of the vectors $\{T_\pi|\pi\in P(k)\}$ is
$$G_{kN}(\pi,\nu)=N^{|\pi\vee\nu|}$$
where $\vee$ is the superposition operation for the partitions, and we have the formula
$$\det(G_{kN})=\prod_{\pi\in P(k)}\frac{N!}{(N-|\pi|)!}$$
at $N\geq k$. In particular, the vectors $\{T_\pi|\pi\in P(k)\}$ are linearly independent at $N\geq k$.
\end{theorem}

\begin{proof}
Regarding the first formula, this comes from the following computation:
\begin{eqnarray*}
G_{kN}(\pi,\nu)
&=&\sum_{j_1\ldots j_l}
\delta_\pi\begin{pmatrix}j_1&\ldots& j_l\end{pmatrix}
\delta_\nu\begin{pmatrix}j_1&\ldots& j_l\end{pmatrix}\\
&=&\sum_{j_1\ldots j_l}
\delta_{\pi\vee\nu}\begin{pmatrix}j_1&\ldots& j_l\end{pmatrix}\\
&=&N^{|\pi\vee\nu|}
\end{eqnarray*}

Next, in order to compute $\det(G_{kN})$, observe that we have the following formula:
\begin{eqnarray*}
N^{|\pi\vee\nu|}
&=&\#\left\{i_1,\ldots,i_k\in\{1,\ldots,N\}\Big|\ker i\geq\pi\vee\nu\right\}\\
&=&\sum_{\tau\geq\pi\vee\nu}\#\left\{i_1,\ldots,i_k\in\{1,\ldots,N\}\Big|\ker i=\tau\right\}\\
&=&\sum_{\tau\geq\pi\vee\nu}N(N-1)\ldots(N-|\tau|+1)
\end{eqnarray*}

Now in view of this latter formula, consider the following matrix:
$$L_{kN}(\pi,\nu)=
\begin{cases}
N(N-1)\ldots(N-|\pi|+1)&{\rm if}\ \nu\leq\pi\\
0&{\rm otherwise}
\end{cases}$$

If we denote by $A_k$ the adjacency matrix of $P(k)$, we conclude that we have:
$$G_{kN}(\pi,\nu)
=\sum_{\tau\geq\pi}L_{kN}(\tau,\nu)
=\sum_\tau A_k(\pi,\tau)L_{kN}(\tau,\nu)
=(A_kL_{kN})(\pi,\nu)$$

Summarizing, we have decomposed our Gram matrix as $G_{kN}=A_kL_{kN}$. Now if we order $P(k)$ as usual, with respect to the number of blocks, and then lexicographically, the matrix $A_k$ is upper triangular, and the matrix $L_{kN}$ is lower triangular, and so:
$$\det(G_{kN})
=\det(A_k)\det(L_{kN})
=\det(L_{kN})
=\prod_{\pi\in P(k)}\frac{N!}{(N-|\pi|)!}$$

And with this done, we are led to the conclusions in the statement. 
\end{proof}

Getting back to $S_N$, as explained before Theorem 4.13, we are done with $t=1$. In order to discuss now the general case $t\in(0,1]$, we will need the Weingarten formula:

\index{Weingarten formula}

\begin{theorem}
For an easy group $G\subset_uO_N$, coming from $D\subset P$, we have
$$\int_Gu_{i_1j_1}\ldots u_{i_kj_k}=\sum_{\pi,\nu\in D(k)}\delta_\pi(i)\delta_\nu(j)W_{kN}(\pi,\nu)$$
with $W_{kN}=G_{kN}^{-1}$ being the inverse of the Gram matrix $G_{kN}(\pi,\nu)=N^{|\pi\vee\nu|}$.
\end{theorem}

\begin{proof}
Consider the integrals in the statement, denoted as follows:
$$P_{i_1\ldots i_k,j_1\ldots j_k}=\int_Gu_{i_1j_1}\ldots u_{i_kj_k}$$

We know from Theorem 4.9 that the matrix $P=(P_{ij})$ is the projection on $Fix(u^{\otimes k})$. On the other hand, by easiness, this space that we are projecting on is given by:
$$Fix(u^{\otimes k})=span\left(T_\pi\Big|\pi\in D(k)\right)$$

In order now to explicitly compute $P$, consider the following linear map:
$$E(x)=\sum_{\pi\in D(k)}<x,T_\pi>T_\pi$$

By linear algebra we have then $P=WE$, where $W$ is the inverse on $Fix(u^{\otimes k})$ of the restriction of $E$. But this latter restriction is the linear map given by the Gram matrix $G_{kN}$, and we are led in this way to the formula in the statement.
\end{proof}

We can go back now to our symmetric group considerations, and we have the following new and final upgrade of our series of results, from the previous section:

\begin{theorem}[again]
The symmetric group $S_N\subset_uO_N$ is easy,
$$Hom(u^{\otimes k},u^{\otimes l})=span\left(T_\pi\Big|\pi\in P(k,l)\right)$$
and this shows that we have $\chi_t\sim p_t$ with $N\to\infty$.
\end{theorem}

\begin{proof}
In the case $t=1$ we are already done, as explained before Theorem 4.13. In the general case now, $t\in(0,1]$, we have the following computation:
\begin{eqnarray*}
\int_{S_N}(u_{11}+\ldots +u_{ss})^k
&=&\sum_{i_1=1}^{s}\ldots\sum_{i_k=1}^s\int_{S_N}u_{i_1i_1}\ldots u_{i_ki_k}\\
&=&\sum_{\pi,\nu\in D(k)}W_{kN}(\pi,\nu)\sum_{i_1=1}^{s}\ldots\sum_{i_k=1}^s\delta_\pi(i)\delta_\nu(i)\\
&=&\sum_{\pi,\nu\in D(k)}W_{kN}(\pi,\nu)G_{ks}(\nu,\pi)\\
&=&Tr(W_{kN}G_{ks})
\end{eqnarray*}

Now since the Gram matrix is asymptotically diagonal, $G_{kN}\simeq diag(N^{|\pi|})$, its inverse is asymptotically diagonal too, $W_{kN}\simeq diag(N^{-|\pi|})$, and the above computation gives:
$$\lim_{N\to\infty}\int_{S_N}\chi_t^k=\sum_{\pi\in P(k)}t^{|\pi|}$$

We conclude that we have $\chi_t\sim p_t$ in the $N\to\infty$ limit, as stated.
\end{proof}

Getting now to the case of the orthogonal group $O_N$, we have here:

\begin{theorem}
The orthogonal group $O_N$ is easy too, according to
$$Hom(u^{\otimes k},u^{\otimes l})=span\left(T_\pi\Big|\pi\in P_2(k,l)\right)$$
with $P_2$ standing for pairings, and with this, we have $\chi_t\sim g_t$ in the $N\to\infty$ limit.
\end{theorem}

\begin{proof}
The study here is similar to the one for $S_N$, with a few new twists:

\medskip

(1) In what regards easiness, consider the semicircle pairing $\cap$. The linear map associated to it being $T_\cap=\sum_ie_i\otimes e_i$, for an arbitrary matrix $U\in M_N(\mathbb R)$, we have: 
\begin{eqnarray*}
U^{\otimes 2}T_\cap=T_\cap
&\iff&\sum_{ijk}U_{ji}U_{ki}e_j\otimes e_k=\sum_ie_i\otimes e_i\\
&\iff&\sum_iU_{ji}U_{ki}=\delta_{jk},\ \forall j,k\\
&\iff&U\in O_N
\end{eqnarray*}

We conclude that for $O_N$ we have $T_\cap\in Fix(u^{\otimes 2})$, and with a bit more work, this gives the inclusion $\supset$ in the statement. To be more precise, you can argue here that the category of representations of $O_N$ containing $T_\pi$, it must contain the whole category $<T_\pi>=span(T_\pi|\pi\in P_2)$ that this arrow generates. Alternatively, you can of course check $\supset$ via a direct computation, variation of the above check of $T_\cap\in Fix(u^{\otimes 2})$.

\medskip

(2) Getting now to $\subset$, this is harder to prove with bare hands, but at the advanced level we can say that this is trivial, using Tannakian duality. Indeed, the spaces on the right in the statement form a category, and more specifically satisfy the axioms of Tannakian categories, so by duality, they must correspond to a certain subgroup $G\subset O_N$. But the computation in (1) shows that this subgroup must be $G=O_N$ itself, as desired.

\medskip

(3) Finally, the probability computations are as before for $S_N$, with, at the end:
$$\int_{O_N}\chi_t^k=Tr(W_{kN}G_{k[tN]})\simeq\sum_{\pi\in P_2(k)}t^{|\pi|}$$

We conclude from this that we have $\chi_t\sim g_t$ with $N\to\infty$, as stated.
\end{proof}

Summarizing, with respect to our goals in this chapter, as formulated in Principles 4.1 and 4.2, we have now a good understanding of $p_t,g_t$. Our plan in what follows will be to discuss $p^2_t,g_t^t$, and then leave $\pi_t,\gamma_t,\pi^2_t,\gamma_t^t$ for later, when doing free probability.

\section*{4c. Reflections, Bessel} 

As a continuation of the above material, which was quite exciting, let us do some more computations, for other interesting groups of matrices $G\subset O_N$ that we know. And here, I don't know about you, but my favorite example is the hyperoctahedral group:

\index{hyperoctahedral group}

\begin{theorem}
Consider the hyperoctahedral group $H_N$, which appears as the symmetry group of the $N$-dimensional cube, viewed as graph, or as metric space,
$$
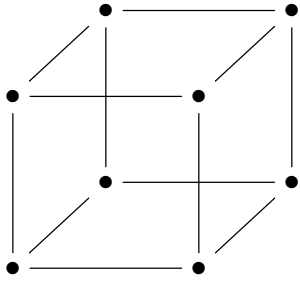
\xymatrix@R=20pt@C=22pt{
&\bullet\ar@{-}[rr]&&\bullet\\
\bullet\ar@{-}[rr]\ar@{-}[ur]&&\bullet\ar@{-}[ur]\\
&\bullet\ar@{-}[rr]\ar@{-}[uu]&&\bullet\ar@{-}[uu]\\
\bullet\ar@{-}[uu]\ar@{-}[ur]\ar@{-}[rr]&&\bullet\ar@{-}[uu]\ar@{-}[ur]
}$$
or equivalently, which is the symmetry group of the $N$ coordinate axes of $\mathbb R^N$:
$$S_N\subset H_N\subset O_N$$
In matrix terms, $H_N$ consists of the permutation-type matrices having $\pm1$ as nonzero entries, and we have a wreath product decomposition as follows:
$$H_N=\mathbb Z_2\wr S_N$$
In this picture, the main character counts the signed number of fixed points, among the coordinate axes, and its truncations count the truncations of such numbers.
\end{theorem}

\begin{proof}
Many things going on here, the idea being as follows:

\medskip

(1) To start with, we can certainly talk about the hyperoctahedral group $H_N$, as being the symmetry group of the $N$-cube, and in small dimensions we have $H_1=\mathbb Z_2$, obvious, then $H_2=D_4=\mathbb Z_4\rtimes\mathbb Z_2$, the symmetry group of the square, and then $H_3=S_4\times\mathbb Z_2$, with the copy of $S_4$ best viewed as permuting the main diagonals of the cube.

\medskip

(2) Next, by centering the $N$-cube at the origin, its symmetry group $H_N$ appears as the symmetry group of the $N$ coordinate axes of $\mathbb R^N$. Which is a quite fruitful viewpoint, proving right away $H_N=M_N(-1,0,1)\cap O_N$, and giving as well the cardinality formula $|H_N|=2^NN!$, which after a bit more work gives $H_N=\mathbb Z_2\wr S_N$, as stated.

\medskip

(3) Finally, in what regards the assertions at the end, concerning the main character and its truncations, these are clear from definitions, exactly as previously for $S_N$. Let us also mention that $H_N$ is easy, but more on this later, when we will really need it.
\end{proof}

Regarding now the character laws, we can compute them by using the same method as for the symmetric group $S_N$, namely inclusion-exclusion, and we have:

\index{main character}

\begin{theorem}
For the hyperoctahedral group $H_N\subset O_N$, the law of the variable
$$\chi_t=\sum_{i=1}^{[tN]}u_{ii}$$
with $u_{ij}$ being the standard coordinates, becomes in the $N\to\infty$ limit the measure
$$p^2_t=e^{-t}\sum_{k=-\infty}^\infty\delta_k\sum_{p=0}^\infty \frac{(t/2)^{|k|+2p}}{(|k|+p)!p!}$$ 
called Bessel law of parameter $t\in(0,1]$.
\end{theorem}

\begin{proof}
This is something very standard, the idea being as follows:

\medskip

(1) We can regard $H_N$ as being the symmetry group of the graph $I_N=\{I^1,\ldots ,I^N\}$ formed by $N$ segments. In this picture, the diagonal coefficients are given by:
$$u_{ii}(g)=\begin{cases}
\ 0\ \mbox{ if $g$ moves $I^i$}\\
\ 1\ \mbox{ if $g$ fixes $I^i$}\\
-1\mbox{ if $g$ returns $I^i$}
\end{cases}$$

We denote by $\uparrow g,\downarrow g$ the number of segments among $\{I^1,\ldots ,I^s\}$ which are fixed, respectively returned by an element $g\in H_N$. With this notation, we have:
$$u_{11}+\ldots+u_{ss}=\uparrow g-\downarrow g$$

Let us denote by $P_N$ probabilities computed over the group $H_N$. The density of the law of $u_{11}+\ldots+u_{ss}$ at a point $k\geq 0$ is then given by the following formula:
\begin{eqnarray*}
D(k)
&=&P_N(\uparrow g-\downarrow g=k)\\
&=&\sum_{p=0}^\infty P_N(\uparrow g=k+p, \downarrow g=p)
\end{eqnarray*}

Also, we have $D(-k)=D(k)$, so it is enough to do the computation for $k\geq0$.

\medskip

(2) Let us first discuss the case $t=1$. We use the fact, that we know well from before, that the probability of $\sigma\in S_N$ to have $m$ fixed points is asymptotically given by $P_m=\frac{1}{em!}$. In terms of probabilities over $H_N$, we obtain from this, as desired:
\begin{eqnarray*}
\lim_{N\to\infty}D(k)
&=&\lim_{N\to\infty}\sum_{p=0}^\infty(1/2)^{k+2p}\begin{pmatrix}k+2p\\ k+p\end{pmatrix} P_N(\uparrow g+\downarrow g=k+2p)\\ 
&=&\sum_{p=0}^\infty(1/2)^{k+2p}\begin{pmatrix}k+2p\\
k+p\end{pmatrix}\frac{1}{e(k+2p)!}\\
&=&\frac{1}{e}\sum_{p=0}^\infty \frac{(1/2)^{k+2p}}{(k+p)!p!}
\end{eqnarray*}

(3) As for the general case $t\in(0,1]$, here the result follows by performing some modifications in the above computation, the density being computed as follows:
\begin{eqnarray*}
\lim_{N\to\infty}D(k)
&=&\lim_{N\to\infty}\sum_{p=0}^\infty(1/2)^{k+2p}\begin{pmatrix}k+2p\\ k+p\end{pmatrix} P_N(\uparrow g+\downarrow g=k+2p)\\
&=&\sum_{p=0}^\infty(1/2)^{k+2p}\begin{pmatrix}k+2p\\
k+p\end{pmatrix}\frac{t^{k+2p}}{e^t(k+2p)!}\\
&=&e^{-t}\sum_{p=0}^\infty \frac{(t/2)^{k+2p}}{(k+p)!p!}
\end{eqnarray*}

Thus, we are led to the conclusion in the statement.
\end{proof}

The above result is quite interesting, and based on it, let us formulate:

\index{Bessel function}
\index{Bessel law}

\begin{definition}
The Bessel law of parameter $t>0$ is the measure
$$p^2_t=e^{-t}\sum_{k=-\infty}^\infty\delta_k\,f_k(t/2)\quad,\quad 
f_k(t)=\sum_{p=0}^\infty \frac{t^{|k|+2p}}{(|k|+p)!p!}$$
with $f_k$ being the Bessel function of the first kind, from classical analysis.
\end{definition}

Let us study now these Bessel laws that we found. In analogy with what we know about the Poisson laws $p_t$, from chapter 1, we first have the following result:

\index{convolution}
\index{convolution semigroup}

\begin{theorem}
The Bessel laws $p^2_t$ have the following properties:
\begin{enumerate}
\item The Fourier transform is $F(x)=\exp((\cos x-1)t)$.

\item We have the semigroup formula $p^2_s*p^2_t=p^2_{s+t}$.
\end{enumerate}
\end{theorem}

\begin{proof}
We use the formula in Definition 4.19. The Fourier transform is:
$$F(x)=e^{-t}\sum_{k=-\infty}^\infty e^{ikx}\,f_k(t/2)$$

We can compute the derivative of $F$ with respect to $t$, as follows:
\begin{eqnarray*}
F(x)'
&=&-e^{-t}\sum_{k=-\infty}^\infty e^{ikx}\,f_k(t/2)+\frac{e^{-t}}{2}\sum_{k=-\infty}^\infty e^{ikx}\,f_k'(t/2)\\
&=&-F(x)+\frac{e^{-t}}{2}\sum_{k=-\infty}^\infty e^{ikx}\,f_k'(t/2)
\end{eqnarray*}

On the other hand, the derivative of $f_k$ with $k\geq 1$ is given by:
\begin{eqnarray*}
f_k'(t)
&=&\sum_{p=0}^\infty \frac{(k+2p)t^{k+2p-1}}{(k+p)!p!}\\
&=&\sum_{p=0}^\infty \frac{(k+p)t^{k+2p-1}}{(k+p)!p!}+\sum_{p=0}^\infty\frac{p\,t^{k+2p-1}}{(k+p)!p!}\\
&=&\sum_{p=0}^\infty \frac{t^{k+2p-1}}{(k+p-1)!p!}+\sum_{p=1}^\infty\frac{t^{k+2p-1}}{(k+p)!(p-1)!}\\
&=&f_{k-1}(t)+f_{k+1}(t)
\end{eqnarray*}

This computation works in fact for any $k\in\mathbb Z$, so we get:
\begin{eqnarray*}
F(x)'
&=&-F(x)+\frac{e^{-t}}{2}
\sum_{k=-\infty}^\infty e^{ikx} (f_{k-1}(t/2)+f_{k+1}(t/2))\\
&=&-F(x)+\frac{e^{-t}}{2} \sum_{k=-\infty}^\infty
e^{i(k+1)x}f_{k}(t/2)+e^{i(k-1)x}f_{k}(t/2)\\
&=&-F(x)+\frac{e^{ix}+e^{-ix}}{2}\,F(x)\\
&=&(\cos x-1)F(x)
\end{eqnarray*}

Thus, by integrating, we are led to the conclusions in the statement.
\end{proof}

In order to unify now the theory of the Poisson and Bessel laws, we have the following key notion, extending the Poisson limit theory from chapter 1:

\index{compound Poisson law}

\begin{definition}
Associated to any compactly supported positive measure $\nu$ on $\mathbb R$ is the probability measure
$$p_\nu=\lim_{n\to\infty}\left(\left(1-\frac{c}{n}\right)\delta_0+\frac{1}{n}\,\nu\right)^{*n}$$
where $c=mass(\nu)$, called compound Poisson law. As example, $p_t=p_{t\delta_1}$.
\end{definition}

In other words, what we are doing here is to generalize the construction in the Poisson Limit Theorem, by allowing the only parameter there, which was the positive real number $t>0$, to be replaced by a certain probability measure $\nu$, of arbitrary mass $c>0$. We will be actually mostly interested in the case where $\nu$ is discrete, and often assume so.

\bigskip

Getting now to some theory, for the above laws, we first have the following result:

\index{Fourier transform}

\begin{proposition}
For a discrete measure, $\nu=\sum_{i=1}^sc_i\delta_{z_i}$ with $c_i>0$ and $z_i\in\mathbb R$, we have the formula
$$F_{p_\nu}(x)=\exp\left(\sum_{i=1}^sc_i(e^{ixz_i}-1)\right)$$
where $F$ denotes as usual the Fourier transform.
\end{proposition}

\begin{proof}
In order to prove the formula in the statement, consider the measure $\mu_n$ appearing in Definition 4.21, under the convolution sign, namely:
$$\mu_n=\left(1-\frac{c}{n}\right)\delta_0+\frac{1}{n}\,\nu$$

We have the following computation, in the context of Definition 4.21:
\begin{eqnarray*}
F_{\mu_n}(x)=\left(1-\frac{c}{n}\right)+\frac{1}{n}\sum_{i=1}^sc_ie^{ixz_i}
&\implies&F_{\mu_n^{*n}}(x)=\left(\left(1-\frac{c}{n}\right)+\frac{1}{n}\sum_{i=1}^sc_ie^{ixz_i}\right)^n\\
&\implies&F_{p_\nu}(x)=\exp\left(\sum_{i=1}^sc_i(e^{ixz_i}-1)\right)
\end{eqnarray*}

Thus, we have obtained the formula in the statement.
\end{proof}

Next, we have the Compound Poisson Limit Theorem, as follows:

\index{compound Poisson Limit theorem}
\index{CPLT}
\index{compound PLT}
\index{compound Poisson limit}

\begin{theorem}[CPLT]
For a discrete measure, $\nu=\sum_{i=1}^sc_i\delta_{z_i}$ with $c_i>0$ and $z_i\in\mathbb R$, we have the following formula,
$$p_\nu=law\left(\sum_{i=1}^sz_i\alpha_i\right)$$
with the variables $\alpha_i$ being Poisson $(c_i)$, and independent.
\end{theorem}

\begin{proof}
Let $\alpha$ be the sum of Poisson variables in the statement. We have:
\begin{eqnarray*}
F_{\alpha_i}(x)=\exp(c_i(e^{ix}-1))
&\implies&F_{z_i\alpha_i}(x)=\exp(c_i(e^{ixz_i}-1))\\
&\implies&F_\alpha(x)=\exp\left(\sum_{i=1}^sc_i(e^{ixz_i}-1)\right)
\end{eqnarray*}

Thus we have the same formula as in Proposition 4.22, as desired.
\end{proof}

Getting back now to the Bessel laws, we have the following result:

\index{Bessel law}

\begin{theorem}
The Bessel laws $p^2_t$ are compound Poisson laws, given by
$$p^2_t=p_{t\varepsilon}$$
where $\varepsilon=\frac{1}{2}(\delta_{-1}+\delta_1)$ is the uniform measure on $\pm1$.
\end{theorem}

\begin{proof}
We recall from Theorem 4.20 that the Fourier transform of the Bessel law $p^2_t$, with respect to a variable $x$, is given by the following formula:
$$F(x)=\exp((\cos x-1)t)=\exp\left(\left(\frac{e^{ix}+e^{-ix}}{2}-1\right)t\right)$$

But, according to Proposition 4.22, this is exactly the Fourier transform of the compound Poisson law $p_{t\varepsilon}$ in the statement. Thus, we have $p^2_t=p_{t\varepsilon}$, as claimed.
\end{proof}

The above result is quite interesting, and together with $p_t=p_{t\delta_1}$ suggests looking, more generally, at the laws $p^s_t=p_{t\varepsilon_s}$, with $\varepsilon_s$ being the uniform measure on the $s$-th roots of unity. We will do this later, in Part II, when discussing the complex measures.

\bigskip

Getting now to the moments of $p^2_t$, it is possible to recapture them from $F(x)=\exp((\cos x-1)t)$, via some work. But as an alternative to this, and standing as an excellent illustration for Principle 4.2, we can do this via group theory power:

\begin{theorem}
The hyperoctadedral group $H_N\subset O_N$ is easy, according to
$$Hom(u^{\otimes k},u^{\otimes l})=span\left(T_\pi\Big|\pi\in P_{even}(k,l)\right)$$
with $P_{even}$ being the partitions all whose blocks have even size. Thus, we have
$$M_k(p^2_t)=\sum_{\pi\in P_{even}(k)}t^{|\pi|}$$
and in particular at $t=1$ we have $M_k(p^2_1)=|P_{even}(k)|$.
\end{theorem}

\begin{proof}
We have two possible proofs for this, which are both elementary:

\medskip

(1) We can either follow the proof of Theorem 4.12, for the symmetric group $S_N$, by adding some signs there, and since at the very end the product of these signs must be 1, the blocks of the corresponding partition must be of even size, as claimed.

\medskip

(2) Or we can follow the proof of Theorem 4.16, for the orthogonal group $O_N$, by replacing the semicircle $\cap\in P(0,2)$ used there by the partition $H\in P(2,2)$, the point being that $T_H\in End(u^{\otimes 2})$ corresponds to the relations defining $H_N\subset O_N$.

\medskip

(3) In short, easiness property proved, one way or another, and with this in hand, the moment formula in the statement is clear, exactly as for $S_N,O_N$ before.
\end{proof}

Let us end this discussion with a grand result about $p^2_t$, in the spirit of what we did before in this book, for other probability measures that we met, as follows:

\begin{theorem}
The Bessel law $p^2_t$ has the following properties:
\begin{enumerate}
\item It is a compound Poisson law, $p^2_t=p_{t\varepsilon}$, with $\varepsilon=(\delta_{-1}+\delta_1)/2$.

\item $p^2_t=e^{-t}\sum_{k=-\infty}^\infty\delta_k\,f_k(t/2)$, with 
$f_k(t)=\sum_{p=0}^\infty \frac{t^{|k|+2p}}{(|k|+p)!p!}$, Bessel.

\item The moments are $M_{2k}=M_{2k}'=\sum_{\pi\in P_{even}(2k)}t^{|\pi|}$.

\item The normalized central moments are $M_{2k}''=\sum_{\pi\in P_{even}(2k)}t^{|\pi|-1}$.

\item We have $E=0$, $V=t$, $\gamma=0$, $\kappa=3+1/t$.

\item The Fourier transform is $F(x)=\exp((\cos x-1)t)$.

\item We have the semigroup property $p^2_s*p^2_t=p^2_{s+t}$.
\end{enumerate}
\end{theorem}

\begin{proof}
We already know all this, save for the formula of the kurtosis, which comes from $P_{even}(4)=\{\cap\,\cap,\cap\!\!\cap,\Cap,\sqcap\hskip-1.6mm\sqcap\hskip-1.6mm\sqcap\}$, which gives $M_4=3t^2+t$, and so $\kappa=3+1/t$.
\end{proof}

Regarding now more advanced questions, in relation with Hankel determinants and orthogonal polynomials, things here are quite tricky. The even moments at $t=1$ are:
$$1\ ,\ 1\ , \ 4\ ,\ 31\ ,\ 379\ ,\ 6556\ ,\ 150349\ ,\ \ldots$$

The first few computations of Hankel determinants are encouraging, as follows:
$$H_1=1\ ,\quad\frac{H_2}{H_1}=1\ ,\quad\frac{H_3}{H_2}=3\ ,\quad\frac{H_4}{H_3}=15\ ,\quad\frac{H_5}{H_4}=120\ ,\quad\frac{H_6}{H_5}=1260$$

However, a bad surprise awaits us at $k=7$, and then at $k=8$ too, as follows:
$$\frac{H_7}{H_6}=2\cdot 9\cdot5\cdot7\cdot29\quad,\quad 
\frac{H_8}{H_7}=4\cdot 27\cdot5\cdot7\cdot83$$

So, where do these numbers, 29 and then 83, come from? Good question, I guess. As for the orthogonal polynomials, nor do I know them. Interesting questions, here.

\section*{4d. Poisson and Gauss}

What is next? More easiness I guess, I mean looking for more easy groups, $G=(G_N)$ with  $G_N\subset O_N$, then computing their character laws in the $N\to\infty$ limit, and then further studying these laws, via various methods, as we just did for $G=H$.

\bigskip

In order to discuss this, let us start with the following key definition, from \cite{bsp}:

\index{category of partitions}

\begin{definition}
A category of partitions is a collection $D=(D(k,l))$ of sets of partitions $D(k,l)\subset P(k,l)$, satisfying the following conditions:
\begin{enumerate}
\item Stability under the horizontal concatenation, $(\pi,\sigma)\to[\pi\sigma]$.

\item Stability under vertical concatenation $(\pi,\sigma)\to[^\sigma_\pi]$.

\item Stability under upside-down turning.

\item $P(k,k)$ contains the identity partition $||\ldots||$.

\item $P(0,2)$ contains the semicircle $\cap$.

\item $P(2,2)$ contains the crossing partition $\slash\hskip-2.0mm\backslash$.
\end{enumerate}
\end{definition}

As examples, we have the categories $P,P_2,P_{even}$, that we met before in relation with the groups $S_N,O_N,H_N$. In fact, any easy group must come from such a category:

\index{easy group}

\begin{theorem}
A subgroup $G\subset_uO_N$ is easy precisely when
$$Hom(u^{\otimes k},u^{\otimes l})=span\left(T_\pi\Big|\pi\in D(k,l)\right)$$
for any $k,l\in\mathbb N$, for a certain category of partitions $D\subset P$.
\end{theorem}  

\begin{proof}
This is something categorical, based on Tannakian duality, the idea being that given an easy group $G\subset_uO_N$, coming from a family of subsets $D(k,l)\subset P(k,l)$, via the formulae in the statement, we can always enlarge the collection $D=(D(k,l))$ into a full category, by using the various operations in Definition 4.27, and in this process, the associated subgroup of $O_N$ will not shrink, I mean will remain $G$ itself.
\end{proof}

The above result is quite interesting, among others because it leads to:

\begin{fact}
We can talk about easiness by staying away from groups, if we want to, by having Definition 4.27 as axiomatics, and then the formulae
$$M_k(\mu_1)=|D(k)|\quad,\quad M_k(\mu_t)=\sum_{\pi\in D(k)}t^{|\pi|}$$
with some care for the second one, as definitions for $\mu_1$, and then for $\mu_t$ with $t>0$.
\end{fact}

To be more precise here, the claim regarding $\mu_1$ is clear, because we have an easy group coming from our category $D\subset P$, and the asymptotic moments of the main character are indeed the numbers $M_k=|D(k)|$. As for the claim regarding $\mu_t$ with $t>0$, this is something more tricky, involving the Weingarten formula, and technically, we must assume here that $D$ is stable under removing blocks. More on this in a moment.

\bigskip

In practice now, shall we say goodbye to groups, and use Fact 4.29 as stated? Or even worse, say goodbye to both groups and probability, and study the categories of partitions only? Not clear at all, but with the orange cat passing by, I know whom to ask:

\begin{tiger}
Power, speed and wisdom come from quantum, which is a mixture of group theory and probability. Stay away from abstractions, these are no good.
\end{tiger}

Okay, thanks cat, so we will keep Fact 4.29 somewhere in our minds, that is not that bad, after all, but resume with some ad-hoc hunting for easiness, using group theory methods. So, looking for a 4th easy group $S_N\subset G_N,H_N\subset O_N$, here is our puzzle:
$$\xymatrix@R=50pt@C=50pt{
G_N\ar[r]&O_N\\
S_N\ar[u]\ar[r]&H_N\ar[u]}$$

But in view of this diagram, it makes sense to look for $G_N$ satisfying:
$$G_N\cap H_N=S_N\quad,\quad <G_N,H_N>=O_N$$

And with this, done, because the first condition suggests taking $G_N=B_N$, the orthogonal bistochastic group, I mean what produces $S_N\subset H_N$ is the fact that rows and columns sum up to 1. As for the second condition, a bit of thinking suggests taking $G_N=B_N$ too. Thus, puzzle solved, and in practice, this leads to the following result:

\index{bistochastic group}
\index{singletons and pairings}

\begin{theorem}
We have easy groups as follows,
$$\xymatrix@R=50pt@C=50pt{
B_N\ar[r]&O_N\\
S_N\ar[u]\ar[r]&H_N\ar[u]}\qquad\xymatrix@R=25pt@C=10pt{\\ :\\ }
\qquad\xymatrix@R=50pt@C50pt{
P_{12}\ar[d]&P_2\ar[d]\ar[l]\\
P&P_{even}\ar[l]}$$
with $B_N\subset O_N$ being the bistochastic group, and $P_{12}$, the singletons and pairings.
\end{theorem}

\begin{proof}
Many things can be said here, the idea being as follows:

\medskip

(1) To start with, $B_N\subset O_N$ is by definition the group of matrices whose rows and columns sum up to 1. Now observe that the bistochasticity condition reads $U\xi=\xi$, with $\xi$ being the all-one vector, and so $\xi\in Fix(u)$. But $\xi=T_|$ being the vector associated to the singleton \!\!${\ }_|\in P(0,1)$, we conclude that $B_N$ is indeed easy, coming from $P_{12}$.

\medskip

(2) Next, in relation with the considerations made before the statement, observe that both the diagrams in the statement $A\subset B,C\subset D$ are intersection and generation diagrams, in the sense that $B\cap C=A$, $<B,C>=D$. Which is nice, and useful.

\medskip

(3) Finally, $\xi\in Fix(u)$ tells us to look at $u-1$, and some routine study here, that we will leave as an exercise, shows that we have in fact an isomorphism $B_N\simeq O_{N-1}$. However, and we insist, from our present easy viewpoint, $B_N$ is a new solution.
\end{proof}

Before going ahead with probabilistic aspects for $B_N$, let us try to have our hunting work finished. And here, in order to cut from complexity, it is convenient to impose an extra axiom, called uniformity, coming from the following result:

\index{uniform group}
\index{removing blocks}

\begin{proposition}
For an easy group $G=(G_N)$, coming from a category of partitions $D\subset P$, the following conditions are equivalent:
\begin{enumerate}
\item $G_{N-1}=G_N\cap O_{N-1}$, via the embedding $O_{N-1}\subset O_N$ given by $u\to diag(u,1)$.

\item $G_{N-1}=G_N\cap O_{N-1}$, via the $N$ possible diagonal embeddings $O_{N-1}\subset O_N$.

\item $D$ is stable under the operation which consists in removing blocks.
\end{enumerate}
If these conditions are satisfied, we say that $G=(G_N)$ is uniform.
\end{proposition}

\begin{proof}
This is something very standard, that we will leave as an exercise. Let us also mention that our condition is something very natural in the probabilistic context, in relation with Weingarten integration, as mentioned in the comments after Fact 4.29.
\end{proof}

Following \cite{bsp}, we can now formulate a nice classification result, as follows:

\index{sizes of blocks}

\begin{theorem}
The uniform easy groups are those that we know, namely
$$\xymatrix@R=50pt@C=50pt{
B_N\ar[r]&O_N\\
S_N\ar[u]\ar[r]&H_N\ar[u]}\qquad\xymatrix@R=25pt@C=10pt{\\ :\\ }
\qquad\xymatrix@R=50pt@C50pt{
P_{12}\ar[d]&P_2\ar[d]\ar[l]\\
P&P_{even}\ar[l]}$$
and with this ultimately coming from only $4$ possible choices for the sizes of blocks.
\end{theorem}

\begin{proof}
Consider an arbitrary easy group $S_N\subset G_N\subset O_N$. This group must then come from a certain category of partitions, as follows:
$$P_2\subset D\subset P$$

Now if we assume $G=(G_N)$ to be uniform, this category of partitions $D$ is uniquely determined by the subset $L\subset\mathbb N$ consisting of the sizes of the blocks of the partitions in $D$. And, our claim is that there are only 4 admissible sets, as follows:

\smallskip

\begin{enumerate}
\item $L=\{2\}$, producing $O_N$.

\medskip

\item $L=\{1,2\}$, producing $B_N$.

\medskip

\item $L=\{2,4,6,\ldots\}$, producing $H_N$.

\medskip

\item $L=\{1,2,3,\ldots\}$, producing $S_N$.
\end{enumerate}

\smallskip

Indeed, assume that $L\subset\mathbb N$ is such that the set $P_L$ consisting of partitions whose sizes of blocks belong to $L$ is a category of partitions. We know from the axioms of the categories of partitions that the semicircle $\cap$ must be in the category, so $2\in L$. With a bit more work, we can see that the following conditions must be satisfied as well:
$$k,l\in L,\,k>l\implies k-l\in L$$
$$k\in L,\,k\geq 2\implies 2k-2\in L$$

But with these two formulae in hand, it is elementary to conclude. Indeed, as explained in \cite{bsp}, the situation $1\in L$ leads to the solutions $L=\{1,2\}$ and $L=\{1,2,3,\ldots\}$, and the situation $1\notin L$ leads to the solutions $L=\{2\}$ and $L=\{2,4,6,\ldots\}$.
\end{proof}

Getting now to probabilistic aspects, the isomorphism $B_N\simeq O_{N-1}$ suggests that we should obtain certain modifications of the normal laws $g_t$. And this is indeed the case:

\index{shifted normal law}

\begin{theorem}
For the group $B_N$ we have $\chi_t\sim g_t^t$ with $N\to\infty$, where
$$g_t^t=\frac{1}{\sqrt{2\pi t}}\,e^{-(x-t)^2/2t}\,dx$$
and the basic properties of this shifted normal law $g_t^t$ are as follows:
\begin{enumerate}
\item The moments are $M_k=\sum_{\pi\in P_{12}(k)}t^{|\pi|}$.

\item The central moments are $M_k'=\delta_{2|k}t^{k/2}k!!$.

\item The normalized central moments are $M_k''=\delta_{2|k}k!!$.

\item We have $E=t$, $V=t$, $\gamma=0$, $\kappa=3$.

\item The Fourier transform is $F(x)=\exp(itx-tx^2/2)$.

\item We have the convolution formula $g_s^s*g_t^t=g_{s+t}^{s+t}$.

\item The Hankel determinant formula at $t=1$ is $H_k=\prod_{i=1}^{k-1}i!$.

\item The orthogonal polynomials are shifted Hermite polynomials.
\end{enumerate}
\end{theorem}

\begin{proof}
By easiness we have the following formula, with on the right $r$ standing for the number of pairs, and so $r+(k-2r)=k-r$ standing for the number of blocks:
$$M_k(\chi_t)\simeq\sum_{\pi\in P_{12}(k)}t^{|\pi|}=\sum_{r=0}^{[k/2]}\binom{k}{2r}(2r)!!t^{k-r}$$

On the other hand, the moments of the shifted normal law $g_t^t$ are given by:
\begin{eqnarray*}
M_k(g_t^t)
&=&\frac{1}{\sqrt{2\pi t}}\int_\mathbb R(x+t)^k\,e^{-x^2/2t}\,dx\\
&=&\sum_{r=0}^{[k/2]}\binom{k}{2r}t^{k-2r}\frac{1}{\sqrt{2\pi t}}\int_\mathbb Rx^{2r}\,e^{-x^2/2t}\,dx\\
&=&\sum_{r=0}^{[k/2]}\binom{k}{2r}(2r)!!t^{k-r}
\end{eqnarray*}

Thus we have the first assertion, and the rest are things that we already know.
\end{proof}

And with this, end of our easiness study, job done, time to light a cigar, and formulate a few philosophical comments. To start with, Theorem 4.33 reformulates as:

\begin{fact}
In the real, classical, uniform easy setting, the laws are
$$\xymatrix@R=45pt@C=50pt{
g_t^t\ar@{-}[r]&g_t\\
p_t\ar@{-}[u]\ar@{-}[r]&p^2_t\ar@{-}[u]}$$
and so, at the ground level in easiness, Poisson and Gauss rule.
\end{fact}

Which sounds quite good, when compared to the rival theory, which is the Askey scheme \cite{awi}. This being said, we are cheating here a bit, because for Poisson and Gauss to truly rule, we would need an extreme familiarity with $p^2_t$, which in turn brings us to the Askey scheme. Quite fun all this, and we will therefore declare a tie.

\bigskip

And more on this, later in this book. We will be back to such things twice, first in Part II, with complex laws $p^s_t,G_t,G_t^t$, and then later in Part IV, with free versions $\pi^s_t,\gamma_t,\gamma_t^t,\Gamma_t,\Gamma_t^t$, and all sorts of intermediate beasts that can be constructed too.

\section*{4e. Exercises}

This was a quite exciting chapter, and as exercises on this, we have:

\begin{exercise}
Do the computations, for the exact law formula for $S_N$.
\end{exercise}

\begin{exercise}
Learn more about the representations of compact groups.
\end{exercise}

\begin{exercise}
Learn more about the uniform measure on compact groups.
\end{exercise}

\begin{exercise}
Learn as well about Tannakian duality for compact groups.
\end{exercise}

\begin{exercise}
Learn more about the hyperoctahedral group $H_N$.
\end{exercise}

\begin{exercise}
Learn as well about $A_N$, and other reflection groups.
\end{exercise}

\begin{exercise}
Do some further computations for the Bessel laws $p^2_t$.
\end{exercise}

\begin{exercise}
Do some further computations for the shifted normal laws $g_t^t$.
'\end{exercise}

As bonus exercise, and no surprise here, read some group theory.

\part{Complex variables}

\ \vskip50mm

\begin{center}
{\em E la pioggia che va, e ritorna il sereno

E col tempo sopra il mondo

Come il sole del mattino

Un amore universale sorgera}
\end{center}

\chapter{Complex variables}

\section*{5a. Complex variables}

Time to upgrade our theory of real random variables $f:X\to\mathbb R$, into a theory of complex random variables $f:X\to\mathbb C$. We will be doing this in the present Part II, with as usual a main focus on the notion of normality, our goal being that of finding and studying the complex analogues of the Poisson law $p_t$, and of the Gauss law $\gamma_t$.

\bigskip

Before anything, however, where do these complex variables $f:X\to\mathbb C$ come from? I mean, certainly not from usual measuring devices, such as thermometers and manometers, I have yet to see such an engineering machine displaying complex numbers. Although, I am pretty much sure that such things have already been invented, and are already for sale on the internet, and making good money from the math nerds buying them.

\bigskip

In answer, the complex variables $f:X\to\mathbb C$ are something quite natural, mathematically speaking, and can appear in several ways. To start with, recall that the complex numbers originally appear as roots of polynomials $P\in\mathbb R[X]$, and think here $i^2=-1$. But this leads us to the following scenario, which looks quite plausible:

\begin{scenario}
Imagine that your problem involves certain polynomials associated to your events, $x\to P_x$, and more specifically, that you are chasing one of the roots $r_i\in\mathbb C$ of these polynomials. Then, what you have is a complex variable, $f(x)=r_i(P_x)$.
\end{scenario}

As a comment here, you can be of course interested, more generally, in a certain function of the roots $\varphi(P)\in\mathbb C$, and in this case your complex variable is $f(z)=\varphi(P_x)$. However, when your function $\varphi$ is symmetric, this function must be a combination of the coefficients of the polynomial, and so must be real, $\varphi(P)\in\mathbb R$. So, we must assume that there is some lack of symmetry in $\varphi$, and with this in mind, the simplest example remains the one above, $\varphi$ being one of the roots. Say the one having biggest absolute value, $|r_i|>|r_j|$ for $j\neq i$, assumed to be unique, in the context of your problem.

\bigskip

Along the same lines, we can have situations involving matrix eigenvalues:

\begin{scenario}
Imagine that your problem involves certain matrices associated to your events, $x\to A_x$, and more specifically, that you are chasing one of the eigenvalues $\lambda_i\in\mathbb C$ of these matrices. Then, what you have is a complex variable, $f(x)=\lambda_i(A_x)$.
\end{scenario}

In this context, the same comments as before apply, with a symmetric function $\varphi(A)$ of the eigenvalues not doing it, due to $\varphi(A)\in\mathbb R$, and with a quite plausible situation, when chasing one particular eigenvalue, as stated, being that when you are interested in the biggest eigenvalue, $|\lambda_i|>|\lambda_j|$ for $j\neq i$, assumed to be unique, in the context of your problem. In fact, thinking a bit, since the eigenvalues appear as roots of the characteristic polynomial, Scenario 5.2 is in fact a particular case of Scenario 5.1.

\bigskip

Moving on, complex variables naturally appear from complex geometry:

\begin{scenario}
Imagine that you are doing geometry, being interested in this compact Lie group $X=G\subset U_N$, or that quotient space $X=G/H\subset\mathbb C^M$, and so on. In this case the coordinate functions $z_i:X\to\mathbb C$ are your complex variables.
\end{scenario}

As before, there are several comments that can be made here. In the group situation, quite remarkably, any compact Lie group, and with this being a real notion, not involving $\mathbb C$, can be shown to appear as a closed subgroup of a unitary group, $G\subset U_N$. In a word, continuous symmetries naturally produce complex numbers, and complex variables. As for the quotient space situation, many things can be said here too. Finally, you can be interested, more generally, in a manifold $X\subset\mathbb C^M$, and here again you have complex variables $z_i:X\to\mathbb C$, provided that you know how to integrate over $X$.

\bigskip

Getting back now to matrices, as a variation of Scenario 5.2, we have:

\begin{scenario}
Imagine that you are doing quantum, and are interested in a certain matrix $A\in M_n(\mathbb C)$. In this case, with the convention $C(M_n)=M_n(\mathbb C)$, your matrix is a function $A:M_n\to\mathbb C$, and in the normal case, $AA^*=A^*A$, its law is $\mu=\sum_i\delta_{\lambda_i}$.
\end{scenario}

Obviously, this is something quite speculative, because that beast $M_n$ is not a probability space in the usual sense. However, the algebra of random variables is there, according to $C(M_n)=M_n(\mathbb C)$, and the expectation is there too, positive and of mass 1, according to $E(A)=tr(A)$, with $tr=Tr/n$ being the normalized trace of matrices. Moreover, in the case where the matrix is normal, $AA^*=A^*A$, we are in fact into classical probability, due to the fact that the algebra $<A>\subset M_n(\mathbb C)$ is commutative, consisting of functions on the spectrum $X=\{\lambda_1,\ldots,\lambda_n\}$, and $\mu=\sum_i\delta_{\lambda_i}$. And, more on such things later.

\bigskip

As a variation of Scenario 5.4, which is of particular interest, we have:

\begin{scenario}
Imagine that you are doing advanced probability, and are interested in a certain random matrix $Z\in M_n(L^\infty(X))$. In this case, you can talk about the law of $Z$, and in the  normal case, $ZZ^*=Z^*Z$, this law is a probability measure on $\mathbb C$.
\end{scenario}

As before, many things can be said here, the idea being that the underlying probability type space is $M_n\times X$ in this case, with the algebra of random variables being there, according to $L^\infty(M_n\times X)=M_n(L^\infty(X))$, and with the expectation being there too, given by $Z\to E(tr(Z))$, with $tr=Tr/n$ being the normalized trace of matrices. And in the case where the matrix is normal, $ZZ^*=Z^*Z$, we are in fact into classical probability, due to the fact that the algebra $<Z>\subset M_n(L^\infty(X))$ is commutative, consisting of functions on the spectrum $\sigma(Z)\subset\mathbb C$. And, more on all this later.

\bigskip

Finally, here is one more scenario where complex variables naturally appear:

\begin{scenario}
Imagine that you are doing usual, real probability, but the real measure $\mu$ that you found is something quite complicated. In this case, one trick to be tried, of Fourier type, would be to blow up $\mu$ into a measure $\varepsilon$ on the unit circle $\mathbb T\subset\mathbb C$.
\end{scenario}

Which sounds quite plausible too, I mean complex numbers coming to the rescue when stuck with complicated formulae involving reals, examples of this phenomenon abound, and we will discuss concrete illustrations of this in probability, as above, a bit later. 

\bigskip

As a conclusion now, with no less than 6 plausible scenarii, we definitely need to have a look at complex probability, which looks like very useful technology. Let us start with some generalities. In analogy with our previous axioms in the real case, we have:

\begin{definition}
The axioms for complex probability theory are as follows:
\begin{enumerate}
\item A probability space is a measured space $(X,M,\nu)$ of mass one, $\nu(X)=1$.

\item A random variable on $X$ is a complex measurable function $f:X\to\mathbb C$.

\item The expectation of such a variable is $E(f)=\int_Xf(x)d\nu(x)\in\mathbb C$.

\item The variance of such a variable is $V(f)=E(|f-E(f)|^2)\geq0$.

\item The law of $f$ is the complex probability measure $\mu=f_*\nu$, push-forward of $\nu$ by $f$.
\end{enumerate}
\end{definition}

Obviously, all this is quite straightforward, generalizing what we previously had in the real case. In relation with the variance, we have the following alternative formula:
\begin{eqnarray*}
V(f)
&=&E\left[(f-E(f))(\overline{f-E(f)})\right]\\
&=&E(f\bar{f})-E(f)\overline{E(f)}\\
&=&E(|f|^2)-|E(f)|^2
\end{eqnarray*}

In relation with the law, observe that we have the following formula, valid for any  $k,l\in\mathbb N$, allowing us in particular, with $k=l=1$, to compute the variance:
$$E(f^k\bar{f}^l)=\int_\mathbb Cz^k\bar{z}^ld\mu(z)$$

In fact, the above quantities $E(f^k\bar{f}^l)$, which uniquely determine the law, can be thought of as being the moments of $f$. And with this showing that $V(f)$ does not really encapsulate the information contained in the second moment, because in the present complex case we have 3 such second moments, namely $E(f\bar{f})$ used by $V(f)$, but also $E(f^2)$ and $E(\bar{f}^2)$, which are conjugate, and whose values are not captured by $V(f)$.

\bigskip

As an illustration for all this, in the discrete case the situation is as follows:

\index{complex variable}
\index{random variable}

\begin{proposition}
For a discrete variable $f:X\to\mathbb C$, the law is the measure
$$\mu=\sum_i c_i\delta_{z_i}\quad,\quad c_i>0\quad,\quad\sum_ic_i=1\quad,\quad z_i\in\mathbb C$$
given by the following formula, with $P$ being the probability over $X$,
$$\mu=\sum_{z\in\mathbb C}P(f=z)\delta_z$$
with the sum being finite or countable, as per our discretness assumption. We have
$$E(f)=\sum_ic_iz_i\quad,\quad V(f)=\sum_ic_i|z_i|^2-\Big|\sum_ic_iz_i\Big|^2$$
and the second moments are $E(f^2)=\sum_ic_iz_i^2$, $E(f\bar{f})=\sum_ic_i|z_i|^2$, $E(\bar{f}^2)=\sum_ic_i\bar{z}_i^2$.
\end{proposition}

\begin{proof}
The first assertion is self-explanatory, then we have the formulae of $E(f)$ and of the three order 2 moments, and then we have the formula of $V(f)$.
\end{proof}

Getting now to independence, things are again similar to those from the real case:

\begin{definition}
We say that two variables $f,g:X\to\mathbb C$ are independent when
$$P(f\in I,g\in J)=P(f\in I)P(g\in J)$$
happens, for any two measurable sets $I,J\subset\mathbb C$. Equivalently, we must have
$$E(f^k\bar{f}^lg^m\bar{g}^n)=E(f^k\bar{f}^l)E(g^m\bar{g}^n)$$
for any $k,l,m,n\in\mathbb N$, that is, the mixed moments must factorize.
\end{definition}

To be more precise here, the equivalence between the above two conditions is standard, exactly as in the real case. As an illustration for this, in the discrete case we have the following computation, which works in both senses, and gives the equivalence:
\begin{eqnarray*}
E(f^k\bar{f}^lg^m\bar{g}^n)
&=&\sum_{y,z\in\mathbb C}P(f=y,g=z)y^k\bar{y}^lz^m\bar{z}^n\\
&=&\sum_{y,z\in\mathbb C}P(f=y)P(g=z)y^k\bar{y}^lz^m\bar{z}^n\\
&=&\sum_{y\in\mathbb C}P(f=y)y^k\bar{y}^l\sum_{z\in\mathbb C}P(g=z)z^m\bar{z}^n\\
&=&E(f^k\bar{f}^l)E(g^m\bar{g}^n)
\end{eqnarray*}

Still in relation with independence, let us record as well the following fact:

\begin{theorem}
Assuming that $f,g:X\to\mathbb C$ are independent, we have
$$\mu_{f+g}=\mu_f*\mu_g$$
where $*$ is the convolution operation for the complex probability measures.
\end{theorem}

\begin{proof}
Again, this is something standard, exactly as in the real case. As an illustration, in the discrete case the formula in the statement follows from:
\begin{eqnarray*}
\mu_{f+g}
&=&\sum_{z\in\mathbb C}P(f+g=z)\delta_z\\
&=&\sum_{x,y\in\mathbb C}P(f=x,g=y)\delta_{x+y}\\
&=&\sum_{x,y\in\mathbb C}P(f=x)P(g=y)\delta_x*\delta_y\\
&=&\left(\sum_{x\in\mathbb C}P(f=x)\delta_x\right)*\left(\sum_{y\in\mathbb C}P(f=y)\delta_y\right)\\
&=&\mu_f*\mu_g
\end{eqnarray*}

Thus, one way or another, we are led to the formula in the statement.
\end{proof}

In order to understand now the relation with the real theory, as developed in Part I, we can decompose any complex variable $f:X\to\mathbb C$ as a sum, as follows:
$$f=g+ih\quad,\quad g=Re(f),\ h=Im(f)$$

In practice, this method leads, among others, to the following useful result:

\index{real part}
\index{imaginary part}
\index{rotation of variable}

\begin{theorem}
For a complex random variable $f:X\to\mathbb C$, decomposed into real and imaginary parts as $f=g+ih$, and with $g,h$ assumed independent, we have
$$\mu_f=\mu_g*i\mu_h$$
with $*$ being the usual convolution operation, $\delta_z*\delta_t=\delta_{z+t}$, and with $\mu\to i\mu$ denoting the rotated version, $\mathbb R\to i\mathbb R$. If $g,h$ are not independent, this formula does not hold.
\end{theorem}

\begin{proof}
This is something very standard, the idea being as follows:

\medskip

(1) As usual, I will do the computations in the discrete case, and leave the discussion of the general continuous case to you, either by using a density argument, or by redoing the computation, with sums replaced by integrals. Assuming $f=g+ih$, we have:
\begin{eqnarray*}
\mu_f
&=&\sum_{z\in\mathbb C}P(f=z)\delta_z\\
&=&\sum_{x,y\in\mathbb R}P(f=x+iy)\delta_{x+iy}\\
&=&\sum_{x,y\in\mathbb R}P(g+ih=x+iy)\delta_{x+iy}\\
&=&\sum_{x,y\in\mathbb R}P(g=x,h=y)\delta_{x+iy}
\end{eqnarray*}

(2) In the case where the real and imaginary parts $g,h:X\to\mathbb R$ are independent, we can say more, with the above computation having the following continuation:
\begin{eqnarray*}
\mu_f
&=&\sum_{x,y\in\mathbb R}P(g=x,h=y)\delta_{x+iy}\\
&=&\sum_{x,y\in\mathbb R}P(g=x)P(h=y)\delta_{x+iy}\\
&=&\sum_{x,y\in\mathbb R}P(g=x)P(h=y)\delta_x*\delta_{iy}\\
&=&\left(\sum_{x\in\mathbb R}P(g=x)\delta_x\right)*\left(\sum_{y\in\mathbb R}P(h=y)\delta_{iy}\right)\\
&=&\mu_g*i\mu_h
\end{eqnarray*}

To be more precise, we have used here in the beginning the independence of the variables $h,g:X\to\mathbb R$, and at the end we have denoted the measure on the right, which is obtained from $\mu_h$ by putting this measure on the imaginary axis, by $i\mu_h$.

\medskip

(3) As for the second assertion, this follows by carefully examining the above computation. Indeed, we used only at one point the independence of $g,h$, so for the formula $\mu_f=\mu_g*i\mu_h$ to hold, the equality used at that point, which is as follows, must hold:
$$\sum_{x,y\in\mathbb R}P(g=x,h=y)\delta_{x+iy}
=\sum_{x,y\in\mathbb R}P(g=x)P(h=y)\delta_{x+iy}$$

But this is the same as saying that the following must hold, for any $x,y$:
$$P(g=x,h=y)=P(g=x)P(h=y)$$

We conclude that, in order for the decomposition formula $\mu_f=\mu_g*i\mu_h$ to hold, the real and imaginary parts $g,h:X\to\mathbb R$ must be independent, as stated.
\end{proof}

As an alternative to the above decomposition into real and imaginary parts, we can write our  complex variable $f:X\to\mathbb C$ in polar form, as follows:
$$f=ge^{ih}\quad,\quad g=|f|,\ e^{ih}=f/|f|$$

In practice, this method leads, among others, to the following useful result:

\index{real part}
\index{imaginary part}
\index{rotation of variable}

\begin{theorem}
For a complex random variable $f:X\to\mathbb C$, written in polar form as $f=ge^{ih}$, and with $g,h$ assumed to be independent, we have
$$\mu_f=\mu_g\times e^{i\mu_h}$$
with $\times$ being the multiplicative convolution, $\delta_z\times\delta_t=\delta_{zt}$, and with $\mu\to e^{i\mu}$ denoting the wrapped version, $\mathbb R\to\mathbb T$. If $g,h$ are not independent, this formula does not hold.
\end{theorem}

\begin{proof}
This is very similar to the proof of Theorem 5.11. In the discrete case, assuming that the polar decomposition is $f=ge^{ih}$, with the components $g:X\to[0,\infty)$ and $h:X\to[0,2\pi)$ being independent, we have the following computation:
\begin{eqnarray*}
\mu_f
&=&\sum_{z\in\mathbb C}P(f=z)\delta_z\\
&=&\sum_{r,t}P(ge^{ih}=re^{it})\delta_{re^{it}}\\
&=&\sum_{r,t}P(g=r,h=t)\delta_{re^{it}}\\
&=&\sum_{r,t}P(g=r)P(h=t)\delta_r\times\delta_{e^{it}}\\
&=&\left(\sum_rP(g=r)\delta_r\right)\times
\left(\sum_tP(h=t)\delta_{e^{it}}\right)\\
&=&\mu_g\times e^{i\mu_h}
\end{eqnarray*}

Conversely, this computation shows as well that, in order to have $\mu_f=\mu_g\times e^{i\mu_h}$, the variables $g,h$ must assumed to be independent. As for the discussion in the general, continuous case, this is again something straightforward, say good exercise for you.
\end{proof}

\section*{5b. Central limits}

Getting now to central limits, in analogy with what we did in chapter 2 in the real case, we would like to investigate the following fundamental question:

\begin{question}
Given complex variables $f_1,f_2,f_3,\ldots:X\to\mathbb C$, assumed to be i.i.d. and centered, what can we say about their rescaled partial sums
$$S_n=\frac{f_1+\ldots+f_n}{\sqrt{n}}$$
in the $n\to\infty$ limit? 
\end{question}

And good question this is. In answer, we already have some good input for this from chapter 2, because the CLT established there has the following consequences:

\bigskip

(1) When our variables are assumed to be real, $f_1,f_2,f_3,\ldots:X\to\mathbb R$, we have $S_n\sim g_t$ in the $n\to\infty$ limit, with $t=V(f_i)$ being the common variance of our variables. 

\bigskip

(2) Similarly, in the purely imaginary case, $f_1,f_2,f_3,\ldots:X\to i\mathbb R$, we have $S_n\sim ig_t$ in the $n\to\infty$ limit, with $t=V(f_i)$ being the common variance of our variables.

\bigskip

(3) More generally, in the radial case, $f_1,f_2,f_3,\ldots:X\to w\mathbb R$ with $w\in\mathbb T$, we have $S_n\sim wg_t$ in the $n\to\infty$ limit, with $t=V(f_i)$ being again the common variance.

\bigskip

In order to deal now with Question 5.13 as stated, without further assumptions on our variables $f_i$, the idea will be very simple, namely decomposing them as follows:
$$f_i=h_i+ik_i$$

In practice, this needs a bit of preliminary discussion, first in regards with independence, and then in regards with the variance. Regarding independence, we have:

\begin{proposition}
Assuming that $f,f':X\to\mathbb C$ are independent, if we write
$$f=h+ik\quad,\quad f'=h'+ik'$$
then the variables $h,k:X\to\mathbb R$ are independent from $h',k':X\to\mathbb R$.
\end{proposition}

\begin{proof}
In order to prove the result in the discrete case, to start with, observe that with $f=h+ik$, $f'=h'+ik'$ as above, we have the following formulae:
$$P(f=a+ib,f'=a'+ib')=P(h=a,k=b,h'=a',k'=b')$$
$$P(f=a+ib)P(f'=a'+ib')=P(h=a,k=b)P(h'=a',k'=b')$$

Now assuming that the variables $f,f':X\to\mathbb C$ are independent, the quantities on the left must be equal, so the quantities on the right must be equal too:
$$P(h=a,k=b,h'=a',k'=b')=P(h=a,k=b)P(h'=a',k'=b')$$

But by summing over $b,b'$ we obtain that $h,h'$ are independent, as shown by:
\begin{eqnarray*}
P(h=a,h'=a')
&=&\sum_{b,b'}P(h=a,k=b,h'=a',k'=b')\\
&=&\sum_{b,b'}P(h=a,k=b)P(h'=a',k'=b')\\
&=&P(h=a)P(h'=a')
\end{eqnarray*}

Similarly, by summing over $b,a'$ we obtain that $h,k'$ are independent, by summing over $a,b'$ we obtain that $k,h'$ are independent, and by summing over $a,a'$ we obtain that $k,k'$ are independent. Thus we proved the result in the discrete case, and the extension to the general, continuous case is straightforward, say good exercise for you.
\end{proof}

The second piece of discussion now, regarding the variance, is as follows:

\begin{proposition}
Given a variable $f:X\to\mathbb C$, decomposed as $f=h+ik$, we have
$$E(f)=E(h)+iE(k)$$
and assuming that this quantity vanishes, we have the following formulae,
$$V(h)=\frac{V(f)+Re(E(f^2))}{2}\quad,\quad V(k)=\frac{V(f)-Re(E(f^2))}{2}$$
and in particular, we have the formula $V(f)=V(h)+V(k)$.
\end{proposition}

\begin{proof}
The first assertion is clear. In order to compute now the variances of the real and imaginary parts, we can use the standard formulae for these parts, namely:
$$h=\frac{f+\bar{f}}{2}\quad,\quad k=\frac{f-\bar{f}}{2i}$$

Indeed, assuming $E(f)=0$, and so in particular $E(h)=0$, we have:
\begin{eqnarray*}
V(h)
&=&E\left(\left(\frac{f+\bar{f}}{2}\right)^2\right)\\
&=&\frac{E(f^2+\bar{f}^2+2f\bar{f})}{4}\\
&=&\frac{E(|f|^2)+(E(f^2)+E(\bar{f}^2))/2}{2}\\
&=&\frac{V(f)+Re(E(f^2))}{2}
\end{eqnarray*}

Similarly, again by assuming $E(f)=0$, which gives $E(k)=0$, we have:
\begin{eqnarray*}
V(k)
&=&E\left(\left(\frac{f-\bar{f}}{2i}\right)^2\right)\\
&=&\frac{E(-f^2-\bar{f}^2+2f\bar{f})}{4}\\
&=&\frac{E(|f|^2)-(E(f^2)+E(\bar{f}^2))/2}{2}\\
&=&\frac{V(f)-Re(E(f^2))}{2}
\end{eqnarray*}

As for the last assertion, this comes from this, or just directly, as follows:
$$V(f)=E(|h+ik|^2)=E(h^2+k^2)=V(h)+V(k)$$

Thus, we are led to the conclusions in the statement.
\end{proof}

We can come back now to the complex CLT problematics, and we first have:

\index{CCLT}
\index{Complex CLT}

\begin{theorem}[CCLT 1]
Given complex variables $f_1,f_2,f_3,\ldots:X\to\mathbb C$, assumed to be i.i.d. and centered, if we write their rescaled partial sums as
$$\frac{f_1+\ldots+f_n}{\sqrt{n}}=h+ik$$
then $h\sim g_s$ and $k\sim g_t$ with $n\to\infty$, where $s=V(h_i)$ and $t=V(k_i)$, with $f_i=h_i+ik_i$. 
\end{theorem}

\begin{proof}
This comes indeed from the real CLT, from chapter 2, as follows:

\medskip

(1) Let us decompose indeed our variables as $f_i=h_i+ik_i$, and consider the following two sequences, which by Proposition 5.14 consist of independent variables:
$$h_1,h_2,h_3,\ldots:X\to\mathbb R$$
$$k_1,k_2,k_3,\ldots:X\to\mathbb R$$

(2) Our claim is that both the above sequences are i.i.d. Indeed, recall from the proof of Proposition 5.15 that for $f=h+ik$, we have the following formulae:
$$h=\frac{f+\bar{f}}{2}\quad,\quad k=\frac{f-\bar{f}}{2i}$$

Thus the laws of both $h,k$ can be recaptured from the law of $f$, and getting back now to the sequences in (1), since the variables $f_i=h_i+ik_i$ were assumed to be identically distributed, it follows that the variables $h_i$ are identically distributed, and that the variables $k_i$ are identically distributed too. Thus, we have two i.i.d. sequences, as claimed.

\medskip

(3) But with this, done, because the real CLT from chapter 2 applies to both sequences $h_1,h_2,h_3,\ldots:X\to\mathbb R$ and $k_1,k_2,k_3,\ldots:X\to\mathbb R$ from (1), and gives the result.

\medskip

(4) Next, observe that in what regards the values of the parameters $s,t$, we have the following aternative formulae for them, coming from Proposition 5.15, which can be quite convenient, making no explicit reference to the decompositions $f_i=h_i+ik_i$:
$$s=\frac{V(f_i)+Re(E(f_i^2))}{2}\quad,\quad t=\frac{V(f_i)-Re(E(f_i^2))}{2}$$

(5) Finally, observe that we have is not a true CLT, because we know nothing about the relation between $h,k$, so we cannot say something about $h+ik$ in the $n\to\infty$ limit. However, our result belongs to the CLT galaxy, and we chose to call it CCLT 1.
\end{proof}

In order to have now convergence, we need to know more about the relation between $h_i,k_i$, when decomposing $f_i=h_i+ik_i$. And with the notion of independence being the simplest such relation, we are led in this way to the following statement:

\begin{theorem}[CCLT 2]
Given i.i.d. centered variables $f_1,f_2,f_3,\ldots:X\to\mathbb C$, with their real parts assumed to be independent from their imaginary parts, we have
$$\frac{f_1+\ldots+f_n}{\sqrt{n}}\sim C_{st}$$
with $n\to\infty$, with the limiting law on the right being given by the formula
$$C_{st}=law(h+ik)$$
with $h\sim g_s$ and $k\sim g_t$, independent, where $s=V(h_i)$ and $t=V(k_i)$, with $f_i=h_i+ik_i$. 
\end{theorem}

\begin{proof}
This comes as a continuation of Theorem 5.16. Indeed, we know from there that if we write $f_i=h_i+ik_i$, then we have, in the $n\to\infty$ limit:
$$\frac{h_1+\ldots+h_n}{\sqrt{n}}\sim g_s\quad,\quad \frac{k_1+\ldots+k_n}{\sqrt{n}}\sim g_t$$

Now since in the present context, the whole family $\{h_i,k_i\}$ was assumed to be independent, the above partial sums are independent. Thus, their limits $h,k$ are independent too, and so $h+ik$ converges with $n\to\infty$ to the law $C_{st}$ in the statement.
\end{proof}

Coming next, still with $f_i=h_i+ik_i$ as usual, let us examine what happens when imposing the strongest condition, namely that the family $\{h_i,k_i\}$ is i.i.d. Here we basically have to set $s=t$ in Theorem 5.17, and we are done. However, before doing so, remember from Proposition 5.15 the following formula, for a centered variable $f=h+ik$:
$$V(f)=V(h)+V(k)$$ 

Now the point is that if we set $s=t$ in Theorem 5.17, as suggested above, this would lead to a limiting law having variance $2t$. And the 2 factor might bring some troubles later, when investigating the limiting law, independently of the CLT.

\bigskip

Thus, we must divide everything by $\sqrt{2}$, as to have in the end a variable having complex variance $t$, in the sense of Definition 5.7. And, we are led in this way to:

\index{complex normal law}
\index{complex Gaussian law}
\index{CCLT}
\index{Complex CLT}
\index{Complex Central Limit Theorem}
\index{complex variables}

\begin{theorem}[CCLT 3]
Given complex variables $f_1,f_2,f_3,\ldots$ whose real and imaginary parts are i.i.d. and centered, and having complex variance $t>0$, we have
$$\frac{f_1+\ldots+f_n}{\sqrt{n}}\sim G_t$$
with the limiting law, called complex normal, or Gaussian law, being given by
$$G_t=law\left(\frac{h+ik}{\sqrt{2}}\right)$$
where $h,k$ are real and independent, each following the law $g_t$.
\end{theorem}

\begin{proof}
This follows indeed from our previous CCLT result, from Theorem 5.17, by setting $s=t$ and dividing everything by $\sqrt{2}$, as explained in the above.
\end{proof}

And with this, end of our discussion regarding the CCLT. We should mention that all this was just an introduction to the subject, and that far more things can be said, notably with some improvements in the exact independence assumptions, and also with the fact that the convergence is in fact in probability. For more on all this, we refer to a solid probability book. In what concerns us, what we have in the above will do.

\section*{5c. Wick formula}

What is next? Systematic study of the complex normal laws $G_t$ that we found, in analogy with what we know from Part I, regarding the real normal laws $g_t$. And as a first result here, exactly as in the real case, our laws form convolution semigroups:

\index{convolution semigroup}

\begin{theorem}
The complex Gaussian laws have the property
$$G_s*G_t=G_{s+t}$$
for any $s,t>0$, and so they form a convolution semigroup.
\end{theorem}

\begin{proof}
This can be seen in two possible ways, both instructive, as follows:

\medskip

(1) Consider a sequence of variables $f_1,f_2,f_3,\ldots\,$, as in the third CCLT, with common complex variance $s>0$, and consider as well a sequence of variables $f_1',f_2',f_3',\ldots\,$, again as in the third CCLT, and independent from the first ones, with common complex variance $t>0$. By the CCLT we have the following convergences, with $n\to\infty$:
$$\frac{f_1+\ldots+f_n}{\sqrt{n}}\sim G_s\quad,\quad 
\frac{f'_1+\ldots+f'_n}{\sqrt{n}}\sim G_t$$

On the other hand, the CCLT applies as well to the sums $f_i+f_i'$, which according to Proposition 5.15 have common complex variance $s+t$, and gives:
$$\frac{f_1+f_1'+\ldots+f_n+f_n'}{\sqrt{n}}\sim G_{s+t}$$

Now by comparing, we are led to the formula $G_s*G_t=G_{s+t}$ in the statement.

\medskip

(2) Alternatively, and without reference to the CCLT, consider real variables $h,k\sim g_s$ and $h',k'\sim g_t$, taken independent. We have then the following computation:
\begin{eqnarray*}
G_s*G_t
&=&law\left(\frac{h+ik}{\sqrt{2}}\right)*\,law\left(\frac{h'+ik'}{\sqrt{2}}\right)\\
&=&law\left(\frac{h+ik}{\sqrt{2}}+\frac{h'+ik'}{\sqrt{2}}\right)\\
&=&law\left(\frac{(h+h')+i(k+k')}{\sqrt{2}}\right)\\
&=&G_{s+t}
\end{eqnarray*}

Thus, one way or another, we are led to the formula in the statement.
\end{proof}

Regarding now the moments, the situation here is more complicated than in the real case, because in order to have good results, we have to deal with both the complex variables, and their conjugates. Let us formulate the following definition:

\index{colored integers}
\index{colored moments}

\begin{definition}
The moments a complex variable $f:X\to\mathbb C$ are the numbers
$$M_k=E(f^k)$$
depending on colored integers $k=\circ\bullet\bullet\circ\ldots\,$, with the conventions
$$f^\emptyset=1\quad,\quad f^\circ=f\quad,\quad f^\bullet=\bar{f}$$
and multiplicativity, in order to define the colored powers $f^k$.
\end{definition}

As an illustration for this notion, which is something very intuitive, here are the formulae of the four possible order 2 moments of a complex variable $f$:
$$M_{\circ\circ}=E(f^2)\quad,\quad M_{\circ\bullet}=E(f\bar{f})$$
$$M_{\bullet\circ}=E(\bar{f}f)\quad,\quad M_{\bullet\bullet}=E(\bar{f}^2)$$

Observe that, since $f,\bar{f}$ commute, we have the following identity, which shows that there is a bit of redundancy in our above definition, as formulated:
$$M_{\circ\bullet}=M_{\bullet\circ}$$

In fact, again since $f,\bar{f}$ commute, we can permute terms, in the general context of Definition 5.20, and restrict the attention to exponents of the following type:
$$k=\ldots\circ\circ\circ\bullet\bullet\bullet\bullet\ldots$$

However, our results about the complex Gaussian laws, and other complex laws, later on, not to talk about laws of matrices, random matrices and other noncommuting variables, that will appear later too, will look better without doing this. So, we will use Definition 5.20 as stated. Getting to work now, we first have the following result:

\index{moments of normal law}

\begin{theorem}
The moments of the complex normal law are given by
$$M_k(G_t)=\begin{cases}
t^pp!&(k\ {\rm uniform, of\ length}\ 2p)\\
0&(k\ {\rm not\ uniform})
\end{cases}$$
where $k=\circ\bullet\bullet\circ\ldots$ is called uniform when it contains the same number of $\circ$ and $\bullet$.
\end{theorem}

\begin{proof}
We must compute the moments, with respect to colored integer exponents $k=\circ\bullet\bullet\circ\ldots$ as above, of the variable from Theorem 5.18, namely:
$$f=\frac{g+ih}{\sqrt{2}}$$

We can assume that we are in the case $t=1$, and the proof here goes as follows:

\medskip

(1) As a first observation, in the case where our exponent $k=\circ\bullet\bullet\circ\ldots$ is not uniform, a standard rotation argument shows that the corresponding moment of $f$ vanishes. To be more precise, the variable $f'=wf$ is complex Gaussian too, for any complex number $w\in\mathbb T$, and from $M_k(f)=M_k(f')$ we obtain $M_k(f)=0$, in this case.

\medskip

(2) In the uniform case now, where the exponent $k=\circ\bullet\bullet\circ\ldots$ consists of $p$ copies of $\circ$ and $p$ copies of $\bullet$\,, the corresponding moment can be computed as follows:
\begin{eqnarray*}
M_k
&=&E[(f\bar{f})^p]\\
&=&\frac{1}{2^p}E[(g^2+h^2)^p]\\
&=&\frac{1}{2^p}\sum_{r=0}^p\binom{p}{r}E(g^{2r})E(h^{2p-2r})\\
&=&\frac{1}{2^p}\sum_{r=0}^p\binom{p}{r}(2r)!!(2p-2r)!!\\
&=&\frac{1}{2^p}\sum_{r=0}^p\frac{p!}{r!(p-r)!}\cdot\frac{(2r)!}{2^rr!}\cdot\frac{(2p-2r)!}{2^{p-r}(p-r)!}\\
&=&\frac{p!}{4^p}\sum_{r=0}^p\binom{2r}{r}\binom{2p-2r}{p-r}
\end{eqnarray*}

(3) In order to finish now the computation, let us recall that we have the following formula, coming from the generalized binomial formula, with exponent $-1/2$:
$$\frac{1}{\sqrt{1+t}}=\sum_{q=0}^\infty\binom{2q}{q}\left(\frac{-t}{4}\right)^q$$

By taking the square of this series, we obtain the following formula:
$$\frac{1}{1+t}
=\sum_{p=0}^\infty\left(\frac{-t}{4}\right)^p\sum_{r=0}^p\binom{2r}{r}\binom{2p-2r}{p-r}$$

Now by looking at the coefficient of $t^p$ on both sides, we conclude that the sum on the right equals $4^p$. Thus, we can finish the moment computation in (2), as follows:
$$M_k=\frac{p!}{4^p}\times 4^p=p!$$

We are therefore led to the conclusion in the statement.
\end{proof}

As before with the real normal laws, or even before with the Poisson laws, a better-looking statement regarding the moments is in terms of partitions, as follows:

\index{matching pairings}

\begin{theorem}
The moments of the complex normal law are the numbers
$$M_k(G_t)=\sum_{\pi\in\mathcal P_2(k)}t^{|\pi|}$$
with $\mathcal P_2(k)$ being the pairings of $\{1,\ldots,k\}$ which are matching, pairing $\circ-\bullet$ symbols.
\end{theorem}

\begin{proof}
This is a reformulation of Theorem 5.21. Indeed, we can assume that we are in the case $t=1$, and here we know from Theorem 5.21 that the moments are:
$$M_k=\begin{cases}
(|k|/2)!&(k\ {\rm uniform})\\
0&(k\ {\rm not\ uniform})
\end{cases}$$

On the other hand, the numbers $|\mathcal P_2(k)|$ are given by exactly the same formula. Indeed, in order to have a matching pairing of $k$, our exponent $k=\circ\bullet\bullet\circ\ldots$ must be uniform, consisting of $p$ copies of $\circ$ and $p$ copies of $\bullet$, with $p=|k|/2$. But then the matching pairings of $k$ correspond to the permutations of the $\bullet$ symbols, as to be matched with $\circ$ symbols, and so we have $p!$ such pairings. Thus, we have the same formula as for the moments of $f$, and we are led to the conclusion in the statement.
\end{proof}

In practice, we also need to know how to compute joint moments. We have here:

\index{joint moments}
\index{Wick formula}

\begin{theorem}[Wick formula]
Given independent variables $f_i$, each following the complex normal law $G_t$, with $t>0$ being a fixed parameter, we have the formula
$$E\left(f_{i_1}^{k_1}\ldots f_{i_s}^{k_s}\right)=t^{s/2}\#\left\{\pi\in\mathcal P_2(k)\Big|\pi\leq\ker i\right\}$$
where $k=k_1\ldots k_s$ and $i=i_1\ldots i_s$, for the joint moments of these variables, where $\pi\leq\ker i$ means that the indices of $i$ must fit into the blocks of $\pi$, in the obvious way.
\end{theorem}

\begin{proof}
This is something well-known, which can be proved as follows:

\medskip

(1) Let us first discuss the case where we have a single variable $f$, which amounts in taking $f_i=f$ for any $i$ in the formula in the statement. What we have to compute here are the moments of $f$, with respect to colored integer exponents $k=\circ\bullet\bullet\circ\ldots\,$, and the formula in the statement tells us that these moments must be:
$$E(f^k)=t^{|k|/2}|\mathcal P_2(k)|$$

But this is the formula in Theorem 5.22, so we are done with this case.

\medskip

(2) In general now, when expanding the product $f_{i_1}^{k_1}\ldots f_{i_s}^{k_s}$ and rearranging the terms, we are left with doing a number of computations as in (1), and then making the product of the expectations that we found. But this amounts in counting the partitions in the statement, with the condition $\pi\leq\ker i$ there standing for the fact that we are doing the various type (1) computations independently, and then making the product.
\end{proof}

The above statement is one of the possible formulations of the Wick formula, and there are many more formulations, which are all useful. For instance, we have:

\index{Wick formula}

\begin{theorem}[Wick formula 2]
Given independent variables $f_i$, each following the complex normal law $G_t$, with $t>0$ being a fixed parameter, we have the formula
$$E\left(f_{i_1}\ldots f_{i_k}\bar{f}_{j_1}\ldots\bar{f}_{j_k}\right)=t^k\#\left\{\pi\in S_k\Big|i_{\pi(r)}=j_r,\forall r\right\}$$
for the non-vanishing joint moments of these variables.
\end{theorem}

\begin{proof}
This follows from the usual Wick formula, from Theorem 5.23. With some changes in the indices and notations, the formula there reads:
$$E\left(f_{I_1}^{K_1}\ldots f_{I_s}^{K_s}\right)=t^{s/2}\#\left\{\sigma\in\mathcal P_2(K)\Big|\sigma\leq\ker I\right\}$$

Now observe that we have $\mathcal P_2(K)=\emptyset$, unless the colored integer $K=K_1\ldots K_s$ is uniform, in the sense that it contains the same number of $\circ$ and $\bullet$ symbols. Up to permutations, the non-trivial case, where the moment is non-vanishing, is the case where the colored integer $K=K_1\ldots K_s$ is of the following special form:
$$K=\underbrace{\circ\circ\ldots\circ}_k\ \underbrace{\bullet\bullet\ldots\bullet}_k$$

So, let us focus on this case, which is the non-trivial one. Here we have $s=2k$, and we can write the multi-index $I=I_1\ldots I_s$ in the following way:
$$I=i_1\ldots i_k\ j_1\ldots j_k$$

With these changes made, the above usual Wick formula reads:
$$E\left(f_{i_1}\ldots f_{i_k}\bar{f}_{j_1}\ldots\bar{f}_{j_k}\right)=t^k\#\left\{\sigma\in\mathcal P_2(K)\Big|\sigma\leq\ker(ij)\right\}$$

The point now is that the matching pairings $\sigma\in\mathcal P_2(K)$, with $K=\circ\ldots\circ\bullet\ldots\bullet\,$, of length $2k$, as above, correspond to the permutations $\pi\in S_k$, in the obvious way. With this identification made, the above modified usual Wick formula becomes:
$$E\left(f_{i_1}\ldots f_{i_k}\bar{f}_{j_1}\ldots\bar{f}_{j_k}\right)=t^k\#\left\{\pi\in S_k\Big|i_{\pi(r)}=j_r,\forall r\right\}$$

Thus, we have reached to the formula in the statement, and we are done.
\end{proof}

Finally, here is one more formulation of the Wick formula, useful as well:

\index{Wick formula}

\begin{theorem}[Wick formula 3]
Given independent variables $f_i$, each following the complex normal law $G_t$, with $t>0$ being a fixed parameter, we have the formula
$$E\left(f_{i_1}\bar{f}_{j_1}\ldots f_{i_k}\bar{f}_{j_k}\right)=t^k\#\left\{\pi\in S_k\Big|i_{\pi(r)}=j_r,\forall r\right\}$$
for the non-vanishing joint moments of these variables.
\end{theorem}

\begin{proof}
This follows from our second Wick formula, from Theorem 5.24, simply by permuting the terms, as to have an alternating sequence of plain and conjugate variables. Alternatively, we can start with Theorem 5.23, and then perform the same manipulations as in the proof of Theorem 5.24, but with the exponent being this time as follows: 
$$K=\underbrace{\circ\bullet\circ\bullet\ldots\ldots\circ\bullet}_{2k}$$

Thus, we are led to the conclusion in the statement.
\end{proof}

And with this, end of our study regarding the complex normal variables. We will heavily use all this in chapters 7-8 below, when doing groups and random matrices.

\section*{5d. Poisson, Bessel}

There is no Gauss without Poisson, and vice versa, as we say in probability, so let us get now into the question of finding the complex analogues of the Poisson laws $p_t$. But here, as we perfectly know from chapter 4, the Poisson laws can spread to the whole real line $\mathbb R$ by using the compound Poisson procedure. And the point is that exactly the same procedure works in the complex case, with the essentials being as follows:

\index{CPLT}
\index{Compound PLT}

\begin{theorem}[CPLT]
Given a compactly supported complex measure $\nu$, of mass $c>0$, the following compound Poisson limit converges:
$$p_\nu=\lim_{n\to\infty}\left(\left(1-\frac{c}{n}\right)\delta_0+\frac{1}{n}\,\nu\right)^{*n}$$
In the discrete case, $\nu=\sum_{k=1}^sc_k\delta_{z_k}$ with $c_k>0$ and $z_k\in\mathbb C$, the Fourier transform is
$$F_{p_\nu}(x)=\exp\left(\sum_{k=1}^sc_k(e^{ixz_k}-1)\right)$$
and the measure itself is given by the following formula, with $f_k\sim p_{c_k}$, independent:
$$p_\nu=law\left(\sum_{k=1}^sz_kf_k\right)$$
As main examples for this construction, $p_t=p_{t\delta_1}$ and $p_t^2=p_{t(\delta_{-1}+\delta_1)/2}$.
\end{theorem}

\begin{proof}
This is something that we know from chapter 4 in the real case, and the proof in the complex case is similar. Here are are details, reproduced for convenience:

\medskip

(1) Regarding the convergence, in the discrete case, which is the one that we are interested in, this follows from the alternative formula for $p_\nu$ in the statement, that we will prove below. As for the general case, this comes from the discrete case by using standard limiting arguments, and we will leave clarifying this as an exercise. 

\medskip

(2) So, assume that we are in the discrete case, $\nu=\sum_{k=1}^sc_k\delta_{z_k}$ with $c_k>0$ and $z_k\in\mathbb C$. If we denote by $\mu_n$ the measure under the convolution sign, we have:
\begin{eqnarray*}
F_{\mu_n}(x)=\left(1-\frac{c}{n}\right)+\frac{1}{n}\sum_{k=1}^sc_ke^{ixz_k}
&\implies&F_{\mu_n^{*n}}(x)=\left(\left(1-\frac{c}{n}\right)+\frac{1}{n}\sum_{k=1}^sc_ke^{ixz_k}\right)^n\\
&\implies&F_{p_\nu}(x)=\exp\left(\sum_{k=1}^sc_k(e^{ixz_k}-1)\right)
\end{eqnarray*}

Thus, we have indeed the Fourier transform formula in the statement.

\medskip

(3) Next, if we denote by $f$ the sum of Poisson variables in the statement, we have:
\begin{eqnarray*}
F_{f_k}(x)=\exp(c_k(e^{ix}-1))
&\implies&F_{z_kf_k}(x)=\exp(c_k(e^{ixz_k}-1))\\
&\implies&F_f(x)=\exp\left(\sum_{k=1}^sc_k(e^{ixz_k}-1)\right)
\end{eqnarray*}

Thus we have the same Fourier transform formula as in (2), as desired.

\medskip

(4) Finally, the formula $p_t=p_{t\delta_1}$ is the usual PLT, which in our context comes from the formula of $p_\nu$ proved in (3), and the formula $p_t^2=p_{t(\delta_{-1}+\delta_1)/2}$ is something that we discussed in chapter 4, coming there via the Fourier formula established in (2).
\end{proof}

Moving on, and still with the idea in mind of finding the correct complex analogues of the Poisson law $p_t$, the two examples that we have, which are both in the real case, namely $p_t=p_{t\delta_1}$ and $p_t^2=p_{t(\delta_{-1}+\delta_1)/2}$, suggest formulating the following definition:

\index{generalized Bessel laws}
\index{complex Bessel laws}

\begin{definition}
The Bessel law of level $s\in\mathbb N\cup\{\infty\}$ and parameter $t>0$ is
$$p_t^s=p_{t\varepsilon_s}$$
with $\varepsilon_s$ being the uniform measure on the $s$-th roots of unity.
\end{definition}

Which sounds quite good and conceptual, hope you agree with me. As a first observation, at $s=1,2$ we obtain the Poisson laws, and the previous real Bessel laws:
$$p^1_t=p_t\quad,\quad 
p^2_t=p_t^2$$

Another important particular case is $s=\infty$, where we obtain a measure which is actually not discrete, that we will call complex Bessel law, and denote as follows:
$$p^\infty_t=P_t$$

Let us develop now some theory, for these Bessel laws that we found. To start with, save for the Poisson limits, what we have in Theorem 5.26 translates into:

\begin{theorem}
The Bessel laws $p_t^s=p_{t\varepsilon_s}$ are given by the following formula, with $f_1,\ldots,f_s$ being Poisson $(t/s)$ and independent, and $w=e^{2\pi i/s}$:
$$p^s_t=law\left(\sum_{k=1}^sw^kf_k\right)$$
The Fourier transform of $p^s_t$ is as follows, in terms of $\exp_sz=\sum_{n=0}^\infty z^{sn}/(sn)!$:
$$F(x)=\exp\left(t(\exp_s(ix)-1)\right)$$
Also, we have the convolution formula $p^s_r*p^s_t=p^s_{r+t}$, for any $r,t>0$.
\end{theorem}

\begin{proof}
The first formula comes from Theorem 5.26, according to our definition of $p^s_t$. Regarding the Fourier transform, again by using Theorem 5.26, we have:
\begin{eqnarray*}
\log F(x)
&=&\sum_{k=1}^s\frac{t}{s}\left(\exp(iw^kx)-1\right)\\
&=&t\left[\left(\frac{1}{s}\sum_{k=1}^s\exp(iw^kx)\right)-1\right]\\
&=&t\left[\left(\frac{1}{s}\sum_{m=0}^\infty\frac{(ix)^m}{m!}\sum_{k=1}^sw^{km}\right)-1\right]\\
&=&t\left(\sum_{n=0}^\infty\frac{(ix)^{sn}}{(sn)!}-1\right)\\
&=&t\left(\exp_s(ix)-1\right)
\end{eqnarray*}

Observe by the way that the function $\exp_sz$ is given at $s=1,2$ by $\exp_1=\exp$ and $\exp_2=\cosh$, so our formula above is compatible with those from chapter 4, at $s=1,2$. Finally, the formula $p^s_r*p^s_t=p^s_{r+t}$ comes from the fact that $\log F$ is linear in $t$.
\end{proof}

In what regards now the density of the Bessel laws, the result here is as follows:

\begin{theorem}
The density of $p^s_t$ is given by the following formula,
$$p^s_t=e^{-t}\sum_{p_1=0}^\infty\ldots\sum_{p_s=0}^\infty\frac{1}{p_1!\ldots p_s!}\,\left(\frac{t}{s}\right)^{p_1+\ldots+p_s}\delta\left(\sum_{k=1}^sw^kp_k\right)$$
where $w=e^{2\pi i/s}$, and with the $\delta$ symbol standing as usual for a Dirac mass.
\end{theorem}

\begin{proof}
In order to prove the result, we compute the Fourier transform of the measure on the right. This is given by the following formula:
\begin{eqnarray*}
F(x)
&=&e^{-t}\sum_{p_1=0}^\infty\ldots\sum_{p_s=0}^\infty\frac{1}{p_1!\ldots p_s!}\left(\frac{t}{s}\right)^{p_1+\ldots+p_s}F\delta\left(\sum_{k=1}^sw^kp_k\right)(x)\\
&=&e^{-t}\sum_{p_1=0}^\infty\ldots\sum_{p_s=0}^\infty\frac{1}{p_1!\ldots p_s!}\left(\frac{t}{s}\right)^{p_1+\ldots+p_s}\exp\left(\sum_{k=1}^siw^kp_kx\right)\\
&=&e^{-t}\sum_{r=0}^\infty\left(\frac{t}{s}\right)^r\sum_{\sum p_i=r}\frac{\exp\left(\sum_{k=1}^siw^kp_kx\right)}{p_1!\ldots p_s!}
\end{eqnarray*}

We multiply by $e^t$, and we compute the derivative with respect to $t$:
\begin{eqnarray*}
(e^tF(x))'
&=&\sum_{r=1}^\infty\frac{r}{s}\left(\frac{t}{s}\right)^{r-1}\sum_{\sum p_i=r}\frac{\exp\left(\sum_{k=1}^siw^kp_kx\right)}{p_1!\ldots p_s!}\\
&=&\frac{1}{s}\sum_{r=1}^\infty\left(\frac{t}{s}\right)^{r-1}\sum_{\sum p_i=r}\left(\sum_{l=1}^sp_l\right)\frac{\exp\left(\sum_{k=1}^siw^kp_kx\right)}{p_1!\ldots p_s!}\\
&=&\frac{1}{s}\sum_{r=1}^\infty\left(\frac{t}{s}\right)^{r-1}\sum_{\sum p_i=r}\sum_{l=1}^s\frac{\exp\left(\sum_{k=1}^siw^kp_kx\right)}{p_1!\ldots p_{l-1}!(p_l-1)!p_{l+1}!\ldots p_s!}\\
\end{eqnarray*}

By using the variable $u=r-1$, we obtain from this the following formula:
\begin{eqnarray*}
(e^tF(x))'
&=&\frac{1}{s}\sum_{u=0}^\infty\left(\frac{t}{s}\right)^u\sum_{\sum q_i=u}\sum_{l=1}^s\frac{\exp\left(iw^lx+\sum_{k=1}^siw^kq_kx\right)}{q_1!\ldots q_s!}\\
&=&\left(\frac{1}{s}\sum_{l=1}^s\exp(iw^lx)\right)\left(\sum_{u=0}^\infty\left(\frac{t}{s}\right)^u\sum_{\sum q_i=u}\frac{\exp\left(\sum_{k=1}^siw^kq_kx\right)}{q_1!\ldots q_s!}\right)\\
&=&\exp_s(ix)\cdot e^tF(x)
\end{eqnarray*}

But $\Phi(t)=\exp(t\exp_s(ix))$ satisfies as well $\Phi'(t)=(\exp_s(ix))\Phi(t)$, so we have:
$$e^tF(x)=\exp(t\exp_s(ix))$$

Thus $\log F(x)=t(\exp_s(ix)-1)$, so the measure in statement is indeed $p^s_t$.
\end{proof}

In order to deal now with moments, the most elegant is to use group theory:

\index{complex reflection group}
\index{roots of unity}

\begin{definition}
The complex reflection group $H_N^s\subset U_N$, depending on parameters
$$N\in\mathbb N\quad,\quad s\in\mathbb N\cup\{\infty\}$$
is the group of permutation-type matrices with $s$-th roots of unity as entries,
$$H_N^s=M_N(\mathbb Z_s\cup\{0\})\cap U_N$$
with the convention $\mathbb Z_\infty=\mathbb T$, at $s=\infty$.
\end{definition}

Observe that at $s=1,2$ we obtain the symmetric and hyperoctahedral groups:
$$H_N^1=S_N\quad,\quad 
H_N^2=H_N$$

Another important particular case is $s=\infty$, where we obtain a compact group which is actually not finite, but is of key importance, that we will denote as follows:
$$H_N^\infty=K_N$$

The relation with the Bessel laws $p^s_t$ comes from the following result:

\begin{theorem}
The following happen:
\begin{enumerate}
\item For the group $H_N^s$, the truncated characters satisfy $\chi_t\sim p^s_t$, with $N\to\infty$.

\item $H_N^s$ is easy, coming from $P^s$, the partitions satisfying $\#\circ=\#\bullet(s)$, in each block.

\item The moments of the Bessel law $p^s_1$ are the numbers $M_k(p^s_1)=|P^s(k)|$.

\item The moments of the Bessel law $p^s_t$ are the numbers $M_k(p^s_t)=\sum_{\pi\in P^s(k)}t^{|\pi|}$.
\end{enumerate}
\end{theorem}

\begin{proof}
This is very standard, as before in chapter 4 at $s=1,2$, as follows:

\medskip

(1) Let us first work out the case $t=1$. Since the limit probability for a random permutation to have exactly $k$ fixed points is $e^{-1}/k!$, we get, as before for $S_N,H_N$:
$$\lim_{N\to\infty}law(\chi_1)=e^{-1}\sum_{k=0}^\infty \frac{1}{k!}\,\varepsilon_s^{*k}$$

On the other hand, by the very definition of the Bessel law $p^s_1$, we have:
$$p^s_1
=\lim_{N\to\infty}\left(\left(1-\frac{1}{N}\right)\delta_0+\frac{1}{N}\,\varepsilon_s\right)^{*N}
=e^{-1}\sum_{k=0}^\infty\frac{1}{k!}\,\varepsilon_s^{*k}$$

In the case where $t\in(0,1]$ is arbitrary, the same method works. We have:
$$\lim_{N\to\infty}law(\chi_t)=e^{-t}\sum_{k=0}^\infty\frac{t^k}{k!}\,\varepsilon_s^{*k}$$

On the other hand, by definition of the Bessel law $p^s_t$, we have as well:
$$p^s_t=\lim_{N\to\infty}\left(\left(1-\frac{1}{N}\right)\delta_0+\frac{1}{N}\,\varepsilon_s\right)^{*[tN]}=e^{-t}\sum_{k=0}^\infty\frac{t^k}{k!}\,\varepsilon_s^{*k}$$

(2) Again, this follows as before, for $S_N,H_N$. Consider indeed an arbitrary linear map $T:(\mathbb C^N)^{\otimes k}\to(\mathbb C^N)^{\otimes l}$, written as follows, with $\lambda\binom{i_1\ldots i_k}{j_1\ldots j_l}\in\mathbb C$ being certain scalars:
$$T(e_{i_1}\otimes\ldots\otimes e_{i_k})=\sum_{j_1\ldots j_l}\lambda
\begin{pmatrix}i_1&\ldots& i_k\\ j_1&\ldots& j_l\end{pmatrix}e_{j_1}\otimes\ldots\otimes e_{j_l}$$

Now pick a group element $g=\sigma^w\in H_N^s$. We have the following computation:
\begin{eqnarray*}
Tg^{\otimes k}(e_{i_1}\otimes\ldots\otimes e_{i_k})
&=&w_{i_1}\ldots w_{i_k}T(e_{\sigma(i_1)}\otimes\ldots\otimes e_{\sigma(i_k)})\\
&=&w_{i_1}\ldots w_{i_k}\sum_{j_1\ldots j_l}\lambda
\begin{pmatrix}\sigma(i_1)&\ldots&\sigma(i_k)\\ j_1&\ldots& j_l\end{pmatrix}e_{j_1}\otimes\ldots\otimes e_{j_l}\\
&=&w_{i_1}\ldots w_{i_k}\sum_{j_1\ldots j_l}\lambda
\begin{pmatrix}\sigma(i_1)&\ldots&\sigma(i_k)\\ \sigma(j_1)&\ldots&\sigma(j_l)\end{pmatrix}e_{\sigma(j_1)}\otimes\ldots\otimes e_{\sigma(j_l)}
\end{eqnarray*}

On the other hand, we have as well the following computation:
\begin{eqnarray*}
g^{\otimes l}T(e_{i_1}\otimes\ldots\otimes e_{i_k})
&=&g^{\otimes l}\sum_{j_1\ldots j_l}\lambda
\begin{pmatrix}i_1&\ldots& i_k\\ j_1&\ldots& j_l\end{pmatrix}e_{j_1}\otimes\ldots\otimes e_{j_l}\\
&=&\sum_{j_1\ldots j_l}w_{j_1}\ldots w_{j_l}\lambda
\begin{pmatrix}i_1&\ldots&i_k\\ j_1&\ldots&j_l\end{pmatrix}e_{\sigma(j_1)}\otimes\ldots\otimes e_{\sigma(j_l)}
\end{eqnarray*}

Now let us examine the equality $Tg^{\otimes k}=g^{\otimes l}T$. Leaving aside the $w$ factors, having this equality for any $g\in H_N^s$ amounts in saying that the function $\lambda$ must be of the form $\lambda\binom{i}{j}=\varphi(\ker\binom{i}{j})$, so that the linear map $T$ itself must be a linear combination of maps $T_\pi$, with $\pi\in P(k,l)$. But then, we have to take into account the equality $w_{i_1}\ldots w_{i_k}=w_{j_1}\ldots w_{j_l}$ too, which tells us that we must have $\pi\in P^s(k,l)$, as claimed.

\medskip

(3) This follows from (2), since the integrals of characters count the fixed points.

\medskip

(4) This follows again from (2), via the integration technology from chapter 4. 
\end{proof}

And for more on the Bessel laws, and on their free analogues too, we refer to \cite{bb+}.

\section*{5e. Exercises}

This was a quite standard chapter, and as exercises on this, we have:

\begin{exercise}
Clarify what we said regarding complex independence.
\end{exercise}

\begin{exercise}
Clarify what we said regarding various complex decompositions.
\end{exercise}

\begin{exercise}
Learn about the various technical versions of the CCLT.
\end{exercise}

\begin{exercise}
Learn more about the Wick formula, and its applications.
\end{exercise}

\begin{exercise}
Prove the convergence of compound Poisson limits, in general.
\end{exercise}

\begin{exercise}
Work out all the details, in the proof of easiness of $H_N^s$.
\end{exercise}

\begin{exercise}
Learn about the complex reflection groups, and their classification.
\end{exercise}

\begin{exercise}
In the Bessel law context, discuss what happens at $s=\infty$.
\end{exercise}

As bonus exercise, and no surprise here, start learning some complex analysis.

\chapter{Rayleigh laws}

\section*{6a. Chi variables}

We have seen so far the basic theory of the real and complex normal laws $g_t$ and $G_t$, and of their discrete analogues $p^s_t$ too, and this is certainly what is needed, for dealing with various probabilistic phenomena, originally appearing from central limits.

\bigskip

At a more advanced level, however, far more things can be said, the idea being that our laws $g_t,G_t,p^s_t$ are surrounded by a myriad of other interesting laws, appearing as variations of them, which must be mastered as well. We will explore this ``normal law galaxy'' in this chapter, and then later too, our guiding principle being as follows:

\index{chi law}
\index{chi squared law}

\begin{principle}
Given a normal variable $f:X\to\mathbb C$, we can talk about
$$\chi=law{(|f|)}\quad,\quad\chi^2=law(|f|^2)$$
both real probability measures, whose knowledge helps in understanding $law(f)$.
\end{principle}

In fact, we already met this principle in chapter 5, when computing the moments of variables $f\sim G_t$, the point being that, by standard symmetry arguments, the moments vanish, unless we are dealing with the moments of $|f|^2$. And as explained there, the moments of $|f|^2$ can be computed indeed, leading to the solution of the $f$ problem.

\bigskip

Before starting, let us specify as well the precise normal variables $f:X\to\mathbb C$ that we are interested in. We would like of course to deal with the main cases $f\sim g_t,G^t,p^s_t$, but also with shifts, and other multiparametric versions of our laws, which can be of interest too. And thinking a bit, there are many of those, with the list being as follows:

\begin{liste}
We can develop our $\chi,\chi^2$ philosophy for $f\sim g_t$ and $f\sim g_t^a$, then for the various types of complex normal variables, which are as follows, with $a=(b+ic)/\sqrt{2}$,
\begin{enumerate}
\item $f\sim G_t$, meaning $f=(h+ik)/\sqrt{2}$ with $h,k\sim g_t$, independent.

\item $f\sim G_t^a$, meaning $f=(h+ik)/\sqrt{2}$ with $h\sim g_t^b,k\sim g_t^c$, independent. 

\item $f\sim G_{st}$, meaning $f=(h+ik)/\sqrt{2}$ with $h\sim g_s,k\sim g_t$, independent.

\item $f\sim G_{st}^a$, meaning $f=(h+ik)/\sqrt{2}$ with $h\sim g_t^b,k\sim g_t^c$, independent. 
\end{enumerate}
and then for the discrete variables $f\sim p^s_t$, that we can further fine-tune too.
\end{liste}

And good list that we have here, 2 questions concerning 7 types of variables, making it for a total of 14 situations to be investigated. And in case this will prove not to be enough, we can still modify $p^s_t=p_{t\varepsilon_s}$ into $p^\varepsilon_t=p_{t\varepsilon}$, with $\varepsilon$ being an arbitrary measure on the unit circle, that we can further shift at $a\in\mathbb C$, and so on, if we want to. 

\bigskip

Getting to work now, let us first talk about $\chi,\chi^2$ for $f\sim g_t$ and $f\sim g_t^a$. We have here 4 straigthforward statements to be formulated, with the first one being:

\begin{theorem}
Assuming $f\sim g_t$, the law $\chi_{1t}=law(|f|)$ has density
$$\chi_{1t}=\sqrt{\frac{2}{\pi t}}\,e^{-x^2/2t}dx$$
on $[0,\infty)$, and the basic properties of this law $\chi_{1t}$ are as follows:
\begin{enumerate}
\item The mean is $E=\sqrt{2t/\pi}$, the variance is $V=t-2t/\pi$.

\item The even moments are $M_{2l}=t^l(2l)!!$.

\item The odd moments are $M_{2l+1}=\sqrt{(2t)^{2l+1}/\pi}\cdot l!$.

\item Equivalently, we have $M_k=\sqrt{(2t)^k/\pi}\cdot\Gamma((k+1)/2)$.

\item We have $\gamma=(4-\pi)\sqrt{2/(\pi-2)^3}$ and $\kappa=(3\pi^2-4\pi-12)/(\pi-2)^2$.
\end{enumerate}
\end{theorem}

\begin{proof}
The density of $f$ is $e^{-x^2/2t}/\sqrt{2\pi t}$ on the whole $\mathbb R$, symmetric with respect to 0, so when taking the absolute value the density will restrict to $[0,\infty)$, and double, leading to the formula in the statement. And with the comment that the 1 subscript in $\chi_{1t}$ stands for the dimensionality of our problem, which is 1. As for the rest:

\medskip

(1) The mean can be computed with the change of variables $x=\sqrt{y}$, as follows:
\begin{eqnarray*}
E
&=&\sqrt{\frac{2}{\pi t}}\int_0^\infty xe^{-x^2/2t}dx\\
&=&\sqrt{\frac{2}{\pi t}}\int_0^\infty \sqrt{y}e^{-y/2t}\,\frac{1}{2\sqrt{y}}\,dy\\
&=&\sqrt{\frac{1}{2\pi t}}\int_0^\infty e^{-y/2t}dy\\
&=&\sqrt{\frac{1}{2\pi t}}\Big[-2te^{-y/2t}\Big]_0^\infty\\
&=&\sqrt{\frac{2t}{\pi}}
\end{eqnarray*}

Regarding now the second moment, by symmetry this is the same as for $f$ itself, namely $M_2=t$. We conclude that the variance is $V=t-2t/\pi$, as stated.

\medskip

(2) By symmetry the even moments are as for $f$ itself, namely $M_{2l}=t^l(2l)!!$.

\medskip

(3) We can compute the odd moments a bit like we did for $E$ before, as follows:
\begin{eqnarray*}
M_{2l+1}
&=&\sqrt{\frac{2}{\pi t}}\int_0^\infty x^{2l+1}e^{-x^2/2t}dx\\
&=&\sqrt{\frac{2}{\pi t}}\int_0^\infty(\sqrt{2ty})^{2l+1}e^{-y}\,\sqrt{\frac{t}{2y}}\,dy\\
&=&\sqrt{\frac{(2t)^{2l+1}}{\pi}}\int_0^\infty y^le^{-y}dy\\
&=&\sqrt{\frac{(2t)^{2l+1}}{\pi}}\cdot l!
\end{eqnarray*}

Here we have used $\int_0^\infty y^le^{-y}dy=l!$, obtained by iterated partial integration.

\medskip

(4) Obviously, we are here into a gamma function business. The idea is that the gamma function is defined as follows, for any $s>0$, with the integral converging:
$$\Gamma(s)=\int_0^\infty x^{s-1}e^{-x}dx$$

By partial integration we have $\Gamma(s+1)=s\Gamma(s)$. When coupled with $\Gamma(1)=1$, which is clear, and with $\Gamma(1/2)=\sqrt{\pi}$, obtained with $x=y^2$, this gives, by recurrence:
$$\Gamma(l+1)=l!\quad,\quad \Gamma\left(l+\frac{1}{2}\right)=\frac{(2l)!!\sqrt{\pi}}{2^l}$$

Now observe that this fits with the formulae in (2,3), and gives the formula in the statement. Alternatively, we can establish this formula directly, as follows:
\begin{eqnarray*}
M_k
&=&\sqrt{\frac{2}{\pi t}}\int_0^\infty x^ke^{-x^2/2t}dx\\
&=&\sqrt{\frac{2}{\pi t}}\int_0^\infty(\sqrt{2ty})^ke^{-y}\,\sqrt{\frac{t}{2y}}\,dy\\
&=&\sqrt{\frac{(2t)^k}{\pi}}\int_0^\infty y^{(k-1)/2}e^{-y}dy\\
&=&\sqrt{\frac{(2t)^k}{\pi}}\cdot\Gamma\left(\frac{k+1}{2}\right)
\end{eqnarray*}

(5) Getting now to the computation of $\gamma,\kappa$, the low order moments are as follows:
$$M_1=\sqrt{\frac{2t}{\pi}}\quad,\quad M_2=t\quad,\quad 
M_3=\sqrt{\frac{(2t)^3}{\pi}}\quad,\quad M_4=3t^2$$

The low order central moments are $M_1'=0$, $M_2'=t-2t/\pi$, and then:
\begin{eqnarray*}
M_3'
&=&M_3-3EM_2+2E^3\\
&=&2t\sqrt{\frac{2t}{\pi}}-3t\sqrt{\frac{2t}{\pi}}+\frac{4t}{\pi}\sqrt{\frac{2t}{\pi}}\\
&=&\left(\frac{4t}{\pi}-t\right)\sqrt{\frac{2t}{\pi}}
\end{eqnarray*} 

As for the next central moment, the fourth one, this is given by:
\begin{eqnarray*}
M_4'
&=&M_4-4EM_3+6E^2M_2-3E^4\\
&=&3t^2-\frac{16t^2}{\pi}+\frac{12t^2}{\pi}-\frac{12t^2}{\pi^2}\\
&=&3t^2-\frac{4t^2}{\pi}-\frac{12t^2}{\pi^2}
\end{eqnarray*} 

Now by normalizing, we obtain the following formula for the skewness $\gamma$:
\begin{eqnarray*}
\gamma
&=&\left(\frac{4t}{\pi}-t\right)\sqrt{\frac{2t}{\pi}}\left(t-\frac{2t}{\pi}\right)^{-3/2}\\
&=&\left(\frac{4}{\pi}-1\right)\sqrt{\frac{2}{\pi}}\left(1-\frac{2}{\pi}\right)^{-3/2}\\
&=&\frac{4-\pi}{\pi}\sqrt{\frac{2}{\pi}}\cdot\frac{\pi^{3/2}}{(\pi-2)^{3/2}}\\
&=&(4-\pi)\sqrt{\frac{2}{(\pi-2)^3}}
\end{eqnarray*}

Similarly, by normalizing we obtain the following formula for the kurtosis $\kappa$:
\begin{eqnarray*}
\kappa
&=&\left(3t^2-\frac{4t^2}{\pi}-\frac{12t^2}{\pi^2}\right)\left(t-\frac{2t}{\pi}\right)^{-2}\\
&=&\left(3-\frac{4}{\pi}-\frac{12}{\pi^2}\right)\left(1-\frac{2}{\pi}\right)^{-2}\\
&=&\left(3-\frac{4}{\pi}-\frac{12}{\pi^2}\right)\frac{\pi^2}{(\pi-2)^2}\\
&=&\frac{3\pi^2-4\pi-12}{(\pi-2)^2}
\end{eqnarray*}

Thus, we are led to the conclusions in the statement.
\end{proof}

As our second statement now, still dealing with $f\sim g_t$, we have:

\begin{theorem}
Assuming $f\sim g_t$, the law $\chi_{1t}^2=law(f^2)$ has density
$$\chi^2_{1t}=\frac{1}{\sqrt{2\pi t}}\,x^{-1/2}e^{-x/2t}dx$$
on $[0,\infty)$, and the basic properties of this law $\chi_{1t}^2$ are as follows:
\begin{enumerate}
\item The mean is $E=t$, the variance is $V=2t^2$.

\item The moments are $M_k=t^k(2k)!!$.

\item The skewness is $\gamma=2\sqrt{2}$, the kurtosis is $\kappa=15$.

\item The Fourier transform is $F(x)=1/\sqrt{1-2ixt}$.
\end{enumerate}
\end{theorem}

\begin{proof}
In order to compute the density, assume $f\sim g_t$, and consider an arbitrary function $\varphi:[0,\infty)\to\mathbb R$. We have then the following computation:
\begin{eqnarray*}
E(\varphi(f^2))
&=&\sqrt{\frac{2}{\pi t}}\int_0^\infty\varphi(x^2)e^{-x^2/2t}\,dx\\
&=&\sqrt{\frac{2}{\pi t}}\int_0^\infty\varphi(y)e^{-y/2t}\,\frac{dy}{2\sqrt{y}}\\
&=&\frac{1}{\sqrt{2\pi t}}\int_0^\infty\varphi(y)y^{-1/2}e^{-y/2t}\,dy
\end{eqnarray*}

Thus the density of $f^2\sim\chi_{1t}^2$ is the one in the statement. As for the rest:

\medskip

(1) The moments being $M_k(f^2)=M_{2k}(f)=t^k(2k)!!$, the mean is $E=M_1=t$, and the variance is $V=M_2-M_1^2=3t^2-t^2=2t^2$, as stated.

\medskip

(2) As already mentioned, the moments are $M_k(f^2)=M_{2k}(f)=t^k(2k)!!$.

\medskip

(3) Getting now to the computation of $\gamma,\kappa$, the low order moments are as follows:
$$M_1=t\quad,\quad M_2=3t^2\quad,\quad 
M_3=15t^3\quad,\quad M_4=105t^4$$

The low order central moments are $M_1'=0$, $M_2'=2t^2$, and then:
\begin{eqnarray*}
M_3'
&=&M_3-3EM_2+2E^3\\
&=&15t^3-9t^3+2t^3\\
&=&8t^3
\end{eqnarray*} 

As for the next central moment, the fourth one, this is given by:
\begin{eqnarray*}
M_4'
&=&M_4-4EM_3+6E^2M_2-3E^4\\
&=&105t^4-60t^4+18t^4-3t^4\\
&=&60t^4
\end{eqnarray*} 

Now by normalizing, we obtain the formulae of $\gamma,\kappa$ in the statement.

\medskip

(4) The Fourier transform computation is straighforward, as follows:
\begin{eqnarray*}
F(x)
&=&\frac{1}{\sqrt{2\pi t}}\int_0^\infty y^{-1/2}e^{-y/2t+ixy}\,dy\\
&=&\frac{1}{\sqrt{2\pi t}}\int_0^\infty y^{-1/2}e^{(2ix-1/t)y/2}\,dy\\
&=&\frac{1}{\sqrt{2\pi t}}\int_0^\infty\left(\frac{z}{1/t-2ix}\right)^{-1/2}e^{-z/2}\,\frac{dz}{1/t-2ix}\\
&=&\frac{1}{\sqrt{1/t-2ix}}\cdot\frac{1}{\sqrt{2\pi t}}\int_0^\infty z^{-1/2}e^{-z/2}\,dz\\
&=&\frac{1}{\sqrt{1/t-2ix}}\cdot\frac{1}{\sqrt{2\pi t}}\cdot\sqrt{2\pi}\\
&=&\frac{1}{\sqrt{1-2ixt}}
\end{eqnarray*}

Thus, we are led to the conclusions in the statement.
\end{proof}

\section*{6b. Rayleigh laws}

Well, I don't know about you, but personally I find Theorems 6.3 and 6.4, dealing with the $\chi,\chi^2$ problematics for $f\sim g_t$, a bit harder than expected, and I don't really feel like getting into shifts $f\sim g_t^a$ right away. Remember indeed from chapter 2 that the moments of these shifts were not really computable, so when putting on top of this some gamma function combinatorics, things will most likely become truly complicated.

\bigskip

In short, giving up, or rather change of plan, my proposal would be to investigate next the $\chi,\chi^2$ problematics for $f\sim G_t$, and for shifts and further parameters that can be added, we can see later. In addition, a quick comparison between Theorems 6.3 and 6.4 shows that the $\chi^2$ problematics is simpler than the $\chi$ one, and so, as final updated plan, we will investigate the $\chi^2$ problematics, and then the $\chi$ problematics, for $f\sim G_t$.

\bigskip

Getting to work now, regarding the $\chi^2$ problematics for $f\sim G_t$, we have:

\begin{theorem}
Assuming $f\sim G_t$, the law $\chi_{2t}^2=law(2|f|^2)$ has density
$$\chi^2_{2t}=\frac{1}{2t}\,e^{-x/2t}dx$$
on $[0,\infty)$, and the basic properties of this law $\chi_{2t}^2$ are as follows:
\begin{enumerate}
\item The mean is $E=2t$, the variance is $V=4t^2$.

\item The moments are $M_k=(2t)^kk!$.

\item The skewness is $\gamma=2$, the kurtosis is $\kappa=9$.

\item The Fourier transform is $F(x)=1/(1-2ixt)$.
\end{enumerate}
\end{theorem}

\begin{proof}
This is something quite subtle, the idea being as follows:

\medskip

(1) To start with, the density in the statement is exponential, $\chi_{2t}^2=e_{1/2t}$. However, what the theorem says is not the usual story of the exponential law.

\medskip

(2) Next, according to our convention for the law $G_t$, we have $\chi_{2t}^2=law(2|f|^2)$, where $f=(a+ib)/\sqrt{2}$, with $a,b\sim g_t$ independent. Thus, $\chi_{2t}^2$ can be thought of as well as being something real, 2-dimensional, and with this justifying the subscript 2:
$$\chi_{2t}^2=law(a^2+b^2)$$

(3) Let us also mention that $\chi_{2t}^2$ is called as well squared Rayleigh law, and more on this later in this chapter, when discussing the Rayleigh law $\chi_{2t}=law(\sqrt{2}|f|)$. 

\medskip

(4) Getting now to the proof, we already computed the moments in chapter 5, and it is most convenient to start by recalling that. The computation was as follows:
\begin{eqnarray*}
M_k
&=&E\left((a^2+b^2)^k\right)\\
&=&\sum_{r=0}^k\binom{k}{r}E(a^{2r})E(b^{2k-2r})\\
&=&t^k\sum_{r=0}^k\binom{k}{r}(2r)!!(2k-2r)!!\\
&=&t^k\sum_{r=0}^k\frac{k!}{r!(k-r)!}\cdot\frac{(2r)!}{2^rr!}\cdot\frac{(2k-2r)!}{2^{k-r}(k-r)!}\\
&=&\frac{t^kk!}{2^k}\sum_{r=0}^k\binom{2r}{r}\binom{2k-2r}{k-r}
\end{eqnarray*}

On the other hand, the generalized binomial formula with exponent $-1/2$ gives:
$$\frac{1}{\sqrt{1+s}}=\sum_{q=0}^\infty\binom{2q}{q}\left(\frac{-s}{4}\right)^q$$

Now by taking the square of this series, we obtain the following formula:
$$\frac{1}{1+s}
=\sum_{k=0}^\infty\left(\frac{-s}{4}\right)^k\sum_{r=0}^k\binom{2r}{r}\binom{2k-2r}{k-r}$$

Thus the sum on the right is $4^k$, and we can finish our moment computation:
$$M_k=\frac{t^kk!}{2^k}\times 4^k=(2t)^kk!$$

(5) The mean is therefore $E=2t$, and the variance is $V=8t^2-4t^2=4t^2$.

\medskip

(6) Let us prove now that the density is the one in the statement, namely:
$$\chi_{2t}^2=\frac{1}{2t}\,e^{-x/2t}dx$$

By partial integration, the moments of this predicted density are given by:
\begin{eqnarray*}
I_k
&=&\frac{1}{2t}\int_0^\infty x^ke^{-x/2t}dx\\
&=&\frac{1}{2t}\int_0^\infty kx^{k-1}\cdot 2te^{-x/2t}dx\\
&=&2tk\cdot\frac{1}{2t}\int_0^\infty x^{k-1}e^{-x/2t}dx\\
&=&2tk\cdot I_{k-1}
\end{eqnarray*}

As for the initial value, for this recurrence, this is given by the following formula, which by the way shows that we have indeed a probability measure, as we should:
$$I_0
=\frac{1}{2t}\int_0^\infty e^{-x/2t}dx
=\frac{1}{2t}\Big[-2te^{-x/2t}\Big]_0^\infty
=1$$

Thus $I_k=(2t)^kk!$, which matches with $M_k(\chi_{2t}^2)=(2t)^kk!$, so job done.

\medskip

(7) The computation of the Fourier transform is very standard, as follows:
\begin{eqnarray*}
F(x)
&=&\frac{1}{2t}\int_0^\infty e^{-y/2t+ixy}\,dy\\
&=&\frac{1}{2t}\int_0^\infty e^{(2ix-1/t)y/2}\,dy\\
&=&\frac{1}{2t}\int_0^\infty e^{-z/2}\,\frac{dz}{1/t-2ix}\\
&=&\frac{1}{1-2ixt}\cdot\frac{1}{2}\int_0^\infty e^{-z/2}\,dz\\
&=&\frac{1}{1-2ixt}
\end{eqnarray*}

(8) According to our moment formula in (4), the low order moments are:
$$M_1=2t\quad,\quad 
M_2=8t^2\quad,\quad 
M_3=48t^3\quad,\quad 
M_4=384t^4$$

Thus the low order central moments are $M_1'=0$, $M_2'=4t^2$, and then:
$$M_3'=M_3-3EM_2+2E^3=16t^3$$
$$M_4'=M_4-4EM_3+6E^2M_2-3E^4=144t^4$$

Now by normalizing, we obtain from this $\gamma=2$ and $\kappa=9$, as stated. 
\end{proof}

Regarding now the $\chi$ laws, called Rayleigh laws, their basic theory is as follows:

\index{Rayleigh laws}

\begin{theorem}
Assuming $f\sim G_t$, the law $\chi_{2t}=law(\sqrt{2}|f|)$ has density
$$\chi_{2t}=\frac{1}{t}\,xe^{-x^2/2t}dx$$
on $[0,\infty)$, and the basic properties of this law $\chi_{2t}$ are as follows:
\begin{enumerate}
\item The mean is $E=\sqrt{t\pi/2}$, the variance is $V=2t-t\pi/2$.

\item The even moments are $M_{2l}=(2t)^ll!$

\item The odd moments are $M_{2l-1}=\sqrt{(2t)^{2l-1}\pi}\cdot(2l)!!/2^l$.

\item Equivalently, the moments are $M_k=\sqrt{(2t)^k}\cdot\Gamma(k/2+1)$.

\item We have $\gamma=2(\pi-3)\sqrt{\pi/(4-\pi)^3}$ and $\kappa=(32-3\pi^2)/(4-\pi)^2$.
\end{enumerate}
\end{theorem}

\begin{proof}
According to our convention for the law $G_t$, we have $\chi_{2t}=law(\sqrt{2}|f|)$, where $f=(a+ib)/\sqrt{2}$, with $a,b\sim g_t$ independent. Thus, the Rayleigh law $\chi_{2t}$ can be thought of as well as being something real, and 2-dimensional, and with this justifying the subscript 2, and more specifically, appearing as follows, with $a,b\sim g_t$ independent:
$$\chi_{2t}=law(\sqrt{a^2+b^2})$$

Let us also mention that at the level of the phenomenology, the Rayleigh law appears in many interesting situations, as for instance in the modeling of wind. And with exercise for you to learn more about this. As for the proof, this is routine, as follows:

\medskip

(1) In order to compute the density, the simplest is to use Theorem 6.5, along with the following fact, coming from the definition of our chi and chi squared laws:
$$g\geq0\quad,\quad g\sim\chi_{2t}\iff g^2\sim\chi_{2t}^2$$

Indeed, assume $g\sim\chi_{2t}$, and consider an arbitrary function $\varphi:[0,\infty)\to\mathbb R$. If we define $\psi:[0,\infty)\to\mathbb R$ by $\psi(x)=\varphi(\sqrt{x})$, we have the following computation:
\begin{eqnarray*}
E(\varphi(g))
&=&E(\psi(g^2))\\
&=&\frac{1}{2t}\int_0^\infty\psi(x)e^{-x/2t}\,dx\\
&=&\frac{1}{2t}\int_0^\infty\varphi(y)e^{-y^2/2t}\,2ydy\\
&=&\frac{1}{t}\int_0^\infty\varphi(y)ye^{-y^2/2t}\,dy
\end{eqnarray*}

Thus the density of $g\sim\chi_{2t}$ is the one in the statement. As for the other assertions:

\medskip

(2) Regarding the moment computation, this is straightforward, as follows:
\begin{eqnarray*}
M_k
&=&\frac{1}{t}\int_0^\infty x^{k+1}e^{-x^2/2t}\,dx\\
&=&\frac{1}{t}\int_0^\infty\left(\sqrt{2ty}\right)^{k+1}e^{-y}\,{\sqrt\frac{t}{2y}}\,dy\\
&=&\sqrt{(2t)^k}\int_0^\infty y^{k/2}e^{-y}\,dy\\
&=&\sqrt{(2t)^k}\cdot\Gamma\left(\frac{k}{2}+1\right)
\end{eqnarray*}

In particular the even moments are given by the following formula:
$$M_{2l}=(2t)^l\Gamma(l+1)=(2t)^ll!$$

As for the odd moments, these are given by the following formula:
$$M_{2l-1}=\sqrt{(2t)^{2l-1}}\cdot\Gamma\left(l+\frac{1}{2}\right)
=\sqrt{(2t)^{2l-1}\pi}\cdot\frac{(2l)!!}{2^l}$$

(3) Getting now to the computation of the parameters $(E,V,\gamma,\kappa)$, by using the moment formulae found above, the first four moments of our law are given by:
$$M_1=\sqrt{(2t)\pi}\cdot\frac{1}{2}=\sqrt{\frac{t\pi}{2}}\quad,\quad M_2=2t\cdot1=2t$$
$$M_3=\sqrt{(2t)^3\pi}\cdot\frac{3}{4}=3t\sqrt{\frac{t\pi}{2}}\quad,\quad M_4=4t^2\cdot2=8t^2$$

We can see that the mean is $E=\sqrt{t\pi/2}$, and the variance is $V=2t-t\pi/2$. Next, the third central moment is given by the following formula:
\begin{eqnarray*}
M_3'
&=&M_3-3EM_2+2E^3\\
&=&3t\sqrt{\frac{t\pi}{2}}-6t\sqrt{\frac{t\pi}{2}}+t\pi\sqrt{\frac{t\pi}{2}}\\
&=&t(\pi-3)\sqrt{\frac{t\pi}{2}}
\end{eqnarray*} 

By normalizing, we conclude that the skewness is given by the following formula:
\begin{eqnarray*}
\gamma
&=&t(\pi-3)\sqrt{\frac{t\pi}{2}}\cdot\frac{1}{\sqrt{(2t-t\pi/2)^3}}\\
&=&(\pi-3)\sqrt{\frac{\pi}{2}}\cdot\frac{1}{\sqrt{(2-\pi/2)^3}}\\
&=&2(\pi-3)\sqrt{\frac{\pi}{(4-\pi)^3}}
\end{eqnarray*}

Next, the fourth central moment is given by the following formula:
\begin{eqnarray*}
M_4'
&=&M_4-4EM_3+6E^2M_2-3E^4\\
&=&8t^2-12t\cdot\frac{t\pi}{2}+6\cdot\frac{t\pi}{2}\cdot 2t-3\,\frac{t^2\pi^2}{4}\\
&=&t^2\left(8-\frac{3\pi^2}{4}\right)
\end{eqnarray*} 

By normalizing, we conclude that the kurtosis is given by the following formula:
\begin{eqnarray*}
\kappa
&=&t^2\left(8-\frac{3\pi^2}{4}\right)\frac{1}{(2t-t\pi/2)^2}\\
&=&\left(8-\frac{3\pi^2}{4}\right)\frac{1}{(2-\pi/2)^2}\\
&=&\frac{32-3\pi^2}{(4-\pi)^2}
\end{eqnarray*}

Thus, we are led to the conclusions in the statement.
\end{proof}

\section*{6c. Parametric sums}

Good work that we did, but with my apologies, we are behind schedule, with $\chi,\chi^2$ done for $g_t,G_t$, but with $g_t^a,G_t^a,G_{st},G_{st}^a$ plus $p^s_t$ still in need to be investigated. Nevermind. So we will speed up a bit, and focus on the essentials, namely density, moments and $E,V,\gamma,\kappa$. Regarding the simplest question, namely $\chi^2$ for $g_t^a$, the result is:

\begin{theorem}
Assuming $f\sim g_t^a$ with $a\geq0$, the law $\chi_{1t}^{2a}=law(f^2)$ has density
$$\chi^{2a}_{1t}=\frac{e^{-a^2/2t}}{\sqrt{2\pi t}}\cdot\frac{e^{-x/2t}}{\sqrt{x}}\cosh\left(\frac{a\sqrt{x}}{t}\right)dx$$
on $[0,\infty)$, and the basic properties of this law $\chi_{1t}^{2a}$ are as follows:
\begin{enumerate}
\item The mean is $E=t+a^2$, the variance is $V=2t^2+4a^2t$.

\item The moments are $M_k=M_{2k}(g_t^a)$, computable by recurrence.

\item $\gamma=(t+3a^2)\sqrt{8t/(t+2a^2)^3}$ and $\kappa=(15t^2+60a^2t+12a^4)/(t+2a^2)^2$.
\end{enumerate}
\end{theorem}

\begin{proof}
Assuming $f\sim g_t^a$ with $a\geq0$, we have the following computation:
\begin{eqnarray*}
E(\varphi(f^2))
&=&\frac{1}{\sqrt{2\pi t}}\int_\mathbb R\varphi(x^2)e^{-(x-a)^2/2t}dx\\
&=&\frac{1}{\sqrt{2\pi t}}\left(\int_{-\infty}^0\varphi(x^2)e^{-(x-a)^2/2t}dx
+\int_0^\infty\varphi(x^2)e^{-(x-a)^2/2t}dx\right)\\
&=&\frac{1}{\sqrt{2\pi t}}\left(\int_0^\infty\varphi(y)e^{-(-\sqrt{y}-a)^2/2t}\frac{dy}{2\sqrt{y}}
+\int_0^\infty\varphi(y)e^{-(\sqrt{y}-a)^2/2t}\frac{dy}{2\sqrt{y}}\right)\\
&=&\frac{1}{\sqrt{2\pi t}}\int_0^\infty\frac{\varphi(y)}{\sqrt{y}}
\cdot\frac{1}{2}\left(e^{-(\sqrt{y}+a)^2/2t}+e^{-(\sqrt{y}-a)^2/2t}\right)dy\\
&=&\frac{e^{-a^2/2t}}{\sqrt{2\pi t}}\int_0^\infty\frac{\varphi(y)}{\sqrt{y}}\cdot e^{-y/2t}
\cdot\frac{1}{2}\left(e^{-a\sqrt{y}/t}+e^{a\sqrt{y}/t}\right)dy\\
&=&\frac{e^{-a^2/2t}}{\sqrt{2\pi t}}\int_0^\infty\frac{\varphi(y)}{\sqrt{y}}\cdot e^{-y/2t}
\cdot\cosh\left(\frac{a\sqrt{y}}{t}\right)dy
\end{eqnarray*}

Thus, the density is the one in the statement. Regarding now the moments, these are the numbers $M_k=M_{2k}(g_t^a)$, and as explained in chapter 2, the whole sequence $M_k(g_t^a)$ is computable by recurrence. In particular, as explained there, we have:
$$M_2(g_t^a)=t+a^2\quad,\quad M_4(g_t^a)=3t^2+6a^2t+a^4$$
$$M_6(g_t^a)=15t^3+45a^2t^2+15a^4t+a^6$$
$$M_8(g_t^a)=105t^4+420a^2t^3+210a^4t^2+28a^6t+a^8$$

Thus we have $E=t+a^2$, and the variance is given by the following formula:
$$V=3t^2+6a^2t+a^4-(t+a^2)^2=2t^2+4a^2t$$

Next, the third and fourth central moments of our law are given by:
$$M_3'=M_3-3EM_2+2E^3=8t^2(t+3a^2)$$
$$M_4'=M_4-4EM_3+6E^2M_2-3E^4=12t^2(5t^2+20a^2t+4a^4)$$

Thus, we are led to the formulae of $\gamma$ and $\kappa$ in the statement.
\end{proof}

Still talking shifts of the real normal laws, we have as well the following result:

\begin{theorem}
Assuming $f\sim g_t^a$ with $a\geq0$, the law $\chi_{1t}^a=law(|f|)$ has density
$$\chi_{1t}^a=\sqrt{\frac{2}{\pi t}}\,e^{-a^2/2t}\cdot e^{-x^2/2t}\cosh\left(\frac{ax}{t}\right)dx$$
on $[0,\infty)$, and the basic properties of this law $\chi_{1t}^a$ are as follows:
\begin{enumerate}
\item The even moments are $M_{2l}=M_{2l}(g_t^a)$, computable by recurrence. 

\item The odd moments are given by $M_{2l+1}=\sqrt{(2t)^{2l+1}\pi}\cdot l!$ at $a=0$.
\end{enumerate}
\end{theorem}

\begin{proof}
Assuming $f\sim g_t^a$ with $a\geq0$, we have the following computation:
\begin{eqnarray*}
E(\varphi(|f|))
&=&\frac{1}{\sqrt{2\pi t}}\int_\mathbb R\varphi(|x|)e^{-(x-a)^2/2t}dx\\
&=&\frac{1}{\sqrt{2\pi t}}\left(\int_{-\infty}^0\varphi(-x)e^{-(x-a)^2/2t}dx
+\int_0^\infty\varphi(x)e^{-(x-a)^2/2t}dx\right)\\
&=&\frac{1}{\sqrt{2\pi t}}\left(\int_0^\infty\varphi(x)e^{-(-x-a)^2/2t}dx
+\int_0^\infty\varphi(x)e^{-(x-a)^2/2t}dx\right)\\
&=&\frac{1}{\sqrt{2\pi t}}\int_0^\infty\varphi(x)\left(e^{-(x+a)^2/2t}+e^{-(x-a)^2/2t}\right)dx\\
&=&\sqrt{\frac{2}{\pi t}}\,e^{-a^2/2t}\int_0^\infty\varphi(x)\cdot e^{-x^2/2t}
\cdot\frac{1}{2}\left(e^{-ax/t}+e^{ax/t}\right)dy\\
&=&\sqrt{\frac{2}{\pi t}}\,e^{-a^2/2t}\int_0^\infty\varphi(x)\cdot e^{-x^2/2t}\cosh\left(\frac{ax}{t}\right)dx
\end{eqnarray*}

Thus, the density is the one in the statement, and with the remark that at $a=0$ we obtain indeed the density in Theorem 6.3. Regarding now the moments, these are:
$$M_k=\sqrt{\frac{2}{\pi t}}\,e^{-a^2/2t}\int_0^\infty x^ke^{-x^2/2t}\cosh\left(\frac{ax}{t}\right)dx$$

And with the above integral being not very inviting, we will stop here, and formulate instead about moments the two things in the statement, that we both know well.
\end{proof}

Quite interesting all this, we just talked about random variables whose expectation $E$ we are unable to compute, guess this means that we are at the advanced level. So very nice, and talking advanced level, perhaps time to ask the cat. Which cat declares:

\begin{cat}
Hypergeometric functions are as fun as quantum groups, the first human to master them both will be awarded an honorary cat title. 
\end{cat}

Humm, interesting what you say, katie. And I guess that's good advice for you reader, too, shall you ever manage to master one of these main virtual objects in mathematics, such as hypergeometric functions, or quantum groups, don't stop there. Being advanced with respect to most humans means nothing, compared to an honorary cat title.

\bigskip

Getting back now to work, time to upgrade to two dimensions, and investigate the $\chi,\chi^2$ problematics for the more specialized normal laws there, namely $G_t^a,G_{st},G_{st}^a$. And with myself lacking I guess the needed titles and diplomas, for making something fun out of shifts, the plan will be investigate first $G_{st}$, and then, hopefully, $G_t^a$, and even $G_{st}^a$.

\bigskip

So, let us start with a discussion regarding $G_{st}$. We already met this law in the context of the various CCLT from chapter 5, but as a matter of having a new start, for all this, here is what is to be known about $G_{st}$, and in fact about $G_t,G_t^a,G_{st}^a$ as well:

\begin{theorem}
The various normal laws $G=law((h+ik)/\sqrt{2})$ are as follows:
\begin{enumerate}
\item The usual normal law $G_t$, coming from $h,k\sim g_t$ independent, and its shift at $a=(b+ic)/\sqrt{2}$, coming from $h\sim g_t^b,k\sim g_t^c$ independent, are
$$G_t=\frac{1}{\pi t}\,e^{-|z|^2/t}dz\quad\ ,\ \quad G_t^a=\frac{1}{\pi t}\,e^{-|z-a|^2/t}dz$$

\item The two-parameter normal law $G_{st}$, coming from $h\sim g_s,k\sim g_t$ independent, and its shift at $a=(b+ic)/\sqrt{2}$, coming from $h\sim g_s^b,k\sim g_t^c$ independent, are
$$G_{st}=\frac{1}{\pi\sqrt{st}}\,e^{-x^2/s-y^2/t}dxdy\quad,\quad G_{st}^a=\frac{1}{\pi\sqrt{st}}e^{-(x-a_x)^2/s-(y-a_y)^2/t}dxdy$$
\end{enumerate}
with the usual convention $z=x+iy$ for the numbers in the complex plane.
\end{theorem}

\begin{proof}
This is something very standard, based on the real formulae, namely:
$$g_t=\frac{1}{\sqrt{2\pi t}}\,e^{-x^2/2t}dx\quad,\quad 
g_t^a=\frac{1}{\sqrt{2\pi t}}\,e^{-(x-a)^2/2t}dx$$

(1) Our first claim is that we have the following dilation formula, over the reals:
$$h\sim\phi(x)dx\ \implies\ \lambda h\sim\frac{1}{\lambda}\phi\left(\frac{x}{\lambda}\right)dx$$

Indeed, with $h\sim\phi(x)dx$, we have the following computation, proving our claim:
$$E(\psi(\lambda h))
=\int_\mathbb R\psi(\lambda x)\phi(x)dx
=\int_\mathbb R\psi(y)\phi\left(\frac{y}{\lambda}\right)\frac{dy}{\lambda}$$

(2) Assume now that we have independent variables $h,k\sim g_t$. According to our dilation formula above, applied with $\lambda=1/\sqrt{2}$, we have:
$$\frac{h}{\sqrt{2}}\sim\frac{1}{\sqrt{\pi t}}\,e^{-x^2/t}dx
\quad,\quad \frac{k}{\sqrt{2}}\sim\frac{1}{\sqrt{\pi t}}\,e^{-y^2/t}dy$$

We conclude that, with $z=x+iy$ as usual, the law of $(h+ik)/\sqrt{2}$ is:
$$G_t=\frac{1}{\pi t}\,e^{-(x^2+y^2)/t}dxdy=\frac{1}{\pi t}\,e^{-|z|^2/t}dz$$

(3) Next, more generally, assume that we have independent variables $h\sim g_t^b,k\sim g_t^c$. According to our dilation formula above, again applied with $\lambda=1/\sqrt{2}$, we have:
$$\frac{h}{\sqrt{2}}\sim\frac{1}{\sqrt{\pi t}}\,e^{-(\sqrt{2}x-b)^2/2t}dx
\quad,\quad \frac{k}{\sqrt{2}}\sim\frac{1}{\sqrt{\pi t}}\,e^{-(\sqrt{2}y-c)^2/2t}dy$$

Thus, with $z=x+iy$ as usual, and $a=(b+ic)/\sqrt{2}$, the law of $(h+ik)/\sqrt{2}$ is:
\begin{eqnarray*}
G_t^a
&=&\frac{1}{\pi t}\,\exp\left(-\frac{1}{2t}\left[(\sqrt{2}x-b)^2+(\sqrt{2}y-c)^2\right]\right)dxdy\\
&=&\frac{1}{\pi t}\,\exp\left(-\frac{1}{t}\left[(x-b/\sqrt{2})^2+(y-c/\sqrt{2})^2\right]\right)dxdy\\
&=&\frac{1}{\pi t}\,\exp\left(-\frac{1}{t}\Big|\Big|(x,y)-\frac{(b,c)}{\sqrt{2}}\Big|\Big|^2\right)dxdy\\
&=&\frac{1}{\pi t}\,e^{-|z-a|^2/t}dz
\end{eqnarray*}

(4) Next, as another generalization of (2), assume that we have independent variables $h\sim g_s,k\sim g_t$. According to our dilation formula, applied with $\lambda=1/\sqrt{2}$, we have:
$$\frac{h}{\sqrt{2}}\sim\frac{1}{\sqrt{\pi s}}\,e^{-x^2/2s}dx
\quad,\quad \frac{k}{\sqrt{2}}\sim\frac{1}{\sqrt{\pi t}}\,e^{-y^2/2t}dy$$

We conclude that the law of $(h+ik)/\sqrt{2}$ is given by the following formula:
$$G_t=\frac{1}{\pi\sqrt{st}}\,e^{-x^2/2s-y^2/2t}dxdy$$

(5) Finally, with full house parameters, meaning independent variables $h\sim g_s^b,k\sim g_t^c$, according to our dilation formula, applied as usual with $\lambda=1/\sqrt{2}$, we have:
$$\frac{h}{\sqrt{2}}\sim\frac{1}{\sqrt{\pi s}}\,e^{-(\sqrt{2}x-b)^2/2s}dx
\quad,\quad \frac{k}{\sqrt{2}}\sim\frac{1}{\sqrt{\pi t}}\,e^{-(\sqrt{2}y-c)^2/2t}dy$$

Thus, with $z=x+iy$ as usual, and $a=(b+ic)/\sqrt{2}$, the law of $(h+ik)/\sqrt{2}$ is:
\begin{eqnarray*}
G_{st}^a
&=&\frac{1}{\pi\sqrt{st}}\,\exp\left(-\frac{(\sqrt{2}x-b)^2}{2s}-\frac{(\sqrt{2}y-c)^2}{2t}\right)dxdy\\
&=&\frac{1}{\pi\sqrt{st}}\,\exp\left(-\frac{(x-b/\sqrt{2})^2}{s}-\frac{(y-c/\sqrt{2})^2}{t}\right)dxdy\\
&=&\frac{1}{\pi\sqrt{st}}e^{-(x-a_x)^2/s-(y-a_y)^2/t}dxdy
\end{eqnarray*}

Thus, we are led to the conclusions in the statement.
\end{proof}

According to our updated plan made above, let us first investigate $G_{st}$, given by:
$$G_{st}=\frac{1}{\pi\sqrt{st}}\,e^{-x^2/s-y^2/t}dxdy$$

However, when looking at moments, we are right away in trouble, with respect to the case $s=t$, previously investigated in this chapter, and in fact going back to chapter 5. Indeed, when $s=t$ our distribution is rotationally invariant, and this allowed us to say that all moments of a variable $f$ following this law vanish, except for those of $|f|^2$.

\bigskip

In the case $s\neq t$ the rotational invariance obviously fails, and so fails the vanishing of the generic moments of $f$. So, what to do? I would say that the wisest would be to restrict the attention to $|f|^2$, that is, to study the $\chi^2$ problematics for $G_{st}$, which after all was on our official plan. And, we are led in this way to the following result:

\begin{theorem}
Assuming $f\sim G_{st}$, the law $\chi_{2st}^2=law(2|f|^2)$ has moments
$$M_k=\frac{k!}{2^k}\sum_{k=a+b}\binom{2a}{a}\binom{2b}{b}s^at^b$$
the Fourier transform is given by the following formula,
$$F(x)=\frac{1}{\sqrt{(1-2six)(1-2tix)}}$$
the expectation is $E=s+t$, the variance is $V=2s^2+2t^2$, and
$$\gamma=\frac{s^3+t^3}{s^2+t^2}\sqrt{\frac{8}{s^2+t^2}}\quad,\quad
\kappa=\frac{15s^4+15t^4+6s^2t^2}{(s^2+t^2)^2}$$
are the skewness and kurtosis.
\end{theorem}

\begin{proof}
This is similar to the proof of Theorem 6.5, the idea being as follows:

\medskip

(1) With $f=(d+ie)/\sqrt{2}$, the moments are given by the following formula:
\begin{eqnarray*}
M_k
&=&E\left((d^2+e^2)^k\right)\\
&=&\sum_{k=a+b}\binom{k}{a}E(d^{2a})E(e^{2b})\\
&=&\sum_{k=a+b}\binom{k}{a}s^a(2a)!!t^b(2b)!!\\
&=&\sum_{k=a+b}\frac{k!}{a!b!}\cdot\frac{(2a)!}{2^aa!}\cdot\frac{(2b)!}{2^bb!}\,s^at^b\\
&=&\frac{k!}{2^k}\sum_{k=a+b}\binom{2a}{a}\binom{2b}{b}s^at^b
\end{eqnarray*}

(2) The generalized binomial formula with exponent $-1/2$ shows that we have:
$$\frac{1}{\sqrt{1-4sz}}=\sum_{a=0}^\infty\binom{2a}{a}(sz)^a
\quad,\quad\frac{1}{\sqrt{1-4tz}}=\sum_{b=0}^\infty\binom{2b}{b}(tz)^b$$

By multiplying these two formulae, and then setting $k=a+b$, we obtain:
\begin{eqnarray*}
\frac{1}{\sqrt{(1-4sz)(1-4tz)}}
&=&\sum_{a=0}^\infty\sum_{b=0}^\infty\binom{2a}{a}\binom{2b}{b}(sz)^a(tz)^b\\
&=&\sum_{k=0}^\infty z^k\sum_{k=a+b}\binom{2a}{a}\binom{2b}{b}s^at^b\\
&=&\sum_{k=0}^\infty z^k\cdot\frac{2^kM_k}{k!}\\
&=&\sum_{k=0}^\infty\frac{(2z)^k}{k!}\cdot M_k
\end{eqnarray*}

Now by setting $z=ix/2$, this latter formula takes the following form:
$$\frac{1}{\sqrt{(1-2six)(1-2tix)}}=\sum_{k=0}^\infty\frac{(ix)^k}{k!}\cdot M_k$$

But on the right we have the Fourier transform $F(x)$, so done with this too.

\medskip

(3) It remains to compute $E,V,\gamma,\kappa$. According to (1), the low order moments are:
$$M_1=s+t$$
$$M_2=3s^2+3t^2+2st$$
$$M_3=15s^3+9s^2t+9st^2+15t^3$$
$$M_4=105s^4+60s^3t+54s^2t^2+60st^3+105t^4$$

Thus the mean is $E=s+t$, and the variance is given by the following formula:
$$V=(3s^2+3t^2+2st)-(s+t)^2=2s^2+2t^2$$

(4) Next, the third central moment of our law is given by the following formula:
\begin{eqnarray*}
M_3'
&=&M_3-3EM_2+2E^3\\
&=&15s^3+9s^2t+9st^2+15t^3\\
&&-3(s+t)(3s^2+3t^2+2st)+2(s+t)^3\\
&=&8s^3+8t^3
\end{eqnarray*} 

By normalizing, we conclude that the skewness is given by the following formula:
$$\gamma=\frac{8s^3+8t^3}{(2s^2+2t^2)\sqrt{2s^2+2t^2}}
=\frac{s^3+t^3}{s^2+t^2}\sqrt{\frac{8}{s^2+t^2}}$$

(5) Next, the fourth central moment is given by the following formula:
\begin{eqnarray*}
M_4'
&=&M_4-4EM_3+6E^2M_2-3E^4\\
&=&105s^4+60s^3t+54s^2t^2+60st^3+105t^4\\
&&-4(s+t)(15s^3+9s^2t+9st^2+15t^3)\\
&&+6(s+t)^2(3s^2+3t^2+2st)-3(s+t)^4\\
&=&60s^4+60t^4+24s^2t^2
\end{eqnarray*} 

By normalizing, we conclude that the kurtosis is given by the following formula:
$$\kappa=\frac{60s^4+60t^4+24s^2t^2}{(2s^2+2t^2)^2}=\frac{15s^4+15t^4+6s^2t^2}{(s^2+t^2)^2}$$

Thus, we are led to the conclusions in the statement.
\end{proof}

Getting now to the $\chi$ problematics for $G_{st}$, things here are quite complicated, again due to the lack of rotational invariance, when $s\neq t$. With bare hands we can only formulate a modest remark, based on the density formula for $G_{st}$ from Theorem 6.10, and on our previous moment formulae from Theorems 6.6 and 6.11, as follows:

\begin{remark}
Assuming $f\sim G_{st}$, the law $\chi_{2st}=law(\sqrt{2}|f|)$ has moments:
$$M_k=\frac{1}{\pi\sqrt{st}}\int_\mathbb R\int_\mathbb R\sqrt{(x^2+y^2)^k}\,e^{-x^2/s-y^2/t}dxdy$$
The moments of even order are given by a quite explicit formula, namely:
$$M_{2l}=\frac{(2l)!}{4^l}\sum_{2k=a+b}\binom{2a}{a}\binom{2b}{b}s^at^b$$
Also, at $s=t$ the moment formula is $M_k=\sqrt{(2t)^k}\cdot\Gamma(k/2+1)$.
\end{remark}

And with this, end of our study. Good work that we did, and time I guess to stop here, and enjoy a good meal, I mean I don't know about you, but I find that these computations burn calories, I definitely need some. Of course we still have the shifts left, and the discrete measures too, but these can wait until tomorrow. Curtain.

\section*{6d. Poisson sums}

Welcome back, and in the hope that we agree on this, time to discuss now the discrete case, looking for some mathematical fun there. Our plan will be as follows:

\begin{plan}
We would like to investigate the $\chi,\chi^2$ problematics for:
\begin{enumerate}
\item Compound Poisson laws $p_{t\nu}$, with $\nu$ discrete probability measure on $\mathbb C$.

\item Compound Poisson laws $p_{t\nu}$, with $\nu$ discrete probability measure on $\mathbb T$.

\item Compound Poisson laws $p_{t\nu}$, with $\nu$ discrete uniform measure on $X\subset\mathbb T$.

\item Compound Poisson laws $p_{t\nu}$, with $\nu=\varepsilon_s$, uniform measure on $\mathbb Z_s\subset\mathbb T$.
\end{enumerate}
\end{plan}

In short, we plan to do things slowly, following the above hierarchy, starting with (1) which is the general case, and ending with (4) which is the Bessel law case, $p_{t\nu}=p^s_t$. So, going now with (1), here is what we need to know, extracted from chapter 5:

\begin{theorem}
Given a complex probability measure $\nu$, and $t>0$, let us set:
$$p_{t\nu}=\lim_{n\to\infty}\left(\left(1-\frac{t}{n}\right)\delta_0+\frac{t}{n}\,\nu\right)^{*n}$$
In the discrete case, $\nu=\sum_{k=1}^sc_k\delta_{w_k}$ with $c_k>0$ and $w_k\in\mathbb C$, the Fourier transform is
$$F_{p_{t\nu}}(x)=\exp\left(t\sum_{k=1}^sc_k(e^{ixw_k}-1)\right)$$
and the measure itself is given by the following formula, with $f_k\sim p_{tc_k}$, independent:
$$p_\nu=law\left(\sum_{k=1}^sw_kf_k\right)$$
As examples, $p_t^s=p_{t\varepsilon_s}$, with $\varepsilon_s$ being the uniform measure on the $s$-th roots of unity.
\end{theorem}

\begin{proof}
This is indeed something that we know well from chapter 5.
\end{proof}

At the above level of generality, we have the following result:

\begin{proposition}
For a measure $p_{t\nu}$ with $\nu=\sum_{k=1}^sc_k\delta_{w_k}$, we have:
$$E=t\sum_{k=1}^sc_kw_k\quad,\quad V=t\sum_{k=1}^sc_k|w_k|^2$$
When $\nu$ is supported on the unit circle $\mathbb T$, the variance formula reads $V=t$.
\end{proposition}

\begin{proof}
We know that $p_{t\nu}$ is the law of the variable $f=\sum_kw_kf_k$, with $f_k\sim p_{tc_k}$, independent. But with this, the formula of the expectation is clear, and we have:
\begin{eqnarray*}
E(|f|^2)
&=&\sum_{kl}w_k\bar{w}_lE(f_kf_l)\\
&=&\sum_kw_k\bar{w}_kE(f_k^2)+\sum_{k\neq l}w_k\bar{w}_lE(f_k)E(f_l)\\
&=&\sum_kw_k\bar{w}_k(tc_k+t^2c_k^2)+\sum_{k\neq l}w_k\bar{w}_l(t^2c_kc_l)\\
&=&t\sum_kw_k\bar{w}_kc_k+t^2\sum_{kl}w_k\bar{w}_lc_kc_l\\
&=&t\sum_k|w_k|^2c_k+t^2\Big|\sum_kw_kc_k\Big|^2
\end{eqnarray*}

Thus, we are led to the conclusions in the statement.
\end{proof}

The above result, with the nice formula $V=t$ at the end, looks quite encouraging. However, for higher moments things become quite complicated, as shown by:

\begin{remark}
The second moment of the variable $|f|^2=|\sum_kw_kf_k|^2$ is
$$E(|f|^4)=\sum_{klmn}w_kw_l\bar{w}_m\bar{w}_nE(f_kf_lf_mf_n)$$
which will spilt into 15 sums, according to the values of $\ker(klmn)\in P(4)$.
\end{remark}

In short, quite tough combinatorics going on here, and with respect to our original Plan 6.13, done I guess with (1) and (2) there, with our hunting trophy being something quite modest, namely Proposition 6.15, computing the mean of $|f|^2$. Guess for a fresh dinner tonight, we will have to stop by the market, or make a deal with the cat.

\bigskip

Getting now to (3) in our Plan 6.13, let us further assume that $\nu$ is the uniform measure on a certain finite subset $X\subset\mathbb T$. In practice, this means $\nu=\frac{1}{s}\sum_{k=1}^s\delta_{w_k}$ with $w_k\in\mathbb T$, so our variable is $f=\sum_kw_kf_k$, with $f_k\sim p_{t/s}$, independent, and we get:
\begin{eqnarray*}
E(|f|^4)
&=&\sum_{klmn}\frac{w_kw_l}{w_mw_n}\,E(f_kf_lf_mf_n)\\
&=&\sum_{\pi\in P(4)}\sum_{\ker{(klmn)}=\pi}\frac{w_kw_l}{w_mw_n}\,E(f_kf_lf_mf_n)\\
&=&\sum_{\pi\in P(4)}M_\pi(p_{t/s})\sum_{\ker{(klmn)}=\pi}\frac{w_kw_l}{w_mw_n}
\end{eqnarray*}

To be more precise, what we called here $M_\pi$ means product of moments, computed with respect to the blocks of $\pi$. More generally, we have the following result:

\begin{proposition}
For a measure $p_{t\nu}$ with $\nu=\frac{1}{s}\sum_{k=1}^s\delta_{w_k}$ with $w_k\in\mathbb T$, we have
$$E(|f|^{2r})=\sum_{\pi\in P(2r)}M_\pi(p_{t/s})\sum_{\ker{(km)}=\pi}\frac{w_{k_1}\ldots w_{k_r}}{w_{m_1}\ldots w_{m_r}}$$
with $M_\pi$ meaning product of moments, with respect to the blocks of $\pi$.
\end{proposition}

\begin{proof}
This follows by redoing the computation above, with an arbitrary exponent:
\begin{eqnarray*}
E(|f|^{2r})
&=&\sum_{km}\frac{w_{k_1}\ldots w_{k_r}}{w_{m_1}\ldots w_{m_r}}\,E(f_{k_1}\ldots f_{k_r}f_{m_1}\ldots f_{m_r})\\
&=&\sum_{\pi\in P(2r)}\sum_{\ker{(km)}=\pi}\frac{w_{k_1}\ldots w_{k_r}}{w_{m_1}\ldots w_{m_r}}\,E(f_{k_1}\ldots f_{k_r}f_{m_1}\ldots f_{m_r})\\
&=&\sum_{\pi\in P(2r)}M_\pi(p_{t/s})\sum_{\ker{(km)}=\pi}\frac{w_{k_1}\ldots w_{k_r}}{w_{m_1}\ldots w_{m_r}}
\end{eqnarray*}

Thus, we are led to the formula in the statement.
\end{proof}

And with this, done I guess with (3) in Plan 6.13 too. Cats are laughing at me, but let me nevertheless formulate a comment on all this, that I was keeping secret:

\begin{comment}
There is something conceptual to be said in general, namely
$$k_n(p_\alpha)=M_n(\alpha)$$
the cumulants of a compound Poisson law are the moments of the input measure.
\end{comment}

And more on such things later, in Part IV, when discussing cumulants. Moving on now, (4) in Plan 6.13, that is about the Bessel laws from chapter 5. So, let us recall from there that we have the following result, and in the hope that I did't forget anything:

\begin{theorem}
The Bessel laws $p_t^s=p_{t\varepsilon_s}$ are given by the following formula, with $f_1,\ldots,f_s$ being Poisson $(t/s)$ and independent, and $w=e^{2\pi i/s}$:
$$p^s_t=law\left(\sum_{k=1}^sw^kf_k\right)$$
The Fourier transform of $p^s_t$ is as follows, in terms of $\exp_sz=\sum_{n=0}^\infty z^{sn}/(sn)!$,
$$F(x)=\exp\left(t(\exp_s(ix)-1)\right)$$
the density is given by the following formula, with $\delta$ being a Dirac mass,
$$p^s_t=e^{-t}\sum_{p_1=0}^\infty\ldots\sum_{p_s=0}^\infty\frac{1}{p_1!\ldots p_s!}\,\left(\frac{t}{s}\right)^{p_1+\ldots+p_s}\delta\left(\sum_{k=1}^sw^kp_k\right)$$
and the moments are given by the following formula,
$$M_r(p^s_t)=\sum_{\pi\in P^s(r)}t^{|\pi|}$$
with $P^s$ standing for the partitions satisfying $\#\circ=\#\bullet(s)$ in each block.
\end{theorem}

\begin{proof}
This is something that we know well from chapter 5, and as illustrations for this, let us briefly discuss a few particular cases, which are of interest:

\medskip

(1) Case $s=1$. Here we have the usual Poisson law $p_t$, the Fourier transform is $F(x)=\exp(t(e^{ix}-1))$, the density is $p_t=e^{-t}\sum_{p\geq0}t^p\delta_p/p!$, and the partitions, which can be thought of as being uncolored, since we are in the real case, are $P^1=P$.

\medskip

(2) Case $s=2$. Here we have the real Bessel law $p_t^2$, the Fourier transform is $F(x)=\exp(t(\cos x-1))$, the density is $p_t^2=\sum_{k\in\mathbb Z}f_k(t/2)\delta_k$, with $f_k(t)=\sum_{p\geq0}t^{|k|+2p}/(|k|+p)!p!$ being the Bessel function of the first kind, and the partitions are $P^2=P_{even}$.

\medskip

(3) Case $s=3$. This is the first interesting case, with respect to what we knew well from Part I, regarding $p_t,p_t^2$, and obviously, all the formulae in the statement are about some crazy aspects of the combinatorics of the 3rd roots of unity. Good to know.

\medskip

(4) Case $s=4$. This is something quite manageable, because the root of unity is $w=i$, which brings us into more familiar territory. For instance, the density is supported on $\mathbb Z[i]$, which is nice. So, good to know, when stuck with $p^3_t$, try instead $p^4_t$.

\medskip

(5) Case $s=\infty$. Here we have the purely complex Bessel law $P_t$. And with the formulae of the Fourier transform and of the density being a bit unclear, the main tool here are the partitions, $P^\infty=\mathcal P_{even}$, the matching partitions with even blocks. See \cite{bb+}.
\end{proof}

With this discussed, time now to get into the $\chi,\chi^2$ problematics for the Bessel laws. In what regards the corresponding densities, things are quickly done, as follows:

\begin{theorem}
Given a variable following a Bessel law, $f\sim p^s_t$, we have
$$|f|\sim e^{-t}\sum_{p_1=0}^\infty\ldots\sum_{p_s=0}^\infty\frac{1}{p_1!\ldots p_s!}\,\left(\frac{t}{s}\right)^{p_1+\ldots+p_s}\delta\left(\,\left|\sum_{k=1}^sw^kp_k\right|\,\right)$$
$$|f|^2\sim e^{-t}\sum_{p_1=0}^\infty\ldots\sum_{p_s=0}^\infty\frac{1}{p_1!\ldots p_s!}\,\left(\frac{t}{s}\right)^{p_1+\ldots+p_s}\delta\left(\left|\sum_{k=1}^sw^kp_k\right|^2\right)$$
with the $\delta$ symbols standing as usual for Dirac masses.
\end{theorem}

\begin{proof}
This follows from the density formula from Theorem 6.19, namely:
$$f\sim e^{-t}\sum_{p_1=0}^\infty\ldots\sum_{p_s=0}^\infty\frac{1}{p_1!\ldots p_s!}\,\left(\frac{t}{s}\right)^{p_1+\ldots+p_s}\delta\left(\sum_{k=1}^sw^kp_k\right)$$

Indeed, with $\delta_z\to\delta_{|z|}$ and $\delta_z\to\delta_{|z|^2}$, we obtain the formulae in the statement.
\end{proof}

Getting now to more specialized questions, going beyond the density knowledge, things are quite complicated for $|f|$. However, in what regards $|f|^2$, we have good results:

\begin{theorem}
Given a variable following a Bessel law, $f\sim p^s_t$, we have
$$E(|f|^{2r})=\sum_{\pi\in P^s(\circ^r\bullet^r)}t^{|\pi|}$$
with $P^s$ standing for the partitions satisfying $\#\circ=\#\bullet(s)$ in each block. Also
$$E=t\quad,\quad V=t^2+t\quad,\quad\gamma=\frac{2t^2+6t+1}{(t+1)\sqrt{t^2+t}}
\quad,\quad \kappa=\frac{9t^3+42t^2+30t+1}{t(t+1)^2}$$
for the variable $|f|^2$, with these formulae being valid respectively at $s\geq 2,3,4,5$.
\end{theorem}

\begin{proof}
According to the general moment formula in Theorem 6.19, we have the following computation, with the dots standing everywhere for a total of $r$ copies:
\begin{eqnarray*}
E(|f|^{2r})
&=&E(f\ldots f\bar{f}\ldots\bar{f})\\
&=&E(f^{\circ\ldots\circ\bullet\ldots\bullet})\\
&=&\sum_{\pi\in P^s(\circ\ldots\circ\bullet\ldots\bullet)}t^{|\pi|}
\end{eqnarray*}

But with this, we can compute the low order moments of $|f|^2$, and then $(E,V,\gamma,\kappa)$:

\medskip

(1) The first moment $M_1$ comes from $P^s(\circ\,\bullet)$, which at $s\geq2$ consists precisely of the partition $\sqcap\,$, with the other 2-element partition, namely $|\ |\,$, being dismissed, since not satisfying the condition $\#\circ=\#\bullet(s)$ in each block. Thus at $s\geq2$ we have:
$$M_1=t$$

Similarly, the second moment $M_2$ comes from $P^s(\circ\circ\bullet\,\bullet)$, which at $s\geq3$ consists of the 2 double pairings which are matching, and the 4-block. Thus at $s\geq3$ we have:
$$M_2=2t^2+t$$

Next, the third moment $M_3$ comes from $P^s(\circ\circ\circ\bullet\bullet\,\bullet)$, which at $s\geq4$ consists of 6 partitions of type $2+2+2$, meaning the 6 possible matching triple pairings, then 9 partitions of type $4+2$, and finally 1 partition of type 6. Thus at $s\geq4$ we have:
$$M_3=6t^3+9t^2+t$$

Finally, $M_4$ comes from $P^s(\circ\circ\circ\circ\bullet\bullet\bullet\,\bullet)$, which at $s\geq5$ consists of 24 partitions of type $2+2+2+2$, then 72 partitions of type $4+2+2$, then 16 partitions of type $6+2$ and 18 partitions of type $4+4$, and finally 1 partition of type 8. Thus at $s\geq5$:
$$M_4=24t^4+72t^3+34t^2+t$$

(2) But with this, we can compute $(E,V,\gamma,\kappa)$, under suitable assumptions on $s$. Indeed, at $s\geq2$ the expectation is $E=t$, and then at $s\geq3$, the variance is:
$$V=(2t^2+t)-t^2=t^2+t$$

Next, assuming $s\geq4$, the third central moment can be computed as follows:
\begin{eqnarray*}
M_3'
&=&M_3-3EM_2+2E^3\\
&=&(6t^3+9t^2+t)-3t(2t^2+t)+2t^3\\
&=&2t^3+6t^2+t
\end{eqnarray*}

Now by normalizing, we obtain the following formula for the skewness, at $s\geq4$:
$$\gamma=\frac{2t^3+6t^2+t}{(t^2+t)\sqrt{t^2+t}}
=\frac{2t^2+6t+1}{(t+1)\sqrt{t^2+t}}$$

Finally, assuming $s\geq5$, the fourth central moment can be computed as follows:
\begin{eqnarray*}
M_4'
&=&M_4-4EM_3+6E^2M_2-3E^4\\
&=&(24t^4+72t^3+34t^2+t)-4t(6t^3+9t^2+t)+6t^2(2t^2+t)-3t^4\\
&=&9t^4+42t^3+30t^2+t
\end{eqnarray*}

Now by normalizing, we obtain the following formula for the kurtosis, at $s\geq5$:
$$\kappa=\frac{9t^4+42t^3+30t^2+t}{(t^2+t)^2}
=\frac{9t^3+42t^2+30t+1}{t(t+1)^2}$$

Thus, we are led to the conclusions in the statement.
\end{proof}

And with this, end of our study, but far more things can be said. Indeed, as mentioned in Comment 6.18, there is some cumulant work to be done, and more on this later. Then, as mentioned in the proof of Theorem 6.19, there is some work to be done at $s=3,4,\infty$, and when $s$ is even too. Next, as explained in chapter 5, one can pass from Fourier to the density via some computations, and a differential equation trick, and the same can be done in reverse, with $|f|^2$ instead of $f$, as to compute the Fourier transform of $|f|^2$. Then, as explained in \cite{bb+}, a number of things can be said about $f^s$ and $|f|^s$. And finally, as explained in \cite{bcu}, regarding more advanced questions, a number of conjectures can be made. We will be hopefully back to some of these questions, later in this book.

\section*{6e. Exercises}

This was a quite computational chapter, and as exercises on this, we have:

\begin{exercise}
Try more advanced computations, for $\chi_{1t}$ and $\chi_{1t}^2$.
\end{exercise}

\begin{exercise}
Try more advanced computations, for $\chi_{2t}$ and $\chi_{2t}^2$.
\end{exercise}

\begin{exercise}
Learn about the occurrences of the Rayleigh law $\chi_{2t}$.
\end{exercise}

\begin{exercise}
Do some further computations for $g_t^a$.
\end{exercise}

\begin{exercise}
Do the computations for $G_t^a$.
\end{exercise}

\begin{exercise}
Do some further computations for $G_{st}$.
\end{exercise}

\begin{exercise}
Do the computations for $G_{st}^a$.
\end{exercise}

\begin{exercise}
Do some further computations for $p^s_t$.
\end{exercise}

As bonus exercise, learn more about the compound Poisson laws, in general.

\chapter{Unitary groups}

\section*{7a. Representations}

We already met compact groups $G\subset U_N$ in the above, the idea being that under suitable assumptions, namely easiness and uniformity, the asymptotic laws $\mu_t$ of the truncations $\chi_t$ of the main character produce interesting semigroups $\{\mu_t\}$. To be more precise, we first saw in chapter 4 that in the orthogonal case, the relevant subgroups $G\subset O_N$, and the real probability measures $\mu_t$ that they produce, are as follows:
$$\xymatrix@R=48pt@C=45pt{
H_N\ar[r]&O_N\\
S_N\ar[u]\ar[r]&B_N\ar[u]}\qquad\xymatrix@R=25pt@C=10pt{\\ :\\ }
\qquad\xymatrix@R=45pt@C=50pt{
p_t^2\ar@{-}[r]&g_t\\
p_t\ar@{-}[u]\ar@{-}[r]&g^t_t\ar@{-}[u]}$$

We also saw in chapter 5, as a generalization of the results on the left edge, that for the reflection groups $H_N^s\subset U_N$ with $s\in\{1,2,\ldots,\infty\}$ we obtain the Bessel laws $p^s_t$. All this is quite interesting, and in this chapter we intend to systematically investigate such phenomena, with our various motivations, and precise plan, being as follows:

\begin{plan}
The easy groups $G\subset U_N$ are worth a detailed look, as follows:
\begin{enumerate}
\item We need to understand well Haar integration, Peter-Weyl and Tannaka.

\item Then have a more conceptual look at $S_N,O_N,H_N,B_N$ and $H_N^s$.

\item Then add the unitary group $U_N$ and the complex bistochastic group $C_N$.

\item Do the probabilistic computations for $G_t^t$, coming from $C_N$.

\item Have a look at Weingarten functions, this being advanced probability.

\item Generalize the Lindst\"om formula for $S_N$, leading to advanced theory too.
\end{enumerate}
\end{plan}

In short, expect a mixture of abstract algebra, to start with, followed by our regular probability business, namely brave computations for $G_t^a$, that we were afraid of, in the previous chapter, motivated by $C_N$ and helped by the choice $a=t$, and then all sorts of advanced probability considerations, featuring orthogonal polynomials and more.

\bigskip

Getting started now, our first goal will be that of explaining the Peter-Weyl theory for the compact groups $G\subset U_N$. At the beginning of everything, we have:

\begin{definition}
A unitary representation of a compact group $G$ is a continuous group morphism into a unitary group
$$u:G\to U_N\quad,\quad g\to u_g$$
which can be faithful or not. The character of such a representation is the function
$$\chi:G\to\mathbb C\quad,\quad g\to Tr(u_g)$$
where $Tr$ is the usual, unnormalized trace of the $N\times N$ matrices.
\end{definition}

At the level of examples, most of the compact groups that we met so far, finite or continuous, naturally appear as closed subgroups $G\subset U_N$. In this case, the embedding $G\subset U_N$ is of course a representation, called fundamental representation.

\bigskip

In general now, let us first discuss the various operations on the representations. We have here the following elementary result, coming from definitions:

\begin{proposition}
The representations of a compact group $G$ are subject to:
\begin{enumerate}
\item Making sums. Given representations $u,v$, of dimensions $N,M$, 
their sum is the $N+M$-dimensional representation $u+v=diag(u,v)$.

\item Making products. Given representations $u,v$, of dimensions $N,M$, their product is the $NM$-dimensional representation $(u\otimes v)_{ia,jb}=u_{ij}v_{ab}$.

\item Taking conjugates. Given a $N$-dimensional representation $u$, its conjugate is the $N$-dimensional representation $(\bar{u})_{ij}=\bar{u}_{ij}$.

\item Spinning by unitaries. Given a $N$-dimensional representation $u$, and a unitary $V\in U_N$, we can spin $u$ by this unitary, $u\to VuV^*$.
\end{enumerate}
\end{proposition}

\begin{proof}
The fact that the operations in the statement are indeed well-defined, among morphisms from $G$ to unitary groups, is indeed clear from definitions.
\end{proof}

In relation now with characters, we have the following result:

\begin{proposition}
We have the following formulae, regarding characters
$$\chi_{u+v}=\chi_u+\chi_v\quad,\quad 
\chi_{u\otimes v}=\chi_u\chi_v\quad,\quad 
\chi_{\bar{u}}=\bar{\chi}_u\quad,\quad
\chi_{VuV^*}=\chi_u$$
in relation with the basic operations for the representations.
\end{proposition}

\begin{proof}
All these assertions are elementary, coming from the following formulae:
$$Tr(diag(U,V))=Tr(U)+Tr(V)\quad,\quad 
Tr(U\otimes V)=Tr(U)Tr(V)$$
$$Tr(\bar{U})=\overline{Tr(U)}\quad,\quad 
Tr(VUV^*)=Tr(U)$$

Thus, we are led to the various character formulae in the statement.
\end{proof}

Assume now that we are given a closed subgroup $G\subset U_N$. By using the above operations, we can construct a whole family of representations of $G$, as follows:

\begin{definition}
Given a closed subgroup $G\subset U_N$, its Peter-Weyl representations are the tensor products between the fundamental representation and its conjugate:
$$u:G\subset U_N\quad,\quad 
\bar{u}:G\subset U_N$$ 
We denote these tensor products $u^{\otimes k}$, with $k=\circ\bullet\bullet\circ\ldots$ being a colored integer, with the colored tensor powers being defined according to the rules 
$$u^{\otimes\circ}=u\quad,\quad
u^{\otimes\bullet}=\bar{u}\quad,\quad
u^{\otimes kl}=u^{\otimes k}\otimes u^{\otimes l}$$
and with the convention that $u^{\otimes\emptyset}$ is the trivial representation $1:G\to U_1$.
\end{definition}

Here are a few examples of such representations, namely those coming from the colored integers of length 2, which will often appear in what follows:
$$u^{\otimes\circ\circ}=u\otimes u\quad,\quad 
u^{\otimes\circ\bullet}=u\otimes\bar{u}$$
$$u^{\otimes\bullet\circ}=\bar{u}\otimes u\quad,\quad
u^{\otimes\bullet\bullet}=\bar{u}\otimes\bar{u}$$

In order to advance, we must develop some general theory. Let us start with:

\begin{definition}
Given a compact group $G$, and two of its representations,
$$u:G\to U_N\quad,\quad 
v:G\to U_M$$
we define the space of intertwiners between these representations as being 
$$Hom(u,v)=\left\{T\in M_{M\times N}(\mathbb C)\Big|Tu_g=v_gT,\forall g\in G\right\}$$
and we use the following conventions:
\begin{enumerate}
\item We use the notations $Fix(u)=Hom(1,u)$, and $End(u)=Hom(u,u)$.

\item We write $u\sim v$ when $Hom(u,v)$ contains an invertible element.

\item We say that $u$ is irreducible, and write $u\in Irr(G)$, when $End(u)=\mathbb C1$.
\end{enumerate}
\end{definition}

To be more precise, the terminology here is very standard, with Fix, Hom, End standing respectively for fixed points, homomorphisms and endomorphisms. We will see later that irreducible means indecomposable, in a suitable sense. 

\bigskip

Here are now a few basic results, regarding the above linear spaces:

\begin{proposition}
The spaces of intertwiners have the following properties:
\begin{enumerate}
\item $T\in Hom(u,v),S\in Hom(v,w)\implies ST\in Hom(u,w)$.

\item $S\in Hom(u,v),T\in Hom(w,z)\implies S\otimes T\in Hom(u\otimes w,v\otimes z)$.

\item $T\in Hom(u,v)\implies T^*\in Hom(v,u)$.
\end{enumerate}
In abstract terms, we say that the Hom spaces form a tensor $*$-category.
\end{proposition}

\begin{proof}
All the formulae in the statement are clear from definitions, via elementary computations. As for the last assertion, this is something coming from (1,2,3). We will be back to tensor categories later on, with more details on this latter fact.
\end{proof}

In order to advance, we will need the following standard linear algebra fact:

\begin{proposition}
Let $A\subset M_N(\mathbb C)$ be a $*$-algebra.
\begin{enumerate}
\item We have $1=p_1+\ldots+p_k$, with $p_i\in A$ being central minimal projections.

\item Each of the spaces $A_i=p_iAp_i$ is a non-unital $*$-subalgebra of $A$.

\item We have a non-unital $*$-algebra sum decomposition $A=A_1\oplus\ldots\oplus A_k$.

\item We have unital $*$-algebra isomorphisms $A_i\simeq M_{n_i}(\mathbb C)$, with $n_i=rank(p_i)$.

\item Thus, we have a $*$-algebra isomorphism $A\simeq M_{n_1}(\mathbb C)\oplus\ldots\oplus M_{n_k}(\mathbb C)$.
\end{enumerate}
\end{proposition}

\begin{proof}
Consider indeed an arbitrary $*$-algebra of the $N\times N$ matrices, $A\subset M_N(\mathbb C)$. Let us first look at the center of this algebra, $Z(A)=A\cap A'$. It is elementary to prove that this center, as an algebra, is of the following form:
$$Z(A)\simeq\mathbb C^k$$

Consider now the standard basis $e_1,\ldots,e_k\in\mathbb C^k$, and let  $p_1,\ldots,p_k\in Z(A)$ be the images of these vectors via the above identification. In other words, these elements $p_1,\ldots,p_k\in A$ are central minimal projections, summing up to 1:
$$p_1+\ldots+p_k=1$$

The idea is then that this partition of the unity eventually leads to the block decomposition of $A$, as in the statement, and we leave the details here as an exercise.
\end{proof}

We can now formulate our first Peter-Weyl theorem, as follows:

\index{Peter-Weyl}

\begin{theorem}[PW1]
Let $u:G\to U_N$ be a representation, consider the algebra $A=End(u)$, and write its unit $1=p_1+\ldots+p_k$ as above. We have then 
$$u=v_1+\ldots+v_k$$
with each $v_i$ being an irreducible representation, obtained by restricting $u$ to $Im(p_i)$.
\end{theorem}

\begin{proof}
This follows indeed from Proposition 7.7 and Proposition 7.8:

\medskip

(1) We first associate to our representation $u:G\to U_N$ the corresponding action map on $\mathbb C^N$. If a linear subspace $V\subset\mathbb C^N$ is invariant, the restriction of the action map to $V$ is an action map too, which must come from a subrepresentation $v\subset u$.

\medskip

(2) Consider now a projection $p\in End(u)$. From $pu=up$ we obtain that the linear space $V=Im(p)$ is invariant under $u$, and so this space must come from a subrepresentation $v\subset u$. It is routine to check that the operation $p\to v$ maps subprojections to subrepresentations, and minimal projections to irreducible representations.

\medskip

(3) With these preliminaries in hand, let us decompose the algebra $End(u)$ as above, by using the decomposition $1=p_1+\ldots+p_k$ into central minimal projections. If we denote by $v_i\subset u$ the subrepresentation coming from the vector space $V_i=Im(p_i)$, then we obtain in this way a decomposition $u=v_1+\ldots+v_k$, as in the statement.
\end{proof}

Here is now our second Peter-Weyl theorem, complementing Theorem 7.9:

\begin{theorem}[PW2]
Given a closed subgroup $G\subset_uU_N$, any of its irreducible smooth representations 
$$v:G\to U_M$$
appears inside a tensor product of the fundamental representation $u$ and its adjoint $\bar{u}$.
\end{theorem}

\begin{proof}
This basically follows from Theorem 7.9, by reasoning as follows:

\medskip

(1) Given $v:G\to U_M$, consider its space of coefficients $C_v\subset C(G)$. The operation $w\to C_w$ is then functorial, mapping subrepresentations into linear subspaces.

\medskip

(2) A closed subgroup $G\subset_uU_N$ is a Lie group, and a representation $v:G\to U_M$ is smooth when we have an inclusion $C_v\subset<C_u>$. This is indeed well-known.

\medskip

(3) By definition of the Peter-Weyl representations, as arbitrary tensor products between the fundamental representation $u$ and its conjugate $\bar{u}$, we have:
$$<C_u>=\sum_kC_{u^{\otimes k}}$$

(4) Now by putting together the above observations (2,3) we conclude that we must have an inclusion as follows, for certain exponents $k_1,\ldots,k_p\in\mathbb N$:
$$C_v\subset C_{u^{\otimes k_1}\oplus\ldots\oplus u^{\otimes k_p}}$$

(5) By using now (1), we deduce that we have an inclusion $v\subset u^{\otimes k_1}\oplus\ldots\oplus u^{\otimes k_p}$, and by applying Theorem 7.9, this leads to the conclusion in the statement.
\end{proof}

In order to further advance, we will need the following standard fact:

\index{Haar measure}
\index{Ces\`aro limit}

\begin{theorem}
The integration over a compact group $G\subset U_N$ can be constructed by starting with any faithful positive unital linear form $\varphi\in C(G)^*$, and setting:
$$\int_G=\lim_{n\to\infty}\frac{1}{n}\sum_{k=1}^n\varphi^{*k}$$
Moreover, for any representation $v:G\to U_n$ we have the following formula,
$$\left(\int_Gv_{ij}\right)_{ij}=P$$
where $P$ is the orthogonal projection onto the linear space $Fix(v)$.
\end{theorem}

\begin{proof}
This is something that we already discussed in chapter 4, in the orthogonal case, and the proof in general is similar. To start with, given a positive unital linear form $\varphi\in C(G)^*$, our claim is that the following limit converges, for any $f\in C(G)$:
$$\int_\varphi f=\lim_{n\to\infty}\frac{1}{n}\sum_{k=1}^n\varphi^{*k}(f)$$

Indeed, by linearity we can assume that $f$ is the coefficient of certain representation $v:G\to U_n$. But in this case, an elementary computation gives the following formula, with $P_\varphi\geq P$ being the orthogonal projection onto the $1$-eigenspace of $[\varphi(v_{ij})]_{ij}$:
$$\left(\int_\varphi v_{ij}\right)_{ij}=P_\varphi$$

The point now is that when $\varphi\in C(G)^*$ is faithful, by using a standard positivity trick, we can prove that we have $P_\varphi=P$. Assume indeed $P_\varphi\xi=\xi$, and let us set:
$$f=\sum_i\left(\sum_jv_{ij}\xi_j-\xi_i\right)\overline{\left(\sum_kv_{ik}\xi_k-\xi_i\right)}$$

A straightforward computation shows then that $\varphi(f)=0$, and so $f=0$, as desired. Thus, we proved our claim, and with this in hand, the left and right invariance of $\int_G=\int_\varphi$ is clear on coefficients, and so in general, and this gives all the assertions. For details here you can check for instance my book \cite{ba1}, where all this is explained in detail.
\end{proof}

We will need as well an algebraic ingredient for our study, as follows:

\index{Frobenius isomorphism}

\begin{proposition}
We have a Frobenius type isomorphism
$$Hom(v,w)\simeq Fix(v\otimes\bar{w})$$
valid for any two representations $v,w$.
\end{proposition}

\begin{proof}
According to the definitions, we have the following equivalences:
\begin{eqnarray*}
T\in Hom(v,w)
&\iff&Tv=wT\\
&\iff&\sum_jT_{aj}v_{ji}=\sum_bw_{ab}T_{bi},\forall a,i
\end{eqnarray*}

On the other hand, we have as well the following equivalences:
\begin{eqnarray*}
T\in Fix(v\otimes\bar{w})
&\iff&(v\otimes\bar{w})T=\xi\\
&\iff&\sum_{jb}v_{ij}w_{ab}^*T_{bj}=T_{ai}\forall a,i
\end{eqnarray*}

But with this in hand, both inclusions follow from the unitarity of $v,w$.
\end{proof}

Good news, we can now formulate our third Peter-Weyl theorem, as follows:

\index{Peter-Weyl}

\begin{theorem}[PW3]
The dense subalgebra $\mathcal C(G)\subset C(G)$ generated by the coefficients of the fundamental representation decomposes as a direct sum 
$$\mathcal C(G)=\bigoplus_{v\in Irr(G)}M_{\dim(v)}(\mathbb C)$$
with the summands being pairwise orthogonal with respect to $<f,g>=\int_Gf\bar{g}$.
\end{theorem}

\begin{proof}
By combining the previous two Peter-Weyl results, we deduce that we have a linear space decomposition as follows:
$$\mathcal C(G)
=\sum_{v\in Irr(G)}C_v
=\sum_{v\in Irr(G)}M_{\dim(v)}(\mathbb C)$$

Thus, in order to conclude, it is enough to prove that for any two irreducible representations $v,w\in Irr(G)$, the corresponding spaces of coefficients are orthogonal:
$$v\not\sim w\implies C_v\perp C_w$$ 

But this follows by Frobenius duality, by integrating. Let us set indeed:
$$P_{ia,jb}=\int_Gv_{ij}\bar{w}_{ab}$$

Then $P$ is the orthogonal projection onto the following vector space:
$$Fix(v\otimes\bar{w})
\simeq Hom(v,w)
=\{0\}$$

Thus we have $P=0$, and this gives the result.
\end{proof}

Finally, we have the following result, completing the Peter-Weyl theory:

\index{Peter-Weyl}

\begin{theorem}[PW4]
The characters of irreducible representations belong to
$$\mathcal C(G)_{central}=\left\{f\in\mathcal C(G)\Big|f(gh)=f(hg),\forall g,h\in G\right\}$$
called algebra of central functions on $G$, and form an orthonormal basis of it.
\end{theorem}

\begin{proof}
Observe first that $\mathcal C(G)_{central}$ is indeed an algebra, which contains all the characters. Conversely, consider a function $f\in\mathcal C(G)$, written as follows:
$$f=\sum_{v\in Irr(G)}f_v$$

The condition $f\in\mathcal C(G)_{central}$ states then that for any $v\in Irr(G)$, we must have:
$$f_v\in\mathcal C(G)_{central}$$

But this means that $f_v$ must be a scalar multiple of $\chi_v$, so the characters form a basis of $\mathcal C(G)_{central}$, as stated. Also, the fact that we have an orthogonal basis follows from Theorem 7.13. As for the fact that the characters have norm 1, this follows from:
$$\int_G\chi_v\bar{\chi}_v
=\sum_{ij}\int_Gv_{ii}\bar{v}_{jj}
=1$$

Here we have used the fact, coming from Frobenius duality, that the various integrals $\int_Gv_{ij}\bar{v}_{kl}$ form altogether the orthogonal projection onto the following vector space:
$$Fix(v\otimes\bar{v})\simeq End(v)=\mathbb C1$$

Thus, good news, we just proved all the Peter-Weyl theorems.
\end{proof}

\section*{7b. Tannaka, easiness}

Getting now to more concrete things, which will turn to be very useful, for our probabilistic purposes, we have the following principle, coming on top of Peter-Weyl:

\begin{principle}
Any compact group $G\subset U_N$ appears as the symmetry group of its corresponding Tannakian category $C_G$,
$$G=G(C_G)$$
and by suitably delinearizing $C_G$, via a Brauer theorem of type $C_G=span(D_G)$, we can view $G$ as symmetry group of a certain combinatorial object $D_G$.
\end{principle}

Excited about this? Does not look easy, all this material, with both Tannaka and Brauer being quite scary names, in the context of algebra. But, believe me, all this is worth learning, and it is good to have in your bag some cutting-edge technology regarding the groups, such as the results of Tannaka and Brauer. So, we will go for this.

\bigskip

Getting started now, we first have a categorical definition, as follows:

\begin{definition}
A tensor category over $H=\mathbb C^N$ is a collection $C=(C_{kl})$ of linear spaces $C_{kl}\subset\mathcal L(H^{\otimes k},H^{\otimes l})$ satisfying the following conditions:
\begin{enumerate}
\item $S,T\in C$ implies $S\otimes T\in C$.

\item If $S,T\in C$ are composable, then $ST\in C$.

\item $T\in C$ implies $T^*\in C$.

\item Each $C_{kk}$ contains the identity operator.

\item $C_{\emptyset k}$ with $k=\circ\bullet,\bullet\circ$ contain the operator $R:1\to\sum_ie_i\otimes e_i$.

\item $C_{kl,lk}$ with $k,l=\circ,\bullet$ contain the flip operator $\Sigma:a\otimes b\to b\otimes a$.
\end{enumerate}
\end{definition}

Here, as usual, the tensor powers $H^{\otimes k}$, which are Hilbert spaces depending on a colored integer $k=\circ\bullet\bullet\circ\ldots\,$, are defined by the following formulae, and multiplicativity:
$$H^{\otimes\emptyset}=\mathbb C\quad,\quad 
H^{\otimes\circ}=H\quad,\quad
H^{\otimes\bullet}=\bar{H}\simeq H$$

We have already met such categories, when dealing with the Tannakian categories of the closed subgroups $G\subset U_N$, and our knowledge can be summarized as follows:

\begin{proposition}
Given a closed subgroup $G\subset U_N$, its Tannakian category
$$C_{kl}=\left\{T\in\mathcal L(H^{\otimes k},H^{\otimes l})\Big|Tg^{\otimes k}=g^{\otimes l}T,\forall g\in G\right\}$$
is a tensor category over $H=\mathbb C^N$. Conversely, given a tensor category $C$ over $\mathbb C^N$,
$$G=\left\{g\in U_N\Big|Tg^{\otimes k}=g^{\otimes l}T,\forall k,l,\forall T\in C_{kl}\right\}$$
is a closed subgroup of $U_N$.
\end{proposition}

\begin{proof}
This is something that we basically know, the idea being as follows:

\medskip

(1) Regarding the first assertion, we have to check here the axioms (1-6) in Definition 7.16. The axioms (1-4) being all clear from definitions, let us establish (5). But this follows from the fact that each element $g\in G$ is a unitary, which can be reformulated as follows, with $R:1\to\sum_ie_i\otimes e_i$ being the map in Definition 7.16:
$$R\in Hom(1,g\otimes\bar{g})\quad,\quad 
R\in Hom(1,\bar{g}\otimes g)$$

Regarding now the condition in Definition 7.16 (6), this comes from the fact that the matrix coefficients $g\to g_{ij}$ and their conjugates $g\to\bar{g}_{ij}$ commute with each other.

\medskip

(2) Regarding the second assertion, we have to check that the subset $G\subset U_N$ constructed in the statement is a closed subgroup. But this is clear from definitions.
\end{proof}

Summarizing, we have so far precise axioms for the tensor categories $C=(C_{kl})$, given in Definition 7.16, as well as correspondences as follows:
$$G\to C_G\quad,\quad 
C\to G_C$$

We will prove in what follows that these correspondences are inverse to each other. In order to get started, we first have the following technical result:

\begin{proposition}
Consider the following conditions:
\begin{enumerate}
\item $C=C_{G_C}$, for any tensor category $C$.

\item $G=G_{C_G}$, for any closed subgroup $G\subset U_N$.
\end{enumerate}
We have then $(1)\implies(2)$. Also, $C\subset C_{G_C}$ is automatic.
\end{proposition}

\begin{proof}
Given $G\subset U_N$, we have $G\subset G_{C_G}$. On the other hand, by using (1) we have $C_G=C_{G_{C_G}}$. Thus, we have an inclusion of closed subgroups of $U_N$, which becomes an isomorphism at the level of the associated Tannakian categories, so $G=G_{C_G}$. Finally, the fact that we have an inclusion $C\subset C_{G_C}$ is clear from definitions.
\end{proof}

The point now is that it is possible to prove that we have $C_{G_C}\subset C$, by doing some abstract algebra, and we are led in this way to the following conclusion:

\begin{theorem}
The Tannakian duality constructions 
$$C\to G_C\quad,\quad 
G\to C_G$$
are inverse to each other.
\end{theorem}

\begin{proof}
This is something quite tricky, the idea being as follows:

\medskip

(1) According to Proposition 7.18, we must prove $C_{G_C}\subset C$. For this purpose, given a tensor category $C=(C_{kl})$ over a Hilbert space $H$, consider the following $*$-algebra:
$$E_C
=\bigoplus_{k,l}C_{kl}
\subset\bigoplus_{k,l}B(H^{\otimes k},H^{\otimes l})
\subset B\left(\bigoplus_kH^{\otimes k}\right)$$

Consider also, inside this $*$-algebra, the following $*$-subalgebra:
$$E_C^{(s)}
=\bigoplus_{|k|,|l|\leq s}C_{kl}
\subset\bigoplus_{|k|,|l|\leq s}B(H^{\otimes k},H^{\otimes l})
=B\left(\bigoplus_{|k|\leq s}H^{\otimes k}\right)$$

(2) It is then routine to check that we have equivalences as follows:
\begin{eqnarray*}
C_{G_C}\subset C
&\iff&E_{C_{G_C}}\subset E_C\\
&\iff&E_{C_{G_C}}^{(s)}\subset E_C^{(s)},\forall s\\
&\iff&E_{C_{G_C}}^{(s)'}\supset E_C^{(s)'},\forall s
\end{eqnarray*}

(3) Summarizing, we would like to prove that we have inclusions $E_C^{(s)'}\subset E_{C_{G_C}}^{(s)'}$. But this can be done by doing some algebra, as explained by Malacarne in \cite{mal}. For an introduction to all this, you can have a look at any of my algebra books.
\end{proof}

With this discussed, let us go back now to Principle 7.15, and develop the second idea there, namely delinearization and Brauer theorems. Following \cite{bsp}, we have the following unitary update of our previous notion of category of partitions, from chapter 4:

\begin{definition}[update]
A category of partitions is a collection $D=\bigsqcup_{k,l}D(k,l)$ of subsets $D(k,l)\subset P(k,l)$, having the following properties:
\begin{enumerate}
\item Stability under the horizontal concatenation, $(\pi,\sigma)\to[\pi\sigma]$.

\item Stability under vertical concatenation $(\pi,\sigma)\to[^\sigma_\pi]$, with matching middle symbols.

\item Stability under the upside-down turning $*$, with switching of colors, $\circ\leftrightarrow\bullet$.

\item Each set $P(k,k)$ contains the identity partition $||\ldots||$.

\item The sets $P(\emptyset,\circ\bullet)$ and $P(\emptyset,\bullet\circ)$ both contain the semicircle $\cap$.

\item The sets $P(k,\bar{k})$ with $|k|=2$ contain the crossing partition $\slash\hskip-2.0mm\backslash$.
\end{enumerate}
\end{definition} 

Observe the similarity with Definition 7.16, and more on this in a moment. In order now to construct a Tannakian category out of such a category, we will need:

\begin{proposition}
Each partition $\pi\in P(k,l)$ produces a linear map
$$T_\pi:(\mathbb C^N)^{\otimes k}\to(\mathbb C^N)^{\otimes l}$$
given by the following formula, with $e_1,\ldots,e_N$ being the standard basis of $\mathbb C^N$,
$$T_\pi(e_{i_1}\otimes\ldots\otimes e_{i_k})=\sum_{j_1\ldots j_l}\delta_\pi\begin{pmatrix}i_1&\ldots&i_k\\ j_1&\ldots&j_l\end{pmatrix}e_{j_1}\otimes\ldots\otimes e_{j_l}$$
and with the Kronecker type symbols $\delta_\pi\in\{0,1\}$ depending on whether the indices fit or not. The assignement $\pi\to T_\pi$ is categorical, in the sense that we have
$$T_\pi\otimes T_\sigma=T_{[\pi\sigma]}\quad,\quad 
T_\pi T_\sigma=N^{c(\pi,\sigma)}T_{[^\sigma_\pi]}\quad,\quad 
T_\pi^*=T_{\pi^*}$$
where $c(\pi,\sigma)$ are certain integers, coming from the erased components in the middle.
\end{proposition}

\begin{proof}
This is something elementary, the computations being as follows:

\medskip

(1) The concatenation axiom follows from the following computation:
\begin{eqnarray*}
&&(T_\pi\otimes T_\sigma)(e_{i_1}\otimes\ldots\otimes e_{i_p}\otimes e_{k_1}\otimes\ldots\otimes e_{k_r})\\
&=&\sum_{j_1\ldots j_q}\sum_{l_1\ldots l_s}\delta_\pi\begin{pmatrix}i_1&\ldots&i_p\\j_1&\ldots&j_q\end{pmatrix}\delta_\sigma\begin{pmatrix}k_1&\ldots&k_r\\l_1&\ldots&l_s\end{pmatrix}e_{j_1}\otimes\ldots\otimes e_{j_q}\otimes e_{l_1}\otimes\ldots\otimes e_{l_s}\\
&=&\sum_{j_1\ldots j_q}\sum_{l_1\ldots l_s}\delta_{[\pi\sigma]}\begin{pmatrix}i_1&\ldots&i_p&k_1&\ldots&k_r\\j_1&\ldots&j_q&l_1&\ldots&l_s\end{pmatrix}e_{j_1}\otimes\ldots\otimes e_{j_q}\otimes e_{l_1}\otimes\ldots\otimes e_{l_s}\\
&=&T_{[\pi\sigma]}(e_{i_1}\otimes\ldots\otimes e_{i_p}\otimes e_{k_1}\otimes\ldots\otimes e_{k_r})
\end{eqnarray*}

(2) The composition axiom follows from the following computation:
\begin{eqnarray*}
&&T_\pi T_\sigma(e_{i_1}\otimes\ldots\otimes e_{i_p})\\
&=&\sum_{j_1\ldots j_q}\delta_\sigma\begin{pmatrix}i_1&\ldots&i_p\\j_1&\ldots&j_q\end{pmatrix}
\sum_{k_1\ldots k_r}\delta_\pi\begin{pmatrix}j_1&\ldots&j_q\\k_1&\ldots&k_r\end{pmatrix}e_{k_1}\otimes\ldots\otimes e_{k_r}\\
&=&\sum_{k_1\ldots k_r}N^{c(\pi,\sigma)}\delta_{[^\sigma_\pi]}\begin{pmatrix}i_1&\ldots&i_p\\k_1&\ldots&k_r\end{pmatrix}e_{k_1}\otimes\ldots\otimes e_{k_r}\\
&=&N^{c(\pi,\sigma)}T_{[^\sigma_\pi]}(e_{i_1}\otimes\ldots\otimes e_{i_p})
\end{eqnarray*}

(3) Finally, the involution axiom follows from the following computation:
\begin{eqnarray*}
&&T_\pi^*(e_{j_1}\otimes\ldots\otimes e_{j_q})\\
&=&\sum_{i_1\ldots i_p}<T_\pi^*(e_{j_1}\otimes\ldots\otimes e_{j_q}),e_{i_1}\otimes\ldots\otimes e_{i_p}>e_{i_1}\otimes\ldots\otimes e_{i_p}\\
&=&\sum_{i_1\ldots i_p}\delta_\pi\begin{pmatrix}i_1&\ldots&i_p\\ j_1&\ldots& j_q\end{pmatrix}e_{i_1}\otimes\ldots\otimes e_{i_p}\\
&=&T_{\pi^*}(e_{j_1}\otimes\ldots\otimes e_{j_q})
\end{eqnarray*}

Summarizing, our correspondence is indeed categorical.
\end{proof}

Following \cite{bsp}, we can now formulate a key theoretical result, as follows:

\begin{theorem}
Any category of partitions $D\subset P$ produces a series of compact groups $G=(G_N)$, with $G_N\subset U_N$ for any $N\in\mathbb N$, via the formula
$$C_{kl}=span\left(T_\pi\Big|\pi\in D(k,l)\right)$$
for any $k,l$, and Tannakian duality. We call such groups easy.
\end{theorem}

\begin{proof}
Indeed, once we fix an integer $N\in\mathbb N$, the various axioms in Definition 7.20 show, via Proposition 7.21, that the following spaces form a Tannakian category:
$$span\left(T_\pi\Big|\pi\in D(k,l)\right)$$

Thus, Tannakian duality applies, and provides us with a closed subgroup $G_N\subset U_N$ such that the following equalities are satisfied, for any colored integers $k,l$:
$$C_{kl}=span\left(T_\pi\Big|\pi\in D(k,l)\right)$$

We are therefore led to the conclusion in the statement.
\end{proof}

And with this, good news, done with the general theory. At the level of examples, and generalizing what we previously knew from chapters 4 and 5, we have:

\begin{theorem}
We have easy groups as follows, with calligraphic standing for matching, meaning satisfying $\#\circ=\#\bullet$, as a weighted equality, in each block: 
$$
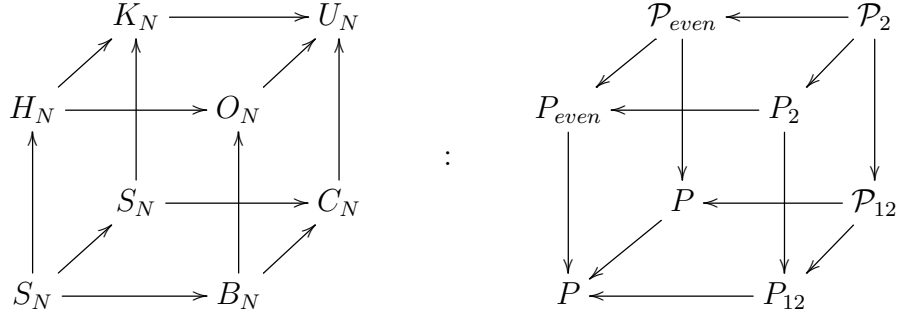
\xymatrix@R=18pt@C=14pt{
&K_N\ar[rr]&&U_N\\
H_N\ar[rr]\ar[ur]&&O_N\ar[ur]\\
&S_N\ar[rr]\ar[uu]&&C_N\ar[uu]\\
S_N\ar[uu]\ar[ur]\ar[rr]&&B_N\ar[uu]\ar[ur]
}\qquad\xymatrix@R=18pt@C=14pt{\\ \\ :}
\qquad\xymatrix@R=18pt@C=10pt{
&\mathcal P_{even}\ar[dl]\ar[dd]&&\mathcal P_2\ar[dl]\ar[ll]\ar[dd]\\
P_{even}\ar[dd]&&P_2\ar[dd]\ar[ll]\\
&P\ar[dl]&&\mathcal P_{12}\ar[dl]\ar[ll]\\
P&&P_{12}\ar[ll]}$$
Moreover, the reflection groups $H_N^s$ with $s\in\{1,2,\ldots,\infty\}$ fit on the diagonal of the left face of the left cube, corresponding to the categories $P^s$, fitting on the right cube. 
\end{theorem}

\begin{proof}
We already know most of these results, from chapters 4-5, but with our discussion there being a bit amateurish anyway, time to redo everything:

\medskip

(1) The unitary group $U_N$ being defined via the relations $u^*=u^{-1}$, $u^t=\bar{u}^{-1}$, the associated Tannakian category is $C=span(T_\pi|\pi\in D)$, with:
$$D
=<{\ }^{\,\cap}_{\circ\bullet}\,\,,{\ }^{\,\cap}_{\bullet\circ}>
=\mathcal P_2$$

(2) The orthogonal group $O_N\subset U_N$ being defined by imposing the relations $u_{ij}=\bar{u}_{ij}$, the associated Tannakian category is $C=span(T_\pi|\pi\in D)$, with:
$$D
=<\mathcal P_2,|^{\hskip-1.32mm\circ}_{\hskip-1.32mm\bullet},|_{\hskip-1.32mm\circ}^{\hskip-1.32mm\bullet}>
=P_2$$

(3) The unitary bistochastic group $C_N\subset U_N$ being defined by imposing the relations $u\xi=\xi$, $\bar{u}\xi=\xi$, the associated Tannakian category is $C=span(T_\pi|\pi\in D)$, with:
$$D=<\mathcal P_2,|_{\hskip-1.32mm\circ},|_{\hskip-1.32mm\bullet}>
=\mathcal P_{12}$$
  
(4) The orthogonal bistochastic group appears as $B_N=C_N\cap O_N$, and we conclude that the associated Tannakian category is $C=span(T_\pi|\pi\in D)$, with:
$$D=<\mathcal P_{12},P_2>=P_{12}$$

(5) In order to discuss now $S_N$, consider the one-block ``fork'' partition, namely:
$$\xymatrix@R=1mm@C=2mm{\\ \\ \mu\ \ =\\ \\ }\ \ \ 
\xymatrix@R=2mm@C=3mm{
\circ\ar@/_/@{-}[dr]&&\circ\\
&\ar@/_/@{-}[ur]\ar@{-}[dd]\\
&&&\\
&\circ}$$

We have then $T_\mu(e_i\otimes e_j)=\delta_{ij}e_i$, and by using this formula, we obtain:
$$T_\mu\in Hom(u^{\otimes 2},u)\iff u_{ij}u_{ik}=\delta_{jk}u_{ij},\forall i,j,k$$

But on the right we have the relations defining $S_N\subset O_N$, and we conclude that $S_N$ is indeed easy, coming from the following category of partitions:
$$D=<\mu>=P$$

(6) In order to discuss $H_N$, consider the following one-block partition:
$$\xymatrix@R=0.5mm@C=2mm{\\ \\ \\ \chi\ \ =\\ \\ }\ \ \ 
\xymatrix@R=2mm@C=3mm{
\circ\ar@/_/@{-}[dr]&&\circ\\
&\ar@/_/@{-}[ur]\ar@{-}[dd]\\
&&&\\
&\ar@/^/@{-}[dr]\ar@/_/@{-}[dl]\\
\circ&&\circ}$$

We have then $T_\chi(e_i\otimes e_j)=\delta_{ij}e_i\otimes e_i$, and by using this formula, we obtain:
$$T_\chi\in End(u^{\otimes 2})\iff \delta_{ik}u_{ia}u_{ib}=\delta_{ab}u_{ia}u_{ka},\forall i,k,a,b$$

But on the right we have the relations defining $H_N\subset O_N$, and we conclude that $H_N$ is indeed easy, coming from the following category of partitions:
$$D=<\chi>=P_{even}$$

(7) In order to discuss now $K_N\subset U_N$, consider the following colored partition:
$$\xymatrix@R=0.5mm@C=2mm{\\ \\ \\ \chi'\ \ =\\ \\ }\ \ \ 
\xymatrix@R=2mm@C=3mm{
\circ\ar@/_/@{-}[dr]&&\bullet\\
&\ar@/_/@{-}[ur]\ar@{-}[dd]\\
&&&\\
&\ar@/^/@{-}[dr]\ar@/_/@{-}[dl]\\
\bullet&&\circ}$$

Our computations from the previous proof, for the group $H_N$, modify into:
$$T_{\chi'}\in Hom(u\otimes\bar{u},\bar{u}\otimes u)\iff \delta_{ik}u_{ia}\bar{u}_{ib}=\delta_{ab}\bar{u}_{ia}u_{ka},\forall i,k,a,b$$

But on the right we have the relations defining $K_N\subset U_N$, and we conclude that $K_N$ is indeed easy, coming from the following category of partitions:
$$D=<\chi'>=\mathcal P_{even}$$

(8) Finally, regarding $H_N^s\subset K_N$, consider the following partition, with $s+2$ legs:
$$\xymatrix@R=6pt@C=12pt{
&&\ar@{-}[dd]\ar@{-}[rrrrr]&\ar@{-}[dd]&&\ar@{-}[dd]&\ar@{-}[dd]&\ar@{-}[dd]\\
\xi&=&&&\ldots\\
&&\circ&\circ&&\circ&\circ&\bullet}$$

We have then $T_\xi=\sum_je_j^{\otimes s+2}$, and by using this formula, inside $K_N$, we obtain:
$$T_\xi\in Fix(u^{\otimes s+1}\otimes\bar{u})\iff\sum_ju_{ij}^{s+1}\bar{u}_{ij}=1$$

But on the right we have the relations defining $H_N^s\subset K_N$, and we conclude that $H_N^s$ is indeed easy, coming from the following category of partitions:
$$D=<\mathcal P_{even},\xi>=P^s$$

Summarizing, theorem proved, in exactly 2 pages. Tannakian duality rules.
\end{proof}

\section*{7c. Laws of characters}

Getting now towards probability and applications, we will need the following key notion, that we already met in chapter 4, in the orthogonal case:

\index{uniform group}

\begin{proposition}
For an easy group $G=(G_N)$, coming from a category of partitions $D\subset P$, the following conditions are equivalent:
\begin{enumerate}
\item $G_{N-1}=G_N\cap U_{N-1}$, via the embedding $U_{N-1}\subset U_N$ given by $u\to diag(u,1)$.

\item $G_{N-1}=G_N\cap U_{N-1}$, via the $N$ possible diagonal embeddings $U_{N-1}\subset U_N$.

\item $D$ is stable under the operation which consists in removing blocks.
\end{enumerate}
If these conditions are satisfied, we say that $G=(G_N)$ is uniform.
\end{proposition}

\begin{proof}
This is indeed something very standard, already discussed in chapter 4 in the orthogonal case, and that we will leave now, as back then, as an exercise.
\end{proof}

We can now formulate a nice probabilistic result, as follows:

\begin{theorem}
For a uniform easy group $G=(G_N)$, we have the formula
$$\lim_{N\to\infty}\int_{G_N}\chi_t^k=\sum_{\pi\in D(k)}t^{|\pi|}$$
with $D\subset P$ being the corresponding category of partitions.
\end{theorem}

\begin{proof}
This is again something very standard, already discussed in chapter 4 in the orthogonal case, and the proof in the general unitary case is similar:

\medskip

(1) At $t=1$ the result, which does not need the uniformity assumption, follows straight from easiness, and from the Lindst\"om linear independence result from chapter 4:
\begin{eqnarray*}
\lim_{N\to\infty}\int_{G_N}\chi^k
&=&\lim_{N\to\infty}\int_{G_N}\chi_{u^{\otimes k}}\\
&=&\lim_{N\to\infty}\dim\left[Fix(u^{\otimes k})\right]\\
&=&\lim_{N\to\infty}\dim\left[span\left(T_\pi\Big|\pi\in D(k)\right)\right]\\
&=&|D(k)|
\end{eqnarray*}

(2) At $t\in(0,1)$ things are more tricky, requiring the use of the Weingarten formula, which is as follows, similar to the one from the orthogonal case, from chapter 4:
$$\int_{G_N}u_{i_1j_1}^{e_1}\ldots u_{i_kj_k}^{e_k}=\sum_{\pi,\nu\in D(k)}\delta_\pi(i)\delta_\nu(j)W_{kN}(\pi,\nu)$$

To be more precise, here $k=(e_1,\ldots,e_k)$ is a colored integer, the $\delta$ signs are usual Kronecker type symbols, checking whether the indices match, and $W_{kN}=G_{kN}^{-1}$ is the inverse of the Gram matrix $G_{kN}(\pi,\nu)=N^{|\pi\vee\nu|}$, with $|.|$ being the number of blocks. 

\medskip

(3) As for the proof of this formula, this is similar to the proof in the orthogonal case, from chapter 4. Consider indeed the above integrals, denoted as follows:
$$P_{i_1\ldots i_k,j_1\ldots j_k}=\int_{G_N}u_{i_1j_1}^{e_1}\ldots u_{i_kj_k}^{e_k}$$

We know from Theorem 7.11 that the matrix $P=(P_{ij})$ is the projection on $Fix(u^{\otimes k})$. On the other hand, by easiness, this space that we are projecting on is given by:
$$Fix(u^{\otimes k})=span\left(T_\pi\Big|\pi\in D(k)\right)$$

In order now to explicitly compute $P$, consider the following linear map:
$$E(x)=\sum_{\pi\in D(k)}<x,T_\pi>T_\pi$$

By linear algebra we have then $P=WE$, where $W$ is the inverse on $Fix(u^{\otimes k})$ of the restriction of $E$. But this latter restriction is the linear map given by the Gram matrix $G_{kN}$, and we are led in this way to the Weingarten formula in (2).

\medskip

(4) Still following the material from chapter 4, by applying this formula, we obtain:
\begin{eqnarray*}
\int_{G_N}(u_{11}+\ldots +u_{ss})^k
&=&\sum_{i_1=1}^{s}\ldots\sum_{i_k=1}^s\int_{G_N}u_{i_1i_1}\ldots u_{i_ki_k}\\
&=&\sum_{\pi,\nu\in D(k)}W_{kN}(\pi,\nu)\sum_{i_1=1}^{s}\ldots\sum_{i_k=1}^s\delta_\pi(i)\delta_\nu(i)\\
&=&\sum_{\pi,\nu\in D(k)}W_{kN}(\pi,\nu)G_{ks}(\nu,\pi)\\
&=&Tr(W_{kN}G_{ks})
\end{eqnarray*}

(5) The point now is that in the uniform case the Gram matrix is asymptotically diagonal, $G_{kN}\simeq diag(N^{|\pi|})$, with this being something clear for all the examples, and true in general too, say via the classification results in \cite{twe}. Now the inverse of $G_{kN}$ being asymptotically diagonal too, $W_{kN}\simeq diag(N^{-|\pi|})$, the above computation gives:
$$\lim_{N\to\infty}\int_{G_N}\chi_t^k=\sum_{\pi\in D(k)}t^{|\pi|}$$

Thus, we have indeed the asymptotic moment formula in the statement.
\end{proof}

At the level of the examples, upgrading our previous knowledge, we have:

\begin{theorem}
For the main easy groups we obtain the following laws,
$$
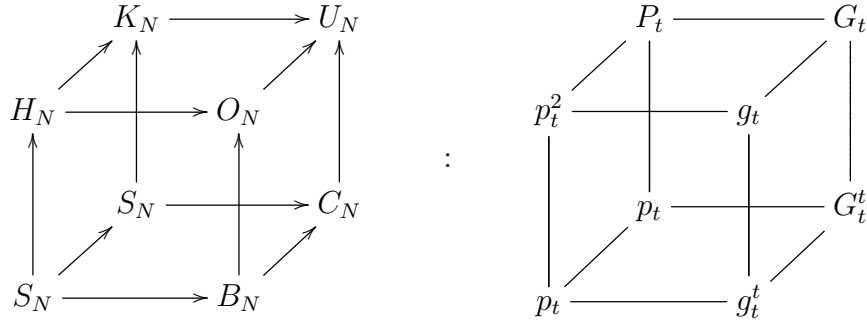
\xymatrix@R=18pt@C=14pt{
&K_N\ar[rr]&&U_N\\
H_N\ar[rr]\ar[ur]&&O_N\ar[ur]\\
&S_N\ar[rr]\ar[uu]&&C_N\ar[uu]\\
S_N\ar[uu]\ar[ur]\ar[rr]&&B_N\ar[uu]\ar[ur]
}\qquad\xymatrix@R=18pt@C=14pt{\\ \\ :}
\qquad\xymatrix@R=16.5pt@C=20pt{
&P_t\ar@{-}[rr]\ar@{-}[dd]&&G_t\ar@{-}[dd]\\
p_t^2\ar@{-}[rr]\ar@{-}[dd]\ar@{-}[ur]&&g_t\ar@{-}[dd]\ar@{-}[ur]\\
&p_t\ar@{-}[rr]\ar@{-}[uu]&&G_t^t\ar@{.}[uu]\\
p_t\ar@{-}[uu]\ar@{-}[ur]\ar@{-}[rr]&&g_t^t\ar@{-}[uu]\ar@{-}[ur]
}$$
and on the left we can insert $H_N^s$ with $s\in\{1,2,\ldots,\infty\}$, corresponding to $p^s_t$.
\end{theorem}

\begin{proof}
These results all follow from the moment formula in Theorem 7.25, and with this being actually something that we already know, from chapters 4-5, save for $U_N$ where the result is clear, and for $C_N$, which needs a bit of work, as follows:

\medskip

(1) According to Theorems 7.23 and 7.25, for the group $C_N$, we have:
$$\lim_{N\to\infty}\int_{C_N}\chi_t^k=\sum_{\pi\in\mathcal P_{12}(k)}t^{|\pi|}
=\sum_{l\subset_u k}t^{|k|-|l|/2}(|l|/2)!$$

To be more precise here, in all this $k=\circ\bullet\bullet\circ\ldots$ is a colored integer, and on the right $l\subset_uk$ stands for the choice of a colored subinteger, which is uniform in the sense of chapter 5, meaning containing the same number of $\circ$ and $\bullet$ symbols. This uniform subinteger $l\subset_uk$ corresponds to the matching pairings inside $\pi$, and there are $(|l|/2)!$ choices of such matching pairings, and the $|k|-|l|$ points left are singletons.

\medskip

(2) On the other hand, according to the formulae in chapter 6, we have:
\begin{eqnarray*}
M_k(G_t^a)
&=&\frac{1}{\pi t}\int_\mathbb Cz^ke^{-|z-a|^2/t}dz\\
&=&\frac{1}{\pi t}\int_\mathbb C(z+a)^ke^{-|z|^2/t}dz\\
&=&\frac{1}{\pi t}\int_\mathbb C\sum_{l\subset k}z^la^{k-l}e^{-|z|^2/t}dz\\
&=&\sum_{l\subset k}a^{k-l}\frac{1}{\pi t}\int_\mathbb Cz^le^{-|z|^2/t}dz\\
&=&\sum_{l\subset k}a^{k-l}M_l(g_t)
\end{eqnarray*}

But we know from chapter 5 that $M_l(g_t)$ vanishes when $l$ is not uniform, and is given by $M_l(g_t)=t^{|l|/2}(|l|/2)!$ in the uniform case. Thus, our moment formula reads:
$$M_k(G_t^a)=\sum_{l\subset_u k}a^{k-l}t^{|l|/2}(|l|/2)!$$

(3) Now in the case where $a$ is a real number, the quantity $a^{k-l}$ is a usual power of $a$, the exponent being the length $|k-l|=|k|-|l|$, and our formula becomes:
$$M_k(G_t^a)=\sum_{l\subset_u k}a^{|k|-|l|}t^{|l|/2}(|l|/2)!$$

But with $a=t$ we obtain precisely the moment formula in (1), and we are done.
\end{proof}

Quite nice all this, and as a continuation, regarding the newcomer $G_t^t$, we have:

\begin{theorem}
For a variable following a shifted normal law, $f\sim G_t^t$, we have
$$E(|f|^{2k})=\sum_{a=0}^k\binom{k}{a}^2a!t^{2k-a}$$
the expectation of $|f|^2$ is given by $E=t^2+t$, the variance is $V=2t^3+t^2$, and 
$$\gamma=\frac{6t+12}{\sqrt{(2t+1)^3}}
\quad,\quad \kappa=\frac{12t^2+36t+11}{(2t+1)^2}$$
are the corresponding skewness and kurtosis.
\end{theorem}

\begin{proof}
According to Theorems 7.23 and 7.25 we have the following formula, with $a$ standing for the number of pairings, and $2k-2a$ for the number of singletons:
$$E(|f|^{2k})=\sum_{\pi\in\mathcal P_{12}(\circ^k\bullet^k)}t^{|\pi|}
=\sum_{a=0}^k\binom{k}{a}^2a!t^{2k-a}$$

In practice, this gives the following formulae, for the low order moments:
$$M_1=t^2+t$$
$$M_2=t^4+4t^3+2t^2$$
$$M_3=t^6+9t^5+18t^4+6t^3$$
$$M_4=t^8+16t^7+72t^6+96t^5+24t^4$$

Thus $E,V$ are given by the formulae in the statement, and then:
$$M_3'=M_3-3EM_2+2E^3=6t^4+12t^3$$
$$M_4'=M_4-4EM_3+6E^2M_2-3E^4=12t^6+36t^5+11t^4$$

Now by normalizing, we obtain the above formulae for the skewness and kurtosis.
\end{proof}

With this discussed, what is next? Looking for more examples I guess, but a bit as previously in the orthogonal case, in chapter 4, what we have is pretty much everything. To be more precise, the situation with the classification work is as follows:

\bigskip

(1) As explained in chapter 4, in the orthogonal case the only uniform easy groups are $S_N,O_N,H_N,B_N$. When lifting the uniformity assumption we have two more examples, namely $S_N'=S_N\times\mathbb Z_2$ and $B_N'=B_N\times\mathbb Z_2$, which are not very interesting. See \cite{bsp}.

\bigskip

(2) In the unitary case the situation is substantially more complex, due to the joint present of three operations, $G\times\mathbb Z_s$, $G\cdot\mathbb Z_s$, $\mathbb Z_s\wr G$. However, the classification can be done, and probabilistically, the main examples remain those in Theorem 7.26. See \cite{twe}.

\bigskip

Well, nevermind. What we have in Theorem 7.26 is not that bad, providing an alternative to the 7 laws scheme from chapter 3, so let us keep building on that:

\begin{plan}
What we can do, in the remainder is this chapter, is to:
\begin{enumerate}
\item Find a clever way, for fixing the lower left edge of the cube.

\item Further study the Weingarten functions.

\item Have a look as well at the Gram determinants.

\item And at Hankel determinants and orthogonal polynomials too. 
\end{enumerate}
\end{plan}

Getting started now, with (1), the lower left edge of our cube looks bad indeed, and what is the fix? Not clear at all, hope you agree with me, and as usual in such difficult situations, we will have to ask for advice. And with the cats being gone, we will have to ask the rat. And here is what Diogenes says, from his hideout under the stove:

\begin{rat}
Peace upon you, and always aim up. Adding new layers of knowledge will gradually erase the previous levels of knowledge, be them correct or flawed.
\end{rat}

Well, looks that this rat fellow knows a thing or two, and no wonder here, with such a modest lifestyle, plenty of time for thinking at mathematics and philosophy. So thanks a lot, and following now your advice, let us have a new look at our cubes from Theorems 7.23 and 7.26, and try to build upwards. And, here is where this leads us:

\begin{theorem}
We can further build on the cube, with formal easy objects as follows,
$$
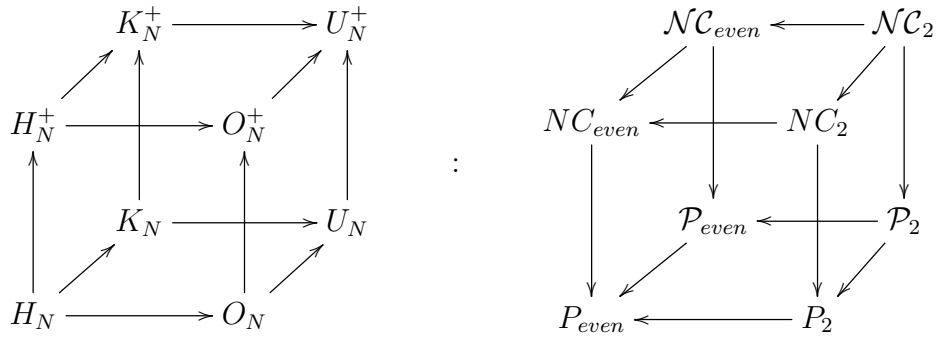
\xymatrix@R=18pt@C=15pt{
&K_N^+\ar[rr]&&U_N^+\\
H_N^+\ar[rr]\ar[ur]&&O_N^+\ar[ur]\\
&K_N\ar[rr]\ar[uu]&&U_N\ar[uu]\\
H_N\ar[uu]\ar[ur]\ar[rr]&&O_N\ar[uu]\ar[ur]
}\qquad\xymatrix@R=18pt@C=14pt{\\ \\ :}
\qquad\xymatrix@R=20pt@C=2pt{
&\mathcal{NC}_{even}\ar[dl]\ar[dd]&&\mathcal{NC}_2\ar[dl]\ar[ll]\ar[dd]\\
NC_{even}\ar[dd]&&NC_2\ar[dd]\ar[ll]\\
&\mathcal P_{even}\ar[dl]&&\mathcal P_2\ar[dl]\ar[ll]\\
P_{even}&&P_2\ar[ll]}$$
and the family of corresponding asymptotic character laws becomes
$$
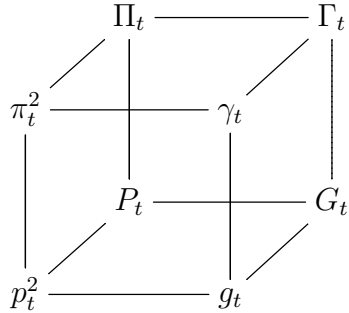
\xymatrix@R=16.5pt@C=20pt{
&\Pi_t\ar@{-}[rr]\ar@{-}[dd]&&\Gamma_t\ar@{-}[dd]\\
\pi_t^2\ar@{-}[rr]\ar@{-}[dd]\ar@{-}[ur]&&\gamma_t\ar@{-}[dd]\ar@{-}[ur]\\
&P_t\ar@{-}[rr]\ar@{-}[uu]&&G_t\ar@{.}[uu]\\
p_t^2\ar@{-}[uu]\ar@{-}[ur]\ar@{-}[rr]&&g_t\ar@{-}[uu]\ar@{-}[ur]
}$$
with the new laws $\pi^2_t,\Pi_t,\Gamma_t$ being formally defined via $M_k=\sum_{\pi\in D(k)}t^{|\pi|}$.
\end{theorem}

\begin{proof}
This is something a bit speculatory, but definitely worth the attention, see how good the above cubes look, the idea with all this being as follows:

\medskip

(1) As starting point, let us have the moment formulae for the Marchenko-Pastur and Wigner laws that we learned in chapter 3, which were as follows:
$$M_k(\pi_t)=\sum_{\pi\in NC(k)}t^{|\pi|}\quad,\quad M_k(\gamma_t)=\sum_{\pi\in NC_2(k)}t^{|\pi|}$$

Obviously, these two formulae don't fit into the framework of Definition 7.20, due to the lack of the basic crossing $\slash\hskip-2.0mm\backslash$ in $NC$ and $NC_2$, involved in axiom (6) there.

\medskip

(2) This being said, recall that abstract algebra is something quite flexible, with in fact Definition 7.20 being itself an upgrade of a previous definition from chapter 4. So, let us further upgrade Definition 7.20, by removing the crossing axiom (6) there.

\medskip

(3) Which is obviously something fruitful, because we have now $\pi_t,\gamma_t$ covered by our formalism. As question, however, what are the group-type objects that we obtain, via Tannakian duality? And in answer, let us not worry with this. We will see later that these are quantum groups, but for the moment, we will not really need this. 

\medskip

(4) So, we have our formalism, and as a matter of having some notations too, given an easy group $G$ coming from an old-style category of partitions $D$, let us denote by $G^+$ the beast coming from the new-style category of partitions $D\cap NC$. Observe that, by functoriality, we have an inclusion $G\subset G^+$. We will call this inclusion ``liberation''.

\medskip

(5) Now getting back to Theorem 7.23, we can liberate things there, leading to:
$$
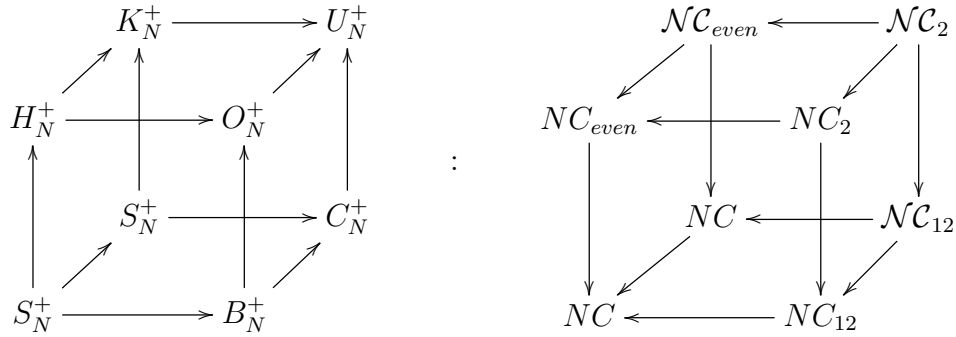
\xymatrix@R=16.6pt@C=15pt{
&K_N^+\ar[rr]&&U_N^+\\
H_N^+\ar[rr]\ar[ur]&&O_N^+\ar[ur]\\
&S_N^+\ar[rr]\ar[uu]&&C_N^+\ar[uu]\\
S_N^+\ar[uu]\ar[ur]\ar[rr]&&B_N^+\ar[uu]\ar[ur]
}\qquad\xymatrix@R=18pt@C=14pt{\\ \\ :}
\qquad\xymatrix@R=20pt@C=2pt{
&\mathcal{NC}_{even}\ar[dl]\ar[dd]&&\mathcal{NC}_2\ar[dl]\ar[ll]\ar[dd]\\
NC_{even}\ar[dd]&&NC_2\ar[dd]\ar[ll]\\
&NC\ar[dl]&&\mathcal{NC}_{12}\ar[dl]\ar[ll]\\
NC&&NC_{12}\ar[ll]}$$

Which looks very nice, and with the corresponding formal asymptotic character laws being as follows, with $\pi_t,\gamma_t$ being the Marchenko-Pastur and Wigner laws, and with the remaining laws $\pi^2_t,\Pi_t,\gamma_t^t,\Gamma_t,\Gamma_t^t$ being formally defined via $M_k=\sum_{\pi\in D(k)}t^{|\pi|}$:
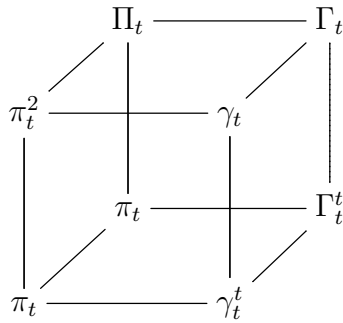
$$\xymatrix@R=16.5pt@C=20pt{
&\Pi_t\ar@{-}[rr]\ar@{-}[dd]&&\Gamma_t\ar@{-}[dd]\\
\pi_t^2\ar@{-}[rr]\ar@{-}[dd]\ar@{-}[ur]&&\gamma_t\ar@{-}[dd]\ar@{-}[ur]\\
&\pi_t\ar@{-}[rr]\ar@{-}[uu]&&\Gamma_t^t\ar@{.}[uu]\\
\pi_t\ar@{-}[uu]\ar@{-}[ur]\ar@{-}[rr]&&\gamma_t^t\ar@{-}[uu]\ar@{-}[ur]}$$

(6) And here comes the point. Now that we have two systems of cubes, the classical ones and the free ones, both suffering from unreliable bottom faces, well, we can follow Rat 7.29, stack the upper faces of the free cubes on the upper faces of the classical cubes, and we are led in this way to the flawless cubes in the statement. Wonderful.
\end{proof}

\section*{7d. Gram determinants}

Getting back now to Plan 7.28, with (1) there discussed, we are left with (2,3,4), plenty of things to be done, and in addition, our formalism and cubes have changed, we would like to deal with both the easy groups, and their liberations.

\bigskip

In short, here we are in a fairytale landscape, where many things can be done. Following \cite{bcu}, we would like to focus now on the Gram determinant question, which is quite central, in all this. By restricting the attention to the orthogonal case, we have:

\begin{question}
In the context of the main orthogonal easy quantum groups
$$\xymatrix@R=17pt@C=17pt{
&H_N^+\ar[rr]&&O_N^+\\
S_N^+\ar[rr]\ar[ur]&&B_N^+\ar[ur]\\
&H_N\ar[rr]\ar[uu]&&O_N\ar[uu]\\
S_N\ar[uu]\ar[ur]\ar[rr]&&B_N\ar[uu]\ar[ur]
}$$
how do the Gram determinants behave, on the vertical?
\end{question}

In answer now, we first need to know how to compute these determinants. In the classical discrete case, the answer is something quite simple, as follows:

\begin{theorem}
For the groups $G=S_N,H_N$ we have the Lindst\"om formula
$$\det(G_{kN})=\prod_{\pi\in D(k)}\frac{N!}{(N-|\pi|)!}$$
with $D=P,P_{even}$, and with $|.|$ being the number of blocks.
\end{theorem}

\begin{proof}
This is something that we already know for $G=S_N$, from chapter 4, and the proof for $G=H_N$ is similar, based in the fact that the corresponding categories of partitions $D=P,P_{even}$ have the property of forming semilattices under $\vee$. Consider indeed the following matrix, obtained by making determinant-preserving operations:
$$G_{kN}'(\pi,\sigma)=\sum_{\pi\leq\tau}\mu(\pi,\tau)N^{|\tau\vee\sigma|}$$

It follows then from the M\"obius inversion formula that we have:
$$G_{kN}'(\pi,\sigma)=
\begin{cases}
N(N-1)\ldots(N-|\sigma|+1)&{\rm if}\ \pi\leq\sigma\\
0&{\rm otherwise}
\end{cases}$$

Thus the matrix is upper triangular, and by computing the product on the diagonal we obtain the formula in the statement, exactly as in the $G=S_N$ case before.
\end{proof}

Next, let us discuss the case of the orthogonal group $O_N$. Here things are more complicated, the combinatorics being that of the Young diagrams. In relation with these, let us denote by $|.|$ the number of boxes, and use quantity $f^\lambda$, which gives the number of standard Young tableaux of shape $\lambda$. With these conventions, we have:

\index{Young tableaux}

\begin{theorem}
The determinant of the Gram matrix of $O_N$ is given by
$$\det(G_{kN})=\prod_{|\lambda|=k/2}f_N(\lambda)^{f^{2\lambda}}$$
where the quantities on the right are $f_N(\lambda)=\prod_{(i,j)\in\lambda}(N+2j-i-1)$.
\end{theorem}

\begin{proof}
For the group $O_N$ the Gram matrix is diagonalizable, as follows:
$$G_{kN}=\sum_{|\lambda|=k/2}f_N(\lambda)P_{2\lambda}$$

To be more precise, here $1=\sum P_{2\lambda}$ is the standard partition of unity associated to the Young diagrams having $k/2$ boxes, and the coefficients $f_N(\lambda)$ are those in the statement. Now since we have $Tr(P_{2\lambda})=f^{2\lambda}$, this gives the formula in the statement.
\end{proof}

As a variation of the above result, for the bistochastic group $B_N$, we have:

\index{bistochastic group}

\begin{theorem}
For the bistochastic group $B_N$ we have
$$\det(G_{kN})=N^{a_k}\prod_{|\lambda|\leq k/2}f_N(\lambda)^{\binom{k}{2|\lambda|}f^{2\lambda}}$$
where $a_k=\sum_{\pi\in P_{12}(k)}(2|\pi|-k)$, and $f_N(\lambda)=\prod_{(i,j)\in\lambda}(N+2j-i-2)$.
\end{theorem}

\begin{proof}
We recall that we have an isomorphism $B_N\simeq O_{N-1}$, given by $u=v+1$, where $u,v$ are the fundamental representations of $B_N,O_{N-1}$. But this gives:
$$Fix(u^{\otimes k})
=Fix\left((v+1)^{\otimes k}\right)
=Fix\left(\sum_{r=0}^k\binom{k}{r}v^{\otimes r}\right)$$

Now if we denote by $\det',f'$ the objects in Theorem 7.33, we obtain:
$$\det(G_{kN})
=N^{a_k}\prod_{r=1}^k{\rm det}'(G_{r,N-1})^{\binom{k}{r}}
=N^{a_k}\prod_{r=1}^k\left(\prod_{|\lambda|=r/2}f'_{N-1}(\lambda)^{f^{2\lambda}}\right)^{\binom{k}{r}}$$

Thus, we are led to the formula in the statement.
\end{proof}

Getting now to the free case, the simplest object here is $O_N^+$, coming from the category of noncrossing pairings $NC_2$. The associated Gram determinant, known as ``meander determinant'', was computed by Di Francesco in \cite{dif}, the result being as follows:

\index{meander determinant}
\index{Gram determinant}
\index{Di Francesco formula}

\begin{theorem}
The determinant of the Gram matrix for $O_N^+$ is given by
$$\det(G_{kN})=\prod_{r=1}^{[k/2]}P_r(N)^{d_{k/2,r}}$$
where $P_r$ are the Chebycheff polynomials, given by
$$P_0=1\quad,\quad 
P_1=X\quad,\quad 
P_{r+1}=XP_r-P_{r-1}$$
and $d_{kr}=f_{kr}-f_{k,r+1}$, with $f_{kr}$ being the following numbers, depending on $k,r\in\mathbb Z$,
$$f_{kr}=\binom{2k}{k-r}-\binom{2k}{k-r-1}$$
with the conventions $f_{kr}=0$ for $k\notin\mathbb Z$, and $\binom{p}{q}=0$ for $q<0$. 
\end{theorem}

\begin{proof}
This is indeed something quite standard, but long and technical. For more on this we refer to \cite{dif}, and to the more recent paper \cite{bcu} too, containing a short proof of this formula, based on advanced planar algebra technology of Jones \cite{jo3}.
\end{proof}

Regarding now $S_N^+$, coming from the category of all noncrossing partitions $NC$, we have here the following formula, also established by Di Francesco in \cite{dif}:

\index{meander determinant}
\index{Gram determinant}
\index{Chebycheff polynomials}
\index{Di Francesco formula}

\begin{theorem}
The determinant of the Gram matrix for $S_N^+$ is given by
$$\det(G_{kN})=(\sqrt{N})^{a_k}\prod_{r=1}^kP_r(\sqrt{N})^{d_{kr}}$$
where $P_r$, $d_{kr}$ are as before, and where $a_k=\sum_{\pi\in NC(k)}(2|\pi|-k)$.
\end{theorem}

\begin{proof}
This comes from Theorem 7.35, via the following standard bijection, obtained via fattening the pairings, and shrinking the partitions:
$$NC(k)\simeq NC_2(2k)$$

Indeed, if we denote by $G'$ the Gram matrix for $O_N^+$, we have the following formula, coming from the above bijection, with $D_{kN}=diag(N^{|\widetilde\pi|/2-k/4})$:
$$G_{kN}=D_{kN}G'_{2k,\sqrt{N}}D_{kN}$$

But with this formula in hand, the result follows from Theorem 7.35.
\end{proof}

As yet another version of Theorem 7.35, we can formulate as well:

\begin{theorem}
For the quantum group $B_N^+$ we have
$$\det(G_{kN})=N^{a_k}\prod_{r=1}^{[k/2]}P_r(N-1)^{\sum_{l=1}^{[k/2]}\binom{k}{2l}d_{lr}}$$
with $P_r$ and $d_{kr}$ being as before, and with $a_k=\sum_{\pi\in NC_{12}(k)}(2|\pi|-k)$.
\end{theorem}

\begin{proof}
The passage $O_N^+\to B_N^+$ is quite similar to the passage $O_N\to B_N$, and by using prime exponents for the various $O_N^+$-related objects, we obtain:
$$\det(G_{kN})
=N^{a_k}\prod_{l=1}^{[k/2]}{\rm det}'(G_{2l,N-1})^{\binom{k}{2l}}
=N^{a_k}\prod_{l=1}^{[k/2]}\left(\prod_{r=1}^lP_r(N-1)^{d_{lr}}\right)^{\binom{k}{2l}}$$

Together with Theorem 7.35, this gives the formula in the statement.
\end{proof}

Finally, in what regards the quantum group $H_N^+$, the result here is as follows:

\begin{theorem}
For the quantum group $H_N^+$ we have the formula
$$\det(G_{kN})=(\sqrt{N})^{a_k}\prod_{r=1}^{[k/2]}P_r(\sqrt{N})^{2d_{k/2,r}'}$$
with $d_{sr}'=f_{sr}'-f_{s,r+1}'$, where 
$f_{sr}'=\binom{3s}{s-r}-\binom{3s}{s-r-1}$ for $s\in\mathbb Z$, $f_{sr}'=0$ for $s\notin\mathbb Z$.
\end{theorem}

\begin{proof}
This can be viewed as yet another variation of Theorem 7.35, involving this time some colors on the strings, and for details here we refer to \cite{dif} and \cite{bcu}.
\end{proof}

Getting now to our original Question 7.31, we have formulae for all determinants there, and the remaining puzzle is that of putting everything together. See \cite{bcu}.

\section*{7e. Exercises}

Welcome to easiness, and no way back. As exercises, many things to be done:

\begin{exercise}
Learn about compact Lie groups, and their embeddings $G\subset U_N$.
\end{exercise}

\begin{exercise}
Learn more about the Haar measure, on various types of groups.
\end{exercise}

\begin{exercise}
Learn the various possible versions of Tannakian duality.
\end{exercise}

\begin{exercise}
Do the intertwining computations, for the main easy groups.
\end{exercise}

\begin{exercise}
Clarify everything that we said, in relation with uniformity.
\end{exercise}

\begin{exercise}
Further study the shifted complex normal laws $G_t^t$.
\end{exercise}

\begin{exercise}
Learn the proof of the Di Francesco determinant formula.
\end{exercise}

\begin{exercise}
Write the Gram determinants in terms of orthogonal polynomials.
\end{exercise}

As bonus exercise, have a look at free quantum groups. But we will be back to this.

\chapter{Random matrices}

\section*{8a. Matrices, laws}

With quantum groups discussed, a complementary discussion, regarding the random matrices, seems unavoidable. Indeed, these matrices are very interesting objects, providing as well models for the Wigner laws $\gamma_t$ and Marchenko-Pastur laws $\pi_t$, which are somewhat simpler. In fact, it is in this way that the laws $\gamma_t,\pi_t$ were discovered, some time ago, by Wigner in the 50s \cite{wig}, and by Marchenko-Pastur in the 60s \cite{mpa}.

\bigskip

This being said, I'm not exactly sure that it is the right time to do this, sort of a random matrix remake of what we did in chapter 7. I mean, this is a modest introduction to probability and normal variables, we have not hit yet the middle of the book, and still have many basic things to be learned. So, we will take a soft approach to this, with random matrices coming as a continuation of the basic material from chapter 5, regarding the complex normal variables, our official motivations being as follows:

\begin{fact}
A random matrix is a matrix with random variables as entries,
$$Z\in M_N(L^\infty(X))$$
with $X$ being a probability space, and the philosophy of these matrices is as follows:
\begin{enumerate}
\item When the usual random variables $f\in L^\infty(X)$ are not enough, for the modeling of your problem, the random matrices with $N\geq2$ are the answer.

\item With suitable definitions, each random matrix has a law, and in relation with your initial question, the problem is that of computing this law.

\item Typically most questions involve random matrices having i.i.d. entries, usually complex normal, $Z_{ij}\sim G_t$, taken up to some symmetry constraints.

\item Most questions involve random matrices coming in series, one for each $N\in\mathbb N$, and the main problem is that of computing the law in the $N\to\infty$ limit.
\end{enumerate}
\end{fact}

So, this was for the story and motivations, and lacking of course remains an explicit example of an application. However, and here comes the difficulty, although applications abound at the advanced level, to questions ranging from quantum physics to stock markets, at the beginner level, which is ours for the moment, there is none.

\bigskip

So, we will have to live with this, I mean just trust me, random matrices are not nonsense, and we will learn this later. At a philosophical level, however, in the lack of a concrete illustration for Fact 8.1, let us point out the obvious similarity with what we did in chapter 7. Indeed, what we have is, a bit as there, advanced beasts depending on $N\in\mathbb N$, constructed out of nothing or almost, using symmetry constraints, and having a law, that we are interested to compute in the $N\to\infty$ limit. And with the story being that we can reach in this way to $\gamma_t,\pi_t$, and many other interesting laws.

\bigskip

Let us record this as a poetical complement to Fact 8.1, as follows:

\begin{poetry}
In the same way as an ocean of independent Bernoulli variables can produce via sums $p_t,g_t$ and other interesting laws, an ocean of independent complex normal variables, arranged into matrices, can produce $\pi_t,\gamma_t$ and other interesting laws.
\end{poetry}

But probably enough talking, let us get to work, and we'll understand later. As a first job to be done, we must talk about the laws of random matrices $Z\in M_N(L^\infty(X))$. And here, in the simplest case, namely $X=\{.\}$, this can be done as follows:

\index{law of matrix}
\index{eigenvalues}
\index{Jordan block}

\begin{theorem}
Each matrix $A\in M_N(\mathbb C)$ has a law, appearing as the following abstract functional, on the algebra of noncommuting polynomials in $2$ variables,
$$\mu_A:\mathbb C<X,X^*>\to\mathbb C\quad,\quad P\to tr(P(A,A^*))$$
with $tr=Tr/N$ being the normalized trace. In the normal case, $AA^*=A^*A$, we have
$$\mu_A(P)=\int_\mathbb CP(z,\bar{z})d\mu_A(z)\quad,\quad\mu_A=\frac{1}{N}\sum_i\delta_{\lambda_i}$$
with $\lambda_1,\ldots,\lambda_N\in\mathbb C$ being the eigenvalues, so the law is a usual complex measure $\mu_A$.
\end{theorem}

\begin{proof}
This is something quite fundamental, based on the spectral theorem for the normal matrices, worth discussing in some detail, the idea being as follows:

\medskip

(1) To start with, as you surely know from linear algebra, at the advanced level the matrices $A\in M_N(\mathbb C)$ do not come alone, but rather in pairs $(A,A^*)$, with the adjoint matrix $A^*\in M_N(\mathbb C)$ being given by one of the following equivalent formulae:
$$<A^*x,y>=<x,Ay>\quad,\quad (A^*)_{ij}=\bar{A}_{ji}$$

To be more precise, the equivalence between the above two formulae is clear, with the second one corresponding to the first one at $x=e_j,y=e_i$. As for our claim that the matrices come in pairs $(A,A^*)$, this is something notorious, some illustrations being:

\medskip

-- Want to study the orthogonal projections? In this case, the adjoint matrices are what you need, the equation of the projections being $P^2=P^*=P$.

\medskip

-- Want to study isometries, also called unitaries? In this case, the adjoint matrices are again what you need, the equation of the unitaries being $U^*=U^{-1}$.

\medskip

-- And so on, and with the remark that restricting the attention to the real case $A\in M_N(\mathbb R)$ won't change the situation, cf. $P^2=P^t=P$ and $U^t=U^{-1}$ there.

\medskip

(2) Next, and getting now to probability, and more specifially to things that we briefly discussed in the beginning of chapter 5, let us make the following speculation:
$$A:M_N\to\mathbb C\ \implies\ E(A)=tr(A)$$

To be more precise, let us make the convention $C(M_N)=M_N(\mathbb C)$, allowing us to view the matrices $A\in M_N(\mathbb C)$ as variables $A:M_N\to\mathbb C$. Of course, that beast $M_N$ is something abstract, not a probability space in the usual sense. However, the algebra of random variables is there, according to $C(M_N)=M_N(\mathbb C)$, and with this being something rock-solid. And the expectation is there too, being of mass 1 as it should, according to the formula $E(A)=tr(A)$, with $tr=Tr/N$ being the normalized trace of matrices. 

\medskip

(3) Next, observe that the expectation, besides being of mass 1, is also positive, as it should. Indeed, in analogy with the fact that a variable $f:X\to\mathbb C$ is positive when it is real, $f:X\to\mathbb R$, and taking positive values, $f(x)\geq0$, we can say that a matrix $A\in M_N(\mathbb C)$ is positive when it is self-adjoint, $A=A^*$, and with positive eigenvalues, $\lambda_i\geq0$. But in this case the expectation is indeed positive, as shown by:
$$A\geq 0\ \implies\ E(A)=\frac{1}{N}\sum_i\lambda_i\geq0$$

(4) Alternatively, in case you are not convinced, in analogy with the fact that the positive variables $f:X\to\mathbb C$ are those of type $f=g\bar{g}$, we can say that the positive matrices $A\in M_N(\mathbb C)$ are those of type $A=BB^*$, and with this definition, we have:
$$A=BB^*\geq 0\ \implies\ E(A)=\frac{1}{N}\sum_{ij}|B_{ij}|^2\geq0$$

(5) Summarizing, one way or another, our speculation in (2) perfectly makes sense, in the hope that we agree on this. And in view of this, the problem appears, what is the law of a complex matrix $A\in M_N(\mathbb C)$, viewed as random variable $A:M_N\to\mathbb C$? 

\medskip

(6) And good question this is. In answer now, since the familiar formula $\mu_f=f_*\nu$ for usual variables $f:(X,\nu)\to\mathbb C$ won't apply to our setting, our space $X=M_N$ being not a probability space in the usual sense, we must trick, and recapture $\mu_A$ via moments.

\medskip

(7) So, our question becomes, what are the moments of a matrix $A\in M_N(\mathbb C)$, viewed as random variable $A:M_N\to\mathbb C$? And here, the answer is straightforward, namely:
$$M_k(A)=tr(A^k)$$

To be more precise, as usual when dealing with complex variables, we are using here colored integers $k=\circ\bullet\bullet\circ\ldots$ as exponents, with the conventions being:
$$A^\emptyset=1\quad,\quad A^\circ=A\quad,\quad A^\bullet=A^*\quad,\quad A^{kl}=A^kA^l$$

(8) Observe now that, unlike in the case of the usual variables $f:X\to\mathbb C$, which commute, our matrices $A$ generically won't commute with their adjoints $A^*$, and so:
$$A^{\circ\bullet}=AA^*\neq A^*A=A^{\bullet\circ}$$

In relation with this, observe however that we have the following computation:
$$M_{\circ\bullet}(A)=tr(AA^*)=tr(A^*A)=M_{\bullet\circ}(A)$$

However, such tricks won't work at higher order, 4 and more, as shown by:
\begin{eqnarray*}
AA^*-A^*A\neq0
&\implies&(AA^*-A^*A)^2>0\\
&\implies&AA^*AA^*-AA^*A^*A-A^*AAA^*+A^*AA^*A>0\\
&\implies&tr(AA^*AA^*-AA^*A^*A-A^*AAA^*+A^*AA^*A)>0\\
&\implies&tr(AA^*AA^*+A^*AA^*A)>tr(AA^*A^*A+A^*AAA^*)\\
&\implies&tr(AA^*AA^*)>tr(AAA^*A^*)
\end{eqnarray*}

(9) Actually, in relation with this, just in case, let us work out as well an explicit example. Playing the role of the bad guy will be a basic Jordan block, as follows:
$$J=\begin{pmatrix}0&1\\0&0\end{pmatrix}$$

We have the following formulae, which show that $J$ is indeed not normal:
$$JJ^*=\begin{pmatrix}1&0\\0&0\end{pmatrix}\quad,\quad 
J^*J=\begin{pmatrix}0&0\\0&1\end{pmatrix}$$

Now let us compute some order 4 moments. We first have:
$$tr(JJ^*JJ^*)
=tr\left(\begin{pmatrix}1&0\\0&0\end{pmatrix}\begin{pmatrix}1&0\\0&0\end{pmatrix}\right)
=tr\begin{pmatrix}1&0\\0&0\end{pmatrix}
=\frac{1}{2}$$

On the other hand, we have as well the following formula:
$$tr(JJJ^*J^*)
=tr\left(\begin{pmatrix}0&0\\0&0\end{pmatrix}
\begin{pmatrix}0&0\\0&0\end{pmatrix}
\right)
=tr\begin{pmatrix}0&0\\0&0\end{pmatrix}
=0$$

Thus, for our Jordan block $J$, there is absolutely no way for $M_{\circ\bullet\circ\bullet}=M_{\circ\circ\bullet\bullet}$ to hold.

\medskip

(10) Summarizing, we have to deal with a massive quantity of moments, basically coming from all colored integers $k=\circ\bullet\bullet\circ\ldots$ as exponents, without many rules for these. Alternatively, by linearizing, the moment combinatorics of $A$ is captured by the following abstract functional, on the algebra of noncommuting polynomials in $2$ variables:
$$\mu_A:\mathbb C<X,X^*>\to\mathbb C\quad,\quad P\to tr(P(A,A^*))$$

(11) And with this discussed, what is next? Well, giving up, I mean enough talking, let's just call the above abstract functional law of $A$, and end of the story.

\medskip

(12) Now with the foundations of our theory laid, let us compute some matrix laws. And here, in view of (8), we definitely want to avoid the non-normal matrices, $AA^*\neq A^*A$, because these will produce a massive amount of moments, which cannot be encoded by a reasonable measure-type object $\mu_a$. That is, for the generic non-normal matrices, the law $\mu_a$ will stay something abstract, as originally constructed in (10) above.

\medskip

(13) So, let us get into the normal case, $AA^*=A^*A$. Here the situation is quite similar to that of the usual variables $f:X\to\mathbb C$, due to the fact that the algebra of random variables that we are truly interested in, namely $<A,A^*>$, is commutative. As a first consequence of this observation, the law constructed in (10) becomes something simpler, namely the following functional, on the algebra of usual polynomials in $2$ variables:
$$\mu_A:\mathbb C[X,X^*]\to\mathbb C\quad,\quad P\to tr(P(A,A^*))$$

(14) So, what is this latter law? In answer, we can use here the spectral theorem for the normal matrices, telling us that we can diagonalize our matrix as $A=UDU^*$, with $U\in U_N$, and $D$ diagonal. Thus, our variable and its conjugate are given by:
$$A=UDU^*\quad,\quad A^*=UD^*U^*$$

But with this in hand, we have the following computation, for any polynomial $P$:
\begin{eqnarray*}
\mu_A(P)
&=&tr(P(A,A^*))\\
&=&tr(P(UDU^*,UD^*U^*))\\
&=&tr(UP(D,D^*)U^*)\\
&=&tr(P(D,D^*))\\
&=&\mu_D(P)
\end{eqnarray*}

(15) Thus, we are left with computing the law of diagonal matrices. But here, given a diagonal matrix $D=diag(\lambda_1,\ldots,\lambda_N)$, we have the following computation:
\begin{eqnarray*}
\mu_D(P)
&=&tr(P(D,D^*))\\
&=&tr(diag(P(\lambda_i,\bar{\lambda}_i)))\\
&=&\frac{1}{N}\sum_iP(\lambda_i,\bar{\lambda}_i)
\end{eqnarray*}

Thus, we are led to the conclusion that the law of a diagonal matrix $D$ is the average of the Dirac masses at the eigenvalues, according to the following formula:
$$\mu_D(P)=\int_\mathbb CP(z,\bar{z})d\mu_D(z)\quad,\quad\mu_D=\frac{1}{N}\sum_i\delta_{\lambda_i}$$

(16) Now by getting back to our normal matrices from (14), we can see that same final conclusion still holds there, namely that the law is the average of the Dirac masses at the eigenvalues. And with this, theorem proved, and end of our discussion.
\end{proof}

\section*{8b. Spectral measures}

With the above discussed, corresponding to the case $X=\{.\}$ in the random matrix context, $Z\in M_N(L^\infty(X))$, let us turn now to the general case, where $X$ is arbitrary. And here, things are obviously quite tricky, because we would need a spectral theorem for the matrices $Z\in M_N(L^\infty(X))$, going well beyond the linear algebra that we know.

\bigskip

Fortunately, there is a clever approach to such questions, as follows:

\begin{principle}
The random matrix algebras are best viewed as operator algebras,
$$M_N(L^\infty(X))=B(H)\quad,\quad H=\mathbb C^N\otimes L^2(X)$$
and the random matrices themselves, as operators $Z\in B(H)$.
\end{principle}

In short, what I propose here is to generalize the linear algebra that we know, including the spectral theorem for the normal matrices, to the case of infinite dimensions, over an arbitrary complex Hilbert space $H$. And afterwards, armed with this knowledge, simply take $H=\mathbb C^N\otimes L^2(X)$, and easily deal with the random matrices $Z\in B(H)$.

\bigskip

Getting started now, we already know about Hilbert spaces from chapter 3, but a few brief reminders won't hurt. At the beginning of everything, we have:

\index{Hilbert space}
\index{scalar product}

\begin{definition}
A complex Hilbert space is a complex vector space $H$ with a scalar product $<x,y>$, taken linear at left and antilinear at right,
$$<\lambda x,y>=\lambda<x,y>\quad,\quad <x,\lambda y>=\bar{\lambda}<x,y>$$
which is complete with respect to corresponding norm
$$||x||=\sqrt{<x,x>}$$
in the sense that any sequence $\{x_n\}$ which is a Cauchy sequence, having the property $||x_n-x_m||\to0$ with $n,m\to\infty$, has a limit, $x_n\to x$.
\end{definition}

Here the fact that $||x||=\sqrt{<x,x>}$ is indeed a norm, satisfying $||x+y||\leq||x||+||y||$, follows from Cauchy-Schwarz, $|\!<x,y>\!|\leq||x||\cdot||y||$, which itself follows as in the $H=\mathbb C^N$ case, by looking at the discriminant of the following degree 2 polynomial:
$$f(t)=||wx+ty||^2$$

As for our convention for the scalar products, written $<x,y>$ and taken linear at left, this is the so-called feline convention. As basic examples now, we have:

\index{square-summable}
\index{Cauchy-Schwarz}

\begin{theorem}
The following are Hilbert spaces:
\begin{enumerate}
\item The space $H=\mathbb C^N$, with scalar product $<x,y>=\sum_ix_i\bar{y}_i$.

\item More generally, $H=l^2(I)$, with scalar product $<x,y>=\sum_ix_i\bar{y}_i$.

\item Even more generally, $H=L^2(X)$, with $<f,g>=\int_Xf(x)\overline{g(x)}dx$.
\end{enumerate}
However, there is only one separable Hilbert space, namely $l^2(\mathbb N)$, or $L^2[0,1]$.
\end{theorem}

\begin{proof}
The first three assertions are very standard, with (1) being clear, (2) being an easy remake of (1), and then (3) being an easy remake of (2). The subtlety is at the end, with separable meaning to have a countable orthonormal basis $\{e_i\}_{i\in I}$, and with the space in question being therefore $H=l^2(I)\simeq l^2(\mathbb N)$. Next, by Weierstrass we have as well as basic example here $H=L^2[0,1]$, and more generally, what happens is that $L^2(X)$ is separable, provided that $X$ itself is separable. See e.g. Rudin \cite{rud}.
\end{proof}

Moving ahead, now that we know what our vector spaces are, we can talk about matrices with respect to them. And the situation here is as follows:

\index{infinite matrix}
\index{linear operator}
\index{bounded operator}
\index{continuous operator}

\begin{theorem}
Given a Hilbert space $H$, consider the linear operators $T:H\to H$, and for each such operator define its norm by the following formula:
$$||T||=\sup_{||x||=1}||Tx||$$
The operators which are bounded, $||T||<\infty$, form then a complex algebra $B(H)$, which is complete with respect to $||.||$. When $H$ comes with a basis $\{e_i\}_{i\in I}$, we have
$$B(H)\subset M_I(\mathbb C)$$
with the correspondence $T\to M$ coming via the usual linear algebra formulae, namely:
$$T(x)=Mx\quad,\quad M_{ij}=<Te_j,e_i>$$
In infinite dimensions, the inclusion $B(H)\subset M_I(\mathbb C)$ is not an equality.
\end{theorem}

\begin{proof}
This is something straightforward, the idea being as follows:

\medskip

(1) The fact that we have indeed an algebra, satisfying the product condition in the statement, follows from the following estimates, which are all elementary:
$$||S+T||\leq||S||+||T||\quad,\quad 
||\lambda T||=|\lambda|\cdot||T||\quad,\quad 
||ST||\leq||S||\cdot||T||$$

(2) Regarding now the completness assertion, if $\{T_n\}\subset B(H)$ is Cauchy then $\{T_nx\}$ is Cauchy for any $x\in H$, so we can define the limit $T=\lim_{n\to\infty}T_n$ by setting:
$$Tx=\lim_{n\to\infty}T_nx$$

Let us first check that the application $x\to Tx$ is linear. We have:
\begin{eqnarray*}
T(x+y)
&=&\lim_{n\to\infty}T_n(x+y)\\
&=&\lim_{n\to\infty}T_n(x)+T_n(y)\\
&=&\lim_{n\to\infty}T_n(x)+\lim_{n\to\infty}T_n(y)\\
&=&T(x)+T(y)
\end{eqnarray*}

Similarly, we have $T(\lambda x)=\lambda T(x)$, and we conclude that $x\to Tx$ is linear.

\medskip

(3) With this done, it remains to prove now that we have $T\in B(H)$, and that $T_n\to T$ in norm. For this purpose, observe that we have:
\begin{eqnarray*}
||T_n-T_m||\leq\varepsilon\ ,\ \forall n,m\geq N
&\implies&||T_nx-T_mx||\leq\varepsilon\ ,\ \forall||x||=1\ ,\ \forall n,m\geq N\\
&\implies&||T_nx-Tx||\leq\varepsilon\ ,\ \forall||x||=1\ ,\ \forall n\geq N\\
&\implies&||T_Nx-Tx||\leq\varepsilon\ ,\ \forall||x||=1\\
&\implies&||T_N-T||\leq\varepsilon
\end{eqnarray*}

But this gives both $T\in B(H)$, and $T_N\to T$ in norm, and we are done.

\medskip

(4) Regarding the embedding, the correspondence $T\to M$ in the statement is indeed linear, and its kernel is $\{0\}$, so we have indeed an embedding as follows, as claimed:
$$B(H)\subset M_I(\mathbb C)$$

In finite dimensions we have an isomorphism, because any matrix $M\in M_N(\mathbb C)$ determines a linear operator $T:\mathbb C^N\to\mathbb C^N$, given by the following formula:
$$<Te_j,e_i>=M_{ij}$$

However, in infinite dimensions we have matrices not producing operators, as for instance the all-one matrix, so the embedding $B(H)\subset M_I(\mathbb C)$ is not an isomorphism.
\end{proof}

As a second and last basic result regarding the operators, we will need:

\index{adjoint operator}
\index{adjoint matrix}

\begin{theorem}
Each operator $T\in B(H)$ has an adjoint $T^*\in B(H)$, given by: 
$$<Tx,y>=<x,T^*y>$$
The operation $T\to T^*$ is antilinear, antimultiplicative, involutive, and satisfies:
$$||T||=||T^*||\quad,\quad ||TT^*||=||T||^2$$
When $H$ comes with a basis $\{e_i\}_{i\in I}$, the operation $T\to T^*$ corresponds to
$$(M^*)_{ij}=\overline{M}_{ji}$$ 
at the level of the associated matrices $M\in M_I(\mathbb C)$.
\end{theorem}

\begin{proof}
This is standard too, and can be proved in 3 steps, as follows:

\medskip

(1) The existence of the adjoint operator $T^*$, given by the formula in the statement, comes from the fact that the function $\varphi(x)=<Tx,y>$ being a linear map $H\to\mathbb C$, we must have a formula as follows, for a certain vector $T^*y\in H$:
$$\varphi(x)=<x,T^*y>$$

Moreover, since this vector is unique, $T^*$ is unique too, and we have as well:
$$(S+T)^*=S^*+T^*\quad,\quad
(\lambda T)^*=\bar{\lambda}T^*\quad,\quad 
(ST)^*=T^*S^*\quad,\quad 
(T^*)^*=T$$

Observe that we have indeed $T^*\in B(H)$, as shown by the following computation:
\begin{eqnarray*}
||T||
&=&\sup_{||x||=1}\sup_{||y||=1}<Tx,y>\\
&=&\sup_{||y||=1}\sup_{||x||=1}<x,T^*y>\\
&=&||T^*||
\end{eqnarray*}

(2) Regarding $||TT^*||=||T||^2$, which is a key formula, observe that we have:
$$||TT^*||
\leq||T||\cdot||T^*||
=||T||^2$$

In the other sense now, observe that we have the following estimate:
\begin{eqnarray*}
||T||^2
&=&\sup_{||x||=1}|<Tx,Tx>|\\
&=&\sup_{||x||=1}|<x,T^*Tx>|\\
&\leq&||T^*T||
\end{eqnarray*}

But by replacing $T\to T^*$ we obtain from this $||T||^2\leq||TT^*||$, as desired.

\medskip

(3) Finally, when $H$ comes with a basis, the formula $<Tx,y>=<x,T^*y>$ applied with $x=e_i$, $y=e_j$ translates into the formula $(M^*)_{ij}=\overline{M}_{ji}$, as desired.
\end{proof}

Getting now a bit abstract, we need to talk about operator algebras, as to suitably cover the random matrix algebras $M_N(L^\infty(X))$. Which can be something quite tricky, but fortunately, we have the following clever definition, due to Gelfand:

\begin{definition}
An abstract operator algebra, or $C^*$-algebra, is a complex algebra $A$ having a norm $||.||$ and an involution $*$, subject to the following conditions:
\begin{enumerate}
\item $A$ is closed with respect to the norm.

\item We have $||aa^*||=||a||^2$, for any $a\in A$.
\end{enumerate}
\end{definition}

In other words, what we did here is to axiomatize the abstract properties of the operator algebras $A\subset B(H)$, chosen closed under the norm and under the adjoint operation, both very natural conditions, without reference to the ambient Hilbert space $H$.

\bigskip

As very basic examples, we have the usual matrix algebras $M_N(\mathbb C)$, with the norm and the involution being the usual matrix norm and involution, given by:
$$||A||=\sup_{||x||=1}||Ax||\quad,\quad 
(A^*)_{ij}=\overline{A}_{ji}$$

Some other basic examples are the algebras $L^\infty(X)$ of essentially bounded functions $f:X\to\mathbb C$ on a measured space $X$, with the usual norm and involution, namely:
$$||f||=\sup_{x\in X}|f(x)|\quad,\quad 
f^*(x)=\overline{f(x)}$$

We can put these two basic classes of examples together, as follows:

\index{random matrix algebra}

\begin{theorem}
The random matrix algebras $A=M_N(L^\infty(X))$ are $C^*$-algebras, with their usual norm and involution, given by:
$$||Z||=\sup_{x\in X}||Z_x||\quad,\quad 
(Z^*)_{ij}=\overline{Z}_{ij}$$
These algebras generalize both the algebras $M_N(\mathbb C)$, and the algebras $L^\infty(X)$.
\end{theorem}

\begin{proof}
The fact that the $C^*$-algebra axioms are satisfied is clear from definitions. As for the last assertion, this follows by taking $X=\{.\}$ and $N=1$, respectively.
\end{proof}

The above result is quite interesting, philosophically, showing that the Gelfand axioms in Definition 8.9 favor the random matrix algebras. Moving on, we can in fact say more about the operator algebra nature of random matrix algebras, as follows:

\begin{theorem}
Any random variable algebra is an operator algebra, as follows:
$$L^\infty(X)\subset B(L^2(X))\quad,\quad 
f\to(g\to fg)$$
More generally, any random matrix algebra is an operator algebra, as follows,
$$M_N(L^\infty(X))\subset B\left(\mathbb C^N\otimes L^2(X)\right)$$
with the embedding being the above one, tensored with the identity.
\end{theorem}

\begin{proof}
We have two assertions to be proved, the idea being as follows:

\medskip

(1) Given $f\in L^\infty(X)$, consider the following operator, acting on $H=L^2(X)$:
$$T_f(g)=fg$$

Observe that $T_f$ is indeed well-defined, and bounded as well, because:
$$||fg||_2
=\sqrt{\int_X|f(x)|^2|g(x)|^2d\mu(x)}
\leq||f||_\infty||g||_2$$

The application $f\to T_f$ being linear, involutive, continuous, and injective as well, we obtain in this way a $C^*$-algebra embedding $L^\infty(X)\subset B(H)$, as desired.

\medskip

(2) Regarding the second assertion, this is best viewed in the following way:
\begin{eqnarray*}
M_N(L^\infty(X))
&=&M_N(\mathbb C)\otimes L^\infty(X)\\
&\subset&M_N(\mathbb C)\otimes B(L^2(X))\\
&=&B\left(\mathbb C^N\otimes L^2(X)\right)
\end{eqnarray*}

Here we have used (1), and some standard tensor product identifications.
\end{proof}

Our purpose in what follows will be to develop the spectral theory of the $C^*$-algebras, and in particular that of the random matrix algebras $M_N(L^\infty(X))$ that we are interested in, our main objective being that of talking about spectral measures, in the normal case, in analogy with what we know about the usual matrices. Let us start with:

\index{spectrum}
\index{functional calculus}
\index{rational function}
\index{self-adjoint}
\index{unitary}
\index{spectral radius}
\index{normal element}

\begin{theorem}
Given an element $a\in A$ of a $C^*$-algebra, define its spectrum as:
$$\sigma(a)=\left\{\lambda\in\mathbb C\Big|a-\lambda\notin A^{-1}\right\}$$
The following spectral theory results hold, exactly as in the $A=B(H)$ case:
\begin{enumerate}
\item We have $\sigma(ab)\cup\{0\}=\sigma(ba)\cup\{0\}$.

\item We have $\sigma(f(a))=f(\sigma(a))$, for any $f\in\mathbb C(X)$ having poles outside $\sigma(a)$.

\item The spectrum $\sigma(a)$ is compact, non-empty, and contained in $D_0(||a||)$.

\item The spectra of unitaries $(u^*=u^{-1})$ and self-adjoints $(a=a^*)$ are on $\mathbb T,\mathbb R$.

\item The spectral radius of normal elements $(aa^*=a^*a)$ is given by $\rho(a)=||a||$.
\end{enumerate}
In addition, assuming $a\in A\subset B$, the spectra of $a$ with respect to $A$ and to $B$ coincide.
\end{theorem}

\begin{proof}
Here the assertions (1-5), which are of course formulated a bit informally, are well-known for the full operator algebra $A=B(H)$, and the proof in general is similar:

\medskip

(1) Assuming that $1-ab$ is invertible, with inverse $c$, we have $abc=cab=c-1$, and it follows that $1-ba$ is invertible too, with inverse $1+bca$. Thus $\sigma(ab),\sigma(ba)$ agree on $1\in\mathbb C$, and by linearity, it follows that $\sigma(ab),\sigma(ba)$ agree on any point $\lambda\in\mathbb C^*$.

\medskip

(2) The formula $\sigma(f(a))=f(\sigma(a))$ is clear for polynomials, $f\in\mathbb C[X]$, by factorizing $f-\lambda$, with $\lambda\in\mathbb C$. Then, the extension to the rational functions is straightforward, because $P(a)/Q(a)-\lambda$ is invertible precisely when $P(a)-\lambda Q(a)$ is.

\medskip

(3) By using $1/(1-b)=1+b+b^2+\ldots$ for $||b||<1$ we obtain that $a-\lambda$ is invertible for $|\lambda|>||a||$, and so $\sigma(a)\subset D_0(||a||)$. It is also clear that $\sigma(a)$ is closed, so what we have is a compact set. Finally, assuming $\sigma(a)=\emptyset$ the function $f(\lambda)=\varphi((a-\lambda)^{-1})$ is well-defined, for any $\varphi\in A^*$, and by Liouville we get $f=0$, contradiction.

\medskip

(4) Assuming $u^*=u^{-1}$ we have $||u||=1$, and so $\sigma(u)\subset D_0(1)$. But with $f(z)=z^{-1}$ we obtain via (2) that we have as well $\sigma(u)\subset f(D_0(1))$, and this gives $\sigma(u)\subset\mathbb T$. As for the result regarding the self-adjoints, this can be obtained from the result for the unitaries, by using (2) with functions of type $f(z)=(z+it)/(z-it)$, with $t\in\mathbb R$.

\medskip

(5) It is routine to check, by integrating quantities of type $z^n/(z-a)$ over circles centered at the origin, and estimating, that the spectral radius is given by $\rho(a)=\lim||a^n||^{1/n}$. But in the self-adjoint case, $a=a^*$, this gives $\rho(a)=||a||$, by using exponents of type $n=2^k$, and then the extension to the general normal case is straightforward.

\medskip 

(6) Regarding now the last assertion, the inclusion $\sigma_B(a)\subset\sigma_A(a)$ is clear. For the converse, assume $a-\lambda\in B^{-1}$, and set $b=(a-\lambda )^*(a-\lambda )$. We have then:
$$\sigma_A(b)-\sigma_B(b)=\left\{\mu\in\mathbb C-\sigma_B(b)\Big|(b-\mu)^{-1}\in B-A\right\}$$

Thus this difference in an open subset of $\mathbb C$. On the other hand $b$ being self-adjoint, its two spectra are both real, and so is their difference. Thus the two spectra of $b$ are equal, and in particular $b$ is invertible in $A$, and so $a-\lambda\in A^{-1}$, as desired.
\end{proof}

We can now prove a key result about the operator algebras, as follows:

\index{Gelfand theorem}
\index{commutative algebra}

\begin{theorem}[Gelfand]
If $X$ is a compact space,  the algebra $C(X)$ of continuous functions on it $f:X\to\mathbb C$ is a $C^*$-algebra, with usual norm and involution, namely:
$$||f||=\sup_{x\in X}|f(x)|\quad,\quad 
f^*(x)=\overline{f(x)}$$
Conversely, any commutative $C^*$-algebra is of this form, $A=C(X)$, with 
$$X=\Big\{\chi:A\to\mathbb C\ ,\ {\rm normed\ algebra\ character}\Big\}$$
with topology making continuous the evaluation maps $ev_a:\chi\to\chi(a)$.
\end{theorem}

\begin{proof}
Given a commutative $C^*$-algebra $A$, let us define $X$ as in the statement. Then $X$ is compact, and $a\to ev_a$ is a morphism of algebras, as follows:
$$ev:A\to C(X)$$

We first prove that $ev$ is involutive. We use the following formula, which is similar to the $z=Re(z)+iIm(z)$ decomposition formula for usual complex numbers:
$$a=\frac{a+a^*}{2}+i\cdot\frac{a-a^*}{2i}$$

Thus it is enough to prove $ev_{a^*}=ev_a^*$ for the self-adjoint elements $a$. But this is the same as proving that $a=a^*$ implies that $ev_a$ is a real function, which is in turn true, by Theorem 8.12, because $ev_a(\chi)=\chi(a)$ is an element of $\sigma(a)$, contained in $\mathbb R$. Next, since $A$ is commutative, each element is normal, so $ev$ is isometric, as shown by:
$$||ev_a||
=\rho(a)
=||a||$$

It remains to prove that $ev$ is surjective. But this follows from the Stone-Weierstrass theorem, because $ev(A)$ is a closed subalgebra of $C(X)$, which separates the points.
\end{proof}

As a main consequence of the Gelfand theorem, we have:

\begin{theorem}
For any normal element $a\in A$ we have $<a>=C(\sigma(a))$. In addition, given a function $f\in C(\sigma(a))$, we can apply it to $a$, and we have
$$\sigma(f(a))=f(\sigma(a))$$
which generalizes the previous rational calculus formula, in the normal case.
\end{theorem}

\begin{proof}
Since $a$ is normal, the $C^*$-algebra $<a>$ that is generates is commutative, so if we denote by $X$ the space of the characters $\chi:<a>\to\mathbb C$, we have $<a>=C(X)$. Now since the map $X\to\sigma(a)$ given by evaluation at $a$ is bijective, we obtain $<a>=C(\sigma(a))$. Thus, we are dealing here with usual functions, and this gives all the assertions.
\end{proof}

In order to get now towards noncommutative probability, we first have to develop the theory of positive elements, and linear forms. First, we have the following result:

\index{positive element}

\begin{proposition}
For an element $a\in A$, the following are equivalent:
\begin{enumerate}
\item $a$ is positive, in the sense that $\sigma(a)\subset[0,\infty)$.

\item $a=b^2$, for some $b\in A$ satisfying $b=b^*$.

\item $a=cc^*$, for some $c\in A$.
\end{enumerate}
\end{proposition}

\begin{proof}
This is something very standard, as follows:

\medskip

$(1)\implies(2)$ Observe first that $\sigma(a)\subset\mathbb R$ implies $a=a^*$. Thus the algebra $<a>$ is commutative, and by using Theorem 8.14, we can set $b=\sqrt{a}$.

\medskip

$(2)\implies(3)$ This is trivial, because we can simply set $c=b$. 

\medskip

$(3)\implies(1)$ We can proceed here by contradiction. Indeed, by multiplying $c$ by a suitable element of the algebra $<cc^*>$, we are led to the existence of an element $d\neq0$ satisfying $-dd^*\geq0$. By writing now $d=x+iy$ with $x=x^*,y=y^*$ we have:
$$dd^*+d^*d
=2(x^2+y^2)
\geq0$$

Thus $d^*d\geq0$, which is easily seen to contradict the condition $-dd^*\geq0$.
\end{proof}

We can talk as well about positive linear forms, as follows:

\index{positive linear form}

\begin{definition}
Consider a linear map $\varphi:A\to\mathbb C$.
\begin{enumerate}
\item $\varphi$ is called positive when $a\geq0\implies\varphi(a)\geq0$.

\item $\varphi$ is called faithful and positive when $a\geq0,a\neq0\implies\varphi(a)>0$.
\end{enumerate}
\end{definition}

In the commutative case, $A=C(X)$, the positive linear forms appear as follows, with $\mu$ being positive, and strictly positive if we want $\varphi$ to be faithful and positive:
$$\varphi(f)=\int_Xf(x)d\mu(x)$$

In general, the positive linear forms can be thought of as being integration functionals with respect to some underlying ``positive measures''. Based on this, we can formulate:

\index{colored moments}

\begin{definition}
Let $A$ be a $C^*$-algebra, given with a positive trace $tr:A\to\mathbb C$.
\begin{enumerate}
\item The elements $a\in A$ are called random variables.

\item The moments of such a variable are the numbers $M_k(a)=tr(a^k)$.

\item The law of such a variable is the functional $\mu_a:P\to tr(P(a))$.
\end{enumerate}
\end{definition}

To be more precise, here the exponent $k=\circ\bullet\bullet\circ\ldots$ is by definition a colored integer, and the powers $a^k$ are defined by the following formulae, and multiplicativity: 
$$a^\emptyset=1\quad,\quad
a^\circ=a\quad,\quad
a^\bullet=a^*$$ 

As for the polynomial $P$, this is a noncommuting $*$-polynomial in one variable: 
$$P\in\mathbb C<X,X^*>$$

At the level of the general theory, we have the following key result, extending the various results from linear algebra, regarding the self-adjoint and normal matrices:

\index{normal element}
\index{spectral measure}

\begin{theorem}
Let $A$ be a $C^*$-algebra, with a trace $tr$, and consider an element $a\in A$ which is normal, in the sense that $aa^*=a^*a$.
\begin{enumerate}
\item $\mu_a$ is a complex probability measure, satisfying $supp(\mu_a)\subset\sigma(a)$.

\item In the self-adjoint case, $a=a^*$, this measure $\mu_a$ is real.

\item Assuming that $tr$ is faithful, we have $supp(\mu_a)=\sigma(a)$.
\end{enumerate}
\end{theorem}

\begin{proof}
According to Theorem 8.14, we have an identification as follows: 
$$<a>=C(\sigma(a))$$

Thus the functional $f(a)\to tr(f(a))$ can be regarded as an integration functional on the algebra $C(\sigma(a))$, and by the Riesz theorem, this latter functional must come from a probability measure $\mu_a$ on the spectrum $\sigma(a)$, in the sense that we must have:
$$tr(f(a))=\int_{\sigma(a)}f(z)d\mu_a(z)$$

We are therefore led to the various conclusions in the statement.
\end{proof}

Getting back now to the random matrices, as a main application, we have:

\begin{theorem}
Given a random matrix $Z\in M_N(L^\infty(X))$, assumed to be normal, $ZZ^*=Z^*Z$, its law, when restricted to the usual polynomials in two variables,
$$\mu_Z:\mathbb C[X,X^*]\to\mathbb C\quad,\quad 
P\to\frac{1}{N}\int_Xtr(P(Z))$$
must come from a probability measure on the spectrum $\sigma(Z)\subset\mathbb C$, as follows:
$$\mu_Z(P)=\int_{\sigma(T)}P(x)d\mu_Z(x)$$
We agree to use the symbol $\mu_Z$ for all these notions.
\end{theorem}

\begin{proof}
This follows indeed from what we know from Theorem 8.18, applied to the normal element $a=Z$, belonging to the $C^*$-algebra $A=M_N(L^\infty(X))$. 
\end{proof}

\section*{8c. Wigner matrices}

We recall that a random matrix algebra is an algebra of type $A=M_N(L^\infty(X))$, and that we are interested in the computation of the laws of the operators $Z\in A$, called random matrices. Regarding the precise classes of random matrices that we are interested in, first we have the complex Gaussian matrices, which are constructed as follows:

\index{Gaussian matrix}

\begin{definition}
A complex Gaussian matrix is a random matrix of type
$$Z\in M_N(L^\infty(X))$$
which has i.i.d. complex normal entries.
\end{definition}

We will see that the above matrices have an interesting, and ``central'' combinatorics, among all kinds of random matrices, with the study of the other random matrices being usually obtained as a modification of the study of the Gaussian matrices.

\bigskip

As a somewhat surprising remark, using real normal variables in Definition 8.20, instead of the complex ones appearing there, leads nowhere. The correct real versions of the Gaussian matrices are the Wigner random matrices, constructed as follows: 

\index{Wigner matrix}

\begin{definition}
A Wigner matrix is a random matrix of type
$$Z\in M_N(L^\infty(X))$$
which has i.i.d. complex normal entries, up to the constraint $Z=Z^*$.
\end{definition}

In other words, a Wigner matrix must be as follows, with the diagonal entries being real normal variables, $a_i\sim g_t$, for some $t>0$, the upper diagonal entries being complex normal variables, $b_{ij}\sim G_t$, the lower diagonal entries being the conjugates of the upper diagonal entries, as indicated, and with all the variables $a_i,b_{ij}$ being independent: 
$$Z=\begin{pmatrix}
a_1&b_{12}&\ldots&\ldots&b_{1N}\\
\bar{b}_{12}&a_2&\ddots&&\vdots\\
\vdots&\ddots&\ddots&\ddots&\vdots\\
\vdots&&\ddots&a_{N-1}&b_{N-1,N}\\
\bar{b}_{1N}&\ldots&\ldots&\bar{b}_{N-1,N}&a_N
\end{pmatrix}$$

As a comment here, for many concrete applications the Wigner matrices are in fact the central objects in random matrix theory, and in particular, they are often more important than the Gaussian matrices. In fact, these are the random matrices which were first considered and investigated, a long time ago, by Wigner himself \cite{wig}.

\bigskip

Finally, we will be interested as well in the complex Wishart matrices, which are the positive versions of the above random matrices, constructed as follows: 

\index{Wishart matrix}

\begin{definition}
A complex Wishart matrix is a random matrix of type
$$Z=YY^*\in M_N(L^\infty(X))$$
with $Y$ being a complex Gaussian matrix.
\end{definition}

As before with the Gaussian and Wigner matrices, there are many possible comments that can be made here, of technical or historical nature. First, using in the above real Gaussian variables instead of complex ones leads to a less interesting combinatorics. Also, these matrices were introduced and studied by Marchenko-Pastur not long after Wigner, and so historically came second \cite{mpa}. Finally, in what regards their combinatorics and applications, these matrices quite often come first, before both the Gaussian and the Wigner ones, with all this being of course a matter of knowledge and taste.

\bigskip

Summarizing, we have three main types of random matrices, which can be somehow designated as ``complex'', ``real'' and ``positive'', and that we will study in what follows. Let us also mention that there are many other interesting classes of random matrices, usually appearing as modifications of the above. More on these later.

\bigskip

Getting to work now, we first have the following result:

\begin{theorem}
Given a sequence of Gaussian random matrices
$$Z_N\in M_N(L^\infty(X))$$
having independent $G_t$ variables as entries, for some fixed $t>0$, we have
$$M_k\left(\frac{Z_N}{\sqrt{N}}\right)\simeq t^{|k|/2}|\mathcal{NC}_2(k)|$$
for any colored integer $k=\circ\bullet\bullet\circ\ldots\,$, in the $N\to\infty$ limit.
\end{theorem}

\begin{proof}
This is something standard, which can be done as follows:

\medskip

(1) We fix $N\in\mathbb N$, and we let $Z=Z_N$. Let us first compute the trace of $Z^k$. With $k=k_1\ldots k_s$, and with the convention $(ij)^\circ=ij,(ij)^\bullet=ji$, we have:
\begin{eqnarray*}
Tr(Z^k)
&=&Tr(Z^{k_1}\ldots Z^{k_s})\\
&=&\sum_{i_1=1}^N\ldots\sum_{i_s=1}^N(Z^{k_1})_{i_1i_2}(Z^{k_2})_{i_2i_3}\ldots(Z^{k_s})_{i_si_1}\\
&=&\sum_{i_1=1}^N\ldots\sum_{i_s=1}^N(Z_{(i_1i_2)^{k_1}})^{k_1}(Z_{(i_2i_3)^{k_2}})^{k_2}\ldots(Z_{(i_si_1)^{k_s}})^{k_s}
\end{eqnarray*}

(2) Next, we rescale our variable $Z$ by a $\sqrt{N}$ factor, as in the statement, and we also replace the usual trace by its normalized version, $tr=Tr/N$. Our formula becomes:
$$tr\left(\left(\frac{Z}{\sqrt{N}}\right)^k\right)=\frac{1}{N^{s/2+1}}\sum_{i_1=1}^N\ldots\sum_{i_s=1}^N(Z_{(i_1i_2)^{k_1}})^{k_1}(Z_{(i_2i_3)^{k_2}})^{k_2}\ldots(Z_{(i_si_1)^{k_s}})^{k_s}$$

Thus, the moment that we are interested in is given by:
$$M_k\left(\frac{Z}{\sqrt{N}}\right)=\frac{1}{N^{s/2+1}}\sum_{i_1=1}^N\ldots\sum_{i_s=1}^N\int_X(Z_{(i_1i_2)^{k_1}})^{k_1}(Z_{(i_2i_3)^{k_2}})^{k_2}\ldots(Z_{(i_si_1)^{k_s}})^{k_s}$$

(3) Let us apply now the Wick formula, from chapter 5. We conclude that the moment that we are interested in is given by the following formula:
\begin{eqnarray*}
&&M_k\left(\frac{Z}{\sqrt{N}}\right)\\
&=&\frac{t^{s/2}}{N^{s/2+1}}\sum_{i_1=1}^N\ldots\sum_{i_s=1}^N\#\left\{\pi\in\mathcal P_2(k)\Big|\pi\leq\ker\left((i_1i_2)^{k_1},(i_2i_3)^{k_2},\ldots,(i_si_1)^{k_s}\right)\right\}\\
&=&t^{s/2}\sum_{\pi\in\mathcal P_2(k)}\frac{1}{N^{s/2+1}}\#\left\{i\in\{1,\ldots,N\}^s\Big|\pi\leq\ker\left((i_1i_2)^{k_1},(i_2i_3)^{k_2},\ldots,(i_si_1)^{k_s}\right)\right\}
\end{eqnarray*}

(4) Our claim now is that in the $N\to\infty$ limit the combinatorics of the above sum simplifies, with only the noncrossing partitions contributing to the sum, and with each of them contributing precisely with a 1 factor, so that we will have, as desired:
\begin{eqnarray*}
M_k\left(\frac{Z}{\sqrt{N}}\right)
&=&t^{s/2}\sum_{\pi\in\mathcal P_2(k)}\Big(\delta_{\pi\in NC_2(k)}+O(N^{-1})\Big)\\
&\simeq&t^{s/2}\sum_{\pi\in\mathcal P_2(k)}\delta_{\pi\in NC_2(k)}\\
&=&t^{s/2}|\mathcal{NC}_2(k)|
\end{eqnarray*}

(5) In order to prove this, the first observation is that when $k$ is not uniform, in the sense that it contains a different number of $\circ$, $\bullet$ symbols, we have $\mathcal P_2(k)=\emptyset$, and so:
$$M_k\left(\frac{Z}{\sqrt{N}}\right)=t^{s/2}|\mathcal{NC}_2(k)|=0$$

(6) Thus, we are left with the case where $k$ is uniform. Let us examine first the case where $k$ consists of an alternating sequence of $\circ$ and $\bullet$ symbols, as follows:
$$k=\underbrace{\circ\bullet\circ\bullet\ldots\ldots\circ\bullet}_{2p}$$

In this case it is convenient to relabel our multi-index $i=(i_1,\ldots,i_s)$, with $s=2p$, in the form $(j_1,l_1,j_2,l_2,\ldots,j_p,l_p)$. With this done, our moment formula becomes:
$$M_k\left(\frac{Z}{\sqrt{N}}\right)
=t^p\sum_{\pi\in\mathcal P_2(k)}\frac{1}{N^{p+1}}\#\left\{j,l\in\{1,\ldots,N\}^p\Big|\pi\leq\ker\left(j_1l_1,j_2l_1,j_2l_2,\ldots,j_1l_p\right)\right\}$$

Now observe that, with $k$ being as above, we have an identification $\mathcal P_2(k)\simeq S_p$, obtained in the obvious way. With this done too, our moment formula becomes:
$$M_k\left(\frac{Z}{\sqrt{N}}\right)
=t^p\sum_{\pi\in S_p}\frac{1}{N^{p+1}}\#\left\{j,l\in\{1,\ldots,N\}^p\Big|j_r=j_{\pi(r)+1},l_r=l_{\pi(r)},\forall r\right\}$$

(7) We are now ready to do our asymptotic study, and prove the claim in (4). Let indeed $\gamma\in S_p$ be the full cycle, which is by definition the following permutation:
$$\gamma=(1 \, 2 \, \ldots \, p)$$

In terms of $\gamma$, the conditions $j_r=j_{\pi(r)+1}$ and $l_r=l_{\pi(r)}$ found above read:
$$\gamma\pi\leq\ker j\quad,\quad 
\pi\leq\ker l$$

Counting the number of free parameters in our moment formula, we obtain:
$$M_k\left(\frac{Z}{\sqrt{N}}\right)
=\frac{t^p}{N^{p+1}}\sum_{\pi\in S_p}N^{|\pi|+|\gamma\pi|}
=t^p\sum_{\pi\in S_p}N^{|\pi|+|\gamma\pi|-p-1}$$

(8) The point now is that the last exponent is well-known to be $\leq 0$, with equality precisely when the permutation $\pi\in S_p$ is geodesic, which in practice means that $\pi$ must come from a noncrossing partition. Thus we obtain, in the $N\to\infty$ limit, as desired:
$$M_k\left(\frac{Z}{\sqrt{N}}\right)\simeq t^p|\mathcal{NC}_2(k)|$$

This finishes the proof in the case of the exponents $k$ which are alternating, and the case where $k$ is an arbitrary uniform exponent is similar, by permuting everything.
\end{proof}

As a comment now on our findings from Theorem 8.23, we have:

\begin{comment}
Given a sequence of complex Gaussian matrices $Z_N\in M_N(L^\infty(X))$, having independent $G_t$ variables as entries, with $t>0$, we have
$$\frac{Z_N}{\sqrt{N}}\sim\Gamma_t$$
in the $N\to\infty$ limit, with the limiting law being defined via its colored moments by
$$M_k(\Gamma_t)=\sum_{\pi\in\mathcal{NC}_2(k)}t^{|\pi|}$$
and being a law which is not a usual measure, called Voiculescu circular law.
\end{comment}

To be more precise here, $\Gamma_t$ is in fact a beast that we already met in chapter 7, when talking easiness, appearing as the asymptotic character law for the quantum group $U_N^+$. We will see later, in chapter 13, that free probability provides a conceptual approach to all this, both Gaussian matrices and quantum groups, and with the central law there, appearing as free analogue of the normal law $G_t$, or as complex analogue of the Wigner semicircle law $\gamma_t$, being the Voiculescu circular law $\Gamma_t$. More on this later.

\bigskip
 
Regarding now the Wigner matrices, which are something more understandable, because they are self-adjoint, we have the following result, about them:

\begin{theorem}
Given a sequence of Wigner random matrices
$$Z_N\in M_N(L^\infty(X))$$
having independent $G_t$ variables as entries, with $t>0$, up to $Z_N=Z_N^*$, we have
$$M_k\left(\frac{Z_N}{\sqrt{N}}\right)\simeq t^{k/2}|NC_2(k)|$$
for any integer $k\in\mathbb N$, in the $N\to\infty$ limit. Thus we have the convergence
$$\frac{Z_N}{\sqrt{N}}\sim\gamma_t$$
in the $N\to\infty$ limit, where $\gamma_t$ is the Wigner semicircle law. 

\end{theorem}

\begin{proof}
We have two possible proofs here, as follows:

\medskip

(1) The formula in the statement can be certainly established via a direct computation based on the Wick formula, similar to that from the proof of Theorem 8.23. We will leave the calculations here, featuring no particular difficulties, as an exercise.

\medskip

(2) However, the best is to deduce this result from Theorem 8.23 itself. Indeed, we know from there that for Gaussian matrices $Y_N\in M_N(L^\infty(X))$ we have the following formula, valid for any colored integer $K=\circ\bullet\bullet\circ\ldots\,$, in the $N\to\infty$ limit:
$$M_K\left(\frac{Y_N}{\sqrt{N}}\right)\simeq t^{|K|/2}|\mathcal{NC}_2(K)|$$

By doing now some combinatorics, we deduce that we have the following formula for the moments of the matrices $Re(Y_N)$, with respect to usual exponents, $k\in\mathbb N$:
\begin{eqnarray*}
M_k\left(\frac{Re(Y_N)}{\sqrt{N}}\right)
&=&2^{-k}\cdot M_k\left(\frac{Y_N}{\sqrt{N}}+\frac{Y_N^*}{\sqrt{N}}\right)\\
&=&2^{-k}\sum_{|K|=k}M_K\left(\frac{Y_N}{\sqrt{N}}\right)\\
&\simeq&2^{-k}\sum_{|K|=k}t^{k/2}|\mathcal{NC}_2(K)|\\
&=&2^{-k}\cdot t^{k/2}\cdot 2^{k/2}|\mathcal{NC}_2(k)|\\
&=&2^{-k/2}\cdot t^{k/2}|NC_2(k)|
\end{eqnarray*}

In terms of the matrices $Z_N=\sqrt{2}Re(Y_N)$, which are of Wigner type, this reads:
$$M_k\left(\frac{Z_N}{\sqrt{2N}}\right)\simeq 2^{-k/2}\cdot t^{k/2}|NC_2(k)|$$

Thus, we are led to the moment formula in the statement, which itself leads to $\gamma_t$.
\end{proof}

So long for the Wigner theorem. Of course, this was just a beginning, and far more things can be said, regarding the precise nature of the convergence, and the fluctuations, with the problematics here being highly non-trivial, and still subject to research. For more on all this, you can consult the various random matrix books and articles referenced at the end of this book. And by the way, do not forget to learn, either from Wigner himself \cite{wig}, or from a book like Mehta \cite{meh}, how Wigner came upon these matrices, in relation with the physics of the heavy atoms. First class science, all this.

\bigskip

In relation now with free probability, that we will learn later in this book, there is no comment in the spirit of Comment 8.24 to be made, because Theorem 8.25 as stated is just fine. However, something interesting happens when looking at families of Wigner matrices, with the finding here, due to Voiculescu \cite{vo3}, being as follows:

\begin{fact}
Given a family of sequences of Wigner matrices, 
$$Z^i_N\in M_N(L^\infty(X))\quad,\quad i\in I$$
with pairwise independent entries, each following the complex normal law $G_t$, with $t>0$, up to the constraint $Z_N^i=(Z_N^i)^*$, the rescaled sequences of matrices
$$\frac{Z^i_N}{\sqrt{N}}\in M_N(L^\infty(X))\quad,\quad i\in I$$
become with $N\to\infty$ semicircular, each following the Wigner law $\gamma_t$, and free.
\end{fact}

To be more precise here, what is new with respect to Theorem 8.25, taking into account the fact that we are now dealing with a family of Wigner matrices, is the word ``free'' at the end. And, more on this later in this book, when doing free probability.

\section*{8d. Wishart matrices}

Let us discuss now the Wishart matrices, which are the positive analogues of the Wigner matrices. Quite surprisingly, the computation here leads to the Catalan numbers, but not in the same way as for the Wigner matrices, the result being as follows:

\begin{theorem}
Given a sequence of complex Wishart matrices
$$W_N=Y_NY_N^*\in M_N(L^\infty(X))$$
with $Y_N$ being $N\times N$ complex Gaussian of parameter $t>0$, we have
$$M_k\left(\frac{W_N}{N}\right)\simeq t^kC_k$$
for any exponent $k\in\mathbb N$, in the $N\to\infty$ limit. Thus we have the convergence
$$\frac{W_N}{tN}\sim\frac{1}{2\pi}\sqrt{4x^{-1}-1}\,dx$$
with $N\to\infty$, with the limiting measure being the Marchenko-Pastur law $\pi_1$.
\end{theorem}

\begin{proof}
There are several possible proofs for this result, as follows:

\medskip

(1) A first method is by using the formula that we have in Theorem 8.23, for the Gaussian matrices $Y_N$. Indeed, we know from there that we have the following formula, valid for any colored integer $K=\circ\bullet\bullet\circ\ldots\,$, in the $N\to\infty$ limit:
$$M_K\left(\frac{Y_N}{\sqrt{N}}\right)\simeq t^{|K|/2}|\mathcal{NC}_2(K)|$$

With $K=\circ\bullet\circ\bullet\ldots\,$, alternating word of length $2k$, with $k\in\mathbb N$, this gives:
$$M_k\left(\frac{Y_NY_N^*}{N}\right)\simeq t^k|\mathcal{NC}_2(K)|$$

Thus, in terms of the Wishart matrix $W_N=Y_NY_N^*$ we have, for any $k\in\mathbb N$:
$$M_k\left(\frac{W_N}{N}\right)\simeq t^k|\mathcal{NC}_2(K)|$$

The point now is that, by doing some combinatorics, we have:
$$|\mathcal{NC}_2(K)|=|NC_2(2k)|=C_k$$

Thus, we are led to the formula in the statement.

\medskip

(2) A second method, that we will explain now as well, is by proving the result directly, starting from definitions. The matrix entries of our matrix $W=W_N$ are given by:
$$W_{ij}=\sum_{r=1}^NY_{ir}\bar{Y}_{jr}$$

Thus, the normalized traces of powers of $W$ are given by the following formula:
\begin{eqnarray*}
tr(W^k)
&=&\frac{1}{N}\sum_{i_1=1}^N\ldots\sum_{i_k=1}^NW_{i_1i_2}W_{i_2i_3}\ldots W_{i_ki_1}\\
&=&\frac{1}{N}\sum_{i_1=1}^N\ldots\sum_{i_k=1}^N\sum_{r_1=1}^N\ldots\sum_{r_k=1}^NY_{i_1r_1}\bar{Y}_{i_2r_1}Y_{i_2r_2}\bar{Y}_{i_3r_2}\ldots Y_{i_kr_k}\bar{Y}_{i_1r_k}
\end{eqnarray*}

By rescaling now $W$ by a $1/N$ factor, as in the statement, we obtain:
$$tr\left(\left(\frac{W}{N}\right)^k\right)=\frac{1}{N^{k+1}}\sum_{i_1=1}^N\ldots\sum_{i_k=1}^N\sum_{r_1=1}^N\ldots\sum_{r_k=1}^NY_{i_1r_1}\bar{Y}_{i_2r_1}Y_{i_2r_2}\bar{Y}_{i_3r_2}\ldots Y_{i_kr_k}\bar{Y}_{i_1r_k}$$

By using now the Wick rule, we obtain the following formula for the moments, with $K=\circ\bullet\circ\bullet\ldots\,$, alternating word of lenght $2k$, and with $I=(i_1r_1,i_2r_1,\ldots,i_kr_k,i_1r_k)$:
\begin{eqnarray*}
M_k\left(\frac{W}{N}\right)
&=&\frac{t^k}{N^{k+1}}\sum_{i_1=1}^N\ldots\sum_{i_k=1}^N\sum_{r_1=1}^N\ldots\sum_{r_k=1}^N\#\left\{\pi\in\mathcal P_2(K)\Big|\pi\leq\ker(I)\right\}\\
&=&\frac{t^k}{N^{k+1}}\sum_{\pi\in\mathcal P_2(K)}\#\left\{i,r\in\{1,\ldots,N\}^k\Big|\pi\leq\ker(I)\right\}
\end{eqnarray*}

In order to compute this quantity, we use the standard bijection $\mathcal P_2(K)\simeq S_k$. By identifying the pairings $\pi\in\mathcal P_2(K)$ with their counterparts $\pi\in S_k$, we obtain:
\begin{eqnarray*}
M_k\left(\frac{W}{N}\right)
&=&\frac{t^k}{N^{k+1}}\sum_{\pi\in S_k}\#\left\{i,r\in\{1,\ldots,N\}^k\Big|i_s=i_{\pi(s)+1},r_s=r_{\pi(s)},\forall s\right\}
\end{eqnarray*}

Now let $\gamma\in S_k$ be the full cycle, which is by definition the following permutation:
$$\gamma=(1 \, 2 \, \ldots \, k)$$

The general factor in the product computed above is then 1 precisely when following two conditions are simultaneously satisfied:
$$\gamma\pi\leq\ker i\quad,\quad 
\pi\leq\ker r$$

Counting the number of free parameters in our moment formula, we obtain:
$$M_k\left(\frac{W}{N}\right)
=t^k\sum_{\pi\in S_k}N^{|\pi|+|\gamma\pi|-k-1}$$

The point now is that the last exponent is well-known to be $\leq 0$, with equality precisely when the permutation $\pi\in S_k$ is geodesic, which in practice means that $\pi$ must come from a noncrossing partition. Thus we obtain, in the $N\to\infty$ limit:
$$M_k\left(\frac{W}{N}\right)\simeq t^kC_k$$

Thus, we are led to the conclusion in the statement.
\end{proof}

As a comment now, while the above result is definitely something interesting at $t=1$, at general $t>0$ this looks more like a ``fake'' generalization of the $t=1$ result, because the law $\pi_1$ stays the same, modulo a trivial rescaling. Which is certainly not good.

\bigskip

The reasons behind this phenomenon are quite subtle, and skipping some discussion here, the point is that Theorem 8.27 is indeed something ``fake'' at general values $t>0$, and the correct generalization of the $t=1$ computation, involving more general classes of complex Wishart matrices, is something a bit more complicated, as follows:

\begin{theorem}
Given a sequence of general complex Wishart matrices
$$W_N=Y_NY_N^*\in M_N(L^\infty(X))$$
with $Y_N$ being $N\times M$ complex Gaussian of parameter $1$, we have
$$\frac{W_N}{N}\sim\max(1-t,0)\delta_0+\frac{\sqrt{4t-(x-1-t)^2}}{2\pi x}\,dx$$
with $M=tN\to\infty$, with the limiting measure being the Marchenko-Pastur law $\pi_t$.
\end{theorem}

\begin{proof}
This is again something which is very standard, as follows:

\medskip

(1) In order to prove the formula in the statement, we can proceed as usual, by using the Wick formula. The matrix entries of our Wishart matrix $W=W_N$ are given by:
$$W_{ij}=\sum_{r=1}^MY_{ir}\bar{Y}_{jr}$$

Thus, the normalized traces of powers of $W$ are given by the following formula:
\begin{eqnarray*}
tr(W^k)
&=&\frac{1}{N}\sum_{i_1=1}^N\ldots\sum_{i_k=1}^NW_{i_1i_2}W_{i_2i_3}\ldots W_{i_ki_1}\\
&=&\frac{1}{N}\sum_{i_1=1}^N\ldots\sum_{i_k=1}^N\sum_{r_1=1}^M\ldots\sum_{r_k=1}^MY_{i_1r_1}\bar{Y}_{i_2r_1}Y_{i_2r_2}\bar{Y}_{i_3r_2}\ldots Y_{i_kr_k}\bar{Y}_{i_1r_k}
\end{eqnarray*}

By rescaling now $W$ by a $1/N$ factor, as in the statement, we obtain:
$$tr\left(\left(\frac{W}{N}\right)^k\right)=\frac{1}{N^{k+1}}\sum_{i_1=1}^N\ldots\sum_{i_k=1}^N\sum_{r_1=1}^M\ldots\sum_{r_k=1}^MY_{i_1r_1}\bar{Y}_{i_2r_1}Y_{i_2r_2}\bar{Y}_{i_3r_2}\ldots Y_{i_kr_k}\bar{Y}_{i_1r_k}$$

(2) By using now the Wick rule, we obtain the following formula for the moments, with $K=\circ\bullet\circ\bullet\ldots\,$, alternating word of lenght $2k$, and $I=(i_1r_1,i_2r_1,\ldots,i_kr_k,i_1r_k)$:
\begin{eqnarray*}
M_k\left(\frac{W}{N}\right)
&=&\frac{1}{N^{k+1}}\sum_{i_1=1}^N\ldots\sum_{i_k=1}^N\sum_{r_1=1}^M\ldots\sum_{r_k=1}^M\#\left\{\pi\in\mathcal P_2(K)\Big|\pi\leq\ker I\right\}\\
&=&\frac{1}{N^{k+1}}\sum_{\pi\in\mathcal P_2(K)}\#\left\{i\in\{1,\ldots,N\}^k,r\in\{1,\ldots,M\}^k\Big|\pi\leq\ker I\right\}
\end{eqnarray*}

(3) In order to compute this quantity, we use the standard bijection $\mathcal P_2(K)\simeq S_k$. By identifying the pairings $\pi\in\mathcal P_2(K)$ with their counterparts $\pi\in S_k$, we obtain:
\begin{eqnarray*}
M_k\left(\frac{W}{N}\right)
&=&\frac{1}{N^{k+1}}\sum_{\pi\in S_k}\#\left\{i\in\{1,\ldots,N\}^k,r\in\{1,\ldots,M\}^k\Big|i_s=i_{\pi(s)+1},r_s=r_{\pi(s)}\right\}
\end{eqnarray*}

Now let $\gamma\in S_k$ be the full cycle, which is by definition the following permutation:
$$\gamma=(1 \, 2 \, \ldots \, k)$$

The general factor in the product computed above is then 1 precisely when following two conditions are simultaneously satisfied:
$$\gamma\pi\leq\ker i\quad,\quad 
\pi\leq\ker r$$

Counting the number of free parameters in our expectation formula, we obtain:
$$M_k\left(\frac{W}{N}\right)
=\frac{1}{N^{k+1}}\sum_{\pi\in S_k}N^{|\gamma\pi|}M^{|\pi|}
=\sum_{\pi\in S_k}N^{|\gamma\pi|-k-1}M^{|\pi|}$$

(4) Now by using the same arguments as in the case $M=N$, from the proof of Theorem 8.27, we conclude that in the $M=tN\to\infty$ limit the permutations $\pi\in S_k$ which matter are those coming from noncrossing partitions, and so that we have:
$$M_k\left(\frac{W}{N}\right)
\simeq\sum_{\pi\in NC(k)}N^{-|\pi|}M^{|\pi|}
=\sum_{\pi\in NC(k)}t^{|\pi|}$$

(5) But these numbers are the moments of the Marchenko-Pastur law $\pi_t$, which in addition has the density given by the formula in the statement, as desired.
\end{proof}

\section*{8e. Exercises}

This was a standard introduction to random matrices, and as exercises, we have:

\begin{exercise}
Learn the spectral theorem for the normal matrices.
\end{exercise}

\begin{exercise}
Further experiment with the moments of the Jordan block $J$.
\end{exercise}

\begin{exercise}
Learn about the spectral theorem for the normal operators.
\end{exercise}

\begin{exercise}
Learn more about $C^*$-algebras, notably with the GNS theorem.
\end{exercise}

\begin{exercise}
Experiment a bit with that mysterious circular law $\Gamma_t$.
\end{exercise}

\begin{exercise}
What happens in the Wigner theorem, when using $\gamma_t$ variables?
\end{exercise}

\begin{exercise}
Learn more about positive matrices, in particular about Wishart.
\end{exercise}

\begin{exercise}
Learn also about the real Wishart matrices, and their properties.
\end{exercise}

As bonus exercise, study the joint distributions of random matrices, with $N\to\infty$.

\part{Higher dimensions}

\ \vskip50mm

\begin{center}
{\em Son io che corro nella tua mente

Non spaventarti, voglio solo un po' giocare

Son tanto sola e sto sognando

Io sto inseguendo la farfalla dell'amor}
\end{center}

\chapter{Gaussian vectors}

\section*{9a. Gaussian vectors}

Welcome to higher dimensions, and our plan for the present Part III will be, a bit like the plan for Parts I and II was, and like the plan for Part IV will be too, first half basic theory, and second half advanced questions. In short, and in the hope that you got it, this book tells two stories, with chapters 1-2, 5-6, 9-10, 13-14 discussing the basics, and with chapters 3-4, 7-8, 11-12, 15-16 discussing more advanced questions.

\bigskip

The question that we will be interested in, in the present chapter, is as follows:

\begin{question}
What are the reasonably general analogues of the normal variables, in $n$ dimensions? What about the basic applications of these?
\end{question} 

And good question this is, with the first part meant to bring us into some interesting mathematics, generalizing what we already know at $n=1,2$, and with the second part destined to bring us into some interesting physics, namely heat diffusion.

\bigskip

Let us start with the case $n=2$. Here we already have some good knowledge from our study of complex probability, from chapters 5-6, but since we would like now to have everything real, this is worth a rediscussion, with real notations. Let us start with:

\begin{definition}
The axioms for 2D probability theory are as follows:
\begin{enumerate}
\item A probability space is a measured space $(X,M,\nu)$ of mass one, $\nu(X)=1$.

\item A random variable on $X$ is a measurable function $f:X\to\mathbb R^2$.

\item The expectation of such a variable is $E(f)=\int_Xf(x)d\nu(x)\in\mathbb R^2$.

\item The variance of such a variable is $V(f)=E(||f-E(f)||^2)\geq0$.

\item The law of $f$ is the 2D probability measure $\mu=f_*\nu$, push-forward of $\nu$ by $f$.
\end{enumerate}
\end{definition}

Obviously, all this is quite straightforward, inspired by the theory in chapter 5. In relation with the variance, observe that we have 3 equivalent formulae for it:
\begin{eqnarray*}
V(f)
&=&E(||f-E(f)||^2)\\
&=&V(f_1)+V(f_2)\\
&=&E(||f||^2)-||E(f)||^2
\end{eqnarray*}

Indeed, the second formula is clear, and the third one follows from it.

\bigskip

In relation now with the law, observe that we have the following formula, valid for any  $k,l\in\mathbb N$, allowing us in particular to compute the expectation and variance:
$$E(f_1^kf_2^l)=\int_{\mathbb R^2}x_1^kx_2^ld\mu(x)$$

In fact, the above quantities $E(f_1^kf_2^l)$, which uniquely determine the law, can be thought of as being the moments of $f$. As an illustration for this, we have:

\index{complex variable}
\index{random variable}

\begin{proposition}
For a discrete variable $f:X\to\mathbb R^2$, the law is the measure
$$\mu=\sum_i c_i\delta_{a_i}\quad,\quad c_i>0\quad,\quad\sum_ic_i=1\quad,\quad a_i\in\mathbb R^2$$
given by the following formula, with $P$ being the probability over $X$,
$$\mu=\sum_{a\in\mathbb R^2}P(f=a)\delta_a$$
with the sum being finite or countable, as per our discretness assumption. We have
$$E(f)=\sum_ic_ia_i\quad,\quad V(f)=\sum_ic_i||a_i||^2-\Big|\Big|\sum_ic_ia_i\Big|\Big|^2$$
and with the notation $a_i=\binom{p_i}{q_i}$, the moments are $E(f_1^kf_2^l)=\sum_ic_ip_i^kq_i^l$.
\end{proposition}

\begin{proof}
The first assertion is self-explanatory. Regarding the moments, we have:
\begin{eqnarray*}
E(f_1^kf_2^l)
&=&\sum_{x\in\mathbb R}P(f_1^kf_2^l=x)x\\
&=&\sum_{p,q\in\mathbb R}P(f_1=p,f_2=q)p^kq^l\\
&=&\sum_{a\in\mathbb R^2}P(f=a)a_1^ka_2^l\\
&=&\int_{\mathbb R^2}x_1^kx_2^ld\mu(x)
\end{eqnarray*}

In terms of $\mu=\sum_i c_i\delta_{a_i}$ with $a_i=\binom{p_i}{q_i}$, this gives the formula in the statement:
$$E(f_1^kf_2^l)=\sum_ic_ip_i^kq_i^l$$

Finally, the formulae of $E(f)$ and $V(f)$ are both clear from definitions.
\end{proof}

In order to get now beyond what we previously knew, from chapter 5, let us have a closer look at the low order moments. At order 2 there are the numbers $E(f_1^2)$, $E(f_2^2)$, that we are certainly familiar with, these being the usual second moments of the variables $f_1,f_2$, and then the mixed moment $E(f_1f_2)$, which is in need of some study.

\bigskip

So, let us temporary pause our 2D study, and formulate, as a basic 1D fact:

\index{covariance}

\begin{theorem}
Given variables $f,g:X\to\mathbb R$, we can talk about their covariance
\begin{eqnarray*}
cov(f,g)
&=&E[(f-E(f))(g-E(g))]\\
&=&E(fg)-E(f)E(g)
\end{eqnarray*}
with this quantity having the following properties:
\begin{enumerate}
\item $cov(f,f)=V(f)$.

\item If $f,g$ are independent, $cov(f,g)=0$.

\item The converse of this latter fact is not true.

\item $cov(f,g)^2\leq V(f)V(g)$.
\end{enumerate}
\end{theorem}

\begin{proof}
This is something very standard, the idea being as follows:

\medskip

(1) This is indeed clear, with either of the definitions for the variance $V(f)$.

\medskip

(2) This is clear too, coming from the multiplicativity formula for mixed moments.

\medskip

(3) To start with, the converse of (2) has obviously zero chances to hold, because we perfectly know from chapter 1 that the independence of $f,g:X\to\mathbb R$ amounts in having the following formula, for the mixed moments, meaning for any $k,l\in\mathbb N$:
$$E(f^kg^l)=E(f^k)E(g^l)$$

This being said, an explicit counterexample will certainly not hurt. Consider two Bernoulli variables $f,g:X\to\mathbb R$, taking their values according to the following table:
\vskip-7mm
$$\xymatrix@R=0pt@C=4pt{
&&\ar@{-}[ddddddd]&&&\ar@{-}[ddddddd]\\
&g\backslash f&\ar@{-}[dd]&a&b\\
\ar@{-}[rrrrrrr]&&&&&&&\\
&c&&pq+x&(1-p)q-x&&q\\
&d&&p(1-q)-x&(1-p)(1-q)+x&&1-q\\
\ar@{-}[rrrrrrr]&&&&&&&\\
&&&p&1-p&\\
&&&&&}$$
\vskip-2mm

The covariance of $f,g$ can be then computed as follows:
\begin{eqnarray*}
cov(f,g)
&=&E(fg)-E(f)E(g)\\
&=&[ac(pq+x)+ad(p(1-q)-x)+bc((1-p)q-x)\\
&&+bd((1-p)(1-q)+x)]-[ap+b(1-p)][cq+d(1-q)]\\
&=&acx-adx-bcx+bdx\\
&=&(a-b)(c-d)x
\end{eqnarray*}

Which is bad news, because $cov(f,g)=0$ means $x=0$, which in turn means that $f,g$ must be independent. Damn, so instead of our counterexample, we have a theorem here, stating that the converse of (2) holds for the Bernoulli variables. Good to know.

\medskip

(3') Counterexample, take two. In view of the above, consider the simplest possible situation which is not (Bernoulli, Bernoulli), which is as follows, with $f$ following a uniform trivalent distribution on $[-1,1]$, and with $g$ being taken to be its square, $g=f^2$: 
\vskip-7mm
$$\xymatrix@R=0pt@C=6pt{
&&\ar@{-}[ddddddd]&&&&\ar@{-}[ddddddd]\\
&g\backslash f&\ar@{-}[dd]&-1&0&1\\
\ar@{-}[rrrrrrrr]&&&&&&&&\\
&0&&0&1/3&0&&1/3\\
&1&&1/3&0&1/3&&2/3\\
\ar@{-}[rrrrrrrr]&&&&&&&&\\
&&&1/3&1/3&1/3\\
&&&&&&}$$
\vskip-2mm

Our variables are then not independent, because the upper left 0 does not come as $0=1/3\cdot1/3$. On the other hand we have $E(f)=0$, $E(g)=1/3$, $E(fg)=0$, so the covariance is $cov(f,g)=0-0=0$, and we have our counterexample, as desired.

\medskip

(4) In order to estimate $cov(f,g)$ we can use the Cauchy-Schwarz inequality, applied for the scalar product $<f,g>=E(fg)$. Indeed, with $a=E(f)$, $b=E(g)$, we have:
\begin{eqnarray*}
cov(f,g)^2
&=&E\big((f-a)(g-b)\big)^2\\
&=&<f-a,g-b>^2\\
&\leq&||f-a||^2||g-b||^2\\
&=&E((f-a)^2)E((g-b)^2)\\
&=&V(f)V(g)
\end{eqnarray*}

Thus, we have indeed the covariance inequality in the statement.
\end{proof}

As a complement to Theorem 9.4, in analogy with the fact that the variance $V(f)$ is usefully complemented by the standard deviation $\sigma_f=\sqrt{V(f)}$, let us formulate:

\index{correlation}

\begin{theorem}[continuation]
Given $f,g:X\to\mathbb R$, we can equally talk about
$$\rho_{fg}=\frac{cov(f,g)}{\sigma_f\sigma_g}$$
called correlation coefficient, which satisfies $\rho_{fg}\in[-1,1]$.
\end{theorem}

\begin{proof}
This is something self-explanatory, with $\rho_{fg}\in[-1,1]$ coming by extracting the square root of the inequality $cov(f,g)^2\leq V(f)V(g)$ from Theorem 9.4.
\end{proof}

Getting back now to our 2D questions, we recall that a given variable $f:X\to\mathbb R^2$ has three order 2 moments, namely $E(f_1^2)$, $E(f_2^2)$, which are familiar objects, namely usual second moments of the variables $f_1,f_2$, and then the mixed moment $E(f_1f_2)$. Now the above results suggest replacing this second moment $E(f_1f_2)$ by the more familiar quantity $cov(f_1,f_2)$, or equivalently, by $\rho_{f_1f_2}$. And the point is that all this leads us to:

\index{covariance matrix}

\begin{definition}
Associated to $f:X\to\mathbb R^2$ is its covariance matrix, given by:
$$\Sigma=\Big[cov(f_i,f_j)\Big]_{ij}$$
Equivalently, and more explicitly, this covariance matrix is given by:
$$\Sigma=\begin{pmatrix}V(f_1)&cov(f_1,f_2)\\ cov(f_1,f_2)&V(f_2)\end{pmatrix}$$
Alternatively, in terms of standard deviations and correlations, we have
$$\Sigma=\begin{pmatrix}\sigma_1^2&\rho\sigma_1\sigma_2\\ \rho\sigma_1\sigma_2&\sigma_2^2\end{pmatrix}$$
with $\sigma_i=\sigma_{f_i}$ and $\rho=\rho_{f_1f_2}$. This matrix is positive, and $Tr(\Sigma)=V(f)$.
\end{definition}

To be more precise, the equivalence between the above 3 definitions of $\Sigma$ comes from the various formulae in Theorems 9.4 and 9.5, and the positivity comes from:
$$Tr(\Sigma)=V(f)\geq 0\quad,\quad \det(\Sigma)=(1-\rho^2)\sigma_1^2\sigma_2^2\geq0$$

Observe that, despite the simple form of $\Sigma$ in terms of standard deviations, the eigenvalues remain elusive, the discriminant of the characteristic polynomial being:
\begin{eqnarray*}
\Delta
&=&Tr(\Sigma)^2-4\det(\Sigma)\\
&=&(\sigma_1^2+\sigma_2^2)^2-4(1-\rho^2)\sigma_1^2\sigma_2^2\\
&=&(\sigma_1^2-\sigma_2^2)^2+8\rho^2\sigma_1^2\sigma_2^2
\end{eqnarray*}

With this discussed, time to get back now to Question 9.1, at $n=2$? Let us start with a real remake of some things that we knew from chapter 6, namely:

\begin{theorem}
Assuming that $f_1\sim g_s^b$ and $f_2\sim g_t^c$ are independent, the variable $f=\binom{f_1}{f_2}$ follows the following law:
$$g_{st}^{bc}=\frac{1}{2\pi\sqrt{st}}\,\exp\left(-\frac{(x-b)^2}{2s}-\frac{(y-c)^2}{2t}\right)dxdy$$
In this situation the covariance matrix is diagonal, $\Sigma=diag(s,t)$.
\end{theorem}

\begin{proof}
According to our formulae for shifted bells, from chapter 2, we have:
$$g_s^b=\frac{1}{\sqrt{2\pi s}}\,e^{-(x-b)^2/2s}\quad,\quad 
g_t^c=\frac{1}{\sqrt{2\pi t}}\,e^{-(x-c)^2/2t}$$

Now by multiplying, we are led to the density formula in the statement. As for $\Sigma=diag(s,t)$, this is clear, since we have $f_1\sim g_s^b$ and $f_2\sim g_t^c$, independent.
\end{proof}

The problem is now, can we get beyond this? And in answer, yes we can, by formulating the following nice answer to the first part of Question 9.1, at $n=2$:

\index{Gaussian vector}

\begin{definition}
We call 2D Gaussian vector any variable 
$$f:X\to\mathbb R^2$$
having the property that $\alpha f_1+\beta f_2$ is normal, for any $\alpha,\beta\in\mathbb R$.
\end{definition}

Obviously, this is something quite conceptual, generalizing the situation in Theorem 9.7. The motivations come from the 2D versions of the CLT, in the spirit of those from chapter 5, leading to such variables. We will leave some learning here as an exercise.

\bigskip

As our most pressing task, we must check that Definition 9.8 brings indeed something new, with respect to what we had in Theorem 9.7. And here, we have:

\begin{theorem}
Given variables $h_1,\ldots,h_m:X\to\mathbb R$ which are normal, $h_i\sim g_{t_i}^{a_i}$, and independent, and scalars $\lambda_1,\ldots,\lambda_m\in\mathbb R$ and $\mu_1,\ldots,\mu_m\in\mathbb R$, the variable
$$f=\binom{\lambda_1h_1+\ldots+\lambda_mh_m}{\mu_1h_1+\ldots+\mu_mh_m}$$
is a Gaussian vector in the above sense, and the corresponding covariance matrix is:
$$\Sigma=\begin{pmatrix}
\sum_i\lambda_i^2t_i&\sum_i\lambda_i\mu_it_i\\
\sum_i\lambda_i\mu_it_i&\sum_i\mu_i^2t_i\end{pmatrix}$$
At $m=2$ and $\lambda=(1,0)$, $\mu=(0,1)$, we recover the situation from Theorem 9.7.
\end{theorem}

\begin{proof}
The first assertion is clear, since a linear combination of independent normal variables is normal, as you can see for instance by looking at the corresponding Fourier transform. Regarding now the computation of $\Sigma$, observe that we have: 
\begin{eqnarray*}
cov(f_1,f_2)
&=&cov\left(\sum_i\lambda_ih_i\,,\ \sum_j\mu_jh_j\right)\\
&=&E\left(\sum_{ij}\lambda_i\mu_jh_ih_j\right)-E\left(\sum_i\lambda_ih_i\right)E\left(\sum_j\mu_jh_j\right)\\
&=&\sum_i\lambda_i\mu_i(a_i^2+t_i)+\sum_{i\neq j}\lambda_i\mu_ja_ia_j-\left(\sum_i\lambda_ia_i\right)\left(\sum_j\mu_ja_j\right)\\
&=&\sum_i\lambda_i\mu_it_i+\sum_{ij}\lambda_i\mu_ja_ia_j-\left(\sum_i\lambda_ia_i\right)\left(\sum_j\mu_ja_j\right)\\
&=&\sum_i\lambda_i\mu_it_i
\end{eqnarray*}

Thus we have the off-diagonal terms of $\Sigma$, and in fact the diagonal terms too, by setting $\lambda_i=\mu_i$ in the above computation. Finally, the last assertion is clear.
\end{proof}

The above result is quite interesting, and suggests the use of some linear algebra, in order to better understand what is going on. Consider the following matrices:
$$L=\begin{pmatrix}\lambda_1&\ldots&\lambda_m\\
\mu_1&\ldots&\mu_m\end{pmatrix}
\quad,\quad 
h=\begin{pmatrix}
h_1\\
\vdots\\
h_m\end{pmatrix}\quad,\quad 
T=\begin{pmatrix}
t_1\\
&\ddots\\
&&t_m
\end{pmatrix}$$

We have then $f=Lh$ and $\Sigma=LTL^t$, and Theorem 9.9 reformulates as follows:

\begin{theorem}
Given $h:X\to\mathbb R^m$ with normal independent entries, and $L\in M_{2\times m}(\mathbb R)$, the variable $Lh:X\to\mathbb R^2$ is Gaussian, with $\Sigma=LTL^t$, $T=diag(V(h_i))$.
\end{theorem}

\begin{proof}
This is linear algebra magic, where all that indices have gone? Mystery.
\end{proof}

With this discussed, let us get now into the computation of the density of variables $Lh$ as above, or more generally, of Gaussian vectors in the sense of Definition 9.8. We first have a technical result, in the general framework of Definition 9.8, as follows:

\begin{proposition}
Given a Gaussian vector $f:X\to\mathbb R^2$, its law depends only on
$$(s,t,\rho,b,c)$$
the parameters being $s=V(f_1)$, $t=V(f_2)$, $\rho=\rho_{f_1f_2}$ and the position $\binom{b}{c}\in\mathbb R^2$.
\end{proposition}

\begin{proof}
According to Definition 9.8, our vector $f:X\to\mathbb R^2$ being Gaussian means that $\beta f_1+\gamma f_2$ is normal, for any $\beta,\gamma\in\mathbb R$. But, we can exploit this as follows:

\medskip

(1) To start with, by using our normality assumption with $\beta=1,\gamma=0$, and then with $\beta=0,\gamma=1$, we must have $f_1\sim g_s^b$ and $f_2\sim g_t^c$, for some $s,t,b,c$.

\medskip

(2) Next, let us study the variables $\beta f_1+\gamma f_2$. The expectation is given by:
$$E(\beta f_1+\gamma f_2)=\beta b+\gamma c$$

As for the variance of these variables, this is given by the following formula:
\begin{eqnarray*}
V(\beta f_1+\gamma f_2)
&=&E\left[(\beta f_1+\gamma f_2)^2\right]-\left[E(\beta f_1+\gamma f_2)\right]^2\\
&=&\beta^2E(f_1^2)+\gamma^2E(f_2^2)+2\beta\gamma E(f_1f_2)-\left[\beta E(f_1)+\gamma E(f_2)\right]^2\\
&=&\beta^2(b^2+s)+\gamma^2(c^2+t)+2\beta\gamma E(f_1f_2)-\beta^2b^2-\gamma^2c^2-2\beta\gamma bc\\
&=&\beta^2s+\gamma^2t+2\beta\gamma(E(f_1f_2)-bc)\\
&=&\beta^2s+\gamma^2t+2\beta\gamma cov(f_1,f_2)\\
&=&\beta^2s+\gamma^2t+2\beta\gamma\sqrt{st}\rho
\end{eqnarray*}

We conclude that our abstract normality assumption on $\beta f_1+\gamma f_2$ reformulates into the following more precise condition, telling us what its normal law exactly is:
$$\beta f_1+\gamma f_2\sim g^{\beta b+\gamma c}_{\beta^2s+\gamma^2t+2\beta\gamma\sqrt{st}\rho}$$

(3) Let us look now at the moments of $\beta f_1+\gamma f_2$. These are given by:
$$M_k(\beta f_1+\gamma f_2)=\sum_{r=0}^k\binom{k}{r}\beta^r\gamma^{k-r}E(f_1^rf_2^{k-r})$$

On the other hand, according to the formula in (2), these moments are given by the following formula, with $\varphi$ being a certain 2-variable polynomial function:
$$M_k(\beta f_1+\gamma f_2)=\varphi(\beta b+\gamma c\,,\,\beta^2s+\gamma^2t+2\beta\gamma\sqrt{st}\rho)$$

We conclude from this, by recurrence, even without knowing the formula of $\varphi$, that we can compute all the moments $E(f_1^kf_2^l)$ in terms of $(s,t,\rho,b,c)$, as desired.
\end{proof}

Good news, we can compute now the density of Gaussian vectors, as follows:

\begin{theorem}
Given a Gaussian vector $f:X\to\mathbb R^2$, its density is
$$g_{st\rho}^{bc}=\frac{1}{2\pi\sqrt{\det\Sigma}}\,\exp\left(-\frac{<\Sigma^{-1}(z-a),z-a>}{2}\right)dz$$
the parameters being $s=V(f_1)$, $t=V(f_2)$, $\rho=\rho_{f_1f_2}$ and the position $a=\binom{b}{c}\in\mathbb R^2$.
\end{theorem}

\begin{proof}
This is something quite tricky, the idea being as follows:

\medskip

(1) Let us first see what happens in the context of Theorem 9.7, where $f=\binom{f_1}{f_2}$ with $f_1\sim g_s^b$ and $f_2\sim g_t^c$, independent. In this case, as mentioned in Theorem 9.7, the covariance matrix is $\Sigma=diag(s,t)$, and our result predicts the correct density:
\begin{eqnarray*}
g_{st0}^{bc}
&=&\frac{1}{2\pi\sqrt{st}}\,\exp\left[-\frac{1}{2}\left<\begin{pmatrix}s^{-1}&0\\0&t^{-1}\end{pmatrix}\binom{x-b}{y-c},\binom{x-b}{y-x}\right>\right]dxdy\\
&=&\frac{1}{2\pi\sqrt{st}}\,\exp\left(-\frac{(x-b)^2}{2s}-\frac{(y-c)^2}{2t}\right)dxdy
\end{eqnarray*}

(2) In order to deal now with the general case, by using Proposition 9.11, we can assume that we are in the framework of Theorem 9.10, and in addition with $m=2$, with this being the smallest $m\in\mathbb N$ allowing us to model the correlation parameter $\rho$.

\medskip

(3) So, assume this, $f=Lh$ with $h_1\sim g_s^{b'}$ and $h_2\sim g_t^{c'}$ independent, and $L\in M_2(\mathbb R)$. According to Theorem 9.7, processed as in (1), we have, with $T=diag(s,t)$:
$$h\sim\frac{1}{2\pi\sqrt{\det T}}\,\exp\left(-\frac{<T^{-1}(z-a'),z-a'>}{2}\right)dz$$

Now observe that for any function $\varphi:\mathbb R^2\to\mathbb R$, we have the following formula:
\begin{eqnarray*}
E(\varphi(f))
&=&E(\varphi(Lh))\\
&=&\frac{1}{2\pi\sqrt{\det T}}\int_{\mathbb R^2}\varphi(Lz)\exp\left(-\frac{<T^{-1}(z-a'),z-a'>}{2}\right)dz\\
&=&\frac{1}{2\pi\sqrt{\det T}}\int_{\mathbb R^2}\varphi(w)\exp\left(-\frac{<T^{-1}(L^{-1}w-a'),L^{-1}w-a'>}{2}\right)\frac{dw}{|\det L|}\\
&=&\frac{1}{2\pi\sqrt{\det(LTL^t)}}\int_{\mathbb R^2}\varphi(w)\exp\left(-\frac{<(LTL^t)^{-1}(w-a),w-a>}{2}\right)dw\\
&=&\frac{1}{2\pi\sqrt{\det\Sigma}}\int_{\mathbb R^2}\varphi(w)\exp\left(-\frac{<\Sigma^{-1}(w-a),w-a>}{2}\right)dw
\end{eqnarray*}

Thus, we are led to the density formula in the statement.
\end{proof}

The above result is quite interesting, and as a first consequence, we have:

\begin{theorem}
For $f=Lh$, with $h:X\to\mathbb R^m$ having normal independent entries and $L\in M_{2\times m}(\mathbb R)$, we have $cov(f_1,f_2)=0$ precisely when $f_1,f_2$ are independent.
\end{theorem}

\begin{proof}
This follows from Theorem 9.12, because $cov(f_1,f_2)=0$ means $\rho=0$, so the law found there is $g_{st0}^{bc}=g_{st}^{bc}$. However, this is overkill, and alternatively, we have:

\medskip

(1) As starting point, we can use the standard fact, that we will leave as an exercise, that given a vector $h:X\to\mathbb R^m$ with entries $h_i\sim g_1$, independent, and an orthogonal matrix $U\in O_m$, the vector $k=Uh:X\to\mathbb R^m$ has also entries $k_i\sim g_1$, independent.

\medskip

(2) Now let prove the result, assuming $h_i\sim g_1$, and with the matrix $L\in M_{2\times m}(\mathbb R)$ being normalized too, as for both its rows to have norm 1. In order to apply the trick in (1), we would like to find an orthogonal matrix $U\in O_m$ which looks as follows:
$$U=\binom{L}{*}$$

(3) But, according to the Gram-Schmidt procedure, this is possible precisely when the rows of $L$ are orthogonal, which in turn corresponds to the condition $cov(f_1,f_2)=0$.

\medskip

(4) Summarizing, we proved our result, under the normalization assumptions in (2). As for the general case, this follows easily from this, say exercise for you.
\end{proof}

\section*{9b. Higher dimensions}

Let us discuss now what happens in arbitrary $n$ dimensions. This will be basically a straightforward remake of what we did in the above at $n=2$. But, we will use this occasion for rearranging a bit the whole theory, upgrading it too, by insisting on the underlying linear algebra, and also, hopefully, adding a few new things as well.

\bigskip

Let us start our discussion with the following straightforward definition:

\begin{definition}
The axioms for $n$-dimensional probability are as follows:
\begin{enumerate}
\item A probability space is a measured space $(X,M,\nu)$ of mass one, $\nu(X)=1$.

\item A random variable on $X$ is a measurable function $f:X\to\mathbb R^n$.

\item The expectation of such a variable is $E(f)=\int_Xf(x)d\nu(x)\in\mathbb R^n$.

\item The variance of such a variable is $V(f)=E(||f-E(f)||^2)\geq0$.

\item The law of $f$ is the probability measure $\mu=f_*\nu$, push-forward of $\nu$ by $f$.
\end{enumerate}
\end{definition}

Obviously, all this is quite straightforward, generalizing what we previously had in 2D. In relation with the variance, we have in fact 3 equivalent formulae, as follows:
\begin{eqnarray*}
V(f)
&=&E(||f-E(f)||^2)\\
&=&V(f_1)+\ldots+V(f_n)\\
&=&E(||f||^2)-||E(f)||^2
\end{eqnarray*}

In relation now with the law, observe that we have the following formula, valid for any  $k_1,\ldots,k_n\in\mathbb N$, allowing us in particular to compute the expectation and variance:
$$E(f_1^{k_1}\ldots f_n^{k_n})=\int_{\mathbb R^n}x_1^{k_1}\ldots x_n^{k_n}d\mu(x)$$

In fact, the above quantities $E(f_1^{k_1}\ldots f_n^{k_n})$, which uniquely determine the law, can be thought of as being the moments of $f$. As an illustration for this, we have:

\index{complex variable}
\index{random variable}

\begin{proposition}
For a discrete variable $f:X\to\mathbb R^n$, the law is the measure
$$\mu=\sum_i c_i\delta_{a_i}\quad,\quad c_i>0\quad,\quad\sum_ic_i=1\quad,\quad a_i\in\mathbb R^n$$
given by the following formula, with $P$ being the probability over $X$,
$$\mu=\sum_{a\in\mathbb R^n}P(f=a)\delta_a$$
with the sum being finite or countable, as per our discretness assumption. We have
$$E(f)=\sum_ic_ia_i\quad,\quad V(f)=\sum_ic_i||a_i||^2-\Big|\Big|\sum_ic_ia_i\Big|\Big|^2$$
and the moments are given by the following formula,
$$E(f_1^{k_1}\ldots f_n^{k_n})=\sum_ic_i(a_i)_1^{k_1}\ldots(a_i)_n^{k_n}$$
valid for any exponents $k_1,\ldots,k_n\in\mathbb N$.
\end{proposition}

\begin{proof}
The first assertion is self-explanatory. Regarding the moments, we have:
\begin{eqnarray*}
E(f_1^{k_1}\ldots f_n^{k_n})
&=&\sum_{x\in\mathbb R}P(f_1^{k_1}\ldots f_n^{k_n}=x)x\\
&=&\sum_{a_1,\ldots,a_n\in\mathbb R}P(f_1=a_1,\ldots,f_n=a_n)a_1^{k_1}\ldots a_n^{k_n}\\
&=&\sum_{a\in\mathbb R^n}P(f=a)a_1^{k_1}\ldots a_n^{k_n}\\
&=&\int_{\mathbb R^n}x_1^{k_1}\ldots x_n^{k_n}d\mu(x)
\end{eqnarray*}

In terms of $\mu=\sum_i c_i\delta_{a_i}$, this gives the formula in the statement, namely:
$$E(f_1^{k_1}\ldots f_n^{k_n})=\sum_ic_i(a_i)_1^{k_1}\ldots(a_i)_n^{k_n}$$

Finally, the formulae of $E(f)$ and $V(f)$ are both clear from definitions.
\end{proof}

Still in analogy with what we did in 2D, let us introduce the following key notion:

\begin{definition}
Associated to $f:X\to\mathbb R^n$ is its covariance matrix, given by:
$$\Sigma=\Big[cov(f_i,f_j)\Big]_{ij}$$
Equivalently, in terms of standard deviations and correlations, we have
$$\Sigma=\Big[\rho_{ij}\sigma_i\sigma_j\Big]_{ij}$$
with $\sigma_i=\sigma_{f_i}$ and $\rho_{ij}=\rho_{f_if_j}$. This matrix is positive, and $Tr(\Sigma)=V(f)$.
\end{definition}

To be more precise, the equivalence between the above two definitions of $\Sigma$ comes from the formulae in Theorems 9.4 and 9.5. Observe also that $\Sigma$ factorizes as:
$$\Sigma=\begin{pmatrix}\sigma_1\\ &\ddots\\ &&\sigma_n\end{pmatrix}
\begin{pmatrix}\rho_{11}&\ldots&\rho_{1n}\\ \vdots&&\vdots\\ \rho_{n1}&\ldots&\rho_{nn}\end{pmatrix}
\begin{pmatrix}\sigma_1\\ &\ddots\\ &&\sigma_n\end{pmatrix}$$

As for the positivity property of $\Sigma$, this comes from the following computation, with the convention $<f,g>=E(fg)$, and with the notation $a_i=E(f_i)$:
\begin{eqnarray*}
\Sigma_{ij}
&=&cov(f_i,f_j)\\
&=&E\left((f_i-a_i)(f_j-a_j)\right)\\
&=&<f_i-a_i,f_j-a_j>
\end{eqnarray*}

Indeed, this shows that $\Sigma$ is the Gram matrix of the vectors $f_i-a_i$, and such Gram matrices being always positive, we have positivity. Finally, $Tr(\Sigma)=V(f)$ is clear.

\bigskip

Moving now towards the generalization of our main results at $n=2$, and with the idea in mind of rearranging a bit that material, retrospectively, let us start with the following answer to the first part of Question 9.1, in the general, $n$-dimensional case:

\begin{definition}
We call Gaussian vector any variable 
$$f:X\to\mathbb R^n$$
having the property that $\sum_i\alpha_if_i$ is normal, for any $\alpha_i\in\mathbb R$.
\end{definition}

As a first illustration for this notion, dealing with the simplest case, let us start with the following straightforward result, coming from things that we know well:

\begin{theorem}
Assuming that we have independent variables $f_i\sim g_{t_i}^{a_i}$, the variable $f=(f_i)$ is Gaussian, and follows the following law:
$$g_T^a=\frac{1}{\sqrt{(2\pi)^nt_1\ldots t_n}}\,\exp\left(-\frac{(x_1-a_1)^2}{2t_1}-\ldots-\frac{(x_n-a_n)^2}{2t_n}\right)dx$$
In this situation the covariance matrix is diagonal, $\Sigma=diag(t_i)$.
\end{theorem}

\begin{proof}
The first assertion is clear, since a linear combination of independent normal variables is normal. Next, according to our formulae from chapter 2, we have:
$$g_{t_i}^{a_i}=\frac{1}{\sqrt{2\pi t_i}}e^{-(x-a_i)^2/2t_i}$$

Now by multiplying, we are led to the density formula in the statement. As for $\Sigma=diag(t_i)$, this is clear, since we have $f_i\sim g_{t_i}^{a_i}$, independent.
\end{proof}

Next, as a key generalization of the construction in Theorem 9.18, we have:

\begin{theorem}
Given a vector $h:X\to\mathbb R^m$ having normal independent entries, and a matrix $L\in M_{n\times m}(\mathbb R)$, the following variable is Gaussian,
$$f=Lh:X\to\mathbb R^n$$
with covariance matrix $\Sigma=LTL^t$, where $T=diag(V(h_i))$. In the case where $m=n$ and $L$ is the identity matrix, we recover the construction in Theorem 9.18. 
\end{theorem}

\begin{proof}
As before, the first assertion is clear, since a linear combination of independent normal variables is normal. Regarding the covariance matrix, we have:
\begin{eqnarray*}
E(f_if_j)
&=&E\left((Lh)_i(Lh)_j\right)\\
&=&\sum_{kl}L_{ik}L_{jl}E(h_kh_l)\\
&=&\sum_kL_{ik}L_{jk}E(h_k^2)+\sum_{k\neq l}L_{ik}L_{jl}E(h_k)E(h_l)
\end{eqnarray*}

By using now our formulae for shifted bells from chapter 2, we obtain:
\begin{eqnarray*}
E(f_if_j)
&=&\sum_kL_{ik}L_{jk}\left[T_{kk}+E(h_k)^2\right]+\sum_{k\neq l}L_{ik}L_{jl}E(h_k)E(h_l)\\
&=&\sum_kL_{ik}L_{jk}T_{kk}+\sum_{kl}L_{ik}L_{jl}E(h_k)E(h_l)\\
&=&\sum_kL_{ik}T_{kk}L_{jk}+\sum_kL_{ik}E(h_k)\sum_lL_{jl}E(h_l)\\
&=&(LTL^t)_{ij}+E(f_i)E(f_j)
\end{eqnarray*}

Thus $cov(f_i,f_j)=(LTL^t)_{ij}$, which shows that we have $\Sigma=LTL^t$, as claimed.
\end{proof}

Let us compute now the density of the variables in Theorem 9.19. In the simplest case, that where the covariance matrix is diagonal, the result is as follows:

\begin{theorem}
In the context of Theorem 9.19, assuming that $\Sigma$ is diagonal, the variables $f_1,\ldots,f_n$ follow to be independent, as in Theorem 9.18, and the density is
$$g_\Sigma^a=\frac{1}{\sqrt{(2\pi)^n\det\Sigma}}\,\exp\left(-\frac{<\Sigma^{-1}(x-a),x-a>}{2}\right)dx$$
the parameters being the covariance matrix $\Sigma\in M_n(\mathbb R)$, and the position $a\in\mathbb R^n$.
\end{theorem}

\begin{proof}
This generalizes again things that we know in 2D, as follows:

\medskip

(1) As starting point, we can use the standard fact, that we will leave as an exercise, that given a vector $h:X\to\mathbb R^m$ with entries $h_i\sim g_1$, independent, and an orthogonal matrix $U\in O_m$, the vector $k=Uh:X\to\mathbb R^m$ has also entries $k_i\sim g_1$, independent.

\medskip

(2) Now let prove the result, assuming $h_i\sim g_1$, and with the matrix $L\in M_{n\times m}(\mathbb R)$ being normalized too, as for all its rows to have norm 1. In order to apply the trick in (1), we would like to find an orthogonal matrix $U\in O_m$ which looks as follows:
$$U=\binom{L}{*}$$

(3) But, according to the Gram-Schmidt procedure, this is possible precisely when the rows of $L$ are orthogonal, which in turn corresponds to the condition $cov(f_i,f_j)=0$.

\medskip

(4) Summarizing, we proved our result, under the normalization assumptions in (2). As for the general case, this follows easily from this, say exercise for you.

\medskip

(5) Finally, the density formula in the statement is the one from Theorem 9.18, written a bit more conceptually, in terms of the covariance matrix $\Sigma\in M_n(\mathbb R)$.
\end{proof}

Quite remarkably, the density formula in Theorem 9.20 holds in fact without assumptions on $\Sigma$, and in fact for any Gaussian vector, the result being as follows:

\begin{theorem}
The density of a Gaussian vector $f:X\to\mathbb R^n$ is given by
$$g_\Sigma^a=\frac{1}{\sqrt{(2\pi)^n\det\Sigma}}\,\exp\left(-\frac{<\Sigma^{-1}(x-a),x-a>}{2}\right)dx$$
the parameters being the covariance matrix $\Sigma\in M_n(\mathbb R)$, and the position $a\in\mathbb R^n$.
\end{theorem}

\begin{proof}
This generalizes again things that we know in 2D, as follows:

\medskip

(1) According to Definition 9.17, our vector $f:X\to\mathbb R^n$ being Gaussian means that $\sum_i\alpha_if_i$ is normal, for any $\alpha_i\in\mathbb R$. In particular, $f_i\sim g_{t_i}^{a_i}$, for some $t_i,a_i$.

\medskip

(2) Next, let us study the variables $\sum_i\alpha_if_i$. Their variance is given by:
\begin{eqnarray*}
V\left(\sum_i\alpha_if_i\right)
&=&E\left[\left(\sum_i\alpha_if_i\right)^2\right]-\left[E\left(\sum_i\alpha_if_i\right)\right]^2\\
&=&\sum_{ij}\alpha_i\alpha_jE(f_if_j)-\sum_{ij}\alpha_i\alpha_jE(f_i)E(f_j)\\
&=&\sum_{ij}\alpha_i\alpha_jcov(f_i,f_j)\\
&=&\sum_{ij}\alpha_i\alpha_j\Sigma_{ij}
\end{eqnarray*}

We conclude that our abstract normality assumption on $\sum_i\alpha_if_i$ reformulates into the following more precise condition, telling us what its normal law exactly is:
$$\sum_i\alpha_if_i\sim g_s^b\quad:\quad s=\sum_{ij}\alpha_i\alpha_j\Sigma_{ij}\quad,\quad  b=\sum_i\alpha_ia_i$$

(3) Let us look now at the moments of $\sum_i\alpha_if_i$. These are given by:
$$M_k\left(\sum_i\alpha_if_i\right)=\sum_{r_1+\ldots+r_n=k}\binom{k}{r_1,\ldots,r_n}\alpha_1^{r_1}\ldots\alpha_n^{r_n}E(f_1^{r_1}\ldots f_n^{r_n})$$

On the other hand, according to the formula in (2), these moments are given by the following formula, with $\varphi$ being a certain 2-variable polynomial function:
$$M_k\left(\sum_i\alpha_if_i\right)=\varphi\left(\sum_i\alpha_ia_i\,,\,\sum_{ij}\alpha_i\alpha_j\Sigma_{ij}\right)$$

We conclude from this, by recurrence, even without knowing the precise formula of $\varphi$, that we can compute all the moments $E(f_1^{k_1}\ldots f_n^{k_n})$ in terms of $(\Sigma,a)$.

\medskip

(4) Thus the density that we are looking for depends only on $(\Sigma,a)$. But this means that we can assume that we are in the framework of Theorem 9.19, and with $m=n$, this being the smallest $m\in\mathbb N$ allowing us to model the covariance matrix $\Sigma$.

\medskip

(5) So, let us assume this, $f=Lh$ with $h_i\sim g_{t_i}^{a_i'}$ independent, and $L\in M_n(\mathbb R)$. According to Theorem 9.18 we have the following formula, with $T=diag(t_i)$:
$$h\sim\frac{1}{\sqrt{(2\pi)^n\det T}}\,\exp\left(-\frac{<T^{-1}(x-a'),x-a'>}{2}\right)dx$$

Now observe that for any function $\varphi:\mathbb R^n\to\mathbb R$, we have the following formula:
\begin{eqnarray*}
E(\varphi(f))
&=&E(\varphi(Lh))\\
&=&\frac{1}{\sqrt{(2\pi)^n\det T}}\int_{\mathbb R^n}\varphi(Lx)\exp\left(-\frac{<T^{-1}(x-a'),x-a'>}{2}\right)dx\\
&=&\frac{1}{\sqrt{(2\pi)^n\det T}}\int_{\mathbb R^n}\varphi(y)\exp\left(-\frac{<T^{-1}(L^{-1}y-a'),L^{-1}y-a'>}{2}\right)\frac{dy}{|\det L|}\\
&=&\frac{1}{\sqrt{(2\pi)^n\det(LTL^t)}}\int_{\mathbb R^n}\varphi(y)\exp\left(-\frac{<(LTL^t)^{-1}(y-a),y-a>}{2}\right)dy\\
&=&\frac{1}{\sqrt{(2\pi)^n\det\Sigma}}\int_{\mathbb R^n}\varphi(y)\exp\left(-\frac{<\Sigma^{-1}(y-a),y-a>}{2}\right)dy
\end{eqnarray*}

Thus, we are led to the density formula in the statement.
\end{proof}

And with this, end of our study of the Gaussian vectors. Good work that we did, we learned the basics, but for the efficient use of these beasts, say in questions from statistics, there are of course far more things to be learned. Exercise for you, I guess.

\section*{9c. Laplace operator}

Moving on, I don't know about you, but personally I find that too much multivariable calculus, to the point of reaching the Jacobian, as we just did, can be a bit frustrating, in the lack of some applications. So, browsing now through various physics and engineering manuals, which keep piling in my living room, here is what I got, for the two of us:

\begin{question}
We would like to solve the main equations in physics, namely:
\begin{enumerate}
\item Laplace: $\Delta\varphi=0$.

\item Heat: $\dot{\varphi}=\alpha\Delta\varphi$.

\item Waves: $\ddot{\varphi}=v^2\Delta\varphi$.

\item Schr\"odinger: $ih\dot{\varphi}=-\frac{h^2}{2m}\Delta\varphi+V\varphi$.
\end{enumerate}
\end{question}

As a first observation, it is not even clear what all these equations are about. I mean, what is that $\Delta$ symbol, then what are $\alpha,v,h,m,V$, and finally, once these equations mathematically making sense, what phenomenon each equation precisely describes.

\bigskip

In answer, and getting back now to my living room, for some further consultation of all that physics and engineering books, $\Delta=\sum_id^2/dx_i^2$ is the Laplace operator, sort of a numeric second derivative, in several variables, and regarding the equations:

\bigskip

(1) The Laplace equation $\Delta\varphi=0$ originally comes from electrostatics, describing steady-state equilibrium and potential fields in source-free regions.

\bigskip

(2) The heat equation $\dot{\varphi}=\alpha\Delta\varphi$ describes heat diffusion, through a medium having thermal diffusivity, or ``speed of diffusion'', if you prefer, $\alpha>0$. 

\bigskip

(3) The wave equation $\ddot{\varphi}=v^2\Delta\varphi$ describes the propagation of various waves, which can be mechanical or electromagnetic, at speed $v>0$.

\bigskip

(4) The Schr\"odinger equation $ih\dot{\varphi}=-\frac{h^2}{2m}\Delta\varphi+V\varphi$ describes the evolution of a wave function under a potential $V$, $m$ being the mass, and $h$ the reduced Planck constant.

\bigskip

Which sounds quite good, and getting to work now, we first need to talk about the Laplace operator $\Delta$, and the Laplace equation $\Delta\varphi=0$. As starting point, we have:

\index{second derivative}
\index{Hessian matrix}
\index{Taylor formula}

\begin{theorem}
The first two derivatives of $\varphi:\mathbb R^n\to\mathbb R$, making the formula
$$\varphi(x+t)\simeq\varphi(x)+\varphi'(x)t+\frac{<\varphi''(x)t,t>}{2}$$
work, with $t\simeq0$, are the following row vector and square matrix,
$$\varphi'(x)=\left(\frac{d\varphi}{dx_i}\right)_i\quad,\quad 
\varphi''(x)=\left(\frac{d^2\varphi}{dx_idx_j}\right)_{ij}$$
called respectively transposed gradient, and Hessian matrix.
\end{theorem}

\begin{proof}
This is indeed something very standard, that you surely know, and I even have proof for this, the Jacobian technology that we used in the proof of Theorem 9.21, and even ages ago, in the proof of Theorem 9.12, being more advanced than this. But of course, if in need of some help, you can have a look for instance at my continuous variable book \cite{ba2}, where all this, gradient, Hessian and Jacobian too, is explained in detail.
\end{proof}

Getting now to the Laplace operator, here is how this appears:

\index{Laplacian}
\index{Laplace operator}

\begin{theorem}
When needing in your problems a numeric second derivative, the trace of the Hessian matrix, called Laplacian, given by the formula
$$\Delta\varphi=\sum_{i=1}^n\frac{d^2\varphi}{dx_i^2}$$
is often what you need. Indeed, $\Delta\varphi$ measures how much different is $\varphi(x)$, compared to the average of $\varphi(y)$, with $y\simeq x$, with this being something very useful.
\end{theorem}

\begin{proof}
This is obviously something a bit heuristic, but good to know. Let us write the formula in Theorem 9.23 as such, and with $t\to-t$ too:
$$\varphi(x+t)\simeq\varphi(x)+\varphi'(x)t+\frac{<\varphi''(x)t,t>}{2}$$
$$\varphi(x-t)\simeq\varphi(x)-\varphi'(x)t+\frac{<\varphi''(x)t,t>}{2}$$

By making the average, we obtain from this the following formula:
$$\frac{\varphi(x+t)+\varphi(x-t)}{2}\simeq\varphi(x)+\frac{<\varphi''(x)t,t>}{2}$$

Now by thinking a bit, we are led to the conclusions in the statement.
\end{proof} 

With this understood, the problem is now, what can we say about the mathematics of $\Delta$? And we have here the following standard question, inspired by linear algebra: 

\begin{question}
The Laplace operator being linear, namely
$$\Delta(a\varphi+b\psi)=a\Delta\varphi+b\Delta\psi$$
what can we say about it, inspired by usual linear algebra?
\end{question}

In answer now, the space of functions $\varphi:\mathbb R^n\to\mathbb R$ on which $\Delta$ acts being infinite dimensional, the usual tools from linear algebra do not apply as such, and we must be extremely careful. In addition, if a function $\varphi:\mathbb R^n\to\mathbb R$ is twice differentiable, nothing will guarantee that the function $\Delta\varphi:\mathbb R^n\to\mathbb R$ is twice differentiable too. Thus, we have some issues with the domain and range of $\Delta$, regarded as linear operator. Not good.

\bigskip

So, shall we trash Question 9.25? Not so quick, because, remarkably, some magic comes at $n=2$ and higher in relation with complex analysis, according to:

\index{holomorphic function}
\index{harmonic function}

\begin{principle}
The functions $\varphi:\mathbb R^n\to\mathbb R$ which are $0$-eigenvectors of $\Delta$,
$$\Delta\varphi=0$$
called harmonic functions, have the following properties:
\begin{enumerate}
\item At $n=1$, nothing spectacular, these are just the linear functions.

\item At $n=2$, these are, locally, the real parts of holomorphic functions.

\item At $n\geq 3$, these still share many properties with the holomorphic functions. 
\end{enumerate}
\end{principle}

In order to understand this, or at least get introduced to it, let us first look at the case $n=2$. Here, any $\varphi:\mathbb R^2\to\mathbb R$ can be regarded as function $\varphi:\mathbb C\to\mathbb R$, depending on $z=x+iy$. Thus, it is natural to enlarge the attention to the functions $\varphi:\mathbb C\to\mathbb C$, and ask which of these functions are harmonic, $\Delta\varphi=0$. And here, we have:

\begin{theorem}
Any holomorphic function $\varphi:\mathbb C\to\mathbb C$, when regarded as function
$$\varphi:\mathbb R^2\to\mathbb C$$
is harmonic. Moreover, the conjugates $\bar{\varphi}$ of holomorphic functions are harmonic too.
\end{theorem}

\begin{proof}
The first assertion comes from the following computation, with $z=x+iy$:
\begin{eqnarray*}
\Delta z^k
&=&\frac{d^2z^k}{dx^2}+\frac{d^2z^k}{dy^2}\\
&=&\frac{d(kz^{k-1})}{dx}+\frac{d(ikz^{k-1})}{dy}\\
&=&k(k-1)z^{k-2}-k(k-1)z^{k-2}\\
&=&0
\end{eqnarray*}

As for the second assertion, this is clear, coming from $\Delta\bar{\varphi}=\overline{\Delta\varphi}$.
\end{proof}

Many more things can be said, along these lines, notably a proof of the assertion (2) in Principle 9.26, which is however a quite tough piece of mathematics, and then with a clarification of the assertion (3) too, from that same principle. See Rudin \cite{rud}.

\section*{9d. Heat diffusion}

Getting now to physics, and to Question 9.22, I would say to skip the Laplace equation, (1) there, on the grounds that we already did some mathematics for that, enough work, and go with a fresh topic, namely (2), heat diffusion. And here, we have:

\index{heat diffusion}
\index{heat propagation}
\index{heat equation}
\index{thermal diffusivity}
\index{lattice model}

\begin{theorem}
Heat diffusion is described by the heat equation
$$\dot{\varphi}=\alpha\Delta\varphi$$
where $\alpha>0$ is a constant, called thermal diffusivity of the medium.
\end{theorem}

\begin{proof}
The study here can be done by using a lattice model, as follows:

\medskip

(1) To start with, as an intuitive explanation, the equation formulated above tells us that the rate of change $\dot{\varphi}$ of the temperature of the material at any given point must be proportional, with proportionality factor $\alpha>0$, to the average difference of temperature between that given point and the surrounding material. Which sounds good.

\medskip

(2) Now the point is that we can more rigorously recover this equation by using a lattice model. Indeed, let us first assume that we are in the one-dimensional case, $n=1$. Here our model can be taken as follows, with distance $l>0$ between neighbors:
$$\xymatrix@R=10pt@C=20pt{
\ar@{-}[r]&\circ_{x-l}\ar@{-}[r]^l&\circ_x\ar@{-}[r]^l&\circ_{x+l}\ar@{-}[r]&
}$$

(3) In order to model heat diffusion, we have to implement the intuitive mechanism explained above, namely ``the rate of change of the temperature of the material at any given point must be proportional, with proportionality factor $\alpha>0$, to the average difference of temperature between that given point and the surrounding material''.

\medskip

(4) In practice, this leads to a condition as follows, expressing the change of the temperature $\varphi$, over a small period of time $\delta>0$:
$$\varphi(x,t+\delta)=\varphi(x,t)+\frac{\alpha\delta}{l^2}\sum_{x\sim y}\left[\varphi(y,t)-\varphi(x,t)\right]$$

To be more precise, we have made several assumptions here, as follows:

\medskip

-- General heat diffusion assumption: the change of temperature at any given point $x$ is proportional to the average over neighbors, $y\sim x$, of the differences $\varphi(y,t)-\varphi(x,t)$ between the temperatures at $x$, and at these neighbors $y$.

\medskip

-- Infinitesimal time and length conditions: in our model, the change of temperature at a given point $x$ is proportional to small period of time involved, $\delta>0$, and is inverse proportional to the square of the distance between neighbors, $l^2$.

\medskip

(5) Regarding these latter assumptions, the one regarding the proportionality with the time elapsed $\delta>0$ is something quite natural, physically speaking, and mathematically speaking too, because we can rewrite our equation as follows, making it clear that we have here an equation regarding the rate of change of temperature at $x$:
$$\frac{\varphi(x,t+\delta)-\varphi(x,t)}{\delta}=\frac{\alpha}{l^2}\sum_{x\sim y}\left[\varphi(y,t)-\varphi(x,t)\right]$$

As for the second assumption that we made above, namely inverse proportionality with $l^2$, this can be justified on physical grounds too, but again, perhaps the best is to do the math, which will show right away where this proportionality comes from. 

\medskip

(6) So, let us do the math. In the context of our 1D model the neighbors of $x$ are the points $x\pm l$, and so the equation that we wrote above takes the following form:
$$\frac{\varphi(x,t+\delta)-\varphi(x,t)}{\delta}=\frac{\alpha}{l^2}\Big[(\varphi(x+l,t)-\varphi(x,t))+(\varphi(x-l,t)-\varphi(x,t))\Big]$$

Now observe that we can write this equation as follows:
$$\frac{\varphi(x,t+\delta)-\varphi(x,t)}{\delta}
=\alpha\cdot\frac{\varphi(x+l,t)-2\varphi(x,t)+\varphi(x-l,t)}{l^2}$$

(7) We recognize on the right the usual approximation of the second derivative, coming from calculus, and more specifically from the Taylor formula, at order 2. Thus, when taking the continuous limit of our model, $l\to 0$, we obtain the following equation:
$$\frac{\varphi(x,t+\delta)-\varphi(x,t)}{\delta}
=\alpha\cdot\varphi''(x,t)$$

Now with $t\to0$, we are led in this way to $\dot{\varphi}(x,t)=\alpha\cdot\varphi''(x,t)$, as desired.

\medskip

(8) With this done, let us discuss now 2 dimensions. Here we can use a similar lattice model, as follows, with all lengths being taken $l>0$, for simplifying:
$$
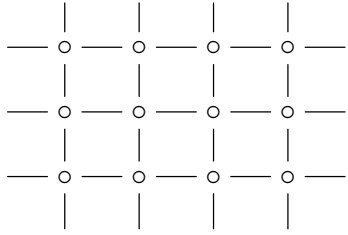
\xymatrix@R=12pt@C=15pt{
&\ar@{-}[d]&\ar@{-}[d]&\ar@{-}[d]&\ar@{-}[d]\\
\ar@{-}[r]&\circ\ar@{-}[r]\ar@{-}[d]&\circ\ar@{-}[r]\ar@{-}[d]&\circ\ar@{-}[r]\ar@{-}[d]&\circ\ar@{-}[r]\ar@{-}[d]&\\
\ar@{-}[r]&\circ\ar@{-}[r]\ar@{-}[d]&\circ\ar@{-}[r]\ar@{-}[d]&\circ\ar@{-}[r]\ar@{-}[d]&\circ\ar@{-}[r]\ar@{-}[d]&\\
\ar@{-}[r]&\circ\ar@{-}[r]\ar@{-}[d]&\circ\ar@{-}[r]\ar@{-}[d]&\circ\ar@{-}[r]\ar@{-}[d]&\circ\ar@{-}[r]\ar@{-}[d]&\\
&&&&&
}$$

(9) We have to implement now the physical heat diffusion mechanism, namely ``the rate of change of the temperature of the material at any given point must be proportional, with proportionality factor $\alpha>0$, to the average difference of temperature between that given point and the surrounding material''. In practice, this leads to a condition as follows, expressing the change of the temperature $\varphi$, over a small period of time $\delta>0$:
$$\varphi(x,y,t+\delta)=\varphi(x,y,t)+\frac{\alpha\delta}{l^2}\sum_{(x,y)\sim(u,v)}\left[\varphi(u,v,t)-\varphi(x,y,t)\right]$$

In fact, we can rewrite our equation as follows, making it clear that we have here an equation regarding the rate of change of temperature at $x$:
$$\frac{\varphi(x,y,t+\delta)-\varphi(x,y,t)}{\delta}=\frac{\alpha}{l^2}\sum_{(x,y)\sim(u,v)}\left[\varphi(u,v,t)-\varphi(x,y,t)\right]$$

(10) So, let us do the math. In the context of our 2D model the neighbors of $x$ are the points $(x\pm l,y\pm l)$, so the equation above takes the following form:
\begin{eqnarray*}
&&\frac{\varphi(x,y,t+\delta)-\varphi(x,y,t)}{\delta}\\
&=&\frac{\alpha}{l^2}\Big[(\varphi(x+l,y,t)-\varphi(x,y,t))+(\varphi(x-l,y,t)-\varphi(x,y,t))\Big]\\
&+&\frac{\alpha}{l^2}\Big[(\varphi(x,y+l,t)-\varphi(x,y,t))+(\varphi(x,y-l,t)-\varphi(x,y,t))\Big]
\end{eqnarray*}

Now observe that we can write this equation as follows:
\begin{eqnarray*}
\frac{\varphi(x,y,t+\delta)-\varphi(x,y,t)}{\delta}
&=&\alpha\cdot\frac{\varphi(x+l,y,t)-2\varphi(x,y,t)+\varphi(x-l,y,t)}{l^2}\\
&+&\alpha\cdot\frac{\varphi(x,y+l,t)-2\varphi(x,y,t)+\varphi(x,y-l,t)}{l^2}
\end{eqnarray*}

(11) As it was the case before in one dimension, we recognize on the right the usual approximation of the second derivative, coming from calculus. Thus, when taking the continuous limit of our model, $l\to 0$, we obtain the following equation:
$$\frac{\varphi(x,y,t+\delta)-\varphi(x,y,t)}{\delta}
=\alpha\left(\frac{d^2\varphi}{dx^2}+\frac{d^2\varphi}{dy^2}\right)(x,y,t)$$

Now with $t\to0$, this gives $\dot{\varphi}(x,y,t)=\alpha\Delta\varphi(x,y,t)$, as desired. 

\medskip

(12) Finally, in arbitrary $n$ dimensions the same argument carries over, namely a straightforward lattice model, and we will leave this as an exercise.
\end{proof}

Regarding now the resolution of the heat equation, we have here:

\index{heat kernel}
\index{initial condition}

\begin{theorem}
The heat diffusion equation, $\dot{\varphi}=\alpha\Delta\varphi$ with $\alpha>0$, with initial condition $\varphi(x,0)=f(x)$, has as solution the function
$$\varphi(x,t)=\int_{\mathbb R^n}K_t(x-y)f(y)dy$$
where the function $K_t:\mathbb R^n\to\mathbb R$, called heat kernel, given by
$$K_t(x)=\frac{1}{\sqrt{(4\pi\alpha t)^n}}\,e^{-||x||^2/4\alpha t}$$
is the standard solution, coming from the initial data $f=\delta_0$, Dirac mass at $0$.
\end{theorem}

\begin{proof}
Observe first that $K_t$ is the standard Gaussian density on $\mathbb R^n$, but we will keep calling it heat kernel, like everyone does. As for the proof, this goes as follows:

\medskip

(1) Let us first discuss what happens in 1 dimension. We first have to check that the heat kernel $K_t$ is indeed a solution. The time derivative is:
$$\dot{K_t}=-\frac{1}{2t\sqrt{4\pi\alpha t}}\,e^{-x^2/4\alpha t}+\frac{x^2}{4\alpha t^2\sqrt{4\pi\alpha t}}\,e^{-x^2/4\alpha t}$$

Regarding the first space derivative, this is given by the following formula:
$$K_t'=-\frac{x}{2\alpha t\sqrt{4\pi\alpha t}}\,e^{-x^2/4\alpha t}$$

As for the second space derivative, this is given by the following formula:
$$K_t''=-\frac{1}{2\alpha t\sqrt{4\pi\alpha t}}\,e^{-x^2/4\alpha t}+\frac{x^2}{4\alpha^2t^2\sqrt{4\pi\alpha t}}\,e^{-x^2/4\alpha t}$$

Thus, we can see that the heat equation $\dot{\varphi}=\alpha\varphi''$ is indeed satisfied.

\medskip

(2) Next, let us convolve $K_t$ with an arbitrary function $f$, as follows:
$$\varphi(x,t)=\int_{\mathbb R^N}K_t(x-y)f(y)dy$$

The point now is that, when doing so, all the computations from (1) will perturb well, according to the following formulae, which are both clear:
$$\dot{\varphi}(x,t)=\int_\mathbb R\dot{K}_t(x-y)f(y)dy$$
$$\varphi''(x,t)=\int_\mathbb RK''_t(x-y)f(y)dy$$

Thus, we can see that the heat equation $\dot{\varphi}=\alpha\varphi''$ is again satisfied. 

\medskip

(3) Next, by using the fact that $K_t$ has mass 1, it is routine to check that $K_t$ comes indeed from the simplest situation, that of a radiating point body placed at $0$, and then that the solution found in (2) by convolving comes indeed from the initial data $f$.

\medskip

(4) Finally, everything extends well to arbitrary dimensions, and we will leave this as an exercise. As for the uniqueness issues, we will leave this as an exercise too.
\end{proof}

As a conclusion to this, the 3D normal laws, also called heat kernel, are smart enough for solving the heat diffusion equation, in our usual 3D. Which is quite remarkable.

\section*{9e. Exercises}

This was a standard introduction to Gaussian vectors, and as exercises, we have:

\begin{exercise}
Do some computations of covariances and correlations.
\end{exercise}

\begin{exercise}
Learn more about technical CLT, producing Gaussian vectors.
\end{exercise}

\begin{exercise}
Learn more about positive matrices, and their properties.
\end{exercise}

\begin{exercise}
Recapture the density of Gaussian vectors, by other means.
\end{exercise}

\begin{exercise}
Clarify what we said, regarding the philosophy of $\Delta$.
\end{exercise}

\begin{exercise}
Learn about harmonic functions, as much as you can.
\end{exercise}

\begin{exercise}
Work out the final details, for the heat equation in 1D.
\end{exercise}

\begin{exercise}
Work out the details for the heat equation, in $n$ dimensions.
\end{exercise}

As bonus exercise, have a look as well at the wave and Schr\"odinger equations.

\chapter{Chi and gamma}

\section*{10a. Chi squared}

We have seen quite a deal of multivariable calculus and probability in chapter 9, in relation with the variables $f:X\to\mathbb R^n$ studied there. In this chapter we get back to 1 dimension, by using a simple idea, namely looking at the laws of the corresponding real variables $||f||^2:X\to\mathbb R$ and $||f||:X\to\mathbb R$, called chi squared and chi laws.

\bigskip

In fact, we already met this idea in chapter 6, at the particular value $n=2$, when studying the Rayleigh laws. So, this chapter will be basically a remake of chapter 6, by adding as much parameters as possible, coming from our work from chapter 9, notably with the dimension parameter $n\in\mathbb N$. And in the end, we will talk more in detail about the case $n=3$, where many interesting things happen. This will be our plan.

\bigskip

Getting started, our previous experience from chapters 6 and 9 shows that many parameters can be a quite tricky business. So, we will do things slowly, by introducing them gradually. Let us start with the following definition, simple and bright:

\index{chi}
\index{chi squared}
\index{Rayleigh law}

\begin{definition}
The chi squared law, depending on $n\in\mathbb N$, is given by
$$\chi_n^2=law\Big(f_1^2+\ldots+f_n^2\Big)$$
with $f_1,\ldots,f_n$ being independent, each following $g_1$. We can talk as well about
$$\chi_n=law\left(\sqrt{f_1^2+\ldots+f_n^2}\right)$$
called chi law, whose simplest $n>1$ instance, $\chi_2$, is the Rayleigh law.
\end{definition}

So, these laws, and their various parametric versions, will be our objects of interest, in this chapter. As in chapter 6, we will first study $\chi_n^2$, which is simpler, and leave $\chi_n$ for later. Getting to work now, as a first and main result regarding $\chi_n^2$, we have:

\index{gamma function}

\begin{theorem}
The density of the chi squared law is given by
$$\chi_n^2=\frac{x^{n/2-1}e^{-x/2}}{2^{n/2}\Gamma(n/2)}\,dx$$
on $[0,\infty)$, with $\Gamma$ being as usual the gamma function.
\end{theorem}

\begin{proof}
This is something quite tricky, the idea being as follows:

\medskip

(1) To start with, and standing as a useful complement to what the statement says, let us recall from chapter 6 that the gamma function is constructed as follows:
$$\Gamma(s)=\int_0^\infty x^{s-1}e^{-x}\,dx$$

By partial integration we have $\Gamma(s+1)=s\Gamma(s)$, and since $\Gamma(1)=1$, trivially, and $\Gamma(1/2)=\sqrt{\pi}$, obtained via $x=y^2$, we have the following formulae, for $N\in\mathbb N$:
$$\Gamma(N)=(N-1)!\quad,\quad\Gamma\left(N+\frac{1}{2}\right)=\frac{(2N)!!}{2^N}\,\sqrt{\pi}$$

Alternatively, we have the following formula, with $c=\sqrt{2},\sqrt{\pi}$ for $n$ even, odd:
$$\Gamma\left(\frac{n}{2}\right)=\frac{(n-1)!!}{2^{(n-1)/2}}\,c$$

(2) Getting now to what our statement says, with $f_1,\ldots,f_n$ being independent variables, all following the normal law $g_1$, we have:
\begin{eqnarray*}
M_k(\chi_n^2)
&=&E((f_1^2+\ldots+f_n^2)^k)\\
&=&\sum_{i_1+\ldots+i_n=k}\binom{k}{i_1,\ldots,i_n}E(f_1^{2i_1})\ldots E(f_n^{2i_n})\\
&=&\sum_{i_1+\ldots+i_n=k}\frac{k!}{i_1!\ldots i_n!}(2i_1)!!\ldots(2i_n)!!\\
&=&\sum_{i_1+\ldots+i_n=k}\frac{k!}{i_1!\ldots i_n!}\cdot\frac{(2i_1)!}{2^{i_1}i_1!}\cdots\frac{(2i_n)!}{2^{i_n}i_n!}\\
&=&\frac{k!}{2^k}\sum_{i_1+\ldots+i_n=k}\frac{(2i_1)!}{i_1!i_1!}\cdots\frac{(2i_n)!}{i_n!i_n!}\\
&=&\frac{k!}{2^k}\sum_{i_1+\ldots+i_n=k}\binom{2i_1}{i_1}\ldots\binom{2i_n}{i_n}
\end{eqnarray*}

In order to finish this computation, we can use the following formula:
$$\frac{1}{\sqrt{1-4z}}=\sum_{i=0}^\infty\binom{2i}{i}z^i$$

Indeed, by taking the $n$-th power of this series, we obtain the following formula:
$$\frac{1}{\sqrt{(1-4z)^n}}
=\sum_kz^k\sum_{i_1+\ldots+i_n=k}\binom{2i_1}{i_1}\ldots\binom{2i_n}{i_n}$$

Now by putting everything together, the conclusion is that we have:
$$\frac{1}{\sqrt{(1-4z)^n}}
=\sum_{k\geq0}\frac{(2z)^k}{k!}\,M_k(\chi_n^2)$$

But the function on the left can be computed by using the generalized binomial formula, with exponent $p=-n/2$, and we obtain in this way:
\begin{eqnarray*}
\frac{1}{\sqrt{(1-4z)^n}}
&=&\sum_{k\geq0}\binom{-n/2}{k}(-4z)^k\\
&=&\sum_{k\geq0}\frac{(-n/2)(-n/2-1)\ldots(-n/2-k+1)}{k!}(-4z)^k\\
&=&\sum_{k\geq0}\frac{n(n+2)\ldots(n+2k-2)}{k!}(2z)^k
\end{eqnarray*}

And with this, good news, done with our moment computation, which yields:
$$M_k(\chi_n^2)=n(n+2)\ldots(n+2k-2)$$

(3) However, not time to celebrate yet, because we still have to do the second part of the work, namely computation of the moments for the predicted density, which is:
$$\chi_n^2=\frac{x^{n/2-1}e^{-x/2}}{2^{n/2}\Gamma(n/2)}\,dx$$

By using the definition of the gamma function, as formulated in (1), and then several times the formula $\Gamma(s+1)=s\Gamma(s)$, the moments of this predicted density are:
\begin{eqnarray*}
I_k
&=&\frac{1}{2^{n/2}\Gamma(n/2)}\int_0^\infty x^{n/2+k-1}e^{-x/2}\,dx\\
&=&\frac{1}{2^{n/2}\Gamma(n/2)}\int_0^\infty (2y)^{n/2+k-1}e^{-y}\,2dy\\
&=&\frac{2^k}{\Gamma(n/2)}\int_0^\infty y^{n/2+k-1}e^{-y}\,dy\\
&=&2^k\cdot\frac{\Gamma(n/2+k)}{\Gamma(n/2)}\\
&=&2^k\cdot\frac{\Gamma(n/2+k)}{\Gamma(n/2+k-1)}\cdot\frac{\Gamma(n/2+k-1)}{\Gamma(n/2+k-2)}\ldots\frac{\Gamma(n/2+1)}{\Gamma(n/2)}\\
&=&2^k\left(\frac{n}{2}+k-1\right)\left(\frac{n}{2}+k-2\right)\ldots\left(\frac{n}{2}\right)\\
&=&n(n+2)\ldots(n+2k-2)
\end{eqnarray*}

We conclude that we have $I_k=M_k(\chi_n^2)$, so job done, and theorem proved.
\end{proof}

As a comment now, you might perhaps say, why doing the above 2 pages of moment computations, instead of simply getting the density of $||f||^2$ from that of $f$, that we know well from chapter 9. Good point, so let us develop now this idea of yours:

\begin{theorem}[again]
The density of the chi squared law is given by
$$\chi_n^2=\frac{x^{n/2-1}e^{-x/2}}{2^{n/2}\Gamma(n/2)}\,dx$$
on $[0,\infty)$, this time coming via a direct density transport.
\end{theorem}

\begin{proof}
Buckle up, the computations will be a bit harder than expected:

\medskip

(1) Consider a vector $f=(f_i)$ whose components follow $g_1$ and are independent, so that $||f||^2\sim\chi_n^2$. We have the following formula, for any function $\varphi:[0,\infty)\to\mathbb R$:
$$E\left(\varphi(||f||^2)\right)
=\frac{1}{\sqrt{(2\pi)^n}}\int_{\mathbb R^n}\varphi(||x||^2)e^{-||x||^2/2}dx$$

Thus, we have our density, provided that we can compute the integral on the right.

\medskip

(2) In order to do so, we can use spherical coordinates. There are many possible conventions here, and we will use in this book the following nice-looking formulae:
$$\begin{cases}
x_1\!\!\!&=\ r\cos t_1\\
x_2\!\!\!&=\ r\sin t_1\cos t_2\\
\vdots\\
x_{n-1}\!\!\!&=\ r\sin t_1\sin t_2\ldots\sin t_{n-2}\cos t_{n-1}\\
x_n\!\!\!&=\ r\sin t_1\sin t_2\ldots\sin t_{n-2}\sin t_{n-1}
\end{cases}$$

We will also need the Jacobian of these coordinates. By developing over the last column, we obtain the following formula for this Jacobian $J_n$:
\begin{eqnarray*}
J_n
&=&r\sin t_1\ldots\sin t_{n-2}\sin t_{n-1}\times \sin t_{n-1}J_{n-1}\\
&+&r\sin t_1\ldots \sin t_{n-2}\cos t_{n-1}\times\cos t_{n-1}J_{n-1}\\
&=&r\sin t_1\ldots\sin t_{n-2}(\sin^2 t_{n-1}+\cos^2 t_{n-1})J_{n-1}\\
&=&r\sin t_1\ldots\sin t_{n-2}J_{n-1}
\end{eqnarray*}

Thus, by recurrence, we conclude that the Jacobian is given by:
$$J_n=r^{n-1}\sin^{n-2}t_1\sin^{n-3}t_2\,\ldots\,\sin^2t_{n-3}\sin t_{n-2}$$

Observe that this formula fits with the well-known formulae at $n=1,2,3$, that you surely know from physics class, namely $J_1=1$, $J_2=r$, $J_3=r^2\sin t_1$.

\medskip

(3) Getting back to our computation (1), this can be continued as follows, with the $2^n$ factor coming from the fact that we restricted attention to $x=(x_i)$ with $x_i\geq0$:
\begin{eqnarray*}
E\left(\varphi(||f||^2)\right)
&=&\frac{1}{\sqrt{(2\pi)^n}}\int_{\mathbb R^n}\varphi(||x||^2)e^{-||x||^2/2}dx\\
&=&\frac{2^n}{\sqrt{(2\pi)^n}}\int_0^{\pi/2}\ldots\int_0^{\pi/2}\int_0^\infty\varphi(r^2)e^{-r^2/2}\times r^{n-1}\\
&&\times\ \sin^{n-2}t_1\sin^{n-3}t_2\,\ldots\,\sin^2t_{n-3}\sin t_{n-2}\,drdt_1\ldots dt_{n-1}\\
&=&\left(\frac{2}{\pi}\right)^{n/2}\int_0^{\pi/2}\sin^{n-2}t_1dt_1\ldots\ldots\int_0^{\pi/2}\sin^0t_{n-1}dt_{n-1}\\
&&\times\int_0^\infty\varphi(r^2)r^{n-1}e^{-r^2/2}dr
\end{eqnarray*}

(4) Now if we denote by $C$ the above product of trigonometric integrals, with the normalization factor in front of it included, we have:
\begin{eqnarray*}
E\left(\varphi(||f||^2)\right)
&=&C\int_0^\infty\varphi(r^2)r^{n-1}e^{-r^2/2}dr\\
&=&C\int_0^\infty\varphi(y)y^{(n-1)/2}e^{-y/2}\frac{dy}{2\sqrt{y}}\\
&=&\frac{C}{2}\int_0^\infty\varphi(y)y^{n/2-1}e^{-y/2}dy
\end{eqnarray*}

(5) As a comment now, the trigonometric factor $C$ is more or less the same thing as the volume $V$ of the unit sphere, as shown by the following computation:
\begin{eqnarray*}
V
&=&2^n\int_0^{\pi/2}\ldots\int_0^{\pi/2}\int_0^1r^{n-1}\\
&&\times\ \sin^{n-2}t_1\sin^{n-3}t_2\,\ldots\,\sin^2t_{n-3}\sin t_{n-2}\,drdt_1\ldots dt_{n-1}\\
&=&\frac{2^n}{n}\int_0^{\pi/2}\sin^{n-2}t_1dt_1\ldots\ldots\int_0^{\pi/2}\sin^0t_{n-1}dt_{n-1}
\end{eqnarray*}

Even better, observe that the area $A$ of the unit sphere is given by:
$$A=nV=2^n\int_0^{\pi/2}\sin^{n-2}t_1dt_1\ldots\ldots\int_0^{\pi/2}\sin^0t_{n-1}dt_{n-1}$$

But with this, our computation in (3) simply reads, with $\psi(r)=\varphi(r^2)e^{-r^2/2}$:
$$\int_{\mathbb R^n}\psi(||x||)dx=A\int_0^\infty\psi(r)r^{n-1}dr$$

In fact, and we will leave this as an instructive geometry exercise, it is possible to come upon this latter formula, valid for any $\psi$, just by thinking. Enjoy.

\medskip

(6) Back to our business now, regardless of the method used, we need to compute that integrals of sines. I mean, even when short-circuiting the spherical coordinates, via some direct geometry, as suggested in (5), this would lead us into the computation of $V$. And here, regardless of the method used, which can be spherical coordinates, or recurrence by slicing, you are led in practice to computing the above integrals of sines. 

\medskip

(7) So, trigonometric integrals. Following Wallis, our claim is that we have the following formulae, where $\varepsilon(p)=1$ if $p$ is even, and $\varepsilon(p)=0$ if $p$ is odd:
$$\int_0^{\pi/2}\cos^pt\,dt=\int_0^{\pi/2}\sin^pt\,dt=\left(\frac{\pi}{2}\right)^{\varepsilon(p)}\frac{p!!}{(p+1)!!}$$

Indeed, our integrals being equal, let us look at the one on the left $I_p$. We have:
\begin{eqnarray*}
(\cos^pt\sin t)'
&=&p\cos^{p-1}t(-\sin t)\sin t+\cos^pt\cos t\\
&=&p\cos^{p+1}t-p\cos^{p-1}t+\cos^{p+1}t\\
&=&(p+1)\cos^{p+1}t-p\cos^{p-1}t
\end{eqnarray*}

By integrating now between $0$ and $\pi/2$, we obtain the following formula:
$$(p+1)I_{p+1}=pI_{p-1}$$

Thus we can compute $I_p$ by recurrence, and we obtain in this way:
\begin{eqnarray*}
I_p
&=&\frac{p-1}{p}\,I_{p-2}\\
&=&\frac{p-1}{p}\cdot\frac{p-3}{p-2}\,I_{p-4}\\
&\vdots&\\
&=&\frac{p!!}{(p+1)!!}\,I_{1-\varepsilon(p)}
\end{eqnarray*}

Now the initial data being $I_0=\frac{\pi}{2}$ and $I_1=1$, we obtain the result.

\medskip

(8) And with this, good news, we can finish. The factor $C$ from (4) is given by:
\begin{eqnarray*}
C
&=&\left(\frac{2}{\pi}\right)^{n/2}\int_0^{\pi/2}\sin^{n-2}t_1dt_1\ldots\ldots\int_0^{\pi/2}\sin^0t_{n-1}dt_{n-1}\\
&=&\left(\frac{2}{\pi}\right)^{n/2}\left(\frac{\pi}{2}\right)^{\varepsilon(n-2)+\ldots+\varepsilon(0)}\frac{(n-2)!!}{(n-1)!!}\cdot\frac{(n-3)!!}{(n-2)!!}\ldots\ldots\frac{0!!}{1!!}\\
&=&\left(\frac{2}{\pi}\right)^{n/2}\left(\frac{\pi}{2}\right)^{[n/2]}\frac{1}{(n-1)!!}\\
&=&\left(\frac{2}{\pi}\right)^{\varepsilon(n)/2}\frac{1}{(n-1)!!}
\end{eqnarray*}

(9) Thus, by getting back now to our formula in (4), that reads:
$$E\left(\varphi(||f||^2)\right)
=\frac{1}{2}\left(\frac{2}{\pi}\right)^{\varepsilon(n)/2}\frac{1}{(n-1)!!}\int_0^\infty\varphi(y)y^{n/2-1}e^{-y/2}dy$$

In order to further process the normalization factor, recall from the proof of Theorem 10.2 that we have the following formula, with $c=\sqrt{2},\sqrt{\pi}$ for $n$ even, odd:
$$\Gamma\left(\frac{n}{2}\right)=\frac{(n-1)!!}{2^{(n-1)/2}}\,c$$

Equivalently, with $\varepsilon(n)=1$ if $n$ is even, and $\varepsilon(n)=0$ if $n$ is odd, we have:
$$2^{n/2}\Gamma\left(\frac{n}{2}\right)=(n-1)!!\,\sqrt{2}\,c=(n-1)!!\,2\left(\frac{\pi}{2}\right)^{\varepsilon(n)/2}$$

Thus, in terms of the gamma function, our expectation formula reads:
$$E\left(\varphi(||f||^2)\right)=\frac{1}{2^{n/2}\Gamma(n/2)}\int_0^\infty\varphi(y)y^{n/2-1}e^{-y/2}dy$$

We are therefore led to the density in the statement.

\medskip

(10) Finally, as a byproduct of our study, let us record some useful formulae, for the volume $V$ and area $A$ of the unit sphere in $\mathbb R^n$, obtained via (5), thanks to the trigonometric integral computed in (8). According to our computations, these are:
$$V=\left(\frac{\pi}{2}\right)^{[n/2]}\frac{2^n}{(n+1)!!}\quad,\quad 
A=\left(\frac{\pi}{2}\right)^{[n/2]}\frac{2^n}{(n-1)!!}$$

Equivalently, in terms of the gamma function, the formulae are a bit more complicated, as follows, with our usual convention $c=\sqrt{2},\sqrt{\pi}$ for $n$ even, odd:
$$V=\left(\frac{\pi}{2}\right)^{[n/2]}\frac{2^{(n-1)/2}c}{\Gamma(n/2+1)}\quad,\quad 
A=\left(\frac{\pi}{2}\right)^{[n/2]}\frac{2^{(n+1)/2}c}{\Gamma(n/2)}$$

Finally, for the story to be complete, let us record as well the Stirling asymptotics:
$$V\simeq\left(\frac{2\pi e}{n}\right)^{n/2}\frac{1}{\sqrt{\pi n}}
\quad,\quad A\simeq\left(\frac{2\pi e}{n}\right)^{n/2}\sqrt{\frac{n}{\pi}}$$

And with of course exercise for you, to learn more about all this, spheres.
\end{proof}

Quite interesting all this, and looking back at what we did, the question comes, in the end, which of Theorems 10.2 and 10.3 is the best? Not an easy question, and as usual in such philosophical situations, time to ask the rat. Who advices the following:

\begin{rat}
Remember, combinatorics and analysis are the same thing, due to the binomial formula $(a+b)^n=\sum_k\binom{n}{k}a^kb^{n-k}$, which produces them both.
\end{rat}

Well, thanks rat, but I must admit that this sounds to me a bit too advanced. I mean, I can certainly understand that at your level of knowledge, Theorems 10.2 and 10.3 are exactly the same thing. But do I have as much time available as you do, for meditating at all this. Remember after all that I am feeding us both, and it would be impossible for me to do that, if spending my days in tranquility, hidden under the stove.

\bigskip

So, going now with cat, for a piece of complementary advice, here is what she says:

\begin{cat}
Efficiency. With combinatorics you got the density and moments in 2 pages, with analysis you got the density only, in 4 pages.
\end{cat}

Which sounds good, thanks cat, cannot argue with this, indeed, and Theorem 10.2 will remain our official one, on the subject. Moving on, based on that theorem, and by adding a few missing things, here is our grand result, about the chi squared laws:

\begin{theorem}
The chi squared laws have the following properties: 
\begin{enumerate}
\item $\chi_n^2=law(f_1^2+\ldots+f_n^2)$, with $f_i\sim g_1$, independent.

\item The density is given by $\chi_n^2=\frac{x^{n/2-1}e^{-x/2}}{2^{n/2}\Gamma(n/2)}\,dx$.

\item The moments are $M_k=n(n+2)\ldots(n+2k-2)$.

\item We have $E=n$, $V=2n$, $\gamma=\sqrt{8/n}$, $\kappa=3+12/n$.

\item The Fourier transform is $F(x)=(1-2ix)^{-n/2}$.
\end{enumerate}
\end{theorem}

\begin{proof}
Here (1) is the definition of the chi squared laws, as formulated before, and (2) and (3) are things that we know well, from Theorem 10.2 and its proof. Regarding now (4), according to our general moment formula in (3), we have:
$$M_1=n$$
$$M_2=n(n+2)$$
$$M_3=n(n+2)(n+4)$$
$$M_4=n(n+2)(n+4)(n+6)$$

Thus we have the expectation $E=n$, and then the variance is given by:
$$V=(n^2+2n)-n^2=2n$$

In order to compute now $\gamma$ and $\kappa$, observe that the needed central moments are:
$$M_3'=M_3-3EM_2+2E^3=8n$$
$$M_4'=M_4-4EM_3+6E^2M_2-3E^4=12n(n+4)$$

Thus, the skewness and kurtotis are given by the following formulae:
$$\gamma=\frac{8n}{2n\sqrt{2n}}=\sqrt{\frac{8}{n}}\quad,\quad\kappa=\frac{12n(n+4)}{4n^2}=3+\frac{12}{n}$$

Finally, regarding (5), the Fourier transform computation is as follows:
\begin{eqnarray*}
F(x)
&=&\frac{1}{2^{n/2}\Gamma(n/2)}\int_0^\infty y^{n/2-1}e^{-y/2+ixy}\,dy\\
&=&\frac{1}{2^{n/2}\Gamma(n/2)}\int_0^\infty y^{n/2-1}e^{(2ix-1)y/2}\,dy\\
&=&\frac{1}{2^{n/2}\Gamma(n/2)}\int_0^\infty\left(\frac{z}{1-2ix}\right)^{n/2-1}e^{-z/2}\,\frac{dz}{1-2ix}\\
&=&(1-2ix)^{-n/2}\cdot\frac{1}{2^{n/2}\Gamma(n/2)}\int_0^\infty z^{n/2-1}e^{-z/2}\,dz\\
&=&(1-2ix)^{-n/2}
\end{eqnarray*}

Thus, we are led to the conclusions in the statement.
\end{proof}

\section*{10b. Parametric versions}

Getting now to the parametric versions of Theorem 10.6, many things to be done. Let us start with a straightforward generalization, obtained via a basic rescaling:

\begin{theorem}
Assuming that $f_1,\ldots,f_n\sim g_t$ are independent, $||f||^2$ follows
$$\chi_{nt}^2=\frac{x^{n/2-1}e^{-x/2t}}{(2t)^{n/2}\Gamma(n/2)}\,dx$$
and the basic properties of this law are as follows:
\begin{enumerate}
\item The moments are $M_k=t^kn(n+2)\ldots(n+2k-2)$.

\item We have $E=nt$, $V=2nt^2$, $\gamma=\sqrt{8/n}$, $\kappa=3+12/n$.

\item The Fourier transform is $F(x)=(1-2ixt)^{-n/2}$.
\end{enumerate}
\end{theorem}

\begin{proof}
This is indeed a basic rescaling of Theorem 10.6, will all the assertions coming via $h_i\sim g_1\implies\sqrt{t}h_i\sim g_t$. Indeed, the density formula comes from:
\begin{eqnarray*}
E(\varphi(th_1^2+\ldots+th_n^2))
&=&\frac{1}{2^{n/2}\Gamma(n/2)}\int_0^\infty\varphi(tx)x^{n/2-1}e^{-x/2}\,dx\\
&=&\frac{1}{2^{n/2}\Gamma(n/2)}\int_0^\infty\varphi(y)(y/t)^{n/2-1}e^{-y/2t}\,\frac{dy}{t}\\
&=&\frac{1}{(2t)^{n/2}\Gamma(n/2)}\int_0^\infty\varphi(y)y^{n/2-1}e^{-y/2t}\,dy
\end{eqnarray*}

As for the other formulae, moments and Fourier, these are clear as well. Finally, observe that all this fits with our previous results from chapter 6, at $n=1,2$.
\end{proof}

More interestingly now, motivated by the fact that the gamma function notoriously works over the whole $(0,\infty)$, a quick examination of what we did, at the end of the proof of Theorem 10.2 and afterwards, shows that our computations work in fact for any $n>0$ real. To be more precise, with the notations $n=2m$ and $t=s/2$, which are both standard, we are led in this way to the following extension of Theorem 10.7:

\index{gamma distribution}

\begin{theorem}
The gamma distribution of parameters $m,s>0$, given by
$$\mu_{ms}=\frac{x^{m-1}e^{-x/s}}{s^m\Gamma(m)}\,dx$$
generalizes the chi squared laws, $\chi^2_{nt}=\mu_{n/2,2t}$, and has the following properties: 
\begin{enumerate}
\item The moments are $M_k=s^km(m+1)\ldots(m+k-1)$.

\item We have $E=sm$, $V=s^2m$, $\gamma=2/\sqrt{m}$, $\kappa=3+6/m$.

\item The Fourier transform is $F(x)=(1-six)^{-m}$.
\end{enumerate}
\end{theorem}

\begin{proof}
This is a straightforward remake of Theorem 10.7, save of course for the key assertion there, which is the Gaussian vector one, the idea being as follows:

\medskip

(1) The moment computation is as follows, using the various properties of the gamma function, and with $M_0=1$ showing that we have indeed a probability measure:
\begin{eqnarray*}
M_k
&=&\frac{1}{s^m\Gamma(m)}\int_0^\infty x^{m+k-1}e^{-x/s}\,dx\\
&=&\frac{1}{s^m\Gamma(m)}\int_0^\infty (sy)^{m+k-1}e^{-y}\,sdy\\
&=&\frac{s^k}{\Gamma(m)}\int_0^\infty y^{m+k-1}e^{-y}\,dy\\
&=&s^k\cdot\frac{\Gamma(m+k)}{\Gamma(m)}\\
&=&s^km(m+1)\ldots(m+k-1)
\end{eqnarray*}

(2) According to our moment formula found above, the first few moments are:
$$M_1=sm$$
$$M_2=s^2m(m+1)$$
$$M_3=s^3m(m+1)(m+2)$$
$$M_4=s^4m(m+1)(m+2)(m+3)$$

Thus the mean is $E=sm$, and the variance is $V=s^2m(m+1)-(sm)^2=s^2m$. In order to compute $\gamma$ and $\kappa$, observe that the needed central moments are given by:
$$M_3'=M_3-3EM_2+2E^3=2s^3m$$
$$M_4'=M_4-4EM_3+6E^2M_2-3E^4=3s^4m(m+2)$$

Thus, the skewness and kurtotis are given by the following formulae:
$$\gamma=\frac{2s^3m}{s^3m\sqrt{m}}=\frac{2}{\sqrt{m}}\quad,\quad\kappa=\frac{3s^4m(m+2)}{s^4m^2}=3+\frac{6}{m}$$

(3) The Fourier transform computation is, as before, as follows:
\begin{eqnarray*}
F(x)
&=&\frac{1}{s^m\Gamma(m)}\int_0^\infty y^{m-1}e^{-y/s+ixy}\,dy\\
&=&\frac{1}{s^m\Gamma(m)}\int_0^\infty y^{m-1}e^{(six-1)y/s}\,dy\\
&=&\frac{1}{2^m\Gamma(m)}\int_0^\infty\left(\frac{z}{1-six}\right)^{m-1}e^{-z/s}\,\frac{dz}{1-six}\\
&=&(1-six)^{-m}\cdot\frac{1}{s^m\Gamma(m)}\int_0^\infty z^{m-1}e^{-z/s}\,dz\\
&=&(1-six)^{-m}
\end{eqnarray*}

Thus, we are led to the conclusions in the statement.
\end{proof}

As another extension now of Theorem 10.7, dealing this time with more general Gaussian vectors, assumed to be centered and with independent entries, we have:

\begin{theorem}
Given variables $f_1,\ldots,f_n$ which are centered normal, $f_i\sim g_{t_i}$, and independent, the corresponding chi squared law $\chi_{nT}^2=law(||f||^2)$ has moments
$$M_k=\frac{k!}{2^k}\sum_{a_1+\ldots+a_n=k}\binom{2a_1}{a_1}\ldots\binom{2a_n}{a_n}t_1^{a_1}\ldots t_n^{a_n}$$
the Fourier transform is given by the following formula,
$$F(x)=\frac{1}{\sqrt{(1-2t_1ix)\ldots(1-2t_nix)}}$$
the expectation is $E=\sum_it_i$, the variance is $V=2\sum_it_i^2$, and
$$\gamma=\frac{\sum_it_i^3}{\sum_it_i^2}\sqrt{\frac{8}{\sum_it_i^2}}\quad\ ,\ \quad
\kappa=3+\frac{12\sum_it_i^4}{(\sum_it_i^2)^2}$$
are the corresponding skewness and kurtosis.
\end{theorem}

\begin{proof}
This is a straightforward generalization of our previous results, as follows:

\medskip

(1) With $f_1,\ldots,f_n$ being as in the statement, $f_i\sim g_{t_i}$, independent, we have:
\begin{eqnarray*}
M_k(\chi_{nT}^2)
&=&E((f_1^2+\ldots+f_n^2)^k)\\
&=&\sum_{a_1+\ldots+a_n=k}\binom{k}{a_1,\ldots,a_n}E(f_1^{2a_1})\ldots E(f_n^{2a_n})\\
&=&\sum_{a_1+\ldots+a_n=k}\frac{k!}{a_1!\ldots a_n!}\,t_1^{a_1}(2a_1)!!\ldots t_n^{a_n}(2a_n)!!\\
&=&\sum_{a_1+\ldots+a_n=k}\frac{k!}{a_1!\ldots a_n!}\cdot\frac{(2a_1)!}{2^{a_1}a_1!}\cdots\frac{(2a_n)!}{2^{a_n}a_n!}\,t_1^{a_1}\ldots t_n^{a_n}\\
&=&\frac{k!}{2^k}\sum_{a_1+\ldots+a_n=k}\frac{(2a_1)!}{a_1!a_1!}\cdots\frac{(2a_n)!}{a_n!a_n!}\,t_1^{a_1}\ldots t_n^{a_n}\\
&=&\frac{k!}{2^k}\sum_{a_1+\ldots+a_n=k}\binom{2a_1}{a_1}\ldots\binom{2a_n}{a_n}\,t_1^{a_1}\ldots t_n^{a_n}
\end{eqnarray*}

(2) The generalized binomial formula with exponent $-1/2$ shows that we have:
$$\frac{1}{\sqrt{1-4t_iz}}=\sum_{a_i=0}^\infty\binom{2a_i}{a_i}(t_iz)^{a_i}$$

By multiplying these formulae, and then setting $k=a_1+\ldots+a_n$, we obtain:
\begin{eqnarray*}
\frac{1}{\sqrt{(1-4t_1z)\ldots(1-4t_nz)}}
&=&\sum_{a_1=0}^\infty\ldots\sum_{a_n=0}^\infty\binom{2a_1}{a_1}\ldots\binom{2a_n}{a_n}(t_1z)^{a_1}\ldots(t_nz)^{a_n}\\
&=&\sum_{k=0}^\infty z^k\sum_{k=a_1+\ldots+a_n}\binom{2a_1}{a_1}\ldots\binom{2a_n}{a_n}t_1^{a_1}\ldots t_n^{a_n}\\
&=&\sum_{k=0}^\infty z^k\cdot\frac{2^kM_k}{k!}\\
&=&\sum_{k=0}^\infty\frac{(2z)^k}{k!}\cdot M_k
\end{eqnarray*}

Now by setting $z=ix/2$, this latter formula takes the following form:
$$\frac{1}{\sqrt{(1-2t_1ix)\ldots(1-2t_nix)}}=\sum_{k=0}^\infty\frac{(ix)^k}{k!}\cdot M_k$$

But on the right we have the Fourier transform $F(x)$, so done with this too.

\medskip

(3) Next, according to the general formula in (1), the first moment is:
$$M_1
=\frac{1}{2}\sum_{a_1+\ldots+a_n=1}\binom{2a_1}{a_1}\ldots\binom{2a_n}{a_n}t_1^{a_1}\ldots t_n^{a_n}
=\frac{1}{2}\sum_i\binom{2}{1}t_i
=\sum_it_i$$

Thus, we have our mean. In order to compute now the variance, we have:
\begin{eqnarray*}
M_2
&=&\frac{1}{2}\sum_{a_1+\ldots+a_n=2}\binom{2a_1}{a_1}\ldots\binom{2a_n}{a_n}t_1^{a_1}\ldots t_n^{a_n}\\
&=&\frac{1}{2}\left[\sum_i\binom{4}{2}t_i^2+\sum_{i<j}\binom{2}{1}\binom{2}{1}t_it_j\right]\\
&=&3\sum_it_i^2+2\sum_{i<j}t_it_j\\
&=&2\sum_it_i^2+\left(\sum_it_i\right)^2
\end{eqnarray*}

We conclude from this that the variance is $V=2\sum_it_i^2$, as stated.

\medskip

(4) Getting now to the third moment, the computation becomes more complex:
\begin{eqnarray*}
M_3
&=&\frac{3}{4}\sum_{a_1+\ldots+a_n=3}\binom{2a_1}{a_1}\ldots\binom{2a_n}{a_n}t_1^{a_1}\ldots t_n^{a_n}\\
&=&\frac{3}{4}\left[\sum_i\binom{6}{3}t_i^3+\sum_{i\neq j}\binom{4}{2}\binom{2}{1}t_i^2t_j+\sum_{i<j<k}\binom{2}{1}\binom{2}{1}\binom{2}{1}t_it_jt_k\right]\\
&=&15\sum_it_i^3+9\sum_{i\neq j}t_i^2t_j+6\sum_{i<j<k}t_it_jt_k\\
&=&8\sum_it_i^3+6\left(\sum_it_i^2\right)\left(\sum_it_i\right)+\left(\sum_it_i\right)^3
\end{eqnarray*}

(5) In order to deal with this, it is convenient to introduce the sums $S_k=\sum_it_i^k$. Indeed, in terms of these sums, the various moment formulae that we found read:
$$M_1=S_1\quad,\quad M_2=2S_2+S_1^2\quad,\quad M_3=8S_3+6S_2S_1+S_1^3$$

We can now compute the third central moment, which follows to be:
$$M_3'=M_3-3EM_2+2E^3=8S_3$$

Thus, we are led to the skewness formula in the statement, namely:
$$\gamma=\frac{8S_3}{2S_2\sqrt{2S_2}}=\frac{S_3}{S_2}\sqrt{\frac{8}{S_2}}$$

(6) Getting now to the fourth moment, the computation is even more complex:
\begin{eqnarray*}
M_4
&=&\frac{3}{2}\sum_{a_1+\ldots+a_n=4}\binom{2a_1}{a_1}\ldots\binom{2a_n}{a_n}t_1^{a_1}\ldots t_n^{a_n}\\
&=&\frac{3}{2}\Big[\sum_i\binom{8}{4}t_i^4+\sum_{i\neq j}\binom{6}{3}\binom{2}{1}t_i^3t_j+\sum_{i<j}\binom{4}{2}\binom{4}{2}t_i^2t_j^2\\
&&+\sum_{j<k}\sum_{i\neq j,k}\binom{4}{2}\binom{2}{1}\binom{2}{1}t_i^2t_jt_k+\sum_{i<j<k<l}\binom{2}{1}\binom{2}{1}\binom{2}{1}\binom{2}{1}t_it_jt_kt_l\Big]\\
&=&105\sum_it_i^4+60\sum_{i\neq j}t_i^3t_j+54\sum_{i<j}t_i^2t_j^2+36\sum_{j<k}\sum_{i\neq j,k}t_i^2t_jt_k+24\sum_{i<j<k<l}t_it_jt_kt_l
\end{eqnarray*}

(7) Now by using tricks as in (3,4) we obtain, in terms of the sums $S_k=\sum_it_i^k$:
$$M_4=48S_4+32S_3S_1+12S_2^2+12S_2S_1^2+S_1^4$$

Thus, the fourth central moment is given by the following formula:
$$M_4'
=M_4-4EM_3+6E^2M_2-3E^4
=48S_4+12S_2^2$$ 

It follows that the kurtosis is given by the following formula:
$$\kappa=\frac{48S_4+12S_2^2}{(2S_2)^4}=3+\frac{12S_4}{S_2^2}$$

Thus, we are led to the various conclusions in the statement.
\end{proof}

\section*{10c. Gamma variables}

What is next? Many potential things to be done, including on one hand unifying Theorems 10.8 and 10.9, and on the other hand, throwing into the picture some shifts and covariances. Let us start with the unification question. We have here:

\index{gamma variables}

\begin{theorem}
Given variables $f_1,\ldots,f_r$ following gamma laws, $f_i\sim \mu_{m_is_i}$, and which are independent, $\mu_{MS}=law(f_1+\ldots+f_r)$ has the following properties:
\begin{enumerate}
\item When $s_i=s$ we have $\mu_{MS}=\mu_{ms}$, with $m=\sum_im_i$.

\item When $m_i=n/2$ with $n\in\mathbb N$ we have $\mu_{MS}=\chi_{nrT}^2$, with $T=S/2$.

\item The Fourier transform is $F(x)=\prod_i(1-s_iix)^{-m_i}$.

\item The moments are $M_k=k!\sum_{a_1+\ldots+a_r=k}\binom{m_1+a_1-1}{a_1}\ldots\binom{m_r+a_r-1}{a_r}s_1^{a_1}\ldots s_r^{a_r}$.

\item $E=S_1$, $V=S_2$, $\gamma=2S_3/S_2\sqrt{S_2}$, $\kappa=3+6S_4/S_2^2$, where $S_k=\sum_im_is_i^k$.
\end{enumerate} 
\end{theorem}

\begin{proof}
This is a straightforward unification of Theorems 10.8 and 10.9, with the idea of the proof, along with the details regarding the statement, being as follows:

\medskip

(1) In the case $s_i=s$ what we have are variables $f_i\sim\mu_{m_is}$, independent, and it follows that the Fourier transform of their sum $f_1+\ldots+f_r$ is given by:
$$F(x)=\prod_iF_i(x)=\prod_i(1-six)^{-m_i}=(1-six)^{-\sum_im_i}$$

But on the right we have the Fourier transform of $\mu_{ms}$ with $m=\sum_im_i$, as desired.

\medskip

(2) In the case $m_i=n/2$ with $n\in\mathbb N$, our variables as follows, with $t_i=s_i/2$:
$$f_i\sim\mu_{m_is_i}=\mu_{n/2,2t_i}=\chi_{nt_i}^2$$

Our variables being independent, the Fourier transform of $f_1+\ldots+f_r$ is:
$$F(x)=\prod_i(1-2ixt_i)^{-n/2}=\left[\prod_i(1-2ixt_i)^{-1/2}\right]^n=\left[F_{\chi_{rT}^2}(x)\right]^n$$

Now observe that the function on the right is the Fourier transform of the following law, with on the right the $T$ vector being repeated $n$ times:
$$\chi_{nrT}^2=\chi_{nr,T\ldots T}^2$$

We conclude that, with this convention, we have indeed $\mu_{MS}=\chi_{nrT}^2$, as stated.

\medskip

(3) The Fourier transform formula in the statement is clear, coming from:
$$F(x)=\prod_iF_i(x)=\prod_i(1-s_iix)^{-m_i}$$

(4) Getting now to the moments, the computation here is as follows, using the moment formula $M_k(\mu_{ms})=s^k\Gamma(m+k)/\Gamma(m)$ from Theorem 10.8, or rather from its proof, and with the outcome being in terms of generalized binomial coefficients $\binom{m_i+a_i-1}{a_i}$:
\begin{eqnarray*}
M_k
&=&E((f_1+\ldots+f_r)^k)\\
&=&\sum_{a_1+\ldots+a_r=k}\binom{k}{a_1,\ldots,a_r}E(f_1^{a_1})\ldots E(f_r^{a_r})\\
&=&\sum_{a_1+\ldots+a_r=k}\frac{k!}{a_1!\ldots a_r!}\,s_1^{a_1}\frac{\Gamma(m_1+a_1)}{\Gamma(m_1)}\ldots s_r^{a_r}\frac{\Gamma(m_r+a_r)}{\Gamma(m_r)}\\
&=&k!\sum_{a_1+\ldots+a_r=k}\frac{\Gamma(m_1+a_1)}{a_1!\Gamma(m_1)}\ldots \frac{\Gamma(m_r+a_r)}{a_r!\Gamma(m_r)}\,s_1^{a_1}\ldots s_r^{a_r}\\
&=&k!\sum_{a_1+\ldots+a_r=k}\binom{m_1+a_1-1}{a_1}\ldots\binom{m_r+a_r-1}{a_r}s_1^{a_1}\ldots s_r^{a_r}
\end{eqnarray*}

(5) Next, according to the above general moment formula, the first moment is:
$$M_1=\sum_i\binom{m_i}{1}s_i
=\sum_im_is_i$$

Thus, we have our mean. In order to compute now the variance, we have:
\begin{eqnarray*}
M_2
&=&2\left[\sum_i\binom{m_i+1}{2}s_i^2+\sum_{i<j}\binom{m_i}{1}\binom{m_j}{1}s_is_j\right]\\
&=&\sum_i(m_i+1)m_is_i^2+\sum_{i\neq j}m_im_js_is_j\\
&=&\left(\sum_im_is_i\right)^2+\sum_im_is_i^2
\end{eqnarray*}

We conclude from this that the variance is $V=\sum_im_is_i^2$, as stated.

\medskip

(6) Getting now to the third moment, the computation becomes more complex, with some standard symmetric function tricks leading to the following formula:
\begin{eqnarray*}
&&M_3\\
&=&6\left[\sum_i\binom{m_i+2}{3}s_i^3+\sum_{i\neq j}\binom{m_i+1}{2}\binom{m_j}{1}s_i^2s_j+\sum_{i<j<k}\binom{m_i}{1}\binom{m_j}{1}\binom{m_k}{1}s_is_js_k\right]\\
&=&\sum_i(m_i+2)(m_i+1)m_is_i^3+3\sum_{i\neq j}(m_i+1)m_im_js_i^2s_j+6\sum_{i<j<k}m_im_jm_ks_is_js_k\\
&=&\left(\sum_im_is_i\right)^3+3\left(\sum_im_is_i^2\right)\left(\sum_im_is_i\right)+2\sum_im_is_i^3
\end{eqnarray*}

(7) In order to deal with this, it is convenient to introduce the sums $S_k=\sum_im_is_i^k$. Indeed, in terms of these sums, the various moment formulae that we found read:
$$M_1=S_1\quad,\quad M_2=S_1^2+S_2\quad,\quad M_3=S_1^3+3S_2S_1+2S_3$$

We can now compute the third central moment, which follows to be:
$$M_3'=M_3-3EM_2+2E^3=2S_3$$

Thus, we are led to the skewness formula in the statement, namely:
$$\gamma=\frac{2S_3}{S_2\sqrt{S_2}}$$

(8) Getting now to the fourth moment, the formula is quite complex, as follows:
\begin{eqnarray*}
&&M_4\\
&=&24\Big[\sum_i\binom{m_i+3}{4}s_i^4+\sum_{i\neq j}\binom{m_i+2}{3}\binom{m_j}{1}s_i^3s_j+\sum_{i<j}\binom{m_i+1}{2}\binom{m_j+1}{2}s_i^2s_j^2\\
&+&\sum_{j<k}\sum_{i\neq j,k}\binom{m_i+1}{2}\binom{m_j}{1}\binom{m_k}{1}s_i^2s_js_k+\sum_{i<j<k<l}\binom{m_i}{1}\binom{m_j}{1}\binom{m_k}{1}\binom{m_l}{1}s_is_js_ks_l\Big]\\
&=&\sum_i(m_i+3)(m_i+2)(m_i+1)m_is_i^4+4\sum_{i\neq j}(m_i+2)(m_i+1)m_im_js_i^3s_j\\
&&+6\sum_{i<j}(m_i+1)m_i(m_j+1)m_js_i^2s_j^2+12\sum_{j<k}\sum_{i\neq j,k}(m_i+1)m_im_jm_ks_i^2s_js_k\\
&&+24\sum_{i<j<k}m_im_jm_km_ls_is_js_ks_l
\end{eqnarray*}

(9) However, by using tricks involving partitions and symmetric functions, as before in (5,6), we obtain, after some work, in terms of the sums $S_k=\sum_im_is_i^k$:
$$M_4=S_1^4+6S_1^2S_2+3S_2^2+8S_3S_1+6S_4$$

Thus, the fourth central moment is given by the following formula:
$$M_4'
=M_4-4EM_3+6E^2M_2-3E^4
=6S_4+3S_2^2$$ 

It follows that the kurtosis is given by the following formula:
$$\kappa=\frac{6S_4+3S_2^2}{S_2^2}=3+\frac{6S_4}{S_2^2}$$

Thus, we are led to the various conclusions in the statement.
\end{proof}

Summarizing, good work that we did, on the unification question, and for more on all this, sums of independent gamma variables, we refer to \cite{mat}, \cite{mos}. The continuation, however, involving shifts and covariances, is a far more complicated business, because we perfectly know from chapters 6 and 9 that the whole territory is a minefield. 

\bigskip

Indeed, getting back to the usual Gaussian setting, and to the law $\chi_{nT}^2$ studied in Theorem 10.9, we would like to generalize that into a law $\chi_\Sigma^{2a}$, involving a covariance matrix $\Sigma\in M_n(\mathbb R)$ and a position parameter $a\in\mathbb R^n$, the precise question being:

\begin{question}
What can we say about the law $\chi_\Sigma^{2a}$ given by
$$\int_0^\infty\varphi(x)d\chi_\Sigma^{2a}(x)
=\frac{1}{\sqrt{(2\pi)^n\det\Sigma}}\int_{\mathbb R^n}\varphi(||x||^2)\exp\left(-\frac{<\Sigma^{-1}(x-a),x-a>}{2}\right)dx$$
the parameters being a covariance matrix $\Sigma\in M_n(\mathbb R)$, and a position vector $a\in\mathbb R^n$?
\end{question}

In answer, we studied this question in chapter 6, at $n=1$. Our computation there was as follows, with the $1\times 1$ covariance matrix being denoted $\Sigma=(t)$, the position vector being $a\in\mathbb R$, and with the law itself $\chi_\Sigma^{2a}$ being denoted in this case $\chi_{1t}^{2a}$:
\begin{eqnarray*}
\int_0^\infty\varphi(x)d\chi_{1t}^{2a}(x)
&=&\frac{1}{\sqrt{2\pi t}}\int_\mathbb R\varphi(x^2)e^{-(x-a)^2/2t}dx\\
&=&\frac{1}{\sqrt{2\pi t}}\int_0^\infty\varphi(y)\left(e^{-(\sqrt{y}+a)^2/2t}+e^{-(\sqrt{y}-a)^2/2t}\right)\frac{dy}{2\sqrt{y}}\\
&=&\frac{e^{-a^2/2t}}{\sqrt{2\pi t}}\int_0^\infty\varphi(y)\cdot\frac{e^{-y/2t}}{\sqrt{y}}
\cdot\cosh\left(\frac{a\sqrt{y}}{t}\right)dy
\end{eqnarray*}

Thus, at $n=1$ our law is computable, given by the following formula:
$$\chi^{2a}_{1t}=\frac{e^{-a^2/2t}}{\sqrt{2\pi t}}\cdot\frac{e^{-x/2t}}{\sqrt{x}}\cosh\left(\frac{a\sqrt{x}}{t}\right)dx$$

Moreover, as also explained in chapter 6, the moments of this law are something quite familiar, as follows, allowing us to easily compute the parameters $E,V,\gamma,\kappa$:
$$M_k(\chi^{2a}_{1t})=M_{2k}(g_t^a)$$

In higher dimensions such things will not work, due to both the covariance matrix and position vector, the problem coming from $<\Sigma^{-1}(x-a),x-a>$, which will not simplify, when plugging in $x$ written in spherical coordinates, unless we are in the case $a=0$ and $\Sigma$ diagonal, which corresponds to the original situation, that from Theorem 10.9.

\bigskip

In short, nice try, but we ended up in a loop. As for the more general case of gamma variables, as in Theorem 10.10, it is not even clear how to talk about shifts and covariances, in that setting. We will leave some thinking here as an interesting exercise.

\section*{10d. Maxwell-Boltzmann}

Not much time left in this chapter, and we still have to discuss the second family of laws introduced in Definition 10.1, namely the chi ones. These are a bit more complicated than the chi squared ones, with their mathematics appearing by extracting the ``square root'' of what we did. Fortunately, this can be done, the essentials being as follows:

\begin{theorem}
The chi law of parameter $n\in\mathbb N$, appearing as
$$\chi_n=law\left(\sqrt{f_1^2+\ldots+f_n^2}\right)$$
with $f_1,\ldots,f_n$ being independent, each following $g_1$, has density and moments given by
$$\chi_n=\frac{x^{n-1}e^{-x^2/2}}{2^{n/2-1}\Gamma(n/2)}\,dx\quad,\quad 
M_k=2^{k/2}\frac{\Gamma\left(\frac{n+k}{2}\right)}{\Gamma(n/2)}$$
with $\Gamma$ being as usual the gamma function.
\end{theorem}

\begin{proof}
This is something very standard, the idea being as follows:

\medskip

(1) In order to compute the density, we can use Theorem 10.2, which gives:
\begin{eqnarray*}
E(\varphi(||f||))
&=&\frac{1}{2^{n/2}\Gamma(n/2)}\int_0^\infty\varphi(\sqrt{x})x^{n/2-1}e^{-x/2}\,dx\\
&=&\frac{1}{2^{n/2}\Gamma(n/2)}\int_0^\infty\varphi(y)y^{n-2}e^{-y^2/2}\,2ydy\\
&=&\frac{1}{2^{n/2-1}\Gamma(n/2)}\int_0^\infty\varphi(y)y^{n-1}e^{-y^2/2}\,dy
\end{eqnarray*}

We conclude from this that the density of $\chi_n$ is the one in the statement.

\medskip

(2) As for the moment computation, this is straightforward too, as follows:
\begin{eqnarray*}
M_k
&=&\frac{1}{2^{n/2-1}\Gamma(n/2)}\int_0^\infty x^{n+k-1}e^{-x^2/2}\,dx\\
&=&\frac{1}{2^{n/2-1}\Gamma(n/2)}\int_0^\infty\left(\sqrt{2y}\right)^{n+k-1}e^{-y}\,\frac{dy}{\sqrt{2y}}\\
&=&\frac{2^{(n+k-1)/2}}{2^{n/2-1/2}\Gamma(n/2)}\int_0^\infty y^{(n+k)/2-1}e^{-y}\,dy\\
&=&\frac{2^{k/2}}{\Gamma(n/2)}\,\Gamma\left(\frac{n+k}{2}\right)
\end{eqnarray*}

(3) Finally, in what regards the moments, whose formula is something quite tricky, let us comment a bit on this. In what regards the even moments, these are given by:
$$M_2=n$$
$$M_4=n(n+2)$$
$$M_6=n(n+2)(n+4)$$
$$M_8=n(n+2)(n+4)(n+6)$$
$$\vdots$$

(4) In what regards the odd moments, however, things are more complicated. According to the general moment formula that we found, at $k=1$, the mean is:
$$E=\sqrt{2}\cdot\frac{\Gamma\left(\frac{n+1}{2}\right)}{\Gamma(n/2)}$$

Now recall from before that the gamma function is given at half-integers by the following almost uniform formula, with $c_N=\sqrt{2},\sqrt{\pi}$ for $N$ even, odd:
$$\Gamma\left(\frac{N}{2}\right)=\frac{(N-1)!!}{2^{(N-1)/2}}\,c_N$$

By using this, the above formula of the mean becomes more explicit, as follows:
$$E
=\sqrt{2}\cdot\frac{n!!}{2^{n/2}}\,c_{n+1}\cdot\frac{2^{(n-1)/2}}{(n-1)!!}\cdot\frac{1}{c_n}
=\frac{n!!}{(n-1)!!}\cdot\frac{c_{n+1}}{c_n}$$

Finally, observe that we have the following formula, for the term on the right:
$$\frac{c_{n+1}}{c_n}=\begin{cases}
\sqrt{\frac{\pi}{2}}&(n\ {\rm even})\\
\sqrt{\frac{2}{\pi}}&(n\ {\rm odd})
\end{cases}$$

(5) As for the higher odd moments, no need for new computations, because the formula $\Gamma(s+1)=s\Gamma(s)$ does the job, and shows that the sequence of odd moments is:
$$M_1=E$$
$$M_3=(n+1)E$$
$$M_5=(n+1)(n+3)E$$
$$M_7=(n+1)(n+3)(n+5)E$$
$$\vdots$$

(6) So, this is the situation, and while we can certainly compute $(E,V,\gamma,\kappa)$, all formulae will involve the above alien quantity $E$. And more on this, a bit later.
\end{proof}

Next, let us add a scaling parameter $t>0$ to what we have. We are led to:

\begin{theorem}
The chi law of parameters $n\in\mathbb N$ and $t>0$, appearing as
$$\chi_{nt}=law\left(\sqrt{f_1^2+\ldots+f_n^2}\right)$$
with $f_1,\ldots,f_n$ being independent, each following $g_t$, has density and moments given by
$$\chi_{nt}=\frac{2x^{n-1}e^{-x^2/2t}}{(2t)^{n/2}\Gamma(n/2)}\,dx
\quad,\quad 
M_k=(2t)^{k/2}\frac{\Gamma\left(\frac{n+k}{2}\right)}{\Gamma(n/2)}$$
with $\Gamma$ being as usual the gamma function.
\end{theorem}

\begin{proof}
This is a straightforward parametric remake of what we have:

\medskip

(1) In order to compute the density, we can use use Theorem 10.7, which gives:
\begin{eqnarray*}
E(\varphi(||f||))
&=&\frac{1}{(2t)^{n/2}\Gamma(n/2)}\int_0^\infty\varphi(\sqrt{x})x^{n/2-1}e^{-x/2t}\,dx\\
&=&\frac{1}{(2t)^{n/2}\Gamma(n/2)}\int_0^\infty\varphi(y)y^{n-2}e^{-y^2/2t}\,2ydy\\
&=&\frac{2}{(2t)^{n/2}\Gamma(n/2)}\int_0^\infty\varphi(y)y^{n-1}e^{-y^2/2t}\,dy
\end{eqnarray*}

We conclude from this that the density of $\chi_n$ is the one in the statement.

\medskip

(2) As for the moment computation, since the moments of the normal laws obey to $M_k(g_t)=t^{k/2}M_k(g_1)$, we conclude from this that we have:
$$M_k(\chi_{nt})=t^{k/2}M_k(\chi_n)=(2t)^{k/2}\frac{\Gamma\left(\frac{n+k}{2}\right)}{\Gamma(n/2)}$$

(3) Finally, in what regards $(E,V,\gamma,\kappa)$, the same comments as before apply, the point being that all formulae will look quite bad, involving $E$. More on this later.
\end{proof}

As a further generalization of what we have, again straightforward, we have:

\index{square root of gamma}

\begin{theorem}
The square root of the gamma law of parameters $m,s>0$, given by
$$f\sim\mu_{ms}\implies\sqrt{f}\sim\mu_{ms}^{1/2}$$
has density and moments given by the following formulae,
$$\chi_{ms}^{1/2}=\frac{2x^{2m-1}e^{-x^2/s}}{s^m\Gamma(m)}\,dx
\quad,\quad 
M_k=s^{k/2}\,\frac{\Gamma(m+k/2)}{\Gamma(m)}$$
with $\Gamma$ being as usual the gamma function.
\end{theorem}

\begin{proof}
This is again a straightforward remake of what we have:

\medskip

(1) In order to compute the density, we can use Theorem 10.8, which gives:
\begin{eqnarray*}
E(\varphi(\sqrt{f}))
&=&\frac{1}{s^m\Gamma(m)}\int_0^\infty\varphi(\sqrt{x})x^{m-1}e^{-x/s}\,dx\\
&=&\frac{1}{s^m\Gamma(m)}\int_0^\infty\varphi(y)y^{2m-2}e^{-y^2/s}\,2ydy\\
&=&\frac{2}{s^m\Gamma(m)}\int_0^\infty\varphi(y)y^{2m-1}e^{-y^2/s}\,dy
\end{eqnarray*}

We conclude from this that the density of $\mu_{ms}^{1/2}$ is the one in the statement.

\medskip

(2) As for the moment computation, this is straigtforward too, as follows:
\begin{eqnarray*}
M_k
&=&\frac{2}{s^m\Gamma(m)}\int_0^\infty x^{2m+k-1}e^{-x^2/s}\,dx\\
&=&\frac{2}{s^m\Gamma(m)}\int_0^\infty\left(\sqrt{sy}\right)^{2m+k-1}e^{-y}\,\frac{\sqrt{s}\,dy}{2\sqrt{y}}\\
&=&\frac{s^{k/2}}{\Gamma(m)}\int_0^\infty y^{m+k/2-1}e^{-y}\,dy\\
&=&\frac{s^{k/2}}{\Gamma(m)}\,\Gamma\left(m+\frac{k}{2}\right)
\end{eqnarray*}

(3) Finally, in what regards $(E,V,\gamma,\kappa)$, the same comments as before apply, the point being that all formulae will look quite bad, involving $E$. More on this later.
\end{proof}

As a last topic of discussion, in our familiar $n=3$ dimensions, we have:

\index{Maxwell-Boltzmann}

\begin{theorem}
The Maxwell-Boltzmann law $\chi_{3t}$ and its moments are given by
$$\chi_{3t}=\sqrt{\frac{2}{t^3\pi}}\,x^2e^{-x^2/2t}dx\quad,\quad 
M_k=2\sqrt{\frac{(2t)^k}{\pi}}\,\Gamma\left(\frac{k+3}{2}\right)$$
with $\Gamma$ being as usual the gamma function, and 
$$E=\sqrt{\frac{8t}{\pi}}\ ,\ 
V=\left(3-\frac{8}{\pi}\right)t\ ,\ 
\gamma=\frac{32-10\pi}{3\pi-8}\sqrt{\frac{2}{3\pi-8}}\ ,\ 
\kappa=\frac{15\pi^2+16\pi-192}{(3\pi-8)^2}$$
are the corresponding mean, variance, skewness and kurtosis.
\end{theorem}

\begin{proof}
We use the general formulae that we found in Theorem 10.13, namely:
$$\chi_{nt}=\frac{2x^{n-1}e^{-x^2/2t}}{(2t)^{n/2}\Gamma(n/2)}\,dx
\quad,\quad 
M_k=(2t)^{k/2}\frac{\Gamma\left(\frac{n+k}{2}\right)}{\Gamma(n/2)}$$

We can see that at $n=3$ the density is as follows, using $\Gamma(3/2)=\sqrt{\pi}/2$:
$$\chi_{3t}=\frac{2x^2e^{-x^2/2t}}{(2t)^{3/2}\sqrt{\pi}/2}\,dx
=\sqrt{\frac{2}{t^3\pi}}\,x^2e^{-x^2/2t}dx$$

As for the moment formula, this is as follows, again by using $\Gamma(3/2)=\sqrt{\pi}/2$:
$$M_k=(2t)^{k/2}\frac{\Gamma\left(\frac{k+3}{2}\right)}{\sqrt{\pi}/2}
=2\sqrt{\frac{(2t)^k}{\pi}}\,\Gamma\left(\frac{k+3}{2}\right)$$

Getting now to the low order moments, these are given by the following formulae:
$$M_1=2\sqrt{\frac{2t}{\pi}}\,\Gamma(2)=\sqrt{\frac{8t}{\pi}}$$
$$M_2=2\sqrt{\frac{(2t)^2}{\pi}}\,\Gamma\left(\frac{5}{2}\right)=3t$$
$$M_3=2\sqrt{\frac{(2t)^3}{\pi}}\,\Gamma(3)=8t\sqrt{\frac{2t}{\pi}}$$
$$M_4=2\sqrt{\frac{(2t)^4}{\pi}}\,\Gamma\left(\frac{7}{2}\right)=15t^2$$

Thus $E,V$ are given by the formulae in the statement. Next, we have:
$$M_3'=M_3-3EM_2+2E^3=\left(\frac{32}{\pi}-10\right)t\sqrt{\frac{2t}{\pi}}$$
$$M_4'=M_4-4EM_3+6E^2M_2-3E^4=15t^2+16\frac{t^2}{\pi}-192\frac{t^2}{\pi^2}$$

Thus, we are led to the formulae of $\gamma$ and $\kappa$ in the statement.
\end{proof}

At the level of the phenomenology, the Maxwell-Boltzmann law appears as follows:

\index{ideal gas}
\index{thermal equilibrium}
\index{Boltzmann constant}
\index{Maxwell-Boltzmann law}
\index{molecular speeds}
\index{kinetic energy}
\index{pressure}
\index{temperature}

\begin{theorem}
For an ideal gas in thermal equilibrium at temperature $T$, the molecular speeds $v\in\mathbb R^3$ are subject to the Maxwell-Boltzmann distribution formula
$$P(v)=\left(\frac{m}{2\pi bT}\right)^{3/2}\exp\left(-\frac{m||v||^2}{2bT}\right)$$
that is, are Gaussian of parameter $t=bT/m$, with $m$ being the mass of the molecules, and $b=1.380\,649\times 10^{-23}$ being the Boltzmann constant. Thus $||v||\sim\chi_{3t}$, with $t=bT/m$.
\end{theorem}

\begin{proof}
Welcome to thermodynamics, the idea being as follows:

\medskip

(1) First comes pressure. For a gas having point molecules, with no collisions between them, the pressure $P$ is easy to compute, by doing the manometer math, given by the following formula, $V$ being the volume of the gas, and $K$ the total kinetic energy:
$$PV=\frac{2K}{3}$$

(2) Next comes temperature. This is something more tricky, and we will take as definition for temperature $T$ the following formula, relating it to the pressure $P$ and volume $V$, where $k=Nb$, with $N$ being the number of molecules present:
$$PV=kT$$

(3) In relation to what we want to do, what we need to know is the following formula for the average molecular energy $K_0=K/N$, coming by combining (1) and (2):
$$\frac{2NK_0}{3}=Nbt\ \implies\ K_0=\frac{3bT}{2}$$

(4) Getting to work now, following Maxwell, we are looking for the probability distribution of the molecular speeds $v\in\mathbb R^3$. Intuition tells us that there are no correlations between the $x,y,z$ directions of space, so we must have a formula as follows:
$$P(v)=f(v_1)f(v_2)f(v_3)$$

(5) On the other hand, by rotational symmetry, $P(v)$ must depend only on the magnitude $||v||$ of the velocity, and not on the direction. Thus, we must have as well:
$$P(v)=\varphi(||v||^2)$$

(6) Now by comparing the requirements in (4) and (5), we are led via some math to the conclusion that $\varphi$ must be an exponential, which amounts in saying that:
$$P(v)=\lambda\exp\left(-C||v||^2\right)$$

(7) Obviously we must have $C>0$, for things to be bounded, and then by integrating we can obtain $\lambda$ as function of $C$, and our formula becomes:
$$P(v)=\left(\frac{C}{\pi}\right)^{3/2}\exp\left(-C||v||^2\right)$$

(8) It remains to find $C$. But for this purpose, observe that, now that we have our distribution, we can compute everything. In particular, we find that on average:
$$v_1^2=v_2^2=v_3^2=\frac{1}{2C}$$

Thus on average $||v||=3/(2C)$, so the average kinetic energy of the molecules is:
$$K_0=\frac{m||v||^2}{2}=\frac{3m}{4C}$$

Now by comparing with (3) we obtain from this $C=m/(2bT)$, as desired.

\medskip

(9) So this was for the story of the formula, following Maxwell, and later Boltzmann came with a fully rigorous proof. And exercise for you, to learn more, about all this.
\end{proof}

\section*{10e. Exercises}

Welcome to statistics I guess, and as exercises on all this, we have:

\begin{exercise}
Study more advanced aspects of $\chi_n^2$.
\end{exercise}

\begin{exercise}
Study more advanced aspects of $\chi_n^{2t}$.
\end{exercise}

\begin{exercise}
Study more advanced aspects of $\mu_{ms}$.
\end{exercise}

\begin{exercise}
Further study the basic aspects of $\chi_{nT}^2$.
\end{exercise}

\begin{exercise}
Further study the basic aspects of $\mu_{MS}$.
\end{exercise}

\begin{exercise}
Can you say something nice about $\chi_\Sigma^{2a}$?
\end{exercise}

\begin{exercise}
Can you talk about something of type $\mu_\Sigma^a$?
\end{exercise}

\begin{exercise}
Study the square roots of the above laws.
\end{exercise}

As bonus exercise, buy a computer, and keep it offline, for computations.

\chapter{Hyperspherical laws}

\section*{11a. Circles, ellipses}

We have seen a lot of multivariable analysis already, in the present Part III, in relation with the normal laws. Our aim in this chapter will be to have a more systematic geometric discussion, regarding the spaces $\mathbb R^N$ and their unit spheres $S^{N-1}_\mathbb R\subset\mathbb R^N$. Which will eventually lead us, and no surprise here, into probability and normal laws.

\bigskip

Let us start with a general mathematical principle, which is something notorious:

\begin{principle}
The interplay between the ambient space and its unit sphere
$$S^{N-1}_\mathbb R\subset\mathbb R^N$$
is the alpha and omega of everything, in $N$-dimensional mathematics.
\end{principle}

To be more precise, this is something that we are very familiar with, and this since high school, in 2 dimensions, where the unit circle $C\subset\mathbb R^2$ produces trigonometry, polar coordinates, complex numbers, and many other things. By the way, speaking high school, we are quite familar as well with 3 dimensions, from geography, where the unit sphere $S\subset\mathbb R^3$, called there Earth, produces latitude, longitude, parallels, meridians, time zones, all sorts of hemispheres, equator, tropics and polar circles, and many more.

\bigskip

As for the more advanced level, we both know from physics classes about the ubiquitous spherical coordinates, I mean any question has a ``radiating center'', inviting for spherical coordinates there. We also know about the stereographic projection, and about the basic theory of curves and surfaces, heavily inspired from that of the spheres. In fact, advanced mathematics and physics quite often mean analysis on manifolds, and with the basic example of manifold, source of inspiration for everything, being $S^{N-1}_\mathbb R\subset\mathbb R^N$.

\bigskip 

So, can we have some probability saying, on all this? Here is our question:

\index{sphere coordinates}

\begin{question}
We would like to compute the joint law of the sphere coordinates
$$x_1:S^{N-1}_\mathbb R\to[-1,1]$$
$$\vdots$$
$$x_N:S^{N-1}_\mathbb R\to[-1,1]$$
and then, have some applications of this probabilistic technology.
\end{question}

Getting started now, let us first investigate the $N=2$ case. Here our sphere is the unit circle $C\subset\mathbb R^2$, the circle coordinates are $x=\cos t,y=\sin t$, and the first question is, how to integrate over $C$ the arbitrary products of these quantities $\cos t,\sin t$. 

\bigskip

In answer, we first have the following formula, for the arbitrary powers of $\cos t,\sin t$, that we already met in chapter 10, and that we reproduce here for convenience:

\index{double factorials}
\index{trigonometric integral}

\begin{theorem}[Wallis]
We have the following formulae,
$$\int_0^{\pi/2}\cos^pt\,dt=\int_0^{\pi/2}\sin^pt\,dt=\left(\frac{\pi}{2}\right)^{\varepsilon(p)}\frac{p!!}{(p+1)!!}$$
where $\varepsilon(p)=1$ if $p$ is even, and $\varepsilon(p)=0$ if $p$ is odd.
\end{theorem}

\begin{proof}
Let us first compute the integral on the left $I_p$. We have:
\begin{eqnarray*}
(\cos^pt\sin t)'
&=&p\cos^{p-1}t(-\sin t)\sin t+\cos^pt\cos t\\
&=&p\cos^{p+1}t-p\cos^{p-1}t+\cos^{p+1}t\\
&=&(p+1)\cos^{p+1}t-p\cos^{p-1}t
\end{eqnarray*}

By integrating now between $0$ and $\pi/2$, we obtain the following formula:
$$(p+1)I_{p+1}=pI_{p-1}$$

Thus we can compute $I_p$ by recurrence, and we obtain in this way:
\begin{eqnarray*}
I_p
&=&\frac{p-1}{p}\,I_{p-2}\\
&=&\frac{p-1}{p}\cdot\frac{p-3}{p-2}\,I_{p-4}\\
&\vdots&\\
&=&\frac{p!!}{(p+1)!!}\,I_{1-\varepsilon(p)}
\end{eqnarray*}

But $I_0=\frac{\pi}{2}$ and $I_1=1$, so we get the result. As for the second formula, this follows from the first one, with $t=\frac{\pi}{2}-s$. Thus, we proved both formulae in the statement.
\end{proof}

More generally now, we have the following result, solving our integration problem:

\index{trigonometric integral}

\begin{theorem}[Wallis 2]
We have the following formula,
$$\int_0^{\pi/2}\cos^pt\sin^qt\,dt=\left(\frac{\pi}{2}\right)^{\varepsilon(p)\varepsilon(q)}\frac{p!!q!!}{(p+q+1)!!}$$
where $\varepsilon(p)=1$ if $p$ is even, and $\varepsilon(p)=0$ if $p$ is odd.
\end{theorem}

\begin{proof}
Let $I_{pq}$ be the integral in the statement. We have the following formula:
\begin{eqnarray*}
(\cos^pt\sin^qt)'
&=&p\cos^{p-1}t(-\sin t)\sin^qt\\
&+&\cos^pt\cdot q\sin^{q-1}t\cos t\\
&=&-p\cos^{p-1}t\sin^{q+1}t+q\cos^{p+1}t\sin^{q-1}t
\end{eqnarray*}

By integrating this between $0$ and $\pi/2$, we obtain, for any $p,q>0$:
$$pI_{p-1,q+1}=qI_{p+1,q-1}$$

Thus, we can compute $I_{pq}$ by recurrence. When $q$ is even we obtain in this way:
\begin{eqnarray*}
I_{pq}
&=&\frac{q-1}{p+1}\,I_{p+2,q-2}\\
&=&\frac{q-1}{p+1}\cdot\frac{q-3}{p+3}\,I_{p+4,q-4}\\
&\vdots&\\
&=&\frac{p!!q!!}{(p+q)!!}\,I_{p+q}
\end{eqnarray*}

But the last term comes from Theorem 11.3, and we obtain the result, for $q$ even:
$$I_{pq}
=\frac{p!!q!!}{(p+q)!!}\left(\frac{\pi}{2}\right)^{\varepsilon(p+q)}\frac{(p+q)!!}{(p+q+1)!!}
=\left(\frac{\pi}{2}\right)^{\varepsilon(p)\varepsilon(q)}\frac{p!!q!!}{(p+q+1)!!}$$

Observe that this gives the result for $p$ even as well, by symmetry. In the remaining case, where both $p,q$ are odd, we can use again $pI_{p-1,q+1}=qI_{p+1,q-1}$, which gives:
\begin{eqnarray*}
I_{pq}
&=&\frac{q-1}{p+1}\,I_{p+2,q-2}\\
&=&\frac{q-1}{p+1}\cdot\frac{q-3}{p+3}\,I_{p+4,q-4}\\
&\vdots&\\
&=&\frac{p!!q!!}{(p+q-1)!!}\,I_{p+q-1,1}
\end{eqnarray*}

But the last term is easy to compute, by using the primitive, as follows:
$$I_{p1}
=\int_0^{\pi/2}\cos^pt\sin t\,dt
=-\frac{1}{p+1}\int_0^{\pi/2}(\cos^{p+1}t)'\,dt
=\frac{1}{p+1}$$

We can therefore finish our computation in the case $p,q$ odd, and we obtain:
$$I_{pq}
=\frac{p!!q!!}{(p+q-1)!!}\cdot\frac{1}{p+q}
=\frac{p!!q!!}{(p+q+1)!!}$$

Thus, we are led to the formula in the statement, the exponent of $\pi/2$ appearing there being $\varepsilon(p)\varepsilon(q)=0\cdot 0=0$ in the present case, and this finishes the proof.
\end{proof}

Good news, we can answer now Question 11.2 at $N=2$, as follows:

\index{arcsine law}

\begin{theorem}
With the notation $z=(x,y)$ for the points $z\in C$, and with respect to the uniform measure on $C$, both variables $x,y:C\to[-1,1]$ follow the law
$$\alpha_{1/4}=\frac{1}{\pi\sqrt{1-x^2}}\,dx$$
that is, the arcsine law on $[-1,1]$. The even the moments of this law are as follows,
$$M_{2r}=\frac{1}{4^r}\binom{2r}{r}$$
and whose odd moments vanish. The expectation, variance, skewness are kurtosis are:
$$E=0\quad,\quad V=\frac{1}{2}\quad,\quad\gamma=0\quad,\quad\kappa=\frac{3}{2}$$
Also, the joint moments of $(x,y)$ are given by the following formula,
$$E(x^{2r}y^{2s})=\frac{1}{4^{r+s}}\cdot\frac{(2r)!(2s)!}{r!s!(r+s)!}$$
with the other joint moments, containing an odd exponent, vanishing.
\end{theorem}

\begin{proof}
We certainly have variables $x,y:C\to[-1,1]$ as above, whose individual and joint behavior, meaning law and moments, we can investigate. But this can be done by using the Wallis formula from Theorem 11.4, the details being as follows:

\medskip

(1) With the convention in the statement, the Wallis formula from Theorem 11.4 takes the following form, in the non-vanishing case, where both exponents are even:
\begin{eqnarray*}
E(x^{2r}y^{2s})
&=&\frac{(2r)!!(2s)!!}{(2r+2s+1)!!}\\
&=&\frac{(2r)!/(2^rr!)\cdot (2s)!/(2^ss!)}{2^{r+s}(r+s)!}\\
&=&\frac{1}{4^{r+s}}\cdot\frac{(2r)!(2s)!}{r!s!(r+s)!}
\end{eqnarray*}

Thus, we are led to the various moment assertions in the statement. 

\medskip

(2) Regarding the formula of the density, this comes either by thinking, or from our previous knowledge of the arcsine law, from chapter 3, or simply as follows:
$$\frac{1}{\pi}\int_{-1}^1\frac{x^{2r}}{\sqrt{1-x^2}}\,dx=\frac{2}{\pi}\int_0^{\pi/2}\frac{\sin^{2r}t}{\cos t}\,\cos tdt=\frac{2}{\pi}\int_0^{\pi/2}\sin^{2r}t\,dt$$

Thus, the moments match indeed, and we have our density, as desired.

\medskip

(3) Finally, the formulae of $E,V,\gamma,\kappa$ are something that we know, from chapter 3.
\end{proof}

Let us record as well a version of Theorem 11.5, focusing on the variables $x^2,y^2$:

\begin{theorem}[version]
With the notation $z=(x,y)$ for the points $z\in C$, and with respect to the uniform measure, both variables $x^2,y^2:C\to[0,1]$ follow the law
$$\alpha_{1/4}'=\frac{1}{\pi\sqrt{x-x^2}}\,dx$$
that is, the arcsine law on $[0,1]$. The moments of this law are as follows,
$$M_r=\frac{1}{4^r}\binom{2r}{r}$$
and the corresponding expectation, variance, skewness are kurtosis are given by:
$$E=\frac{1}{2}\quad,\quad V=\frac{1}{8}\quad,\quad\gamma=0\quad,\quad\kappa=\frac{3}{2}$$
Also, the joint moments of $(x^2,y^2)$ are given by the following formula,
$$E(x^{2r}y^{2s})=\frac{1}{4^{r+s}}\cdot\frac{(2r)!(2s)!}{r!s!(r+s)!}$$
the covariance is $cov(x^2,y^2)=-1/8$, and the correlation is $\rho_{x^2y^2}=-1$.
\end{theorem}

\begin{proof}
This is a straightforward remake of Theorem 11.5, as follows:

\medskip

(1) To start with, the joint moment formula is the one from Theorem 11.5, left unchanged, and in particular we have the individual moments, which give the above density. Alternatively, we can use the fact, that we know well from chapter 3, that when a variable $f$ follows the arcsine law $\alpha_t$, its square $f^2$ follows the rival arcsine law $\alpha_t'$. 

\medskip

(2) Next, according to our general moment formula, at low order we have:
$$M_1=\frac{1}{2}\quad,\quad M_2=\frac{3}{8}\quad,\quad M_3=\frac{5}{16}\quad,\quad M_4=\frac{35}{128}$$

Thus the expectation is $E=1/2$, the variance is $V=1/8$, and then, we have:
$$M_3'=M_3-3EM_2+2E^3=0$$
$$M_4'=M_4-4EM_3+6E^2M_2-3E^4=\frac{3}{128}$$

We are therefore led to the formulae of $\gamma$ and $\kappa$ in the statement.

\medskip

(3) Finally, the assertions at the end are trivial, but for the fun, let me report my own experience with this. So, when writing this book, Theorem 11.5 typed, good, then the present Theorem typed too, also routine, and then I say to myself, eventually, time for some interesting stuff, computing some covariances. So, here is what I did:
$$cov(x^2,y^2)=E(x^2y^2)-E(x^2)E(y^2)=\frac{1}{16}\cdot\frac{2!2!}{1!1!2!}-\frac{1}{4}=\frac{1}{8}-\frac{1}{4}=-\frac{1}{8}$$

Which sounds a bit bizarre, because this implies $\rho_{x^2y^2}=-1$, so equality in the Cauchy-Schwarz inequality leading to $\rho\in[-1,1]$, what is going on here. And after 1 hour of thinking and paranoid doublechecks, of pretty much everything, I realized that there is no contradiction, because $x^2=\cos^2t$ and $y^2=\sin^2t$ satisfy $x^2+y^2=1$, and so $\rho_{x^2y^2}=-1$, the Cauchy-Schwarz vectors being $x^2-1/2=-(y^2-1/2)$. Simple like that.
\end{proof}

Moving on, now that we understood the circle, probabilistically speaking, before getting into higher dimensions, let us have as well a look at ellipses. Regarding these, in the hope that you are familiar with them, and here is what is to be known about them:

\index{ellipse}
\index{conic}

\begin{theorem}
We can talk about ellipses, in several possible ways, as follows:
\begin{enumerate}
\item As the curves appearing via $d(z,p)+d(z,q)=l$, with $p,q\in\mathbb R^2$, and $l>0$.

\item Up to a translation and rotation, as the curves given by $(x/a)^2+(y/b)^2=1$.

\item Or equivalently, in polar coordinates, given by $x=a\cos t$, $y=b\sin t$.

\item As the conics, that is, the plane curves of degree $2$, which are compact. 

\item As the compact curves obtained by cutting a 3D cone with a plane.

\item As the pictures of a circle, taken from various 3D perspectives.

\item As the trajectories of the planets around the Sun.
\end{enumerate}
\end{theorem}

\begin{proof}
All this is first-class mathematics, the idea being as follows:

\medskip

(1) This is something practical, allowing you to draw ellipses, armed with a string.

\medskip

(2) The equivalence with (1) comes from some computations, good exercise for you.

\medskip

(3) This is indeed something which is clearly equivalent to (2).

\medskip

(4) This is indeed equivalent to (2), via some standard algebraic computations.

\medskip

(5) This is indeed equivalent to (4), via some standard 3D computations.

\medskip

(6) This is the same thing as (5), for the attention of Humanities students.

\medskip

(7) This is tough, due to Kepler and Newton, check any of my calculus books.
\end{proof}

Getting now to integration matters, over the ellipses, things are quite interesting, because we have at least 3 different versions of Question 11.2, as follows:

\begin{question}
Given an ellipse $F$, appearing via $(x/a)^2+(y/b)^2=1$, what is the joint law of the coordinates $x:F\to[-a,a]$ and $y:F\to[-b,b]$, with respect to:
\begin{enumerate}
\item The measure on $F$ obtained by transporting the uniform measure on $C$.

\item The measure on $F$ obtained by walking on $F$, at constant speed.

\item The measure on $F$ obtained by traveling on $F$, planet style.
\end{enumerate}
\end{question}

To be more precise here, the measure in (1) is the ``dumb'' one, push-forward of the uniform measure on $C$, via the map $(x,y)\to(ax,by)$. The measure in (2), which is more subtle, is the usual path integral one, and more on this in a moment. As for the measure in (3), this is yet another beast, with the uniform parameter here being the time $t$, when traveling on $F$ gravitationally, according to the laws of Kepler and Newton.

\bigskip

In answer now, in what regards the first question, we have here:

\begin{theorem}
For an ellipse $F$ coming via $(x/a)^2+(y/b)^2=1$, and with respect to the measure transported from $C$ in the obvious way, that is, with 
$$x=a\cos t\quad,\quad y=b\sin t$$
with $t$ being assumed to be uniform, the joint moments of the coordinates are
$$E(x^{2r}y^{2s})=\frac{a^{2r}b^{2s}}{4^{r+s}}\cdot\frac{(2r)!(2s)!}{r!s!(r+s)!}$$
and these coordinates follow the arcsine laws of $[-a,a]$ and $[-b,b]$, respectively.
\end{theorem}

\begin{proof}
This is something self-explanatory, with the main formula coming from:
$$E(x^{2r}y^{2s})=\int_0^{2\pi}(a\cos t)^{2r}(b\sin t)^{2s}dt=a^{2r}b^{2s}\cdot\frac{1}{4^{r+s}}\cdot\frac{(2r)!(2s)!}{r!s!(r+s)!}$$

As for the other assertion, regarding the arcsine laws, that is clear too.
\end{proof}

Getting now to our second question about ellipses, what we can say is as follows:

\index{length of ellipse}
\index{elliptic integral}

\begin{theorem}
For an ellipse $F$ coming via $(x/a)^2+(y/b)^2=1$, its length is
$$L=\int_0^{2\pi}\sqrt{a^2\sin^2t+b^2\cos^2t}\,dt$$
with this integral being generically not computable. Thus, with respect to the measure on $F$ obtained via the path integral over it, the joint moments of the coordinates are
$$E(x^{2r}y^{2s})=\frac{1}{L}\int_0^{2\pi}(a\cos t)^{2r}(b\sin t)^{2s}\sqrt{a^2\sin^2t+b^2\cos^2t}\,dt$$
and with the integral on the right being generically not computable, either.
\end{theorem}

\begin{proof}
This is something quite standard, and a bit surprising, as follows:

\medskip

(1) To start with, what is the length of a curve $\gamma:[a,b]\to\mathbb R^2$? Good question, and in answer, a physicist would say that this is the quantity obtained by integrating the magnitude of the velocity vector over the curve, with respect to time, leading to:
$$L(\gamma)=\int_a^b||\gamma'(t)||dt$$

(2) Regarding now mathematicians, these would say that the length of a curve is the following quantity, with $(t_0=a,t_1,\ldots,t_{n-1},t_n=b)$ being a uniform division of $(a,b)$:
$$L(\gamma)=\lim_{n\to\infty}\sum_{i=1}^n||\gamma(t_i)-\gamma(t_{i-1})||$$

But, by using the fundamental theorem of calculus, this is the same as (1).

\medskip

(3) Getting now to the ellipses, we can compute their length as follows:
\begin{eqnarray*}
L
&=&\int_0^{2\pi}\sqrt{\left(\frac{dx}{dt}\right)^2+\left(\frac{dy}{dt}\right)^2}\,dt\\
&=&\int_0^{2\pi}\sqrt{\left(\frac{da\cos t}{dt}\right)^2+\left(\frac{db\sin t}{dt}\right)^2}\,dt\\
&=&\int_0^{2\pi}\sqrt{a^2\sin^2t+b^2\cos^2t}\,dt
\end{eqnarray*}

(4) Now the point is that when $a=b=R$ we get of course $L=2\pi R$, as we should, but in general, when $a\neq b$, there is no trick for computing the above integral.

\medskip

(5) As for the last assertion, this is something self-explanatory, and we will leave this as an exercise, further meditating at all this, elliptic integrals.
\end{proof}

Finally, getting to our third question about ellipses, things are quite hard here:

\index{Kepler 2}

\begin{comment}
In order to deal with our third question about ellipses, regarding the gravitational travel integration, what we need to know is the Kepler 2 law, stating that the radius vector from the Sun to a planet sweeps equal areas in equal times.
\end{comment}

In short, what we have here is a quite interesting mechanics problem, which is not exactly of trivial type, and we will leave some exploration here, after of course some due learning of the findings of Kepler and Newton, as an interesting exercise.

\bigskip

And with this, end of our study of the unit circle. Good learning this was, and in what follows we will try to understand what happens to all this, in higher dimensions.

\section*{11b. Spherical integrals}

Getting now to higher dimensions, we potentially have many things to be done. Let us start with the spherical coordinates, that we will heavily use, in what follows:

\index{spherical coordinates}

\begin{theorem}
We have spherical coordinates in $N$ dimensions,
$$\begin{cases}
x_1\!\!\!&=\ r\cos t_1\\
x_2\!\!\!&=\ r\sin t_1\cos t_2\\
\vdots\\
x_{N-1}\!\!\!&=\ r\sin t_1\sin t_2\ldots\sin t_{N-2}\cos t_{N-1}\\
x_N\!\!\!&=\ r\sin t_1\sin t_2\ldots\sin t_{N-2}\sin t_{N-1}
\end{cases}$$
the corresponding Jacobian being given by the following formula,
$$J(r,t)=r^{N-1}\sin^{N-2}t_1\sin^{N-3}t_2\,\ldots\,\sin^2t_{N-3}\sin t_{N-2}$$
and with this generalizing the known formulae at $N=2,3$.
\end{theorem}

\begin{proof}
This is something that we met in chapter 10, but always good to talk about it again. To start with, there are several possible conventions for the spherical coordinates, and we will use here the above ones, designed to look good in arbitrary $N$ dimensions. Next, regarding the Jacobian, by developing over the last column, we obtain:
\begin{eqnarray*}
J_N
&=&r\sin t_1\ldots\sin t_{N-2}\sin t_{N-1}\times \sin t_{N-1}J_{N-1}\\
&+&r\sin t_1\ldots \sin t_{N-2}\cos t_{N-1}\times\cos t_{N-1}J_{N-1}\\
&=&r\sin t_1\ldots\sin t_{N-2}(\sin^2 t_{N-1}+\cos^2 t_{N-1})J_{N-1}\\
&=&r\sin t_1\ldots\sin t_{N-2}J_{N-1}
\end{eqnarray*}

Thus, by recurrence, we are led to the formula in the statement, which by the way generalizes the well-known formulae at $N=2,3$, namely $J_2=r$ and $J_3=r^2\sin t_1$.
\end{proof}

Good news, we can now integrate over the spheres, as follows:

\index{spherical integral}
\index{double factorials}

\begin{theorem}[Wallis 3]
The polynomial integrals over the unit sphere $S^{N-1}_\mathbb R\subset\mathbb R^N$, with respect to the uniform mass $1$ measure, are given by the formula
$$\int_{S^{N-1}_\mathbb R}x_1^{k_1}\ldots x_N^{k_N}\,dx=\frac{(N-1)!!k_1!!\ldots k_N!!}{(N+\sum k_i-1)!!}$$
valid when all exponents $k_i$ are even. If an exponent is odd, the integral vanishes.
\end{theorem}

\begin{proof}
This is something routine, generalizing Wallis 2, as follows:

\medskip

(1) Assume first that one of the exponents $k_i$ is odd. We can make then the following change of variables, which shows that the integral in the statement vanishes:
$$x_i\to-x_i$$

(2) Assume now that all the exponents $k_i$ are even. As a first observation, the result holds indeed at $N=2$, due to the Wallis formula from Theorem 11.4, which reads:
$$\int_0^{\pi/2}\cos^pt\sin^qt\,dt
=\left(\frac{\pi}{2}\right)^{\varepsilon(p)\varepsilon(q)}\frac{p!!q!!}{(p+q+1)!!}
=\frac{p!!q!!}{(p+q+1)!!}$$

(3) In the general case now, where the dimension $N\in\mathbb N$ is arbitrary, the integral in the statement can be written in spherical coordinates, as follows:
$$I=\frac{2^N}{A}\int_0^{\pi/2}\ldots\int_0^{\pi/2}x_1^{k_1}\ldots x_N^{k_N}J\,dt_1\ldots dt_{N-1}$$

To be more precise, here $A$ is the area of the sphere, $J$ is the Jacobian, and the $2^N$ factor comes from the restriction to the $1/2^N$ part of the sphere where all coordinates are positive. According to our formulae from chapter 10, the normalization constant is:
$$\frac{2^N}{A}=\left(\frac{2}{\pi}\right)^{[N/2]}(N-1)!!$$

As for the unnormalized integral, by using Theorem 11.12, this is given by:
\begin{eqnarray*}
I'=\int_0^{\pi/2}\ldots\int_0^{\pi/2}
&&(\cos t_1)^{k_1}
(\sin t_1\cos t_2)^{k_2}\\
&&\vdots\\
&&(\sin t_1\sin t_2\ldots\sin t_{N-2}\cos t_{N-1})^{k_{N-1}}\\
&&(\sin t_1\sin t_2\ldots\sin t_{N-2}\sin t_{N-1})^{k_N}\\
&&\sin^{N-2}t_1\sin^{N-3}t_2\ldots\sin^2t_{N-3}\sin t_{N-2}\,dt_1\ldots dt_{N-1}
\end{eqnarray*}

(4) By rearranging the terms, we obtain the following formula:
\begin{eqnarray*}
I'
&=&\int_0^{\pi/2}\cos^{k_1}t_1\sin^{k_2+\ldots+k_N+N-2}t_1\,dt_1\\
&&\int_0^{\pi/2}\cos^{k_2}t_2\sin^{k_3+\ldots+k_N+N-3}t_2\,dt_2\\
&&\vdots\\
&&\int_0^{\pi/2}\cos^{k_{N-2}}t_{N-2}\sin^{k_{N-1}+k_N+1}t_{N-2}\,dt_{N-2}\\
&&\int_0^{\pi/2}\cos^{k_{N-1}}t_{N-1}\sin^{k_N}t_{N-1}\,dt_{N-1}
\end{eqnarray*}

Now by using the above-mentioned Wallis formula at $N=2$, this gives:
\begin{eqnarray*}
I'
&=&\frac{k_1!!(k_2+\ldots+k_N+N-2)!!}{(k_1+\ldots+k_N+N-1)!!}\left(\frac{\pi}{2}\right)^{\varepsilon(N-2)}\\
&&\frac{k_2!!(k_3+\ldots+k_N+N-3)!!}{(k_2+\ldots+k_N+N-2)!!}\left(\frac{\pi}{2}\right)^{\varepsilon(N-3)}\\
&&\vdots\\
&&\frac{k_{N-2}!!(k_{N-1}+k_N+1)!!}{(k_{N-2}+k_{N-1}+l_N+2)!!}\left(\frac{\pi}{2}\right)^{\varepsilon(1)}\\
&&\frac{k_{N-1}!!k_N!!}{(k_{N-1}+k_N+1)!!}\left(\frac{\pi}{2}\right)^{\varepsilon(0)}
\end{eqnarray*}

(5) But since the double factorials telescope, we obtain, after simplification:
$$I'=\frac{k_1!!\ldots k_N!!}{(\sum k_i+N-1)!!}\left(\frac{\pi}{2}\right)^{[N/2]}$$

Now by multiplying by $2^N/A$, we are led to the formula in the statement.
\end{proof}

Next, we have the following useful generalization of the above formula:

\index{spherical integral}

\begin{theorem}[Wallis 3']
We have the following integration formula over the sphere $S^{N-1}_\mathbb R\subset\mathbb R^N$, with respect to the uniform measure, valid for any exponents $k_i\in\mathbb N$,
$$\int_{S^{N-1}_\mathbb R}|x_1^{k_1}\ldots x_N^{k_N}|\,dx=\left(\frac{2}{\pi}\right)^{\Sigma(k_1,\ldots,k_N)}\frac{(N-1)!!k_1!!\ldots k_N!!}{(N+\sum k_i-1)!!}$$
with $\Sigma=[odds/2]$ if $N$ is odd and $\Sigma=[(odds+1)/2]$ if $N$ is even, where ``odds'' denotes the number of odd numbers in the sequence $k_1,\ldots,k_N$.
\end{theorem}

\begin{proof}
As before, the formula holds at $N=2$, due to Theorem 11.4. In general, the integral in the statement can be written in spherical coordinates, as follows:
$$I=\left(\frac{2}{\pi}\right)^{[N/2]}(N-1)!!\int_0^{\pi/2}\ldots\int_0^{\pi/2}x_1^{k_1}\ldots x_N^{k_N}J\,dt_1\ldots dt_{N-1}$$

The unnormalized integral on the right can be written, as before, as follows:
\begin{eqnarray*}
I'
&=&\int_0^{\pi/2}\cos^{k_1}t_1\sin^{k_2+\ldots+k_N+N-2}t_1\,dt_1\\
&&\int_0^{\pi/2}\cos^{k_2}t_2\sin^{k_3+\ldots+k_N+N-3}t_2\,dt_2\\
&&\vdots\\
&&\int_0^{\pi/2}\cos^{k_{N-2}}t_{N-2}\sin^{k_{N-1}+k_N+1}t_{N-2}\,dt_{N-2}\\
&&\int_0^{\pi/2}\cos^{k_{N-1}}t_{N-1}\sin^{k_N}t_{N-1}\,dt_{N-1}
\end{eqnarray*}

Now by using the Wallis 2 formula, from Theorem 11.4, we obtain:
\begin{eqnarray*}
I'
&=&\frac{\pi}{2}\cdot\frac{k_1!!(k_2+\ldots+k_N+N-2)!!}{(k_1+\ldots+k_N+N-1)!!}\left(\frac{2}{\pi}\right)^{\delta(k_1,k_2+\ldots+k_N+N-2)}\\
&&\frac{\pi}{2}\cdot\frac{k_2!!(k_3+\ldots+k_N+N-3)!!}{(k_2+\ldots+k_N+N-2)!!}\left(\frac{2}{\pi}\right)^{\delta(k_2,k_3+\ldots+k_N+N-3)}\\
&&\vdots\\
&&\frac{\pi}{2}\cdot\frac{k_{N-2}!!(k_{N-1}+k_N+1)!!}{(k_{N-2}+k_{N-1}+k_N+2)!!}\left(\frac{2}{\pi}\right)^{\delta(k_{N-2},k_{N-1}+k_N+1)}\\
&&\frac{\pi}{2}\cdot\frac{k_{N-1}!!k_N!!}{(k_{N-1}+k_N+1)!!}\left(\frac{2}{\pi}\right)^{\delta(k_{N-1},k_N)}
\end{eqnarray*}

In order to compute this quantity, let us denote by $F$ the part involving the double factorials, and by $P$ the part involving the powers of $\pi/2$, so that we have:
$$I'=F\cdot P$$

Regarding $F$, there are many cancellations there, and we end up with:
$$F=\frac{k_1!!\ldots k_N!!}{(\sum k_i+N-1)!!}$$

As in what regards $P$, the $\delta$ exponents on the right sum up to the following number:
$$\Delta(k_1,\ldots,k_N)=\sum_{i=1}^{N-1}\delta(k_i,k_{i+1}+\ldots+k_N+N-i-1)$$

In other words, with this notation, the above formula reads:
\begin{eqnarray*}
I'
&=&\left(\frac{\pi}{2}\right)^{N-1}\frac{k_1!!k_2!!\ldots k_N!!}{(k_1+\ldots+k_N+N-1)!!}\left(\frac{2}{\pi}\right)^{\Delta(k_1,\ldots,k_N)}\\
&=&\left(\frac{2}{\pi}\right)^{\Delta(k_1,\ldots,k_N)-N+1}\frac{k_1!!k_2!!\ldots k_N!!}{(k_1+\ldots+k_N+N-1)!!}\\
&=&\left(\frac{2}{\pi}\right)^{\Sigma(k_1,\ldots,k_N)-[N/2]}\frac{k_1!!k_2!!\ldots k_N!!}{(k_1+\ldots+k_N+N-1)!!}
\end{eqnarray*}

To be more precise, the formula relating $\Delta$ to $\Sigma$ follows from a number of simple observations, the first of which being the fact that, due to obvious parity reasons, the sequence of $\delta$ numbers appearing in the definition of $\Delta$ cannot contain two consecutive zeroes. Now together with $I=(2^N/V)I'$, this gives the formula in the statement.
\end{proof}

Finally, we have the following complex version of Theorems 11.13 and 11.14:

\index{spherical integral}

\begin{theorem}[Wallis 4]
We have the following integration formula over the complex sphere $S^{N-1}_\mathbb C\subset\mathbb C^N$, with respect to the uniform mass $1$ measure, 
$$\int_{S^{N-1}_\mathbb C}|z_1|^{2k_1}\ldots|z_N|^{2k_N}\,dz=\frac{(N-1)!k_1!\ldots k_n!}{(N+\sum k_i-1)!}$$
valid for any exponents $k_i\in\mathbb N$. As for the other polynomial integrals in $z_1,\ldots,z_N$ and their conjugates $\bar{z}_1,\ldots,\bar{z}_N$, these all vanish.
\end{theorem}

\begin{proof}
Consider a polynomial integral over $S^{N-1}_\mathbb C$, containing the same number of plain and conjugated variables, as to not vanish trivially, written as follows:
$$I=\int_{S^{N-1}_\mathbb C}z_{i_1}\bar{z}_{i_2}\ldots z_{i_{2k-1}}\bar{z}_{i_{2k}}\,dz$$

By using transformations of type $p\to\lambda p$ with $|\lambda|=1$, we see that this integral $I$ vanishes, unless each $z_a$ appears as many times as $\bar{z}_a$ does, and this gives the last assertion. So, assume now that we are in the non-vanishing case. Then the $k_a$ copies of $z_a$ and the $k_a$ copies of $\bar{z}_a$ produce by multiplication a factor $|z_a|^{2k_a}$, so we have:
$$I=\int_{S^{N-1}_\mathbb C}|z_1|^{2k_1}\ldots|z_N|^{2k_N}\,dz$$

Now by using the standard identification $S^{N-1}_\mathbb C\simeq S^{2N-1}_\mathbb R$, we obtain:
\begin{eqnarray*}
I
&=&\int_{S^{2N-1}_\mathbb R}(x_1^2+y_1^2)^{k_1}\ldots(x_N^2+y_N^2)^{k_N}\,d(x,y)\\
&=&\sum_{r_1\ldots r_N}\binom{k_1}{r_1}\ldots\binom{k_N}{r_N}\int_{S^{2N-1}_\mathbb R}x_1^{2k_1-2r_1}y_1^{2r_1}\ldots x_N^{2k_N-2r_N}y_N^{2r_N}\,d(x,y)
\end{eqnarray*}

By using the formula in Theorem 11.13 for integrating, we obtain:
\begin{eqnarray*}
&&I\\
&=&\sum_{r_1\ldots r_N}\binom{k_1}{r_1}\ldots\binom{k_N}{r_N}\frac{(2N-1)!!(2r_1)!!\ldots(2r_N)!!(2k_1-2r_1)!!\ldots (2k_N-2r_N)!!}{(2N+2\sum k_i-1)!!}\\
&=&\sum_{r_1\ldots r_N}\binom{k_1}{r_1}\ldots\binom{k_N}{r_N}\frac{2^{N-1}(N-1)!\prod(2r_i)!/(2^{r_i}r_i!)\prod(2k_i-2r_i)!/(2^{k_i-r_i}(k_i-r_i)!)}{2^{N+\sum k_i-1}(N+\sum k_i-1)!}\\
&=&\sum_{r_1\ldots r_N}\binom{k_1}{r_1}\ldots\binom{k_N}{r_N}
\frac{(N-1)!(2r_1)!\ldots (2r_N)!(2k_1-2r_1)!\ldots (2k_N-2r_N)!}{4^{\sum k_i}(N+\sum k_i-1)!r_1!\ldots r_N!(k_1-r_1)!\ldots (k_N-r_N)!}
\end{eqnarray*}

Now observe that can rewrite this quantity in the following way:
\begin{eqnarray*}
&&I\\
&=&\sum_{r_1\ldots r_N}\frac{k_1!\ldots k_N!(N-1)!(2r_1)!\ldots (2r_N)!(2k_1-2r_1)!\ldots (2k_N-2r_N)!}{4^{\sum k_i}(N+\sum k_i-1)!(r_1!\ldots r_N!(k_1-r_1)!\ldots (k_N-r_N)!)^2}\\
&=&\sum_{r_1}\binom{2r_1}{r_1}\binom{2k_1-2r_1}{k_1-r_1}\ldots\sum_{r_N}\binom{2r_N}{r_N}\binom{2k_N-2r_N}{k_N-r_N}\frac{(N-1)!k_1!\ldots k_N!}{4^{\sum k_i}(N+\sum k_i-1)!}\\
&=&4^{k_1}\times\ldots\times 4^{k_N}\times\frac{(N-1)!k_1!\ldots k_N!}{4^{\sum k_i}(N+\sum k_i-1)!}\\
&=&\frac{(N-1)!k_1!\ldots k_N!}{(N+\sum k_i-1)!}
\end{eqnarray*}

To be more precise, we have used here the following formula, that we know well from chapter 5, from our study there of the complex normal variables:
$$\sum_r\binom{2r}{r}\binom{2l-2r}{l-r}=4^l$$

Thus, we have obtained the formula in the statement, as desired.
\end{proof}

\section*{11c. Hyperspherical laws}

Good news, we can do now some probability over the spheres. Our source of inspiration will be Theorems 11.5 and 11.6, dealing with the 2D real case, and as technical tools, we will use the Wallis formulae from Theorems 11.13, 11.14 and 11.15. We first have:

\index{hyperspherical law}
\index{normal law}

\begin{theorem}
The moments and density of the hyperspherical laws, that of the coordinates $x_i:S^{N-1}_\mathbb R\to[-1,1]$, are as follows, with $e_N=2,\pi$ for $N$ odd, even,
$$M_{2r}=\frac{(N-1)!!(2r)!!}{(N+2r-1)!!}\quad,\quad\mu_N=\frac{(N-1)!!}{e_N(N-2)!!}(1-x^2)^{(N-3)/2}dx$$
with the corresponding expectation, variance, skewness and kurtosis being given by:
$$E=0\quad,\quad V=\frac{1}{N}\quad,\quad\gamma=0\quad,\quad\kappa=\frac{3N}{N+2}$$
At $N=1,2,3,4$ we obtain the following well-known laws, all supported by $[-1,1]$:
$$\mu_1=\frac{\delta_{-1}+\delta_1}{2}\quad,\quad\mu_2=\frac{1}{\pi\sqrt{1-x^2}}\,dx
\quad,\quad\mu_3=\frac{dx}{2}\quad,\quad\mu_4=\frac{2}{\pi}\sqrt{1-x^2}\,dx$$
Also, the joint moments of the variables $x_i:S^{N-1}_\mathbb R\to[-1,1]$ are given by
$$E(x_1^{2r_1}\ldots x_N^{2r_N})=\frac{(N-1)!!(2r_1)!!\ldots(2r_N)!!}{(N+\sum 2r_i-1)!!}$$
and $y_i=\sqrt{N}x_i$ become standard normal, following $g_1$, and independent with $N\to\infty$.
\end{theorem}

\begin{proof}
This is something standard, based on Theorem 11.13, as follows:

\medskip

(1) To start with, as before in Theorem 11.5, we can certainly denote the points of the sphere $x\in S^{N-1}_\mathbb R$ as $x=(x_1,\ldots,x_N)$, and therefore talk about the $N$ random variables $x_i:S^{N-1}_\mathbb R\to[-1,1]$, with respect to the uniform, mass 1 measure on the sphere.

\medskip

(2) Next, the joint moment formula, which implies the individual moment formula, is the Wallis 3 formula, from Theorem 11.13. So, we got our laws, and what is left is to compute the density, discuss $E,V,\gamma,\kappa$, discuss $N=1,2,3,4$, and discuss $N=\infty$ too.

\medskip

(3) The simplest question is the asymptotic one. Indeed, with $N\to\infty$ we have:
$$E(x_i^{2r})
=\frac{(N-1)!!}{(N+2r-1)!!}\times(2r)!!
\simeq N^{-r}(2r)!!
=N^{-r}M_{2r}(g_1)$$

Thus the variables $y_i=\sqrt{N}x_i$ become standard normal with $N\to\infty$, as stated.

\medskip

(4) Getting now to the joint moments of $y_i=\sqrt{N}x_i$, with $N\to\infty$ we have:
\begin{eqnarray*}
E(y_1^{2r_1}\ldots y_N^{2r_N})
&=&N^{\sum r_i}E(x_1^{2r_1}\ldots x_N^{2r_N})\\
&=&N^{\sum r_i}\frac{(N-1)!!(2r_1)!!\ldots(2r_N)!!}{(N+\sum 2r_i-1)!!}\\
&=&N^{\sum r_i}\frac{(N-1)!!}{(N+\sum 2r_i-1)!!}\times(2r_1)!!\ldots(2r_N)!!\\
&\simeq&N^{\sum r_i}N^{-\sum r_i}(2r_1)!!\ldots(2r_N)!!\\
&=&(2r_1)!!\ldots(2r_N)!!
\end{eqnarray*}

Thus the asymptotic joint moments factorize, and by using the moment characterization of independence we conclude that $y_i$ are asymptotically independent, as stated.

\medskip

(5) Still talking moments, let us do now the $E,V,\gamma,\kappa$ computations. Our measure being centered we have $E=\gamma=0$, and for computing $V,\kappa$ we will need:
$$M_2=\frac{(N-1)!!}{(N+1)!!}=\frac{1}{N}\quad,\quad M_4=\frac{(N-1)!!3}{(N+3)!!}=\frac{3}{N(N+2)}$$

Thus the variance and kurtosis are given by the following formulae:
$$V=M_2=\frac{1}{N}\quad,\quad\kappa=\frac{M_4}{V^2}=\frac{3N}{N+2}$$

(6) As another thing that we can do with moments, we can quickly discuss the cases $N=1,2,3,4$ as well. Indeed, at these particular values of $N$, the moments become:
$$M_{2r}^{(1)}=\frac{(2r)!!}{(2r)!!}=1$$
$$M_{2r}^{(2)}=\frac{(2r)!!}{(2r+1)!!}=\frac{(2r)!/2^rr!}{2^rr!}=\frac{1}{4^r}\binom{2r}{r}$$
$$M_{2r}^{(3)}=\frac{(2r)!!}{(2r+2)!!}=\frac{1}{2r+1}$$
$$M_{2r}^{(4)}=\frac{2\cdot(2r)!!}{(2r+3)!!}=\frac{2\cdot(2r)!/2^rr!}{2^{r+1}(r+1)!}=\frac{1}{4^r}\cdot\frac{1}{r+1}\binom{2r}{r}$$

But we recognize here the moments of the Bernoulli, arcsine, uniform and semicircle laws on $[-1,1]$. Of course, there are far more things that can be said here, and we will leave some geometric thinking here, short-circuiting the above, as an exercise.

\medskip

(7) Summarizing, job done, and we are left with finding the density. But this can be done again via moments, because the moments of the measure in the statement are:
\begin{eqnarray*}
I_{2r}
&=&\frac{(N-1)!!}{e_N(N-2)!!}\int_{-1}^1x^{2r}(1-x^2)^{(N-3)/2}dx\\
&=&\frac{2(N-1)!!}{e_N(N-2)!!}\int_0^1x^{2r}(1-x^2)^{(N-3)/2}dx\\
&=&\frac{2(N-1)!!}{e_N(N-2)!!}\int_0^{\pi/2}\cos^{2r}t\sin^{N-3}t\cdot\sin tdt\\
&=&\frac{2(N-1)!!}{e_N(N-2)!!}\int_0^{\pi/2}\cos^{2r}t\sin^{N-2}t\,dt\\
&=&\frac{2(N-1)!!}{e_N(N-2)!!}\left(\frac{\pi}{2}\right)^{\varepsilon(N)}\frac{(2r)!!(N-2)!!}{(N+2r-1)!!}\\
&=&\frac{2}{e_N}\left(\frac{\pi}{2}\right)^{\varepsilon(N)}\frac{(N-1)!!(2r)!!}{(N+2r-1)!!}\\
&=&\frac{(N-1)!!(2r)!!}{(N+2r-1)!!}
\end{eqnarray*}

Thus, we are led to the conclusions in the statement. 
\end{proof}

Regarding now the squares of the hyperspherical variables, we have:

\index{beta laws}

\begin{theorem}[version]
The squared coordinates $x_i^2:S^{N-1}_\mathbb R\to[0,1]$ follow beta laws, having moments and density as follows, with $e_N=2,\pi$ for $N$ odd, even,
$$M_r=\frac{(N-1)!!(2r)!!}{(N+2r-1)!!}\quad,\quad\nu_N=\frac{(N-1)!!}{e_N(N-2)!!}\,x^{-1/2}(1-x)^{(N-3)/2}dx$$
the expectation being $E=1/N$, and the variance, skewness and kurtosis being:
$$V=\frac{2N-2}{N^2(N+2)}\ ,\quad 
\gamma=\frac{N-2}{N+4}\sqrt{\frac{8(N+2)}{N-1}}\ ,\quad 
\kappa=\frac{3(N+2)(5N^2-13N+12)}{(N-1)(N+4)(N+6)}$$
At $N=1,2,3,4$ we obtain the following well-known laws, all supported by $[0,1]$:
$$\nu_1=\delta_1\quad,\quad\nu_2=\frac{1}{\pi\sqrt{x-x^2}}\,dx
\quad,\quad\nu_3=\frac{dx}{2\sqrt{x}}\quad,\quad\nu_4=\frac{2}{\pi}\sqrt{x^{-1}-1}\,dx$$
Also, the joint moments of the variables $x_i^2:S^{N-1}_\mathbb R\to[0,1]$ are given by
$$E(x_1^{2r_1}\ldots x_N^{2r_N})=\frac{(N-1)!!(2r_1)!!\ldots(2r_N)!!}{(N+\sum 2r_i-1)!!}$$
and $y_i=Nx_i^2$ become chi squared, following $\chi_1^2$, and independent with $N\to\infty$.
\end{theorem}

\begin{proof}
This is a straightforward remake of Theorem 11.16, as follows:

\medskip

(1) The joint moment formula is the one in Theorem 11.16, unchanged, so we have this, and as consequences, the individual moment formula, and the $N\to\infty$ result too.

\medskip

(2) Regarding the particular cases $N=1,2,3,4$, the measures in the statement come from those in Theorem 11.16, by squaring, or via moments, or via direct geometry.

\medskip

(3) In general now, the density formula formula in the statement comes from the density formula in Theorem 11.16, by squaring:
\begin{eqnarray*}
E(\varphi(x_i^2))
&=&\frac{2(N-1)!!}{e_N(N-2)!!}\int_0^1\varphi(x^2)(1-x^2)^{(N-3)/2}dx\\
&=&\frac{2(N-1)!!}{e_N(N-2)!!}\int_0^1\varphi(y)(1-y)^{(N-3)/2}\frac{dy}{2\sqrt{y}}\\
&=&\frac{(N-1)!!}{e_N(N-2)!!}\int_0^1\varphi(y)y^{-1/2}(1-y)^{(N-3)/2}dy
\end{eqnarray*}

(4) Regarding now the computation of $E,V,\gamma,\kappa$, the low order moments are:
$$M_1=\frac{1}{N}$$
$$M_2=\frac{3}{N(N+2)}$$
$$M_3=\frac{15}{N(N+2)(N+4)}$$
$$M_4=\frac{105}{N(N+2)(N+4)(N+6)}$$

We conclude that the mean is $E=1/N$, and that the variance is given by:
$$V=\frac{3}{N(N+2)}-\frac{1}{N^2}=\frac{3N-N-2}{N^2(N+2)}
=\frac{2N-2}{N^2(N+2)}$$

Next, the third central moment is given by the following formula:
$$M_3'=M_3-3EM_2+2E^3=\frac{8(N-1)(N-2)}{N^3(N+2)(N+4)}$$

Thus, we obtain the following formula, for the corresponding skewness:
$$\gamma=\frac{8(N-1)(N-2)}{N^3(N+2)(N+4)}\cdot\frac{N^2(N+2)}{2(N-1)}\sqrt{\frac{N^2(N+2)}{2(N-1)}}=\frac{N-2}{N+4}\sqrt{\frac{8(N+2)}{N-1}}$$

Next, the fourth central moment is given by the following formula:
$$M_4'=M_4-4EM_3+6E^2M_2-3E^4=\frac{12(N-1)(5N^2-13N+12)}{N^4(N+2)(N+4)(N+6)}$$

Thus, we obtain the following formula, for the corresponding kurtosis:
\begin{eqnarray*}
\kappa
&=&\frac{12(N-1)(5N^2-13N+12)}{N^4(N+2)(N+4)(N+6)}\cdot\frac{N^4(N+2)^2}{4(N-1)^2}\\
&=&\frac{3(N+2)(5N^2-13N+12)}{(N-1)(N+4)(N+6)}
\end{eqnarray*}

(5) Finally, as a bonus result, the covariances of our variables are given by:
$$cov(x_i^2,x_j^2)=\frac{(N-1)!!}{(N+3)!!}-\frac{1}{N^2}=\frac{1}{N(N+2)}-\frac{1}{N^2}=-\frac{2}{N^2(N+2)}$$

Thus, the corresponding correlations are given by the following formula:
$$\rho_{x_i^2x_j^2}=-\frac{2}{N^2(N+2)}\cdot\frac{N^2(N+2)}{2N-2}=-\frac{1}{N-1}$$

And with the worries at $N=2$ already discussed, in the proof of Theorem 11.6, I will leave the $N<2$ discussion to you. I mean, that $\rho<-1$ is not normal, right.
\end{proof}

Finally, regarding the absolute values of the hyperspherical variables, we have:

\begin{theorem}[version]
The moments and density of the absolute values of the coordinates $|x_i|:S^{N-1}_\mathbb R\to[0,1]$ are as follows, with $e_N=2,\pi$ for $N$ odd, even,
$$M_k=\left(\frac{2}{\pi}\right)^{\delta_{2|N}\delta_{2\not\:|\,k}}\frac{(N-1)!!k!!}{(N+k-1)!!}\quad,\quad\eta_N=\frac{2(N-1)!!}{e_N(N-2)!!}(1-x^2)^{(N-3)/2}dx$$
and at $N=1,2,3,4$ we obtain the following well-known laws, all supported by $[0,1]$:
$$\mu_1=\delta_1\quad,\quad\mu_2=\frac{2}{\pi\sqrt{1-x^2}}\,dx
\quad,\quad\mu_3=dx\quad,\quad\mu_4=\frac{4}{\pi}\sqrt{1-x^2}\,dx$$
The joint moments are as follows, with $\Sigma=[odds/2]$ if $N$ is odd and $\Sigma=[(odds+1)/2]$ if $N$ is even, where ``odds'' is the number of odd numbers in the sequence $k_1,\ldots,k_N$,
$$E(|x_1|^{k_1}\ldots|x_N|^{k_N})=\left(\frac{2}{\pi}\right)^{\Sigma(k_1,\ldots,k_N)}\frac{(N-1)!!k_1!!\ldots k_N!!}{(N+\sum k_i-1)!!}$$
and $y_i=\sqrt{N}|x_i|$ become standard chi, following $\chi_1$, and independent with $N\to\infty$.
\end{theorem}

\begin{proof}
This follows indeed from what we have, the idea being as follows:

\medskip

(1) The density comes from the density in Theorem 11.16, by doubling on $[0,1]$, and the $N=1,2,3,4$ particular cases come from Theorem 11.16 too, again by doubling.

\medskip

(2) We have the joint moment formula from Theorem 11.14, which produces in particular the individual moment formula, and the asymptotic result as well.

\medskip

(3) Observe that we have skipped saying something about $E,V,\gamma,\kappa$, and this because, as in the usual chi case, the low order moments do not look very good:
$$M_1=\left(\frac{2}{\pi}\right)^{\delta_{2|N}}\frac{(N-1)!!}{N!!}\quad,\quad M_2=\frac{1}{N}$$
$$M_3=\left(\frac{2}{\pi}\right)^{\delta_{2|N}}\frac{(N-1)!!2}{(N+2)!!}\quad,\quad M_4=\frac{3}{N(N+2)}$$

(4) Finally, as a bonus computation, let us compute $cov(|x_i|,|x_j|)$. We have:
$$E(|x_ix_j|)=\left(\frac{2}{\pi}\right)^{\Sigma(110\ldots0)}\frac{(N-1)!!}{(N+1)!!}=\frac{2}{\pi}\cdot\frac{1}{N}$$

Which looks nice, but when passing to $cov$, that alien $M_1$ number will appear:
$$cov(|x_i|,|x_j|)=\frac{2}{\pi}\cdot\frac{1}{N}-\left(\frac{2}{\pi}\right)^{2\delta_{2|N}}\left(\frac{(N-1)!!}{N!!}\right)^2$$

In practice, at $N=2,3,4$ we obtain from this the following figures:
$$cov^{(2)}=\frac{\pi-4}{\pi^2}\quad,\quad cov^{(3)}=\frac{8-3\pi}{12\pi}
\quad,\quad cov^{(4)}=\frac{9\pi-32}{18\pi^2}$$

And with this, end of our study of the absolute values of sphere coordinates.
\end{proof}

In the complex case now, we have a similar result, as follows:

\begin{theorem}
The moments of the complex hyperspherical variables are
$$\int_{S^{N-1}_\mathbb C}|z_i|^{2k}dz=\frac{(N-1)!k!}{(N+k-1)!}$$
with the other moments vanishing. The normalized complex hyperspherical variables 
$$w_i=\sqrt{N}z_i$$
become standard complex normal, following $G_1$, and independent with $N\to\infty$.
\end{theorem}

\begin{proof}
This is a standard complex remake of Theorem 11.16, as follows:

\medskip

(1) To start with, we can certainly denote the points of the complex sphere $z\in S^{N-1}_\mathbb C$ as $z=(z_1,\ldots,z_N)$, and therefore talk about the $N$ random variables $z_i:S^{N-1}_\mathbb C\to D$, with $D\subset\mathbb C$ being the unit disk, with respect to the uniform measure on the sphere.

\medskip

(2) In order to study the individual and joint laws of these variables $z_1,\ldots,z_N$, and their asymptotics, we can use the Wallis 4 formula, from Theorem 11.15, namely:
$$E\left(|z_1|^{2k_1}\ldots|z_N|^{2k_N}\right)=\frac{(N-1)!k_1!\ldots k_n!}{(N+\sum k_i-1)!}$$

(3) As a first observation, by assuming that all exponents are 0, except for the $i$-th one, we obtain the individual moment formula in the statement. Next, observe that by using this, we have the following estimate, in the $N\to\infty$ limit:
$$E(|z_i|^{2k})
=\frac{(N-1)!}{(N+k-1)!}\times k!
\simeq N^{-k}k!
=N^{-k}M_{2k}(G_1)$$

Thus $w_i=\sqrt{N}z_i$ become standard complex normal with $N\to\infty$, as stated.

\medskip

(4) Getting now to the joint moments of $w_i=\sqrt{N}z_i$, with $N\to\infty$ we have:
\begin{eqnarray*}
\int_{S^{N-1}_\mathbb C}|w_1|^{2k_1}\ldots |w_N|^{2k_N}\,dz
&=&N^{\sum k_i}\int_{S^{N-1}_\mathbb C}|z_1|^{2k_1}\ldots|z_N|^{2k_N}\,dz\\
&=&N^{\sum k_i}\frac{(N-1)!k_1!\ldots k_N!}{(N+\sum k_i-1)!}\\
&=&N^{\sum k_i}\frac{(N-1)!}{(N+\sum k_i-1)!}\times k_1!\ldots k_N!\\
&\simeq&N^{\sum k_i}N^{-\sum k_i}k_1!\ldots k_N!\\
&=&k_1!\ldots k_N!
\end{eqnarray*}

Thus the joint moments factorize, and $w_i$ are asymptotically independent.
\end{proof}

As a comment now on all this, in relation with groups, by using the standard actions $O_N\curvearrowright S^{N-1}_\mathbb R$ and $U_N\curvearrowright S^{N-1}_\mathbb C$ we can see that the standard coordinates $u_{ij}\in C(O_N)$ and $u_{ij}\in C(U_N)$ are real and complex hyperspherical. The asymptotics from Theorems 11.16 and 11.19 can be alternatively recovered from the Weingarten formula, but the precise results, at fixed $N\in\mathbb N$, not exactly, or at least not with bare hands. Good to know.

\section*{11d. Clebsch-Gordan}

It has become customary in this book to do some physics at the end of each chapter, and in relation with what we did in the above, I am afraid that we will do some pure mathematics instead, that you can call of course particle physics too. To be more precise, we will be talking about the Clebsch-Gordan rules, which are essential to the functioning of both our beloved pure mathematics, and of our beloved Standard Model.

\bigskip

Getting started, we first have the following result, regarding the group $SU_2$:

\begin{theorem}
In the standard picture of the group $SU_2$, which is as follows,
$$SU_2=\left\{\begin{pmatrix}x+iy&z+it\\ -z+it&x-iy\end{pmatrix}\Big|(x,y,z,t)\in S^3_\mathbb R\right\}$$ 
the main character is $\chi=2x$, following the Wigner semicircle law $\gamma_1$.
\end{theorem}

\begin{proof}
The starting point is the defining formula for the group $SU_2$, as being the group of $2\times2$ unitary matrices having determinant one:
$$SU_2=\left\{U\in M_2(\mathbb C)\Big|U^*=U^{-1},\,\det U=1\right\}$$

In practice, if we pick $U=\binom{a\ b}{c\ d}$ of determinant 1, the equation $U^*=U^{-1}$ reads:
$$\begin{pmatrix}\bar{a}&\bar{c}\\ \bar{b}&\bar{d}\end{pmatrix}
=\begin{pmatrix}d&-b\\ -c&a\end{pmatrix}$$

Thus we must have $d=\bar{a}$, $c=-\bar{b}$, and since with these conditions imposed, the condition $\det U=1$ reads $|a|^2+|b|^2=1$, we are led to the following conclusion:
$$SU_2=\left\{\begin{pmatrix}a&b\\ -\bar{b}&\bar{a}\end{pmatrix}\Big|\ |a|^2+|b|^2=1\right\}$$

But with $a=x+iy$, $b=z+it$ this gives the formula in the statement, namely:
$$SU_2=\left\{\begin{pmatrix}x+iy&z+it\\ -z+it&x-iy\end{pmatrix}\Big|(x,y,z,t)\in S^3_\mathbb R\right\}$$

Summarizing, we have $SU_2\simeq S^3_\mathbb R$, and it is also known that the uniform measure on $SU_2$, in a group-theoretical sense, coincides with the uniform measure on $S^3_\mathbb R$. But with this, done, because if we look at the main character of $SU_2$, this character is $\chi=2x$ in the above $SU_2\simeq S^3_\mathbb R$ picture, which by Theorem 11.16 follows the Wigner law $\gamma_1$.
\end{proof}

Regarding the group $SO_3$, which is the close friend of $SU_2$, we have:

\index{Euler-Rodrigues}
\index{axis of rotation}

\begin{theorem}
In the standard picture of the elements $U\in SO_3$, namely
$$U=\begin{pmatrix}
x^2+y^2-z^2-t^2&2(yz-xt)&2(xz+yt)\\
2(xt+yz)&x^2+z^2-y^2-t^2&2(zt-xy)\\
2(yt-xz)&2(xy+zt)&x^2+t^2-y^2-z^2
\end{pmatrix}$$
with $(x,y,z,t)\in S^3_\mathbb R$, the variable $\chi+1=4x^2$ follows the Marchenko-Pastur law $\pi_1$.
\end{theorem}

\begin{proof}
This is something quite tricky, the idea being as follows:

\medskip

(1) To start with, following Euler, any rotation $U\in SO_3$ has a rotation axis. Indeed, we have the following computation, using some linear algebra magic:
\begin{eqnarray*}
\det(U-1)
&=&\det(U^t-1)\\
&=&\det(U^t(1-U))\\
&=&\det(U^t)\det(1-U)\\
&=&\det(1-U)
\end{eqnarray*}

Thus $\det(U-1)=0$, which tells us that $U$ must have a $1$-eigenvector, as desired. 

\medskip

(2) Next, following Euler-Rodrigues, we have the formula in the statement, with the parameters $(x,y,z,t)\in S^3_\mathbb R$ encapsulating the direction of the rotation axis, and the rotation angle. Alternatively, if you want to trade heavy trigonometry for heavy algebra, the Euler-Rodrigues formula comes from the double cover map $SU_2\to SO_3$, obtained via the adjoint action of $SU_2$ on $\mathbb C^2\simeq\mathbb R^4$, by removing the trivial fixed direction.

\medskip

(3) In short, quite complicated all this, and with the complete picture, needed for properly understanding the double cover map $SU_2\to SO_3$, involving as well the 4 Pauli spin matrices, appearing as matrix coefficients of $x,y,z,t$ in the formula from Theorem 11.20. Have a look at any of my algebra books, where all this is explained.

\medskip

(4) With this discussed, in the Euler-Rodrigues picture, the main character is:
\begin{eqnarray*}
\chi
&=&3x^2-y^2-z^2-t^2\\
&=&4x^2-(x^2+y^2+z^2+t^2)\\
&=&4x^2-1
\end{eqnarray*}

Thus $\chi+1=4x^2$, which by Theorem 11.16 follows the Marchenko-Pastur law $\pi_1$.
\end{proof}

As a continuation of the above, at a more conceptual level, we know from the Peter-Weyl theory from chapter 7 that the law of the main character ultimately comes by decomposing the Peter-Weyl representations $u^{\otimes k}$ into irreducibles. So, let us try to understand this. For the group $SU_2$, further building on Theorem 11.20, we have:

\index{Clebsch-Gordan}
\index{fusion rules}

\begin{theorem}
The irreducible representations of $SU_2$ are all self-adjoint, and can be labeled by positive integers, with their fusion rules being as follows,
$$r_k\otimes r_l=r_{|k-l|}+r_{|k-l|+2}+\ldots+r_{k+l}$$
called Clebsch-Gordan rules. The corresponding dimensions are $\dim r_k=k+1$.
\end{theorem}

\begin{proof}
With the knowledge that we have, we can prove this result as follows:

\medskip

(1) According to Theorem 11.20, the moments of the main character are:
$$\int_{SU_2}\chi^{2k}=C_k$$

(2) Our claim is that we can construct, by recurrence on $k\in\mathbb N$, a sequence $r_k$ of irreducible, self-adjoint and distinct representations of $SU_2$, satisfying:
$$r_0=1\quad,\quad
r_1=u\quad,\quad 
r_k+r_{k-2}=r_{k-1}\otimes r_1$$

Indeed, assume that $r_0,\ldots,r_{k-1}$ are constructed, and let us construct $r_k$. We have:
$$r_{k-1}+r_{k-3}=r_{k-2}\otimes r_1$$

Thus $r_{k-1}\subset r_{k-2}\otimes r_1$, and since $r_{k-2}$ is irreducible, by Frobenius we have:
$$r_{k-2}\subset r_{k-1}\otimes r_1$$

We conclude there exists a certain representation $r_k$ such that:
$$r_k+r_{k-2}=r_{k-1}\otimes r_1$$

(3) By recurrence, $r_k$ is self-adjoint. Now observe that according to our recurrence formula, we can split $u^{\otimes k}$ as a sum of the following type, with positive coefficients:  
$$u^{\otimes k}=c_kr_k+c_{k-2}r_{k-2}+\ldots$$

We conclude by Peter-Weyl that we have an inequality as follows, with equality precisely when $r_k$ is irreducible, and non-equivalent to the other summands $r_i$:
$$\sum_ic_i^2\leq\dim(End(u^{\otimes k}))$$

(4) But according to (1) the number on the right is $C_k$, and some straightforward combinatorics, based on the fusion rules, shows that the number on the left is $C_k$ as well. Thus we have equality in our estimate, so our representation $r_k$ is irreducible, and non-equivalent to $r_{k-2},r_{k-4},\ldots$ Moreover, this representation $r_k$ is not equivalent to $r_{k-1},r_{k-3},\ldots$ either, with this coming from $r_p\subset u^{\otimes p}$ for any $p$, and from:
$$\dim(Fix(u^{\otimes 2s+1}))=\int_{SU_2}\chi^{2s+1}=0$$

(5) Thus, we proved our claim. Now since each irreducible representation of $SU_2$ appears into some $u^{\otimes k}$, and we know how to decompose each $u^{\otimes k}$ into sums of representations $r_k$, these representations $r_k$ are all the irreducible representations of $SU_2$, and we are done with the main assertion. As for the dimension formula, this is clear.
\end{proof}

Regarding now $SO_3$, we have here a similar result, as follows:

\begin{theorem}
The irreducible representations of $SO_3$ are all self-adjoint, and can be labeled by positive integers, with their fusion rules being as follows,
$$r_k\otimes r_l=r_{|k-l|}+r_{|k-l|+1}+\ldots+r_{k+l}$$ 
also called Clebsch-Gordan rules. The corresponding dimensions are $\dim r_k=2k+1$.
\end{theorem}

\begin{proof}
As before with $SU_2$, there are many possible proofs here, which are all instructive. Here is our take on the subject, in the spirit of our proof for $SU_2$:

\medskip

(1) According to Theorem 11.21, the moments of the main character are:
$$\int_{SO_3}\chi^k=C_k$$

(2) Our claim now is that we can construct, by recurrence on $k\in\mathbb N$, a sequence $r_k$ of irreducible, self-adjoint and distinct representations of $SO_3$, satisfying:
$$r_0=1\quad,\quad
r_1=u-1\quad,\quad 
r_k+r_{k-1}+r_{k-2}=r_{k-1}\otimes r_1$$

Indeed, assume that $r_0,\ldots,r_{k-1}$ are constructed, and let us construct $r_k$. The Frobenius trick from the proof for $SU_2$ will no longer work, due to some technical reasons, so we have to invoke (1). To be more precise, by integrating characters we obtain:
$$r_{k-1},r_{k-2}\subset r_{k-1}\otimes r_1$$

We conclude that there exists a representation $r_k$ such that:
$$r_{k-1}\otimes r_1=r_k+r_{k-1}+r_{k-2}$$

(3) Once again by integrating characters, we conclude that $r_k$ is irreducible, and non-equivalent to $r_1,\ldots,r_{k-1}$, and this proves our claim. Also, since any irreducible representation of $SO_3$ must appear in some tensor power of $u$, and we can decompose each $u^{\otimes k}$ into sums of representations $r_p$, we conclude that these representations $r_p$ are all the irreducible representations of $SO_3$. Finally, the dimension formula is clear.
\end{proof}

As a conclusion to this, which is good to know, the advanced combinatorics of the Wigner and Marchenko-Pastur laws $\gamma_1$ and $\pi_1$ comes from the Clebsch-Gordan rules for $SU_2$ and $SO_3$. There are of course far more things that can be said, with for instance an alternative proof for this coming from a certain ``super-easiness'' property of $SU_2$ and $SO_3$. And, for more on this, you can have a look at the quantum algebra literature.

\section*{11e. Exercises}

This was an exciting geometric chapter, and as exercises, we have:

\begin{exercise}
Learn more about ellipses, and their story.
\end{exercise}

\begin{exercise}
Do some computations with elliptic integrals. 
\end{exercise}

\begin{exercise}
Study as well the Kepler 2 integration problem.
\end{exercise}

\begin{exercise}
Do more computations for $x_i:S^{N-1}_{\mathbb R}\to[-1,1]$.
\end{exercise}

\begin{exercise}
Do more computations for $x_i^2:S^{N-1}_{\mathbb R}\to[0,1]$.
\end{exercise}

\begin{exercise}
Do more computations for $|x_i|:S^{N-1}_{\mathbb R}\to[0,1]$.
\end{exercise}

\begin{exercise}
Do more computations for $z_i:S^{N-1}_{\mathbb C}\to D$.
\end{exercise}

\begin{exercise}
Try as well to use $O_N,U_N$ and Weingarten.
\end{exercise}

As bonus exercise, learn about $SU_2$ and $SO_3$, as much as you can.

\chapter{Discrete measures}

\section*{12a. Negative binomials}

I don't know about you, but personally, after all the mutivariable calculus from the last 3 chapters, I miss discrete mathematics and the Poisson laws. So, this chapter will be about this, discrete laws, from an advanced viewpoint, benefiting from our multivariable experience from the previous 3 chapters. We have at least 4 things to be done:

\bigskip

(1) Higher dimensional Bessel type laws, obtained by summing independent Poisson variables suitably distributed over spheres, or other manifolds $S\subset\mathbb R^N$.

\bigskip

(2) Easiness without limits, meaning relaxing all 3 assumptions in what we have so far of discrete type, namely law of $\chi_t$, for the group $H_N^s$, in the $N\to\infty$ limit. 

\bigskip

(3) Heavy statistics, meaning hypergeometric laws, positive and negative, and their beta versions, obtained by replacing the Bernoulli sampling with a beta sampling.

\bigskip

(4) Fourier type manipulations for complicated real discrete measures, by blowing them up on the unit circle $C\subset\mathbb R^2$, or on other manifolds $S\subset\mathbb R^N$.

\bigskip

Which sounds quite exciting, I mean imagine that there is a World War between continuous and discrete mathematics, and the continuous camp has striked hard, with 3 heavy chapters in a row. As discrete mathematicians we must strike back, with something terribly complicated and advanced, and that can be a mixture of (1,2,3,4).

\bigskip 

In practice now, and trying however to remain a bit philosophers, I would propose to leave aside (1) and (2), which are fairly specialized material. So, have a look here at the literature, or at my specialized Poisson law book, that I have in preparation. And, the plan will be to develop (3), and hopefully have a look at (4) too, towards the end. 

\bigskip

Getting to work, as our starting result, which is something elementary, we have:

\begin{theorem}
When flipping a $p$-biased coin until reaching to $m$ heads,
$$P(s\ {\rm tails})=\binom{s+m-1}{s}(1-p)^sp^m$$
with this being called negative binomial law $c_{mp}$ of parameters $m\in\mathbb N$ and $p\in[0,1]$.
\end{theorem}

\begin{proof}
When computing $P(s\ {\rm tails})$, we certainly have the $(1-p)^sp^m$ factor, coming from the $m$ heads and $s$ tails. As for the multiplicity, this is the binomial coefficient in the statement, coming from choosing the positions of the $s$ tails, among the total of $s+m-1$ attempts, up to the last one, which does not count, because it must be heads. 
\end{proof}

The terminology in the above result comes from the following key fact:

\index{Pascal law}
\index{negative binomial law}
\index{generalized binomial coefficient}

\begin{theorem}
The negative binomial law $c_{mp}$ is best viewed as
$$P(s)=\binom{-m}{s}(p-1)^sp^m$$
with $\binom{-m}{s}$ being a generalized binomial coefficient.
\end{theorem}

\begin{proof}
The above formula comes indeed from the following computation:
\begin{eqnarray*}
\binom{-m}{s}
&=&\frac{-m(-m-1)\ldots(-m-s+1)}{s!}\\
&=&(-1)^s\,\frac{m(m+1)\ldots(m+s-1)}{s!}\\
&=&(-1)^s\binom{s+m-1}{s}
\end{eqnarray*}

Now observe that in this picture, the mass 1 property is obvious, coming from:
\begin{eqnarray*}
\sum_{s\geq0}P(s)
&=&\sum_{s\geq0}\binom{-m}{s}(p-1)^sp^m\\
&=&p^m\sum_{s\geq0}\binom{-m}{s}(p-1)^s\\
&=&p^m[1+(p-1)]^{-m}\\
&=&p^mp^{-m}\\
&=&1
\end{eqnarray*}

Thus, we are led to the conclusions in the statement.
\end{proof}

In analogy now with the formula $b_{np}=b_p^{*n}$ for the usual binomial laws, we have:

\begin{theorem}
The negative binomial laws $c_p=c_{1p}$ of exponent $m=1$ are
$$P(s)=(1-p)^sp$$
with these being called geometric laws, and we have $c_{mp}=c_p^{*m}$, for any $m\in\mathbb N$.
\end{theorem}

\begin{proof}
The first assertion is clear from definitions, and with the name geometric laws coming from the geometric series which produces them, namely:
$$\sum_{s\geq0}P(s)=\sum_{s\geq0}(1-p)^sp=\frac{p}{1-(1-p)}=1$$

As for the second assertion, this comes from the following computation:
\begin{eqnarray*}
c_p^{*m}
&=&\left(\sum_{a_1\geq0}(1-p)^{a_1}p\delta_{a_1}\right)*\ldots*\left(\sum_{a_m\geq0}(1-p)^{a_m}p\delta_{a_m}\right)\\
&=&\sum_{a_1,\ldots,a_m\geq0}(1-p)^{a_1+\ldots+a_m}p^m\delta_{a_1+\ldots+a_m}\\
&=&\sum_{s\geq0}\#\left\{a_1,\ldots,a_m\geq0\Big|a_1+\ldots+a_m=s\right\}(1-p)^sp^m\delta_s\\
&=&\sum_{s\geq0}\#\left\{b_1,\ldots,b_m\geq1\Big|b_1+\ldots+b_m=s+m\right\}(1-p)^sp^m\delta_s\\
&=&\sum_{s\geq0}\binom{s+m-1}{m-1}(1-p)^sp^m\delta_s\\
&=&\sum_{s\geq0}\binom{s+m-1}{s}(1-p)^sp^m\delta_s\\
&=&c_{mp}
\end{eqnarray*}

To be more precise here, the count at the end comes from the fact that, in order to partition $s+m$, best thought as a sequence of $s+m$ boxes, as $s+m=b_1+\ldots+b_m$, we have to insert $m-1$ markers between these boxes, at different locations, at the $s+m-1$ positions available. Thus, we are led to the conclusion in the statement.
\end{proof}

We have as well a negative Poisson Limit Theorem, as follows:

\index{PLT}
\index{Poisson Limit Theorem}
\index{Bernoulli laws}
\index{Poisson limit}
\index{negative PLT}
\index{NPLT}
\index{geometric law}

\begin{theorem}[negative PLT]
We have the following convergence, in moments,
$$\left(\sum_{k\geq0}\left(\frac{t}{m+t}\right)^k\frac{m}{m+t}\,\delta_k\right)^{*m}\to p_t$$
for any $t>0$. Equivalently, $c_{mp}\to p_t$, when $p=m/(m+t)$ with $t>0$ fixed.
\end{theorem}

\begin{proof}
Let us denote by $\nu_m$ the geometric law under the convolution sign:
$$\nu_m=\sum_{k\geq0}\left(\frac{t}{m+t}\right)^k\frac{m}{m+t}\,\delta_k$$

We have the following computation, for the Fourier transform of the limit: 
\begin{eqnarray*}
F_{\delta_k}(x)=e^{ikx}
&\implies&F_{\nu_m}(x)=\sum_{k\geq0}\left(\frac{t}{m+t}\right)^k\frac{m}{m+t}\,e^{ikx}\\
&\implies&F_{\nu_m}(x)=\frac{m}{m+t-te^{ix}}\\
&\implies&F_{\nu_m^{*m}}(x)=\left(\frac{m}{m+t-te^{ix}}\right)^m\\
&\implies&F_{\nu_m^{*m}}(x)=\left(1+\frac{t(e^{ix}-1)}{m+t-te^{ix}}\right)^m\\
&\implies&F(x)=\exp\left(t(e^{ix}-1)\right)
\end{eqnarray*}

Thus, we obtain indeed the Fourier transform of $p_t$, as desired.
\end{proof}

Again inspired from what we did before for the usual binomial laws, we have:

\begin{theorem}
The moments of $c_{mp}$ are given by the following formula,
$$M_k=\sum_{\pi\in P(k)}\frac{(m+|\pi|-1)!}{(m-1)!}\left(\frac{1-p}{p}\right)^{|\pi|}$$
with $|.|$ standing as usual for the number of blocks of partitions, and 
$$E=\frac{m(1-p)}{p}\ ,\  V=\frac{m(1-p)}{p^2}\ ,\ \gamma=\frac{2-p}{\sqrt{m(1-p)}}\ ,\ \kappa=3+\frac{6}{m}+\frac{p^2}{m(1-p)}$$
are the corresponding mean, variance, skewness and kurtosis.
\end{theorem}

\begin{proof}
Regarding the mean, this can be computed as follows:
\begin{eqnarray*}
E
&=&\sum_{s\geq1}s\binom{-m}{s}(p-1)^sp^m\\
&=&\sum_{s\geq1}-m\binom{-m-1}{s-1}(p-1)^sp^m\\
&=&-mp^m(p-1)\sum_{s\geq1}\binom{-m-1}{s-1}(p-1)^{s-1}\\
&=&mp^m(1-p)[1+(p-1)]^{-m-1}\\
&=&mp^m(1-p)p^{-m-1}\\
&=&\frac{m(1-p)}{p}
\end{eqnarray*}

With a similar trick, we can compute the difference $M_2-M_1$, as follows:
\begin{eqnarray*}
M_2-M_1
&=&\sum_{s\geq2}s(s-1)\binom{-m}{s}(p-1)^sp^m\\
&=&\sum_{s\geq2}m(m+1)\binom{-m-2}{s-2}(p-1)^sp^m\\
&=&m(m+1)p^m(p-1)^2\sum_{s\geq2}\binom{-m-2}{s-2}(p-1)^{s-2}\\
&=&m(m+1)p^m(1-p)^2[1+(p-1)]^{-m-2}\\
&=&m(m+1)p^m(1-p)^2p^{-m-2}\\
&=&\frac{m(m+1)(1-p)^2}{p^2}
\end{eqnarray*}

We conclude that the second moment is given by the following formula:
$$M_2=\frac{m(1-p)(1+m-mp)}{p^2}$$

Thus, we are led to the variance formula in the statement, namely:
$$V=\frac{m(1-p)(1+m-mp)}{p^2}-\frac{m^2(1-p)^2}{p^2}=\frac{m(1-p)}{p^2}$$

In order to compute now the higher moments, with the same trick, we can use the following formula, that we know well since chapter 1, from our study there of $b_{np}$:
$$s^k=\sum_{\pi\in P(k)}\frac{s!}{(s-|\pi|)!}$$

Indeed, with this as input for our moment computations, we are led to:
$$M_k=\sum_{\pi\in P(k)}\frac{(m+|\pi|-1)!}{(m-1)!}\left(\frac{1-p}{p}\right)^{|\pi|}$$

Getting now to $M_3,M_4$, these do not look very good, but their central versions do:
$$M_3'=\frac{m(1-p)(2-p)}{p^3}$$
$$M_4'=\frac{m(1-p)((3m+6)(1-p)+p^2)}{p^4}$$

And with this, we are led to the formulae of $\gamma,\kappa$ in the statement, as desired.
\end{proof}

\section*{12b. Hypergeometric laws}

Let us go back now to the usual binomial laws, studied in chapter 1. As a variation of that construction, we can talk about hypergeometric laws, as follows:

\index{hypergeometric law}
\index{Vandermonde formula}

\begin{theorem}
Given a population of size $N\in\mathbb N$, with $m\in\{0,\ldots,N\}$ of these objects having a certain feature, the probability of having $s$ successes, when performing $p\in\{0,\ldots,N\}$ draws without replacement, and looking for that feature, is
$$P(s)=\frac{\binom{m}{s}\binom{N-m}{p-s}}{\binom{N}{p}}$$
with this being called hypergeometric law of parameters $(N,m,p)$. When doing the same thing, but with replacement, we obtain a binomial law of parameter $m/N$.
\end{theorem}

\begin{proof}
Many things can be said here, the idea being as follows:

\medskip

(1) To start with, the main assertion is clear, because we have a total of $\binom{N}{p}$ possibilities for our $p$ draws without replacement, and among these draws, those amounting to $s$ successes come by multiplying $\binom{m}{s}$, standing for the $s$ successes, and $\binom{N-m}{p-s}$, standing for the $p-s$ fails. Thus, we are led to the formula of $P(s)$ in the statement.

\medskip

(2) Regarding now the fact that we have indeed a probability measure, this is something which is plainly clear, because when performing our $p$ draws, something must happen, I mean we have a well-defined number of successes $s\in\mathbb N$, and $\sum_sP(s)=1$.

\medskip

(3) However, let us comment more on this. To start with, we must have:
$$0\leq s\leq m\quad,\quad 0\leq p-s\leq N-m$$

In other words, the number of successes $s\in\mathbb N$, where $P(s)>0$, is subject to:
$$s\in\{\max(0,p+m-N),\ldots,\min(m,p)\}$$

(4) Getting now to the mass 1 formula $\sum_sP(s)=1$, this looks in practice as follows, with the bounds for the summing parameter $s$ being those found above:
$$\sum_s\binom{m}{s}\binom{N-m}{p-s}=\binom{N}{p}$$

But this comes from the Vandermonde formula, which is as follows, coming from $(1+x)^m(1+x)^n=(1+x)^{m+n}$, by looking at the coefficient of $x^p$, on both sides:
$$\sum_{s=0}^p\binom{m}{s}\binom{n}{p-s}=\binom{m+n}{p}$$

(5) Finally, in what regards the last assertion of the theorem, assume that we are doing our $p$ draws, but this time with replacement. In this case, for having $s$ successes, we have $\binom{p}{s}$ choices for the occurences of these successes, and then $(m/N)^s$ probability for the successes, and $(1-m/N)^{n-s}$ probability for the fails, so we obtain:
$$P(s)=\binom{p}{s}\left(\frac{m}{N}\right)^s\left(1-\frac{m}{N}\right)^{n-s}$$

But this is a binomial law of parameter $m/N$, as stated.
\end{proof}

Still at the theoretical level, as an interesting observation, we have:

\begin{proposition}
The hypergeometric law of parameters $(N,m,p)$, given by
$$P(s)=\frac{m!p!(N-m)!(N-p)!}{s!(m-s)!(p-s)!N!(N-m-p+s)!}$$
is symmetric in $m,p$, so is equal to the hypergeometric law of parameters $(N,p,m)$.
\end{proposition}

\begin{proof}
This is indeed something self-explanatory, coming from definitions.
\end{proof}

At a more advanced level, we can use some algebra. Consider the standard matrix coordinates over the symmetric group $S_N\subset O_N$, which form a matrix, as follows:
$$u=\begin{pmatrix}
u_{11}&\ldots&u_{1p}&\ldots&\ldots&u_{1N}\\
\vdots&\bigstar&\vdots&&&\vdots\\
u_{m1}&\ldots&u_{mp}&\ldots&\ldots&u_{mN}\\
\vdots&&\vdots&&&\vdots\\
\vdots&&\vdots&&&\vdots\\
u_{N1}&\ldots&u_{Np}&\ldots&\ldots&u_{NN}
\end{pmatrix}$$

The point now is that the sum of the entries of the rectangular matrix $\bigstar$ appearing in the upper left corner is hypergeometric, as shown by the following result:

\index{hypergeometric variable}

\begin{theorem}
The following variable over the symmetric group $S_N\subset O_N$,
$$S_{mp}=\sum_{i=1}^m\sum_{j=1}^pu_{ij}$$
is hypergeometric, having parameters $(N,m,p)$.
\end{theorem}

\begin{proof}
We know that the coordinates of the symmetric group, viewed as group of permutation matrices, $S_N\subset O_N$, are given by the following formula:
$$u_{ij}=\chi\left(\sigma\in S_N\Big|\sigma(j)=i\right)$$

Thus, the variable in the statement is given by the following formula:
$$S_{mp}=\sum_{i=1}^m\sum_{j=1}^p\chi\left(\sigma\in S_N\Big|\sigma(j)=i\right)$$

Let us compute now the law of this variable. We have the following formula:
\begin{eqnarray*}
P(S_{mp}=s)
&=&\frac{1}{N!}\,\#\left(\sigma\in S_N\Big|S_{mp}(\sigma)=s\right)\\
&=&\frac{1}{N!}\,\#\left(\sigma\in S_N\Big|\#\left\{i\leq m,j\leq p\big|\sigma(j)=i\right\}=s\right)\\
&=&\frac{1}{N!}\,\#\left(\sigma\in S_N\Big|\#\left\{j\leq p\big|\sigma(j)\leq m\right\}=s\right)
\end{eqnarray*}

Now let us try to count the permutations on the right. Such a permutation $\sigma\in S_N$ comes from 5 operations involved, as follows:

\medskip

(1) We first have to pick a subset $X\subset\{1,\ldots,p\}$, with $|X|=s$.

\medskip

(2) We also have to pick a subset $Y\subset\{1,\ldots,m\}$, with $|Y|=s$.

\medskip

(3) Then, we have to bijectively map $X\to Y$.

\medskip

(4) Next, we have to map injectively $\{1,\ldots,p\}-X\to\{m+1,\ldots,N\}$.

\medskip

(5) And finally, we have to bijectively map $\{p+1,\ldots,N\}$ to what is left.

\medskip

Regarding now the precise numbers of choices for these 5 operations involved, these are easy to compute, as follows:

\medskip

(1) Here we have $\binom{p}{s}$ choices.

\medskip

(2) Here we have $\binom{m}{s}$ choices.

\medskip

(3) Here we have $s!$ choices.

\medskip

(4) Here we have $\frac{(N-m)!}{(N-m-p+s)!}$ choices.

\medskip

(5) Here we have $(N-p)!$ choices.

\medskip

We can now finish our probability computation started above, as follows:
\begin{eqnarray*}
P(S_{mp}=s)
&=&\frac{1}{N!}\binom{p}{s}\binom{m}{s}s!\,\frac{(N-m)!}{(N-m-p+s)!}(N-p)!\\
&=&\frac{1}{N!}\cdot\frac{p!}{s!(p-s)!}\cdot\frac{m!}{s!(m-s)!}\cdot s!\cdot\frac{(N-m)!}{(N-m-p+s)!}(N-p)!\\
&=&\frac{p!m!(N-m)!(N-p)!}{N!s!(p-s)!(m-s)!(N-m-p+s)!}\\
&=&\frac{m!}{s!(m-s)!}\cdot\frac{(N-m)!}{(p-s)!(N-m-p+s)!}\cdot\frac{p!(N-p)!}{N!}\\
&=&\frac{\binom{m}{s}\binom{N-m}{p-s}}{\binom{N}{p}}
\end{eqnarray*}

We are therefore led to the conclusion in the statement.
\end{proof}

As a first observation, the above trivializes the symmetry observation from Proposition 12.7. As another success of our group theory approach, by some magic, we have:

\index{moments of hypergeometric law}

\begin{theorem}
For the hypergeometric law of parameters $(N,m,p)$ we have
$$M_k=\sum_{\pi\in P(k)}\frac{m!}{(m-|\pi|)!}\cdot\frac{p!}{(p-|\pi|)!}\cdot\frac{(N-|\pi|)!}{N!}$$
where $|.|$ denotes as usual the number of blocks of partitions, and
$$E=\frac{mp}{N}\quad,\quad V=\frac{mp(N-m)(N-p)}{N^2(N-1)}\quad,\quad 
\gamma=\frac{(N-2m)(N-2p)\sqrt{N-1}}{(N-2)\sqrt{mp(N-m)(N-p)}}$$
$$\kappa=3+\frac{\begin{pmatrix}N^2(N-1)[N(N+1)-6m(N-m)-6p(N-p)]\\
+6mp(N-m)(N-p)(5N-6)\end{pmatrix}}{mp(N-m)(N-p)(N-2)(N-3)}$$
and the corresponding mean, variance, skewness and kurtosis.
\end{theorem}

\begin{proof}
We have the following computation, using technology from chapter 4:
\begin{eqnarray*}
M_k
&=&\int_{S_N}S_{mp}^k\\
&=&\int_{S_N}\left(\sum_{i=1}^m\sum_{j=1}^pu_{ij}\right)^k\\
&=&\int_{S_N}\sum_{i_1=1}^m\ldots\sum_{i_k=1}^m\sum_{j_1=1}^p\ldots\sum_{j_k=1}^pu_{i_1j_1}\ldots u_{i_kj_k}\\
&=&\sum_{i_1=1}^m\ldots\sum_{i_k=1}^m\sum_{j_1=1}^p\ldots\sum_{j_k=1}^p\int_{S_N}u_{i_1j_1}\ldots u_{i_kj_k}\\
&=&\sum_{i_1=1}^m\ldots\sum_{i_k=1}^m\sum_{j_1=1}^p\ldots\sum_{j_k=1}^p\delta_{\ker i,\ker j}\frac{(N-|\ker i|)!}{N!}\\
&=&\sum_{\pi\in P(k)}\sum_{i_1=1}^m\ldots\sum_{i_k=1}^m\sum_{j_1=1}^p\ldots\sum_{j_k=1}^p\delta_{\ker i,\ker j,\pi}\frac{(N-|\pi|)!}{N!}\\
&=&\sum_{\pi\in P(k)}\left(\sum_{i_1=1}^m\ldots\sum_{i_k=1}^m\delta_{\ker i,\pi}\right)
\left(\sum_{j_1=1}^p\ldots\sum_{j_k=1}^p\delta_{\ker j,\pi}\right)\frac{(N-|\pi|)!}{N!}\\
&=&\sum_{\pi\in P(k)}\frac{m!}{(m-|\pi|)!}\cdot\frac{p!}{(p-|\pi|)!}\cdot\frac{(N-|\pi|)!}{N!}
\end{eqnarray*}

As for the numerics, I will leave them to you as an exercise, enjoy.
\end{proof}

Next, we can talk as well about negative hypergeometric laws, as follows:

\index{negative hypergeometric law}
\index{negative binomial law}

\begin{theorem}
Given a population of size $N\in\mathbb N$, with $m\in\{0,\ldots,N\}$ of these objects having a certain feature, the probability of having $s$ successes, when performing draws without replacement, and stopping after $r\in\{0,\ldots,N-m\}$ failures, is
$$P(s)=\frac{\binom{s+r-1}{s}\binom{N-r-s}{m-s}}{\binom{N}{m}}$$
with this being called negative hypergeometric law of parameters $(N,m,r)$. When doing this with replacement, we obtain a negative binomial law of parameter $m/N$.
\end{theorem}

\begin{proof}
Many things can be said here, the idea being as follows:

\medskip

(1) Observe first the similarity with the definition of the usual hypergeometric laws, from Theorem 12.6. To be more precise, the common part is that, in both these cases, we have a population of size $N\in\mathbb N$, with $m\in\{0,\ldots,N\}$ of these objects having a certain feature, and we are looking for the probability of having $s$ successes, when performing draws with replacement. The difference comes from the stopping convention, for the hypergeometric laws this being after $p\in\{0,\ldots,N\}$ draws, and for the negative hypergeometric laws this being after $r\in\{0,\ldots,N-m\}$ failures.

\medskip

(2) In practice, based on this, we can develop the theory of negative hypergeometric laws by using the previously developed theory of the usual hypergeometric laws. Indeed, to start with, the probability $P(s)$ that we are interested in comes by multiplying a probability coming from a usual hypergeometric law, corresponding to the $m$ succeses in the first $m+r-1$ draws, and a basic fraction, corresponding to the failure required at the $(m+r)$-th draw. In practice, this gives the following formula, as desired:
\begin{eqnarray*}
P(s)
&=&\frac{\binom{m}{s}\binom{N-m}{s+r-1-s}}{\binom{N}{s+r-1}}\cdot\frac{N-m-(r-1)}{N-(s+r-1)}\\
&=&\frac{\binom{m}{s}\binom{N-m}{r-1}}{\binom{N}{s+r-1}}\cdot\frac{N-m-r+1}{N-s-r+1}\\
&=&\frac{m!(N-m)!(s+r-1)!(N-s-r+1)!(N-m-r+1)}{s!(m-s)!(r-1)!N!(N-m-r+1)!(N-s-r+1)}\\
&=&\frac{m!(N-m)!(s+r-1)!(N-s-r)!}{s!(m-s)!(r-1)!N!(N-m-r)!}\\
&=&\frac{(s+r-1)!(N-r-s)!(N-m)!m!}{s!(r-1)!(m-s)!(N-r-m)!N!}\\
&=&\frac{\binom{s+r-1}{s}\binom{N-r-s}{m-s}}{\binom{N}{m}}
\end{eqnarray*}

(3) Regarding now the fact that we have indeed a probability measure, $\sum_sP(s)=1$, this is something which is plainly clear from definitions, but which can be deduced as well from the Vandermonde formula, along with the identity $\binom{n}{k}=(-1)^k\binom{k-n-1}{k}$:
\begin{eqnarray*}
\sum_s\binom{s+r-1}{s}\binom{N-r-s}{m-s}
&=&\sum_s(-1)^s\binom{-r}{s}(-1)^{m-s}\binom{m+r-N-1}{m-s}\\
&=&(-1)^m\sum_s\binom{-r}{s}\binom{m+r-N-1}{m-s}\\
&=&(-1)^m\binom{m-N-1}{m}\\
&=&\binom{N}{m}
\end{eqnarray*}

(4) Finally, in what regards the second assertion of the theorem, assume that we are doing our draws, but this time with replacement. In this case, for having $s$ successes, we have $\binom{s+r-1}{s}$ choices for the occurences of these successes, and then $(m/N)^s$ probability for the successes, and $(1-m/N)^r$ probability for the fails, so we obtain:
$$P(s)=\binom{s+r-1}{s}\left(\frac{m}{N}\right)^s\left(1-\frac{m}{N}\right)^r$$

But this is a negative binomial law of parameter $m/N$, as stated.
\end{proof}

As a main result now regarding the negative hypergeometric laws, we have:

\begin{theorem}
For the negative hypergeometric law of parameters $(N,m,r)$,
$$M_k=\sum_{\pi\in P(k)}\frac{m!}{(m-|\pi|)!}\cdot\frac{(r+|\pi|-1)!}{(r-1)!}\cdot\frac{(N-m)!}{(N-m+|\pi|)!}$$
where $|.|$ denotes as usual the number of blocks of partitions, we have
$$E=\frac{mr}{N-m+1}\quad,\quad V=\frac{mr(N+1)(N-m-r+1)}{(N-m+1)^2(N-m+2)}$$
$$\gamma=\frac{(N+m+1)(N-m-2r+1)}{N-m+3}\sqrt{\frac{N-m+2}{mr(N+1)(N-m-r+1)}}$$
and with the notations $c=N-m+1$ and $d=r(N-m-r+1)$, 
$$\kappa=\frac{(c+1)(c^2[c^2+(6m-1)c+6m^2]+3d[(m-2)c^2+m(m-6)c-6m^2])}{md(c+2)(c+3)(c+m)}$$
is the corresponding kurtosis.
\end{theorem}

\begin{proof}
This is something quite standard, the idea being that the formula in the statement follows by recurrence, by using our standard trick, namely: 
$$s^k=\sum_{\pi\in P(k)}\frac{s!}{(s-|\pi|)!}$$

We will leave the computations here, as well as the remaining discussion, regarding the mean, variance, skewness and kurtosis, as an instructive exercise. Enjoy.
\end{proof}

\section*{12c. Beta binomials}

With the above, done with our series of generalizations? You must be kidding. Indeed, in advanced statistics it is sometimes useful to have some further generalizations of all the above, with the usual Bernoulli sampling replaced by a beta sampling. And with this being  something quite interesting, for other mathematical purposes as well.

\bigskip

In order to discuss this, we must first have a closer look at the beta laws. We already talked about these, vaguely, in chapter 3, and then we met them, at certain special values of the parameters, in chapter 11, disguised there as squared hyperspherical laws. In order to have a systematic look at this, we will need a piece of fine mathematics, namely:

\begin{theorem}
We have the following formula, valid for any $a,b>0$,
$$\int_0^1x^{a-1}(1-x)^{b-1}dx=\frac{\Gamma(a)\Gamma(b)}{\Gamma(a+b)}$$
with this quantity being denoted $B(a,b)$, and called beta function. 
\end{theorem}

\begin{proof}
This is something quite tricky, the idea being as follows;

\medskip

(1) As a first comment, you might wonder why calling B beta, but the thing is, the letter B is originally Greek, upper-case $\beta$. And with this causing various troubles, because when needing distinct upper-case $b,\beta$ letters, you are tempted to use $B,\mathfrak B$. With by the way mea culpa for having done such things in the past. We won't do them here.

\medskip

(2) As a second comment, in the case where the exponents are integers, $a,b\in\mathbb N$, the formula in the statement takes the following form, in terms of usual factorials:
$$\int_0^1x^{a-1}(1-x)^{b-1}dx=\frac{(a-1)!(b-1)!}{(a+b-1)!}$$

But with $x=\cos^2t$ this is exactly the Wallis 2 formula from chapter 11, at odd values $p=2a-1$ and $q=2b-1$ of the exponents there. In fact, the Wallis 2 formula corresponds precisely to the formula in the statement at $a,b\in\mathbb N/2$, and this is why in chapter 11 we managed to get away with Wallis, without using the beta function.

\medskip

(3) Getting now to the general case, let us attempt to compute $\Gamma(a)\Gamma(b)$. By definition of gamma, this amounts in computing a certain integral over $\mathbb R^2$, as follows:
\begin{eqnarray*}
\Gamma(a)\Gamma(b)
&=&\int_0^\infty x^{a-1}e^{-x}dx\int_0^\infty y^{b-1}e^{-y}dy\\
&=&\int_0^\infty\int_0^\infty x^{a-1}y^{b-1}e^{-x-y}dxdy
\end{eqnarray*}

(4) Now comes the trick. Let us perform the following change of variables:
$$\begin{cases}x=st\\ y=s(1-t)\end{cases}\iff
\begin{cases}s=x+y\\ t=x/(x+y)\end{cases}$$ 

We have then $dxdy=sdsdt$, coming from the following Jacobian computation:
$$\frac{dxdy}{dsdt}
=\begin{vmatrix}dx/ds&dx/dt\\ \\ dy/ds&dy/dt\end{vmatrix}
=\begin{vmatrix}t&s\\ \\1-t&-s\end{vmatrix}=-s$$

(5) We conclude that our integral over $\mathbb R^2$ above takes the following form:
\begin{eqnarray*}
\Gamma(a)\Gamma(b)
&=&\int_0^1\int_0^\infty(st)^{a-1}(s(1-t))^{b-1}e^{-s}sdsdt\\
&=&\int_0^1s^{a+b-1}e^{-s}ds\int_0^\infty t^{a-1}(1-t)^{b-1}dt\\
&=&\Gamma(a+b)\int_0^\infty t^{a-1}(1-t)^{b-1}dt
\end{eqnarray*}

Thus, we are led to the conclusion in the statement.
\end{proof}

Before getting into probability, let us study a bit more the beta function. We know that for integer variables $a,b\in\mathbb N$, we have the following formula:
$$B(a,b)=\frac{(a-1)!(b-1)!}{(a+b-1)!}$$

Thus, the beta function is some sort of generalized binomial coefficient, written upside down, and coming with an extra factor, or a factor missing, upon taste. We have:

\index{Pascal formula}
\index{Vandermonde formula}

\begin{theorem}
We have the following formula,
$$\sum_{s=0}^n\binom{n}{s}B(s+a,n-s+b)=B(a,b)$$
generalizing the usual Vandermonde formula for binomials.
\end{theorem}

\begin{proof}
This is very standard, as for the usual Vandermonde formula, as follows:

\medskip

(1) At $n=0$ the formula to be proved is $B(a,b)=B(a,b)$, true.

\medskip

(2) At $n=1$ the formula, called Pascal formula, can be proved as follows:
\begin{eqnarray*}
B(a,b+1)+B(a+1,b)
&=&\frac{\Gamma(a)\Gamma(b+1)}{\Gamma(a+b+1)}+\frac{\Gamma(a+1)\Gamma(b)}{\Gamma(a+b+1)}\\
&=&\frac{b\Gamma(a)\Gamma(b)}{(a+b)\Gamma(a+b)}+\frac{a\Gamma(a)\Gamma(b)}{(a+b)\Gamma(a+b)}\\
&=&\frac{\Gamma(a)\Gamma(b)}{\Gamma(a+b)}\\
&=&B(a,b)
\end{eqnarray*}

(3) At $n=2$ now, we can use the above Pascal formula three times, as follows:
\begin{eqnarray*}
&&B(a,b+2)+2B(a+1,b+1)+B(a+2,b)\\
&=&[B(a,b+2)+B(a+1,b+1)]+[B(a+1,b+1)+B(a+2,b)]\\
&=&B(a,b+1)+B(a+1,b)\\
&=&B(a,b)
\end{eqnarray*}

(4) And so on, with the general case following by iterating Pascal, as above.
\end{proof}

Time, eventually, to do some probability? Based on what we have, we can now formulate the following key result, which is something quite far-reaching:

\index{beta distribution}

\begin{theorem}
The beta distribution of parameters $a,b>0$, namely
$$\nu_{ab}=\frac{x^{a-1}(1-x)^{b-1}}{B(a,b)}\,dx$$
on $[0,1]$, with $B(a,b)=\Gamma(a)\Gamma(b)/\Gamma(a+b)$ being the beta function, has moments
$$M_k=\frac{a(a+1)\ldots(a+k-1)}{(a+b)(a+b+1)\ldots(a+b+k-1)}$$
and in particular, the mean and variance are given by the following formulae,
$$E=\frac{a}{a+b}\quad,\quad V=\frac{ab}{(a+b)^2(a+b+1)}$$
and the skewness and kurtosis are given by the following formulae,
$$\gamma=\frac{2(b-a)\sqrt{a+b+1}}{(a+b+2)\sqrt{ab}}\quad,\quad
\kappa=3+\frac{6[(a-b)^2(a+b+1)-ab(a+b+2)]}{ab(a+b+2)(a+b+3)}$$
and with all this generalizing many previous formulae, from this book.
\end{theorem}

\begin{proof}
The moments can be indeed computed as follows, using Theorem 12.12:
\begin{eqnarray*}
M_k
&=&\frac{1}{B(a,b)}\int_0^1x^{a+k-1}(1-x)^{b-1}dx\\
&=&\frac{B(a+k,b)}{B(a,b)}\\
&=&\frac{\Gamma(a+k)\Gamma(b)}{\Gamma(a+b+k)}\cdot\frac{\Gamma(a+b)}{\Gamma(a)\Gamma(b)}\\
&=&\frac{\Gamma(a+k)}{\Gamma(a)}\cdot\frac{\Gamma(a+b)}{\Gamma(a+b+k)}\\
&=&\frac{a(a+1)\ldots(a+k-1)}{(a+b)(a+b+1)\ldots(a+b+k-1)}
\end{eqnarray*}

In particular, we have the formula of $E$. Getting now to the variance, this is:
$$V=\frac{a(a+1)}{(a+b)(a+b+1)}-\left(\frac{a}{a+b}\right)^2
=\frac{ab}{(a+b)^2(a+b+1)}$$

As for the computations of the skewness and kurtosis, these are similar.
\end{proof}

Getting back now to our usual business in this chapter, namely discrete laws, we can talk about beta binomial distributions, constructed in the following way:

\index{beta binomial law}
\index{beta function}

\begin{definition}
The beta binomial distributions are the discrete laws given by
$$P(s)=\binom{n}{s}\frac{B(s+a,n-s+b)}{B(a,b)}$$
with $B(a,b)=\Gamma(a)\Gamma(b)/\Gamma(a+b)$ being as usual the beta function.
\end{definition}

As a first job, let us verify that we have indeed probability measures, of mass 1. But this comes indeed from the Vandermonde formula in Theorem 12.13, as follows:
$$\sum_{s=0}^nP(s)
=\frac{1}{B(a,b)}\sum_{s=0}^n\binom{n}{s}B(s+a,n-s+b)
=1$$

Next, we have the following fact, which is the key to understanding our laws:

\begin{theorem}
At integer values of the parameters, $a,b\in\mathbb N$, we have
$$P(s)=\frac{\binom{s+a-1}{s}\binom{n-s+b-1}{n-s}}{\binom{n+a+b-1}{n}}$$
which is a negative hypergeometric law, of parameters $(N,m,r)=(n+a+b-1,n,a)$.
\end{theorem}

\begin{proof}
At integer values of the parameters $a,b$, the beta function is given by:
$$B(a,b)=\frac{\Gamma(a)\Gamma(b)}{\Gamma(a+b)}=\frac{(a-1)!(b-1)!}{(a+b-1)!}$$

By using this, the formula of the beta binomial distribution is as follows:
\begin{eqnarray*}
P(s)
&=&\binom{n}{s}\frac{B(s+a,n-s+b)}{B(a,b)}\\
&=&\frac{n!}{s!(n-s)!}\cdot\frac{(s+a-1)!(n-s+b-1)!}{(n+a+b-1)!}\cdot\frac{(a+b-1)!}{(a-1)!(b-1)!}\\
&=&\frac{n!(a+b-1)!}{(n+a+b-1)!}\cdot\frac{(s+a-1)!}{s!(a-1)!}\cdot\frac{(n-s+b-1)!}{(n-s)!(b-1)!}\\
&=&\frac{\binom{s+a-1}{s}\binom{n-s+b-1}{n-s}}{\binom{n+a+b-1}{n}}
\end{eqnarray*}

Thus, we are led to the conclusion in the statement.
\end{proof}

We can see now the interest in Definition 12.15, the point being that, in certain situations where the negative hypergeometric laws do not produce accurate results, we can use the beta binomial laws. For more here, we refer to a solid statistics book. 

\bigskip

In what follows we will be mostly interested in the abstract mathematics of the beta binomial laws. In order to deal with the moment problem, let us formulate:

\begin{definition}
The factorial moments of $f:X\to\mathbb R$ are the numbers
$$\widetilde{M}_k=E[f(f-1)\ldots(f-k+1)]$$
with the convention $\widetilde{M}_0=1$.
\end{definition}

This definition is something quite subtle, the idea is that the above normalization ``kills the underlying partitions''. Indeed, recall the following formula, that we know well for numbers, and which must therefore hold as well for random variables:
$$f^k=\sum_{\pi\in P(k)}f(f-1)\ldots(f-|\pi|+1)$$

By applying the expectation on both sides, we obtain the following formula:
$$M_k=\sum_{\pi\in P(k)}\widetilde{M}_{|\pi|}$$

As an illustration for this conversion formula, and for the power of Definition 12.17, in general, for the Poisson laws we have the following remarkably simple formula:
$$\widetilde{M_k}(p_t)=t^k$$

Getting now to the beta binomial laws, the result for them is as follows:

\index{factorial moments}

\begin{theorem}
The factorial moments of the beta binomial distribution are
$$\widetilde{M}_k=\frac{n!}{(n-k)!}\cdot\frac{B(a+k,b)}{B(a,b)}$$
with $B$ being as usual the beta function. The mean and variance are
$$E=\frac{na}{a+b}\quad,\quad V=\frac{nab(n+a+b)}{(a+b)^2(a+b+1)}$$
the skewness is given by the following formula,
$$\gamma=\frac{(2n+a+b)(b-a)}{a+b+2}\sqrt{\frac{a+b+1}{nab(n+a+b)}}$$
and with the notations $c=a+b$ and $d=ab$, 
$$\kappa=\frac{(c+1)(c^2[c^2+(6n-1)c+6n^2]+3d[(n-2)c^2+n(n-6)c-6n^2])}{nd(c+2)(c+3)(c+n)}$$
is the corresponding kurtosis.
\end{theorem}

\begin{proof}
To start with, let us record the following more digest formulation of the moment formula in the statement, in terms of the gamma function only:
\begin{eqnarray*}
\widetilde{M}_k
&=&\frac{n!}{(n-k)!}\cdot\frac{B(a+k,b)}{B(a,b)}\\
&=&\frac{n!}{(n-k)!}\cdot\frac{\Gamma(a+k)\Gamma(b)}{\Gamma(a+b+k)}\cdot\frac{\Gamma(a+b)}{\Gamma(a)\Gamma(b)}\\
&=&\frac{n!}{(n-k)!}\cdot\frac{\Gamma(a+k)}{\Gamma(a)}\cdot\frac{\Gamma(a+b)}{\Gamma(a+b+k)}
\end{eqnarray*}

Equivalently, by using the formula $\Gamma(s+1)=s\Gamma(s)$ for both fractions, we have the following formula, directly in terms of the parameters $a,b>0$:
$$\widetilde{M_k}=\frac{n(n-1)\ldots(n-k+1)a(a+1)\ldots(a+k-1)}{(a+b)(a+b+1)\ldots(a+b+k-1)}$$

As for the proof, this is something quite standard, by using the same tricks as for the negative hypergeometric laws. We will leave this as an exercise.
\end{proof}

And with this, done? Not yet, because we have one more thing to talk about:

\index{negative beta binomial law}

\begin{definition}
The negative beta binomial law of parameters $r,a,b>0$ is
$$P(s)=\frac{\Gamma(s+b)}{s!\Gamma(b)}\cdot\frac{B(s+r,a+b)}{B(r,a)}$$
with $B=\Gamma(a)\Gamma(b)/\Gamma(a+b)$ being as usual the beta function.
\end{definition}

Obviously, this is something quite tricky. As a first observation, by writing everything in terms of the gamma function, and rearranging the terms, we have as well:
$$P(s)=\frac{1}{s}\cdot\frac{B(r+a,s+b)}{B(a,b)B(r,s)}$$

In practice, the simplest case is $r=1$, and we have here the following result:

\begin{proposition}
The negative beta binomial distribution at $r=1$ is given by
$$P(s)=\frac{B(a+1,b+s)}{B(a,b)}$$
and with this being called beta geometric law of parameters $a,b>0$.
\end{proposition}

\begin{proof}
This comes from our beta-only formula above, which reads:
$$P(s)=\frac{1}{s}\cdot\frac{B(a+1,s+b)}{B(a,b)B(1,s)}=\frac{\Gamma(s+1)}{s\Gamma(s)}\cdot\frac{B(a+1,b+s)}{B(a,b)}$$

Thus, we are led to the formula in the statement.
\end{proof}

Generally speaking, the negative beta binomial distributions are best thought as coming in 3 stages, namely the beta geometric laws above, corresponding to the case $r=1$, then their generalizations with $r\in\mathbb N$, and then their further generalizations with $r>0$. Many other things can be here, and as before for the usual beta binomial distributions, for more on the phenomenology of these laws, we refer to a solid statistics book.

\bigskip

Getting now to the mathematics of these laws, as a central result, we have:

\begin{theorem}
The factorial moments of the negative beta binomial law are
$$\widetilde{M}_k=\frac{\Gamma(r+k)\Gamma(a-k)\Gamma(b+k)}{\Gamma(r)\Gamma(a)\Gamma(b)}$$
with $r,a,b>0$ being as usual the parameters.
\end{theorem}

\begin{proof}
This is obviously something quite complicated, that we will not attempt to prove here. This being said, a few comments. As a first observation, by using the identity $\Gamma(s+1)=s\Gamma(s)$, the formula in the statement can be written as follows:
$$\widetilde{M}_k=\frac{(r+k-1)(r+k-2)\ldots r(b+k-1)(b+k-2)\ldots b}{(a-1)(a-2)\ldots(a-k)}$$

In practice, at $k=0$ this tells us that we have $\widetilde{M}_0=1$, with this showing that our measures are indeed of mass one, as they should. Next, at $k=1,2$ we get:
$$\widetilde{M}_1=\frac{rb}{a-1}\quad,\quad \widetilde{M}_2=\frac{(r+1)r(b+1)b}{(a-1)(a-2)}$$

And so on, the idea being that the formula of $\widetilde{M}_k$ is something simple and useful.
\end{proof}

Regarding the standard parameters $(E,V,\gamma,\kappa)$ of our distribution, we have:

\begin{theorem}
The mean and variance of the negative beta binomial law are
$$E=\frac{rb}{a-1}\quad,\quad V=\frac{rb(r+a-1)(a+b-1)}{(a-1)^2(a-2)}$$
the corresponding skewness is given by the following formula,
$$\gamma=\frac{(2r+a-1)(2b+a-1)\sqrt{a-2}}{(a-3)\sqrt{rb(r+a-1)(a+b-1)}}$$
and the corresponding kurtosis is given by a formula of the following type,
$$\kappa=\frac{(a-2)P(r,a,b)}{rb(a-3)(a-4)(r+a-1)^2(a+b-1)^2}$$
assuming $a>1,2,3,4$ respectively, with $P$ being a certain polynomial.
\end{theorem}

\begin{proof}
This follows indeed from Theorem 12.21, good exercise for you, and good exercise as well, to study that polynomial $P$, and let me know if you find something.
\end{proof}

\section*{12d. Circular measures}

With the end of this chapter approaching, time for some physics, as per our usual policy in this Part III. However, it is not very clear what type of physics we can do, in the discrete setting, and as usual in such situations, I will have to ask for advice.

\bigskip

And with cats and rats gone, I will turn to neighbor's goats. These fellows are quite fun, and certainly doing well their goat work, but in what regards physics and philosophy, I'm a bit worried, I mean these guys are constantly into ramming, they even taught this to some of the cats, who are now into ramming too. This being said, never knows, so let's see what one of these fellows thinks, about the mysteries of quantum physics: 

\begin{goat}
You should do some math for the ADE graphs, I mean these are like Mike Tyson or Michael Jordan, the greatest of all time.
\end{goat}

Well, quite interesting what you say, goat, so thanks a lot for the suggestion, and enjoy your meal. So, getting now to what goat says, let us start with:

\index{algebraic manifolds}
\index{quantum field theory}
\index{subfactors}
\index{knots and links}
\index{planar algebras}

\begin{fact}
The ADE graphs classify the following:
\begin{enumerate}
\item Basic Lie groups and algebras.

\item Subgroups of $SU_2$ and of $SO_3$.

\item Singularities of algebraic manifolds.

\item Basic invariants of knots and links.

\item Subfactors and planar algebras of small index.

\item Subgroups of the quantum permutation group $S_4^+$.

\item Basic quantum field theories, and other physics beasts.
\end{enumerate}
\end{fact}

Which sounds quite exciting, good mathematics that we will be learning here. In practice now, the A graphs are as follows, with the distinguished vertex being denoted $\bullet$, and with $A_n$ having by definition $n\geq2$ vertices, and $\tilde{A}_{2n}$ having $2n\geq2$ vertices:
$$A_n=\bullet-\circ-\circ\cdots\circ-\circ-\circ\hskip18mm 
A_{\infty}=\bullet-\circ-\circ-\circ\cdots\hskip7mm$$
\vskip-3mm
$$\ \ \ \ \ \ \ \tilde{A}_{2n}=
\begin{matrix}
\circ&\!\!\!\!-\circ-\circ\cdots\circ-\circ-&\!\!\!\!\circ\\
|&&\!\!\!\!|\\
\bullet&\!\!\!\!-\circ-\circ-\circ-\circ-&\!\!\!\!\circ\\
\\
\\
\end{matrix}\hskip20mm 
\tilde{A}_\infty=
\begin{matrix}
\circ&\!\!\!\!-\circ-\circ-\circ\cdots\\
|&\\
\bullet&\!\!\!\!-\circ-\circ-\circ\cdots\\
\\
\\
\end{matrix}
\hskip15mm$$
\vskip-7mm

Next come the D graphs, which are as follows, with $D_n$ having $n\geq3$ vertices, $\tilde{D}_n$ having $n+1\geq5$ vertices, and with $D_\infty$ being the limiting graph of $\tilde{D}_n$:
$$D_n=\bullet-\circ-\circ\dots\circ-
\begin{matrix}\ \circ\\
\ |\\
\ \circ \\
\ \\
\  \end{matrix}-\circ\hskip71mm$$
\vskip-7mm
$$\hskip7mm\tilde{D}_n=\bullet-
\begin{matrix}\circ\\
|\\
\circ\\
\ \\
\ \end{matrix}-\circ\dots\circ-
\begin{matrix}\ \circ\\
\ |\\
\ \circ \\
\ \\
\  \end{matrix}-\circ\hskip18mm$$
\vskip-7mm
$$\hskip50mm D_\infty=\bullet-
\begin{matrix}\circ\\
|\\
\circ\\
\ \\
\ \end{matrix}-\circ-\circ\cdots$$
\vskip-7mm

As a comment now, the labeling conventions for the AD graphs, while very standard, can be a bit confusing. The first graph in each series is by definition as follows:
$$A_2=\bullet-\circ\hskip13mm 
\tilde{A}_2=\begin{matrix}
\circ\\
||\\
\bullet\\
&\\
&\\
\end{matrix}\hskip13mm 
D_3=\begin{matrix}\ \circ\\
\ |\\
\ \bullet \\
\ \\
\  \end{matrix}-\circ \hskip13mm
\tilde{D}_4=\bullet-\!\!\!\!\!\begin{matrix}
\circ\hskip5mm \circ\\
\backslash\ \,\slash\\
\circ\\
&\\
&\\
\end{matrix}\!\!\!\!\!\!\!\!\!\!-\circ$$
\vskip-7mm

Finally, there are also a number of exceptional ADE graphs. First we have:
$$E_6=\bullet-\circ-
\begin{matrix}\circ\\
|\\
\circ\\
\ \\
\ \end{matrix}-
\circ-\circ\hskip71mm$$
\vskip-13mm
$$E_7=\bullet-\circ-\circ-
\begin{matrix}\circ\\
|\\
\circ\\
\ \
\\
\ \end{matrix}-
\circ-\circ\hskip18mm$$
\vskip-15mm
$$\hskip30mm E_8=\bullet-\circ-\circ-\circ-
\begin{matrix}\circ\\
|\\
\circ\\
\ \\
\ \end{matrix}-
\circ-\circ$$
\vskip-5mm

Then, we have extended versions of the above exceptional graphs, as follows:
$$\tilde{E}_6=\bullet-\circ-\begin{matrix}
\circ\\
|
\\
\circ\\
|&\\
\circ&\!\!\!\!-\ \circ\\
\ \\
\   \\
\ \\
\ \end{matrix}-\circ\hskip71mm$$
\vskip-23mm
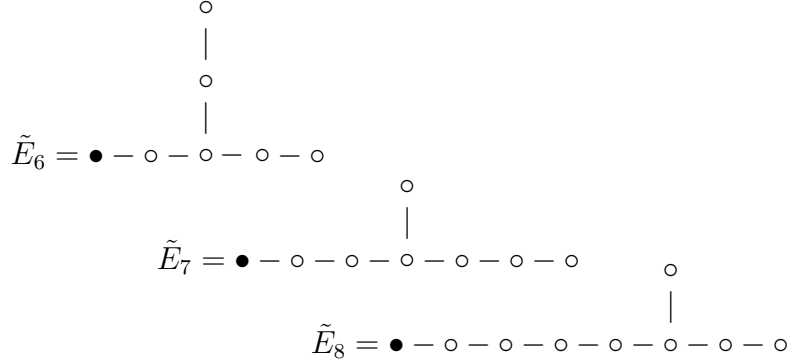
$$\tilde{E}_7=\bullet-\circ-\circ-
\begin{matrix}\circ\\
|\\
\circ\\
\ \\
\ \end{matrix}-
\circ-\circ-\circ\hskip18mm$$
\vskip-15mm
$$\hskip30mm \tilde{E}_8=\bullet-\circ-\circ-\circ-\circ-
\begin{matrix}\circ\\
|\\
\circ\\
\ \\
\ \end{matrix}-
\circ-\circ$$
\vskip-5mm

Getting to work now, as a first observation, we know that $A_\infty$ and $\tilde{A}_\infty$ are the graphs that we previously called $\mathbb N$ and $\mathbb Z$, and studied in some detail, in chapter 3. So, based on our previous work for these graphs, let us formulate the following definition:

\index{Poincar\'e series}
\index{bipartite graph}
\index{positive spectral measure}

\begin{definition}
The positive spectral measure $\nu$ of a rooted bipartite graph $X$ is the real probability measure having as moments
$$\int_\mathbb Rx^kd\nu(x)=L_{2k}$$
where $L_{2k}$ is the number of $2k$-loops based at the root. Equivalently, we must have
$$\int_\mathbb R\frac{1}{1-xz}\,d\nu(x)=f(z)$$
where $f(z)=\sum_{k=0}^\infty L_{2k}z^k$ is the Poincar\'e series of $X$.
\end{definition}

Here the existence of $\nu$, and the fact that this is indeed a positive measure, meaning a measure supported on $[0,\infty)$, comes from the following simple fact:

\begin{proposition}
The positive spectral measure of a rooted bipartite graph $X$ is given by the following formula, with $d$ being the adjacency matrix of the graph,
$$\nu=law(d^2)$$
and with the probabilistic computation being with respect to the expectation 
$$A\to<A>$$
with $<A>$ being the $(*,*)$-entry of a matrix $A$, where $*$ is the root.
\end{proposition}

\begin{proof}
With the above conventions, we have the following computation:
$$f(z)
=\sum_{k=0}^\infty L_{2k}z^k
=\sum_{k=0}^\infty\left<d^{2k}\right>z^k
=\left<\frac{1}{1-d^2z}\right>$$

But this shows that we have $\nu=law(d^2)$, as desired.
\end{proof}

Let us introduce as well the following notion, which is something more subtle, coming from the work of Jones on subfactor theory, and related topics \cite{jo1}, \cite{jo2}, \cite{jo3}, \cite{jo4}:

\index{circular measure}
\index{circular spectral measure}

\begin{definition}
The circular measure $\varepsilon$ of a rooted bipartite graph $X$ is given by
$$d\varepsilon(q)=d\nu((q+q^{-1})^2)$$
where $\nu$ is the associated positive spectral measure.
\end{definition}

To be more precise, we know from Proposition 12.26 that the positive measure $\nu$ is the spectral measure of a certain positive matrix, $d^2\geq0$, and it follows from this, and from basic spectral theory, that this measure is supported by the positive reals:
$$supp(\nu)\subset\mathbb R_+$$

But then, with this observation in hand, we can define indeed the circular measure $\varepsilon$ as above, as being the pullback of $\nu$ via the following map:
$$\mathbb R\cup\mathbb T\to\mathbb R_+\quad,\quad 
q\to (q+q^{-1})^2$$

As a basic example for this, to start with, assume that $\nu$ is a discrete measure, supported by $n$ positive numbers $x_1<\ldots<x_n$, with corresponding densities $p_1,\ldots,p_n$:
$$\nu=\sum_{i=1}^n p_i\delta_{x_i}$$

For each $i\in\{1,\ldots,n\}$ the equation $(q+q^{-1})^2=x_i$ has then four solutions, that we can denote $q_i,q_i^{-1},-q_i,-q_i^{-1}$. And with this notation, we have:
$$\varepsilon=\frac{1}{4}\sum_{i=1}^np_i\left(\delta_{q_i}+\delta_{q_i^{-1}}+\delta_{-q_i}+\delta_{-q_i^{-1}}\right)$$

In general, the basic properties of $\varepsilon$ can be summarized as follows:

\begin{theorem}
The circular measure has the following properties:
\begin{enumerate}
\item $\varepsilon$ has equal density at $q,q^{-1},-q,-q^{-1}$.

\item The odd moments of $\varepsilon$ are $0$.

\item The even moments of $\varepsilon$ are half-integers.

\item When $X$ has norm $\leq 2$, $\varepsilon$ is supported by the unit circle.

\item When $X$ is finite, $\varepsilon$ is discrete.

\item If $K$ is a solution of $d=K+K^{-1}$, then $\varepsilon=law(K)$. 
\end{enumerate}
\end{theorem}

\begin{proof}
These results can be deduced from definitions, the idea being that (1-5) are trivial, and that (6) follows from the formula of $\nu$ from Proposition 12.26.
\end{proof}

Getting now to computations, we first have the following result:

\begin{theorem}
The circular measure of the basic index $4$ graph, namely 
$$\begin{matrix}
&\circ&\!\!\!\!-\circ-\circ\cdots\circ-\circ-&\!\!\!\!\circ\cr
\tilde{A}_{2n}=&|&&\!\!\!\!|\cr
&\bullet&\!\!\!\!-\circ-\circ-\circ-\circ-&\!\!\!\!\circ\cr\cr\cr\end{matrix}$$
\vskip-10mm

\noindent is the uniform measure on the $2n$-roots of unity.
\end{theorem}

\begin{proof}
Let us identify the vertices of $X=\tilde{A}_{2n}$ with the group $\{w^k\}$ formed by the $2n$-th roots of unity in the complex plane, where $w=e^{\pi i/n}$. The adjacency matrix of $X$ acts then on the functions $f\in C(X)$ in the following way:
$$df(w^s)=f(w^{s-1})+f(w^{s+1})$$

But this shows that we have $d=K+K^{-1}$, where $K$ is given by:
$$Kf(w^s)=f(w^{s+1})$$

Thus we can use Proposition 12.26 and Theorem 12.28 (6), and we get:
$$\varepsilon=law(K)$$

But this is the uniform measure on the $2n$-roots of unity, as claimed.
\end{proof}

Observe that Theorem 12.29 is something quite interesting, I mean you can try to do a few loop computations for $\tilde{A}_{2n}$, with bare hands, and these will lead you nowhere. So, what we have here is a good idea, with Definition 12.27 being something quite subtle.

\bigskip

Moving on now, with this circular measure idea of Jones adopted, in order to deal with the other ADE graphs, let us introduce the following densities:
$$\alpha=Re(1-q^2)$$
$$\beta=Re(1-q^4)$$
$$\gamma=Re(1-q^6)$$

We have then the following result, $d_n$ being the uniform measure on the $2n$-th roots of unity, and $d_n'=2d_{2n}-d_n$ being the uniform measure on the odd $4n$-roots of unity:

\index{ADE graph}
\index{circular measure}
\index{cyclotomic measure}

\begin{theorem}
The circular measures of ADE graphs are given by:
\begin{enumerate}
\item $A_{n-1}\to\alpha_n$.

\item $\tilde{A}_{2n}\to d_n$.

\item $D_{n+1}\to\alpha_n'$.

\item $\tilde{D}_{n+2}\to (d_n+d_1')/2$.

\item $E_6\to\alpha_{12}+(d_{12}-d_6-d_4+d_3)/2$.

\item $E_7\to\beta_9'+(d_1'-d_3')/2$.

\item $E_8\to\alpha_{15}'+\gamma_{15}'-(d_5'+d_3')/2$.

\item $\tilde{E}_{n+3}\to (d_n+d_3+d_2-d_1)/2$.
\end{enumerate}
\end{theorem}

\begin{proof}
This is something which can be proved in three steps, as follows:

\medskip

(1) For the simplest graph, namely the circle $\tilde{A}_{2n}$, we already have the result, from Theorem 12.29, with the proof there being something elementary.

\medskip

(2) For the other non-exceptional graphs, that is, of type A and D, the same method works, namely direct loop counting, with some matrix tricks.

\medskip

(3) As for the exceptional graphs, of type E, basically the same method works, but the computations are more complicated, using some root of unity know-how. 

\medskip

(4) So, this was for the idea, and in practice, read Jones \cite{jo1}, \cite{jo2}, \cite{jo3}, \cite{jo4}, and then have a look at my paper with Bisch \cite{bbi}, where all this is explained.
\end{proof}

And we will end our study of the discrete measures with this. As a conclusion, at the advanced level, the study of discrete laws is all mutidimensional analysis.

\section*{12e. Exercises}

Exciting chapter that we had here, and as exercises on this, we have:

\begin{exercise}
Work out the details, for the moment formula for $c_{mp}$.
\end{exercise}

\begin{exercise}
Work out the details, for the formulae of $\gamma,\kappa$ for $c_{mp}$.
\end{exercise}

\begin{exercise}
Fill in the missing details, for the hypergeometric laws.
\end{exercise}

\begin{exercise}
Fill in the missing details, for the negative hypergeometric laws.
\end{exercise}

\begin{exercise}
Learn more about the beta distributions, and their applications.
\end{exercise}

\begin{exercise}
Learn more about the binomial beta laws, positive and negative.
\end{exercise}

\begin{exercise}
Clarify what we said, regarding the circular measures.
\end{exercise}

\begin{exercise}
Learn as well about blowup on more complicated manifolds.
\end{exercise}

As bonus exercise, reiterated, learn some math from Jones \cite{jo1}, \cite{jo2}, \cite{jo3}, \cite{jo4}.

\part{Quantum versions}

\ \vskip50mm

\begin{center}
{\em Ma il cammino di ogni speranza

Si ferma un momento e poi se ne va

La tua voce gia vola nel vento

Non e che un ricordo e non tornera}
\end{center}

\chapter{Free probability}

\section*{13a. Free variables}

Welcome to free probability. We have met some already, in this book, in the context of our various considerations in relation with groups, random matrices, or just moment combinatorics, with our knowledge being summarized in the following principle:

\begin{principle}
The free analogues of the Gaussian and Poisson laws $g_t,p_t$ are the Wigner and Marchenko-Pastur laws $\gamma_t,\pi_t$.
\end{principle}

However, this remains something vague, more of a physics finding, and our goal in this chapter will be that of putting this on firm ground. Following Voiculescu \cite{vo1} we will talk about freeness, free convolution and free limiting theorems, and with $\gamma_t,\pi_t$ appearing via the free analogues of the CLT and PLT, the above principle will become a Theorem.

\bigskip

Then afterwards, once we have our theory, we will start systematically developing it, with the aim of having it on par with the classical theory, developed in Parts I-II-III. Many things to be done here, and chapters 14-15-16 will be an introduction to that.

\bigskip

Getting started now, what is freeness? Nothing left to lose, Janis Joplin would say, although Nelson Mandela might claim the opposite. In short, take it easy, this is something quite subtle, that will take us some time to understand. Mathematically, we have:

\begin{principle}
The scheme for reaching to freeness is as follows,
$${\rm commutativity}\quad\to\quad{\rm noncommutativity}\quad\to\quad {\rm freeness}$$
with some basic illustrations for this being:
\begin{enumerate}
\item Abelian groups $\to$ arbitrary groups $\to$ free groups.

\item Commuting algebras $\to$ noncommuting algebras $\to$ free algebras.

\item Commuting variables $\to$ noncommuting variables $\to$ free variables.

\item Classical geometry $\to$ noncommutative geometry $\to$ free geometry.

\item Classical mechanics $\to$ quantum mechanics $\to$ free mechanics.
\end{enumerate}
\end{principle}

To be more precise here, (1) is something self-explanatory, with commutativity in the group theory setting, $gh=hg$, being called abelianity, and with the free groups $F_N=<g_1,\ldots,g_N|\,\emptyset>$ being, obviously, free. As for (2), this is still algebra, and as basic examples here, we can say that among the group algebras $A=\mathbb C[G]$, those coming from abelian groups are commutative, and those coming from free groups are free.

\bigskip

Regarding now (3), this is something more subtle, that we would like to understand. Formally this can be thought of as coming from (2), but as we perfectly know as probabilists, things are a bit more complicated than this. Indeed, the simplest examples of noncommuting variables are the usual matrices $A\in M_N(\mathbb C)$, followed by the random matrices $Z\in M_N(L^\infty(X))$, followed by the arbitrary operators $T\in B(H)$, and properly talking about freeness in this setting is most likely a non-trivial business. 

\bigskip

As for (4) and (5), these are something even more complicated, with (4) being a potential theory that can be built on top of (3), and we will see later in this book that this can be done indeed, and with (5) being a potential ``theory of everything'' that can be built on top of (4), but no one knows how to do it, for 50 years and counting.

\bigskip

So, this is the situation, layers of knowledge leading to freeness, or rather to us humans understanding the freeness that this world naturally has. In practice now, benefiting from our experience from chapter 8, we can have as entry point, in relation with (3):

\index{random variable}
\index{moments}
\index{colored moments}
\index{law}
\index{distribution}

\begin{definition}
Let $A$ be a $C^*$-algebra, that is, a complex algebra with an involution and norm satisfying $||aa^*||=||a||^2$, given with a positive unital trace $tr:A\to\mathbb C$.
\begin{enumerate}
\item The elements $a\in A$ are called random variables.

\item The moments of such a variable are the numbers $M_k(a)=tr(a^k)$.

\item The law of such a variable is the functional $\mu:P\to tr(P(a))$.
\end{enumerate}
\end{definition}

We refer to the material in chapter 8 for more on all this. In regards with the $C^*$-algebras, the summary of what we know, along with a little bit more, is as follows:

\bigskip

(1) The starting point are the bounded linear operators on a Hilbert space, $T\in B(H)$. Indeed, these operators have an adjoint operation given by $<T^*x,y>=<x,Ty>$, a norm given by $||T||=\sup_{||x||=1}||Tx||$, and the key formula $||TT^*||=||T||^2$ is satisfied. Thus, any norm closed $*$-algebra $A\subset B(H)$ is a $C^*$-algebra, in our sense.

\bigskip

(2) To be known about $C^*$-algebras are 3 theorems. First is the Gelfand theorem, related to the spectral theorem for normal operators, stating that the commutative $C^*$-algebras are precisely the algebras of type $A=C(X)$, with $X$ being a compact space. Which is something truly fundamental, that we will heavily use, in what follows.

\bigskip

(3) The second theorem is the Gelfand-Naimark-Segal one, stating that any $C^*$-algebra appears as closed $*$-subalgebra $A\subset B(H)$, with respect to some suitable Hilbert space, $H=L^2(A)$. Technically, we will not really need this. And finally, as a third theorem, that we know from chapter 7, in finite dimensions we have $A=\oplus_kM_{n_k}(\mathbb C)$.

\bigskip

So, this was for the brief story of the $C^*$-algebras. Regarding now the other things appearing in Definition 13.3, the situation with them is quite straightforward: 

\bigskip

(1) By using the Gelfand theorem we can talk about positive elements $a\geq0$, with these appearing as $a=bb^*$, and then about positive linear forms too, those satisfying $a\geq0\implies\varphi(a)\geq0$. But with this, we have our expectations $tr:A\to\mathbb C$, which must be by definition positive, unital, and satisfying the trace condition $tr(ab)=tr(ba)$.

\bigskip

(2) Next, in the definition of the moments, $k=\circ\bullet\bullet\circ\ldots$ is as usual a colored integer, and the corresponding powers $a^k$ are defined by $a^\circ=a$, $a^\bullet=a^*$ and multiplicativity. As for the definition of the law, we use there noncommuting $*$-polynomials in one variable, $P\in\mathbb C<X,X^*>$. Observe that the law is uniquely determined by the moments.

\bigskip

(3) Finally, the whole picture includes a speculation too. Indeed, according to the Gelfand theorem we can think of any $C^*$-algebra as being of the form $A=C(X)$, with $X$ being a ``compact quantum space''. And with this in hand, we can better understand what we are doing in Definition 13.3, that being, obviously, ``quantum probability''.

\bigskip

Generally speaking, Definition 13.3 is something quite abstract, and there is no other way of doing things, at that level of generality. However, remarkably, we have:

\index{normal element}
\index{spectral measure}
\index{functional calculus}
\index{Riesz theorem}

\begin{theorem}
Assuming that $a\in A$ is normal, $aa^*=a^*a$, its law corresponds to a probability measure on its spectrum $\sigma(a)\subset\mathbb C$, according to the following formula,
$$tr(P(a))=\int_{\sigma(a)}P(z)d\mu(z)$$
and in the self-adjoint case, $a=a^*$, this law is a real probability measure. When the variable is not normal and $tr$ is faithful, such a probability measure does not exist.
\end{theorem}

\begin{proof}
This is something that we know from chapter 8, the idea being as follows:

\medskip

(1) Assuming that $a\in A$ is normal the algebra that it generates is commutative, $<a>=C(\sigma(a))$, so the functional $P\to tr(P(a))$ can be though of as being an integration functional on the spectrum $\sigma(a)\subset\mathbb C$, and the Riesz theorem provides us with the desired probability measure $\mu$, making the formula $tr(P(a))=\int P(z)d\mu(z)$ work.

\medskip

(2) As a technical comment here, when assuming that our positive trace $tr$ is strictly positive, of faithful, in the sense that $a>0\implies tr(a)>0$, the support of the spectral measure $\mu$ is the whole spectrum $\sigma(a)$. Finally, still talking support, when $a$ is self-adjoint, $a=a^*$, its spectrum is real, and so $\mu$ is a real probability measure, as stated.

\medskip

(3) In what regards now the last assertion, this is something that we know too from chapter 8, coming from $aa^*\neq a^*a\implies tr(aa^*aa^*)>tr(aaa^*a^*)$. Indeed, in view of this inequality, it is impossible to obtain the quantities $tr(aa^*aa^*)\neq tr(aaa^*a^*)$ by integrating $z\bar{z}z\bar{z}=zz\bar{z}\bar{z}$ with respect to a certain complex probability measure $\mu$.
\end{proof}

And with this, done with the generalities? Not yet, because there is an important technical comment to be formulated, about all this, as follows:

\begin{comment}
Quantum probability can be done in two possible ways:
\begin{enumerate}
\item Using $C^*$-algebras as we do. With these concretely appearing as norm closed $*$-algebras $A\subset B(H)$, being in the commutative case of the form $A=C(X)$, with $X$ being a compact space, and being in finite dimensions $A=\oplus_kM_{n_k}(\mathbb C)$.

\item Using von Neumann algebras. With these being the weakly closed $*$-subalgebras $A\subset B(H)$, being in the commutative case of the form $A=L^\infty(X)$, with $X$ being a measured space, and being in finite dimensions $A=\oplus_kM_{n_k}(\mathbb C)$.
\end{enumerate}
\end{comment}

So, this is the truth about quantum probability, theory coming in 2 stages, and in this book we will be learning the basics, at stage 1. With the remark of course that the algebras $A=L^\infty(X)$ being technically $C^*$-algebras, we will cover them too, but not in a fully satisfactory way, the problem coming from Gelfand, $A=C(\widehat{X})$ with $\widehat{X}$ being a certain abstract compactification of $X$, that you don't want to hear about. And, we will leave some thinking at all this, pros and cons of the $C^*$-algebra theory, as an exercise.

\bigskip

Getting to work now, as a first straightforward definition, we have:

\index{abstract independence}
\index{independent algebras}

\begin{definition}
We call two subalgebras $A,B\subset C$ independent when the following condition is satisfied, for any $a\in A$ and $b\in B$: 
$$tr(ab)=tr(a)tr(b)$$
Equivalently, the following condition must be satisfied, for any $a\in A$ and $b\in B$: 
$$tr(a)=tr(b)=0\implies tr(ab)=0$$
Also, $a,b\in C$ are called independent when $A=<a>$ and $B=<b>$ are independent.
\end{definition}

Observe that the above two independence conditions are indeed equivalent, with this following from the following computation, with the convention $a'=a-tr(a)$:
\begin{eqnarray*}
tr(ab)
&=&tr[(a'+tr(a))(b'+tr(b))]\\
&=&tr(a'b')+t(a')tr(b)+tr(a)tr(b')+tr(a)tr(b)\\
&=&tr(a'b')+tr(a)tr(b)\\
&=&tr(a)tr(b)
\end{eqnarray*}

It is possible to develop some theory here, but this leads to the usual CLT. As a much more interesting notion now, we have Voiculescu's freeness \cite{vo1}:

\index{freeness}
\index{free algebras}

\begin{definition}
Given a pair $(C,tr)$, we call two subalgebras $A,B\subset C$ free when the following condition is satisfied, for any $a_i\in A$ and $b_i\in B$:
$$tr(a_i)=tr(b_i)=0\implies tr(a_1b_1a_2b_2\ldots)=0$$
Also, $a,b\in C$ are called free when $A=<a>$ and $B=<b>$ are free.
\end{definition}

As a first observation on this notion, which is similar to independence, there is a certain lack of symmetry between Definition 13.6 and Definition 13.7, because the latter does not include an explicit formula for the following type of quantities:
$$tr(a_1b_1a_2b_2\ldots)$$

However, this is not an issue, and is simply due to the fact that the formula in the free case is something more complicated, the precise result being as follows:

\index{free trace}

\begin{proposition}
If $A,B\subset C$ are free, the restriction of $tr$ to $<A,B>$ can be computed in terms of the restrictions of $tr$ to $A,B$. To be more precise, we have
$$tr(a_1b_1a_2b_2\ldots)=P\Big(\{tr(a_{i_1}a_{i_2}\ldots)\}_i,\{tr(b_{j_1}b_{j_2}\ldots)\}_j\Big)$$
where $P$ is certain polynomial, depending on the length of $a_1b_1a_2b_2\ldots\,$, having as variables the traces of products $a_{i_1}a_{i_2}\ldots$ and $b_{j_1}b_{j_2}\ldots\,$, with $i_1<i_2<\ldots$ and $j_1<j_2<\ldots$
\end{proposition}

\begin{proof}
This is something a bit theoretical, so let us begin with an example. The computation after Definition 13.6 perfectly works when $a,b$ are free, and gives:
$$tr(ab)=tr(a)tr(b)$$

In general, the situation is of course more complicated than this, but the same trick applies. To be more precise, we can start our computation as follows:
\begin{eqnarray*}
tr(a_1b_1a_2b_2\ldots)
&=&tr\big[(a_1'+tr(a_1))(b_1'+tr(b_1))(a_2'+tr(a_2))\ldots\big]\\
&=&tr(a_1'b_1'a_2'b_2'\ldots)+{\rm other\ terms}\\
&=&{\rm other\ terms}
\end{eqnarray*}

Now regarding the ``other terms'', those which are left, each of them will consist of a product of traces of type $tr(a_i)$ and $tr(b_i)$, and then a trace of a product still remaining to be computed, which is of the following form, for some elements $\alpha_i\in A$ and $\beta_i\in B$:
$$tr(\alpha_1\beta_1\alpha_2\beta_2\ldots)$$

Now since the length of $\alpha_1\beta_1\alpha_2\beta_2\ldots$ is smaller than the length of the original product $a_1b_1a_2b_2\ldots$, we are led into of recurrence, and this gives the result.
\end{proof}

Let us discuss now some examples of independence and freeness. We first have the following result, from \cite{vo1}, which is something straightforward:

\begin{theorem}
Given two algebras $(A,tr)$ and $(B,tr)$, the following hold:
\begin{enumerate}
\item $A,B$ are independent inside their tensor product $A\otimes B$, endowed with its canonical tensor product trace, given on basic tensors by $tr(a\otimes b)=tr(a)tr(b)$.

\item $A,B$ are free inside their free product $A*B$, endowed with its canonical free product trace, given by the formulae in Proposition 13.8.
\end{enumerate}
\end{theorem}

\begin{proof}
Both the above assertions are clear from definitions, as follows:

\medskip

(1) This is clear with either of the definitions of the independence, from Definition 13.6, because we have, by construction of the product trace:
$$tr(ab)
=tr[(a\otimes1)(1\otimes b)]
=tr(a\otimes b)
=tr(a)tr(b)$$

(2) This is clear too from definitions, the only point being that of showing that the notion of freeness, or the recurrence formulae in Proposition 13.8, can be used in order to construct a canonical free product trace, on the free product of the algebras involved:
$$tr:A*B\to\mathbb C$$

But this can be checked for instance by using a GNS construction. Indeed, consider the GNS constructions for the algebras $(A,tr)$ and $(B,tr)$:
$$A\to B(L^2(A))\quad,\quad
B\to B(L^2(B))$$

By taking the free product of these representations, we obtain a representation as follows, with the $*$ on the right being a free product of pointed Hilbert spaces:
$$A*B\to B(L^2(A)*L^2(B))$$

Now by composing with the linear form $T\to<T\xi,\xi>$, where $\xi=1_A=1_B$ is the common distinguished vector of $L^2(A)$, $L^2(B)$, we obtain a linear form, as follows:
$$tr:A*B\to\mathbb C$$

It is routine then to check that $tr$ is indeed a trace, and this is the ``canonical free product trace'' from the statement. Then, an elementary computation shows that $A,B$ are free inside $A*B$, with respect to this trace, and this finishes the proof.
\end{proof}

More concretely now, we have the following result, also from \cite{vo1}:

\begin{theorem}
We have the following results, valid for group algebras:
\begin{enumerate}
\item $C^*(G),C^*(H)$ are independent inside $C^*(G\times H)$.

\item $C^*(G),C^*(H)$ are free inside $C^*(G*H)$.
\end{enumerate}
\end{theorem}

\begin{proof}
This something very standard, the idea being as follows:

\medskip

(1) To start with, recall that associated to a discrete group $G$ is its group algebra $\mathbb C[G]$, formal linear span of the group elements $g\in G$, with involution $g^*=g^{-1}$. We can then talk about $C^*$-norms on this algebra, subject to $||aa^*||=||a||^2$, and define $C^*(G)$ as being the completion of  $\mathbb C[G]$ with respect to the biggest such norm.

\medskip

(2) And this, coming with the remark that such $C^*$-norms exist indeed, a basic example being the one coming from the embedding $\mathbb C[G]\subset B(l^2(G))$ given by $g(\delta_h)=\delta_{gh}$. And with the extra remark that the biggest $C^*$-norm is bounded, due to the fact that the standard generators $g\in G$ being unitaries, $g^*=g^{-1}$, they must satisfy $||g||=1$.

\medskip

(3) Getting now to what the statement says, we can use the general results in Theorem 13.9, along with the following two isomorphisms, which are both standard:
$$C^*(G\times H)=C^*(G)\otimes C^*(H)\quad,\quad
C^*(G*H)=C^*(G)*C^*(H)$$

Alternatively, we can check the independence and freeness formulae on group elements, which is something trivial, and then conclude by linearity.
\end{proof}

Summarizing, what we have so far is a notion of freeness, appearing as a natural analogue of the notion of independence. We will see in what follows that this analogy can be pushed to remarkably high levels, with free probability being quite often on par with classical probability, and sometimes even beating it, on certain aspects.

\section*{13b. Free convolution}

All the above was quite theoretical, and as a concrete application of the above results, bringing us into probability, we have the following result, from \cite{vo1}:

\index{free convolution}

\begin{theorem}
We have a free convolution operation $\boxplus$ for the distributions
$$\mu:\mathbb C<X,X^*>\to\mathbb C$$
which is well-defined by the following formula, with $a,b$ taken to be free:
$$\mu_a\boxplus\mu_b=\mu_{a+b}$$
This restricts to an operation, still denoted $\boxplus$, on the real probability measures.
\end{theorem}

\begin{proof}
We have several verifications to be performed here, as follows:

\medskip

(1) We first have to check that given two variables $a,b$, which live respectively in certain $C^*$-algebras $A,B$, we can recover them inside some $C^*$-algebra $C$, with exactly the same distributions $\mu_a,\mu_b$, as to be able to sum them and talk about $\mu_{a+b}$. But this comes from Theorem 13.9, because we can set $C=A*B$, as explained there.

\medskip

(2) The other verification which is needed is that of the fact that if two variables $a,b$ are free, then the distribution $\mu_{a+b}$ depends only on the distributions $\mu_a,\mu_b$. But for this purpose, we can use the general formula from Proposition 13.8, namely:
$$tr(a_1b_1a_2b_2\ldots)=P\Big(\{tr(a_{i_1}a_{i_2}\ldots)\}_i,\{tr(b_{j_1}b_{j_2}\ldots)\}_j\Big)$$

Indeed, by plugging in arbitrary powers of $a,b$ as variables $a_i,b_j$, we obtain a family of formulae of the following type, with $Q$ being certain polyomials:
$$tr(a^{k_1}b^{l_1}a^{k_2}b^{l_2}\ldots)=Q\Big(\{tr(a^k)\}_k,\{tr(b^l)\}_l\Big)$$

Thus the moments of $a+b$ depend only on the moments of $a,b$, with of course colored exponents in all this, according to our moment conventions, and this gives the result.

\medskip

(3) Finally, in what regards the last assertion, this is clear from the fact that if the variables $a,b$ are self-adjoint, then so is their sum $a+b$.
\end{proof}

Along the same lines, but with some technical subtleties this time, we can talk as well about multiplicative free convolution, following \cite{vo2}, as follows:

\index{multiplicative free convolution}
\index{self-adjoint variable}

\begin{theorem}
We have a multiplicative free convolution operation $\boxtimes$ on the distributions, which is well-defined by the following formula, with $a,b$ taken to be free:
$$\mu_a\boxtimes\mu_b=\mu_{ab}$$
In the case of positive variables, we can equally set, with the same outcome
$$\mu_a\boxtimes\mu_b=\mu_{\sqrt{a}\,b\,\sqrt{a}}$$
and so we have an operation, still denoted $\boxtimes$, on the probability measures on $\mathbb R_+$.
\end{theorem}

\begin{proof}
We have two statements here, the idea being as follows:

\medskip

(1) The verifications for the fact that $\boxtimes$ as above is indeed well-defined at the general distribution level are identical to those done before for $\boxplus$, with the result basically coming from the formula in Proposition 13.8, and with Theorem 13.9 invoked as well, in order to say that we have a model, and so we can indeed use that formula.

\medskip

(2) Regarding the last assertion, real measures, this was something trivial for $\boxplus$, but is something trickier for $\boxtimes$, because if we take $a,b$ to be self-adjoint, their product $ab$ will be not self-adjoint, in general. However, when $a>0$ we can trick, by setting:
$$c=\sqrt{a}\,b\,\sqrt{a}$$

Indeed, this new variable is then self-adjoint, and its moments are given by:
\begin{eqnarray*}
tr(c^k)
&=&tr[(\sqrt{a}\,b\,\sqrt{a})^k]\\
&=&tr[\sqrt{a}\,baba\ldots bab\,\sqrt{a}]\\
&=&tr[\sqrt{a}\cdot \sqrt{a}\,baba\ldots bab]\\
&=&tr[(ab)^k]
\end{eqnarray*}

We conclude that the $\boxtimes$ operation restricts into an operation as follows:
$$\boxtimes:\mathcal P(\mathbb R_+)\times\mathcal P(\mathbb R)\to\mathcal P(\mathbb R)$$

Moreover, since $b>0\implies c>0$, we can make things symmetric, as stated.
\end{proof}

The problem is now how to linearize the above operations $\boxplus$ and $\boxtimes$. And here, leaving aside the arbitrary distributions, and focusing on the case of the real measures, and leaving aside as well $\boxtimes$, say for later, we are led to the following concrete question:

\begin{question}
What is the free analogue of the logarithm of the Fourier transform $F_f(x)=E(e^{ixf})$, linearizing the operation $\boxplus$ for the real probability measures?
\end{question}

In answer now, this is something quite tricky, and Voiculescu's idea in \cite{vo1} was that of solving this problem by temporarily exiting the class of real probability measures, and even the general tracial framework of Definition 13.3. Let us start with:

\index{shift}

\begin{theorem}
Consider the shift operator on $H=l^2(\mathbb N)$, given by $S(e_i)=e_{i+1}$. The operators of the following type, with $f\in\mathbb C[X]$ being a polynomial, 
$$T=S^*+f(S)$$
model then in uncolored moments, up to finite order, all the distributions $\mu:\mathbb C[X]\to\mathbb C$, when regarded as random variables with respect to the state $\varphi(T)=<Te_0,e_0>$.
\end{theorem}

\begin{proof}
There are two things to be done here, the hard one, understanding what the statement exactly says, and then the easy one, proving what is to be proved:

\medskip

(1) Regarding the statement, as just mentioned, this is something tricky, exiting the framework of Definition 13.3. That definition was something clean, with the expectation being a trace, $tr(ab)=tr(ba)$, and with us being interested in the colored moments $M_k(a)$. What we are doing here is to use the state $\varphi(T)=<Te_0,e_0>$, which is not a trace, and get interested in uncolored moments only, that is, $M_k(a)=\varphi(a^k)$ with $k\in\mathbb N$.

\medskip

(2) Why all this, you would wonder. In answer, not my idea, that is Voiculescu's. And as we will soon discover, such tricks can eventually solve Question 13.13.

\medskip

(3) Good, so with the statement agreed upon, by presidential decision, let us turn now to the proof. The shift and its adjoint are given by the following formulae:
$$S(e_i)=e_{i+1}\quad,\quad S^*(e_i)=
\begin{cases}
e_{i-1}&(i>0)\\
0&(i=0)
\end{cases}$$

Consider now a variable as in the statement, written as follows:
$$T=S^*+a_01+a_1S+a_2S^2+\ldots+a_nS^n$$

The computation of the moments of $T$ amounts then in expanding $T^k$, and looking for terms which simplify up to 1, according to the rule $S^*S=1$. In practice, we get:
$$\varphi(T)=\varphi(a_01)=a_0$$
$$\varphi(T^2)=\varphi(a_0^211+a_1S^*S)=a_0^2+a_1$$
$$\varphi(T^3)=\varphi(a_0^3111+a_0a_11S^*S+a_0a_1S^*1S+a_0a_1S^*S1)=a_0^2+3a_0a_1$$

At order 4 the computation is similar, but a bit more complex, as follows:
\begin{eqnarray*}
\varphi(T^4)
&=&\varphi\big[a_0^41111+a_0^2a_1(11S^*S+1S^*1S+1S^*S1+S^*11S+S^*1S1+S^*S11)\\
&&+a_1^2(S^*S^*SS+S^*SS^*S)+a_2S^*S^*S^2\big]\\
&=&a_0^4+6a_0^2a_1+2a_1^2+a_2
\end{eqnarray*}

And so on, you get the point, the idea being that we are led into a certain recurrence, with the $2k$-th moment of our distribution $\varphi(T^{2k})$ providing us with the needed value of the coefficient $a_k$. Thus, we are led to the conclusion in the statement.
\end{proof}

Before getting further, with free products of models as above, let us work out a basic example, which is something fundamental, that we will need in what follows:

\index{random walk}
\index{Dyck paths}
\index{Catalan numbers}
\index{Wigner law}
\index{semicircle law}

\begin{proposition}
In the context of the above correspondence, the variable
$$T=S+S^*$$
follows the Wigner semicircle law $\gamma_1=\frac{1}{2\pi}\sqrt{4-x^2}\,dx$.
\end{proposition}

\begin{proof}
In order to compute the law of variable $T$ in the statement, we can use the moment method. Indeed, the moments of this variable are as follows:
$$M_k
=\varphi(T^k)
=\varphi((S+S^*)^k)
=\#(1\in(S+S^*)^k)$$

Now since $S$ shifts to the right on $\mathbb N$, and $S^*$ shifts to the left, while remaining positive, we are left with counting the length $k$ paths on $\mathbb N$ starting and ending at 0. But this is something that we are familiar with, the number of such paths being $L_{2k}=C_k$, the Catalan numbers, and with the corresponding measure being the Wigner law $\gamma_1$.
\end{proof}

Getting back now to our linearization program for $\boxplus$, the next step is that of taking a free product of the model found in Theorem 13.14 with itself.  There are two approaches here, one being a bit abstract, and the other one being more concrete. We first have:

\index{semigroup algebra}
\index{shift}

\begin{proposition}
We can talk about semigroup algebras $C^*_{red}(G)\subset B(l^2(G))$, a bit as we did for the group algebras, and at the level of examples:
\begin{enumerate}
\item With $G=\mathbb N$ we recover the shift algebra $A=<S>$ on $H=l^2(\mathbb N)$.

\item With $G=\mathbb N*\mathbb N$, we obtain the algebra $A=<S_1,S_2>$ on $H=l^2(\mathbb N*\mathbb N)$.
\end{enumerate}
\end{proposition}

\begin{proof}
We can talk indeed about semigroup algebras $C^*_{red}(G)\subset B(l^2(G))$, with the semigroup elements $g\in G$ being now isometries, and with this in hand:

\medskip

(1) With $G=\mathbb N$ we recover the shift algebra $A=<S>$ on the Hilbert space $H=l^2(\mathbb N)$, the shift $S$ itself being the isometry associated to the element $1\in\mathbb N$.

\medskip

(2) With $G=\mathbb N*\mathbb N$ we recover the double shift algebra $A=<S_1,S_2>$ on the Hilbert space $H=l^2(\mathbb N*\mathbb N)$, the two shifts $S_1,S_2$ themselves being the isometries associated to two copies of the element $1\in\mathbb N$, one for each of the two copies of $\mathbb N$ which are present.
\end{proof}

In what follows we will rather use an equivalent, second approach to our problem, which is exactly the same thing, but formulated in a less abstract way, as follows:

\index{free Fock space}
\index{creation operator}

\begin{proposition}
We can talk about the algebra of creation operators $S_x:v\to x\otimes v$ on the free Fock space associated to a real Hilbert space $H$, given by 
$$F(H)=\mathbb C\Omega\oplus H\oplus H^{\otimes2}\oplus\ldots$$
and at the level of examples, we have:
\begin{enumerate}
\item With $H=\mathbb C$ we recover the shift algebra $A=<S>$ on $H=l^2(\mathbb N)$.

\item With $H=\mathbb C^2$, we obtain the algebra $A=<S_1,S_2>$ on $H=l^2(\mathbb N*\mathbb N)$.
\end{enumerate}
\end{proposition}

\begin{proof}
We can talk indeed about the algebra $A(H)$ of creation operators on the free Fock space $F(H)$ associated to a real Hilbert space $H$, with the remark that, in terms of the abstract semigroup notions from Proposition 13.16, we have:
$$A(\mathbb C^k)=C^*_{red}(\mathbb N^{*k})\quad,\quad 
F(\mathbb C^k)=l^2(\mathbb N^{*k})$$

As for the assertions (1,2) in the statement, these are both clear, either directly, or by passing via (1,2) from Proposition 13.16, which were both clear as well.
\end{proof}

The advantage with this latter model comes from the following result, from \cite{vo1}, which has a very simple formulation, without linear combinations or anything:

\index{creation operator}
\index{annihilation operator}
\index{vacuum vector}
\index{free Fock space}

\begin{proposition}
Given a real Hilbert space $H$, and two orthogonal vectors $x\perp y$, the corresponding creation operators $S_x$ and $S_y$ are free with respect to
$$\varphi(T)=<T\Omega,\Omega>$$
called state associated to the vacuum vector.
\end{proposition}

\begin{proof}
In standard tensor product notation for the elements of the free Fock space $F(H)$, the formula of a creation operator associated to a vector $x\in H$ is as follows:
$$S_x(y_1\otimes\ldots\otimes y_n)=x\otimes y_1\otimes\ldots\otimes y_n$$

As for the formula of the adjoint of this creation operator, called annihilation operator associated to the vector $x\in H$, this is as follows: 
$$S_x^*(y_1\otimes\ldots\otimes y_n)=<x,y_1>\otimes y_2\otimes\ldots\otimes y_n$$

We obtain from this the following formula, which holds for any two vectors $x,y\in H$:
$$S_x^*S_y=<x,y>id$$

But with these formulae in hand, the result follows by doing some elementary computations, in the spirit of those done for the group algebras, in the above.
\end{proof}

With this technology in hand, let us go back to our linearization program for $\boxplus$. We know from Theorem 13.14 how to model the individual distributions $\mu\in\mathcal P(\mathbb R)$, and by combining this with Proposition 13.17 and Proposition 13.18, we therefore know how to freely model pairs of distributions $\mu,\nu\in\mathcal P(\mathbb R)$, as required by the convolution problem. We are therefore left with doing the sum in the model, and then computing its distribution. And the point here is that, still following \cite{vo1}, we have the following key result:

\index{freeness}

\begin{theorem}
Given two polynomials $f,g\in\mathbb C[X]$, consider the variables 
$$S^*+f(S)\quad,\quad 
T^*+g(T)$$
where $S,T$ are two creation operators, or shifts, associated to a pair of  orthogonal norm $1$ vectors. These variables are then free, and their sum has the same law as
$$R^*+(f+g)(R)$$
with $R$ being the usual shift on $l^2(\mathbb N)$.
\end{theorem}

\begin{proof}
We have two assertions here, the idea being as follows:

\medskip

(1) The freeness assertion comes from the general freeness result from Proposition 13.18, via the various identifications coming from the previous results.

\medskip

(2) Regarding the second assertion, the idea is that this comes from a $45^\circ$ rotation trick. Let us write indeed the two variables in the statement as follows:
$$X=S^*+a_0+a_1S+a_2S^2+\ldots$$
$$Y=T^*+b_0+b_1T+a_2T^2+\ldots$$

Now let us perform the following $45^\circ$ base change, on the real span of the vectors $s,t\in H$ producing our two shifts $S,T$, as follows:
$$r=\frac{s+t}{\sqrt{2}}\quad,\quad 
u=\frac{s-t}{\sqrt{2}}$$

The new shifts, associated to these vectors $r,u\in H$, are then given by:
$$R=\frac{S+T}{\sqrt{2}}\quad,\quad
U=\frac{S-T}{\sqrt{2}}$$

By using now these two new shifts, which are free according to Proposition 13.18, we obtain the following equality of distributions:
\begin{eqnarray*}
X+Y
&=&S^*+T^*+\sum_ka_kS^k+b_kT^k\\
&=&\sqrt{2}R^*+\sum_ka_k\left(\frac{R+U}{\sqrt{2}}\right)^k+b_k\left(\frac{R-U}{\sqrt{2}}\right)^k\\
&\sim&\sqrt{2}R^*+\sum_ka_k\left(\frac{R}{\sqrt{2}}\right)^k+b_k\left(\frac{R}{\sqrt{2}}\right)^k\\
&\sim&R^*+\sum_ka_kR^k+b_kR^k
\end{eqnarray*}

To be more precise, here at the end we have used the freeness property of $R,U$ in order to cut $U$ from the computation, as it cannot bring anything, and then we did a basic rescaling at the very end. Thus, we are led to the conclusion in the statement.
\end{proof}

As a conclusion, the operation $\mu\to f$ from Theorem 13.14 linearizes $\boxplus$. In order to finish, we are left with a routine computation inside $C^*(\mathbb N)$, which leads to:

\index{Cauchy transform}
\index{R-transform}
\index{free convolution}
\index{free Fourier transform}

\begin{theorem}
Given a real probability measure $\mu$, define its $R$-transform as follows:
$$G_\mu(\xi)=\int_\mathbb R\frac{d\mu(x)}{\xi-x}\implies G_\mu\left(R_\mu(z)+\frac{1}{z}\right)=z$$
The free convolution operation is then linearized by this $R$-transform.
\end{theorem}

\begin{proof}
Still following Voiculescu \cite{vo1}, this can be done as follows:

\medskip

(1) Consider a variable as in Theorem 13.14, as follows, with $f$ being a polynomial:
$$X=S^*+f(S)$$

In order to establish the result, we must prove that the $R$-transform of $X$, constructed according to the procedure in the statement, is the polynomial $f$ itself.

\medskip

(2) In order to do so, we fix $|z|<1$ in the complex plane, and we set:
$$q_z=\delta_0+\sum_{k=1}^\infty z_k\delta_k$$

The shift operator and its adjoint act then on this vector as follows:
$$Sq_z=z^{-1}(q_z-\delta_0)\quad,\quad 
S^*q_z=zq_z$$

It follows that the adjoint of our operator $X$ acts on this vector as follows:
\begin{eqnarray*}
X^*q_z
&=&(S+f(S^*))q_z\\
&=&z^{-1}(q_z-\delta_0)+f(z)q_z\\
&=&(z^{-1}+f(z))q_z-z^{-1}\delta_0
\end{eqnarray*}

Now observe that the above formula can be written in the following way:
$$z^{-1}\delta_0=(z^{-1}+f(z)-X^*)q_z$$

The point now is that when $|z|$ is small, the operator appearing on the right is invertible. Thus, with this assumption, we can rewrite the above formula as follows:
$$(z^{-1}+f(z)-X^*)^{-1}\delta_0=zq_z$$

Now by applying the trace, we are led to the following formula:
\begin{eqnarray*}
tr\left[(z^{-1}+f(z)-X^*)^{-1}\right]
&=&\left<(z^{-1}+f(z)-X^*)^{-1}\delta_0,\delta_0\right>\\
&=&<zq_z,\delta_0>\\
&=&z
\end{eqnarray*}

(3) Let us apply now the procedure in the statement to the real probability measure $\mu$ modeled by our variable $X$. The Cauchy transform $G_\mu$ is given by:
$$G_\mu(\xi)
=tr((\xi-X)^{-1})
=\overline{tr\Big((\bar{\xi}-X^*)^{-1}\Big)}
=tr((\xi-X^*)^{-1})$$

Now observe that, with the choice $\xi=z^{-1}+f(z)$ for our complex variable, the trace formula found in (2) above tells us that we have:
$$G_\mu\big(z^{-1}+f(z)\big)=z$$

Thus we have $R_\mu(z)=f(z)$, which finishes the proof, as explained in (1).
\end{proof}

\section*{13c. Limiting theorems}

With the above linearization technology in hand, we can do many things. First, we have the following free analogue of the CLT, due to Voiculescu \cite{vo1}:

\index{FCLT}
\index{Free CLT}
\index{R-transform}

\begin{theorem}[Free CLT]
Given self-adjoint variables $f_1,f_2,f_3,\ldots$ which are f.i.d., centered, with variance $t>0$, we have, with $n\to\infty$, in moments,
$$\frac{f_1+\ldots+f_n}{\sqrt{n}}\sim\gamma_t$$
with the limiting measure being the Wigner semicircle law $\gamma_t$.
\end{theorem}

\begin{proof}
The $R$-transform of the variable on the left in the statement can be computed by using the linearization property from Theorem 13.20, and is given by:
$$R(z)
=nR_f\left(\frac{z}{\sqrt{n}}\right)
\simeq tz$$

Regarding now the right term, as explained in chapter 3, the Cauchy transform is:
$$G_{\gamma_t}(\xi)=\frac{\xi-\sqrt{\xi^2-4t}}{2t}$$

But this gives $R_{\gamma_t}(z)=tz$, according to the following computation:
\begin{eqnarray*}
G_{\gamma_t}(tz+z^{-1})
&=&\frac{tz+z^{-1}-\sqrt{(tz+z^{-1})^2-4t}}{2t}\\
&=&\frac{tz+z^{-1}+tz-z^{-1}}{2t}\\
&=&z
\end{eqnarray*}

Observe that $R_{\gamma_t}(z)=tz$ follows in fact as well from the following formula, coming from Proposition 13.15, and from the technical details of the $R$-transform:
$$S+S^*\sim\gamma_1$$

Thus, the laws in the statement have the same $R$-transforms, so they are equal.
\end{proof}

Next, still following \cite{vo1}, we have the following free analogue of the CCLT:

\index{FCCLT}
\index{Free CCLT}

\begin{theorem}[Free CCLT]
Given random variables $f_1,f_2,f_3,\ldots$ which are f.i.d., centered, with variance $t>0$, we have, with $n\to\infty$, in moments,
$$\frac{f_1+\ldots+f_n}{\sqrt{n}}\sim\Gamma_t$$
where $\Gamma_t$ is the Voiculescu circular law of parameter $t$, given by the formula
$$\Gamma_t=law\left(\frac{a+ib}{\sqrt{2}}\right)$$
with $a,b$ being free, each following the Wigner semicircle law $\gamma_t$.
\end{theorem}

\begin{proof}
Let us decompose our variables into real and imaginary parts, as follows:
$$f_i=\frac{h_i+ik_i}{\sqrt{2}}$$

The variables $h_i$ and $k_i$ satisfy then the assumptions of the free CLT, so we have:
$$\frac{h_1+\ldots+h_n}{\sqrt{n}}\sim\gamma_t\quad,\quad 
\frac{k_1+\ldots+k_n}{\sqrt{n}}\sim\gamma_t$$

Now since the two limiting semicircle laws that we obtain in this way are free, their rescaled sum is circular, in the above sense, and this gives the result.
\end{proof}

Following now Hiai and Petz \cite{hpe}, we have the following free analogue of the PLT:

\index{FPLT}
\index{free PLT}
\index{Marchenko-Pastur law}
\index{free Poisson law}

\begin{theorem}[Free PLT]
We have the following convergence, in moments
$$\left(\left(1-\frac{t}{n}\right)\delta_0+\frac{t}{n}\,\delta_1\right)^{\boxplus n}\to\pi_t$$
the limiting law being the Marchenko-Pastur law of parameter $t>0$, 
$$\pi_t=\max(1-t,0)\delta_0+\frac{\sqrt{4t-(x-1-t)^2}}{2\pi x}\,dx$$
also called free Poisson law of parameter $t$.
\end{theorem}

\begin{proof}
This follows from Theorem 13.20, and from some computations for $\pi_t$:

\medskip

(1) Consider the measure $\eta$ in the statement, appearing under the convolution sign. The Cauchy transform of this measure is easy to compute, and is given by:
$$G_\eta(\xi)=\left(1-\frac{t}{n}\right)\frac{1}{\xi}+\frac{t}{n}\cdot\frac{1}{\xi-1}$$

In order to prove the result, we want to compute the following $R$-transform:
$$R
=R_{\eta^{\boxplus n}}(z)
=nR_\eta(z)$$

According to the formula of $G_\eta$, the equation for this function $R$ is as follows:
$$\left(1-\frac{t}{n}\right)\frac{1}{1/z+R/n}+\frac{t}{n}\cdot\frac{1}{1/z+R/n-1}=z$$

By multiplying both sides by $n/z$, this equation can be written as follows:
$$\frac{t+zR}{1+zR/n}=\frac{t}{1+zR/n-z}$$

With $n\to\infty$ this reads $t+zR=t/(1-z)$, so the limiting $R$-transform is:
$$R(z)=\frac{t}{1-z}$$

(2) Getting now to the $R$-transform of $\pi_t$, as a first observation, things fine at $t=1$, because according to our formulae from chapter 3, we have:
\begin{eqnarray*}
G_{\pi_1}(\xi)=\frac{1-\sqrt{1-4\xi^{-1}}}{2}
&\implies&G\left(\frac{1}{1-z}+\frac{1}{z}\right)=\frac{1-\sqrt{1-4(z-z^2)}}{2}=z\\
&\implies&R_{\pi_1}(z)=\frac{1}{1-z}
\end{eqnarray*}

In order to deal now with the general case, $t>0$, recall from chapter 3 that the law $\pi_t$ appears via Stieltjes inversion from the following Cauchy transform:
$$G_{\pi_t}(\xi)
=\frac{1-t+\xi-\sqrt{(\xi-1-t)^2-4t}}{2\xi}$$

Now observe that, by using this formula, we have the following computation:
\begin{eqnarray*}
&&G_{\pi_t}\left(\frac{t}{1-z}+\frac{1}{z}\right)\\
&=&\frac{z-z^2}{2(1+tz-z)}\left(1-t+\frac{1+tz-z}{z-z^2}
-\sqrt{\left(\frac{1+tz-z}{z-z^2}-1-t\right)^2-4t}\right)\\
&=&\frac{z-z^2}{2(1+tz-z)}\left(\frac{1+tz^2-z^2}{z-z^2}
-\frac{\sqrt{\left(1-2z+(1+t)z^2\right)^2-4t(z-z^2)^2}}{z-z^2}\right)\\
&=&\frac{1+tz^2-z^2-\sqrt{\left(1-(1-\sqrt{t})z\right)^2\left(1-(1+\sqrt{t})z\right)^2}}{2(1+tz-z)}\\
&=&\frac{1+tz^2-z^2-\left(1-(1-\sqrt{t})z\right)\left(1-(1+\sqrt{t})z\right)}{2(1+tz-z)}\\
&=&\frac{1+tz^2-z^2-(1-2z+(1-t)z^2)}{2(1+tz-z)}\\
&=&\frac{2z+(2t-2)z^2}{2(1+tz-z)}\\
&=&z
\end{eqnarray*}

Thus $R_{\pi_t}(z)=t/(1-z)$, and we are led to the conclusion in the statement.
\end{proof}

The above result is quite remarkable, putting an end to the chronic difficulties that we were having in this book, in dealing with the Marchenko-Pastur laws $\pi_t$ at $t\neq1$. Indeed, with our free probability knowledge, we can rewite the story of these laws, as follows:

\begin{theorem}
The Marchenko-Pastur laws, $\pi_t$ with $t>0$, appear from
$$\pi_1=\frac{1}{2\pi}\sqrt{4x^{-1}-1}\,dx$$
which is itself a basic beta law, or the square of $\gamma_1$, via the formula $\pi_t=\pi_1^{\boxplus t}$.
\end{theorem}

\begin{proof}
This is something very simple and beautiful, and self-explanatory, based on Theorem 13.23, or simply on the fact that $R_{\pi_t}(z)=t/(1-z)$ is linear in $t$. By the way, observe that we have as well $\gamma_t=\gamma_1^{\boxplus t}$ and $\Gamma_t=\Gamma_1^{\boxplus t}$, for similar reasons.
\end{proof}

Next, we can talk about compound free Poisson limits and laws, as follows:

\index{compound Poisson law}

\begin{definition}
Associated to any compactly supported positive measure $\nu$ on $\mathbb C$ is the probability measure
$$\pi_\nu=\lim_{n\to\infty}\left(\left(1-\frac{c}{n}\right)\delta_0+\frac{1}{n}\,\nu\right)^{\boxplus n}$$
where $c=mass(\nu)$, called compound free Poisson law.
\end{definition}

In what follows we will be mostly interested in the case where $\nu$ is discrete, as is for instance the case for the measure $\nu=t\delta_1$ with $t>0$, which produces the free Poisson laws. The following result allows one to detect the compound free Poisson laws:

\index{R-transform}

\begin{proposition}
For $\nu=\sum_{i=1}^sc_i\delta_{z_i}$ with $c_i>0$ and $z_i\in\mathbb C$, we have
$$R_{\pi_\nu}(y)=\sum_{i=1}^s\frac{c_iz_i}{1-yz_i}$$
where $R$ denotes as usual the Voiculescu $R$-transform.
\end{proposition}

\begin{proof}
In order to prove this result, let $\eta_n$ be the measure appearing in Definition 13.25, under the free convolution sign. The Cauchy transform of $\eta_n$ is then given by:
$$G_{\eta_n}(\xi)=\left(1-\frac{c}{n}\right)\frac{1}{\xi}+\frac{1}{n}\sum_{i=1}^s\frac{c_i}{\xi-z_i}$$

Consider now the $R$-transform of the measure $\eta_n^{\boxplus n}$, which is given by:
$$R_{\eta_n^{\boxplus n}}(y)=nR_{\eta_n}(y)$$

By using the general theory of the $R$-transform, from Theorem 13.20, the above formula of $G_{\eta_n}$ shows that the equation for $R=R_{\eta_n^{\boxplus n}}$ is as follows:
\begin{eqnarray*}
&&\left(1-\frac{c}{n}\right)\frac{1}{1/y+R/n}+\frac{1}{n}\sum_{i=1}^s\frac{c_i}{1/y+R/n-z_i}=y\\
&\implies&\left(1-\frac{c}{n}\right)\frac{1}{1+yR/n}+\frac{1}{n}\sum_{i=1}^s\frac{c_i}{1+yR/n-yz_i}=1
\end{eqnarray*}

Now multiplying by $n$, then rearranging the terms, and letting $n\to\infty$, we get:
\begin{eqnarray*}
\frac{c+yR}{1+yR/n}=\sum_{i=1}^s\frac{c_i}{1+yR/n-yz_i}
&\implies&c+yR_{\pi_\rho}(y)=\sum_{i=1}^s\frac{c_i}{1-yz_i}\\
&\implies&R_{\pi_\rho}(y)=\sum_{i=1}^s\frac{c_iz_i}{1-yz_i}
\end{eqnarray*}

Thus, we are led to the conclusion in the statement.
\end{proof}

Finally, we have the following result, providing an alternative to Definition 13.25, and which, together with Definition 13.25, can be thought of as being the free CPLT:

\index{CFPLT}
\index{Compound FPLT}

\begin{theorem}[Free CPLT]
For $\nu=\sum_{i=1}^sc_i\delta_{z_i}$ with $c_i>0$ and $z_i\in\mathbb C$, we have
$$\pi_\nu=law\left(\sum_{i=1}^sz_i\alpha_i\right)$$
with the variables $\alpha_i$ being free Poisson$(c_i)$, and free.
\end{theorem}

\begin{proof}
Let $\alpha$ be the sum of free Poisson variables in the statement:
$$\alpha=\sum_{i=1}^sz_i\alpha_i$$

In order to prove the result, we will show that the $R$-transform of $\alpha$ is given by the formula in Proposition 13.26. We have the following computation:
\begin{eqnarray*}
R_{\alpha_i}(y)=\frac{c_i}{1-y}
&\implies&R_{z_i\alpha_i}(y)=\frac{c_iz_i}{1-yz_i}\\
&\implies&R_\alpha(y)=\sum_{i=1}^s\frac{c_iz_i}{1-yz_i}
\end{eqnarray*}

Thus we have the same formula as in Proposition 13.26, and we are done.
\end{proof}

\section*{13d. Random matrices}

We have learned so far the basics of free probability, and obviously, the theory can be indefinitely developed, sky is the limit. Before doing that, however, let us have a bit of thinking, I mean, in relation with concrete questions, is this theory any good?

\bigskip

In answer, let us look at the random matrices. As our first success, we certainly have Theorem 13.24, which provides a more conceptual approach to the laws $\pi_t$ with $t>0$ discovered by Marchenko and Pastur in \cite{mpa}, in relation with the Wishart matrices.

\bigskip

However, the story is far from being over with this, and following Voiculescu \cite{vo3}, we have the following claim, which is obviously first class mathematical discovery:

\begin{claim}
The following happen, in connection with the random matrices:
\begin{enumerate}
\item The Gaussian matrices follow with $N\to\infty$ the Voiculescu circular law $\Gamma_t$.

\item When looking at a family of such matrices, with $N\to\infty$ these become free.

\item The Wigner and Wishart matrices are subject to asymptotic freeness too.
\end{enumerate}
\end{claim}

So, we will discuss this, in the remainder of this chapter. Getting started now, for our computations, we will need explicit models for $\Gamma_t$. Following \cite{vo1}, let us start with:

\index{shift}

\begin{proposition}
Let $H$ be the complex Hilbert space having as basis the colored integers $k=\circ\bullet\bullet\circ\ldots$\,, and consider the shift operators on this space: 
$$S:k\to\circ k\quad,\quad 
T:k\to\bullet k$$
We have then the following equalities of distributions,
$$S+S^*\sim\gamma_1\quad,\quad 
S+T^*\sim\Gamma_1$$
with respect to the state $\varphi(T)=<Te,e>$, where $e$ is the empty word.
\end{proposition}

\begin{proof}
We already know $S+S^*\sim\gamma_1$, from Proposition 13.15, and $S+T^*\sim\Gamma_1$ follows from this, by using the freeness result from Theorem 13.19.
\end{proof}

At the combinatorial level now, we have the following result:

\index{circular variable}
\index{noncrossing pairings}

\begin{theorem}
A variable $a\in A$ follows the law $\Gamma_1$ precisely when
$$tr(a^k)=|\mathcal{NC}_2(k)|$$
for any colored integer $k=\circ\bullet\bullet\circ\ldots$
\end{theorem}

\begin{proof}
By using Proposition 13.29, it is enough to do the computation in the model there. To be more precise, we can use the following explicit formulae for $S,T$:
$$S:k\to\circ k\quad,\quad 
T:k\to\bullet k$$

With these formulae in hand, our claim is that we have the following formula:
$$<(S+T^*)^ke,e>=|\mathcal{NC}_2(k)|$$

In order to prove this formula, let us expand the quantity $(S+T^*)^k$, and then apply the state $\varphi$. With respect to our previous computations, for $S+S^*$, what happens is that the contributions will come this time via the following formulae, which must successively apply, as to collapse the whole product of $S,S^*,T,T^*$ variables into a 1 quantity:
$$S^*S=1\quad,\quad
T^*T=1$$

As before, in the proof for the semicircle laws, these applications of the rules $S^*S=1$, $T^*T=1$ must appear in a noncrossing manner, but what happens now, in contrast with the computation from the proof before, where $S+S^*$ was self-adjoint, is that at each point where the exponent $k$ has a $\circ$ entry we must use $T^*T=1$, and at each point where the exponent $k$ has a $\bullet$ entry we must use $S^*S=1$. Thus the contributions, which are each worth 1, are parametrized by the partitions $\pi\in\mathcal{NC}_2(k)$, as desired.
\end{proof}

More generally now, by rescaling, we have the following result:

\begin{theorem}
A variable $a\in A$ is circular, $a\sim\Gamma_t$, precisely when
$$tr(a^k)=t^{|k|/2}|\mathcal{NC}_2(k)|$$
for any colored integer $k=\circ\bullet\bullet\circ\ldots$
\end{theorem}

\begin{proof}
This follows indeed from Theorem 13.30, by rescaling. Alternatively, we can get this as well directly, by suitably modifying Proposition 13.29 first.
\end{proof}

Even more generally now, we have the following free version of the Wick rule:

\index{free Wick formula}
\index{circular system}

\begin{theorem}
Given free variables $a_i$, each following the Voiculescu circular law $\Gamma_t$, with $t>0$ being a fixed parameter, we have the Wick type formula
$$tr(a_{i_1}^{k_1}\ldots a_{i_s}^{k_s})=t^{s/2}\#\left\{\pi\in\mathcal{NC}_2(k)\Big|\pi\leq\ker i\right\}$$
where $k=k_1\ldots k_s$ and $i=i_1\ldots i_s$, for the joint moments of these variables, with the inequality $\pi\leq\ker i$ on the right being taken in a technical, appropriate sense.
\end{theorem}

\begin{proof}
This follows a bit as in the classical case, the idea being as follows:

\medskip

(1) In the case where we have a single complex normal variable $a$,  we have to compute the moments of $a$, with respect to colored integer exponents $k=\circ\bullet\bullet\circ\ldots\,$, and the formula in the statement coincides with the one in Theorem 13.31, namely:
$$tr(a^k)=t^{|k|/2}|\mathcal{NC}_2(k)|$$

(2) In general now, when expanding the product $a_{i_1}^{k_1}\ldots a_{i_s}^{k_s}$ and rearranging the terms, we are left with doing a number of computations as in (1), and then making the product of the expectations that we found. But this amounts precisely in counting the partitions in the statement, with the condition $\pi\leq\ker i$ there standing precisely for the fact that we are doing the various type (1) computations independently. See \cite{nsp}, \cite{vdn}.
\end{proof}

Getting back now to the case of the single variables, let us record as well:

\begin{theorem}
The moments of the Voiculescu circular law are the numbers
$$M_k(\Gamma_t)=\sum_{\pi\in\mathcal{NC}_2(k)}t^{|\pi|}$$
with $\mathcal{NC}_2$ standing as usual for the noncrossing matching pairings.
\end{theorem}

\begin{proof}
This comes indeed from the formula in Theorem 13.31, which gives:
$$tr(a^k)=t^{|k|/2}|\mathcal{NC}_2(k)|=\sum_{\pi\in\mathcal{NC}_2(k)}t^{|\pi|}$$

Thus, we are led to the conclusion in the statement.
\end{proof}

As an application of this, let us go back to the random matrices. We first have:

\begin{theorem}
Given complex Gaussian matrices $Z_N\in M_N(L^\infty(X))$, having independent $G_t$ variables as entries, with $t>0$, we have
$$\frac{Z_N}{\sqrt{N}}\sim\Gamma_t$$
in the $N\to\infty$ limit, with the limiting measure being Voiculescu's circular law.
\end{theorem}

\begin{proof}
As explained in chapter 8, the asymptotic moments are given by:
$$M_k\left(\frac{Z_N}{\sqrt{N}}\right)\simeq t^{|k|/2}|\mathcal{NC}_2(k)|$$

Now by using Theorem 13.31, we are led to the conclusion in the statement.
\end{proof}

Getting now to asymptotic freeness, following Voiculescu \cite{vo3}, we first have:

\index{Wigner matrix}
\index{asymptotic freeness}

\begin{theorem}
Given a family of sequences of Wigner matrices, 
$$Z^i_N\in M_N(L^\infty(X))\quad,\quad i\in I$$
with pairwise independent entries, each following the complex normal law $G_t$, with $t>0$, up to the constraint $Z_N^i=(Z_N^i)^*$, the rescaled sequences of matrices
$$\frac{Z^i_N}{\sqrt{N}}\in M_N(L^\infty(X))\quad,\quad i\in I$$
become with $N\to\infty$ semicircular, each following the Wigner law $\gamma_t$, and free.
\end{theorem}

\begin{proof}
The first part of the statement is the Wigner theorem from chapter 8, and what is new with respect to that is the ``free'' at the end. Now in order to prove this, asymptotic freeness, we can assume that we are dealing with the case of 2 sequences of matrices, $|I|=2$. So, assume that we have Wigner matrices as follows:
$$Z_N,Z_N'\in M_N(L^\infty(X))$$

We have to prove that these matrices become asymptotically free, with $N\to\infty$. In order to do this, we can use a trick. Consider indeed the following matrix:
$$Y_N=\frac{Z_N+iZ_N'}{\sqrt{2}}$$

This is a complex Gaussian matrix, so by Theorem 13.34 we have, with $N\to\infty$:
$$\frac{Y_N}{\sqrt{N}}\sim\Gamma_t$$

Now recall that the Voiculescu circular law $\Gamma_t$ was by definition the law of the following variable, with $a,b$ being semicircular, each following the law $\gamma_t$, and free:
$$c=\frac{a+ib}{\sqrt{2}}$$

We are therefore in the situation where the variable $(Z_N+iZ_N')/\sqrt{N}$, which has asymptotically semicircular real and imaginary parts, converges to the distribution of $a+ib$, equally having semicircular real and imaginary parts, and with these real and imaginary parts being free. Thus $Z_N,Z_N'$ become asymptotically free, as desired.
\end{proof}

Getting now to the complex case, we have a similar result here, as follows:

\index{Gaussian matrix}
\index{asymptotic freeness}

\begin{theorem}
Given a family of sequences of complex Gaussian matrices, 
$$Z^i_N\in M_N(L^\infty(X))\quad,\quad i\in I$$
with pairwise independent entries, each following the complex normal law $G_t$, with $t>0$, the rescaled sequences of matrices
$$\frac{Z^i_N}{\sqrt{N}}\in M_N(L^\infty(X))\quad,\quad i\in I$$
become with $N\to\infty$ circular, each following the Voiculescu law $\Gamma_t$, and free.
\end{theorem}

\begin{proof}
This follows indeed from Theorem 13.35, which applies to the real and imaginary parts of our complex Gaussian matrices, and gives the result.
\end{proof}

Finally, we have as well an asymptotic freeness result for the Wishart matrices, coming again as a consequence of Theorem 13.35, and completing the proof of Claim 13.28.

\section*{13e. Exercises}

Welcome to freeness, such a pleasure having you here, and as exercises, we have:

\begin{exercise}
Learn more $C^*$-algebra basics, including the GNS theorem.
\end{exercise}

\begin{exercise}
Have a look as well at the von Neumann algebra theory.
\end{exercise}

\begin{exercise}
Further meditate on the multiplicative convolution operation $\boxtimes$.
\end{exercise}

\begin{exercise}
Find a proof for the $R$-transform, by staying in the tracial setting.
\end{exercise}

\begin{exercise}
Compute, with bare hands, the Cauchy transform of $\pi_t$.
\end{exercise}

\begin{exercise}
Study a bit the free compound Poisson laws that we constructed.
\end{exercise}

\begin{exercise}
Clarify what we said, in regards with the free Wick formula.
\end{exercise}

\begin{exercise}
Do some joint moment computations, for random matrices.
\end{exercise}

As bonus exercise, and no surprise here, read the originals \cite{vo1}, \cite{vo2}, \cite{vo3}.

\chapter{The bijection}

\section*{14a. Cumulants}

In order to further advance in our study, we need to better understand the correspondence between classical and free. At the general level this comes of course from ``independence gets replaced by freeness''. But this remains quite vague, I mean go get the formula of $\pi_t$ out of the formula of $p_t$, based on this, this is no easy task.

\bigskip

So, let us have a closer look at what we did in the previous chapter, all the formulae there, regarding $\gamma_t,\Gamma_t,\pi_t,\pi_\nu$, and compare them with our previous formulae from this book, regarding their classical counterparts $g_t,G_t,p_t,p_\nu$. We are led in this way to:

\begin{speculation}
The correspondence $m_t\leftrightarrow\mu_t$ between the classical and free can be understood in several possible ways, as follows,
\begin{enumerate}
\item Basic guideline: ``independence gets replaced by freeness''.

\item Same thing, sounding better: ``the operation $*$ gets replaced by $\boxplus$''.

\item Moments, we have here the following connecting formula, with the categories $D\subset P$ and $D'\subset NC$ being related by $D'=D\cap NC$ and $D=<D',\slash\hskip-2.1mm\backslash>$:
$$M_k(m_t)=\sum_{\pi\in D(k)}t^{|\pi|}\quad \longleftrightarrow\quad M_k(\mu_t)=\sum_{\pi\in D'(k)}t^{|\pi|}$$

\item Functional transforms, we have here the following connecting formula:
$$\log F_{m_t}(-iz)=\sum_{n=1}^\infty c_n\cdot\frac{z^n}{n!}
\quad \longleftrightarrow\quad zR_{\mu_t}(z)=\sum_{n=1}^\infty c_n\cdot z^n$$
\end{enumerate}
and there are certainly more ways, by looking at orthogonal polynomials, Gram or Hankel determinants, or Weingarten functions. The only wrong way is via densities.
\end{speculation}

To be more precise here, (1) and (2) are from Mao's red book, (3) is something that we know well for $g_t\leftrightarrow\gamma_t$, $G_t\leftrightarrow\Gamma_t$, $p_t\leftrightarrow\pi_t$, and I will leave to you as an exercise to establish this for $p^s_t\leftrightarrow\pi^s_t$ too, with the obvious definition for $\pi^s_t$, but we will be back to this with soon, and the assertions at the end are a quote from comrade Stalin.

\bigskip

As for (4), this is something new, coming by closely examining the $R$-transform formulae from the previous chapter. Indeed, for the Poisson laws, this certainly holds:
$$\log F_{p_t}(-iz)=(e^z-1)t=\sum_{n=1}^\infty t\cdot\frac{z^n}{n!}$$
$$zR_{\pi_t}(z)=\frac{zt}{1-z}=\sum_{n=1}^\infty t\cdot z^n$$

Observe that by linearity this holds, more generally, for the compound Poisson laws. Next, for the Gaussian laws this holds too, due to more trivial reasons, as follows:
$$\log F_{g_t}(-iz)=\frac{z^2t}{2}=t\cdot\frac{z^2}{2!}$$
$$zR_{\gamma_t}(z)=z^2t=t\cdot z^2$$

As yet another piece of evidence, for the Dirac mass $\delta_t$, which reigns over probability theory at large, being at the same time classical and free, this works too:
$$\log F_{\delta_t}(-iz)=zt=t\cdot\frac{z}{1!}$$
$$zR_{\delta_t}(z)=zt=t\cdot z$$

Summarizing, we have our speculation up and running, and the problem is now, which of the 2 main ways there shall we follow. In answer, we will follow a middle path:

\begin{principle}
The corrrespondence betwen classical and free is best understood via quantities $k_n,\kappa_n$, called cumulants and free cumulants, appearing as follows:
$$\log F_m(-iz)=\sum_{n=1}^\infty k_n(m)\cdot\frac{z^n}{n!}
\quad,\quad 
zR_\mu(z)=\sum_{n=1}^\infty \kappa_n(\mu)\cdot z^n$$
Indeed, we can say that we have $m\leftrightarrow\mu$ when $k_n(m)=\kappa_n(\mu)$. And with some further combinatorial work, this will clarify our correspondence for the moments too. 
\end{principle} 

So, this is the situation, hope you got it, what we have initially suggests jumping into heavy analysis, by comparing $\log F_m$ and $R_\mu$, and this is indeed what Bercovici and Pata did in \cite{bpa}, when first examining the correspondence question. However, from a more modern perspective, it is actually better to replace analysis by combinatorics, by comparing instead the coefficients of $\log F_m$ and $R_\mu$, with this being equivalent to what Bercovici-Pata did, but coming with a bonus, namely moment results too.

\bigskip

Getting to work now, we first need to talk about cumulants $k_n$. And here, forgetting all the above, and everything advanced in general, and coming as a continuation of what we did in chapter 1, following Rota, let us formulate the following key definition:

\index{cumulant}
\index{cumulant-generating function}
\index{Taylor coefficient}
\index{Fourier transform}
\index{generating series}

\begin{definition}
Associated to any real probability measure $\mu=\mu_f$ is the following modification of the logarithm of the Fourier transform $F_\mu(z)=E(e^{izf})$,
$$K_\mu(z)=\log E(e^{zf})$$
called cumulant-generating function. The Taylor coefficients $k_n(\mu)$ of this series, given by
$$K_\mu(z)=\sum_{n=1}^\infty k_n(\mu)\,\frac{z^n}{n!}$$
are called cumulants of the measure $\mu$. We also use the notations $k_f,K_f$ for these cumulants and their generating series, where $f$ is a variable following the law $\mu$.
\end{definition}

In other words, the cumulants are more or less the coefficients of the logarithm of the Fourier transform $\log F_\mu$, up to some normalizations. To be more precise, we have $K_\mu(z)=\log F_\mu(-iz)$, so the formula relating $\log F_\mu$ to the cumulants $k_n(\mu)$ is:
$$\log F_\mu(-iz)=\sum_{n=1}^\infty k_n(\mu)\,\frac{z^n}{n!}$$

We will see in a moment the reasons for the above normalizations, namely change of variables $z\to -iz$, and Taylor coefficients instead of plain coefficients, the idea being that for simple laws like $p_t,g_t$, we will obtain in this way very simple quantities.

\bigskip

As a first observation, the sequence of cumulants $k_1,k_2,k_3,\ldots$ appears as a modification of the sequence of moments $M_1,M_2,M_3,\ldots\,$, the numerics being as follows:

\index{sequence of moments}
\index{sequence of cumulants}

\begin{proposition}
The sequence of cumulants $k_1,k_2,k_3,\ldots$ appears as a modification of the sequence of moments $M_1,M_2,M_3,\ldots\,$, and uniquely determines $\mu$. We have
$$k_1=M_1$$
$$k_2=-M_1^2+M_2$$
$$k_3=2M_1^3-3M_1M_2+M_3$$
$$k_4=-6M_1^4+12M_1^2M_2-3M_2^2-4M_1M_3+M_4$$
$$\vdots$$
in one sense, and in the other sense we have
$$M_1=k_1$$
$$M_2=k_1^2+k_2$$
$$M_3=k_1^3+3k_1k_2+k_3$$
$$M_4=k_1^4+6k_1^2k_2+3k_2^2+4k_1k_3+k_4$$
$$\vdots$$
with in both cases the correspondence being polynomial, with integer coefficients.
\end{proposition}

\begin{proof}
We know from Definition 14.3 that the cumulants are given by:
$$\log E(e^{zf})=\sum_{s=1}^\infty k_s(f)\,\frac{z^s}{s!}$$

By exponentiating, we obtain from this the following formula:
$$E(e^{zf})=\exp\left(\sum_{s=1}^\infty k_s(f)\,\frac{z^s}{s!}\right)$$

Now by looking at the terms of order $1,2,3,4$, this gives the above formulae.
\end{proof}

The interest in cumulants comes from the fact that $\log F_\mu$, and so the cumulants $k_n(\mu)$ too, linearize the convolution. To be more precise, we have the following result:

\index{linearization of convolution}
\index{additivity of cumulants}

\begin{theorem}
The cumulants have the following properties:
\begin{enumerate}
\item $k_n(cf)=c^nk_n(f)$.

\item $k_1(f+d)=k_1(f)+d$, and $k_n(f+d)=k_n(f)$ for $n>1$.

\item $k_n(f+g)=k_n(f)+k_n(g)$, if $f,g$ are independent.
\end{enumerate}
\end{theorem}

\begin{proof}
Here (1) and (2) are both clear from definitions, because we have:
\begin{eqnarray*}
K_{cf+d}(z)
&=&\log E(e^{z(cf+d)})\\
&=&\log[e^{zd}\cdot E(e^{zcf})]\\
&=&zd+K_f(cz)
\end{eqnarray*}

As for (3), this follows from the fact that the Fourier transform $F_f(x)=E(e^{ixf})$ satisfies the following formula, whenever $f,g$ are independent random variables:
$$F_{f+g}(x)=F_f(x)F_g(x)$$

Indeed, by applying the logarithm, we obtain the following formula:
$$\log F_{f+g}(x)=\log F_f(x)+\log F_g(x)$$

With the change of variables $x=-iz$, we obtain the following formula:
$$K_{f+g}(z)=K_f(z)+K_g(z)$$

Thus, at the level of coefficients, we obtain $k_n(f+g)=k_n(f)+k_n(g)$, as claimed.
\end{proof}

At the level of basic examples now, we have the following result:

\begin{theorem}
The sequence of cumulants $k_1,k_2,k_3,k_4,\ldots$ is as follows:
\begin{enumerate}
\item For $\mu=\delta_t$ the cumulants are $t,0,0,0,\ldots$

\item For $\mu=g_t$ the cumulants are $0,t,0,0,\ldots$

\item For $\mu=p_t$ the cumulants are $t,t,t,t,\ldots$

\item For $\mu=p_t^2$ the cumulants are $0,t,0,t,\ldots$
\end{enumerate}
Also, for the compound Poisson laws the cumulants are $k_n(p_\nu)=M_n(\nu)$.
\end{theorem}

\begin{proof}
We have 5 computations to be done, the idea being as follows:

\medskip

(1) For a Dirac mass, $\mu=\delta_t$, this comes from the following computation:
$$K_\mu(z)
=\log E(e^{tz})
=\log(e^{tz})
=tz$$

(2) For a normal law, $\mu=g_t$, this comes from the following computation:
$$K_\mu(z)
=\log F_\mu(-iz)
=\log\exp\left[-t(-iz)^2/2\right]
=tz^2/2$$

(3) For a Poisson law, $\mu=p_t$, this comes from the following computation:
$$K_\mu(z)
=\log F_\mu(-iz)
=\log\exp\left[(e^{i(-iz)}-1)t\right]
=(e^z-1)t$$

(4) For a Bessel law, $\mu=p_t^2$, this comes from the following computation:
$$K_\mu(z)
=\log F_\mu(-iz)
=\log\exp\left[\left(\frac{e^z+e^{-z}}{2}-1\right)t\right]
=\left(\frac{e^z+e^{-z}}{2}-1\right)t$$

(5) In order to prove the last assertion, generalizing (3,4), by continuity we can assume that our input measure is discrete, $\nu=\sum_i t_i\delta_{z_i}$ with $t_i>0$ and $z_i\in\mathbb R$. We have:
\begin{eqnarray*}
K_{p_\nu}(y)
&=&\log\exp\left[\sum_it_i(e^{yz_i}-1)\right]\\
&=&\sum_it_i\sum_{n\geq1}\frac{(yz_i)^n}{n!}\\
&=&\sum_{n\geq1}\frac{y^n}{n!}\sum_it_iz_i^n\\
&=&\sum_{n\geq1}\frac{y^n}{n!}\,M_n(\nu)
\end{eqnarray*}

Thus, we are led to the conclusion in the statement.
\end{proof}

Getting back to theory now, the sequence of cumulants $k_1,k_2,k_3,\ldots$ appears as a modification of the sequence of moments $M_1,M_2,M_3,\ldots\,$, and understanding the relation between moments and cumulants will be our next task. Let us start with:

\index{M\"obius function}
\index{lattice}
\index{lattice of partitions}

\begin{definition}
The M\"obius function of any lattice, and so of $P$, is given by
$$\mu(\pi,\nu)=\begin{cases}
1&{\rm if}\ \pi=\nu\\
-\sum_{\pi\leq\tau<\nu}\mu(\pi,\tau)&{\rm if}\ \pi<\nu\\
0&{\rm if}\ \pi\not\leq\nu
\end{cases}$$
with the construction being performed by recurrence.
\end{definition}

As an illustration here, for $P(2)=\{||,\sqcap\}$, we have by definition:
$$\mu(||,||)=\mu(\sqcap,\sqcap)=1$$

Also, $||<\sqcap$, with no intermediate partition in between, so we obtain:
$$\mu(||,\sqcap)=-\mu(||,||)=-1$$

Finally, we have $\sqcap\not\leq||$, and so we have as well the following formula:
$$\mu(\sqcap,||)=0$$

Thus, the M\"obius matrix $M_{\pi\nu}=\mu(\pi,\nu)$ of the lattice $P(2)=\{||,\sqcap\}$ is as follows:
$$M=\begin{pmatrix}1&-1\\ 0&1\end{pmatrix}$$

Back to the general case now, the main interest in the M\"obius function comes from the M\"obius inversion formula, which can be formulated as follows:

\index{M\"obius inversion}
\index{M\"obius matrix}

\begin{theorem}
We have the following implication,
$$f(\pi)=\sum_{\nu\leq\pi}g(\nu)
\quad\implies\quad
g(\pi)=\sum_{\nu\leq\pi}\mu(\nu,\pi)f(\nu)$$
valid for any two functions $f,g:P(n)\to\mathbb C$.
\end{theorem}

\begin{proof}
Consider the adjacency matrix of $P$, given by the following formula:
$$A_{\pi\nu}=\begin{cases}
1&{\rm if}\ \pi\leq\nu\\
0&{\rm if}\ \pi\not\leq\nu
\end{cases}$$

Our claim is that the inverse of this matrix is the M\"obius matrix of $P$, given by:
$$M_{\pi\nu}=\mu(\pi,\nu)$$

Indeed, the above matrix $A$ is upper triangular, and when trying to invert it, we are led to the recurrence in Definition 14.7, so to the M\"obius matrix $M$. Thus we have:
$$M=A^{-1}$$

Thus, in practice, we are led to the inversion formula in the statement.
\end{proof}

As an illustration here, for $P(2)$ the formula $M=A^{-1}$ appears as follows:
$$\begin{pmatrix}1&-1\\ 0&1\end{pmatrix}=
\begin{pmatrix}1&1\\ 0&1\end{pmatrix}^{-1}$$

With these ingredients in hand, let us go back to probability. We first have:

\index{cumulant}
\index{classical cumulant}
\index{generalized cumulant}
\index{multiplicativity over blocks}

\begin{definition}
We define quantities $M_\pi(f),k_\pi(f)$, depending on partitions 
$$\pi\in P(k)$$
by starting with $M_n(f),k_n(f)$, and using multiplicativity over the blocks. 
\end{definition}

To be more precise, the convention here is that for the one-block partition $1_n\in P(n)$, the corresponding moment and cumulant are the usual ones, namely:
$$M_{1_n}(f)=M_n(f)\quad,\quad k_{1_n}(f)=k_n(f)$$

Then, for an arbitrary partition $\pi\in P(k)$, we decompose this partition into blocks, having sizes $b_1,\ldots,b_s$, and we set, by multiplicativity over blocks:
$$M_\pi(f)=M_{b_1}(f)\ldots M_{b_s}(f)\quad,\quad k_\pi(f)=k_{b_1}(f)\ldots k_{b_s}(f)$$

With this convention, following Rota and others, we can now formulate a key result, fully clarifying the relation between moments and cumulants, as follows:

\index{moment-cumulant formula}

\begin{theorem}
We have the moment-cumulant formulae
$$M_n(f)=\sum_{\nu\in P(n)}k_\nu(f)\quad,\quad 
k_n(f)=\sum_{\nu\in P(n)}\mu(\nu,1_n)M_\nu(f)$$
or, equivalently, we have the moment-cumulant formulae
$$M_\pi(f)=\sum_{\nu\leq\pi}k_\nu(f)\quad,\quad 
k_\pi(f)=\sum_{\nu\leq\pi}\mu(\nu,\pi)M_\nu(f)$$
where $\mu$ is the M\"obius function of $P(n)$.
\end{theorem}

\begin{proof}
There are several things going on here, the idea being as follows:

\medskip

(1) According to our conventions, and to M\"obius inversion, the 4 formulae in the statement are all equivalent. In what follows we will focus on the first 2 formulae.

\medskip

(2) Before anything, let us first work out some examples. At $n=1,2,3$ the moment formula in the statement gives, in tune with Proposition 14.4:
$$M_1=k_|=k_1$$
$$M_2=k_{|\,|}+k_\sqcap=k_1^2+k_2$$
$$M_3=k_{|\,|\,|}+k_{\sqcap|}+k_{\sqcap\hskip-2.8mm{\ }_|}+k_{|\sqcap}+k_{\sqcap\hskip-0.5mm\sqcap}=k_1^3+3k_1k_2+k_3$$

At $n=4$ now, the computation is as follows, again in tune with Proposition 14.4:
\begin{eqnarray*}
M_4
&=&k_{|\,|\,|}+(\underbrace{k_{\sqcap\,|\,|}+\ldots}_{6\ terms})+(\underbrace{k_{\sqcap\,\sqcap}+\ldots}_{3\ terms})+(\underbrace{k_{\sqcap\hskip-0.5mm\sqcap\,|}+\ldots}_{4\ terms})+k_{\sqcap\hskip-0.5mm\sqcap\hskip-0.5mm\sqcap}\\
&=&k_1^4+6k_1^2k_2+3k_2^2+4k_1k_3+k_4
\end{eqnarray*}

As for the cumulant formula in the statement, at $n=1,2,3$ this gives the following formulae for the low order cumulants, again in tune with Proposition 14.4:
$$k_1=M_|=M_1$$
$$k_2=(-1)M_{|\,|}+M_\sqcap=-M_1^2+M_2$$
$$k_3=2M_{|\,|\,|}+(-1)M_{\sqcap|}+(-1)M_{\sqcap\hskip-2.8mm{\ }_|}+(-1)M_{|\sqcap}+M_{\sqcap\hskip-0.5mm\sqcap}=2M_1^3-3M_1M_2+M_3$$

Finally, at $n=4$, after computing the M\"obius function of $P(4)$, we obtain the following formula for the fourth cumulant, again in tune with Proposition 14.4:
\begin{eqnarray*}
k_4
&=&(-6)M_{|\,|\,|}+2(\underbrace{M_{\sqcap\,|\,|}+\ldots}_{6\ terms})+(-1)(\underbrace{M_{\sqcap\,\sqcap}+\ldots}_{3\ terms})+(-1)(\underbrace{M_{\sqcap\hskip-0.5mm\sqcap\,|}+\ldots}_{4\ terms})+M_{\sqcap\hskip-0.5mm\sqcap\hskip-0.5mm\sqcap}\\
&=&-6M_1^4+12M_1^2M_2-3M_2^2-4M_1M_3+M_4
\end{eqnarray*}

(3) Time now to get to work, and prove the result. As mentioned above, the formulae in the statement are all equivalent, and it is enough to prove the first one, namely:
$$M_n(f)=\sum_{\nu\in P(n)}k_\nu(f)$$

In order to do this, we use the very definition of the cumulants, namely:
$$\log E(e^{zf})=\sum_{s=1}^\infty k_s(f)\,\frac{z^s}{s!}$$

By exponentiating, we obtain from this the following formula:
$$E(e^{zf})=\exp\left(\sum_{s=1}^\infty k_s(f)\,\frac{z^s}{s!}\right)$$

(4) Let us first compute the function on the left. This is easily done, as follows:
$$E(e^{zf})
=E\left(\sum_{n=0}^\infty\frac{(zf)^n}{n!}\right)
=\sum_{n=0}^\infty M_n(f)\,\frac{z^n}{n!}$$

(5) Regarding now the function on the right, this is given by:
\begin{eqnarray*}
\exp\left(\sum_{s=1}^\infty k_s(f)\,\frac{z^s}{s!}\right)
&=&\sum_{p=0}^\infty\frac{\left(\sum_{s=1}^\infty k_s(f)\,\frac{z^s}{s!}\right)^p}{p!}\\
&=&\sum_{p=0}^\infty\frac{1}{p!}\sum_{s_1=1}^\infty k_{s_1}(f)\,\frac{z^{s_1}}{s_1!}\ldots\ldots\sum_{s_p=1}^\infty k_{s_p}(f)\,\frac{z^{s_p}}{s_p!}\\
&=&\sum_{p=0}^\infty\frac{1}{p!}\sum_{s_1=1}^\infty\ldots\sum_{s_p=1}^\infty k_{s_1}(f)\ldots k_{s_p}(f)\,\frac{z^{s_1+\ldots+s_p}}{s_1!\ldots s_p!}
\end{eqnarray*}

But the point now is that all this leads us into partitions. Indeed, we are summing over indices $s_1,\ldots,s_p\in\mathbb N$, which can be thought of as corresponding to a partition of $n=s_1+\ldots+s_p$. So, let us rewrite our sum, as a sum over partitions. For this purpose, recall that the number of partitions $\nu\in P(n)$ having blocks of sizes $s_1,\ldots,s_p$ is:
$$\binom{n}{s_1,\ldots,s_p}=\frac{n!}{p_1!\ldots p_s!}$$

Also, when resumming over partitions, there will be a $p!$ factor as well, coming from the permutations of $s_1,\ldots,s_p$. Thus, our sum can be rewritten as follows:
\begin{eqnarray*}
\exp\left(\sum_{s=1}^\infty k_s(f)\,\frac{z^s}{s!}\right)
&=&\sum_{n=0}^\infty\sum_{p=0}^\infty\frac{1}{p!}\sum_{s_1+\ldots+s_p=n}k_{s_1}(f)\ldots k_{s_p}(f)\,\frac{z^n}{s_1!\ldots s_p!}\\
&=&\sum_{n=0}^\infty\frac{z^n}{n!}\sum_{p=0}^\infty\frac{1}{p!}\sum_{s_1+\ldots+s_p=n}\binom{n}{s_1,\ldots,s_p}k_{s_1}(f)\ldots k_{s_p}(f)\\
&=&\sum_{n=0}^\infty\frac{z^n}{n!}\sum_{\nu\in P(n)}k_\nu(f)
\end{eqnarray*}

(6) We are now in position to conclude. According to (3,4,5), we have:
$$\sum_{n=0}^\infty M_n(f)\,\frac{z^n}{n!}=\sum_{n=0}^\infty\frac{z^n}{n!}\sum_{\nu\in P(n)}k_\nu(f)$$

Thus, we have the following formula, valid for any $n\in\mathbb N$:
$$M_n(f)=\sum_{\nu\in P(n)}k_\nu(f)$$

We are therefore led to the conclusions in the statement.
\end{proof}

The above cumulant technology is quite powerful, and can be applied as well to other measures that we know. For the exponential law $e_t$, the cumulants are:
$$k_n=\frac{(n-1)!}{t^n}$$

For the uniform law $u_t$ the cumulants are quite tricky to compute, and are given by the following formula, with $B_n\in\mathbb Q$ being the $n$-th Bernoulli number:
$$k_1=\frac{t}{2}\quad,\quad k_n=\frac{t^n}{n}\cdot B_n\ (n>1)$$

Finally, there is also a discussion to be made, in relation with the notion of infinite divisibility of measures. We will leave some learning here as an exercise.

\section*{14b. Free cumulants}

Getting back now to free probability, following Speicher \cite{sp1}, let us formulate:

\index{cumulant}
\index{free cumulant}

\begin{definition}
The free cumulants $\kappa_n(a)$ of a variable $a\in A$ are defined by
$$R_a(z)=\sum_{n=1}^\infty\kappa_n(a)z^{n-1}$$
with $R$ being as usual the Voiculescu $R$-transform.
\end{definition}

As before with the classical cumulants, we have a number of basic examples and illustrations, and a number of basic general results. Let us start with some numerics:

\begin{proposition}
The free cumulants $\kappa_1,\kappa_2,\kappa_3,\ldots$ appear as a modification of the moments $M_1,M_2,M_3,\ldots\,$, and uniquely determine $\mu$. We have
$$\kappa_1=M_1$$
$$\kappa_2=-M_1^2+M_2$$
$$\kappa_3=2M_1^3-3M_1M_2+M_3$$
$$\kappa_4=-5M_1^4+10M_1^2M_2-2M_2^2-4M_1M_3+M_4$$
$$\vdots$$
in one sense, and in the other sense we have
$$M_1=\kappa_1$$
$$M_2=\kappa_1^2+\kappa_2$$
$$M_3=\kappa_1^3+3\kappa_1\kappa_2+\kappa_3$$
$$M_4=\kappa_1^4+6\kappa_1^2\kappa_2+2\kappa_2^2+4\kappa_1\kappa_3+\kappa_4$$
$$\vdots$$
with in both cases the correspondence being polynomial, with integer coefficients.
\end{proposition}

\begin{proof}
This is something quite instructive, a bit more complicated than the similar result in the classical case, worth discussing in detail, the idea being as follows:

\medskip

(1) We know that the Cauchy transform is the following function:
$$G(\xi)=\sum_{n=0}^\infty\frac{M_n}{\xi^{n+1}}$$

By inverting, $K(G(\xi))=\xi$, the numbers $\kappa_n$ appear as follows, with $G=G(\xi)$:
$$\frac{1}{G}+\sum_{n=1}^\infty\kappa_nG^{n-1}=\xi$$

With $\xi=z^{-1}$ our equation takes the following form, with $G=G(z^{-1})$:
$$\frac{z}{G}+z\sum_{n=1}^\infty\kappa_nG^{n-1}=1$$

With $\psi=\sum_{n=1}^\infty M_nz^n$ we have $G/z=1+\psi$, and our equation becomes:
$$\frac{1}{1+\psi}+\sum_{n=1}^\infty\kappa_nz^n(1+\psi)^{n-1}=1$$

(2) By expanding the fraction on the left, our equation becomes:
$$\sum_{n=0}^\infty(-\psi)^n+\sum_{n=1}^\infty\kappa_nz^n(1+\psi)^{n-1}=1$$

Moreover, we can cancel the 1 term on both sides, and our equation becomes:
$$\sum_{n=1}^\infty(-\psi)^n+\sum_{n=1}^\infty\kappa_nz^n(1+\psi)^{n-1}=0$$

Alternatively, we can write our equation in the following way:
$$\sum_{n=1}^\infty\kappa_nz^n(1+\psi)^{n-1}=-\sum_{n=1}^\infty(-\psi)^n$$

(3) Good news, this latter equation is something that we are eventually happy with. By remembering that we have $\psi=\sum_{n=1}^\infty M_nz^n$, our equation looks as follows:
\begin{eqnarray*}
&&\kappa_1z+\kappa_2z^2(1+M_1z+M_2z^2+\ldots)+\kappa_3z^3(1+M_1z+M_2z^2+\ldots)^2+\ldots\\
&=&(M_1z+M_2z^2+\ldots)-(M_1z+M_2z^2+\ldots)^2+(M_1z+M_2z^2+\ldots)^3-\ldots
\end{eqnarray*}

Which looks nice and exploitable, and with this, we are basically done.

\medskip

(4) Indeed, by looking at terms of order $1,2,3,4$, we have the following equations:
$$\kappa_1=M_1$$
$$\kappa_2=M_2-M_1^2$$
$$\kappa_2M_1+\kappa_3=M_3-2M_1M_2+M_1^3$$
$$\kappa_4+2\kappa_3M_1+\kappa_2M_2=M_4-2M_1M_3-M_2^2+3M_1^2M_2-M_1^4$$

Thus, we are led to the formulae of $\kappa_1,\kappa_2,\kappa_3,\kappa_4$ in the statement. And then, by inversion, to the formulae of $M_1,M_2,M_3,M_4$ as well, as desired. 
\end{proof}

Observe the similarity with the formulae for classical cumulants. In fact, we have:

\begin{conclusion}
The first three classical and free cumulants coincide, but
$$k_4=-6M_1^4+12M_1^2M_2-3M_2^2-4M_1M_3+M_4$$
$$\kappa_4=-5M_1^4+10M_1^2M_2-2M_2^2-4M_1M_3+M_4$$
are different, and the same happens at higher order as well.
\end{conclusion}

This is something quite interesting, and we will back later with a conceptual explanation for this, via partitions, the idea being that all this comes from:
$$P(n)=NC(n)\iff n\leq 3$$

But more on this later. At the level of basic general results, in analogy with what we know about the classical cumulants, we first have the following result:

\begin{theorem}
The free cumulants have the following properties:
\begin{enumerate}
\item $\kappa_n(\lambda a)=\lambda^n\kappa_n(a)$.

\item $\kappa_n(a+b)=\kappa_n(a)+\kappa_n(b)$, if $a,b$ are free.
\end{enumerate}
\end{theorem}

\begin{proof}
This is something very standard, the idea being as follows:

\medskip

(1) We have the following Cauchy transform computation:
$$G_{\lambda a}(\xi)
=\int_\mathbb R\frac{d\mu_{\lambda a}(x)}{\xi-x}
=\int_\mathbb R\frac{d\mu_a(y)}{\xi-\lambda y}
=\frac{1}{\lambda}\int_\mathbb R\frac{d\mu_a(y)}{\xi/\lambda-y}
=\frac{1}{\lambda}\,G_a\left(\frac{\xi}{\lambda}\right)$$

But this gives the following formula, by using the definition of the $R$-transform:
$$G_{\lambda a}\left(\lambda R_a(\lambda z)+\frac{1}{z}\right)
=\frac{1}{\lambda}\,G_a\left(R_a(\lambda z)+\frac{1}{\lambda z}\right)
=\frac{1}{\lambda}\cdot\lambda z
=z$$

Thus we have the formula $R_{\lambda a}(z)=\lambda R_a(\lambda z)$, which gives (1). 

\medskip

(2) This follows indeed from the fact that the $R$-transform linearizes $\boxplus$.
\end{proof}

Again in analogy with the classical case, at the level of the main examples we have the following result, with $\pi^2_t=\pi_{t\varepsilon}$ with $\varepsilon=(\delta_{-1}+\delta_1)/2$ being a free Bessel law:

\begin{theorem}
The sequence of free cumulants $\kappa_1,\kappa_2,\kappa_3,\kappa_4,\ldots$ is as follows:
\begin{enumerate}
\item For $\mu=\delta_t$ the free cumulants are $t,0,0,0,\ldots$

\item For $\mu=\gamma_t$ the free cumulants are $0,t,0,0,\ldots$

\item For $\mu=\pi_t$ the free cumulants are $t,t,t,t,\ldots$

\item For $\mu=\pi_t^2$ the free cumulants are $0,t,0,t,\ldots$
\end{enumerate}
Also, for the compound free Poisson laws the free cumulants are $\kappa_n(\pi_\nu)=M_n(\nu)$.
\end{theorem}

\begin{proof}
The proofs are analogous to those from the classical case, as follows:

\medskip

(1) For $\mu=\delta_t$ we have $G_\mu(\xi)=1/(\xi-t)$, and so $R_\mu(z)=t$, as desired.

\medskip

(2) For $\mu=\gamma_t$ we have, as computed before, $R_\mu(z)=tz$, as desired.

\medskip

(3) For $\mu=\pi_t$ we have, also from before, $R_\mu(z)=t/(1-z)$, as desired.

\medskip

(4) For $\mu=\pi_t^2$, a free Bessel law, this can be established too, but the best is to prove directly the last assertion, which generalizes (3,4). With $\nu=\sum_ic_i\delta_{z_i}$ we have:
\begin{eqnarray*}
R_{\pi_\nu}(y)
&=&\sum_i\frac{c_iz_i}{1-yz_i}\\
&=&\sum_{n\geq1}y^{n-1}\sum_ic_iz_i^n\\
&=&\sum_{n\geq 1}y^{n-1}\,M_n(\nu)
\end{eqnarray*}

Thus, we are led to the conclusion in the statement.
\end{proof}

As before in the classical case, we can define generalized free cumulants, $\kappa_\pi(a)$ with $\pi\in P(k)$, by starting with the numeric free cumulants $\kappa_n(a)$, as follows:

\begin{definition}
We define free cumulants $\kappa_\pi(a)$, depending on partitions 
$$\pi\in P(k)$$
by starting with $\kappa_n(a)$, and using multiplicativity over the blocks. 
\end{definition}

To be more precise, the convention here is that for the one-block partition $1_n\in P(n)$, the corresponding free cumulant is the usual one, namely:
$$\kappa_{1_n}(a)=\kappa_n(a)$$

Then, for an arbitrary partition $\pi\in P(k)$, we decompose this partition into blocks, having sizes $b_1,\ldots,b_s$, and we set, by multiplicativity over blocks:
$$\kappa_\pi(a)=\kappa_{b_1}(a)\ldots\kappa_{b_s}(a)$$

With this convention, we have the following result, due to Speicher \cite{sp1}:

\begin{theorem}
We have the moment-cumulant formulae
$$M_n(a)=\sum_{\nu\in NC(n)}\kappa_\nu(a)\quad,\quad 
\kappa_n(a)=\sum_{\nu\in NC(n)}\mu(\nu,1_n)M_\nu(a)$$
or, equivalently, we have the moment-cumulant formulae
$$M_\pi(a)=\sum_{\nu\leq\pi}\kappa_\nu(a)\quad,\quad 
\kappa_\pi(a)=\sum_{\nu\leq\pi}\mu(\nu,\pi)M_\nu(a)$$
where $\mu$ is the M\"obius function of $NC(n)$.
\end{theorem}

\begin{proof}
As before in the classical case, the 4 formulae in the statement are equivalent, via M\"obius inversion. Thus, it is enough to prove one of them, and we will discuss the first formula, which in practice is the most useful one. Thus, we must prove that:
$$M_n(a)=\sum_{\nu\in NC(n)}\kappa_\nu(a)$$

(1) In order to prove this, let us get back to the construction of the free cumulants, from Definition 14.11. The Cauchy transform of $a$ is the following function:
$$G_a(\xi)=\sum_{n=0}^\infty\frac{M_n(a)}{\xi^{n+1}}$$

Consider the inverse of this Cauchy transform $G_a$, with respect to composition:
$$K_a(G_a(\xi))=\xi$$

According to Definition 14.11, the free cumulants $\kappa_n(a)$ appear then as follows:
$$K_a(z)=\frac{1}{z}+\sum_{n=1}^\infty\kappa_n(a)z^{n-1}$$

(2) In practice, as seen in the proof of Proposition 14.12, this formula leads to:
\begin{eqnarray*}
&&\kappa_1z+\kappa_2z^2(1+M_1z+M_2z^2+\ldots)+\kappa_3z^3(1+M_1z+M_2z^2+\ldots)^2+\ldots\\
&=&(M_1z+M_2z^2+\ldots)-(M_1z+M_2z^2+\ldots)^2+(M_1z+M_2z^2+\ldots)^3-\ldots
\end{eqnarray*}

But, in case you have followed the proof of Proposition 14.12, you know how to exploit this formula at order $n=1,2,3,4$. The same method works in general, and after some computations, this leads to the formula that we want to establish, namely:
$$M_n(a)=\sum_{\nu\in NC(n)}\kappa_\nu(a)$$

(3) We are therefore led to the conclusions in the statement. All this was of course quite brief, and for details here, we refer for instance to Nica-Speicher \cite{nsp}.
\end{proof}

The above free cumulant technology is quite powerful, and can be applied as well to other measures that we know. For the arcsine law $\alpha_t$, the free cumulants are:
$$\kappa_{2r+2}=(-1)^r2t^{r+1}\frac{1}{r+1}\binom{2r}{r}$$

For the uniform law $u_t$ the free cumulants are more tricky to compute, and are given by the following formula, with $B_n\in\mathbb Q$ being the $n$-th Bernoulli number:
$$\kappa_n=\frac{(-t)^n}{n!}\cdot B_n$$

Finally, there is also a discussion to be made, in relation with the notion of free infinite divisibility of measures. We will leave some learning here as an exercise.

\section*{14c. The bijection}

As a main application now of the theories of Rota cumulants and Speicher free cumulants, following Bercovici and Pata \cite{bpa}, we can formulate the following simple and bright definition, making a rock-solid connection between classical and free:

\index{convolution semigroup}
\index{free convolution semigroup}
\index{Bercovici-Pata bijection}
\index{classical cumulants}
\index{free cumulants}

\begin{definition}
A convolution semigroup of measures
$$\{m_t\}_{t>0}\quad:\quad m_s*m_t=m_{s+t}$$
is in Bercovici-Pata bijection with a free convolution semigroup of measures
$$\{\mu_t\}_{t>0}\quad:\quad \mu_s\boxplus\mu_t=\mu_{s+t}$$
when the classical cumulants of $m_t$ coincide with the free cumulants of $\mu_t$.
\end{definition}

Getting now to the examples, we know from the above cumulant computations that we have correspondences $g_t\leftrightarrow\gamma_t$ and $p_t\leftrightarrow\pi_t$, as well as the trivial correspondence $\delta_t\leftrightarrow\delta_t$. In order to find more examples, we can take some inspiration from the group theory material from chapter 7, which suggests looking at two more situations:
$$g_t^t\leftrightarrow\gamma_t^t\qquad,\qquad p_t^2\leftrightarrow\pi_t^2$$

So, let us do this. Let us first welcome the shifted semicircle law $\gamma_t^t$, with due respect, meaning full probabilistic computations for it, the result being as follows:

\begin{theorem}
The shifted semicircle law $\gamma_t^t$ of parameter $t>0$, namely
$$\gamma_t^t=\frac{1}{2\pi t}\sqrt{4t-(x-t)^2}\,dx$$
on $[t-2\sqrt{t},t+2\sqrt{t}]$, has the following properties:
\begin{enumerate}
\item The mean is $E=t$, the variance is $V=t$.

\item The moments are $M_k=\sum_{s=0}^{[k/2]}\binom{k}{2s}\frac{1}{s+1}\binom{2s}{s}t^{k-s}$.

\item Equivalently, we have $M_k=\sum_{\pi\in NC_{12}(k)}t^{|\pi|}$.

\item The skewness is $\gamma=0$, the kurtosis is $\kappa=2$.

\item $G(\xi)=(\xi-t-\sqrt{(\xi-t)^2-4t})/2t$.

\item $M(z)=(1-tz-\sqrt{(1-tz)^2-4tz^2})/2tz^2$.

\item The $R$-transform is $R(z)=t+tz$.

\item The free cumulants are $t,t,0,0,\ldots$

\item We have the formula $\gamma_s^s\boxplus\gamma_t^t=\gamma_{s+t}^{s+t}$.

\item This measure is the free version of $g_t^t$.

\item The Hankel determinants at $t=1$ are $H_k=1$, exactly as for $\gamma_1$.

\item The orthogonal polynomials are modified Chebycheff polynomials.
\end{enumerate}
\end{theorem}

\begin{proof}
The first moment formula is what comes out of the density, and the second moment formula, which is more practical, comes from this. Numerically, we have:
$$M_1=t$$
$$M_2=t^2+t$$
$$M_3=t^3+3t^2$$
$$M_4=t^4+6t^3+2t^2$$

But this gives the formulae of $E,V,\gamma,\kappa$ in the statement. Next, we have:
$$G(\xi)=G_{\gamma_t}(\xi-t)=\frac{\xi-t-\sqrt{(\xi-t)^2-4t}}{2t}$$ 

But this gives the formula of $M(z)=z^{-1}G(z^{-1})$ in the statement, and we can capture the $R$-transform too, coming from $G(t+tz+z^{-1})=z$. The free cumulant formulae follow, as well as the semigroup assertion, and the liberation claim too. Finally, the discussion regarding the Hankel determinants and orthogonal polynomials is standard.
\end{proof}

Next, let us welcome the free Bessel law $\pi_t^2$, with due respect too. We have:

\begin{theorem}
The free Bessel law $\pi_t^2$ of parameter $t>0$, appearing as
$$\pi_t^2=\pi_{t\varepsilon}\quad,\quad\varepsilon=\frac{\delta_{-1}+\delta_1}{2}$$
has the following properties:
\begin{enumerate}
\item The mean is $E=0$, the variance is $V=t$.

\item The moments are $M_k=\sum_{\pi\in NC_{even}(k)}t^{|\pi|}$.

\item Equivalently, $M_{2r}=\sum_{b=1}^r\frac{1}{b}\binom{r-1}{b-1}\binom{2r}{b-1}t^b$.

\item At $t=1$ the moments are $M_{2r}=\frac{1}{2r+1}\binom{3r}{r}$.

\item The skewness is $\gamma=0$, the kurtosis is $\kappa=2+1/t$.

\item $M=M(z)$ satisfies $M=1+(zM)^2(M+t-1)$.

\item $G=G(\xi)$ satisfies $\xi G=1+\xi G^3+(t-1)G^2$.

\item The $R$-transform is $R(z)=tz/(1-z^2)$.

\item The free cumulants are $0,t,0,t,\ldots$

\item We have the formula $\pi^2_s\boxplus\pi^2_t=\pi^2_{s+t}$.

\item This measure is the free version of $p_t^2$.

\item The density is analytic, with an atom at $0$ when $t<1$.
\end{enumerate}
\end{theorem}

\begin{proof}
The moment formulae are something that can be obtained from Fourier, but we prefer to discuss this in the next chapter, using easiness. Numerically, we get:
$$M_1=0\quad,\quad 
M_2=t\quad,\quad
M_3=0\quad,\quad
M_4=2t^2+t$$

But this gives the above formulae of $E,V,\gamma,\kappa$. As for the rest, this is standard, and we partly know all this, save for the study of the density, for which we refer to \cite{bb+}.
\end{proof}

We can go back now to the Bercovici-Pata bijection, and formulate a nice result:

\begin{theorem}
In the standard cube of basic real probability limiting measures
$$
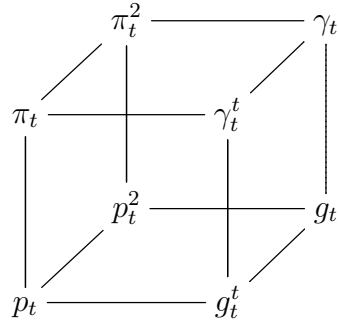
\xymatrix@R=16pt@C=20pt{
&\pi_t^2\ar@{-}[rr]\ar@{-}[dd]&&\gamma_t\ar@{-}[dd]\\
\pi_t\ar@{-}[rr]\ar@{-}[dd]\ar@{-}[ur]&&\gamma_t^t\ar@{-}[dd]\ar@{-}[ur]\\
&p_t^2\ar@{-}[rr]\ar@{-}[uu]&&g_t\ar@{.}[uu]\\
p_t\ar@{-}[uu]\ar@{-}[ur]\ar@{-}[rr]&&g_t^t\ar@{-}[uu]\ar@{-}[ur]
}$$
the upper measures appear as free versions of the lower measures.
\end{theorem}

\begin{proof}
This comes indeed from our various cumulant computations above.
\end{proof}

At a more advanced level now, we have an easiness result too, as follows:

\begin{theorem}[continuation]
The moments are $M_k=\sum_{\pi\in D(k)}t^{|\pi|}$, with
$$
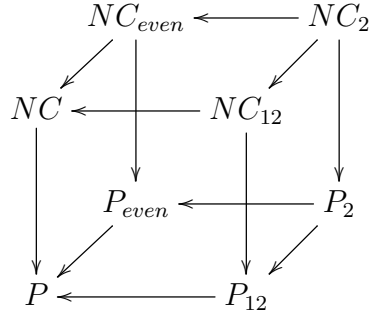
\xymatrix@R=18pt@C=2pt{
&NC_{even}\ar[dl]\ar[dd]&&NC_2\ar[dl]\ar[ll]\ar[dd]\\
NC\ar[dd]&&NC_{12}\ar[dd]\ar[ll]\\
&P_{even}\ar[dl]&&P_2\ar[dl]\ar[ll]\\
P&&P_{12}\ar[ll]
}$$
being the corresponding categories of partitions $D\subset P$.
\end{theorem}

\begin{proof}
This is again self-explanatory, based on our various moment results. Let us also mention that, in the context of the classification results from chapter 4, with $L\subset\mathbb N$ being as there the set of sizes of blocks of the partitions in $D$, the classical cumulants on the bottom are $k_n=t\delta_{n\in L}$, and the free cumulants on top are $\kappa_n=t\delta_{n\in L}$. See \cite{bsp}.
\end{proof}

Getting now to the unitary case, meaning complex measures such as $G_t$, or distributions which are not usual measures such as $\Gamma_t$, things here are more tricky, because the Bercovici-Pata bijection is something of real nature. We can nevertheless persist, and adopt the easiness point of view on liberation, in order to have our say here.

\bigskip

So, let us welcome to our family of basic limiting measures the shifted circular law $\Gamma_t^t$, which is something self-explanatory, along with the full family of free Bessel laws $\{\pi^s_t|s\in\{1,2,\ldots,\infty\}\}$, appearing as $\pi^s_t=\pi_{t\varepsilon_s}$, with $\varepsilon_s$ being the uniform measure on the $s$-th roots of unity, with in particular the purely complex free Bessel law $\Pi_t=\pi_t^\infty$.

\bigskip

With these conventions, we have the following complex analogue of Theorem 14.21, which is formulated with respect to the easiness viewpoint on liberation:

\begin{theorem}
In the standard cube of basic complex limiting measures
$$
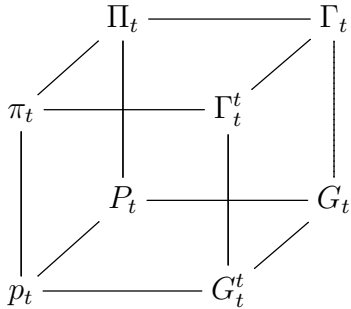
\xymatrix@R=16pt@C=20pt{
&\Pi_t\ar@{-}[rr]\ar@{-}[dd]&&\Gamma_t\ar@{-}[dd]\\
\pi_t\ar@{-}[rr]\ar@{-}[dd]\ar@{-}[ur]&&\Gamma_t^t\ar@{-}[dd]\ar@{-}[ur]\\
&P_t\ar@{-}[rr]\ar@{-}[uu]&&G_t\ar@{.}[uu]\\
p_t\ar@{-}[uu]\ar@{-}[ur]\ar@{-}[rr]&&G_t^t\ar@{-}[uu]\ar@{-}[ur]
}$$
the upper measures appear as free versions of the lower measures.
\end{theorem}

\begin{proof}
This is standard too, with the corresponding categories $D\subset P$ being:
$$
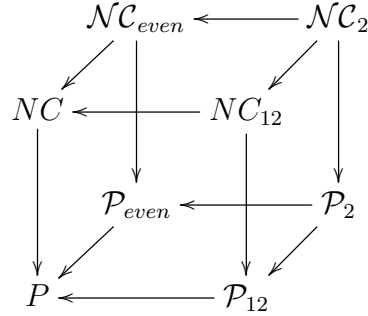
\xymatrix@R=18pt@C=2pt{
&\mathcal{NC}_{even}\ar[dl]\ar[dd]&&\mathcal{NC}_2\ar[dl]\ar[ll]\ar[dd]\\
NC\ar[dd]&&NC_{12}\ar[dd]\ar[ll]\\
&\mathcal P_{even}\ar[dl]&&\mathcal P_2\ar[dl]\ar[ll]\\
P&&\mathcal P_{12}\ar[ll]
}$$

Thus, we are led to the conclusion in the statement. See \cite{bb+}, \cite{bsp}, \cite{twe}.
\end{proof}

As a last manipulation, we can merge our real and complex cubes, and we obtain:

\begin{theorem}
In the standard cube of basic probability limiting measures
$$
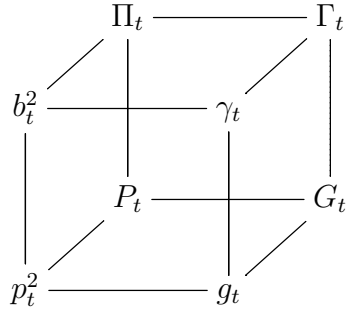
\xymatrix@R=16pt@C=20pt{
&\Pi_t\ar@{-}[rr]\ar@{-}[dd]&&\Gamma_t\ar@{-}[dd]\\
b_t^2\ar@{-}[rr]\ar@{-}[dd]\ar@{-}[ur]&&\gamma_t\ar@{-}[dd]\ar@{-}[ur]\\
&P_t\ar@{-}[rr]\ar@{-}[uu]&&G_t\ar@{.}[uu]\\
p_t^2\ar@{-}[uu]\ar@{-}[ur]\ar@{-}[rr]&&g_t\ar@{-}[uu]\ar@{-}[ur]
}$$
the upper measures appear as free versions of the lower measures.
\end{theorem}

\begin{proof}
We would like to merge our 3D cubes, real and complex, into a 4D cube. Unfortunately, this cannot be properly done, due to the double presence of $p_t,\pi_t$. So, what to do? Well, do some cleanup, leading to the above 3D cube, corresponding to:
$$
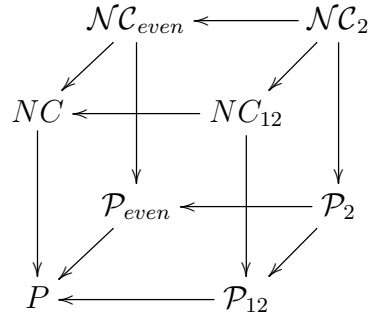
\xymatrix@R=18pt@C=2pt{
&\mathcal{NC}_{even}\ar[dl]\ar[dd]&&\mathcal{NC}_2\ar[dl]\ar[ll]\ar[dd]\\
NC\ar[dd]&&NC_{12}\ar[dd]\ar[ll]\\
&\mathcal P_{even}\ar[dl]&&\mathcal P_2\ar[dl]\ar[ll]\\
P&&\mathcal P_{12}\ar[ll]
}$$

Finally, observe that our cube is something quite subtle, featuring a classical/free correspondence on the vertical, a discrete/continuous correspondence on the horizontal, and a real/complex correspondence on the depth. Nice and conceptual.
\end{proof}

And with this, end of our discussion regarding the Bercovici-Pata bijection. All this was of course basics, and at the advanced level, crazy things can happen, such as:

\begin{question}[Belinschi, Bo\.zejko, Lehner, Speicher]
The normal law $g_1$ being freely infinitely divisible, what is its classical analogue?
\end{question}

And isn't this puzzling, weird question that we have here, obtained by playing with fire, namely applying free transforms to classical measures, although the converse can be certainly attempted too. For more on all this, we refer to their paper \cite{bbl}.

\section*{14d. Chi variables}

Getting now to the chi, chi square problematics in free probability, we would need a linearization result for $\boxtimes$. And we have here the following result of Voiculescu \cite{vo2}:

\index{S-transform}
\index{Stieltjes transform}
\index{multiplicative free convolution}

\begin{theorem}
The operation $\boxtimes$ can be linearized as follows:
\begin{enumerate}
\item Start with $\psi(z)=M_1z+M_2z^2+M_3z^3+\ldots$

\item Compute its inverse, $\psi(\chi(z))=z$.

\item Compute $S(z)=(1+z^{-1})\chi(z)$.

\item Then $\log S$ linearizes the free multiplicative convolution, $S_{\mu\boxtimes\nu}=S_\mu S_\nu$.
\end{enumerate}
\end{theorem}

\begin{proof}
There are several proofs here, and following Haagerup \cite{haa}, we have:

\medskip

(1) According to our conventions from chapter 13, we want to prove that, given noncommutative variables $a,b$ which are free, we have the following formula:
$$S_{\mu_{ab}}(z)=S_{\mu_a}(z)S_{\mu_b}(z)$$

(2) For this purpose, consider the orthogonal shifts $S,T$ on the free Fock space, from chapter 13. By using the algebraic arguments there, from the proof of the $R$-transform theorem, we can assume, a bit as there, that our variables are as follows:
$$a=(1+S)f(S^*)\quad,\quad 
b=(1+T)g(T^*)$$

Our claim, which will prove the theorem, is that we have the following formulae, for the $S$-transforms of the various variables involved:
$$S_{\mu_a}(z)=\frac{1}{f(z)}\quad,\quad 
S_{\mu_b}(z)=\frac{1}{g(z)}\quad,\quad 
S_{\mu_{ab}}(z)=\frac{1}{f(z)g(z)}$$

(3) Let us first compute $S_{\mu_a}$. We know that we have $a=(1+S)f(S^*)$, with $S$ being the shift on $l^2(\mathbb N)$. Given $|z|<1$, consider the following vector:
$$p=\sum_{k\geq0}z^ke_k$$ 

The shift and its adjoint act on this vector in the following way:
$$Sp=\sum_{k\geq0}z^ke_{k+1}=\frac{p-e_0}{z}\quad,\quad 
S^*p=\sum_{k\geq1}z^ke_{k-1}=zp$$

Thus $f(S^*)p=f(z)p$, and we deduce from this that we have:
\begin{eqnarray*}
ap
&=&(1+S)f(z)p\\
&=&f(z)(p+Sp)\\
&=&f(z)\left(p+\frac{p-e_0}{z}\right)\\
&=&\left(1+\frac{1}{z}\right)f(z)p-\frac{f(z)}{z}e_0
\end{eqnarray*}

By dividing everything by $(1+1/z)f(z)$, this formula becomes:
$$\frac{z}{1+z}\cdot\frac{1}{f(z)}\,ap=p-\frac{e_0}{1+z}$$

We can write this latter formula in the following way:
$$\left(1-\frac{z}{1+z}\cdot\frac{1}{f(z)}\,a\right)p=\frac{e_0}{1+z}$$

Now by inverting, we obtain from this the following formula:
$$\left(1-\frac{z}{1+z}\cdot\frac{1}{f(z)}\,a\right)^{-1}e_0=(1+z)p$$

(4) But this gives us the formula of $S_{\mu_a}$. Indeed, consider the following function:
$$\rho(z)=\frac{z}{1+z}\cdot\frac{1}{f(z)}$$

With this notation, the formula that we found in (3) becomes:
$$(1-\rho(z)a)^{-1}e_0=(1+z)p$$

By using this, in terms of $\varphi(T)=<Te_0,e_0>$, we obtain:
\begin{eqnarray*}
\varphi\left((1-\rho(z)a)^{-1}\right)
&=&<(1-\rho(z)a)^{-1}e_0,e_0>\\
&=&<(1+z)p,e_0>\\
&=&1+z
\end{eqnarray*}

Thus the above function $\rho$ is the inverse of the following function:
$$\psi(z)=\varphi\left(\frac{1}{1-za}\right)-1$$

But this latter function is the $\psi$ function from the statement, and so $\rho$ is the function $\chi$ from the statement, and we can finish our computation, as follows:
$$S_{\mu_a}(z)
=\frac{1+z}{z}\cdot\rho(z)
=\frac{1+z}{z}\cdot\frac{z}{1+z}\cdot\frac{1}{f(z)}
=\frac{1}{f(z)}$$

(5) A similar computation, or just a symmetry argument, gives $S_{\mu_b}(z)=1/g(z)$. In order to compute now $S_{\mu_{ab}}(z)$, we use a similar trick. Consider the following vector of $l^2(\mathbb N*\mathbb N)$, with the primes and double primes referring to the two copies of $\mathbb N$:
$$q=e_0+\sum_{k\geq1}(e_1'+e_1''+e_1'\otimes e_1'')^{\otimes k}$$

The adjoints of the shifts $S,T$ act as follows on this vector:
$$S^*q=z(1+T)q\quad,\quad T^*q=zq$$

By using these formulae, we have the following computation:
\begin{eqnarray*}
abq
&=&(1+S)f(S^*)(1+T)g(T^*)q\\
&=&(1+S)f(S^*)(1+T)g(z)q\\
&=&g(z)(1+S)f(S^*)(1+T)q
\end{eqnarray*}

In order to compute the last term, observe that we have:
\begin{eqnarray*}
S^*(1+T)q
&=&(S^*+S^*T)q\\
&=&S^*q\\
&=&z(1+T)q
\end{eqnarray*}

Thus $f(S^*)(1+T)q=f(z)(1+T)q$, and back to our computation, we have:
\begin{eqnarray*}
abq
&=&g(z)(1+S)f(z)(1+T)q\\
&=&f(z)g(z)(1+S)(1+T)q\\
&=&f(z)g(z)\left(\frac{1+z}{z}\cdot q-\frac{e_0}{z}\right)
\end{eqnarray*}

Now observe that we can write this formula as follows:
$$\left(1-\frac{z}{1+z}\cdot\frac{1}{f(z)g(z)}\cdot ab\right)q=\frac{e_0}{1+z}$$

By inverting, we obtain from this the following formula:
$$\left(1-\frac{z}{1+z}\cdot\frac{1}{f(z)g(z)}\cdot ab\right)^{-1}e_0=(1+z)q$$

(6) But this formula that we obtained is similar to the formula that we obtained at the end of (3) above. Thus, we can use the same argument as in (4), and we obtain:
$$S_{\mu_{ab}}(z)=\frac{1}{f(z)g(z)}$$

We are therefore done with the computations, and this finishes the proof.
\end{proof}

With the above result, coupled with the other results from chapter 13, we can potentially do many things, in relation with the chi, chi square problematics in free probability. We will discuss in what follows two such results, one continuous and one discrete. 

\bigskip

The continuous result that we would like to discuss, due to Voiculescu \cite{vo3}, which is something of key importance for the study of circular variables, is as follows:

\begin{theorem}
The polar decomposition of circular variables is $c=uq$, with:
\begin{enumerate}
\item $u$ being a Haar unitary, meaning a unitary with $tr(u^k)=\delta_{k0}$, for $k\in\mathbb Z$.

\item $q$ being a quarter-circular, meaning $q=|s|$ with $s$ semicircular.

\item $u,q$ being free.
\end{enumerate}
\end{theorem}

\begin{proof}
This is something quite tricky, the idea being as follows:

\medskip

(1) To start with, the statement assumes that we are inside a von Neumann algebra, and this in order for the polar decomposition of operators to stay inside the algebra. Also, the trace must be assumed to be faithful, in order to be able to recover things.

\medskip

(2) In order now to get familiar with the polar decomposition, let us first decompose a semicircular $s$. For this purpose, we can use the following simple model for $s$:
$$s=x\in L^\infty\left([-2,2],\gamma_1\right)$$

Indeed, we can see right away that the decomposition is $s=eq$, with $e=sgn(x)$ and $q=|x|$. Thus, probabilistically, we have $s=eq$, with $e$ being a unitary having moments $1,0,1,0,\ldots$, then $q=|s|$ being quarter-circular, and with $e,q$ being independent. 

\medskip

(3) Getting now to circular variables, as in the statement, things here are more tricky, because we don't have available classical models as in (2). However, the known models, using shifts on Fock spaces, or random matrices, both work. Alternatively, the result can be deduced from a variation of Theorem 14.26. For a discussion here, see \cite{nsp}, \cite{vo3}.
\end{proof}

As for the discrete result that we would like to present, from \cite{bb+}, this is:

\begin{theorem}
The modified free Bessel law $\tilde{\pi}^s_t$, which is the real probability measure given by the following moment formula, which must hold for any $k\in\mathbb N$,
$$M_k(\tilde{\pi}^s_t)=M_{sk}(\pi^s_t)$$
is subject to the following formulae, involving $\pi=\pi_1$ only:
\begin{enumerate}
\item $\tilde{\pi}^s_t=\pi^{\boxtimes s-1}\boxtimes\pi^{\boxplus t}$, at $s\geq1,t>0$.

\item $\tilde{\pi}^s_t=((1-t)\delta_0+t\delta_1)\boxtimes\pi^{\boxtimes s}$, at $s>0,t\in(0,1]$.
\end{enumerate}
\end{theorem}

\begin{proof}
This is again something quite tricky, the idea being as follows:

\medskip

(1) To start with, at $s=2$, the free Bessel law $\pi^2_t$ being a real probability measure, we have available a squared version of it, given by the following formula:
$$f\sim\pi^2_t\ \implies\ f^2\sim\tilde{\pi}^2_t$$

Equivalently, in terms of moments, this modified free Bessel law $\tilde{\pi}^2_t$ is given by:
$$M_k(\tilde{\pi}^2_t)=M_{2k}(\pi^2_t)$$

And the point is that many things can be said about $\tilde{\pi}^2_t$, the idea being that its combinatorics is simpler than that of $\pi^2_t$ itself. With as illustrations here being our claims in the statement, which are both crucially based on the doubling operation.

\medskip

(2) In view of this, the problem appears, what is the correct generalization of the above law $\tilde{\pi}^2_t$, in the general case $s\in\mathbb N\cup\{\infty\}$? And here, at $s=1$ to start with, since the Marchenko-Pastur law $\pi_t=\pi_t^1$ will not simplify when squaring, or doing other types of similar manipulations on it, we must leave this law untouched, $\tilde{\pi}^1_t=\pi^1_t$.

\medskip

(3) Next, some further thinking along these lines, and some computations too, suggest a manipulation as follows, which is in tune with what we have above at $s=1,2$:
$$f\sim\pi^s_t\ \implies\ f^s\sim\tilde{\pi}^s_t$$

However, passed $s=\infty$ that we obviously have to leave aside, this won't work well at $s\geq3$, because we will not obtain in this way a real probability measure. So, the trick is the one in the statement, restricting the attention to the uncolored moments of $f^s$, which amounts in defining $\tilde{\pi}^s_t$ by the following formula, which must hold for any $k\in\mathbb N$:
$$M_k(\tilde{\pi}^s_t)=M_{sk}(\pi^s_t)$$

(4) With this understood, getting now to what the statement says, in order to prove that, we first have the following computation, for the $S$-transform of $\pi_t=\pi^{\boxplus t}$:
\begin{eqnarray*}
\xi G^2+1=(\xi+1-t)G
&\implies&z(\psi+1)^2+1=(1+z-zt)(\psi+1)\\
&\implies&\chi(z+1)^2+1=(1+\chi-\chi t)(z+1)\\
&\implies&\chi(z+1)(t+z)=z\\
&\implies&S=1/(t+z)
\end{eqnarray*}

We conclude that the $S$-transform of the measure $\pi^{\boxtimes s-1}\boxtimes\pi^{\boxplus t}$ is given by:
$$S(z)=\frac{1}{(1+z)^{s-1}}\cdot\frac{1}{t+z}$$

(5) Next, the $S$-transform of $(1-t)\delta_0+t\delta_1$ can be computed as follows:
\begin{eqnarray*}
M=1+tz/(1-z)
&\implies&\psi=tz/(1-z)\\
&\implies&z=t\chi/(1-\chi)\\
&\implies&\chi=z/(t+z)\\
&\implies& S=(1+z)/(t+z)
\end{eqnarray*}

Thus, the $S$-transform of the measure $((1-t)\delta_0+t\delta_1)\boxtimes\pi^{\boxtimes s}$ is given by:
$$S(z)=\frac{1}{(1+z)^s}\cdot\frac{1+z}{t+z}$$

(6) Summarizing, done with two measures. In order to deal now with $\tilde{\pi}^s_t$, we first have the following computation, for the $R$-transform of $\pi^s_t$, with $w=e^{2\pi i/s}$:
$$R(z)=\frac{t}{s}\sum_{k=1}^s\frac{w^k}{1-w^kz}=t\cdot\frac{z^{s-1}}{1-z^s}$$ 

Next, we can do some standard manipulations on this formula, as follows:
\begin{eqnarray*}
R=tz^{s-1}/(1-z^s)
&\implies&M=1+(zM)^s(M+t-1)\\
&\implies&\tilde{M}=1+z\tilde{M}^s(\tilde{M}+t-1)\\
&\implies&\tilde{\psi}=z(1+\tilde{\psi})^s(t+\tilde{\psi})\\
&\implies&z=\tilde{\chi}(1+z)^s(t+z)\\
&\implies&(1+z^{-1})z=\tilde{S}(1+z)^s(t+z)\\
&\implies&\tilde{S}=(1+z)^{1-s}/(t+z)
\end{eqnarray*}

Thus the $S$-transform of $\tilde{\pi}^s_t$ coincides with the $S$-transforms computed in (4,5), and we are done. For further details here, and more on the free Bessel laws, we refer to \cite{bb+}.
\end{proof}

\section*{14e. Exercises}

This was an exciting chapter, and as exercises about this, we have:

\begin{exercise}
Compute the cumulants of the exponential law $e_t$.
\end{exercise}

\begin{exercise}
Compute the cumulants of the uniform law $u_t$.
\end{exercise}

\begin{exercise}
Compute the free cumulants of the arcsine law $\alpha_t$.
\end{exercise}

\begin{exercise}
Compute the free cumulants of the uniform law $u_t$.
\end{exercise}

\begin{exercise}
Clarify everything that we said, regarding $\gamma_t^t$.
\end{exercise}

\begin{exercise}
Clarify everything that we said, regarding $\pi^2_t$.
\end{exercise}

\begin{exercise}
Learn more about circular and semicircular systems.
\end{exercise}

\begin{exercise}
Learn more about $\tilde{\pi}^s_t$, and about their classical versions $\tilde{p}^s_t$ too.
\end{exercise}

As bonus exercise, have some thinking at what happens in the unitary case.

\chapter{Quantum groups}

\section*{15a. Quantum groups}

We have seen in the last 2 chapters that a reliable theory of free probability can be developed, essentially based on Voiculescu's $R$-transform, which linearizes the free convolution operation $\boxplus$. Indeed, with this transform, which is the free analogue of $\log F$, linearizing the usual convolution operation $*$, we can do as many things in free probability as you can do in classical probability, using Fourier. And, sky is the limit.

\bigskip

This being said, we sort of reached the limit already, at least at the level of basic things that can be said. So, what is next? We can try here to get more technical, along the lines of the 5-year plan evoked in the beginning of chapter 14. Or upgrade to free geometry and physics, as per the general hierarchy explained in the beginning of chapter 13. Or do something else, so many things that can be potentially tried.

\bigskip

Not clear what to do, and looking for advice, I think I know whom to ask. The lady comes quite often to the garden, picking pretty much everything that can be picked, that I normally leave there to the attention of bigger animals, like cats, chickens or foxes. I bet that, among the million things accumulated at her nest, she has some free probability measures, that she can lend me for this chapter. Let's see if we can make a deal:

\begin{magpie}
There's basic stealing, and advanced stealing.
\end{magpie}

Well, thanks magpie, not that I got what I wanted, but this looks like a good piece of advice, just be free as a bird, and enjoy. So, getting back now to what we were doing, the truth is, we've been shamelessly stealing theorems from classical probability, and it is probably most reasonable to keep doing that, but at a more advanced level:

\begin{plan}
We can upgrade to geometry, with our free manifolds $X$ appearing as abstract spectra of suitably chosen algebras $A$, and do some probability after.
\end{plan} 

Getting to work now, with free geometry, what comes first? Free spheres, you would say, and normally good answer. However, thinking a bit, we would like afterwards to integrate on these spheres, and the rotationally invariant measure that we need is most likely a ``quantum rotationally invariant'' measure, involving the action of a ``free rotation group''. So, we need quantum groups first, and the spheres will come after.

\bigskip

So, let us introduce the quantum groups. The good axioms here, due to Woronowicz \cite{wo1}, and slightly modified for our present purposes, are as follows:

\begin{definition}
A Woronowicz algebra is a $C^*$-algebra $A$, given with a unitary matrix $u\in M_N(A)$ whose coefficients generate $A$, such that the formulae
$$\Delta(u_{ij})=\sum_ku_{ik}\otimes u_{kj}\quad,\quad
\varepsilon(u_{ij})=\delta_{ij}\quad,\quad 
S(u_{ij})=u_{ji}^*$$
define morphisms of $C^*$-algebras $\Delta:A\to A\otimes A$, $\varepsilon:A\to\mathbb C$ and $S:A\to A^{opp}$, called comultiplication, counit and antipode. 
\end{definition}

To be more precise, we are using here the general $C^*$-algebra formalism, that we know well since chapters 8 and 13. The tensor product $\otimes$ appearing above can be any $C^*$-algebra tensor product, and more on this in a moment. As for $A^{opp}$, this is the opposite algebra, having multiplication $a\cdot b=ba$, and more on this in a moment too.

\bigskip

So, let's see how this works. With the convention that $A$ is cocommutative when $\Sigma\Delta=\Delta$, where $\Sigma(a\otimes b)=b\otimes a$ is the flip, we have the following key result:

\begin{proposition}
The following are Woronowicz algebras:
\begin{enumerate}
\item $C(G)$, with $G\subset U_N$ compact Lie group, with structural maps as follows:
$$\Delta(f)=\big[(g,h)\to f(gh)\big]\quad,\quad 
\varepsilon(f)=f(1)\quad,\quad
S(f)=\big[g\to f(g^{-1})\big]$$

\item $C^*(\Gamma)$, with $F_N\to\Gamma$ finitely generated group, with structural maps as follows:
$$\Delta(g)=g\otimes g\quad,\quad 
\varepsilon(g)=1\quad,\quad
S(g)=g^{-1}$$
\end{enumerate}
Moreover, we obtain in this way all the commutative/cocommutative algebras.
\end{proposition}

\begin{proof}
For the first assertion we can use the matrix $u=(u_{ij})$ formed by the standard matrix coordinates of $G$, which is by definition given by:
$$g=\begin{pmatrix}
u_{11}(g)&\ldots&u_{1N}(g)\\
\vdots&&\vdots\\
u_{N1}(g)&\ldots&u_{NN}(g)
\end{pmatrix}$$

Indeed, all axioms are clearly satisfied. Next, for the second assertion we can use the diagonal matrix formed by the generators, with again the axioms being satisfied, and with the remark here that, for things to work fine, the target of $S$ must be indeed $A^{opp}$:
$$u=\begin{pmatrix}
g_1&&0\\
&\ddots&\\
0&&g_N
\end{pmatrix}$$

Finally, regarding the last assertion, in the commutative case this follows from the Gelfand theorem, and in the cocommutative case, we will be back to this.
\end{proof}

The above result is quite interesting, and based on it, and on Gelfand duality in general, we can formulate the following definition, complementing Definition 15.3:

\begin{definition}[continuation]
Given a Woronowicz algebra $(A,u)$, we write it as
$$A=C(G)=C^*(\Gamma)$$
and call $G$ a compact quantum Lie group, and $\Gamma$ a finitely generated discrete quantum group. We say that these quantum groups are dual to each other, and write: 
$$G=\widehat{\Gamma}\quad,\quad\Gamma=\widehat{G}$$
Also, we agree to identify two Woronowicz algebras, $(A,u)=(B,v)$, when we have an isomorphism of $*$-algebras $<u_{ij}>\simeq<v_{ij}>$ mapping $u_{ij}\to v_{ij}$.
\end{definition}

To be more precise, the convention $A=C(G)$ comes from Proposition 15.4 and Gelfand, the convention $A=C^*(\Gamma)$ comes from Proposition 15.4, and the resulting duality between $G,\Gamma$ generalizes the usual Pontrjagin duality, for the abelian groups. 

\bigskip

As for the convention at the end, which by the way clarifies our $\otimes$ comment made after Definition 15.3, this is something more technical, coming from the fact that a non-amenable discrete group $\Gamma$ has several possible completions of $\mathbb C[\Gamma]$, that we must identify, in order to have uniqueness of $\Gamma$, as a discrete quantum group. More on this later.

\bigskip

Let us discuss now some tools for studying our beasts. We first have:

\begin{proposition}
Let $(A,u)$ be a Woronowicz algebra.
\begin{enumerate} 
\item $\Delta,\varepsilon$ satisfy the usual axioms for a comultiplication and a counit, namely:
$$(\Delta\otimes id)\Delta=(id\otimes \Delta)\Delta$$
$$(\varepsilon\otimes id)\Delta=(id\otimes\varepsilon)\Delta=id$$

\item $S$ satisfies the antipode axiom, on the $*$-algebra generated by entries of $u$: 
$$m(S\otimes id)\Delta=m(id\otimes S)\Delta=\varepsilon(.)1$$

\item In addition, the square of the antipode is the identity, $S^2=id$.
\end{enumerate}
\end{proposition}

\begin{proof}
As a first observation, the result holds in the commutative case, $A=C(G)$ with $G\subset U_N$. Indeed, here we know from Proposition 15.4 that $\Delta,\varepsilon,S$ appear as functional analytic transposes of the multiplication, unit and inverse maps $m,u,i$:
$$\Delta=m^t\quad,\quad 
\varepsilon=u^t\quad,\quad 
S=i^t$$

Thus, the various conditions in the statement on $\Delta,\varepsilon,S$ simply come from the group axioms satisfied by $m,u,i$. In general now, the first axiom follows from:
$$(\Delta\otimes id)\Delta(u_{ij})=(id\otimes \Delta)\Delta(u_{ij})=\sum_{kl}u_{ik}\otimes u_{kl}\otimes u_{lj}$$

As for the other axioms, the verifications here are similar.
\end{proof}

In order to reach to more advanced results, the idea will be that of doing representation theory. Following Woronowicz \cite{wo1}, let us start with the following definition:

\begin{definition}
Given $(A,u)$, we call corepresentation of it any unitary matrix $v\in M_n(\mathcal A)$, with $\mathcal A=<u_{ij}>$, satisfying the same conditions as $u$, namely:
$$\Delta(v_{ij})=\sum_kv_{ik}\otimes v_{kj}\quad,\quad
\varepsilon(v_{ij})=\delta_{ij}\quad,\quad
S(v_{ij})=v_{ji}^*$$
We also say that $v$ is a representation of the underlying compact quantum group $G$.
\end{definition}

In the commutative case, $A=C(G)$ with $G\subset U_N$, we obtain in this way the finite dimensional unitary smooth representations $v:G\to U_n$, via the following formula:
$$v(g)=\begin{pmatrix}
v_{11}(g)&\ldots&v_{1n}(g)\\
\vdots&&\vdots\\
v_{n1}(g)&\ldots&v_{nn}(g)
\end{pmatrix}$$

With these conventions, we have the following fundamental result, from \cite{wo1}:

\begin{theorem}
Any Woronowicz algebra has a unique Haar integration, 
$$\left(\int_G\otimes\, id\right)\Delta=\left(id\otimes\int_G\right)\Delta=\int_G(.)1$$
and for any corepresentation $v\in M_n(\mathbb C)\otimes A$ we have the formula
$$\left(id\otimes\int_G\right)v=P$$
where $P$ is the orthogonal projection onto $Fix(v)=\{\xi\in\mathbb C^n|v\xi=\xi\}$.
\end{theorem}

\begin{proof}
This is something that we know well in the classical case, for the algebras of type $A=C(G)$ with $G\subset U_N$, from chapter 7, and the proof in general is similar.
\end{proof}

Next, we can develop the Peter-Weyl theory for the corepresentations of $A$. Consider the dense subalgebra $\mathcal A=<u_{ij}>$, and endow it with the scalar product $<a,b>=\int_Gab^*$. With this convention, we have the following key result, also from \cite{wo1}:

\begin{theorem}
We have the following Peter-Weyl type results:
\begin{enumerate}
\item Any corepresentation decomposes as a sum of irreducible corepresentations.

\item Each irreducible corepresentation appears inside a certain $u^{\otimes k}$.

\item $\mathcal A=\bigoplus_{v\in Irr(A)}M_{\dim(v)}(\mathbb C)$, the summands being pairwise orthogonal.

\item The characters of irreducible corepresentations form an orthonormal system.
\end{enumerate}
\end{theorem}

\begin{proof}
This is something that we know well from chapter 7, when $G\subset U_N$ is a compact group. In general, the proof is quite similar, by using Theorem 15.8.
\end{proof}

Finally, no discussion about compact and discrete quantum groups would be complete without a word on amenability. The result here, again from \cite{wo1}, is as follows:

\begin{theorem}
Let $A_{full}$ be the enveloping $C^*$-algebra of $\mathcal A$, and $A_{red}$ be the quotient of $A$ by the null ideal of the Haar integration. The following are then equivalent:
\begin{enumerate}
\item The Haar functional of $A_{full}$ is faithful.

\item The projection map $A_{full}\to A_{red}$ is an isomorphism.

\item The counit map $\varepsilon:A_{full}\to\mathbb C$ factorizes through $A_{red}$.

\item We have $N\in\sigma(Re(\chi_u))$, the spectrum being taken inside $A_{red}$.
\end{enumerate}
If this is the case, we say that the underlying discrete quantum group $\Gamma$ is amenable.
\end{theorem}

\begin{proof}
This is well-known in the group dual case, $A=C^*(\Gamma)$, with $\Gamma$ being a usual discrete group. In general, the result follows by adapting the group dual case proof:

\medskip

$(1)\iff(2)$ This simply follows from the fact that the GNS construction for the algebra $A_{full}$ with respect to the Haar functional produces the algebra $A_{red}$.

\medskip

$(2)\iff(3)$ Here $\implies$ is trivial, and conversely, a counit $\varepsilon:A_{red}\to\mathbb C$ produces an isomorphism $\Phi:A_{red}\to A_{full}$, by slicing the map $\widetilde{\Delta}:A_{red}\to A_{red}\otimes A_{full}$.

\medskip

$(3)\iff(4)$ Here $\implies$ is clear, coming from $\varepsilon(N-Re(\chi (u)))=0$, and the converse can be proved by doing some functional analysis. See \cite{wo1}.
\end{proof}

The compact quantum groups include the compact Lie groups, $G\subset U_N$, and the abstract duals $G=\widehat{\Gamma}$ of the finitely generated groups $F_N\to\Gamma$. Following Wang \cite{wan}, let us discuss now a number of truly ``new'' quantum groups. We first have:

\index{free orthogonal group}
\index{free unitary group}

\begin{theorem}
The following universal algebras are Woronowicz algebras,
$$C(O_N^+)=C^*\left((u_{ij})_{i,j=1,\ldots,N}\Big|u=\bar{u},u^t=u^{-1}\right)$$
$$C(U_N^+)=C^*\left((u_{ij})_{i,j=1,\ldots,N}\Big|u^*=u^{-1},u^t=\bar{u}^{-1}\right)$$
so the underlying quantum spaces $O_N^+,U_N^+$ are compact quantum groups.
\end{theorem}

\begin{proof}
This comes from the elementary fact that if a matrix $u=(u_{ij})$ is orthogonal or biunitary, then so must be the following matrices:
$$(u^\Delta)_{ij}=\sum_ku_{ik}\otimes u_{kj}\quad,\quad 
(u^\varepsilon)_{ij}=\delta_{ij}\quad,\quad 
(u^S)_{ij}=u_{ji}^*$$

Thus we can define $\Delta,\varepsilon,S$ by using the universal property of $C(O_N^+)$, $C(U_N^+)$.
\end{proof}

Before going further, let us emphasize the fact that $O_N^+,U_N^+$ are really huge, much bigger that anything that you can imagine. Indeed, besides of course containing $O_N,U_N$, these quantum groups contain the duals of $L_N=\mathbb Z_2^{*N}$ and $F_N=\mathbb Z^{*N}$, and this due to the following quotient maps, coming from the universal properties of $C(O_N^+)$, $C(U_N^+)$:
$$C(O_N^+)\to C^*(L_N)\quad,\quad C(U_N^+)\to C^*(F_N)$$

Moving on, still following Wang \cite{wan}, let us discuss now the construction and basic properties of the quantum permutation group $S_N^+$. The construction is as follows:

\begin{theorem}
The following universal $C^*$-algebra, with magic meaning formed by projections $(p^2=p^*=p)$, summing up to $1$ on each row and each column,
$$C(S_N^+)=C^*\left((u_{ij})_{i,j=1,\ldots,N}\Big|u={\rm magic}\right)$$
is a Woronowicz algebra, and the underlying quantum group $S_N^+$, called quantum permutation group, appears as a liberation of the usual permutation group $S_N$.
\end{theorem}

\begin{proof}
The first assertion comes from the standard fact, say easy exercise for you, that if a matrix $u=(u_{ij})$ is magic, then so must be the following matrices:
$$(u^\Delta)_{ij}=\sum_ku_{ik}\otimes u_{kj}\quad,\quad 
(u^\varepsilon)_{ij}=\delta_{ij}\quad,\quad 
(u^S)_{ij}=u_{ji}^*$$

As for the second assertion, what that precisely says is that we have an inclusion $S_N\subset S_N^+$, with the group $S_N$ being the classical version of $S_N^+$, obtained at the algebra level by dividing $C(S_N^+)$ by its commutator ideal. Which, in practice, means:
$$C(S_N)=C^*_{comm}\left((u_{ij})_{i,j=1,\ldots,N}\Big|u={\rm magic}\right)$$

But this is something coming from Gelfand. Indeed, the algebra on the right must be of the form $A=C(X)$, with $X$ being a certain compact space. Now since we have coordinates $u_{ij}:X\to\mathbb C$, our space must appear as $X\subset M_N(\mathbb C)$, and then the magic condition on the matrix $u=(u_{ij})$ tells us that the matrices $g\in X$ must have 0-1 entries, summing up to 1 on each row and column. Thus $X=S_N$, as desired.
\end{proof}

As before with the quantum groups $O_N^+,U_N^+$, what our construction produces is something huge, with the result here, which is worth recording, being as follows:

\begin{theorem}
The quantum permutation groups $S_N^+$ are as follows:
\begin{enumerate}
\item At $N=2,3$ we have $S_N^+=S_N$.

\item At $N=4$ we have $\widehat{D}_\infty\subset S_4^+$, so $S_4^+$ is not finite.

\item At $N=5$ we have $\widehat{\mathbb Z_2*\mathbb Z_3}\subset S_5^+$, so $S_5^+$ is not coamenable.

\item At $N=6$ and higher, $S_N^+$ remains not coamenable.
\end{enumerate}
\end{theorem}

\begin{proof}
This is something quite tricky, the idea being as follows:

\medskip

(1) This comes from the fact that a magic matrix of size $N=2,3$ must have commuting entries, which is obvious at $N=2$, and is a good exercise at $N=3$.

\medskip

(2) This comes from the fact that if $u,v$ are magic, then so is $diag(u,v)$. Indeed, we obtain in this way a quotient map $C(S_4^+)\to C^*(\mathbb Z_2*\mathbb Z_2)$, so an embedding of the dual of $D_\infty=\mathbb Z_2*\mathbb Z_2$ into $S_4^+$. Now $D_\infty$ being obviously infinite, so must be $S_4^+$.

\medskip

(3) This comes as in (2), by using a $diag(u,v)$ trick, and with the group $\mathbb Z_2*\mathbb Z_3$ being notoriously non-amenable, so must be the dual of $S_5^+$, having it as quotient.

\medskip

(4) This comes from (3), because we have $S_5^+\subset S_N^+$, obtained via $diag(u,1_{N-5})$.
\end{proof}

Moving on, we can complete our collection with some further beasts, as follows:

\begin{theorem}
We have compact quantum groups as follows:
\begin{enumerate}
\item Subgroups $B_N^+\subset O_N^+$ and $C_N^+\subset U_N^+$, obtained via the fact that $u=(u_{ij})$ must be bistochastic, that is, with sums $1$ on each row and column.

\item Quantum groups $S_N^+\subset H_N^{s+}\subset U_N^+$ with $s\in\{2,3,\ldots,\infty\}$, constructed according to the formula $H_N^{s+}=\mathbb Z_s\wr_*S_N^+$, with $\wr_*$ being a free wreath product.
\end{enumerate}
\end{theorem}

\begin{proof}
Here (1) is something self-explanatory, with the quantum group property coming as before for $O_N^+,U_N^+,S_N^+$, by looking at the matrices $u^\Delta,u^\varepsilon,u^S$, and (2) is something more tricky and algebraic, and for details here, we refer to \cite{bb+}. 
\end{proof}

And with this, good news, done with the preliminaries, and we can go ahead with easiness and probability, for our free quantum groups. Following the material from chapter 7, which was itself based on the material from chapter 4, let us start with:

\begin{definition}[update]
A category of partitions is a collection $D=\bigsqcup_{k,l}D(k,l)$ of subsets $D(k,l)\subset P(k,l)$, having the following properties:
\begin{enumerate}
\item Stability under the horizontal concatenation, $(\pi,\sigma)\to[\pi\sigma]$.

\item Stability under vertical concatenation $(\pi,\sigma)\to[^\sigma_\pi]$, with matching middle symbols.

\item Stability under the upside-down turning $*$, with switching of colors, $\circ\leftrightarrow\bullet$.

\item Each set $P(k,k)$ contains the identity partition $||\ldots||$.

\item The sets $P(\emptyset,\circ\bullet)$ and $P(\emptyset,\bullet\circ)$ both contain the semicircle $\cap$.
\end{enumerate}
\end{definition} 

In other words, what we have here is our previous definition for the categories of partitions, from chapter 7, with the last axiom there, stating that $P(k,\bar{k})$ with $|k|=2$ must contain the crossing partition $\slash\hskip-2.0mm\backslash\,$, slashed. With this, following \cite{bsp}, we have:

\begin{theorem}
Any category of partitions $D\subset P$ produces a series of compact quantum groups $G=(G_N)$, with $G_N\subset U_N^+$ for any $N\in\mathbb N$, via the formula
$$Hom(u^{\otimes k},u^{\otimes l})=span\left(T_\pi\Big|\pi\in D(k,l)\right)$$
for any $k,l$, and Tannakian duality. We call such quantum groups easy.
\end{theorem}

\begin{proof}
Indeed, the axioms in Definition 15.15 show, via the standard properties of $\pi\to T_\pi$, that the following spaces form a Tannakian category, in the sense of \cite{wo2}:
$$C_{kl}=span\left(T_\pi\Big|\pi\in D(k,l)\right)$$

Thus, Tannakian duality applies, and provides us with a closed subgroup $G_N\subset U_N^+$ such that the following equalities are satisfied, for any colored integers $k,l$:
$$Hom(u^{\otimes k},u^{\otimes l})=C_{kl}$$

We are therefore led to the conclusion in the statement.
\end{proof}

At the level of the examples now, we certainly have the easy groups from chapter 7, and we have as well their free versions constructed above, the result being:

\begin{theorem}
We have easy quantum groups as follows, with calligraphic standing for matching, meaning satisfying $\#\circ=\#\bullet$, as a weighted equality, in each block: 
$$
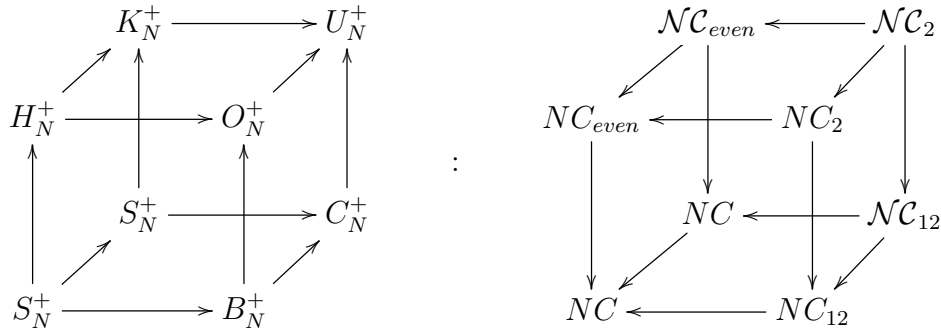
\xymatrix@R=16pt@C=15pt{
&K_N^+\ar[rr]&&U_N^+\\
H_N^+\ar[rr]\ar[ur]&&O_N^+\ar[ur]\\
&S_N^+\ar[rr]\ar[uu]&&C_N^+\ar[uu]\\
S_N^+\ar[uu]\ar[ur]\ar[rr]&&B_N^+\ar[uu]\ar[ur]
}\qquad\xymatrix@R=18pt@C=14pt{\\ \\ :}
\qquad\xymatrix@R=19.2pt@C=0pt{
&\mathcal{NC}_{even}\ar[dl]\ar[dd]&&\mathcal{NC}_2\ar[dl]\ar[ll]\ar[dd]\\
NC_{even}\ar[dd]&&NC_2\ar[dd]\ar[ll]\\
&NC\ar[dl]&&\mathcal{NC}_{12}\ar[dl]\ar[ll]\\
NC&&NC_{12}\ar[ll]}$$
Moreover, the quantum reflection groups $H_N^{s+}$ with $s\in\{1,2,\ldots,\infty\}$ fit on the left diagonal of the left cube, corresponding to the categories $NC^s$, fitting on the right cube. 
\end{theorem}

\begin{proof}
This is identical to the proof of the classical result, from chapter 7, with the diagrams, intertwiners and everything being identical, save for the very end in each case, with the categories generated by the relevant partitions being noncrossing, due to the fact that we removed the crossing axiom, in our updated Definition 15.15.
\end{proof}

In relation now with the easiness considerations from chapter 14, we first have:

\begin{theorem}
The basic orthogonal easy quantum groups are
$$
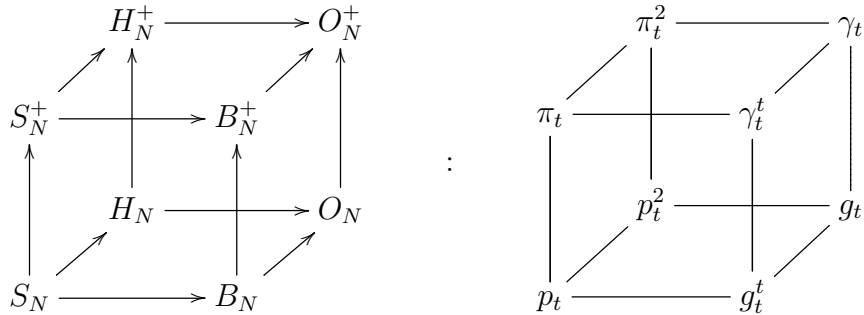
\xymatrix@R=16pt@C=15pt{
&H_N^+\ar[rr]&&O_N^+\\
S_N^+\ar[rr]\ar[ur]&&B_N^+\ar[ur]\\
&H_N\ar[rr]\ar[uu]&&O_N\ar[uu]\\
S_N\ar[uu]\ar[ur]\ar[rr]&&B_N\ar[uu]\ar[ur]
}\qquad\xymatrix@R=18pt@C=14pt{\\ \\ :}
\qquad\xymatrix@R=15pt@C=20pt{
&\pi_t^2\ar@{-}[rr]\ar@{-}[dd]&&\gamma_t\ar@{-}[dd]\\
\pi_t\ar@{-}[rr]\ar@{-}[dd]\ar@{-}[ur]&&\gamma_t^t\ar@{-}[dd]\ar@{-}[ur]\\
&p_t^2\ar@{-}[rr]\ar@{-}[uu]&&g_t\ar@{.}[uu]\\
p_t\ar@{-}[uu]\ar@{-}[ur]\ar@{-}[rr]&&g_t^t\ar@{-}[uu]\ar@{-}[ur]
}$$
with on the right being the asymptotic laws of truncated characters.
\end{theorem}

\begin{proof}
We already know the result for the bottom objects, from chapter 4, and the proof for the top objects is identical, using the Weingarten formula. Let us also mention that the orthogonal easy quantum groups were fully classified by Raum and Weber in \cite{rwe}, with the uniform classical and free objects being the above ones. See \cite{bsp}, \cite{rwe}.
\end{proof}

In the unitary case now, still in relation with the easiness considerations from chapter 14, and further clarifying things there, we have a similar result, as follows:

\begin{theorem}
The basic unitary easy quantum groups are
$$\xymatrix@R=16pt@C=15pt{
&K_N^+\ar[rr]&&U_N^+\\
S_N^+\ar[rr]\ar[ur]&&C_N^+\ar[ur]\\
&K_N\ar[rr]\ar[uu]&&U_N\ar[uu]\\
S_N\ar[uu]\ar[ur]\ar[rr]&&C_N\ar[uu]\ar[ur]
}\qquad\xymatrix@R=18pt@C=14pt{\\ \\ :}
\qquad\xymatrix@R=16pt@C=20pt{
&\Pi_t\ar@{-}[rr]\ar@{-}[dd]&&\Gamma_t\ar@{-}[dd]\\
\pi_t\ar@{-}[rr]\ar@{-}[dd]\ar@{-}[ur]&&\Gamma_t^t\ar@{-}[dd]\ar@{-}[ur]\\
&P_t\ar@{-}[rr]\ar@{-}[uu]&&G_t\ar@{.}[uu]\\
p_t\ar@{-}[uu]\ar@{-}[ur]\ar@{-}[rr]&&G_t^t\ar@{-}[uu]\ar@{-}[ur]
}$$
with on the right being the asymptotic laws of truncated characters.
\end{theorem}

\begin{proof}
The story here is similar to that of Theorem 15.18, with the remark however that at the classification level, things are far more complicated. See \cite{mwe}, \cite{rwe}.
\end{proof}

Finally, as explained in chapter 14, we cannot really merge our pairs of 3D cubes into 4D cubes. Instead, we can get rid of some of the objects, and we are led in this way to:

\begin{theorem}
The very basic easy quantum groups are
$$\xymatrix@R=16.1pt@C=15pt{
&K_N^+\ar[rr]&&U_N^+\\
H_N^+\ar[rr]\ar[ur]&&O_N^+\ar[ur]\\
&K_N\ar[rr]\ar[uu]&&U_N\ar[uu]\\
H_N\ar[uu]\ar[ur]\ar[rr]&&O_N\ar[uu]\ar[ur]
}\qquad\xymatrix@R=18pt@C=14pt{\\ \\ :}
\qquad\xymatrix@R=16pt@C=20pt{
&\Pi_t\ar@{-}[rr]\ar@{-}[dd]&&\Gamma_t\ar@{-}[dd]\\
b_t^2\ar@{-}[rr]\ar@{-}[dd]\ar@{-}[ur]&&\gamma_t\ar@{-}[dd]\ar@{-}[ur]\\
&P_t\ar@{-}[rr]\ar@{-}[uu]&&G_t\ar@{.}[uu]\\
p_t^2\ar@{-}[uu]\ar@{-}[ur]\ar@{-}[rr]&&g_t\ar@{-}[uu]\ar@{-}[ur]
}$$
with on the right being the asymptotic laws of truncated characters.
\end{theorem}

\begin{proof}
This follows indeed from Theorems 15.18 and 15.19.
\end{proof}

And with this, end of our discussion regarding quantum groups and easiness. Good to have this done, we have now a good understanding of the probabilistic easiness phenomena from chapter 14. But of course, many other things can be said, as a continuation, notably with De Finetti theorems, and we refer here to \cite{bcs}, \cite{cur}, \cite{csp}, \cite{dfr}, \cite{ksp}, \cite{liu}.

\section*{15b. Spheres, manifolds}

What is next? According to our Plan 15.2, doing some geometry, and then, guided by Magpie 15.1, we will do some probability, hopefully of truly advanced type. So, just a bit more algebra and abstractions, and after that probability, that is promised.

\bigskip

To start with, following \cite{bgo}, we can talk about free spheres, as follows:

\begin{definition}
We have free real and complex spheres, defined via
$$C(S^{N-1}_{\mathbb R,+})=C^*\left(x_1,\ldots,x_N\Big|x_i=x_i^*,\sum_ix_i^2=1\right)$$
$$C(S^{N-1}_{\mathbb C,+})=C^*\left(x_1,\ldots,x_N\Big|\sum_ix_ix_i^*=\sum_ix_i^*x_i=1\right)$$
with the symbol $C^*$ standing for universal enveloping $C^*$-algebra.
\end{definition}

Here the fact that these algebras are indeed well-defined comes from the following estimate, which shows that the biggest $C^*$-norms on these $*$-algebras are bounded:
$$||x_i||^2
=||x_ix_i^*||
\leq\Big|\Big|\sum_ix_ix_i^*\Big|\Big|
=1$$

As a first result now, regarding the above free spheres, we have:

\begin{theorem}
We have embeddings of compact quantum spaces, as follows,
$$\xymatrix@R=15mm@C=15mm{
S^{N-1}_{\mathbb R,+}\ar[r]&S^{N-1}_{\mathbb C,+}\\
S^{N-1}_\mathbb R\ar[r]\ar[u]&S^{N-1}_\mathbb C\ar[u]
}$$
and the spaces on top appear as liberations of the spaces on the bottom.
\end{theorem}

\begin{proof}
In practice, we must establish the following isomorphisms, where the symbol $C^*_{comm}$ stands for ``universal commutative $C^*$-algebra generated by'':
$$C(S^{N-1}_\mathbb R)=C^*_{comm}\left(x_1,\ldots,x_N\Big|x_i=x_i^*,\sum_ix_i^2=1\right)$$
$$C(S^{N-1}_\mathbb C)=C^*_{comm}\left(x_1,\ldots,x_N\Big|\sum_ix_ix_i^*=\sum_ix_i^*x_i=1\right)$$

It is enough to establish the second isomorphism. So, consider the second universal commutative $C^*$-algebra $A$ constructed above. Since the standard coordinates on $S^{N-1}_\mathbb C$ satisfy the defining relations for $A$, we have a quotient map of as follows:
$$A\to C(S^{N-1}_\mathbb C)$$

Conversely, let us write $A=C(S)$, by using the Gelfand theorem. The variables $x_1,\ldots,x_N$ become in this way true coordinates, providing us with an embedding $S\subset\mathbb C^N$. Also, the quadratic relations become $\sum_i|x_i|^2=1$, so we have $S\subset S^{N-1}_\mathbb C$. Thus, we have a quotient map $C(S^{N-1}_\mathbb C)\to A$, as desired, and this gives the results.
\end{proof}

By using the free spheres constructed above, we can now formulate:

\begin{definition}
A real algebraic manifold $X\subset S^{N-1}_{\mathbb C,+}$ is a closed quantum subspace defined, at the level of the corresponding $C^*$-algebra, by a formula of type
$$C(X)=C(S^{N-1}_{\mathbb C,+})\Big/\Big<f_i(x_1,\ldots,x_N)=0\Big>$$
for certain family of noncommutative polynomials, as follows:
$$f_i\in\mathbb C<x_1,\ldots,x_N>$$
We denote by $\mathcal C(X)$ the $*$-subalgebra of $C(X)$ generated by the coordinates $x_1,\ldots,x_N$. 
\end{definition}

As a basic example here, we have the free real sphere $S^{N-1}_{\mathbb R,+}$. The classical spheres $S^{N-1}_\mathbb C,S^{N-1}_\mathbb R$, and their real submanifolds, are covered as well by this formalism. We have as well as examples the compact quantum groups, embedded as follows:
$$G\subset U_N^+\ \implies \ G\subset U_N^+\subset S^{N^2-1}_{\mathbb C,+}\quad,\quad x_{ij}=\frac{u_{ij}}{\sqrt{N}}$$

At the level of the general theory, we have the following version of the Gelfand theorem:

\begin{theorem}
If $X\subset S^{N-1}_{\mathbb C,+}$ is an algebraic manifold, as above, we have
$$X_{class}=\left\{x\in S^{N-1}_\mathbb C\Big|f_i(x_1,\ldots,x_N)=0\right\}$$
and $X$ appears as a liberation of $X_{class}$.
\end{theorem}

\begin{proof}
This is something that we already met, in the context of the free spheres. In general, the proof is similar, by using the Gelfand theorem. Indeed, if we denote by $X_{class}'$ the manifold constructed in the statement, then we have a quotient map of $C^*$-algebras as follows, mapping standard coordinates to standard coordinates:
$$C(X_{class})\to C(X_{class}')$$

Conversely now, from $X\subset S^{N-1}_{\mathbb C,+}$ we obtain $X_{class}\subset S^{N-1}_\mathbb C$. Now since the relations defining $X_{class}'$ are satisfied by $X_{class}$, we obtain an inclusion $X_{class}\subset X_{class}'$. Thus, at the level of algebras of continuous functions, we have a quotient map of $C^*$-algebras as follows, mapping standard coordinates to standard coordinates:
$$C(X_{class}')\to C(X_{class})$$

Thus, we have constructed a pair of inverse morphisms, and we are done.
\end{proof}

Finally, once again at the level of the general theory, we have:

\begin{definition}
We agree to identify two real algebraic submanifolds $X,Y\subset S^{N-1}_{\mathbb C,+}$ when we have a $*$-algebra isomorphism between $*$-algebras of coordinates
$$f:\mathcal C(Y)\to\mathcal C(X)$$
mapping standard coordinates to standard coordinates.
\end{definition}

Observe that this convention, which is something technical, coming from amenability, fits with the similar convention from Definition 15.5, via the embeddings $G\subset S^{N^2-1}_{\mathbb C,+}$ mentioned after Definition 15.23. Observe also that our whole formalism, Definitions 15.23 and 15.25, fixes some chronic bugs of Gelfand duality, that we preferred not to talk about, in chapters 8 and 13. Well, now you know, and very nice all this.

\bigskip

With this discussed, let us go back now to the spheres, and try to integrate over them. For this purpose, it is convenient to deal at the same time with all 4 spheres in Theorem 15.22, that we will call ``basic spheres'', and generically denote $S$. With the convention that the corresponding unitary quantum groups are denoted $U$, we first have:

\begin{proposition}
For the basic spheres, we have a diagram as follows,
$$\xymatrix@R=50pt@C=50pt{
C(S)\ar[r]^\alpha\ar[d]^\gamma&C(S)\otimes C(U)\ar[d]^{\gamma\otimes id}\\
C(U)\ar[r]^\Delta&C(U)\otimes C(U)
}$$
where on top $\alpha(x_i)=\sum_jx_j\otimes u_{ji}$, and on the left $\gamma(x_i)=u_{1i}$.
\end{proposition}

\begin{proof}
The diagram in the statement commutes indeed on the standard coordinates $x_i$. Thus by linearity and multiplicativity, the whole the diagram commutes.
\end{proof}

In practice now, the map $\alpha$ constructed above is what is called a coaction, and skipping some discussion here, the conclusion is that we have a good relationship $S\leftrightarrow U$, a bit like in the classical case. Now based on this, in analogy with what happens in the classical case, we can introduce the uniform integration over our spheres $S$, as follows:

\index{integration over spheres}

\begin{definition}
We endow each of the algebras $C(S)$ with the functional
$$\int_S:C(S)\to C(U)\to\mathbb C$$
obtained by composing the morphism $x_i\to u_{1i}$ with the Haar integration of $C(U)$.
\end{definition}

In order to efficiently integrate now over our generic spheres $S$, in the lack of some trick like spherical coordinates, we can use the Weingarten formula:

\index{integration over spheres}
 
\begin{theorem}
The integration over the basic spheres is given by
$$\int_Sx_{i_1}^{e_1}\ldots x_{i_k}^{e_k}=\sum_{\pi,\sigma\in D(k)}\delta_\sigma(i)W_{kN}(\pi,\sigma)$$
where $D\subset P$ are the partitions for $U$, and $W_{kN}=G_{kN}^{-1}$, with $G_{kN}(\pi,\sigma)=N^{|\pi\vee\sigma|}$. 
\end{theorem}

\begin{proof}
According to our conventions, the integration over $S$ is a particular case of the integration over $U$, via $x_i=u_{1i}$. By using the Weingarten formula for $U$, we get:
\begin{eqnarray*}
\int_Sx_{i_1}^{e_1}\ldots x_{i_k}^{e_k}
&=&\int_Uu_{1i_1}^{e_1}\ldots u_{1i_k}^{e_k}\\
&=&\sum_{\pi,\sigma\in D(k)}\delta_\pi(1)\delta_\sigma(i)W_{kN}(\pi,\sigma)\\
&=&\sum_{\pi,\sigma\in D(k)}\delta_\sigma(i)W_{kN}(\pi,\sigma)
\end{eqnarray*}

Thus, we are led to the formula in the statement.
\end{proof}

Still following \cite{bgo}, we have the following key result:

\index{ergodicity}

\begin{theorem}
The integration functional of $S$ has the ergodicity property
$$\left(id\otimes\int_U\right)\alpha(x)=\int_Sx$$
where $\alpha:C(S)\to C(S)\otimes C(U)$ is the universal affine coaction map.
\end{theorem}

\begin{proof}
In the real case, $x_i=x_i^*$, it is enough to check the equality in the statement on an arbitrary product of coordinates, $x_{i_1}\ldots x_{i_k}$. The left term is as follows:
\begin{eqnarray*}
\left(id\otimes\int_U\right)\alpha(x_{i_1}\ldots x_{i_k})
&=&\sum_{j_1\ldots j_k}x_{j_1}\ldots x_{j_k}\int_Uu_{j_1i_1}\ldots u_{j_ki_k}\\
&=&\sum_{j_1\ldots j_k}\ \sum_{\pi,\sigma\in D(k)}\delta_\pi(j)\delta_\sigma(i)W_{kN}(\pi,\sigma)x_{j_1}\ldots x_{j_k}\\
&=&\sum_{\pi,\sigma\in D(k)}\delta_\sigma(i)W_{kN}(\pi,\sigma)\sum_{j_1\ldots j_k}\delta_\pi(j)x_{j_1}\ldots x_{j_k}\\
&=&\sum_{\pi,\sigma\in D(k)}\delta_\sigma(i)W_{kN}(\pi,\sigma)
\end{eqnarray*}

Now by comparing with the formula in Theorem 15.28, we obtain the result. As for the proof in the complex case, this is similar, by adding exponents everywhere.
\end{proof} 

Still following \cite{bgo}, we have as well the following useful abstract result:

\index{integration over spheres}

\begin{theorem}
There is a unique positive unital trace $tr:C(S)\to\mathbb C$ satisfying
$$(tr\otimes id)\alpha(x)=tr(x)1$$
where $\alpha$ is the coaction map of the corresponding quantum isometry group,
$$\alpha:C(S)\to C(S)\otimes C(U)$$
and this is the canonical integration, as constructed in Definition 15.27.
\end{theorem}

\begin{proof}
First of all, it follows from the Haar integral invariance condition for $U$ that the canonical integration has indeed the invariance property in the statement, namely:
$$(tr\otimes id)\alpha(x)=tr(x)1$$

In order to prove now the uniqueness, let $tr$ be as in the statement. We have:
$$tr\left(id\otimes\int_U\right)\alpha(x)
=\int_U(tr\otimes id)\alpha(x)
=\int_U(tr(x)1)
=tr(x)$$

On the other hand, according to Theorem 15.29, we have as well:
$$
tr\left(id\otimes\int_U\right)\alpha(x)
=tr\left(\int_Sx\right)
=\int_Sx$$

We therefore conclude that $tr$ equals the standard integration, as claimed.
\end{proof}

Getting now to probabilistic aspects, we know from chapter 11 that for the classical spheres the laws of the standard coordinates, which are called hyperspherical laws, become normal and independent, with $N\to\infty$. Here is the free analogue of this result:

\begin{theorem}
For the free spheres $S^{N-1}_{\mathbb R,+}$, $S^{N-1}_{\mathbb C,+}$, the rescaled coordinates 
$$y_i=\sqrt{N}x_i$$
become semicircular/circular and free, in the $N\to\infty$ limit.
\end{theorem}

\begin{proof}
In the real case, the Weingarten formula from Theorem 15.28, coupled with the standard fact that the Gram and Weingarten matrices are asymptotically diagonal, $G_{kN}\simeq diag(N^{k/2})$ and $W_{kN}\simeq diag(N^{-k/2})$, gives the following estimate:
$$\int_{S^{N-1}_{\mathbb R,+}}x_{i_1}\ldots x_{i_k}
=\sum_{\pi,\sigma\in NC_2(k)}\delta_\sigma(i)W_{kN}(\pi,\sigma)
\simeq N^{-k/2}\sum_{\sigma\in NC_2(k)}\delta_\sigma(i)$$

With this formula in hand, we can compute the asymptotic moments of each coordinate $x_i$. Indeed, by setting $i_1=\ldots=i_k=i$, all Kronecker symbols are 1, and we obtain:
$$\int_{S^{N-1}_{\mathbb R,+}}x_i^k\,dx\simeq N^{-k/2}|NC_2(k)|=N^{-k/2}M_k(\gamma_1)$$

Thus the rescaled coordinates $y_i=\sqrt{N}x_i$ become semicircular in the $N\to\infty$ limit, and the asymptotic freeness result follows as well, from the above general joint moment estimate. As for the proof in the complex case, this is similar. See \cite{bgo}.
\end{proof}

So long for the free spheres, and integration over them. This is of course just the tip of the iceberg, and the next step is to consider quotient spaces of type $X=G_N/G_{N-M}$, with $G$ being an easy quantum group, which unify the spheres with the quantum groups. Then, we can look at more complicated quotient spaces, with a class of fairly general examples here, covering those before with $L=M$, being as follows:
$$X=(G_M\times G_N)\big/(G_L\times G_{M-L}\times G_{N-L})$$

And there is even more, because these latter quotient spaces are at the core of what is called ``affine homogeneous spaces'', and again, a bit as before for the quantum groups and the spheres, a Weingarten formula can be worked out, along with some basic asymptotic consequences, involving the basic limiting laws from classical and free probability. For more on all this we refer to the easy geometry literature, \cite{bgo} and its follow-ups.
 
\section*{15c. Hyperspherical laws}

The problem now, which is highly non-trivial, is that of computing the moments of the coordinates of the free sphere at fixed values of $N\in\mathbb N$. And the answer here, from \cite{bcz}, based on advanced quantum group and calculus techniques, is as follows:

\index{free hyperspherical law}
\index{special functions}
\index{deformed quantum groups}

\begin{theorem}
The moments of the free hyperspherical law are given by
$$\int_{S^{N-1}_{\mathbb R,+}}x_1^{2l}=\frac{1}{(N+2)^l}\cdot\frac{q+1}{q-1}\cdot\frac{1}{l+1}\sum_{r=-l-1}^{l+1}(-1)^r\binom{2l+2}{l+r+1}\frac{r}{1+q^r}$$
where $q\in [-1,0)$ is such that $q+q^{-1}=-N$.
\end{theorem}

\begin{proof}
This is something quite tricky, sort of thing that a physicist does in 20 minutes, and mathematicians spend 2 years after in understanding, the idea being:

\medskip

(1) In order to get into the mood, let us first see what our formula tells us, at $l=1$:
\begin{eqnarray*}
M_2
&=&\frac{1}{N+2}\cdot\frac{q+1}{q-1}\cdot\frac{1}{2}\sum_{r=-2}^2(-1)^r\binom{4}{r+2}\frac{r}{1+q^r}\\
&=&\frac{1}{N+2}\cdot\frac{q+1}{q-1}\cdot\frac{1}{2}\left[-\frac{2}{1+q^{-2}}+\frac{4}{1+q^{-1}}-\frac{4}{1+q}+\frac{2}{1+q^2}\right]\\
&=&\frac{1}{N+2}\cdot\frac{q+1}{q-1}\left[\frac{1-q^2}{1+q^2}+\frac{2(q-1)}{1+q}\right]\\
&=&\frac{q+1}{N+2}\left[\frac{2}{1+q}-\frac{1+q}{1+q^2}\right]
=\frac{q+1}{N+2}\cdot\frac{q^2-2q+1}{(1+q)(1+q^2)}\\
&=&\frac{1}{N+2}\left[1-\frac{2q}{1+q^2}\right]
=\frac{1}{N+2}\left[1-2\left(q+\frac{1}{q}\right)^{-1}\right]\\
&=&\frac{1}{N+2}\left[1+\frac{2}{N}\right]=\frac{1}{N}
\end{eqnarray*}

Which was not particularly fun, hope you agree with me. Note by the way that this is the correct answer, because we have $x_1^2+\ldots+x_N^2=1$, uniformly.

\medskip

(2) At $l=2$ now, the computation becomes particularly tough, but still can be done by hand. Then at $l=3$, things complicated. In any case, here are the formulae:
$$M_4=\frac{2}{N(N+1)}\quad,\quad 
M_6=\frac{5N-7}{N(N+1)(N^2-2)}$$

And, are these correct. In answer, the Weingarten formula in Theorem 15.28 gives:
$$M_4=Tr\left[\begin{pmatrix}1&1\\1&1\end{pmatrix}
\begin{pmatrix}N^2&N\\N&N^2\end{pmatrix}^{-1}\right]
=\frac{2}{N(N+1)}$$
$$M_6=Tr\left[\begin{pmatrix}
1&1&1&1&1\\
1&1&1&1&1\\
1&1&1&1&1\\
1&1&1&1&1\\
1&1&1&1&1
\end{pmatrix}
\begin{pmatrix}
N^3&N^2&N^2&N&N^2\\
N^2&N^3&N&N^2&N\\
N^2&N&N^3&N^2&N\\
N&N^2&N^2&N^3&N^2\\
N^2&N&N&N^2&N^3
\end{pmatrix}^{-1}\right]
=\frac{5N-7}{N(N+1)(N^2-2)}$$

(3) Summarizing, our formula looks good. Observe also that, as a byproduct of our computations, the expectation, variance, skewness and kurtosis are:
$$E=0\quad,\quad V=\frac{1}{N}\quad,\quad\gamma=0\quad,\quad\kappa=\frac{2N}{N+1}$$

As a comment here, $E,V,\gamma$ equal those from chapter 11 for the classical hyperspherical laws, but the kurtosis, which was there $\kappa=3N/(N+2)$, has changed. Which is not surprising, because we should have $\kappa\simeq3$ for classical laws, and $\kappa\simeq2$ for free laws.

\medskip

(4) As yet another verification, this time in the general case, $l\in\mathbb N$, up to some noise, our formula predicts a moment $M_{2l}$ having as denominator the following product:
$$P=\prod_{r=-l-1}^{l+1}(1+q^r)$$ 

And this fits with the product of Chebycheff polynomials from the Di Francesco formula for $\det(G_{kN})$ for the quantum group $O_N^+$, from chapter 7, and we will leave some verifications here, first at small values of $l$, and then in general, as an exercise.

\medskip

(5) Next, the $N\to\infty$ asymptotics of our formula, meant to fit with Theorem 15.31, are less obvious to work out. And more on this, later in the proof.

\medskip

(6) Finally, before getting into the proof, let us mention that \cite{bcz} contains as well a number of useful reformulations and corollaries of the formula in the statement, notably with the fact that the support is as follows, with the density being analytic:
$$S=\left[-\frac{2}{\sqrt{N+2}}\ ,\ \frac{2}{\sqrt{N+2}}\right]$$ 

(7) Getting now to the proof, the idea will be that of transporting our computation from space to space, no less than 4 times, up to reaching a tough computation on the unit circle $\mathbb T$, and solving that with some help from advanced orthogonal polynomials:
$$S^{N-1}_{\mathbb R,+}
\ \longrightarrow\ O_N^+ 
\ \longrightarrow\ O_F^+=SU_2^q
\ \longrightarrow\ B(l^2(\mathbb N))
\ \longrightarrow\ \mathbb T$$

(8) Regarding the first transport, via $x_1=u_{11}$, we know this, from the previous section. Next, to any $F\in GL_N(\mathbb R)$ satisfying $F^2=1$ let us associate the following algebra:
$$C(O_F^+)=C^*\left((u_{ij})_{i,j=1,\ldots,N}\Big|u=F\bar{u}F={\rm unitary}\right)$$

Observe that we have $O_{I_N}^+=O_N^+$. In general, the above algebra satisfies Woronowicz's generalized axioms in \cite{wo1}, which do not include the strong antipode axiom $S^2=id$. 

\medskip

(9) At $N=2$, up to a trivial equivalence relation on the matrices $F$, and on the quantum groups $O_F^+$, we can assume that $F$ is as follows, with $q\in [-1,0)$:
$$F=\begin{pmatrix}0&\sqrt{-q}\\
1/\sqrt{-q}&0\end{pmatrix}$$

Our claim is that for this matrix $F$ we have an identification $O_F^+=SU^q_2$. Indeed, the relations $u=F\bar{u}F$ tell us that $u$ must be of the following form:
$$u=\begin{pmatrix}\alpha&-q\gamma^*\cr \gamma&\alpha^*\end{pmatrix}$$

Thus $C(O_F^+)$ is the universal algebra generated by two elements $\alpha,\gamma$, with the relations making the above matrix $u$ unitary. But these unitarity conditions are as follows:
$$\alpha\gamma=q\gamma\alpha\quad,\quad
\alpha\gamma^*=q\gamma^*\alpha\quad,\quad 
\gamma\gamma^*=\gamma^*\gamma$$
$$\alpha^*\alpha+\gamma^*\gamma=1\quad,\quad 
\alpha\alpha^*+q^2\gamma\gamma^*=1$$

We recognize here the relations in \cite{wo1} defining the algebra $C(SU^q_2)$, and it follows that we have an isomorphism of Hopf $C^*$-algebras, as follows:
$$C(O_F^+)\simeq C(SU^q_2)$$

(10) Now back to the general case, let us try to understand the integration over $O_F^+$. Given a noncrossing pairing $\pi\in NC_2(2k)$ and an index $i=(i_1,\ldots,i_{2k})$, we set:
$$\delta_\pi^F(i)=\prod_{s\in\pi}F_{i_{s_l}i_{s_r}}$$

Here the product is over all strings $s=\{s_l\curvearrowright s_r\}$ of $\pi$. Our claim is that the following family of vectors, with $\pi\in NC_2(2k)$, spans the space of fixed vectors of $u^{\otimes 2k}$:
$$\xi_\pi=\sum_i\delta_\pi^F(i)e_{i_1}\otimes\ldots\otimes e_{i_{2k}}$$ 

Indeed, having $\xi_\cap$ fixed by $u^{\otimes 2}$ is equivalent to assuming that $u=F\bar{u}F$ is unitary. By using now the above vectors, we obtain the following Weingarten formula:
$$\int_{O_F^+}u_{i_1j_1}\ldots u_{i_{2k}j_{2k}}=\sum_{\pi\sigma}\delta_\pi^F(i)\delta_\sigma^F(j)W_{kN}(\pi,\sigma)$$

(11) With these preliminaries in hand, let us start the computation. Let $q\in [-1,0)$ be such that $q+q^{-1}=-N$. Our claim is that we have the following formula:
$$\int_{O_N^+}\varphi(\sqrt{N+2}\,u_{ij})=\int_{SU^q_2}\varphi(\alpha+\alpha^*+\gamma-q\gamma^*)$$

Indeed, the moments of the variable on the left are given by:
$$\int_{O_N^+}u_{ij}^{2k}=\sum_{\pi\sigma}W_{kN}(\pi,\sigma)$$

On the other hand, the moments of the variable on the right, which in terms of the fundamental corepresentation $v=(v_{ij})$ is given by $w=\sum_{ij}v_{ij}$, are given by:
$$\int_{SU^q_2}w^{2k}=\sum_{ij}\sum_{\pi\sigma}\delta_\pi^F(i)\delta_\sigma^F(j)W_{kN}(\pi,\sigma)$$

We deduce that $w/\sqrt{N+2}$ has the same moments as $u_{ij}$, which proves our claim.

\medskip

(12) Next, before jumping into computations over $SU^q_2$, it is convenient to introduce a useful degree of flexibility, coming from an extra parameter. The point indeed is that the following variables over $SU^q_2$ have the same law, whenever $AB=-q$:
$$\alpha+\alpha^*+\gamma-q\gamma^*\sim \alpha+\alpha^*+A\gamma+B\gamma^*$$

Indeed, this comes via Weingarten, applied as in (11). Thus, as a conclusion, the law that we want to compute is that of the following variable, with  $AB=-q$:
$$\sqrt{N+2}\,u_{ij}\sim \alpha+\alpha^*+A\gamma+B\gamma^*$$

(13) In order to do now the computation over $SU^q_2$, we can use a matrix model due to Woronowicz \cite{wo1}, where the standard generators $\alpha,\gamma$ are mapped as follows:
$$\pi_u(\alpha)e_k=\sqrt{1-q^{2k}}e_{k-1}\quad,\quad 
\pi_u(\gamma)e_k=uq^k e_k$$

Here $u\in\mathbb T$ is a parameter, and $(e_k)$ is the standard basis of $l^2(\mathbb N)$. The point with this representation is that it allows the computation of the Haar functional. Indeed, if $D$ is the diagonal operator given by $D(e_k)=q^{2k}e_k$, then the formula is as follows:
$$\int _{SU^q_2}x=(1-q^2)\int_{\mathbb T}tr(D\pi_u(x))\frac{du}{2\pi iu}$$

(14) Thus, the law of the variable that we are interested in is of the following form:
$$\int_{SU^q_2}\varphi(\alpha+\alpha^*+A\gamma+B\gamma^*)=(1-q^2)\int_{\mathbb T}tr(D\varphi(M))\frac{du}{2\pi iu}$$

To be more precise, this formula holds indeed, with $M$ being as follows:
$$M(e_k)=\sqrt{1-q^{2k+2}}\,e_{k+1}+q^k(Au+Bu^{-1})e_k+\sqrt{1-q^{2k}}\,e_{k-1}$$

(15) Next, let us perform a change of basis. Consider another copy of $l^2(\mathbb N)$, with orthonormal basis denoted $\{e_n'\}$, and define a linear map $J:l^2(\mathbb N)\to l^2(\mathbb N)$ as follows:
$$J(e_k')=\sqrt{1-q^2}\sqrt{1-q^4}\ldots\sqrt{1-q^{2k}}\,e_k$$

In terms of $D'=J^{-1}DJ$ and $M'=J^{-1}MJ$, our formula in (14) becomes:
$$\int\varphi(\alpha+\alpha^*+A\gamma+B\gamma^*)=(1-q^2)\int_{\mathbb T}tr(D'\varphi(M'))\frac{du}{2\pi iu}$$

Since $J$ is diagonal, we have $D'(e_k')=q^{2k}e_k'$. As for $M'$, this is given by:
\begin{eqnarray*}
M'(e_k)
&=&J^{-1}MJ(e_k')\\
&=&\sqrt{1-q^2}\sqrt{1-q^4}\ldots\sqrt{1-q^{2k}}J^{-1}M(e_k)\\
&=&e_{k+1}'+q^k(Au+Bu^{-1})e_k'+(1-q^{2k})e_{k-1}'
\end{eqnarray*}

(16) As a conclusion to all this, by removing all prime signs, and by remembering as well from (12) what we wanted to do, we have a formula as follows:
$$\int_{O_N^+}\varphi(\sqrt{N+2}\,u_{ij})=(1-q^2)\int_{\mathbb T}tr(D\varphi(M))\frac{du}{2\pi iu}$$

To be more precise, this formula holds indeed, with $M$ being as follows:
$$M(e_k)=e_{k+1}+q^k(Au+Bu^{-1})e_k+(1-q^{2k})e_{k-1}$$

(17) This was for the easy part. The hard part, inspired from Koelink-Verding \cite{kve}, is that of computing the above integral. First, the $q$-shifted factorial is given by:
$$(a;q)_k=(1-a)(1-qa)\ldots(1-q^{k-1}a)$$

Next, the multivariable version of this $q$-shifted factorial is given by:
$$(a_1,\ldots,a_r;q)_k=(a_1;q)_k\ldots(a_r;q)_k$$

Finally, the $q$-hypergeometric series is given by the following formula:
$${}_r\varphi_s\left(\begin{matrix}{a_1,\ldots,a_r}\cr {b_1,\ldots,b_s}\end{matrix}\,;q,z\right)
=\sum_{k=0}^\infty\left( (-1)^k q^{k(k-1)/2}\right)^{s-r+1}\frac{(a_1,\ldots,a_r;q)_k}{(q,b_1,\ldots,b_s;q)_k}\,z^k$$

(18) With this, we can now introduce the Al-Salam-Chihara polynomials, depending on two parameters $a,b$, which are defined by the following formula:
$$Q_k(x)=\frac{(ab;q)_k}{a^k}\,{}_3\varphi_2\left(\begin{matrix}q^{-k},az,az^{-1}\cr ab,0\end{matrix}\,;q,q\right)$$

Here we use the convenient parameterization $2x=z+z^{-1}$. These polynomials are known to satisfy the following recurrence relation, and we refer to \cite{awi} for this:
$$2xQ_k(x)=Q_{k+1}(x)+q^k(a+b)Q_k(x)+(1-q^k)(1-abq^{k-1})Q_{k-1}(x)$$

(19) Now by getting back to (16), let us set $a=Au$ and $b=Bu^{-1}$. In terms of these new parameters $a,b$, the formula there for our tridiagonal operator $M$ reads:
$$M(e_k)=e_{k+1}+q^k(a+b)e_k+(1-q^{2k})e_{k-1}$$

Now recall from (16) that we have $AB=-q$, and so $ab=-q$. Thus the recurrence formula for the corresponding Al-Salam-Chihara polynomials is:
$$2xQ_k(x)=Q_{k+1}(x)+q^k(a+b)Q_k(x)+(1-q^{2k})Q_{k-1}(x)$$

We deduce that the collection $\{Q_k(x)\}$ of Al-Salam-Chihara polynomials evaluated at $x$ play the role of eigenvectors for the operator $M$, with corresponding eigenvalue $2x$.

\medskip

(20) And with this, we are now into orthogonal polynomials, and some 3 pages of calculations, based on the Poisson kernel for the Al-Salam-Chihara polynomials, and on some remarkable formulae from \cite{ars}, at black belt level ${}_8\varphi_7$, lead to the formula in the statement. For details on all this, and for more, we refer to \cite{bcz}. 
\end{proof}

\section*{15d. Hypergeometric laws}

With the above computation done, end of the world, you would say. Well, not really, because abstract algebra can strike, and there will be a further twist to the story. We will be explaining this in what follows, which is material from the paper \cite{bbc}. 

\bigskip

Let us start with something routine. Still guided by Magpie 15.1, now that we stole the hyperspherical laws, let us steal as well the hypergeometric ones:

\index{free hypergeometric law}

\begin{definition}
The free hypergeometric law of parameters $(N,m,p)$ is the law of
$$\sigma_{mp}=\sum_{i=1}^m\sum_{j=1}^pu_{ij}$$
over the quantum permutation group $S_N^+$.
\end{definition}

Observe the similarity with the formula from chapter 12 for the classical hypergeometric laws, involving the usual permutation group $S_N$. We can actually exploit this similarity for getting some useful pieces of information about our laws, as follows:

\begin{theorem}
The moments of the free hypergeometric laws are given by
$$M_k=\sum_{\pi,\nu\in NC(k)}W_{kN}(\pi,\nu)m^{|\pi|}p^{|\nu|}$$
with $W_{kN}$ being the Weingarten matrix for $S_N^+$. In particular we have
$$E=\frac{mp}{N}\quad,\quad V=\frac{mp(N-m)(N-p)}{N^2(N-1)}\quad,\quad 
\gamma=\frac{(N-2m)(N-2p)\sqrt{N-1}}{(N-2)\sqrt{mp(N-m)(N-p)}}$$
exactly as for the usual hypergeometric laws. As for the kurtosis, this is
$$\kappa=\frac{P(N,m,p)}{N^{10}(N-1)^{11}(N-2)^6(N^2-3N+1)m^2p^2(N-m)^2(N-p)^2}$$
with $P(N,m,p)$ being a certain polynomial, which might absorb some bottom terms.
\end{theorem}

\begin{proof}
According to the Weingarten formula for $S_N^+$, we have:
\begin{eqnarray*}
M_k
&=&\int_{S_N^+}\sum_{i_1=1}^m\ldots\sum_{i_k=1}^m
\sum_{j_1=1}^p\ldots\sum_{j_k=1}^pu_{i_1j_1}\ldots u_{i_pj_p}\\
&=&\sum_{i_1=1}^m\ldots\sum_{i_k=1}^m
\sum_{j_1=1}^p\ldots\sum_{j_k=1}^p
\sum_{\pi,\nu\in NC(k)}\delta_\pi(i)\delta_\nu(j)W_{kN}(\pi,\nu)\\
&=&\sum_{\pi,\nu\in NC(k)}W_{kN}(\pi,\nu)
\sum_{i_1=1}^m\ldots\sum_{i_k=1}^m\delta_\pi(i)
\sum_{j_1=1}^p\ldots\sum_{j_k=1}^p\delta_\nu(j)\\
&=&\sum_{\pi,\nu\in NC(k)}W_{kN}(\pi,\nu)m^{|\pi|}p^{|\nu|}
\end{eqnarray*}

Now the point is that, since we have $P(k)=NC(k)$ at $k=1,2,3$, the Weingarten matrices for $S_N$ and $S_N^+$ coincide at $k=1,2,3$, so we have that $E,V,\gamma$ figures, from chapter 12. Finally, regarding the kurtosis, yes I know that's not professional, that denominator is the Di Francesco determinant of $G_{4N}$, divided by the square of the variance $V$.
\end{proof}

Still following \cite{bbc}, we have the following basic asymptotic study:

\begin{theorem}
We have the following asymptotic results,
\begin{enumerate}
\item With $N,m,p\to\infty$ and $mp/N\to t\in(0,\infty)$, we have:
$$\sigma_{mp}\sim\pi_t$$

\item With $M,m,p\to\infty$ and $m/N\to\nu\in (0,1)$ and $p/N\to 0$, we have:
$$\frac{\sigma_{mp}-p\nu}{\sqrt{p\nu(1-\nu)}}\sim\gamma_1$$
\end{enumerate}
with $\sigma_{mp}$ following the free hypergeometric law of parameters $(N,m,p)$.
\end{theorem}

\begin{proof}
This is quite standard, by using the Weingarten formula for the moments of the free hypergeometric laws, from Theorem 15.34, namely:
$$M_k=\sum_{\pi,\nu\in NC(k)}W_{kN}(\pi,\nu)m^{|\pi|}p^{|\nu|}$$

Indeed, we have the following standard estimate, for the Weingarten matrix:
$$W_{kN}(\pi,\nu)= 
\begin{cases}N^{-|\pi|}+O(N^{-|\pi|-1})&{\rm if}\ \pi=\nu\\ 
O(N^{|\pi\vee\nu|-|\pi|-|\nu|})&{\rm if}\ \pi\neq\nu
\end{cases}$$

Now in the regime $N,m,p\to\infty$ and $mp/N\to t\in(0,\infty)$, it follows that we have:
$$W_{kN}(\pi,\nu)m^{|\pi|}p^{|\nu|} \to 
\begin{cases}
t^{|\pi|}&{\rm if}\ \pi=\nu\\ 
0&{\rm if}\ \pi\neq\nu 
\end{cases}$$

Thus, in this regime we obtain indeed the moments of the Marchenko-Pastur law:
$$M_k=\sum_{\pi\in NC(k)}t^{|\pi|}$$

As for the second assertion, the proof is a bit more technical. See \cite{bbc}.
\end{proof}

With this discussed, time now for abstract algebra to strike. We have indeed the following result from \cite{bbc}, which was at that time something quite unexpected:

\begin{theorem}
The projective orthogonal and unitary quantum groups coincide, $PO_n^+=PU_n^+$, and appear as a twist of the quantum permutation group $S_{n^2}^+$,
$$PO_n^+=PU_n^+\sim S_{n^2}^+$$
and with this being something having no classical counterpart.
\end{theorem}
 
\begin{proof}
As mentioned, this is something quite surprising, the idea being as follows:

\medskip

(1) To start with, the projective version of $(A,u)$ is the pair $(PA,v)$, with $PA=<v_{ij}>$ and $v=u\otimes\bar{u}$. The inclusion $O_n^+\subset U_n^+$ produces then an inclusion $PO_n^+\subset PU_n^+$, and since the main character on both sides follows $\pi_1$, we have an isomorphism:
$$PO_n^+=PU_n^+$$

Observe that this isomorphism has no classical counterpart, because $O_n,U_n$ have nothing much to do with each other, and so do their projective versions $PO_n,PU_n$.

\medskip

(2) Next, the above isomorphism is even better understood in the context of the quantum permutation groups $S_X^+$ of the finite quantum spaces $X$. Indeed, with the convention that $M_n$ is the quantum space given by $C(M_n)=M_n(\mathbb C)$, we have:
$$PO_n^+=PU_n^+=S_{M_n}^+$$

(3) Next, the fact that we have $M_n\sim\{1,\ldots,n^2\}$, twist of quantum spaces, raises the possibility for the corresponding quantum permutation groups to be related by a twisting operation. So, putting everything together, we would reach to a scheme as follows:
$$PO_n^+=PU_n^+=S_{M_n}^+\sim S_{n^2}^+$$ 

(4) And, is is true or not. At $n=2$ things are quite special, with the groups $SU_2$, $SO_3$ and their magic being involved, and it was known, basically since the probabilistic computations in \cite{bco}, that we have isomorphims as follows, confirming the above:
$$PO_2^+=PU_2^+=PU_2=SO_3\sim SO_3^{-1}\simeq S_4^+$$

(5) Now the point is that this holds in general too, as explained in \cite{bbc}, the reason behind this being the fact that the quantum automorphism group of a twisted group algebra is the twist of the quantum automorphism group of the original group algebra.

\medskip

(6) Finally, let us mention that all this is of potential key important in physics, showing that in the context of the Standard Model, in Chamseddine-Connes formulation \cite{cco}, as analysed by Bhowmick, Andrea and Dabrowski in \cite{bdd}, acting on the QED and QCD parts are twists of $S_4^+,S_9^+$. Which is quite interesting, suggesting that QED and QCD, suitably twisted, might be some sort of Yang-Mills theories based on $S_4^+,S_9^+$, respectively.
\end{proof}

Quite interesting all this, we are now deep into quarks. Still following \cite{bbc}, we have the following remarkable probabilistic consequence of the above physics trickery:

\begin{theorem}
The following families of variables have the same joint law,
$$\Big\{v_{ij}^2\Big|i,j=1,\ldots,n\Big\}\in C(O_n^+)\ \sim\ \Big\{\rho_{ij}=\frac{1}{n}\sum_{ab}u_{ia,jb}\Big|i,j=1,\ldots,n\Big\}\in C(S_{n^2}^+)$$
so in particular, the squared free hyperspherical law of parameter $n$ coincides with the $1/n$ rescaling of the free hypergeometric law of parameters $(n^2,n,n)$.
\end{theorem}

\begin{proof}
This comes from Theorem 15.36, but as pointed out in \cite{bbc}, thinking a bit retrospectively, this can be recovered as well directly, via Weingarten, as follows:
$$\int_{O_n^+}v_{ij}^{2k}=\sum_{\pi,\nu\in NC_2(2k)}W_{2k,n}(\pi,\nu)$$
$$\int_{S_{n^2}^+}\rho_{ij}^k=\sum_{\pi,\nu\in NC_2(2k)}n^{|\pi'|+|\sigma'|-k}W_{k,n^2}(\pi',\nu')$$

To be more precise, here $\pi\to\pi'$ is the shrinking operation for the pairings, discussed in chapter 7, in the context of Gram determinants, and this gives the result.
\end{proof}

Observe now that the squared free hyperspherical law of parameter $n$ is something that we know well, having struggled with it for pages in a row, in the previous section. What a comeback. In other words, we are dealing here with a generalization of that.

\bigskip

As a puzzling consequence of the above remarkable coincidence, we have:

\begin{theorem}
The squared hypespherical law of parameter $n$, having
$$E=\frac{1}{n}\quad,\quad V=\frac{2n-2}{n^2(n+2)}\quad,\quad\gamma=\frac{n-2}{n+4}\sqrt{\frac{8(n+2)}{n-1}}$$
and the $1/n$ rescaling of the hypergeometric law of parameters $(n^2,n,n)$, having
$$E=\frac{1}{n}\quad,\quad V=\frac{n-1}{n^2(n+1)}\quad,\quad\gamma=\frac{(n-2)^2}{n^3(n^2-2)}\sqrt{\frac{n+1}{n-1}}$$
have as common free analogue the free law discussed above, having
$$E=\frac{1}{n}\quad,\quad V=\frac{n-1}{n^2(n+1)}\quad,\quad\gamma=\frac{(n-2)^2}{n^3(n^2-2)}\sqrt{\frac{n+1}{n-1}}$$
and with the corresponding $3$ kurtoses being distinct.
\end{theorem}

\begin{proof}
The classical formulae come indeed from our results in chapter 12, and the free formulae come from Theorem 15.32 and its proof, or from Theorem 15.34.
\end{proof}

Finally, still talking violations of Bercovici-Pata, the above remains something quite specialized. At a more basic level, the main characters of the classical and free tori are related by the Meixner/free Meixner correspondence, discussed in \cite{ans}, \cite{bwy}.

\section*{15e. Exercises}

This was a quite technical chapter, and as exercises here, we have:

\begin{exercise}
Learn more about amenability, and related topics.
\end{exercise}

\begin{exercise}
Investigate the free quantum groups, at small values of $N$.
\end{exercise}

\begin{exercise}
Learn about quantum isometries, and related topics.
\end{exercise}

\begin{exercise}
Learn about the various types of quantum quotient spaces.
\end{exercise}

\begin{exercise}
Clarify what we said, regarding the deformations of $O_N^+$.
\end{exercise}

\begin{exercise}
Further explore the free hyperspherical laws.
\end{exercise}

\begin{exercise}
Learn various algebraic twisting techniques.
\end{exercise}

\begin{exercise}
Further explore the free hypergeometric laws.
\end{exercise}

As bonus exercise, and no surprise here, systematically learn quantum groups.

\chapter{Wishart matrices}

\section*{16a. Block transposition}

Time to end this book, and not very clear what to talk about, so many possible choices.  However, looking back at what we did, in the last few chapters, there is something frustrating there, namely the fact that the free probability that we learned, while certainly having the potential of being as wide as classical probability, right, remains something a bit abstract, not really alive, heavily depending on classical probability.

\bigskip

So, can we have some alive theory going on, an independent creature, self-replicating at will, with laws and theorems endlessly producing further laws and theorems, without any help from the outside? That would be very nice, and getting started with this, I think I should first ask the rat, as usual in such philosophical situations:

\begin{rat}
Any strong and healthy rat can produce dozens of baby rats, over the time. The same should happen for the random matrices.
\end{rat}

Well, thanks rat, this sounds like a good idea, manipulating and cutting and gluing the random matrices, until reaching to a Frankenstein. By the way, good luck in finding a lady rat, and there should be certainly place under the stove for all of you guys, unless the cats have their say, I mean we can perhaps talk later, about all this. 

\bigskip

Getting started now, we will be mostly interested in the Wishart matrices, which being ``discrete'', are easier to manipulate. Regarding them, let us recall that we have:

\index{asymptotic freeness}

\begin{theorem}
Given a sequence of complex Wishart matrices
$$W_N=Y_NY_N^*\in M_N(L^\infty(X))$$
with $Y_N$ being $N\times M$ complex Gaussian of parameter $1$, we have
$$\frac{W_N}{N}\sim\max(1-t,0)\delta_0+\frac{\sqrt{4t-(x-1-t)^2}}{2\pi x}\,dx$$
with $M=tN\to\infty$, with the limiting measure being the Marchenko-Pastur law $\pi_t$.
\end{theorem}

\begin{proof}
This is the Marchenko-Pastur theorem \cite{mpa}, that we know well since chapter 8. Let us also mention, at the level of the basics, that, as explained in chapter 13, as a consequence of the results of Voiculescu in \cite{vo3}, when looking at a family of sequences of such matrices, in the $N\to\infty$ limit these become asymptotically free. 
\end{proof}

Regarding now the possible manipulations on the Wishart matrices, there are many of them, gradually explored over the time, and quite systematically in the last 20 years, notably with an  interesting 2003 paper of Graczyk, Letac and Massam \cite{glm}. Of particular interest here is the study of more general products of random matrices, be them of Gaussian, Wigner or Wishart type, and for some recent developments here, in relation with the free Bessel laws, we refer to \cite{bb+} and subsequent papers.

\bigskip

However, it is not about products, which remain something quite technical, that we would like to talk about, in this chapter. As something which is perhaps more fruitful, or at least looks more basic, we have the idea of block-transposing or block-modifying a Wishart matrix, according to the following suprising finding of Aubrun \cite{aub}:

\begin{fact}
When suitably block-transposing the entries of a complex Wishart matrix, we obtain as asymptotic distribution a shifted version of Wigner's semicircle law. 
\end{fact}

So, let us first try to understand this, and for generalizations, we can see later. Following \cite{aub}, and some more recent work in \cite{bn1}, where the Aubrun finding was further studied, and better understood, let us start with the following definition:

\begin{definition}
The partial transpose of a complex Wishart matrix $W$ of parameters $(N,M)=(dn,dm)$ is the matrix
$$\widetilde{W}=(id\otimes t)W$$
where $id$ is the identity of $M_d(\mathbb C)$, and $t$ is the transposition of $M_n(\mathbb C)$. 
\end{definition}

In more familiar terms of bases and indices, the standard decomposition $\mathbb C^{dn}=\mathbb C^d\otimes\mathbb C^n$ induces an algebra decomposition $M_{dn}(\mathbb C)=M_d(\mathbb C)\otimes M_n(\mathbb C)$, and with this convention made, the partial transpose matrix $\widetilde{W}$ constructed above has entries as follows:
$$\widetilde{W}_{ia,jb}=W_{ib,ja}$$

Our goal in what follows will be that of computing the law of $\widetilde{W}$, first when $d,n,m$ are fixed, and then in the $d\to\infty$ regime. For this purpose, we will need a number of standard facts regarding the noncrossing partitions. Let us start with:

\begin{proposition}
For a permutation $\sigma\in S_p$, we have the formula
$$|\sigma|+\#\sigma=p$$
where $|\sigma|$ is the number of cycles of $\sigma$, and $\#\sigma$ is the minimal $k\in\mathbb N$ such that $\sigma$ is a product of $k$ transpositions. Also, the following formula defines a distance on $S_p$, 
$$(\sigma,\pi)\to\#(\sigma^{-1}\pi)$$
and the set of permutations $\sigma\in S_p$ which saturate the triangular inequality 
$$\#\sigma+\#(\sigma^{-1}\gamma)=\#\gamma=p-1$$
where $\gamma\in S_p$ is a full cycle, is in bijection with the set $NC(p)$. 
\end{proposition}

\begin{proof}
All this is standard combinatorics, that can be found for instance, along with many other things, in the paper of Biane \cite{bia}, and that we will leave as an exercise.
\end{proof}

We use the standard bijection $NC(p)\simeq NC_2(2p)$, denoted $\pi\to\widetilde{\pi}$, obtained by fattening the partitions. We have the following formula,  where $\vee$ is the join operation on $NC_2(2p)$, and $\rho_{12} = (12)(34)\ldots(2p-1,2p)$ is the fattened identity permutation:
$$|\pi|=|\widetilde{\pi}\vee\rho_{12}|$$

Similarly, we have the formula $|\pi\gamma|=|\widetilde{\pi}\vee\rho_{14}|$, where $\rho_{14}$ is the pairing corresponding to the fattening of the inverse full cycle $\gamma^{-1}(i) = i-1$, which pairs an element $2i$ with $2(i-1)-1 = 2i-3$, or, equivalently, an element $i \in \{1, \ldots, 2p\}$ with $i+(-1)^{i+1}3$. 

\bigskip

We will need the following well-known result, which is something standard:

\begin{proposition}
The number $||\pi||$ of blocks having even size is given by
$$1+||\pi||=|\pi\gamma|$$
for every noncrossing partition $\pi \in NC(p)$.
\end{proposition}

\begin{proof}
We use a recurrence over the number of blocks of $\pi$. If $\pi$ has just one block, its associated geodesic permutation is $\gamma$ and we have, as desired:
$$|\gamma^2|=1+\delta_{2|p}$$

For the partitions $\pi$ having more than one block, we can assume without loss of generality that $\pi = \hat 1_k \sqcup \pi'$, where $\hat 1_k$ is a contiguous block of size $k$. Recall now that the number of blocks of the permutation $\pi\gamma$ is given by the following formula, where $\rho_{14} \in P_2(2p)$ is the pair partition which pairs an element $i$ with $i+(-1)^{i+1}3$:
$$|\pi\gamma|=|\widetilde{\pi}\vee\rho_{14}|$$

If $k$ is an even number, $k=2r$, in order to prove the result, consider the following partition, which contains the block $(1 \, 4 \, 5 \, 8 \, \ldots 4r-3 \, 4r)$, along with the blocks coming from the elements of the form $4i+2, 4i+3$ from $\{1, \ldots, 4r\}$ and from $\pi'$:
$$\sigma=\widetilde{\hat 1_{2r} \sqcup \pi'}\vee \rho_{14}$$

Now recall that we can count the blocks of the join of two partitions by drawing them one beneath the other, and then counting the number of connected components of the curve, without taking into account the possible crossings. We conclude that we have the following formula, where $\rho'_{14}$ is $\rho_{14}$ restricted to the set $\{2k+1, 2k+2 \ldots, 2p\}$:
$$|\widetilde{\pi}\vee\rho_{14}|=1+|\widetilde{\pi'}\vee\rho'_{14}|$$

As for the case where $k$ is odd, here there is no extra block appearing, and so:
$$|\widetilde{\pi}\vee\rho_{14}|=|\widetilde{\pi'}\vee\rho'_{14}|$$

Summarizing, we are led to the conclusion in the statement.
\end{proof}

We can now investigate the block-transposed Wishart matrices, and we have:

\begin{theorem}
For any $p\geq 1$ we have the formula
$$\lim_{d\to\infty}(E\circ tr)\big(m\widetilde{W}\big)^p
=\sum_{\pi\in NC(p)}m^{|\pi|}n^{||\pi||}$$
where $|.|$ and $||.||$ are the number of blocks, and the number of blocks of even size. 
\end{theorem}

\begin{proof}
The matrix elements of the partial transpose matrix are given by:
$$\widetilde{W}_{ia,jb}=W_{ib,ja}=(dm)^{-1}\sum_{k=1}^d\sum_{c=1}^mG_{ib,kc}\bar{G}_{ja,kc}$$

This gives the following formula, for the corresponding normalized trace:
\begin{eqnarray*}
tr(\widetilde{W}^p)
&=&(dn)^{-1}(dm)^{-p}\sum_{i_1,\ldots,i_p=1}^d\sum_{a_1,\ldots,a_p=1}^n\prod_{s=1}^p W_{i_sa_{s+1},i_{s+1}a_s} \\
&=&(dn)^{-1}(dm)^{-p}\sum_{i_1,\ldots,i_p=1}^d\sum_{a_1,\ldots,a_p=1}^n\prod_{s=1}^p \sum_{j_1,\ldots,j_p=1}^d\sum_{b_1,\ldots,b_p=1}^mG_{i_sa_{s+1},j_sb_s}\bar{G}_{i_{s+1}a_s,j_sb_s}
\end{eqnarray*}

After interchanging the product with the last two sums, the average of the general term can be computed by the Wick rule, namely:
$$E\left(\prod_{s=1}^pG_{i_sa_{s+1},j_sb_s}\bar{G}_{i_{s+1}a_s,j_sb_s}\right)
=\sum_{\pi\in S_p}\prod_{s=1}^p\delta_{i_s,i_{\pi(s)+1}}\delta_{a_{s+1},a_{\pi(s)}}\delta_{j_s,j_{\pi(s)}}\delta_{b_s,b_{\pi(s)}}$$

Let $\gamma\in S_p$ be the full cycle $\gamma=(1 \, 2 \, \ldots \, p)^{-1}$. The general factor in the above product is 1 if and only if the following four conditions are simultaneously satisfied:
$$\gamma^{-1}\pi\leq \ker i\quad,\quad
\pi\gamma \leq \ker a\quad,\quad
\pi \leq \ker j\quad,\quad
\pi \leq \ker b$$

Counting the number of free parameters in the above equation, we obtain:
\begin{eqnarray*}
(E\circ tr)(\widetilde{W}^p)
&=&(dn)^{-1}(dm)^{-p}\sum_{\pi\in S_p}d^{|\pi|+|\gamma^{-1}\pi|}m^{|\pi|}n^{|\pi\gamma|}\\
&=&\sum_{\pi\in S_p}d^{|\pi|+|\gamma^{-1}\pi|-p-1}m^{|\pi|-p}n^{|\pi\gamma|-1}
\end{eqnarray*}

The exponent of $d$ in the last expression on the right is as follows:
\begin{eqnarray*}
N(\pi)
&=&|\pi|+|\gamma^{-1}\pi|-p-1\\
&=&p-1-(\#\pi+\#(\gamma^{-1}\pi))\\
&=&p-1-(\#\pi+\#(\pi^{-1}\gamma))
\end{eqnarray*}

As explained in the beginning of this section, this quantity is known to be $\leq 0$, with equality iff $\pi$ is geodesic, hence associated to a noncrossing partition. Thus:
$$(E\circ tr)(\widetilde{W}^p)=(1+O(d^{-1}))m^{-p}n^{-1}\sum_{\pi\in NC(p)}m^{|\pi|} n^{|\pi\gamma|}$$

But together with $|\pi\gamma|=||\pi||+1$, this gives the result.
\end{proof}

We would like now to find an equation for the moment generating function of the asymptotic law of $m\widetilde{W}$. This moment generating function is defined as:
$$F(z)=\lim_{d\to\infty}(E\circ tr)\left(\frac{1}{1-zm\widetilde{W}}\right)$$

We have the following result, regarding this moment generating function:

\begin{theorem}
The moment generating function of $m\widetilde{W}$ satisfies the equation
$$(F-1)(1-z^2F^2)=mzF(1+nzF)$$
in the $d\to\infty$ limit.
\end{theorem}

\begin{proof}
We use the formula in Theorem 16.7. If we denote by $N(p,b,e)$ the number of partitions in $NC(p)$ having $b$ blocks and $e$ even blocks, we have:
\begin{eqnarray*}
F
&=&1+\sum_{p=1}^\infty\sum_{\pi\in NC(p)} z^pm^{|\pi|}n^{||\pi||}\\
&=&1+\sum_{p=1}^\infty\sum_{b=0}^\infty\sum_{e=0}^\infty z^pm^bn^eN(p,b,e)
\end{eqnarray*}

Let us try to find a recurrence formula for the numbers $N(p,b,e)$. If we look at the block containing $1$, this block must have $r\geq 0$ other legs, and we get:
\begin{eqnarray*}
N(p,b,e) 
&=&\sum_{r\in 2\mathbb N}\sum_{p=\Sigma p_i+r+1}\sum_{b=\Sigma b_i+1}\sum_{e=\Sigma e_i}N(p_1,b_1,e_1)\ldots N(p_{r+1},b_{r+1},e_{r+1})\\
&+&\sum_{r\in 2\mathbb N+1}\sum_{p=\Sigma p_i+r+1}\sum_{b=\Sigma b_i+1}\sum_{e=\Sigma e_i+1}N(p_1,b_1,e_1)\ldots N(p_{r+1},b_{r+1},e_{r+1})
\end{eqnarray*}

Here $p_1,\ldots,p_{r+1}$ are the number of points between the legs of the block containing 1, so that we have $p=(p_1+\ldots+p_{r+1})+r+1$, and the whole sum is split over two cases, $r$ even or odd, because the parity of $r$ affects the number of even blocks of our partition. Now by multiplying everything by a $z^pm^bn^e$ factor, and by carefully distributing the various powers of $z,m,b$ on the right, we obtain the following formula:
\begin{eqnarray*}
z^pm^bn^eN(p,b,e)
&=&m\sum_{r\in 2\mathbb N}z^{r+1}\sum_{p=\Sigma p_i+r+1}\sum_{b=\Sigma b_i+1}\sum_{e=\Sigma e_i}\prod_{i=1}^{r+1}z^{p_i}m^{b_i}n^{e_i}N(p_i,b_i,e_i)\\
&+&mn\sum_{r\in 2\mathbb N+1}z^{r+1}\sum_{p=\Sigma p_i+r+1}\sum_{b=\Sigma b_i+1}\sum_{e=\Sigma e_i+1}\prod_{i=1}^{r+1}z^{p_i}m^{b_i}n^{e_i}N(p_i,b_i,e_i)
\end{eqnarray*}

Let us sum now all these equalities, over all $p\geq 1$ and over all $b,e\geq 0$. According to the definition of $F$, at left we obtain $F-1$. As for the two sums appearing on the right, that is, at right of the two $z^{r+1}$ factors, when summing them over all $p\geq 1$ and over all $b,e\geq 0$, we obtain in both cases $F^{r+1}$. So, we have the following formula:
\begin{eqnarray*}
F-1
&=&m\sum_{r\in 2\mathbb N}(zF)^{r+1}+mn\sum_{r\in 2\mathbb N+1}(zF)^{r+1}\\
&=&m\,\frac{zF}{1-z^2F^2}+mn\,\frac{z^2F^2}{1-z^2F^2}\\
&=&mzF\,\frac{1+nzF}{1-z^2F^2}
\end{eqnarray*}

But this gives the formula in the statement, and we are done.
\end{proof}

Our goal now will be that of further processing the formula in Theorem 16.8, as to reach to a formula for the density of the corresponding law. We first have:

\begin{theorem}
The Cauchy transform of $m\widetilde{W}$ satisfies the equation
$$(\xi G-1)(1-G^2)=mG(1+nG)$$
in the $d\to\infty$ limit. Moreover, this equation simply reads
$$R=\frac{m}{2}\left(\frac{n+1}{1-z}-\frac{n-1}{1+z}\right)$$
with the usual $R$-transform substitutions $G\to z$ and $\xi\to R+z^{-1}$.
\end{theorem}

\begin{proof}
We have two assertions to be proved, the first one being standard, and the second one being something quite magic, the idea being as follows:

\medskip

(1) Consider the equation of $F$, found in Theorem 16.8, namely:
$$(F-1)(1-z^2F^2)=mzF(1+nzF)$$

With $z\to\xi^{-1}$ and $F\to\xi G$, so that $zF\to G$, we obtain, as desired:
$$(\xi G-1)(1-G^2)=mG(1+nG)$$

(2) Next, let us look at the equation of the Cauchy transform that we have. With the substitutions $\xi \to K$ and $G\to z$, this equation becomes:
$$(zK-1)(1-z^2)=mz(1+nz)$$

The point now is that with $K\to R+z^{-1}$ this latter equation becomes:
$$zR(1-z^2)=mz(1+nz)$$

But the solution of this latter equation is trivial to compute, given by:
$$R=m\,\frac{1+nz}{1-z^2}=\frac{m}{2}\left(\frac{n+1}{1-z}-\frac{n-1}{1+z}\right)$$

Thus, we are led to the conclusion in the statement.
\end{proof}

We can now answer our original question, regarding the limiting law, as follows:

\begin{theorem}
With $d\to\infty$ we have the following convergence in law,
$$m\widetilde{W}\sim\mu_{mn}$$
where $\mu_{mn}$ is the free difference of free Poisson laws of parameters $m(n\pm 1)/2$.
\end{theorem}

\begin{proof}
Consider two free variables $a,b$, following free Poisson laws of parameters $s,t$. According to our formulae from chapter 13, the corresponding $R$-transforms are:
$$R_a=\frac{s}{1-z}\quad,\quad R_b=\frac{t}{1-z}$$

By using the general dilation formula $R_{qb}(z)=qR_b(qz)$  at $q=-1$, we deduce that we have $R_{-b}=-t/(1+z)$. Thus, our two formulae above can be written as follows:
$$R_a=\frac{s}{1-z}\quad,\quad R_{-b}=-\frac{t}{1+z}$$

Now when making the sum, the corresponding $R$-transform will be given by:
$$R_{a-b}=\frac{s}{1-z}-\frac{t}{1+z}$$

But with $s=m(n+1)/2$ and $t=m(n-1)/2$, this latter formula reads:
$$R=\frac{m}{2}\left(\frac{n+1}{1-z}-\frac{n-1}{1+z}\right)$$

And since this is the formula found in Theorem 16.9, we obtain the result.
\end{proof}

\section*{16b. Shifted semicircles}

Still following \cite{bn1}, let us have now a closer look at the asymptotic laws that we found in Theorem 16.10, among others with the goal of recovering, as a particular case of that result, the original finding by Aubrun, regarding shifted semicircles \cite{aub}.

\bigskip

It is convenient to enlarge our study to all possible real values of $m,n$, as follows:

\begin{definition}
The real measure $\mu_{mn}$ of parameters $m>0,n\geq 1$ is given by
$$a\sim\pi_{m(n+1)/2}\ ,\ b\sim\pi_{m(n-1)/2}\ ,\ {\rm free}\quad\implies\quad a-b\sim\pi_{mn}$$
that is, is the free difference of free Poisson laws of parameters $m(n\pm 1)/2$.
\end{definition}

As a first remark, the various formulae that we found before, namely moment formula, and various equations for the functional transforms, are valid in this more general setting. Let us collect these formulae, along with a bit more, in a single statement:

\begin{theorem}
The measure $\mu_{mn}$ is a compound free Poisson law,
$$\mu_{mn}=\pi_\nu\quad,\quad\nu=\frac{m(n-1)}{2}\,\delta_{-1}+\frac{m(n+1)}{2}\,\delta_1$$
and has the following properties:
\begin{enumerate}
\item The moments are $M_k=\sum_{\pi\in NC(k)}m^{|\pi|}n^{||\pi||}$.

\item The mean is $E=m$, the variance is $V=mn$.

\item $\gamma=1/(n\sqrt{mn})$ and $\kappa=2+1/(mn)$.

\item $(M-1)(1-z^2M^2)=mzM(1+nzM)$.

\item $(\xi G-1)(1-G^2)=mG(1+nG)$.

\item $R(z)=m(1+nz)/(1-z^2)$.
\end{enumerate}
\end{theorem}

\begin{proof}
Regarding the first assertion, let us recall from chapter 13 that the $R$-transform of a compound free Poisson law is given by the following formula:
$$\nu=\sum_ic_i\delta_{y_i}\quad\implies\quad R_{\pi_{\nu}}(z)=\sum_i\frac{c_iy_i}{1-y_iz}$$

For the measure $\nu$ in the statement, we obtain the following $R$-transform:
$$R(z)=\frac{m(n-1)}{2}\cdot\frac{-1}{1+z}+\frac{m(n+1)}{2}\cdot\frac{1}{1-z}
=\frac{m}{2}\left(\frac{n+1}{1-z}-\frac{n-1}{1+z}\right)$$

But this is the $R$-transform of $\mu_{mn}$ that we computed before. As for the rest, these are basically things that we know from before, coming now in reverse order:

\medskip

(1) We have the formula of the $R$-transform, and by using the manipulations from the proof of Theorem 16.9, in reverse, we obtain the equation for the Cauchy transform $G$. 

\medskip

(2) Next, the same procedure leads to the equation for moment generating function, that we can call now $M$ as usual, since done with matrices and the $M,N$ there. 

\medskip

(3) Next, the same procedure leads to the formula in the statement for the moments, that we can call $M_k$ as usual, now that done with indices and combinatorics.

\medskip

(4) It remains to compute $E,V,\gamma,\kappa$. The low order moments are as follows:
$$M_1=m$$
$$M_2=m^2+mn$$
$$M_3=m^3+3m^2n+m$$
$$M_4=m^4+6m^3n+4m^2+2m^2n^2+mn$$

Thus the mean is $E=m$, and the variance is $V=mn$. Next, we have:
$$M_3'=M_3-3EM_2+2E^3=m$$
$$M_4'=M_4-4EM_3+6E^2M_2-3E^4=mn(2mn+1)$$ 

Thus the skewness and kurtosis are given by the following formulae:
$$\gamma=\frac{m}{mn\sqrt{mn}}
=\frac{1}{n\sqrt{mn}}\quad,\quad  
\kappa=\frac{mn(2mn+1)}{m^2n^2}
=2+\frac{1}{mn}$$

Finally, observe that we have $\kappa\simeq 2$ as we should, our law being of free nature. 
\end{proof}

Regarding now the main particular cases of measures of type $\mu_{mn}$, corresponding to some previously known random matrix computations, the situation is as follows:

\begin{theorem}
The measures $\mu_{mn}$ are as follows:
\begin{enumerate}
\item At $n=1$ we have the Marchenko-Pastur laws, $\mu_{m1}=\pi_m$.

\item With $n=sm\to\infty$ we obtain Aubrun's shifted semicircles.

\item With $m=t/n\to 0$ we obtain the free Bessel law $\pi^2_t$.
\end{enumerate}
\end{theorem}

\begin{proof}
These results come either from Definition 16.11, or from Theorem 16.12:

\medskip

(1) This follows from Definition 16.11, because at $n=1$ we have $m(n\pm 1)/2=m,0$. Observe that this is in tune with the Theorem 16.2, because at $n=1$ the partial transposition is trivial, so Theorem 16.10 computes the asymptotic law of $mW$.

\medskip

(2) In order to explain what is going on here, recall from Theorem 16.10 that, with the notations there, we have $m\widetilde{W}\sim\mu_{mn}$ in the $d\to\infty$ limit. Thus, with the convention that $D$ is the dilation operation, given by $a\sim\mu\implies qa\sim D_q(\mu)$, we have:
$$\widetilde{W}\sim D_{m^{-1}}(\mu_{mn})$$

Now the point is that, what Aubrun discovered in \cite{aub} is the following formula, valid in the regime in the statement, $n=sm\to\infty$, and with $d\to\infty$ too, as usual:
$$\widetilde{W}\sim\gamma_s^1$$

So, let us prove this. In view of what we have, this amounts in proving that:
$$\lim_{m\to\infty}D_{m^{-1}}(\mu_{m,sm})=\gamma_s^1$$

For this purpose, let $\Gamma=\Gamma(\xi)$ be the Cauchy transform of $\nu=D_{m^{-1}}(\mu_{mn})$. The Cauchy transform of $\mu_{mn}=D_m(\nu)$ is then $G(m\xi)=m^{-1}\Gamma$, and by making the replacements $\xi\to m\xi$ and $G\to m^{-1}\Gamma$ in the equation of $G$ in Theorem 16.12, we obtain: 
$$(\xi\Gamma-1)(1-m^{-2}\Gamma^2)=\Gamma(1+nm^{-1}\Gamma)$$

In the limit $n=sm\to\infty$, as indicated above, this equation becomes:
\begin{eqnarray*}
\xi\Gamma-1=\Gamma(1+s\Gamma)
&\implies&s\Gamma^2+(1-\xi)\Gamma+1=0\\
&\implies&\Gamma=\frac{\xi-1\pm\sqrt{(1-\xi)^2-4s}}{2s}
\end{eqnarray*}

By applying now the Stieltjes inversion formula, we obtain the following density:
$$\varphi(x)=\frac{\sqrt{4s-(x-1)^2}}{2s\pi}$$

But this is exactly the density of the semicircle law $\gamma_s^1$, and we are done.

\medskip

(3) This comes from Theorem 16.12, because with $m=t/n\to 0$, the measure $\nu$ appearing there, producing the compound free Poisson law $\mu_{mn}=\pi_\nu$, is given by: 
\begin{eqnarray*}
\nu
&=&\frac{m(n-1)}{2}\,\delta_{-1}+\frac{m(n+1)}{2}\,\delta_1\\
&=&\frac{t-m}{2}\,\delta_{-1}+\frac{t+m}{2}\,\delta_1\\
&\simeq&t\cdot\frac{\delta_{-1}+\delta_1}{2}
\end{eqnarray*}

Thus, we obtain indeed the free Bessel law $\pi^2_t$, as constructed in chapter 14.
\end{proof}

As a comment now on (3), we were talking before the statement of the theorem about ``previous random matrix computations'', and in the case of the free Bessel law $\pi^2_t$, this refers to the computations in \cite{bb+} using products of Gaussian matrices, producing matrix models for $\pi^2_1$, or rather for its modified version $\tilde{\pi}^2_1$, in the sense of chapter 14.

\bigskip

Now the point is that these latter computations have nothing to do with our current block-transposition business for Wishart matrices. Thus, in short, what we have here are different approaches to the modeling on $\pi^2_1$. And in relation with this, an interesting question appears. Indeed, with the methods in \cite{bb+} working for any $\pi^s_1$, we have:

\begin{question}
Can our block-transposition business for Wishart matrices be suitably modified, as to cover the general free Bessel laws $\pi^s_t$?
\end{question}

And good question this is, very natural because, as explained earlier in this book, the free Bessel laws $\pi^s_t$ are at the core of free probability of discrete type, so any candidate for a serious subtheory in free probability should have these laws $\pi^s_t$ well understood.

\bigskip

Well, so this is at least how the general philosophy goes. In practice now, Question 16.14 has led to a lot of work, first with my second paper with Nechita, \cite{bn2}, laying the foundations for a possible answer, and we will talk about this in a moment, and then with two follow-ups to that \cite{bn2} paper, one by myself not solving the problem, but working on that was fun, and one by Nechita with Arizmendi and Vargas \cite{anv}, basically solving the problem, but in a quite complicated way. Somehow, Question 16.14 still stands.

\bigskip

So, this is the situation and.. but wait, someone is meowing. That is the young grey cat, he's been around the house for a few weeks already, and boy, is this young fellow quick and smart, compared to the other fossils. Let's see what he has to say:

\begin{cat}
Come on man, you can do everything with random matrices, even talk about the zeroes of zeta. 
\end{cat}

Okay, thanks cat, shall I understand that there is hope. So, we will add this to our to-do list in free probability and random matrices, find simple models for the laws $\pi^s_t$, along of course with improving the physics that we have, quantum and statistical mechanics, which is  already not that bad, and why not, some thinking about zeta too.

\bigskip

Getting back to work now, still following \cite{bn1}, we have the following key supplementary result, adding to Theorem 16.12, regarding the densities of the measures $\mu_{mn}$:

\begin{theorem}
The density of $\mu_{mn}$ has the following properties:
\begin{enumerate}
\item It has at most one atom, located at $0$, of mass $\max(1-mn,0)$.

\item The support is positive precisely when $n\leq m/4+1/m$ and $m \geq 2$.
\end{enumerate}
\end{theorem}

\begin{proof}
We have two assertions here, the idea being as follows:

\medskip

(1) The first assertion follows from the general characterization of atoms of a free additive convolution, given by  Bercovici and Voiculescu in \cite{bvo}, namely $x$ is an atom for the free additive convolution of two measures $\nu_1$ and $\nu_2$ precisely when:
$$x=x_1+x_2\quad,\quad\nu_1(\{x_1\})+\nu_2(\{x_2\})>1$$

In addition, if this is the case, we have the following formula, also from \cite{bvo}:
$$[\nu_1\boxplus\nu_2](\{x\})=\nu_1(\{x_1\})+\nu_2(\{x_2\})-1$$

In our situation, we can apply this with $\nu_1$ being a free Poisson law of parameter $m(n+1)/2$, and $\nu_2$ being the image of a free Poisson law of parameter $m(n-1)/2$ through the negation map, and we are led to the formula in the statement.

\medskip

(2) This is something more complicated, based on a detailed study of the degree 3 equation  for the Cauchy transform, from Theorem 16.12, namely:
$$(\xi G-1)(1-G^2)=mG(1+nG)$$

We use the standard fact that the support of the absolutely continuous part of a probability measure is a union of intervals defined by the points where the analyticity of the Cauchy transform $G$ breaks. In our case, for the measure $\mu_{mn}$, these points are the roots of the discriminant of the above equation defining $G$, which is as follows:
\begin{eqnarray*}
\Delta(\xi)
&=&4\xi^4-12m\xi^3+(n^2m^2+12m^2-20nm-8)\xi^2\\
&&\ +(-2n^2m^3-4m^3+22nm^2-20m)\xi\\
&&\ +n^2m^4-4n^3m^3-2nm^3+12n^2m^2+m^2-12nm+4
\end{eqnarray*}

Now observe that $\Delta(\xi)=0$ is a degree 4 equation in $\xi$ that has either 2 or 4 real solutions, with 4 complex solutions being ruled out, since the support must be non-empty. In order to decide between these cases, we compute the discriminant of $\Delta$:
$$\Delta_2(m,n)=-256m^2(n-1)(n+1)\left(m^3n^3+15m^2n^2-27m^2+48mn-64\right)^3$$

Indeed, recall that the sign of the discriminant of a quartic equation permits to decide whether the equation has two real and two complex roots or 4 roots of the same type, real or complex, the precise result being that the discriminant is negative precisely when the quartic has two real and two complex roots. In our case, due to our assumption $n\geq1$, and with the case $n=1$ being trivial, already discussed in Theorem 16.13 (1), the sign of $\Delta_2$ is the opposite of the sign of the following two-variable polynomial: 
$$P(m,n) = m^3 n^3+15 m^2 n^2+48 m n-27 m^2-64$$

But this leads, via some calculus, partly computer aided, to the conclusion in the statement. For details and pictures here, taking about 3 pages, we refer to \cite{bn1}.
\end{proof}

Finally, as a comment on all this, the above computations, and particularly the conditions $n\leq m/4+1/m$ and $m\geq 2$ found in \cite{bn1} and discussed above for the positivity of the support, are something quite interesting in quantum information theory, where the block transposition of the Wishart matrices can be a useful tool, and where positivity is a key asset. For more on all this, we refer to \cite{aub}, \cite{bn1}, and their various follow-ups.

\section*{16c. Block modifications}

As a continuation of the above, let us discuss now more general block modifications, following the more recent paper \cite{bn2}. We have the following construction:

\begin{definition}
Given a complex Wishart matrix $W=YY^*\in M_{dn}(L^\infty(X))$ as befoe, with $Y$ being a complex Gaussian $dn\times dm$ matrix, and a linear map
$$\varphi:M_n(\mathbb C)\to M_n(\mathbb C)$$
we consider the following matrix, obtained by applying $\varphi$ to the $n\times n$ blocks of $W$,
$$\widetilde{W}=(id\otimes\varphi)W\in M_{dn}(L^\infty(X))$$
and call it block-modified Wishart matrix.
\end{definition}

Here we are using some standard tensor product identifications, the details being as follows. Let $Y$ be a complex Gaussian $dn\times dm$ matrix, as above:
$$Y\in M_{dn\times dm}(L^\infty(X))$$

We can then form the corresponding complex Wishart matrix, as follows: 
$$W=YY^*\in M_{dn}(L^\infty(X))$$

The size of this matrix being a composite number, $N=dn$, we can regard this matrix as being a $n\times n$ matrix, with random $d\times d$ matrices as entries. Equivalently, by using standard tensor product notations, this amounts in regarding $W$ as follows:
$$W\in M_d(L^\infty(X))\otimes M_n(\mathbb C)$$

With this done, we can come up with our linear map, namely:
$$\varphi:M_n(\mathbb C)\to M_n(\mathbb C)$$

We can apply $\varphi$ to the tensors on the right, and we obtain a matrix as follows:
$$\widetilde{W}=(id\otimes\varphi)W\in M_d(L^\infty(X))\otimes M_n(\mathbb C)$$

Finally, we can forget now about tensors, and as a conclusion to all this, we have constructed a matrix as follows, that we can call block-modified Wishart matrix:
$$\widetilde{W}\in M_{dn}(L^\infty(X))$$

In practice now, what we mostly need for fully understanding Definition 16.17 are examples. Following Aubrun \cite{aub}, then \cite{bn1}, and the paper by Collins and Nechita \cite{cne} too, we have the following basic examples, for our general construction:

\begin{proposition}
We have the following examples of block-modified Wishart matrices $\widetilde{W}=(id\otimes\varphi)W$, coming from various linear maps $\varphi:M_n(\mathbb C)\to M_n(\mathbb C)$:
\begin{enumerate}
\item Wishart matrices: $\widetilde{W}=W$, obtained via $\varphi=id$.

\item Aubrun matrices: $\widetilde{W}=(id\otimes t)W$, with $t$ being the transposition.

\item Collins-Nechita one: $\widetilde{W}=(id\otimes\varphi)W$, with $\varphi=tr(.)1$.

\item Collins-Nechita two: $\widetilde{W}=(id\otimes\varphi)W$, with $\varphi$ erasing the off-diagonal part.
\end{enumerate}
\end{proposition}

\begin{proof}
This is indeed something self-explanatory, standing more as a definition.
\end{proof}

Getting back now to the general case, that of Definition 16.17, the linear map there $\varphi:M_n(\mathbb C)\to M_n(\mathbb C)$ is certainly useful for understanding the construction of the block-modified Wishart matrix $\widetilde{W}=(id\otimes\varphi)W$, as illustrated by the above examples. 

\bigskip

In practice, however, we would like to have as block-modification ``data'' something more concrete, such as a usual matrix. To be more precise, we would like to use:

\begin{proposition}
We have a correspondence between linear maps 
$$\varphi:M_n(\mathbb C)\to M_n(\mathbb C)$$
and square matrices $\Lambda\in M_n(\mathbb C)\otimes M_n(\mathbb C)$, given by the formula
$$\Lambda_{ab,cd}=\varphi(e_{ac})_{bd}$$
where $e_{ab}\in M_n(\mathbb C)$ are the standard generators of the matrix algebra $M_n(\mathbb C)$, given by the formula $e_{ab}:e_b\to e_a$, with $\{e_1,\ldots,e_n\}$ being the standard basis of $\mathbb C^n$.
\end{proposition}

\begin{proof}
This is standard linear algebra. Given a linear map $\varphi:M_n(\mathbb C)\to M_n(\mathbb C)$, we can associated to it numbers $\Lambda_{ab,cd}\in\mathbb C$ by the formula in the statement, namely:
$$\Lambda_{ab,cd}=\varphi(e_{ac})_{bd}$$

Now by using these $n^4$ numbers, we can construct a $n^2\times n^2$ matrix, as follows:
$$\Lambda=\sum_{abcd}\Lambda_{ab,cd}e_{ac}\otimes e_{bd}\in M_n(\mathbb C)\otimes M_n(\mathbb C)$$

Thus, we have constructed a correspondence $\varphi\to\Lambda$, and since this correspondence is injective, and the dimensions match, this correspondence is bijective, as claimed.
\end{proof}

Now by getting back to the block-modified Wishart matrices, we have:

\begin{proposition}
Given a Wishart $dn\times dn$ matrix $W=YY^*$, and a linear map
$$\varphi:M_n(\mathbb C)\to M_n(\mathbb C)$$
the entries of the corresponding block-modified matrix $\widetilde{W}=(id\otimes\varphi)W$ are given by
$$\widetilde{W}_{ia,jb}=\sum_{cd}\Lambda_{ca,db}W_{ic,jd}$$
where $\Lambda\in M_n(\mathbb C)\otimes M_n(\mathbb C)$ is the square matrix associated to $\varphi$, as above.
\end{proposition}

\begin{proof}
Again, this is linear algebra, coming from the following computation:
$$\widetilde{W}_{ia,jb}
=\sum_{cd}W_{ic,jd}\varphi(e_{cd})_{ab}
=\sum_{cd}\Lambda_{ca,db}W_{ic,jd}$$

We are therefore led to the conclusion in the statement.
\end{proof}

At the level of the main examples, from Proposition 16.18, the very basic linear maps $\varphi:M_n(\mathbb C)\to M_n(\mathbb C)$ used there can only correspond to some basic examples of matrices $\Lambda\in M_n(\mathbb C)\otimes M_n(\mathbb C)$, via the correspondence in Proposition 16.19. This is indeed the case, and in order to clarify this, and at a rather conceptual level, let us formulate, inspired by the representation theory material from chapter 7, the following definition:

\begin{definition}
Let $P(k,l)$ be the set of partitions between an upper row of $k$ points, and a lower row of $l$ points. Associated to any $\pi\in P(k,l)$ is the linear map
$$T_\pi(e_{i_1}\otimes\ldots\otimes e_{i_k})=\sum_{j_1\ldots j_l}\delta_\pi\begin{pmatrix}i_1&\ldots&i_k\\ j_1&\ldots&j_l\end{pmatrix}e_{j_1}\otimes\ldots\otimes e_{j_l}$$
between tensor powers of $\mathbb C^N$, called ``easy'', with the Kronecker type symbol on the right being given by $\delta_\pi=1$ when the indices fit, and $\delta_\pi=0$ otherwise.
\end{definition}

Observe the obvious connection with notion of easy group, from chapter 7, the point being that a closed subgroup $G\subset U_N$ is easy precisely when its Tannakian category $C_G=(C_G(k,l))$ with $C_G(k,l)\subset\mathcal L((\mathbb C^N)^k,(\mathbb C^N)^l)$ is spanned by easy maps. 

\bigskip

For our purposes, we will need a slight modification of Definition 16.21, as follows:

\begin{definition}
Associated to any partition $\pi\in P(2s,2s)$ is the linear map
$$\varphi_\pi(e_{a_1\ldots a_s,c_1\ldots c_s})=\sum_{b_1\ldots b_s}\sum_{d_1\ldots d_s}\delta_\pi\begin{pmatrix}a_1&\ldots&a_s&c_1&\ldots&c_s\\ b_1&\ldots&b_s&d_1&\ldots&d_s\end{pmatrix}e_{b_1\ldots b_s,d_1\ldots d_s}$$
obtained from $T_\pi$ by contracting all the tensors, via the operation
$$e_{i_1}\otimes\ldots\otimes e_{i_{2s}}\to e_{i_1\ldots i_s,i_{s+1}\ldots i_{2s}}$$
with $\{e_1,\ldots,e_N\}$ standing as usual for the standard basis of $\mathbb C^N$.
\end{definition}

In relation with our Wishart matrix considerations, the point is that the above linear map $\varphi_\pi$ can be viewed as a ``block-modification'' map, as follows:
$$\varphi_\pi:M_{N^s}(\mathbb C)\to M_{N^s}(\mathbb C)$$

As an illustration, let us discuss the case $s=1$. There are 15 partitions $\pi\in P(2,2)$, and among them, the most ``basic'' are the 4 partitions $\pi\in P_{even}(2,2)$. We have:

\begin{theorem}
The partitions $\pi\in P_{even}(2,2)$ are as follows,
$$\pi_1=\begin{bmatrix}\circ&\bullet\\ \circ&\bullet\end{bmatrix}\quad,\quad
\pi_2=\begin{bmatrix}\circ&\bullet\\ \bullet&\circ\end{bmatrix}\quad,\quad
\pi_3=\begin{bmatrix}\circ&\circ\\ \bullet&\bullet\end{bmatrix}\quad,\quad
\pi_4=\begin{bmatrix}\circ&\circ\\ \circ&\circ\end{bmatrix}$$
with the associated linear maps $\varphi_\pi:M_n(\mathbb C)\to M_n(\mathbb C)$ being as follows,
$$\varphi_1(A)=A\quad,\quad
\varphi_2(A)=A^t\quad,\quad
\varphi_3(A)=Tr(A)1\quad,\quad
\varphi_4(A)=A^\delta$$
and the associated square matrices $\Lambda_\pi\in M_n(\mathbb C)\otimes M_n(\mathbb C)$ being as follows,
$$\Lambda^1_{ab,cd}=\delta_{ab}\delta_{cd}\quad,\quad 
\Lambda^2_{ab,cd}=\delta_{ad}\delta_{bc}\quad,\quad
\Lambda^3_{ab,cd}=\delta_{ac}\delta_{bd}\quad,\quad
\Lambda^4_{ab,cd}=\delta_{abcd}$$
producing the main examples of block-modified matrices, from Proposition 16.18.
\end{theorem}

\begin{proof}
This is something elementary, coming from the formula in Definition 16.21. Indeed, in the case $s=1$, that we are interested in here, that formula becomes:
$$\varphi_\pi(e_{ac})=\sum_{bd}\delta_\pi\begin{pmatrix}a&c\\ b&d\end{pmatrix}e_{bd}$$

Now in the case of the 4 partitions in the statement, such maps are given by:
$$\varphi_1(e_{ac})=e_{ac}\quad,\quad 
\varphi_2(e_{ac})=e_{ca}\quad,\quad 
\varphi_3(e_{ac})=\delta_{ac}\sum_be_{bb}\quad,\quad
\varphi_4(e_{ac})=\delta_{ac}e_{aa}$$

Thus, we obtain the formulae in the statement. Regarding now the associated square matrices, appearing via $\Lambda_{ab,cd}=\varphi(e_{ac})_{bd}$, these are given by:
$$\Lambda^1_{ab,cd}=\delta_{ab}\delta_{cd}\quad,\quad 
\Lambda^2_{ab,cd}=\delta_{ad}\delta_{bc}\quad,\quad
\Lambda^3_{ab,cd}=\delta_{ac}\delta_{bd}\quad,\quad
\Lambda^4_{ab,cd}=\delta_{abcd}$$

Thus, we are led to the conclusions in the statement.
\end{proof}

Moving ahead now, we would like to study the distribution of the arbitrary block-modified Wishart matrices $\widetilde{W}=(id\otimes\varphi)W$. We will use as before the moment method. However, things will be more tricky in the present setting, and we will need:

\begin{definition}
The generalized colored moments of a random matrix 
$$W\in M_N(L^\infty(X))$$
with respect to a colored integer $e=e_1\ldots e_p$, and a permutation $\sigma\in S_p$, are the numbers
$$M^\sigma_e(W)=\frac{1}{N^{|\sigma|}}\,E\left(\sum_{i_1,\ldots,i_p}W^{e_1}_{i_1i_{\sigma(1)}}\ldots W^{e_p}_{i_pi_{\sigma(p)}}\right)$$
where $|\sigma|$ is the number of cycles of $\sigma$.
\end{definition}

This is something quite technical, in the spirit of the free probability and free cumulant theory in \cite{nsp}, that we will need in what follows. In order to understand how these generalized moments work, consider the standard cycle in $S_p$, namely:
$$\gamma=(1\to2\to\ldots\to p\to 1)$$

If we use this cycle $\gamma\in S_p$ as our permutation $\sigma\in S_p$ in the above definition, the corresponding generalized moment of a random matrix $W$ is then the usual moment:
$$M^\gamma_e(W)
=\frac{1}{N}\,E\left(\sum_{i_1,\ldots,i_p}W^{e_1}_{i_1i_2}\ldots W^{e_p}_{i_pi_1}\right)
=(E\circ tr)(W^{e_1}\ldots W^{e_p})$$

In general, we can decompose the computation of $M^\sigma_e(W)$ over the cycles of $\sigma$, and we obtain in this way a certain product of moments of $W$. See \cite{nsp}.

\bigskip

As a second illustration now, in relation with the usual square matrices, and more specifically with the square matrices $\Lambda\in M_n(\mathbb C)\otimes M_n(\mathbb C)$ as in Proposition 16.19, we have the following formula, that we will use many times in what follows:

\begin{proposition}
Given a usual square matrix, of composed size,
$$\Lambda\in M_n(\mathbb C)\otimes M_n(\mathbb C)$$
we have the following generalized moment formula,
$$(M^\sigma_e\otimes M^\tau_e)(\Lambda)=\frac{1}{n^{|\sigma|+|\tau|}}\sum_{i_1,\ldots, i_p}\sum_{j_1,\ldots,j_p}\Lambda_{i_1j_1,i_{\sigma(1)}j_{\tau(1)}}^{e_1}\ldots\ldots\Lambda_{i_pj_p,i_{\sigma(p)}j_{\tau(p)}}^{e_p}$$
valid for any two permutations $\sigma,\tau\in S_p$, and any colored integer $e=e_1\ldots e_p$.
\end{proposition}

\begin{proof}
This is something obvious, applying the construction in Definition 16.24 with $N=n^2$, $X=\{.\}$, $W=\Lambda$, and then making a tensor product of the corresponding moments $M^\sigma_e$, $M^\tau_e$, regarded as linear functionals on $M_n(\mathbb C)\otimes M_n(\mathbb C)$.
\end{proof}

Consider now the embedding $NC(p)\subset S_p$ obtained by ``cycling inside each block''. That is, each block $b=\{b_1,\ldots,b_k\}$ with $b_1<\ldots<b_k$ of a given noncrossing partition $\sigma\in NC(p)$ produces by definition the cycle $(b_1\ldots b_k)$ of the corresponding permutation $\sigma\in S_p$. Observe that the one-block partition $\gamma\in NC(p)$ corresponds in this way to the standard cycle $\gamma\in S_p$. Also, the number of blocks $|\sigma|$ of a partition $\sigma\in NC(p)$ corresponds to the number of cycles $|\sigma|$ of the corresponding permutation $\sigma\in S_p$. 

\bigskip

With these conventions, we have the following result, from \cite{bn2}, generalizing our various Wishart matrix moment computations, that we did so far in this book:

\index{block-modified matrix}

\begin{theorem}
The asymptotic moments of a block-modified Wishart matrix 
$$\widetilde{W}=(id\otimes\varphi)W$$
with parameters $d,m,n\in\mathbb N$ as before, are given by the formula
$$\lim_{d\to\infty}M_e\left(\frac{\widetilde{W}}{d}\right)=\sum_{\sigma\in NC(p)}(mn)^{|\sigma|}(M^\sigma_e\otimes M^\gamma_e)(\Lambda)$$
where $\Lambda\in M_n(\mathbb C)\otimes M_n(\mathbb C)$ is the square matrix associated to $\varphi:M_n(\mathbb C)\to M_n(\mathbb C)$.
\end{theorem}

\begin{proof}
We use the formula for the matrix entries of $\widetilde{W}$, directly in terms of the matrix $\Lambda$ associated to the map $\varphi$, from Proposition 16.20, namely:
$$\widetilde{W}_{ia,jb}=\sum_{cd}\Lambda_{ca,db}W_{ic,jd}$$

By conjugating this formula, we obtain the following formula for the entries of the adjoint matrix $\widetilde{W}^*$, that we will need as well, in what follows:
$$\widetilde{W}_{ia,jb}^*
=\sum_{cd}\bar{\Lambda}_{db,ca}\bar{W}_{jd,ic}
=\sum_{cd}\Lambda^*_{ca,db}W_{ic,jd}$$

Thus, we have the following global formula, valid for any exponent $e\in\{1,*\}$:
$$\widetilde{W}_{ia,jb}^e=\sum_{cd}\Lambda^e_{ca,db}W_{ic,jd}$$

In order to compute the moments of $\widetilde{W}$, observe first that we have:
\begin{eqnarray*}
tr(\widetilde{W}^{e_1}\ldots\widetilde{W}^{e_p})
&=&\frac{1}{dn}\sum_{i_ra_r}\prod_s\widetilde{W}_{i_sa_s,i_{s+1}a_{s+1}}^{e_s}\\
&=&\frac{1}{dn}\sum_{i_ra_rc_rd_r}\prod_s\Lambda_{c_sa_s,d_sa_{s+1}}^{e_s}W_{i_sc_s,i_{s+1}d_s}\\
&=&\frac{1}{dn}\sum_{i_ra_rc_rd_rj_rb_r}\prod_s\Lambda_{c_sa_s,d_sa_{s+1}}^{e_s}Y_{i_sc_s,j_sb_s}\bar{Y}_{i_{s+1}d_s,j_sb_s}
\end{eqnarray*}

The average of the general term can be computed by the Wick rule, which gives:
$$E\left(\prod_sY_{i_sc_s,j_sb_s}\bar{Y}_{i_{s+1}d_s,j_sb_s}\right)
=\#\left\{\sigma\in S_p\Big|i_{\sigma(s)}=i_{s+1},c_{\sigma(s)}=d_s,j_{\sigma(s)}=j_s,b_{\sigma(s)}=b_s\right\}$$

Let us look now at the above sum. The $i,j,b$ indices range over sets having respectively $d,d,m$ elements, and they have to be constant under the action of $\sigma\gamma^{-1},\sigma,\sigma$. Thus when summing over these $i,j,b$ indices we simply obtain a factor as follows:
$$f=d^{|\sigma\gamma^{-1}|}d^{|\sigma|}m^{|\sigma|}$$

Thus, we obtain the following moment formula:
$$(E\circ tr)(\widetilde{W}^{e_1}\ldots\widetilde{W}^{e_p})
=\frac{1}{dn}\sum_{\sigma\in S_p}d^{|\sigma\gamma^{-1}|}(dm)^{|\sigma|}\sum_{a_rc_r}\prod_s\Lambda_{c_sa_s,c_{\sigma(s)}a_{s+1}}^{e_s}$$

On the other hand, we know from Proposition 16.25 that the generalized moments of the matrix $\Lambda\in M_n(\mathbb C)\otimes M_n(\mathbb C)$ are given by the following formula:
$$(M^\sigma_e\otimes M^\tau_e)(\Lambda)=\frac{1}{n^{|\sigma|+|\tau|}}\sum_{i_1\ldots i_p}\sum_{j_1\ldots j_p}\Lambda_{i_1j_1,i_{\sigma(1)}j_{\tau(1)}}^{e_1}\ldots\ldots\Lambda_{i_pj_p,i_{\sigma(p)}j_{\tau(p)}}^{e_p}$$

By combining the above two formulae, we obtain the following moment formula:
$$(E\circ tr)(\widetilde{W}^{e_1}\ldots\widetilde{W}^{e_p})
=\sum_{\sigma\in S_p}d^{|\sigma|+|\sigma\gamma^{-1}|-1}(mn)^{|\sigma|}(M^\sigma_e\otimes M^\gamma_e)(\Lambda)$$

We use now the standard fact, that we know well from before, that for $\sigma\in S_p$ we have an inequality as follows, with equality precisely when $\sigma\in NC(p)$:
$$|\sigma|+|\sigma\gamma^{-1}|\leq p+1$$

Thus with $d\to\infty$ the sum restricts over the partitions $\sigma\in NC(p)$, and we get:
$$\lim_{d\to\infty}M_e\big(\widetilde{W}\big)=d^p\sum_{\sigma\in NC(p)}(mn)^{|\sigma|}(M^\sigma_e\otimes M^\gamma_e)(\Lambda)$$

Thus, we are led to the conclusion in the statement.
\end{proof}

\section*{16d. Poisson laws}

In order to interpret the asymptotic moment formula that we found in Theorem 16.26, things are quite tricky, because at that level of generality, the resulting measure can vary a lot, generically escaping from the class of measures that we can efficiently study.

\bigskip

Still following \cite{bn2}, our purpose in what follows will be that of identifying the situations where Theorem 16.26 produces compound free Poisson laws, or at least the obvious situations when this happens. Let us introduce the following technical notion:

\index{multiplicative matrix}

\begin{definition}
We call a matrix $\Lambda\in M_n(\mathbb C)\otimes M_n(\mathbb C)$ multiplicative when
$$(M^\sigma_e\otimes M^\gamma_e)(\Lambda)=(M^\sigma_e\otimes M^\sigma_e)(\Lambda)$$
holds for any $p\in\mathbb N$, any exponents $e_1,\ldots,e_p\in\{1,*\}$, and any $\sigma\in NC(p)$.
\end{definition}

Obviously, this is something quite technical, but we will see in a moment where the above condition comes from, the idea being that the combinatorics in the formula from Theorem 16.26 will heavily simplify, and we will be led to compound Poisson laws. Also, in what regards the examples, there are many of them, and we will discuss this later.

\bigskip

Getting back now to the block-modified Wishart matrices, and to the asymptotic moment formula that we found in Theorem 16.26, as just mentioned, the above notion of multiplicative matrix is exactly what we need, for further improving all this. Indeed, we have the following result, substantially building on Theorem 16.26:

\begin{theorem}
Consider a block-modified Wishart matrix 
$$\widetilde{W}=(id\otimes\varphi)W$$
and assume that the matrix $\Lambda\in M_n(\mathbb C)\otimes M_n(\mathbb C)$ associated to $\varphi$ is multiplicative. Then
$$\frac{\widetilde{W}}{d}\sim\pi_{mn\rho}$$
holds, in moments, in the $d\to\infty$ limit, where $\rho=law(\Lambda)$.
\end{theorem}

\begin{proof}
This is something quite tricky, using all the above:

\medskip

(1) Our starting point is the asymptotic moment formula found in Theorem 16.26, for an arbitrary block-modified Wishart matrix, namely:
$$\lim_{d\to\infty}M_e\left(\frac{\widetilde{W}}{d}\right)=\sum_{\sigma\in NC_p}(mn)^{|\sigma|}(M^\sigma_e\otimes M^\gamma_e)(\Lambda)$$

(2) Now since our modification matrix $\Lambda\in M_n(\mathbb C)\otimes M_n(\mathbb C)$ was assumed to be multiplicative, in the sense of Definition 16.27, this formula reads:
$$\lim_{d\to\infty}M_e\left(\frac{\widetilde{W}}{d}\right)=\sum_{\sigma\in NC_p}(mn)^{|\sigma|}(M^\sigma_e\otimes M^\sigma_e)(\Lambda)$$

(3) Our claim now, which will lead right away to the result, is that given an arbitrary square matrix $\Lambda\in M_n(\mathbb C)\otimes M_n(\mathbb C)$, multiplicative or not, having eigenvalue distribution $\rho=law(\Lambda)$, the moments of the corresponding compound free Poisson law $\pi_{mn\rho}$ are given by the following formula, for any choice of an extra parameter $m\in\mathbb N$:
$$M_e(\pi_{mn\rho})=\sum_{\sigma\in NC_p}(mn)^{|\sigma|}(M^\sigma_e\otimes M^\sigma_e)(\Lambda)$$

(4) Indeed, we know from chapter 14 that in the real case, the free cumulants of $\pi_{mn\rho}$ are the moments of $mn\rho$. But, as explained in \cite{nsp}, the same happens in the general complex case, in the sense that the free $*$-cumulants of $\pi_{mn\rho}$ are the $*$-moments of $mn\rho$. We conclude that these free $*$-cumulants are given by the following formula:
\begin{eqnarray*}
\kappa_p^e(\pi_{mn\rho})
&=&M^p_e(mn\rho)\\
&=&mn\cdot M^p_e(\Lambda)\\
&=&mn\cdot (M^\gamma_e\otimes M^\gamma_e)(\Lambda)
\end{eqnarray*}

But with this in hand, by using now Speicher's moment-cumulant formula, that we discussed in chapter 14 in the real case, and which holds in the general complex case too, as explained for instance in \cite{nsp}, we are led to the moment formula in (3).

\medskip

(5) The point now is that with that formula from (3) in hand, our previous asymptotic moment formula for the block-modified Wishart matrix $\widetilde{W}$ simply reads:
$$\lim_{d\to\infty}M_e\left(\frac{\widetilde{W}}{d}\right)=M_e(\pi_{mn\rho})$$

Thus we have indeed $\widetilde{W}/d\sim\pi_{mn\rho}$, in the $d\to\infty$ limit, as stated.
\end{proof}

Let us work out now some explicit consequences of Theorem 16.28, by using the modified easy linear maps from Definition 16.22. We recall from there that any modified easy linear map $\varphi_\pi$ can be viewed as a ``block-modification'' map, as follows:
$$\varphi_\pi:M_{N^s}(\mathbb C)\to M_{N^s}(\mathbb C)$$

In order to verify that the corresponding matrices $\Lambda_\pi$ are multiplicative, we will need to check that all the functions $\varphi(\sigma,\tau)=(M_\sigma^e\otimes M_\tau^e)(\Lambda_\pi)$ have the following property:
$$\varphi(\sigma,\gamma)=\varphi(\sigma,\sigma)$$

For this purpose, we can use the following result, coming from \cite{bn2}:

\begin{proposition}
The following functions $\varphi:NC(p)\times NC(p)\to\mathbb R$ are multiplicative, in the sense that they satisfy the condition $\varphi(\sigma,\gamma)=\varphi(\sigma,\sigma)$:
\begin{enumerate}
\item $\varphi(\sigma,\tau)=|\sigma\tau^{-1}|-|\tau|$.

\item $\varphi(\sigma,\tau)=|\sigma\tau|-|\tau|$.

\item $\varphi(\sigma,\tau)=|\sigma\wedge\tau|-|\tau|$.
\end{enumerate}
\end{proposition}

\begin{proof}
All this is elementary, and can be proved as follows:

\medskip

(1) This follows indeed from the following computation:
$$\varphi_1(\sigma,\gamma)
=|\sigma\gamma^{-1}|-1
=p-|\sigma|
=\varphi_1(\sigma,\sigma)$$

(2) This follows indeed from the following computation:
$$\varphi_2(\sigma,\gamma)
=|\sigma\gamma|-1
=|\sigma^2|-|\sigma|
=\varphi_2(\sigma,\sigma)$$

(3) This follows indeed from the following computation:
$$\varphi_3(\sigma,\gamma)
=|\gamma|-|\gamma|
=0
=|\sigma|-|\sigma|
=\varphi_3(\sigma,\sigma)$$

Thus, we are led to the conclusions in the statement.
\end{proof}

We can get back now to the easy modification maps, and we have:

\begin{proposition}
The partitions $\pi\in P_{even}(2,2)$ are as follows,
$$\pi_1=\begin{bmatrix}\circ&\bullet\\ \circ&\bullet\end{bmatrix}\quad,\quad
\pi_2=\begin{bmatrix}\circ&\bullet\\ \bullet&\circ\end{bmatrix}\quad,\quad
\pi_3=\begin{bmatrix}\circ&\circ\\ \bullet&\bullet\end{bmatrix}\quad,\quad
\pi_4=\begin{bmatrix}\circ&\circ\\ \circ&\circ\end{bmatrix}$$
with the associated linear maps $\varphi_\pi:M_n(\mathbb C)\to M_N(\mathbb C)$ being as follows:
$$\varphi_1(A)=A\quad,\quad
\varphi_2(A)=A^t\quad,\quad
\varphi_3(A)=Tr(A)1\quad,\quad
\varphi_4(A)=A^\delta$$
The corresponding matrices $\Lambda_\pi$ are all multiplicative, in the sense of Definition 16.27.
\end{proposition}

\begin{proof}
In order to prove this, recall from Theorem 16.23 that the associated square matrices, appearing via $\Lambda_{ab,cd}=\varphi(e_{ac})_{bd}$, are given by:
$$\Lambda^1_{ab,cd}=\delta_{ab}\delta_{cd}\quad,\quad 
\Lambda^2_{ab,cd}=\delta_{ad}\delta_{bc}\quad,\quad
\Lambda^3_{ab,cd}=\delta_{ac}\delta_{bd}\quad,\quad
\Lambda^4_{ab,cd}=\delta_{abcd}$$

Since these matrices are all self-adjoint, we can assume that all the exponents are 1 in Definition 16.27, and the multiplicativity condition there becomes:
$$(M_\sigma\otimes M_\gamma)(\Lambda)=(M_\sigma\otimes M_\sigma)(\Lambda)$$

In order to check this condition, observe that for the above 4 matrices, we have:
\begin{eqnarray*}
(M^\sigma\otimes M^\tau)(\Lambda_1)&=&\frac{1}{n^{|\sigma|+|\tau|}}\sum_{i_1\ldots i_p}\delta_{i_{\sigma(1)}i_{\tau(1)}}\ldots\delta_{i_{\sigma(p)}i_{\tau(p)}}=n^{|\sigma\tau^{-1}|-|\sigma|-|\tau|}\\
(M^\sigma\otimes M^\tau)(\Lambda_2)&=&\frac{1}{n^{|\sigma|+|\tau|}}\sum_{i_1\ldots i_p}\delta_{i_1i_{\sigma\tau(1)}}\ldots\delta_{i_pi_{\sigma\tau(p)}}=n^{|\sigma\tau|-|\sigma|-|\tau|}\\
(M^\sigma\otimes M^\tau)(\Lambda_3)&=&\frac{1}{n^{|\sigma|+|\tau|}}\sum_{i_1\ldots i_p}\sum_{j_1\ldots j_p}\delta_{i_1i_{\sigma(1)}}\delta_{j_1j_{\tau(1)}}\ldots\delta_{i_pi_{\sigma(p)}}\delta_{j_pj_{\tau(p)}}=1\\
(M^\sigma\otimes M^\tau)(\Lambda_4)&=&\frac{1}{n^{|\sigma|+|\tau|}}\sum_{i_1\ldots i_p}\delta_{i_1i_{\sigma(1)}i_{\tau(1)}}\ldots\delta_{i_pi_{\sigma(p)}i_{\tau(p)}}=n^{|\sigma\wedge\tau|-|\sigma|-|\tau|}
\end{eqnarray*}

By using now the results in Proposition 16.29, this gives the result.
\end{proof}

Summarizing, the partitions $\pi\in P_{even}(2,2)$ provide us with some concrete input for Theorem 16.28. The point now is that, when using this input, we obtain the main known computations for the block-modified Wishart matrices, from \cite{aub}, \cite{cne}, \cite{mpa}:

\begin{theorem}
The asymptotic distribution results for the block-modified Wishart matrices coming from the partitions $\pi_1,\pi_2,\pi_3,\pi_4\in P_{even}(2,2)$ are as follows:
\begin{enumerate}
\item Marchenko-Pastur: $\frac{1}{d}W\sim\pi_t$, where $t=m/n$.

\item Aubrun type: $\frac{1}{d}(id\otimes t)W\sim\pi_\nu$, with $\nu=\frac{m(n-1)}{2}\,\delta_{-1}+\frac{m(n+1)}{2}\,\delta_1$. 

\item Collins-Nechita one: $n(id\otimes tr(.)1)W\sim\pi_t$, where $t=mn$.

\item Collins-Nechita two: $\frac{1}{d}(id\otimes(.)^\delta)W\sim\pi_m$.
\end{enumerate}
\end{theorem}

\begin{proof}
All these results follow from Theorem 16.28, with the maps $\varphi_1,\varphi_2,\varphi_3,\varphi_4$ in Proposition 16.30 producing the 4 matrices in the statement, modulo some rescalings:

\medskip

(1) Here $\Lambda=\sum_{ac}e_{ac}\otimes e_{ac}$, and so $\Lambda=nP$, where $P$ is the rank one projection on $\sum_ae_a\otimes e_a\in\mathbb C^n\otimes\mathbb C^n$. Thus we have the following formula, which gives the result:
$$\rho=\frac{n^2-1}{n^2}\,\delta_0+\frac{1}{n^2}\,\delta_n$$

(2) Here $\Lambda=\sum_{ac}e_{ac}\otimes e_{ca}$ is the flip operator, $\Lambda(e_c\otimes e_a)=e_a\otimes e_c$. Thus $\rho=\frac{n-1}{2n}\,\delta_{-1}+\frac{n+1}{2n}\,\delta_1$, and so we have the following formula, which gives the result:
$$mn\rho=\frac{m(n-1)}{2}\,\delta_{-1}+\frac{m(n+1)}{2}\,\delta_1$$

(3) Here $\Lambda=\sum_{ab}e_{aa}\otimes e_{bb}$ is the identity matrix, $\Lambda=1$. Thus in this case we have the following formula, which gives $\pi_{mn\rho}=\pi_{mn}$, and so $n\widetilde{W}\sim\pi_{mn}$, as claimed:
$$\rho=\delta_1$$

(4) Here $\Lambda=\sum_ae_{aa}\otimes e_{aa}$ is the orthogonal projection on $span(e_a\otimes e_a)\subset\mathbb C^n\otimes\mathbb C^n$. Thus we have the following formula, which gives the result:
$$\rho=\frac{n-1}{n}\,\delta_0+\frac{1}{n}\,\delta_1$$

Summarizing, we have proved all the assertions in the statement.
\end{proof}

And with this, end of our sequence of theorems. Good work that we did, having Theorem 16.31 stated as such is a good thing. Of course, the linear maps used there are quite trivial, with the corresponding matrices having very simple spectral theory, but there is a beginning for everything, and Theorem 16.31 stands for such a beginning.

\bigskip

Getting now to final conclusions, what we have in Theorems 16.26, 16.28 and 16.31 is certainly quite nice, doing some much needed unification work, and providing a solid basis for a more advanced theory of block-modified Wishart matrices. However, going beyond this, and reaching to a true next-level theory, remains a challenging question:

\bigskip

(1) One problem is to stay with Theorem 16.28 as main result, and look for more general classes of multiplicative matrices, generalizing those used in Theorem 16.31, which are, after all, of quite trivial nature. In practice, this leads to delicate combinatorics.

\bigskip

(2) The other problem is that of staying with Theorem 16.26 as main result, and inventing the correct generalizations of the multiplicative matrices, and of the compound free Poisson laws,   allowing for advances. Again, this leads to delicate combinatorics.

\bigskip

We refer to the papers \cite{aub}, \cite{bn1}, \cite{bn2} and their various follow-ups, including the paper \cite{anv}, making some important advances, for more on all this, and with a mention as well to \cite{bb+} and its various follow-ups, concerned with a closely related problematics.

$$***$$

\bigskip

And with this, end of our discussion regarding the various manipulations that can be done on the Wishart matrices. Of course, this remains one of the many questions that can be of interest, in relation with the random matrices, in general.

\bigskip

Talking now random matrices in general, there are many things to be learned, going in many possible directions, sometimes in relation with free probability, and sometimes not. In fact, the random matrices, first discovered by Wigner in the 50s \cite{wig} are much older than free probability, a theory invented by Voiculescu in the 80s \cite{vo1}, and to be more precise, as of now, moment when I type these lines, twice as old.

\bigskip

Thus, many things to be learned. Of particular interest, and going well beyond what we have been doing in this book, which remains something introductory, is the study of the fluctuations. There is an enormous quantity of work which has gone into this direction, and for getting an idea of how things looked like more recently, 30 years ago, have a look at the paper of Tracy and Widom \cite{twi}, which remains a classic.

\bigskip

Finally, don't forget physics. For some reason, all good random matrix theorists know well physics, and all good theoretical physicists know well random matrices. There are a lot of things to be learned here too, in relation with quantum and statistical mechanics, and more than certainly, lots of interesting new things to be discovered.

\section*{16e. Exercises}

Congratulations for having read this book, and no exercises for this final chapter. Instead, you can have a look at the various books and articles referenced below. Which are a bit biased towards free probability and random matrices, yes I know, but recall from what we learned in this book, starting with page $N=1$, probability theory is always something biased, and you can see that even in the $N\to\infty$ limit.

\baselineskip=14pt

\printindex

\end{document}